\documentclass[b5paper]{book}
\usepackage{geometry}
\usepackage[noinfo,width=176truemm,height=246truemm,center]{crop}

\usepackage{ccicons}
\usepackage{cite}
\usepackage{caption}
\usepackage{xfrac}
\usepackage{xparse}
\usepackage{relsize}
\usepackage{colortbl}
\usepackage{bussproofs}
\usepackage{xcolor}
\usepackage{graphicx}
\usepackage{makeidx}
\usepackage{xstring}
\usepackage{braket}
\usepackage{mathtools}
\usepackage{multicol}
\usepackage{marginnote}
\usepackage{ifthen}
\usepackage{xifthen}
\usepackage{paralist}
\usepackage{calc}
\usepackage{tikz}
\usepackage{tikz-cd}
\usetikzlibrary{calc,arrows,matrix,decorations.pathreplacing}
\usepackage{picins}
\usepackage{amssymb}
\usepackage{amsmath}
\usepackage{mathrsfs}
\usepackage{nicefrac}
\usepackage{stmaryrd}
\usepackage{sectsty}
\usepackage{xr}
\usepackage{tikz}
\usepackage[all]{xypic}
\usepackage{hyperref} % <- should be last because it redefined \ref, etc..

\NewDocumentEnvironment{solution}{g}{%
    \leavevmode\unskip%
    \par\vskip.6em\noindent%
    \textPointHeaderI{Solution to~\sref{#1}.}
    }{\leavevmode\unskip\ignorespaces}

\NewDocumentEnvironment{erratum}{g}{%
    \leavevmode\unskip%
    \par\vskip.6em\noindent%
    \textPointHeaderI{Erratum to~\sref{#1}.}
    }{\leavevmode\unskip\ignorespaces}

\definecolor{pitchblack}{cmyk}{1 1 1 1}
\definecolor{pitchgray}{cmyk}{.5 .5 .5 .5}

\colorlet{darkgreen}{green!50!pitchblack}
\colorlet{darkblue}{blue!75!pitchblack}
\colorlet{lightblue}{blue!50!white}
\colorlet{lightgray}{pitchgray!50!white}
\colorlet{darkgray}{pitchgray!75!pitchblack}
\colorlet{darkred}{red!75!pitchblack}

\allsectionsfont{\sffamily\color{darkblue}}

\newcommand{\textPointHeaderI}[1]{\textcolor{darkblue}{\textbf{\textsf{#1}}}}
\newcommand{\textPointHeaderII}[1]{\textcolor{darkblue}{\textsf{#1}}}
\newcommand{\textPointHeaderIII}[1]{\textcolor{lightgray}{\textsf{#1}}}
\newcommand{\textParsecNumber}[1]{\textcolor{darkblue}{\textbf{\textsf{#1}}}}
\newcommand{\textPointNumberI}[1]{\textcolor{darkblue}{\textsf{#1}}}
\newcommand{\textPointNumberII}[1]{\textcolor{lightblue}{\textsf{#1}}}
\newcommand{\textPointNumberIII}[1]{\textcolor{lightgray}{\textsf{#1}}}
\newcommand{\textDefine}[1]{\textcolor{darkblue}{#1}}
\newcommand{\textSref}[1]{\textsf{#1}}

\newcommand{\Define}[1]{\textDefine{#1}}

\renewcommand{\leq}{\leqslant}
\renewcommand{\geq}{\geqslant}
\renewcommand{\nleq}{\not\leqslant}

\newcommand{\op}[1]{{#1}{^{\mathsf{op}}}}

\newcommand{\ketbra}[2]{\left|#1\right>\!\left<#2\right|}

\newcommand{\TR}{\mathop{\mathrm{tr}}}

\NewDocumentCommand{\vmleq}{o}{\mathrel{%
\IfNoValueTF{#1}{\lesssim}{\lesssim_{#1}}}}

\newcommand{\C}{\mathbb{C}}
\newcommand{\N}{\mathbb{N}}
\newcommand{\Z}{\mathbb{Z}}
\newcommand{\R}{\mathbb{R}}
\newcommand{\I}{\mathbb{I}}

\newcommand{\scrA}{\mathscr{A}}
\newcommand{\scrB}{\mathscr{B}}
\newcommand{\scrC}{\mathscr{C}}

\newcommand{\scrH}{\mathscr{H}}

\newcommand{\scrP}{\mathscr{P}}
\newcommand{\scrK}{\mathscr{K}}

\newcommand{\scrS}{\mathscr{S}}
\newcommand{\scrT}{\mathscr{T}}

\newcommand{\scrV}{\mathscr{V}}
\newcommand{\scrW}{\mathscr{W}}

\newcommand{\ceil}[1]{\left\lceil#1\right\rceil}
\newcommand{\cceil}[1]{\left\lceil\!\!\left\lceil#1\right\rceil\!\!\right\rceil}

\newcommand{\floor}[1]{\left\lfloor#1\right\rfloor}

\newcommand{\Real}[1]{#1_{\R}}
\newcommand{\Imag}[1]{#1_{\I}}

\newcommand\T{\mathrm{t}} % for transpose

\DeclareMathOperator{\Stat}{Stat}
\DeclareMathOperator{\Pred}{Pred}
\DeclareMathOperator{\Pure}{Pure}
\DeclareMathOperator{\SPred}{SPred}
\DeclareMathOperator{\IM}{im}
\newcommand\BOX{%
    {\rlap{\kern .5pt\tikz[baseline]{\draw[line width=0.08ex] (0,0) rectangle (0.6ex,0.6ex)}}%
    \phantom{\diamond}}}
\DeclareMathOperator{\IMperp}{im^\perp}
\DeclareMathOperator{\cok}{cok}
\DeclareMathOperator{\Par}{Par}
\DeclareMathOperator{\Tot}{Tot}
\newcommand{\andthen}[2]{#1 \mathbin{\&} #2}
\newcommand\asrt{\mathrm{asrt}}
\newcommand\sef{\mathrm{seff}}
\newcommand\pproj{\mathord{\vartriangleright}}
\newcommand{\cmpr}[2]{\ensuremath{\{#1|{\kern.2ex}#2\}}}
\newcommand{\id}{\mathrm{id}}

\DeclareMathOperator{\Scal}{Scal}

\DeclareMathOperator{\Inv}{Inv}

\DeclareMathOperator{\supp}{supp}

\DeclareMathOperator{\Nsb}{Nsb}
\newenvironment{spmatrix}{%
    \left(\begin{smallmatrix}}{%
    \end{smallmatrix}\right)}

\newcommand\after{\mathop{\circ}}
\newcommand\hafter{\mathop{\hat\circ}}
\newcommand\hotimes{\mathbin{\hat{\otimes}}}
\newcommand\ad{\mathrm{ad}}
\newcommand{\uwlim}{\qopname\relax m{uwlim}}
\newcommand{\unlim}{\qopname\relax m{unlim}}
\newcommand{\uslim}{\qopname\relax m{uslim}}

\newcommand\ncp{\mathrm{ncp}}
\newcommand\pto{\rightharpoonup}
\newcommand\pfrom{\leftharpoonup}
\newcommand{\pullback}[1][dr]{\save*!/#1-1.2pc/#1:(-1,1)@^{|-}\restore}
\newcommand\vN{\mathsf{vN}}
\newcommand\CvN{\mathsf{CvN}}
\newcommand\OUG{\mathsf{OUG}}
\newcommand\EJA{\mathsf{EJA}}
\newcommand\OUS{\mathsf{OUS}}
\newcommand\SET{\mathsf{Set}}
\newcommand\CRng{\mathsf{CRng}}
\newcommand\Rng{\mathsf{Rng}}
\newcommand\bCH{\mathsf{CH}}
\newcommand\AConvM{\mathsf{AConv}_M}
\newcommand\EMod{\mathsf{EMod}}
\newcommand{\bigovee}{\mathop{\vphantom{\sum}\mathchoice%
        {\vcenter{\hbox{\huge $\ovee$}}}%
        {\vcenter{\hbox{\Large $\ovee$}}}%
        {\ovee}{\ovee}}\displaylimits}
\DeclareMathOperator{\Kl}{\mathcal{K}\!\ell}
\newcommand\rfrac[2]{\reflectbox{\nicefrac{\reflectbox{\ensuremath{#1}}}{\reflectbox{\ensuremath{#2}}}}}

\makeatletter
\newcommand{\bigperp}{%
  \mathop{\mathpalette\bigp@rp\relax}%
  \displaylimits
}

\newcommand{\bigp@rp}[2]{%
  \vcenter{
    \m@th\hbox{\scalebox{\ifx#1\displaystyle2.1\else1.5\fi}{$#1\perp$}}
  }%
}
\makeatother

\makeatletter
\def\ourrawref#1{%
    \expandafter\expandafter\expandafter
    \@car\csname r@#1\endcsname\@nil
}
\makeatother

\newcommand\numberthis{\addtocounter{equation}{1}\tag{\theequation}}

\newcommand{\qed}{\hfill\textcolor{darkblue}{\ensuremath{\square}}}

\NewDocumentCommand{\sref}{m}{\textSref{\ref{#1}}}

\newcommand{\doi}[1]{\href{https://doi.org/#1}{doi:#1}}

\newcommand{\purl}[1]{\href{#1}{#1}}

\newcommand{\spacingfix}[0]{\vspace{-1\baselineskip}\ignorespaces}
\newcommand{\IGNORE}[1]{}

\RenewDocumentCommand{\sref}{m}{\textSref{%
	\ifthenelse{%
		\equal{\value{parsec}}{\ourrawref{#1::parsec}}%
	}{%
		\ref{#1::point}%
	}{%
        \ref{#1}%
	}%
}}

\newcommand{\colofon}[5]{%
\vspace*{4em}
\noindent
\begin{center}
    #5
    \vspace{3em}

    \textcolor{darkblue}{\textsf{\textbf{ Identifiers}}}
    \vspace{.5em}

\href{https://hdl.handle.net/#4}{hdl:~#4}
    \vspace{.5em}

\href{https://arxiv.org/abs/#3}{arXiv:~#3}
    \vspace{.5em}

\textsc{isbn:}~#2
\end{center}
\vspace{3em}
\begin{center}
    \textcolor{darkblue}{\textsf{\textbf{ Persistent links}}}

\vspace{.5em}
\purl{https://arxiv.org/abs/#3}

\vspace{.5em}
\purl{https://doi.org/#4}

\vspace{.5em}
\purl{https://hdl.handle.net/#4}
\end{center}
\vspace{3em}

\begin{center}
    \textcolor{darkblue}{\textsf{\textbf{Source code}}}
\vspace{.5em}

\begin{tabular}{rl}
\LaTeX{} & \purl{https://github.com/westerbaan/theses} \\
cover & \purl{https://github.com/westerbaan/ndpt}
\end{tabular}

\end{center}
\vfill{}

\noindent
Printed by GVO drukkers \& vormgevers B.V., Ede, 
\purl{https://proefschriften.nl}.

\vspace{1em}
\noindent
Where applicable,
\ccCopy{}~2019 #1,
\ccLogo{}\,\ccAttribution{}  available under \textsc{cc by}, \cite{ccby40}.
\newpage
}

\newcounter{tmptmp}

\newcounter{parsec} % keeps track of the current parsec number
\newcounter{parsecMajor}
\newcounter{parsecMinor} % parsec = 10 * parsecMajor + parsecMinor

\newcommand\refsforparsec{%
    \setcounter{parsecMajor}{\value{parsec}/10}%
    \setcounter{parsecMinor}{\value{parsec}-10*\value{parsecMajor}}%
    \setcounter{parsec}{\value{parsec}-1}%
    \refstepcounter{parsec}%
}

\NewDocumentEnvironment{parsec}{g o}{%
	\leavevmode\unskip%
	\par\vskip1em\noindent%
    \renewcommand{\theparsec}{\the\value{parsecMajor}\alph{parsecMinor}}%
    \setcounter{parsec}{#1}%
    \refsforparsec{}%
	\setcounter{tmptmp}{2*\value{parsec}}%
	\markboth{\the\value{tmptmp}}{\the\value{tmptmp}}%
	\IfValueT{#2}{\label{#2}}%
    \renewcommand{\thepoint}{\the\value{parsecMajor}\alph{parsecMinor}}%
    \setcounter{point}{\value{point}-1}%
    \refstepcounter{point}% 
    \label{parsec-\the\value{parsec}}%
        \marginnote{\makebox[3em][c]{\textParsecNumber{%
            \the\value{parsecMajor}%
            \alph{parsecMinor}%
            }}}%
	\ignorespaces%
}{%
	\leavevmode\unskip%
	\setcounter{tmptmp}{2*\value{parsec}+1}%
	\markboth{\the\value{tmptmp}}{\the\value{tmptmp}}%
	\ignorespaces%
}

\newcounter{point} % keeps track of the current point
\newcounter{pointMajor}
\newcounter{pointMinor}
\numberwithin{point}{parsec}
\newcounter{pointdepth}  % keeps track of the depth of the current point
\newcommand\refsforpoint{%
    \setcounter{pointMajor}{\value{point}/10}%
    \setcounter{pointMinor}{\value{point}-10*\value{pointMajor}}%
    \setcounter{point}{\value{point}-1}%
    \refstepcounter{point}%
}

\NewDocumentEnvironment{point}{g o g}{%
	\leavevmode\unskip%
	\setcounter{pointdepth}{\value{pointdepth}+1}%
	\refstepcounter{point}% TODO: <- is this needed?
    \setcounter{point}{#1}%
    \IfValueT{#2}{%
        \renewcommand{\thepoint}{\the\value{parsec}}%
        \refsforpoint{}%
        \label{#2::parsec}%
        \renewcommand{\thepoint}{\Roman{pointMajor}\alph{pointMinor}}%
        \refsforpoint{}%
        \label{#2::point}%
        \renewcommand{\thepoint}{\theparsec\,\Roman{pointMajor}\alph{pointMinor}}%
        \refsforpoint{}%
        \label{#2}%
    }%
    \renewcommand{\thepoint}{\theparsec\,\Roman{pointMajor}\alph{pointMinor}}%
    \refsforpoint{}%
    \label{parsec-\the\value{parsec}.\the\value{point}}%
	\ifthenelse{\equal{\value{point}}{10}}{}{%
		\ifthenelse{\equal{\value{pointdepth}}{1}}{%
			\par\penalty-100\vskip.6em\noindent%
		}{%
			\ifthenelse{\equal{\value{pointdepth}}{2}}{%
				\par\penalty-50\vskip.2em\noindent%
			}{%
				\par\penalty-25\vskip.1em\noindent%
			}%
		}%
		\marginnote{\makebox[2em][c]{\small%
                    \ifthenelse{\equal{\value{pointdepth}}{1}}{%
                        \textPointNumberI{\Roman{pointMajor}\alph{pointMinor}}%
                    }{%
                        \ifthenelse{\equal{\value{pointdepth}}{2}}{%
                            \textPointNumberII{\Roman{pointMajor}\alph{pointMinor}}%
                        }{%
                            \textPointNumberIII{\Roman{pointMajor}\alph{pointMinor}}%
                        }%
		}}}%
	}%
	\IfValueT{#3}{%
		\ifthenelse{\equal{\value{pointdepth}}{1}}{%
			\textPointHeaderI{#3}%
		}{%
			\ifthenelse{\equal{\value{pointdepth}}{2}}{%
				\textPointHeaderII{#3}%
			}{%
				\textPointHeaderIII{(#3)}%
			}}%
	\ \ }%	
\ignorespaces%
}{%
\leavevmode\unskip%
\setcounter{pointdepth}{\value{pointdepth}-1}%
\ignorespaces%
}

\usepackage{fancyhdr}
\renewcommand{\chaptermark}[1]{}  
\renewcommand{\sectionmark}[1]{}

\newcounter{firstParsec}
\newcounter{lastParsec}
\newcounter{firstParsecF}
\newcounter{lastParsecF}
\newcounter{firstParsecMajor}
\newcounter{firstParsecMinor}
\newcounter{lastParsecMajor}
\newcounter{lastParsecMinor}

\newcounter{parsecToBeContinued}  
\newcommand\ourfancyfooters{%
    \fancyfoot[CE]{%
        \IfInteger{\rightmark}{%
            \setcounter{tmptmp}{\leftmark+0}%
            \ifthenelse{\equal{\value{tmptmp}}{\value{lastParsecF}}}{%
                \setcounter{firstParsecF}{0}%
            }{%
                \setcounter{firstParsecF}{\rightmark+0}%
            }
        }{%
        }%
    }%
    \fancyfoot[RO]{%
        \textsf{\footnotesize\textcolor{lightgray}{\thepage}}}%
    \fancyfoot[CO]{%
        \ifthenelse{\equal{\value{firstParsecF}}{0}}{%
            \setcounter{firstParsecF}{\rightmark+0}%
        }{% 
        }
        \setcounter{firstParsec}{\value{firstParsecF}/2}%
        \setcounter{lastParsecF}{\leftmark+0}%
        \ifthenelse{\equal{\value{lastParsecF}}{0}}{%
        }{%
            \setcounter{lastParsec}{\value{lastParsecF}/2}%
            \textPointNumberI{%
                \ifthenelse{\equal{\value{parsecToBeContinued}}{1}}{..}{}%
                \setcounter{firstParsecMajor}{\value{firstParsec}/10}%
                \setcounter{firstParsecMinor}{\value{firstParsec}-\value{firstParsecMajor}*10}%
                \the\value{firstParsecMajor}\alph{firstParsecMinor}%
                \ifthenelse{\equal{\value{firstParsec}}{\value{lastParsec}}}{%
                }{%  no footer without parsecs
                    \setcounter{tmptmp}{\value{firstParsec}+1}%
                    \ifthenelse{\equal{\value{tmptmp}}{\value{lastParsec}}}{, }{--}% TODO: fix this, maybe?
                    \setcounter{lastParsecMajor}{\value{lastParsec}/10}%
                    \setcounter{lastParsecMinor}{\value{lastParsec}-\value{lastParsecMajor}*10}%
                    \the\value{lastParsecMajor}\alph{lastParsecMinor}%
                }%
                \setcounter{tmptmp}{\value{lastParsec}*2}%
                \ifthenelse{\equal{\value{tmptmp}}{\value{lastParsecF}}}{%
                    ..\setcounter{parsecToBeContinued}{1}%
                }{%
                    \setcounter{parsecToBeContinued}{0}%
                }%
            }%
        }%
    }%
}

\ourfancyfooters
\fancypagestyle{plain}{\ourfancyfooters}

\makeatletter
\def\@wrindex#1{%
   \protected@write\@indexfile{}%
   {\string\indexentry{#1|parsechyperlink{\the\value{parsec}.\the\value{point}}}{\the\value{parsec}.\the\value{point}}}%
 \endgroup
 \@esphack}%
\makeatother

\let\oldchapter\chapter%
\renewcommand\chapter[1]{%
    \oldchapter{#1}%
    \setcounter{tmptmp}{(\value{parsec}/10)*10+10}
    \addtocontents{parsectoc}{\protect\contentsline{chapter}{\numberline{\thechapter} #1}{\sref{parsec-\the\value{tmptmp}}}{chapter.\thechapter}}%
}
\let\oldsection\section%
\renewcommand\section[1]{%
    \oldsection{#1}%
    \setcounter{tmptmp}{(\value{parsec}/10)*10+10}
    \addtocontents{parsectoc}{\protect\contentsline{section}{\numberline{\thesection} #1}{\sref{parsec-\the\value{tmptmp}}}{section.\thesection}}%
}
\let\oldsubsection\subsection%
\renewcommand\subsection[1]{%
    \oldsubsection{#1}%
    \setcounter{tmptmp}{(\value{parsec}/10)*10+10}
    \addtocontents{parsectoc}{\protect\contentsline{subsection}{\numberline{\thesubsection} #1}{\sref{parsec-\the\value{tmptmp}}}{subsection.\thesubsection}}%
}
\let\oldsubsubsection\subsubsection%
\renewcommand\subsubsection[1]{%
    \oldsubsubsection{#1}%
    \setcounter{tmptmp}{(\value{parsec}/10)*10+10}
    \addtocontents{parsectoc}{\protect\contentsline{subsubsection}{\numberline{\thesubsubsection} #1}{\sref{parsec-\the\value{tmptmp}}}{subsubsection.\thesubsubsection}}%
}

\newcommand\backmattertitle[1]{{%
    \sffamily\color{darkblue}\Huge\bfseries #1\vspace{1em}}}



\usepackage{subfiles}

\makeindex

\begin{document}

% \maketitle
\newcommand{\titelpagina}[8]{%
\vspace*{4em}
\begin{center}
{\sffamily\color{darkblue}\huge Dagger and Dilation\\
    \vskip.4em \large in the Category of Von Neumann Algebras}
\end{center}
\vspace{2em}
\begin{center}
{\par\noindent\sffamily\color{darkblue}\large\textbf{#1}}
\end{center}
\begin{center}
#2
\end{center}
\vspace{1em}
\begin{center}
#3
\end{center}
\begin{center}
#4
\end{center}
\begin{center}
#5
\end{center}
\begin{center}
#6
\end{center}
\vspace{10em}
\begin{center}%
#7
\end{center}%
\begin{center}%
\large\color{darkblue}Bastiaan Eelco \textsc{Westerbaan}
\end{center}%
\begin{center}%
#8
\end{center}
\newpage
}

\newcommand{\achterkanttitelpagina}[3]{
\vspace*{8em}
\begin{center}
\textbf{\large\sffamily\color{darkblue}#1:}
\end{center}
\begin{center}
Prof.~dr.~B.P.F.~\textsc{Jacobs}
\end{center}
\vspace{5em}
\begin{center}
\textbf{\large\sffamily\color{darkblue}#2:}
\end{center}
\begin{center}
Prof.~dr.~J.D.M.~\textsc{Maassen}\\
\vspace{1em}
    Prof.~dr.~P.~\textsc{Panangaden} \\ {\footnotesize(McGill University, Canada)} \\
\vspace{1em}
    Prof.~dr.~P.~\textsc{Selinger} \\ {\footnotesize(Dalhousie University, Canada)} \\
\vspace{1em}
    Dr.~C.J.M.~\textsc{Heunen} \\ {\footnotesize(University of Edinburgh, #3)} \\
\vspace{1em}
    Dr.~S.~\textsc{Staton} \\ {\footnotesize(University of Oxford, #3)}
\end{center}
\newpage
}

\newcommand{\whiteout}[1]{{\color{white}#1}}

\titelpagina{\whiteout{P}$\,$}{$\,$\whiteout{ter}\\
\whiteout{aan}$\,$\\
\whiteout{op}$\,$\\
\whiteout{volgens}$\,$\\
\whiteout{in}$\,$}{\whiteout{op}$\,$}{\whiteout{dinsdag}$\,$}{\whiteout{om}$\,$}{\whiteout{10.30}$\,$}{\whiteout{door}}{$\,$\whiteout{geboren}$\,$\\
\whiteout{te}$\,$}

\colofon{B.E.~Westerbaan}{978-94-6332-485-4}{1803.01911}{2066/201785}{Third corrected version\\(april 26th, 2019)}

\titelpagina{Proefschrift}{ter verkrijging van de graad van doctor\\
aan de Radboud Universiteit Nijmegen\\
op gezag van de rector magnificus prof. dr. J.H.J.M. \textsc{van Krieken},\\
volgens besluit van het college van decanen \\
in het openbaar te verdedigen}{op}{\textbf{dinsdag 14 mei 2019}}{om}{\textbf{12.30 uur precies}}{door}{geboren op 30 augustus 1988\\
te Nijmegen}

\achterkanttitelpagina{Promotor}{Manuscriptcommissie}{Verenigd Koninkrijk}

\titelpagina{Doctoral Thesis}{to obtain the degree of doctor \\
from  Radboud University Nijmegen \\
on the authority of the Rector Magnificus prof. dr. J.H.J.M. \textsc{van Krieken},\\
according to the decision of the Council of Deans \\
to be defended in public}{on}{\textbf{Tuesday, May 14, 2019}}{at}{\textbf{12.30 hours}}{by}{born on August 30, 1988\\
in Nijmegen (the Netherlands)}

\achterkanttitelpagina{Supervisor}{Doctoral Thesis Committee}{United Kingdom}

\vspace*{.1em}

\newpage

\makeatletter\@starttoc{parsectoc}\makeatother

\chapter{Introduction}
\begin{parsec}{1340}[bas-first-parsec]%
\begin{point}{10}%
This thesis is a mathematical study of quantum computing, concentrating
    on two related, but independent topics.
First up are \emph{dilations}, covered in chapter~\ref{chapter1}.
    In chapter \ref{chapter2}
    ``diamond, andthen, dagger''
    we turn to the second topic: \emph{effectus theory}.
Both chapters, or rather parts, can be read separately
    and feature a comprehensive introduction of their own.
For now, we focus on what they have in common.
\begin{point}{20}%
Both parts study the category of von Neumann algebras
    (with ncpsu-maps\footnote{%
        An abbreviation for \Define{n}ormal \Define{c}ompletely \Define{p}ositive contractive (i.e.~\Define{s}ub\Define{u}nital) linear maps.} between them).
A von Neumann algebra is a model for a (possibly infinite-dimensional)
    quantum datatype, for instance:
\begin{center}
    \begin{tabular}{lllll}
        $M_2$ & a qubit &\quad\qquad& $M_3 \otimes M_3$ & an (ordered) pair of qutrits \\
        $M_3$ & a qutrit && $\C^2$ & a bit \\
        $M_2 \oplus M_3$ & a qubit or a qutrit &
        & $\scrB(\ell^2)$ & quantum integers
\end{tabular}
\end{center}
The ncpsu-maps between these algebras represent physical
    (but not necessarily computable) operations
        in the opposite direction
        (as in the `Heisenberg picture'):
\begin{enumerate}
\item
    \emph{Measure a qubit in the computational basis} \quad
            ($\mathrm{qubit} \to \mathrm{bit}$)\\
        $\C_2 \rightarrow M_2,$ \quad 
        $(\lambda,\mu) \ \mapsto\ \lambda \ketbra{0}{0} + \mu \ketbra{1}{1}$,
        \quad using braket notation\footnote{%
            See e.g.~\cite[fig.~2.1]{nielsen2002quantum},
            so here
            $\ket{0}\equiv\begin{spmatrix}1 \\ 0 \end{spmatrix}$,
            $\bra{0}\equiv\begin{spmatrix}1 & 0 \end{spmatrix}$,
            $\ket{1}\equiv\begin{spmatrix}0 \\ 1 \end{spmatrix}$ and
                $\ketbra{0}{0} = \begin{spmatrix} 1 & 0
                \\ 0 & 0 \end{spmatrix}$.}
\item
    \emph{Discard the second qubit} \quad
            ($\mathrm{qubit} \times \mathrm{qubit}\to \mathrm{qubit}$) \\
        $M_2 \rightarrow M_2 \otimes M_2$, \quad
        $a \ \mapsto \ a \otimes 1$
\item
    \emph{Apply a Hadamard-gate to a qubit} \quad
            ($\mathrm{qubit} \to \mathrm{qubit}$) \\
        $M_2 \rightarrow M_2$, \quad
        $a \ \mapsto \ H^*a H $, \quad where~$H = \frac{1}{\sqrt{2}}\begin{spmatrix}
            1 & 1 \\ 1 & -1
        \end{spmatrix}$
\item
    \emph{Initialize a qutrit as 0} \quad
            ($\mathrm{empty} \to \mathrm{qutrit}$) \\
        $M_3 \rightarrow \C$, \quad
        $a \ \mapsto \ \bra{0}a\ket{0}$
\end{enumerate}
\end{point}
\spacingfix{}
\begin{point}{30}%
The first part of this thesis starts with a concrete question:
    does the famous  Stinespring dilation theorem~\cite{stinespring},
    a cornerstone of quantum theory,
    extend to every ncpsu-map between von Neumann algebras?
Roughly, this means that every ncpsu-map  splits
    in two maps of a particularly nice form.
We will see that indeed it does.
The insights attained pondering this question lead naturally,
        as they so often do in mathematics, to other new results:
    for instance, we learn more about  the tensor product of
    dilations (\sref{paschke-tensor})
        and the pure maps associated to them (\sref{paschke-pure}).
But perhaps the mathematically
    most interesting by-catch of pondering the question
    are the contributions (see \sref{overview-dils}) to the theory
    of self-dual Hilbert C$^*$-modules over a von Neumann algebra,
    which underlie these dilations.
\end{point}
\begin{point}{40}%
In contrast to the first part of this thesis,
    the second part starts
    with the abstract (or rather: vague) question:
    what characterizes the category of von Neumann algebras?
The goal is to find a set of reasonable axioms
    which pick out the category of von Neumann algebras
    among all other categories.
In other words: a reconstruction of quantum theory with categorical axioms
    (presupposing von Neumann algebras are the right model of quantum theory).
    Unfortunately, we do not reach this goal here.\footnote{%
        Building on this thesis,
            van de Wetering recently found a reconstruction
            of finite-dimensional quantum theory~\cite{wetering,weteringeffthe}.}
However, by studying von Neumann algebras
    with this self-imposed restriction,
    one is led to, as again is often the case in mathematics,
    upon several new notions and results,
    which are of independent interest,
    such as, for instance,~$\diamond$-adjointness,
    pure maps and the existence of their dagger.
\end{point}
\end{point}
\begin{point}{50}{Twin theses}%
A substantial part of the research presented here
    has been performed in close collaboration with
    my twin brother Abraham Westerbaan.
Our results are interwoven to such a degree
    that separating them is difficult
    and, more importantly, detrimental to the clarity of presentation.
Therefore we decided to write our theses as two consecutive volumes
    (with consecutive numbering, more on that later).
    Abraham's thesis~\cite{bram} is the first volume
    and covers the preliminaries on von Neumann
    algebras and, among other results, an axiomatization of
        the sequential product~$b\mapsto \sqrt{a}b\sqrt{a}$.
The present thesis is the second volume and builds upon
    the results of Abraham's thesis,
            in its study of dilations and effectus theory.
Because of this clear division,
    a significant amount of Abraham's original work appears in my thesis
        and, vice versa, my work appears in his.
However, the large majority of the
        new results in each of our theses is our own.
\end{point}
\begin{point}{60}{Preliminaries}%
For the first part of this thesis,
    only knowledge is presumed of, what is considered to be,
    the basic theory of von Neumann algebras.
Unfortunately, this `basic theory' takes most introductory texts
    close to a thousand pages to develop.
This situation has led my twin brother to ambitiously
     redevelop the basic theory of von Neumann algebras
     tailored to the needs of our theses (and hopefully others
     in the field).
In his thesis, Abraham presents~\cite{bram}
    a clear, comprehensive and self-contained development of the basic
    theory of von Neumann algebras in less than 200 pages.
For a traditional treatment,
    I can wholeheartedly recommend the books of Kadison and Ringrose \cite{kr}.
For the second part of this thesis,
    only elementary category theory is required
    as covered in any introductory text, like e.g.~\cite{awodey}.
\end{point}
\begin{point}{70}[bintr-example]{Writing style}%
This thesis is divided into points,
    each a few paragraphs long.
Points are numbered with Roman numerals in the margin:
    this is point~\sref{bintr-example}.
Points  are organized into groups,
    which are assigned Arabic numerals,
    that appear at the start of the group in the margin.
The group of this point is~\sref{bas-first-parsec}.
Group numbers below~\sref{bas-first-parsec}
    appear in my twin brother's thesis~\cite{bram}.
Page numbers have been replaced by the
    group numbers that appear on the spread.\footnote{%
        A spread are two pages that are visible at the same time.}
This unusual system allows for
    very short references:
    instead of writing `see Theorem 2.3.4 on page 123',
    we simply write `see \sref{existence-paschke}'.
Even better: `see the middle of the proof of Theorem 2.3.4 on page 123'
    becomes the more precise `see \sref{paschke-uniqueness}'.

Somewhat unconventionally, this thesis includes exercises.
Most exercises replace straight-forward, but repetitive, proofs.
Other exercises, marked with an asterisk *,
    are harder and explore tangents.
In both cases, these exercises are not meant to force the reader into the
    role of a pupil, but rather to aid her, by reducing clutter.
    Solutions to these exercises (page \pageref{sols})
    and errata to the printed version of this text (page \pageref{errata}),
    can be found at the end of this PDF.
Possibly further errata will be made available
    on the arXiv under~\href{https://arxiv.org/abs/1803.01911}{1803.01911}.
\end{point}
\begin{point}{80}{Advertisements}%
This thesis is based on the publications~\cite{wwpaschke,westerbaan2016universal,cho2015quotient,statesofconvexsets}.
During my doctoral research,
    I also published the following papers,
    which although somewhat related,
    would've made this thesis too broad in scope if included.
\begin{enumerate}
    \item In \emph{A Kochen--Specker System has at least 22 Vectors} \cite{uijlen2016kochen},
        Uijlen and I justify the title building upon the work
        of Arends, Ouaknine and Wampler.
\item 
    In \emph{Unordered Tuples in Quantum Computation} \cite{bags}
        Furber and I argue that~$M_3 \oplus \C$ is the right algebra
        modeling an unordered pair of qubits.  Then we compute
        the algebra for an unordered~$n$-tuple of~$d$-level systems.
\item
    In \emph{Picture-perfect QKD} \cite{kissinger2017picture}
        Kissinger, Tull and I show how to graphically derive a security
        bound on quantum key distribution using the continuity of Stinespring's dilation.
\item
    In \cite{schwabe2016solving}, Schwabe and I
            work out the details of a quantum circuit
            to break binary~$\mathcal{MQ}$ using Grover's algorithm.
\end{enumerate}
In the second part of this thesis,
    we
 use the notion of effectus (introduced by Jacobs \cite{newdirections}),
    solely for studying quantum theory.
Effectus theory, however, has more breadth to it,
 see e.g.~\cite{jacobs2017quantum,
cho2017disintegration,
adams2015type,
jacobs2016hyper,
jacobs2017channel,
jacobs2017formal,
cho2017efprob,
jacobs2017probability,
jacobs2017recipe,
jacobs2016effectuses,
jacobs2016affine,
jacobs2017distances,
jacobs2015effect}.
For a long, but comprehensive introduction, see \cite{effintro}.
\begin{point}{90}%
Next, I would like to mention some other recent approaches
    to the mathematical study of quantum theory
        which are related in spirit and of which I am aware.
First, there are colleagues who also study particular categories
    related to quantum theory.
Spearheaded by Coecke and Abramsky \cite{abramsky2004categorical},
    dozens of authors around Oxford
    have developed insightful graphical calculi
    by studying the category of Hilbert spaces
    categorically \cite{coecke2017picturing}.
Not explicitly categorically (but certainly in spirit)
    is the work around operational probabilistic theories,
    see e.g.~\cite{DAriano2016}.
Besides these approaches (and that of effectus theory),
    there are many other varied categorical forages
    into quantum theory, e.g.~\cite{kornell2012,rennela2017infinite,staton,furber2013kleisli}.
Then there are authors who pick out different structures from quantum theory
    to study.
Most impressive is the work of Alfsen and Shultz,
    who studied and axiomatized the state
    spaces of many operator algebras, including von Neumann algebras~\cite{alfsen2012}.
Recently, the poset of commutative subalgebras of a fixed operator algebra
    has received attention,
    see e.g.~\cite{heunen2015domains,bert,heunen2012bohrification}.
Finally,
    tracing back to von Neumann himself,
    the lattice of projections (representing sharp predicates)
    and the poset of effects (fuzzy predicates),
    has have been studied by many, see e.g.~\cite{dvurecenskij2013new}.
\end{point}
\end{point}
\begin{point}{91}{About the cover}
One of the most familiar structures
    in operator algebras (the mathematical field of this thesis)
    is perhaps the infinite-dimensional
    Hilbert space~$\ell^2(\N)$.
What would its unit sphere look like?
The cover depicts an approximation:
    a 12-dimensional spherical shell split into two.
On the outside the pieces are completely reflective ---
    on the inside they are dark blue and only mostly reflective.
The pieces are sat between two (11-dimensional) hyperplanes
    with a somewhat reflective checkerboard pattern.
    The planes are cropped at a certain radius from the origin
        turning them into \emph{hyperdisks}.
The picture is rendered using a homebrew\footnote{%
    Source code is available at~\url{https://github.com/westerbaan/ndpt}.}
    \emph{path tracer},
    which simulates how rays of light would bounce to end up
    in the camera.
The calculations involved (such as the one to determine how a ray bounces
    from an~$n$-dimensional sphere)
    are all conveniently expressed and performed using
    inner products.
The primary rays cast by the camera have the biggest presence
    in the first three dimensions, only ten percent in the fourth,
    five percent in the fifth et cetera.
This causes more higher-dimensional oddities to appear in the image
    where rays have travelled further or bounced farther in the
    higher dimensions.
\end{point}
\begin{point}{100}{Funding}
My research has been supported by the
European Research Council under grant agreement \textnumero~320571.
\end{point}
\begin{point}{110}{Acknowledgments}%
Progress in theoretical fields
    is often dramatized
    as the fruit of personal struggles of  \emph{einzelg\"angers}
    grown in deep focus and isolation.
But what is missing in this depiction is that
    often the seeds were sown
    in casual conversation among colleagues.
I would like to acknowledge two seeds I'm aware of.
First, it was Sam Staton who suggested looking at a universal property
    for instruments, which led to quotients and comprehension.
Secondly, Chris Heunen insisted,
    against my initial skepticism,
    that Stinespring's dilation must have a
    universal property.
I've also had the immense pleasure of discussing research with
    Robin Adams,
    Guillaume Allais,
    Henk Barendregt,
    Kenta Cho,
    Joan Daemen,
    Robert Furber,
    Bart Jacobs,
    Bas Joosten,
    Martti Karvonen,
    Aleks Kissinger,
    Mike Koeman,
    Hans Maassen,
    Joshua Moerman,
    Teun van Nuland,
    Dusko Pavlovic,
    Arnoud van Rooij,
    Frank Roumen,
    Marc Schoolderman,
    Peter Schwabe,
    Peter Selinger,
    Alexandra Silva,
    Michael Skeide,
    Pawel Sobocinski,
    Sean Tull,
    Sander Uijlen,
    Marloes Venema,
    Abraham Westerbaan,
    John van de Wetering and
    Fabio Zanasi.
I'm particularly grateful to Arnoud, Abraham and John
    for proofreading this thesis.
I'm honored to have been received by Dusko
    for a two-month research visit,
    on which he, Jason, Jon, Mariam, Liang-Ting, Marine, Muzamil, Pawel
    and Peter-Michael made me feel at home on the other side of the planet.
On this side of the planet,
    it is hard to imagine a friendlier group of colleagues.
I'm much obliged to Marko van Eekelen,
    Wouter Geraedts,
    Engelbert Hubbers,
    Carlo Meijer and
    Joeri de Ruiter
    with whom I've performed security audits as a side-gig,
    which allowed for an extension of my position.
I'm grateful to my supervisor, Bart Jacobs,
    who showed me a realistic optimism,
    which is in such short supply these days.
I'm touched by the patience of my friends, (extended) family and Annelies,
    during the final stretch of writing this thesis.

Last, but not least, I would like to share my deep appreciation
    for the (under)graduate education
    I've been fortunate enough to receive at the Radboud Universiteit.
Only now, a few years later,
    I can see how rare it is,
    not only for its uncompromised treatment of the topics at hand,
    but also for the personal attention to each student.
Of all the great professors,
Henk Barendregt,
    Mai Gehrke, Arnoud van Rooij and Wim Veldman
    have been of particular personal influence to me --- I owe them a lot.
\end{point}
\end{parsec}

% vim: ft=tex.latex

\chapter{Dilations}\label{chapter1}

\begin{parsec}{1350}[dils-intro]%
\begin{point}{10}%
In this chapter we will study dilations.
The common theme is that a complicated map
    is actually the composition of a simpler map
    after a representation into a larger algebra.
    The famous Gel'fand--Naimark--Segal construction (see \sref{gns}
        in \cite{bram})
    is a dilation theorem in disguise:
    it shows that every state is a vector state on a larger algebra:
\end{point}
\begin{point}{20}[dils-gns]{Theorem (GNS')}%
    For each pu-map\footnote{%
        pu-map is an abbreviation for positive unital linear map.
        See~\sref{maps} and~\sref{p-uwcont}
            for other abbreviations, like nmiu-map
            and ncpsu-map.
        }~$\omega\colon \scrA \to \C$
        \index{map!ncp(u)-}
        \index{map!nmiu-}
        from a~C$^*$-algebra~$\scrA$,
    there is a Hilbert space~$\scrH$,
    an miu-map~$\varrho\colon \scrA \to \scrB(\scrH)$
    and a vector~$x \in \scrH$
    such that~$\omega = h \after \varrho$
    where~$h \colon \scrB(\scrH) \to \C$
    is given by~$h(T) \equiv \left<x,Tx\right>$.
\end{point}
\begin{point}{30}%
Probably the most famous dilation theorem is that
of Stinespring~\cite[thm.~1]{stinespring}.
(We will see a detailed proof later, in \sref{dils-proof-stinespring}.)
\end{point}
\begin{point}{40}[stinespring-theorem]{Theorem (Stinespring)}
    For every cp-map~$\varphi\colon \scrA \to \scrB (\scrH)$
        between a C$^*$-algebra~$\scrA$
            and the C$^*$-algebra of
            bounded operators on a Hilbert space~$\scrH$,
    there is a Hilbert space~$\scrK$,
        an miu-map~$\varrho\colon \scrA \to \scrB(\scrK)$
        and a bounded operator~$V\colon \scrH \to \scrK$
        such that~$\varphi = \ad_V \after \varrho$,
        where~$\Define{\ad_V}\index{ad@$\ad_V$!Hilbert spaces} \colon \scrB(\scrK) \to \scrB(\scrH)$
        is given by~$\ad_V(T)\equiv V^*TV$.
Furthermore
\begin{enumerate}
\item
    if~$\varphi$ is unital, then~$V$ is an isometry \emph{and}
\item
    if~$\scrA$ is a von Neumann algebra
    and~$\varphi$ is normal, then~$\varrho$ is normal as well.
\end{enumerate}
\end{point}
\spacingfix{}
\begin{point}{50}%
Stinespring's theorem
is fundamental in the study
of quantum information and quantum computing:
it is used to prove entropy inequalities (e.g.~\cite{lindblad}),
bounds on optimal cloners (e.g.~\cite{werner}),
full completeness of quantum programming languages (e.g.~\cite{staton}),
security of quantum key distribution (e.g.~\cite{werner2,kissinger2017picture}),
to analyze quantum alternation (e.g.~\cite{prakash}),
to categorify quantum processes (e.g.~\cite{selinger}) \emph{and}
as an axiom to single out
quantum theory among information processing theories \cite{chiribella}.
These are only a few examples
    of the many thousands building on Stinespring's discovery.

Stinespring's theorem only applies
to maps of the form~$\scrA \to \scrB(\scrH)$
    and so do most of its useful consequences.
One wonders:
    is there an extension of Stinespring's theorem
    to arbitrary ncp-maps~$\scrA \to \scrB$?
A different and perhaps less frequently asked question is whether
     Stinespring's dilation is categorical in some way?
That is: does it have a defining universal property?
Both questions turn out to have an affirmative answer:
Paschke's generalization of GNS for Hilbert C$^*$-modules~\cite{paschke}
    turns out to have the same universal property
        as Stinespring's dilation and so it extends Stinespring
        to arbitrary ncp-maps.

\begin{center}
    \centerline{\begin{tabular}{l|lll}
    Dilation & Maps & LHS & RHS \\ \hline
        (normal) GNS
    & (n)p-maps
    & (n)miu
        & ncp vector state  \\
          \cite{gelfand1943imbedding}, \sref{gns}, \sref{dils-gns} &
    $\scrA \xrightarrow{\varphi} \C$&
    $\scrA \xrightarrow{\varrho} \scrB(\scrK)$ &
    $\scrB(\scrK) \xrightarrow{\langle x, (\ )x\rangle} \C$ \\
         & \multicolumn{3}{l}{\small (Hilbert space~$\scrK$;
                        von Neumann algebra~$\scrA$; $x \in \scrK$)} \\
     \arrayrulecolor{gray}\hline
        (normal) Stinespring
   & (n)cp-maps
    & (n)miu
        & pure ncp-map\\
        \cite{stinespring},
        \sref{stinespring-theorem}&$\scrA \xrightarrow{\varphi} \scrB(\scrH)$
    &$\scrA \xrightarrow{\varrho} \scrB(\scrK)$ 
   &$\scrB(\scrK) \xrightarrow{\ad_{V}} \scrB(\scrH)$ \\
         & \multicolumn{3}{l}{\small (Hilbert spaces~$\scrH$, $\scrK$;
                        vN alg.~$\scrA$;
                    $V \in \scrB(\scrH,\scrK)$)} \\
     \arrayrulecolor{gray}\hline
        Paschke
   & ncp-maps
    & nmiu
        & pure ncp-map\\
        \cite{paschke,wwpaschke}, \sref{existence-paschke}
        &$\scrA \xrightarrow{\varphi} \scrB$
    &$\scrA \xrightarrow{\varrho} \scrP$
   &$\scrP \xrightarrow{h} \scrB$  \\
         & \multicolumn{3}{l}{\small (Von Neumann algebras~$\scrA$, $\scrB$,
                $\scrP$)} \\
     \arrayrulecolor{gray}\hline
        KSGNS
    & cp-maps
    & miu
        & cp-map\\
        \cite{ksgns}, cf.~\cite{lance}
        &$\scrA\xrightarrow{\varphi} \scrB^a(E)$
    &$\scrA \xrightarrow{\varrho} \scrB^a(F)$
    &$\scrB^a(F) \xrightarrow{\ad_T} \scrB^a(E)$  \\
         & \multicolumn{3}{l}{
             \small (C$^*$-algs~$\scrA$, $\scrB$;
                    Hilb.~$\scrB$-modules~$E$, $F$;
                    $T \in \scrB^a (F,  E)$)}
\end{tabular}}
    \captionof{table}{Various dilation theorems.
    Maps~$\varphi$ of the given type
        split as~$\varphi = h \after \varrho$,
        where~$\varrho$ is as in LHS
        and~$h$ is as in RHS.
    The first three appear prominently in this thesis;
        the last, KSGNS, is only briefly touched upon in \sref{ksgns}.}
\end{center}
\spacingfix{} % ok, wtf?
\spacingfix{}
\begin{point}{60}[overview-dils]%
We start this chapter with a detailed proof of Stinespring's theorem.
Then we show that it obeys a universal property.
Before we move on to Paschke's GNS,
    we need to develop Paschke's theory of self-dual Hilbert C$^*$-modules.
Paschke's work builds on Sakai's characterization of von Neumann algebras,
    which would take considerable effort to develop in detail.
Thus to be as self-contained as possible,
    we avoid Sakai's characterization (in contrast to our
        article~\cite{wwpaschke})
    and give new proofs
    of Paschke's results where required.
One major difference is that we will use
    a completion (see \sref{dils-completion}) of a uniform space
    instead of considering the dual space of a Hilbert C$^*$-module.
This uniformity (which we call the \emph{ultranorm uniformity})
    plays a central r\^ole
    in our development of the theory
    and many of the new results.
We finish the first part of this chapter
    by constructing Paschke's dilation
    and establishing that it indeed extends Stinespring's dilation.
We start the second part of this chapter
    by characterizing when the representation in the Paschke dilation
    of an ncp-map is injective (see \sref{paschke-injective}),
    which is a generalization of our answer \cite{stineinj}
    to the same question about the Stinespring dilation.
Then we allow for an intermezzo
    to prove several new results on (self-dual) Hilbert C$^*$-modules,
    notably:
\begin{enumerate}
    \item
    we generalize the Kaplansky density theorem
    to Hilbert C$^*$-modules, see \sref{kaplansky-hilbmod};
    \item
   we prove that every self-dual Hilbert~C$^*$-module over
            a von Neumann algebra factor
            is of a particularly nice form,
            see \sref{selfdual-normalish-form}, and
    \item
    we greatly simplify the construction of the exterior tensor product
        for Hilbert~C$^*$-modules over von Neumann algebras
        in \sref{univprop-ext-tensor}
        and show that its self-dual completion
        is particularly well-behaved, see e.g.~\sref{ba-ext-tensor-pres}.
\end{enumerate}
We use this new insight into the exterior tensor product
    to compute the tensor product
    of dilations, see~\sref{paschke-tensor}.
In the final sections of this chapter we turn our attention to the maps
    that occur on the right-hand side of the Paschke dilation.
In Stinespring's dilation these maps are always of the form~$\ad_V$,
    but in Paschke's dilation different maps may appear ---
    in any case, they are pure in the sense of~\sref{pure}.
The main result of these last sections
    is an equivalence between this and another
    seemingly unrelated notion of purity:
    we show that an ncp-map~$\varphi$ is pure in the sense of~\sref{pure}
    if and only if the map on the left-hand side of
    the Paschke dilation of~$\varphi$ is surjective.
This is a bridge to the next chapter
    in which we prove in greater generality
    that these pure maps have a~$\dagger$-structure,
    see~\sref{dagger-theorem}.
\end{point}
\end{point}
\end{parsec}

\section{Stinespring's theorem}
\begin{parsec}{1360}[dils-completion-to-hilb]%
\begin{point}{10}%
In the proof of Stinespring's theorem
    (\sref{dils-proof-stinespring}),
    it will be necessary to
    ``complete''\footnote{Note that
        we do not require an inner product to be definite.
    The inner product on a Hilbert space \emph{is} definite.
    Thus the completion will quotient out those vectors with
        zero norm with respect to the inner product.}
        a (complex) vector space with inner product into a Hilbert space.
The details of this classical result were only
    sketched~\sref{completion-inner-product-space},
    where it's used in the GNS-construction~(\sref{gns}).
Here we will work through the details
    as the corresponding completion (\sref{dils-completion}) required
    in Paschke's dilation
    is more complex and its exposition will benefit
    from this familiar analog.
\end{point}

\begin{point}{20}[prop-complete-into-hilbert-space]{Proposition}%
    Let~$V$ be a (complex) vector space with inner
        product~$[\,\cdot\,,\,\cdot\,]$.
        \index{$[\,\cdot\,,\,\cdot\,]$!inner product}
    There is a Hilbert space~$\scrH$
        together with a bounded linear map~$\eta\colon V \to \scrH$
            such that
        \begin{inparaenum}
        \item
        $[v,w] = \left<\eta(v), \eta(w)\right>$
            for all~$v,w \in V$ and
        \item
        the image of~$\eta$ is dense in~$\scrH$.
        \end{inparaenum}
\begin{point}{30}{Proof}%
We will form~$\scrH$ from the set of Cauchy sequences from~$V$
    with a little twist.
Recall that two Cauchy
    sequences~$(v_n)_n$ and~$(w_n)_n$ in~$V$
    are said to be equivalent
    if for every~$\varepsilon > 0$
    there is a~$n_0$
    such that~$\| v_n - w_n \| \leq \varepsilon$
    for all~$n \geq n_0$,
    where as usual $\|v\| \equiv \sqrt{[v,v]}$.
Call a Cauchy sequence~$(v_n)_n$ \Define{fast}\index{Cauchy sequence!fast}
    if for each~$n_0$
    we have~$\| v_n - v_m\| \leq \frac{1}{2^{n_0}}$
    for all~$n,m \geq n_0$.
Clearly every Cauchy sequence has a fast subsequence
    which is (as all subsequences are) equivalent with it.
\begin{point}{40}%
    Define~$\scrH$ to be the set of fast Cauchy sequences modulo
        equivalence.\footnote{We use fast Cauchy sequences
            in the definition of~$\scrH$ only to shorten the text:
            any Cauchy sequence is equivalent to a fast one anyway.}
For brevity we will denote an element of~$\scrH$,
    which is an equivalence class of fast Cauchy sequences, simply by
    a single representative.
Also, we often tersely write~$v$ for the Cauchy sequence~$(v_n)_n$.
The set~$\scrH$ is a metric space
    with the standard
    distance~$d(v, w) \equiv \lim_{n\to\infty} \| v_n - w_n\|$.
It is independent of representatives
    as by the reverse and regular triangle inequality combined,
     we get
\begin{equation*}
    \bigl| \| v_n - w_n \| \,-\, \| v_n' - w_n' \| \bigr|
        \ \leq \ 
    \| v_n - v_n' \| \,+\, \|w_n - w_n'\| \ \rightarrow \ 0
\end{equation*}
for~$v,w$ equivalent with a Cauchy sequence~$v'$, respectively $w'$.
To show~$\scrH$ is a complete metric space, assume
$v^1, v^2, \ldots$
is a fast Cauchy sequence of fast Cauchy sequences.
We will show~$(v_n^n)_n$ is a Cauchy sequence.
Here the restriction to fast sequences bears fruit:
without it, the limit sequence is (notationally) harder to construct.
Assume~$n,m,k \geq N$ for some~$N$.  We have
\begin{equation*}
    \| v^n_n - v^m_m \|
       \  \leq \ 
    \| v^n_n - v^n_k \| \,+\,
    \| v^n_k - v^m_k \| \,+\,
    \| v^m_k - v^m_m \|\  \leq \ 
    \| v^n_k - v^m_k \|+ \frac{2}{2^N}.
\end{equation*}
Because~$\lim_{k\to\infty} \|v_k^n-v_k^m \| =d(v^n,v^m) \leq 2^{-N}$
    we can find a~$k \geq N$
    such that~$\| v^n_k - v^m_k \| \leq \frac{2}{2^N}$
    and so~$\|v^n_n - v^m_m\| \leq \frac{4}{2^N}$.
    Thus~$(v^n_n)_n$ is a Cauchy sequence.
It's easily checked~$v^1, v^2, \ldots$
converges to~$(v^n_n)_n$ with respect to~$d$,
with one caveat: the sequence~$(v^n_n)_n$ might not be fast
    and so, might not be in~$\scrH$.
Like any other Cauchy sequence, $(v^n_n)_n$
    has a fast subsequence.
The equivalence class
    of this subsequence is the limit of~$v^1, v^2, \ldots$ in~$\scrH$.
We have shown~$\scrH$ is a complete metric space.
\end{point}
\begin{point}{50}[prop-hilbert-space-completion-extension]%
Define~$\eta\colon V \to \scrH$
    to be the map which sends~$v$ to the equivalence class of the
    constant sequence~$(v)_n$.
By construction of~$\eta$ and~$\scrH$,
    the image of~$\eta$ is dense in~$\scrH$.
Let~$f\colon V \to X$
    be any uniformly continuous map to a complete metric space~$X$.
We will show there is a unique continuous~$g\colon \scrH \to X$
    such that~$g \after \eta = f$.
From the uniform continuity it easily follows
    that~$f$ maps Cauchy sequences to Cauchy sequences
    and preserves equivalence between them.
Together with the completeness of ~$X$
    there is a unique~$g\colon \scrH \to X$
    fixed by~$g((v_n)_n) = \lim_{n\to\infty}f(v_n)$.
Clearly~$g \after \eta = f$.
Finally, uniqueness of~$g$ follows readily from
    the fact~$g$ is fixed on the image of~$\eta$,
    which is dense in~$\scrH$.
\end{point}
\begin{point}{60}[completion-to-hilb-vect]%
Thus we directly get a scalar multiplication
    on~$\scrH$ by extending~$\eta \after r_z \colon V \to \scrH$
    where~$r_z(v) = zv$ is scalar multiplication by~$z \in \C$,
    which is uniformly continuous.
In fact~$z (v_n)_n \equiv (z v_n)_n$.
Extending addition is more involved, but straightforward:
    given representatives~$v,w \in \scrH$
    we know~$(v_n+w_n)_n$ is Cauchy in~$V$
    by uniform continuity of addition and so
    picking a fast subsequence~$v+w$ of $(v_n+w_n)_n$
        fixes an addition on~$\scrH$.
With this expression for addition, and the similar one for
    scalar multiplication, it is easy to see
    they turn~$\scrH$ into a complex vector space with zero~$\eta(0)$
    for which~$\eta$ is linear.
Extending the inner product is bit trickier.
\end{point}
\begin{point}{70}%
To define the inner product on~$\scrH$,
first note that for Cauchy sequences~$v$ and~$w$ in~$V$
we have (using Cauchy--Schwarz, \sref{chilb-cs}, on the second line):
\begin{align*}
    \bigl|[v_n,w_n] - [v_m,w_m]\bigr|
    & \ =\  \bigl|[v_n,w_n-w_m] + [v_n - v_m,w_m]\bigr| \\
    & \ \leq\  \|v_n\| \|w_n - w_m\| + \|v_n-v_m\|\|w_m\|.
\end{align*}
Thus as~$(\|v_n\|)_n$ and~$(\|w_n\|)_n$ are bounded
we see~$([v_n,w_n])_n$ is Cauchy.
In a similar fashion we see~$\bigl|[v_n,w_n] - [v'_n,w'_n]\bigr|
    \leq \|v_n\| \|w_n - w_n'\| + \|v_n-v_n'\|\|w'_n\|\to 0$ as~$n\to \infty$
    for~$v'$ and~$w'$ Cauchy sequences equivalent
to~$v$ respectively~$w$.
Thus
$([v_n,w_n])_n$ is equivalent to
$([v'_n,w'_n])_n$
    and so
    we define~$\left<v,w\right>
        \equiv \lim_{n\to \infty} [v_n,w_n]$ on~$\scrH$.
With vector space structure from~\sref{completion-to-hilb-vect}
we easily see this is an inner product~$\scrH$.
The metric induced by the inner product coincides with~$d$:
\begin{equation*}
    \left<v-w,v-w\right>
   \ =\ \lim_{n}[v_n-w_n,v_n-w_n]
    \ =\ \lim_{n}\|v_n-w_n\|^2
    \ =\  d(v,w)^2.
\end{equation*}
And so~$\scrH$ is a Hilbert space.
Finally, $\left<\eta(v),\eta(w)\right>=[v,w]$ is direct.\qed
\end{point}
\end{point}
\end{point}
\end{parsec}

\begin{parsec}{1370}[dils-stinespring]%
\begin{point}{10}%
We are ready to prove Stinespring's dilation theorem.
\begin{point}{20}[dils-proof-stinespring]{%
    Proof of Stinespring's theorem, \sref{stinespring-theorem}}%
Let~$\varphi\colon \scrA \to \scrB(\scrH)$
    be any cp-map.
Write~$\scrA \odot \scrH$ for the tensor product of~$\scrA$ and~$\scrH$
    as vector spaces.
    \index{*odot@$\odot$!vector spaces}
We will define~$\scrK$ as the Hilbert space completion
    of~$\scrA \odot \scrH$ with respect to the inner
    product~$[\,\cdot\,,\,\cdot\,]$
    \index{$[\,\cdot\,,\,\cdot\,]$!inner product}
    fixed by
\begin{equation}\label{eq-stinespring-innerprod}
    [a\otimes x,\, b \otimes y] \ \equiv\  \left<x, \,\varphi(a^*b)\, y\right>_{\scrH}.
\end{equation}
To proceed, we have to show that~$[\,\cdot\,,\,\cdot\,]$
    is indeed an inner product.
By linear extension, \eqref{eq-stinespring-innerprod}
    fixes a sesquilinear form~$[\,\cdot\,,\,\cdot\,]$.
As~$\varphi$ preserves involution as a positive map,
see \sref{cstar-p-implies-i},
    the sesquilinear form~$[\,\cdot\,,\,\cdot\,]$
is also symmetric, i.e.~$[t,s]=\overline{[s,t]}$.
    Assume~$\sum^n_{i=1} a_i\otimes x_i$
is an arbitrary element of~$\scrA \odot \scrH$.
To see~$[\,\cdot\,,\,\cdot\,]$ is positive,
we need to show
\begin{equation*}
    0 \ \leq\  \bigl[\sum_i a_i\otimes x_i, \sum_j a_j\otimes x_j\bigr]
        \ \equiv\  \sum_{i,j} \left< x_i, \varphi(a_i^*a_j) x_j \right>.
\end{equation*}
The trick is to consider the matrix algebra~$M_n\scrB(\scrH)$
acting on~$\scrH^{\oplus n}$ in the obvious way:
$(A (x_1,\ldots,x_n)^\T)_i \equiv \sum_j A_{ij} x_j$.
(In fact, this yields an
    isomorphism~$M_n\scrB(\scrH) \cong \scrB(\scrH^{\oplus n})$.)
Writing~$x$ for~$(x_1,\ldots,x_n)^\T$ in~$\scrH^{\oplus n}$,
    we get
\begin{equation}\label{eq-stinespring-norm-tensor}
    \sum_{i,j} \left< x_i, \,\varphi(a_i^*a_j)\, x_j \right>
    \ =\  \bigl<
    x,
    (M_n\varphi) \,\bigl(\, (a_i^*a_j)_{ij} \, \bigr)\,
    x \bigr>_{\scrH^{\oplus n}}.
\end{equation}
The matrix~$(a_i^*a_j)_{ij}$
is positive as~$(C^*C)_{ij} = a_i^*a_j$
    with~$C_{ij} \equiv \frac{1}{\sqrt{n}} a_j$.
    By complete positivity~$M_n\varphi(\,(a_i^*a_j)_{ij}\,)$ is positive
    and so is~\eqref{eq-stinespring-norm-tensor},
    hence~$[\,\cdot\,,\,\cdot\,]$ is an inner product.
Write~$\eta\colon \scrA\odot \scrH \to \scrK$ for the Hilbert space completion
    of~$\scrA \odot \scrH$ with the inner product
    described in~\sref{prop-complete-into-hilbert-space}.
\begin{point}{30}[stinespring-extend-operator]%
We need a lemma before we continue.
Let~$T\colon \scrA \odot \scrH \to \scrA \odot \scrH$
be a bounded linear map.
We show there is an extension~$\hat{T} \colon \scrK \to \scrK$.
The operator~$T$ is
uniformly continuous and so is~$\eta\after T \colon \scrA\odot\scrH \to\scrK$.
Thus by \sref{prop-hilbert-space-completion-extension}
there is a unique continuous extension~$\hat{T} \colon \scrK \to \scrK$
with~$\hat{T}(\eta(t)) = \eta(T(t))$
    for all~$t \in \scrA\odot\scrH$.
Clearly~$\hat{T}$ is linear on the image of~$\eta$,
    which is dense and so~$\hat{T}$ is linear everywhere.
As~$\hat{T}$ is linear and continuous, it must be
    bounded, so~$\hat{T} \in \scrB(\scrK)$.
It is easy to see~$\widehat{T+S}=\hat{T}+\hat{S}$ and
    $\widehat{\lambda T} = \lambda \hat{T}$
    for operators~$S,T$ on~$\scrA \odot \scrH$ and~$\lambda \in \C$.
Also~$\widehat{TS} = \hat{T}\hat{S}$:
indeed~$\hat{T}\hat{S} \eta(t)
            = \hat{T} \eta(St)
            = \eta(TSt)$
    and so by uniqueness~$\hat{T}\hat{S} = \widehat{TS}$.
\end{point}
\begin{point}{40}%
    Assume~$b\in \scrA$.
    Let~$\varrho_0(b)$ be the operator on~$\scrA \odot \scrH$
    fixed by
\begin{equation*}
    \varrho_0(b)\, a\otimes x\, \  \equiv\  (b a) \otimes x.
\end{equation*}
Clearly $\varrho_0$ is linear, unital and multiplicative.
We want to show~$\varrho_0(b)$ is bounded for fixed~$b \in \scrA$,
    so that we can define~$\varrho(b) \equiv \widehat{\varrho_0(b)}
            \in \scrB(\scrK)$ using \sref{stinespring-extend-operator}.
    To show~$\varrho_0(b)$ is bounded, we claim that in~$M_n\scrB(\scrH)$
we have
\begin{equation}\label{stinespring-rho-bound-lemma}
(a_i^*b^*ba_j)_{ij} \ \leq\  \|b\|^2 (a_i^*a_j)_{ij}.
\end{equation}
Indeed, as~$b^*b \leq \|b\|^2$,
    there is some~$c$ with~$c^*c = \|b\|^2 - b^*b$.
Define~$C_{ij} \equiv \frac{1}{\sqrt{n}} ca_j$.
We compute:
$ (C^*C)_{ij} = a_i^*c^*ca_j = a_i^* (\|b\|^2 - b^*b) a_j$
and so~\eqref{stinespring-rho-bound-lemma} holds.
Hence
\begin{alignat*}{2}
    \bigl\| \varrho_0(b) \,\sum_i a_i \otimes x_i \bigr\|^2
    & \ = \ \bigl\| \sum_i (ba_i) \otimes x_i \bigr\|^2 \\
    & \ =\  \bigl< x, (M_n\varphi)(\, (a_i^* b^*b a_j)_{ij}\,)\,x\bigr> \\
    &\ \leq\  \|b\|^2 \, \bigl< x, (M_n\varphi)(\, (a_i^* a_j)_{ij}\,) \, x\bigr>
    &\qquad& \text{by \eqref{stinespring-rho-bound-lemma}}\\
    &\ =\  \|b\|^2\, \bigl\| \sum_i a_i \otimes x_i \bigr\|^2
\end{alignat*}
and so~$\varrho_0(b)$ is bounded.
Define~$\varrho(b) \colon \scrA \to \scrB(\scrH)$
    using~\sref{stinespring-extend-operator}
    by~$\varrho(b) \equiv \widehat{\varrho_0(b)}$.
\end{point}
\begin{point}{50}%
We already know~$\varrho$ is a~mu-map.
We want to show it preserves involution: $\varrho(c^*) = \varrho(c)^*$
for all~$c \in \scrA$.
Indeed, for every~$a,b \in \scrA$ and~$x,y \in \scrH$
    we have
    \begin{equation*}
        [\varrho_0(c^*) \, a\otimes x,b \otimes y]
        \ =\  [(c^* a)\otimes x,b \otimes y]
        \ =\  \left<x, \varphi(a^*cb) y \right>
        \ =\  [a\otimes x,\varrho_0(c)\, b \otimes y].
    \end{equation*}
Hence~$
    \left<\varrho(c^*) \eta(a\otimes x), \eta(b\otimes y)\right>=
    \left< \eta(a\otimes x), \varrho(c)\eta(b\otimes y)\right>$.
    As the linear span of vectors of the form~$\eta(a \otimes x)$
        is dense in~$\scrK$,
        we see~$\varrho(c^*)$ is adjoint to~$\varrho(c)$
        and so~$\varrho(c^*) = \varrho(c)^*$, as desired.
\end{point}
\begin{point}{60}%
    Let~$V_0 \colon \scrH \to \scrA \odot \scrH$
        be given by~$V_0 x = 1 \otimes x$.
        This operator~$V_0$ is bounded:
\begin{equation*}
        \| V_0 x\|\  =\  \left<x, \varphi(1^*1) x\right>^{\frac{1}{2}}
        \ =\  \|\sqrt{\varphi(1)} x\| \ \leq \ \|\sqrt{\varphi(1)}\| \|x\|.
\end{equation*}
    Define~$V \equiv \eta \after V_0$.
    For all~$a \in \scrA$ and~$x,y \in \scrH$, we have
\begin{equation*}
            [a \otimes x, V_0 y]
          \  =\ \left<x, \varphi(a^*)y\right>
          \  =\ \left<x, \varphi(a)^*y\right>
          \  =\ \left<\varphi(a) x, y\right>
\end{equation*}
    and so~$V^*$ satisfies~$V^* \eta(a \otimes x) = \varphi(a)x$.
    Hence for all~$a \in \scrA$ and~$x \in \scrH$
\begin{equation*}
    V^* \varrho(a) V x
        \ =\  V^* \varrho(a) \eta(1 \otimes x)
        \ =\  V^* \eta(a\otimes x)
        \ =\  \varphi(a)x.
\end{equation*}
We have shown~$\ad_V \after \varrho = \varphi$.
    If~$\varphi$ is unital,
    then~$V^*Vx = V^* (1\otimes x) = \varphi(1)x=x$
        for all~$x \in \scrH$ and so~$V$ is an isometry.
    (Or equivalently~$\ad_V$ is unital.)
\end{point}
\begin{point}{70}%
Assume~$\scrA$ is a von Neumann algebra.
It remains to be shown that~$\varrho$
    is normal when~$\varphi$ is.
So assume~$\varphi$ is normal.
To start, note that
    the vector functionals~$\left< a\otimes x, (\,\cdot\,) a\otimes x\right>$
    are together a faithful set of np-functionals.
Thus, by~\sref{normal-faithful},
    it is sufficient
    to show that
\begin{equation*}
    b \ \mapsto \ \left<a \otimes x, \varrho(b) \, a \otimes x \right>
    \ \equiv \ \left<x, \varphi(a^* b a) \,x \right>
\end{equation*}
    is normal for every~$a\in \scrA$ and~$x \in \scrH$.
    This is indeed the case, as normal
    maps are closed under composition and~$\left<x, (\,\cdot\,) x\right>$
    is normal by \sref{hilb-suprema},~$a^* (\,\cdot\,) a$ by \sref{ad-normal}
    and~$\varphi$ by assumption.
        \qed
\end{point}
\end{point}
\begin{point}{80}{Remark}%
Stinespring's theorem generalizes GNS.
It seems tempting to directly prove Stinespring instead of GNS ---
however, for our version with ncp-maps,
we used GNS itself in the proof of Stinespring's theorem to show
that~$\varrho$ is normal.
\end{point}
\end{point}
\end{parsec}

\begin{parsec}{1380}%
\begin{point}{10}%
In the following series of exercises, the reader will develop
    a few consequences of Stinespring's theorem
    which are in some fields better-known than
    Stinespring's theorem itself.
As preparation we study the nmiu-maps between type I von Neumann algebras.
The following result can be derived
    from the more general~\cite[thm.~5.5]{Takesaki1},
    but we will give a direct proof,
    which seems to be new.~\cite{ref1382}
\end{point}
\begin{point}{20}[nmiu-between-type-I]{Proposition}%
Let~$\varrho\colon \scrB(\scrH) \to \scrB(\scrK)$    
    be any non-zero\footnote{
        If~$\scrK$ is 0-dimensional,
            then there zero-map~$\scrB(\scrH) \to \scrB(\scrK)$
            is an nmiu-map.} nmiu-map for Hilbert spaces~$\scrH$ and~$\scrK$.
There is a Hilbert space~$\scrK'$
    and a unitary~$U\colon \scrK \to \scrH \otimes \scrK'$
    (see~\sref{cstar-unitary})
    such that~$\varrho(a) = U^* (a\otimes 1) U $
    for all~$a \in \scrB(\scrH)$.
    \begin{point}{30}[kernel-ultraweak-twosided-ideal-dils]{Proof}%
To start, we will show that~$\varrho$ is injective
    by proving~$\ker\varrho = \{0\}$.
Clearly~$\ker \varrho$ is a two-sided ideal
    --- indeed, if~$\varrho(a) = 0$ for some~$a\in \scrB(\scrH)$,
    then~$\varrho(ab) = \varrho(a)\varrho(b) = 0$
        and~$\varrho(ba) = 0$
        for any~$b \in \scrB(\scrH)$.
As~$\varrho$ is normal,
    we have~$\varrho(\sup D) = \sup_{d \in D} \varrho(d) = 0$
    for any bounded directed~$D\subseteq \ker \varrho$
    of self-adjoint elements.
So by~\sref{prop:weakly-closed-ideal},
    we know~$\ker \varrho = z \scrB(\scrH)$
    for some central projection~$z \in \scrB(\scrH)$.
Either~$z=0$ or~$z=1$.
Clearly~$z\neq 1$, for otherwise~$0=\varrho\neq 0$, which is absurd.
    So~$z=0$ and~$\ker \varrho = \{0\}$, thus~$\varrho $ is injective.
\begin{point}{40}%
Pick an orthonormal basis~$(e_i)_{i \in I}$ of~$\scrH$.
Write~$p_i$ for the projection onto~$\C e_i$,
    i.e.~$p_i \equiv \ketbra{e_i}{e_i}$.
Choose any~$i_0 \in I$
    and write~$\scrK' \equiv \varrho(p_{i_0}) \scrK$.
Pick any orthonormal basis~$(d_j)_{j \in J}$ of~$\scrK'$
    and write~$q_j$ for the projection onto~$\C d_j$,
    i.e.~$q_j \equiv \ketbra{d_j}{d_j}$.
Let~$u_i$ denote the unique partial
        isometry fixed by~$u_i e_i = e_{i_0}$ and
                $u_i e_j = 0$ for~$j\neq i$ ---
                viz.~$u_i \equiv \ketbra{e_{i_0}}{e_i}$.
    Define~$r_{ij} \equiv \varrho(u_i^*) d_j$.
We will show~$(r_{ij})_{(i,j) \in I\times J}$ is an orthonormal
    basis of~$\scrK$.
For this, it is sufficient to
    show~$r_{ij}$ are unit vectors
    and~$\sum_{i,j} \ketbra{r_{ij}}{r_{ij}} = 1$
    (where the sum is defined as supremum over the partial sums),
    using that then~$\ketbra{r_{ij}}{r_{ij}}$
        are pairwise orthogonal projections
        by \sref{orthogonal-tuple-of-projections}.
Clearly these~$r_{ij}$ are unit vectors,
     as by~$u_iu_i^* = p_{i_0}$
        and~$\varrho(p_{i_0}) d_j = d_j$,
        we have~$\|r_{ij}\|^2 = \langle \varrho(u_i u_i^*) d_j,d_j\rangle = 1$.
To show~$1 = \sum_{i,j} \ketbra{r_{ij}}{r_{ij}}$, we compute
\begingroup\allowdisplaybreaks
    \begin{alignat*}{2}
    1 
    & \ = \  \varrho(1)
            &\qquad&\text{as $\varrho$ is unital} \\
    & \ = \  \varrho\bigl( \,\sum_i p_i \,\bigr)
            &&\text{as $e_i$ is an orthn.~basis}\\
        & \ = \  \sum_i \varrho(p_i) &&\text{as $\varrho$ is normal}\\
        & \ = \  \sum_i \varrho(u_i^* p_{i_0} u_i) && \text{by dfn.~$u_i$}\\
    & \ = \  \sum_i \varrho(u_i)^* \varrho(p_{i_0}) \varrho(u_i) 
            && \text{as~$\varrho$ is mi-map}\\
    & \ = \  \sum_i \varrho(u_i)^* \,\bigl( \sum_j q_j \bigr) \, \varrho(u_i)
            && \text{as $d_j$ is an orthn.~basis of~$\scrK'$}
        \\
    & \ = \  \sum_{i,j} \varrho(u_i)^*  \,q_j\, \varrho(u_i) 
            && \text{by \sref{ad-normal}}\\
        & \ = \  \sum_{i,j} \ketbra{r_{ij}}{r_{ij}}.
\end{alignat*}
\endgroup
\end{point}
\begin{point}{50}%
Let~$U \colon \scrK \to \scrH \otimes \scrK'$
    denote the unitary fixed by~$U r_{ij} = e_i \otimes d_j$.
For any~$i\in I$, $j \in J$ and~$a\in\scrB(\scrH)$,
    it is straightforward to show that~$\ketbra{a e_i}{e_{i_0}}
                    = \sum_k \langle e_k, ae_i \rangle \ketbra{e_k}{e_{i_0}}$,
        where the sum converges ultrastrongly
        and consequently
\begingroup\allowdisplaybreaks
\begin{align*}
    \varrho(a)\, r_{ij}
        & \ = \ \varrho(a) \varrho(\ketbra{e_i}{e_{i_0}}) \,d_j 
            &\qquad&\text{by dfn.~$r_{ij}$ and~$u_i$}\\
        & \ = \ \varrho(\ketbra{a e_i}{e_{i_0}}) \,d_j \\
        & \ = \ \varrho\bigl(\, \sum_k \langle e_k, ae_i \rangle
                    \ketbra{e_k}{e_{i_0}} \, \bigr)\,d_j  \\
        & \ = \ \bigl( \sum_k  \langle e_k, ae_i \rangle
                    \, \varrho(\ketbra{e_k}{e_{i_0}})\bigr)\,d_j 
                    &&\text{as $\varrho$ is us-cont.~by~\sref{p-uwcont}} \\
        & \ = \  \sum_k  \langle e_k, ae_i \rangle
                    \, \varrho(\ketbra{e_k}{e_{i_0}})\,d_j \\
        & \ = \ \sum_k  \langle e_k, ae_i \rangle
                    \, U^*\, e_k\otimes d_j 
                    && \text{by dfn.~$U$ and~$u_k$}\\
        & \ = \  U^*\,  \bigl(\, \sum_k e_k \langle e_k, ae_i \rangle 
                   \, \bigr)  \otimes d_j\\
        & \ = \  U^*\,  (ae_i) \otimes d_j 
                   && \text{as $e_k$ is an orthn.~basis}\\
        & \ = \  U^*\,  (a \otimes 1) \, U \, r_{ij} && \text{by dfn.~$r_{ij}$.}
\end{align*}
\endgroup
We have shown~$\varrho(a) = U^* (a \otimes 1) U$. \qed
\end{point}
\end{point}
\begin{point}{60}[typei-inner-auto]{Corollary}%
The nmiu-isomorphisms~$\scrB(\scrH) \to \scrB(\scrK)$
    are precisely of the form~$\ad_U$
    for some unitary~$U\colon \scrK \to \scrH$.
\end{point}
\end{point}
\begin{point}{70}[physics-stinespring]{Exercise*}%
Use \sref{nmiu-between-type-I} and~\sref{stinespring-theorem} to show
    that for every ncp-map~$\varphi\colon \scrB(\scrH) \to \scrB(\scrK)$
    there are a Hilbert space~$\scrK'$
    and a bounded operator~$V\colon \scrK \to \scrH \otimes \scrK'$
    such that~$\varphi(a) = V^* (a \otimes 1)V$.

    Conclude that any quantum channel
        $\Phi$ from~$\scrH$ to itself
        (i.e.~completely positive trace-preserving
            linear map between trace-class operators over~$\scrH$,
            see e.g.~\cite[\S2.6.2]{tomamichel})
        is of the
        form~$\Phi(\varrho) = \TR_{\scrK'}
            [U^* \varrho \otimes \ketbra{v_0}{v_0} U]$
            for some unitary~$U$.
\end{point}
\begin{point}{80}[kraus-exercise]{Exercise* (Kraus' decomposition)}%
Let~$\varphi \colon \scrB(\scrH) \to \scrB(\scrK)$
    be any ncp-map.
By~\sref{physics-stinespring}
    there is a Hilbert space~$\scrK'$
    and bounded operator~$V\colon \scrK \to \scrH \otimes \scrK'$
    such that~$\varphi(a) = V^* (a \otimes 1) V$.
Choose an orthonormal basis~$(e_i)_{i \in I}$ of~$\scrK'$.
Show~$\varphi(a) = \sum_{i \in I} V^* (a \otimes \ketbra{e_i}{e_i}) V$,
where the sum converges ultraweakly.
Deduce there are projections~$P_i \colon \scrH \otimes \scrK' \to \scrH$
with~$\varphi(a) = \sum_{i \in I} V^*P_i^* a P_iV$.
Conclude all ncp-maps~$\varphi\colon \scrB(\scrH) \to \scrB(\scrK)$
are of the form~$\varphi(a) = \sum_{i \in I} V_i^* a V_i$
for some bounded operators~$V_i\colon \scrK \to \scrH$,
for which the partial sums of~$\sum_i V_i^*V_i$ are bounded.
These~$V_i$ are called \Define{Kraus operators}\index{Kraus operators}
    for~$\varphi$.  (For a different approach,
    see e.g.~\cite[thm.~2.3]{davies1976quantum}.)

Now assume~$\scrH$ and~$\scrK$ are finite dimensional.\footnote{%
        For arbitrary~$\scrH$ and~$\scrK$,
            the cardinality of the set of Kraus operators
            can be chosen to be less than or equal
            to~$(\dim \scrH) \cdot (\dim \scrK)$,
            with dimension understood as the cardinality of
                an orthonormal basis and the usual product on cardinals.}
Show that in this case
    the total number of Kraus operators can be chosen
    to be less than or equal to~$(\dim \scrH) \cdot (\dim \scrK)$.
    (There is a similar fact for infinite-dimensional~$\scrH$ and~$\scrK$,
        where the product of numbers is replaced by a maximum of cardinals.)
\end{point}
\end{parsec}

\begin{parsec}{1390}%
\begin{point}{10}{Definition}%
Let~$\varphi\colon \scrA \to \scrB(\scrH)$
    be any ncp-map
    for some Hilbert space~$\scrH$ and von Neumann algebra~$\scrA$.
    A \Define{normal Stinespring dilation}\index{Stinespring dilation}
        (cf.~\cite[ch.~4]{paulsen})
    is a triple~$(\scrK,\varrho, V)$
    with~$\scrK$ a Hilbert space,
    $\varrho\colon \scrA \to \scrB(\scrK)$ an nmiu-map
        (i.e.~normal representation)
        and~$V\colon \scrH \to \scrK$ a bounded operator
        such that~$\varphi$ decomposes as~$\varphi = \ad_V \after \varrho$.
    We just saw every~$\varphi$ has (at least one)
        normal Stinespring dilation.
    A (normal) Stinespring dilation
$(\scrK,\varrho,V)$ is said to be \Define{minimal}\index{Stinespring dilation!minimal}
    if the linear span of~$\varrho(\scrA)V\scrH
        \ \equiv \ \{ \,\varrho(a) Vx; \ a \in \scrA,\ x\in \scrH\,\}$
is dense in~$\scrK$.
By construction the (normal) Stinespring dilation
    from~\sref{dils-proof-stinespring} is minimal.
Every (normal) Stinespring dilation~$(\scrK, \varrho, V)$
    of~$\varphi$
    restricts to a minimal dilation~$(\scrK', \varrho,V)$
    for~$\varphi$, where~$\scrK' \subseteq \scrK$
    is the norm-closure of the linear span of~$\varrho(\scrA)V\scrH$.
\end{point}
\begin{point}{20}%
It is a well-known fact that all minimal normal Stinespring dilations
    for a fixed~$\varphi$ are unitarily equivalent,
    see for instance~\cite[prop.~4.2]{paulsen}.
We will adapt its proof to show that a minimal normal Stinespring
    dilation admits a universal property.
Later we will use this to prove a second universal property \sref{stinespring-is-paschke}
    that allows us to generalize Stinespring's dilation to
    arbitrary ncp-maps.
The modification of the familiar argument
    is mostly straightforward, except for the following lemma,
    which we published earlier as \cite[lem.~11]{wwpaschke}.
\end{point}
\spacingfix{} % wtf?
\parpic[r]{
$\xymatrix{\mathscr{A} \ar[r]^\varrho \ar[rd]_{\varrho'}&
\mathscr{B} \\
& \mathscr{C} \ar[u]_\sigma}$}
\begin{point}{30}[dils-univlemma]{Lemma}%
Let~$\varrho\colon \scrA \to \scrB$ and~$\varrho'\colon \scrA \to \scrC$
    be nmiu-maps between von Neumann algebras,
        and let~$\sigma \colon \scrC \to \scrB$ be an ncp-map
        such that~$\sigma \after \varrho' = \varrho$. Then
\begin{equation*}
\sigma\bigl(\varrho'(a_1)\,c\, \varrho'(a_2)\bigr) \ =\   \varrho(a_1) \,\sigma(c) \,\varrho(a_2)
\end{equation*}
    for all~$a_1,a_2 \in \scrA$ and~$c \in \scrC$.
\begin{point}{40}{Proof}%
By Theorem 3.1 of \cite{choi} (see \sref{choi}),
we know that for all~$c,d\in\scrC$:
\begin{equation}
    \sigma(d^*d) = \sigma(d)^*\sigma(d) \quad \implies \quad
    \sigma(cd) = \sigma(c)\sigma(d).\label{dils-eq-choi}
\end{equation}
We can apply this to~$\varrho'(a)^*\varrho'(a)$
with~$a \in \scrA$ --- indeed we have
\begin{equation*}
\sigma(\varrho'(a)^* \varrho'(a))
  \   =\ \sigma(\varrho'(a^*a))
  \  =\ \varrho(a^*a)
  \  =\ \varrho(a)^*\varrho(a)
  \  =\ \sigma(\varrho'(a))^*\sigma(\varrho'(a))
\end{equation*}
and so by~\eqref{dils-eq-choi},
    we have~$\sigma(c\varrho'(a))
        = \sigma(c)\sigma(\varrho'(a)) = \sigma(c)\varrho(a)$
        for all~$c \in \scrC$.
Thus also (taking adjoints):
$\sigma(\varrho'(a)c) = \varrho(a)\sigma(c)$ for all~$c \in \scrC$.
Hence
\begin{equation*}
\sigma(\varrho'(a_1)c\varrho'(a_2))
            \ =\  \varrho(a_1)\sigma(c \varrho'(a_2))
            \ =\  \varrho(a_1) \sigma(c) \varrho(a_2)
\end{equation*}
    for all~$a_1,a_2 \in \scrA$ and~$c\in \scrC$ as desired.\qed
\end{point}
\end{point}
\begin{point}{50}[dils-univ-stinespring]{Proposition}%
    Assume~$\varphi\colon \scrA \to \scrB(\scrH)$
        is any ncp-map with normal
        Stinespring dilations~$(\scrK, \varrho, V)$
    and~$(\scrK', \varrho', V')$.
        If~$(\scrK,\varrho,V)$ is minimal,
        then there is a unique isometry~$S\colon \scrK \to \scrK'$
        such that~$SV=V'$ and~$\varrho = \ad_S\after \varrho'$.
\begin{point}{60}{Proof}
(This is the same proof we gave in~\cite{wwpaschke}.)
First we will deal with a pathological case.
If~$V = 0$, then~$\varphi = 0$, $V' = 0$, $\scrK = \{0\}$
and~$\varrho = 0$.  Thus the unique linear map~$S\colon \{0\} \to \scrK'$
satisfies the requirements.  Assume~$V \neq 0$.
\begin{point}{70}{Uniqueness}%
Let~$S_1,S_2\colon \scrK \to \scrK'$
be any isometries for which~$S_k V = V'$ and~$\ad_{S_k} \after \varrho' = \varrho$
for~$k =1,2$.
We want to show~$S_1=S_2$.
First we will prove that~$\ad_{S_1} = \ad_{S_2}$.
For~$k=1,2$, $n, m\in \N$, $a_1, \ldots, a_n,\alpha_1, \ldots, \alpha_m \in \scrA$,
$x_1, \ldots, x_n, y_1, \ldots, y_m \in \scrH$ \emph{and}
    $c \in \scrB(\scrK')$,
    we have
\begin{align*}
    \bigl< \ad_{S_k}(c)
        \sum_i \varrho(a_i) V x_i,
        \sum_j \varrho(\alpha_j) V y_j \bigr>
    &\ = \ \sum_{i,j}
        \bigl< V^* \varrho(\alpha_j^*) \,\ad_{S_k}(c)\, \varrho(a_i) V x_i, y_i
            \bigr> \\
            & \ \stackrel{\sref{dils-univlemma}}{=} \ \sum_{i,j}
        \bigl< V^* \ad_{S_k} \bigl(\varrho'(\alpha_j^*) \,c\,
            \varrho'(a_i)\bigr) V x_i, y_j\bigr> \\
    & \ =\  \sum_{i,j}
        \bigl< (V')^* \varrho'(\alpha_j^*) \,c\, \varrho'(a_i) V' x_i, y_j
            \bigr>.
\end{align*}
As the linear span of~$\varrho(\scrA) V\scrH$ is dense in~$\scrK$,
we see~$\ad_{S_1}(c) = \ad_{S_2}(c)$ for all~$c \in \scrB(\scrK')$.
So~$\ad_{S_1} = \ad_{S_2}$
    and thus $\lambda S_1= S_2$ for some~$\lambda \in \C$,
cf.~\cite[lem.~9]{westerbaan2016universal}.
As~$V \neq 0$, there is an~$x \in \scrH$ with~$Vx \neq 0$.
Also~$S_1 V x \neq 0$ as~$S_1$ is an isometry and thus injective.
From this and~$S_1 V x = V'x = S_2Vx = \lambda S_1 Vx$, we get~$\lambda=1$.
Hence~$S_1=S_2$.
\end{point}
\begin{point}{80}{Existence}%
For any~$n \in \N$, $a_1, \ldots, a_n \in \scrA$
    and~$x_1, \ldots, x_n \in \scrH$, we have
\begin{align*}
\bigl\| \sum_i \varrho(a_i)Vx_i \bigr\|^2
&\ =\  \sum_{i,j} \bigl< V^* \varrho(a_j^*a_i) Vx_i, x_j\bigr> \\
&\ =\  \sum_{i,j} \bigl< \varphi(a_j^* a_i)x_i, x_j \bigr> \\
&\ =\  \bigl\| \sum_i \varrho'(a_i) V'x_i \bigr\|^2.
\end{align*}
It follows (again
    using denseness of the linear span of~$\varrho(\scrA) V \scrH$),
that there is a unique isometry~$S\colon \scrK \to \scrK'$
    such that~$S\varrho(a) Vx = \varrho'(a)V' x$
        for all~$a\in \scrA$ and~$x \in\scrH$.
As~$SVx = S\varrho(1)Vx=\varrho'(1)V'x = V'x$ for all~$x\in \scrH$,
    we have~$SV = V'$. Furthermore
we have
\begin{align*}
    S \varrho(a) \sum_i \varrho(a_i)Vx_i
    &\ =\  \sum_i S \varrho(aa_i) Vx_i \\
    &\ =\  \sum_i \varrho'(aa_i) V'x_i \\
    &\  =\ \varrho'(a) \sum_i \varrho'(a_i)V'x_i\\
    &\ =\ \varrho'(a) S \sum_i \varrho(a_i)V x_i.
\end{align*}
for all~$n \in \N$, $a, a_1, \ldots, a_n \in \scrA$ and
    $x_1, \ldots, x_n \in \scrH$.
So~$S \varrho(a) = \varrho'(a)S$.
Thus $S^* \varrho'(a) S = S^*S\varrho(a) = \varrho(a)$ for all~$a\in \scrA$,
whence~$\ad_S \after \varrho' = \varrho$. \qed
\end{point}
\end{point}
\end{point}

\spacingfix{}
\begin{point}{90}[exc-chris-univ-prop]{Exercise*}%
    The statement of the universal property
        \sref{dils-univ-stinespring}
        for Stinespring
        was proposed by Chris Heunen in an unpublished note.
    In this note he shows it is equivalent to the existence
        of a left-adjoint, as we will see in this exercise.
    
    Let~$\mathsf{Rep}_{\mathrm{cp}}$
        denote the category with as objects
            ncpu-maps~$\varphi\colon \scrA \to \scrB(\scrH)$
            and as arrows between~$\varphi\colon \scrA \to \scrB(\scrH)$
            and~$\varphi'\colon \scrA' \to \scrB(\scrH')$
            pairs~$(m, S)$
            where~$m\colon \scrA \to \scrA'$ is an nmiu-map
            and~$S\colon \scrH \to \scrH'$ is an isometry
            such that~$\varphi = \ad_S \after \varphi' \after m$.
    Write~$\mathsf{Rep}$
        for the full subcategory of~$\mathsf{Rep}_{\mathrm{cp}}$
        consisting of objects~$\varphi\colon \scrA \to \scrB(\scrH)$
            that are nmiu.

Show using~\sref{dils-univ-stinespring}
    that the inclusion functor~$U \colon \mathsf{Rep}
    \to \mathsf{Rep}_{\mathrm{cp}}$ has a left adjoint.
\begin{point}{91}{Remarks}%
Recently there have been two interesting categorical
    discoveries concerning Stinespring's dilation theorem,
    which deserve a mention.
Parzygnat found a rather different characterization
    of Stinespring's dilation
    theorem as an adjunction~\cite{parzygnat2018stinespring}.
Secondly, Huot and Staton discovered
    that the category of completely positive maps
    between matrix algebras
    is the universal monoidal category with terminal unit
    with a functor from the category of matrix algebras
    with isometries between them.~\cite{huot}
This fact is surprisingly closely related
    to Stinespring's dilation theorem for matrix algebras.
\end{point}

\spacingfix{}
\begin{point}{100}%
In the special case of ncp-maps
    of the form~$\varphi \colon \scrB(\scrH) \to \scrB(\scrK)$,
    there is a different and noteworthy  property
    concerning dilations `of the same dimension'
    called \emph{essential uniqueness of purification}.
This property is the cornerstone of several publications,
    notably \cite{chiribella}, which popularized it.
\end{point}
\end{point}
\begin{point}{110}[ess-uniq-pur]{Exercise* (Essential uniqueness of purification)}
    Let~$\varphi \colon \scrB(\scrH) \to \scrB(\scrK)$
        be any ncp-map.
    By~\sref{physics-stinespring}
    here is a Hilbert space~$\scrK'$
        and bounded operator~$V\colon \scrK \to \scrH \otimes \scrK'$
        such that~$V^* (a\otimes 1) V = \varphi(a)$.
        Assume there is a bounded
            operator~$W\colon \scrK \to \scrH \otimes \scrK'$
        with~$W^* (a \otimes 1) W = \varphi(a)$ as well.
    We will show there is a unitary~$U \colon \scrK' \to \scrK'$
    such that~$V = (1 \otimes U) W$.

    First show that in the special case the dilations are minimal
    (i.e.~the linear span of~$\{(a \otimes 1) V x;\  a \in \scrB(\scrH), x \in \scrK\}$ is dense
            in $\scrH \otimes \scrK'$),
        there is a (unique) unitary~$U_0$ on $\scrH \otimes \scrK'$
        with~$V = U_0 W$ and~$(a \otimes 1) U_0 = U_0 (a \otimes 1)$
            for all~$a \in \scrB(\scrH)$.
    Derive from the latter
    that for each~$y \in \scrK'$ and unit-vector~$e \in \scrH$,
        there is a~$y' \in \scrK'$
        with~$U_0 (e \otimes y) = e \otimes y'$.
    Conclude~$U_0 = 1 \otimes U$ for some unitary~$U$ on~$\scrK'$.
    Show  that the general case follows from the special case
        in which the dilations are minimal.
\end{point}
\end{parsec}
\begin{parsec}{1400}%
\begin{point}{10}%
    We are ready to state the more general universal property
        that Stinespring's dilation admits,
        which will allow us to generalize it to arbitrary ncp-maps.
\end{point}
\parpic[r]{
$\xymatrix@C=1em@R=1em{
\mathscr{A} 
\ar[rr]^\varphi
\ar[rd]_{\varrho}
\ar@/_1.5em/[rdd]_{\varrho'}
& & \mathscr{B} 
\\ &
\mathscr{P}
\ar[ru]_{h}
& \\ &
\mathscr{P}'
\ar@/_1.5em/[ruu]_{h'}
\ar@{..>}[u]_-{\sigma}
}$}%
\spacingfix{}
\begin{point}{20}[def-paschke]{Definition}%
    Let~$\varphi\colon \scrA \to \scrB$
        be any ncp-map between von Neumann algebras.
    A \Define{Paschke dilation}\index{Paschke dilation} of~$\varphi$
    is a triple~$(\scrP, \varrho, h)$
    with a von Neumann algebra~$\scrP$;
    an nmiu-map $\varrho \colon \scrA \to \scrP$
    and ncp-map~$h \colon \scrP \to \scrB$
    with~$\varphi = h \after \varrho$
    and the following universal property holds.
\begin{quote}
    For every triple~$(\scrP', \varrho', h')$
    with a von Neumann algebra~$\scrP'$;
    an nmiu-map $\varrho'\colon \scrA \to \scrP'$
    and ncp-map~$h' \colon \scrP' \to \scrB$
    with~$\varphi = h' \after \varrho'$,
    there is a unique ncp-map~$\sigma \colon \scrP' \to \scrP$
    (the ``mediating map'')
    with~$\sigma \after \varrho' = \varrho$
    and~$h \after \sigma = h'$.
\end{quote}

\par % XXX hack to get parpic to work
\end{point}
\spacingfix{}
\begin{point}{30}[stinespring-is-paschke]{Theorem}%
Let~$\varphi\colon \scrA \to \scrB(\scrH)$
    be any ncp-map
    together with a minimal normal Stinespring dilation~$(\scrK,\varrho,V)$.
    Then $(\scrB(\scrK), \varrho, \ad_V)$
        is a Paschke dilation of~$\varphi$.
\begin{point}{40}{Proof}%
(This is essentially the same proof
    we gave in \cite{wwpaschke}.)
\parpic[r]{
\xymatrix{
\mathscr{A} \ar@/_/[rdd]_{\varrho'}\ar[rd]^\varrho \ar[rr]^\varphi
&& \scrB(\mathscr{H}) \\
& \scrB(\mathscr{K}) \ar[ru]^{\ad_V} \\
& \mathscr{P}' \ar[ruu]|{h'} \ar[r]_{\tilde\varrho}
& \scrB(\tilde{\mathscr{K}}) \ar[uu]_{\ad_{\tilde V}} \ar[lu]|{\ad_S}
}
}%
\noindent Let~$\scrP'$ be any von Neumann algebra
    together with nmiu-map~$\varrho'\colon \scrA \to \scrP'$
    and ncp-map~$h'\colon \scrP' \to \scrB(\scrH)$
    for which it holds that~$h' \after \varrho' = \varphi$.
We must show that there is a unique
    ncp-map~$\sigma\colon \scrP' \to \scrB(\scrK)$
    such that~$\sigma \after \varrho' = \varrho$
    and~$\ad_V \after \sigma = h'$.
First we will prove uniqueness of this map~$\sigma$.
Then we will show that such a~$\sigma$ exists
    by appealing to the universal property~\sref{dils-univ-stinespring}.

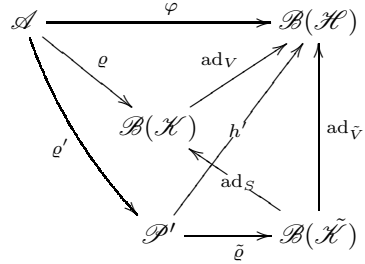
\spacingfix{}
\begin{point}{50}{Uniqueness}%
    Let~$\sigma_1,\sigma_2\colon \scrP' \to \scrB(\scrH)$
        be two such ncp-maps.
    Similarly to the uniqueness proof in \sref{dils-univ-stinespring},
    we compute
        that for all~$n,m \in \N$,
        $a_1, \ldots, a_n,\alpha_1, \ldots, \alpha_m \in \scrA$,
        $x_1, \ldots, x_n,y_1 \ldots, y_m \in \scrH$
        and~$c \in \scrP'$, we have
\begin{align*}
    \bigl< \sigma_k(c)
        \sum_i \varrho(a_i) V x_i,
        \sum_j \varrho(a_j) V x_j \bigr>
    &\ = \ \sum_{i,j}
        \bigl< V^* \varrho(\alpha_j^*) \sigma_k(c) \varrho(a_i) V x_i, y_i
            \bigr> \\
            & \ \stackrel{\mathclap{\sref{dils-univlemma}}}{=} \ \sum_{i,j}
        \bigl< V^* \sigma_k \bigl(\varrho'(\alpha_j^*) c
            \varrho'(a_i)\bigr) V x_i, y_j\bigr> \\
    & \ =\  \sum_{i,j}
            \bigl< h'(\varrho'(\alpha_j^*) c \varrho'(a_i) )x_i, y_j
            \bigr>.
\end{align*}
By definition of minimality of a Stinespring dilation,
the set of vectors of the form~$\sum_i \varrho(a_i) Vx_i$
    is norm dense in~$\scrK$
    and so~$\sigma_1(c)=\sigma_2(c)$ as desired.
\end{point}

\spacingfix{}

\begin{point}{60}{Existence}%
Let~$(\tilde\scrK, \tilde\varrho, \tilde V)$
be a minimal normal Stinespring dilation of~$h'$.
The triple~$(\tilde\scrK, \tilde\varrho \after \varrho', \tilde V)$
is a normal Stinespring dilation of~$\varphi$.
Thus by \sref{dils-univ-stinespring}
    there is a (unique) isometry~$S \colon \scrK \to \tilde\scrK$
    with~$SV = \tilde V$
    and~$\ad_S \after \tilde\varrho \after \varrho' = \varrho$.
Define~$\sigma = \ad_S \after \tilde\varrho$.
Clearly~$\sigma \after \varrho' = \ad_S \after \tilde\varrho \after \varrho'
    = \varrho$
    and~$\ad_V \after \sigma =
        \ad_V \after \ad_S \after \tilde\varrho
        = \ad_{SV} \after \tilde\varrho
        = \ad_{\tilde V} \after \tilde\varrho = h'$
        as desired. \qed
\end{point}
\par
\end{point}
\spacingfix{}
\begin{point}{70}%
The converse holds as well:
    by \sref{paschke-unique-up-to-iso}
    every Paschke dilation of a map~$\varphi \colon \scrA \to \scrB(\scrH)$
    is up to a unique nmiu-isomorphism a minimal Stinespring dilation.
Later we will show that
    every ncp-map~$\varphi\colon \scrA \to \scrB$
    has a Paschke dilation:
    it turns out that Paschke's generalization of
    GNS to Hilbert C$^*$-modules fits the bill.
Before we show this,
    we look at some other examples and basic properties
        of Paschke dilations.
First, as Paschke dilations are defined by a simple universal property,
    they are unique up to a unique nmiu-isomorphism:
\par
\end{point}
\end{point}
\spacingfix{}
\begin{point}{80}[paschke-unique-up-to-iso]{Lemma}%
If $(\scrP_1, \varrho_1, h_1)$
and $(\scrP_2, \varrho_2, h_2)$
    are two Paschke dilations for 
    the same ncp-map $\varphi\colon \scrA \to \scrB$,
    then there is a unique (nmiu)
    isomorphism~$\vartheta \colon \scrP_1 \to \scrP_2$
    with~$\vartheta \after \varrho_1 = \varrho_2$
    and~$h_2 \after \vartheta = h_1$.
\begin{point}{90}{Proof}%
There are unique mediating maps~$\sigma_1\colon \scrP_1 \to \scrP_2$
and~$\sigma_2 \colon \scrP_2 \to \scrP_1$.
It is easy to see~$\sigma_1 \after \sigma_2$
satisfies the same defining property
as the unique mediating map~$\id\colon \scrP_2 \to \scrP_2$
and so~$\sigma_1 \after \sigma_2 = \id$.
Similarly~$\sigma_2 \after \sigma_1 = \id$.
Define~$\vartheta = \sigma_1$.
We just saw~$\vartheta$ is an ncp-isomorphism.
Note~$\vartheta(1) = \vartheta(\varrho_1(1)) = \varrho_2(1) = 1$
and so~$\vartheta$ is unital.
But then
    $\vartheta$ is an nmiu-isomorphism
    by \sref{iso}. \qed
\par
\end{point}
\end{point}
\spacingfix{}
\begin{point}{100}[paschke-basics]{Exercise}%
The following can be shown using only the defining universal property.
\begin{enumerate}
\item
Show~$(\scrA,\varrho,\id)$ is a Paschke dilation of
    an nmiu-map~$\varrho\colon \scrA \to \scrB$.
\item
Let~$(\scrP, \varrho, h)$ be any Paschke dilation.
Prove~$(\scrP, \id, h)$ is a Paschke dilation for~$h$.
(Hint: use~$\id\colon \scrP \to \scrP$
is the unique ncp-map with~$\id \after \varrho = \varrho$
and~$h \after \id = h$.)
\item
    Let~$ \left(\begin{smallmatrix} \varphi_1 \\ \varphi_2
    \end{smallmatrix}\right) \colon \scrA \to \scrB_1 \oplus \scrB_2$
    be an ncp-map.
Show~$(\scrP_1 \oplus \scrP_2,
        \left(\begin{smallmatrix} \varrho_1 \\ \varrho_2 \end{smallmatrix}\right),
        h_1 \oplus h_2)$
        is a Paschke dilation of~$\left(\begin{smallmatrix} \varphi_1 \\ \varphi_2 
        \end{smallmatrix}\right)$
    if~$(\scrP_i, \varrho_i, h_i)$
    is a Paschke dilation of~$\varphi_i$ for~$i=1,2$.
\item
Let~$\varphi\colon \scrA \to \scrB$
    be any ncp-map with Paschke dilation~$(\scrP, \varrho, h)$
    and~$\lambda \in \R$ with $\lambda > 0$.
Prove~$(\scrP, \varrho, \lambda h)$
    is a Paschke dilation of~$\lambda \varphi$.
\end{enumerate}
\spacingfix{}
\end{point}
\begin{point}{110}{Examples}%
Later we will be able to compute more dilations.
\begin{enumerate}
\item
In \sref{paschke-tensor}, we will see
    that if~$(\scrP_i, \varrho_i, h_i)$
    is a Paschke dilation of an ncp-map~$\varphi_i \colon \scrA_i \to \scrB_i$
    for~$i=1,2$,
    then~$(\scrP_1 \otimes \scrP_2, \varrho_1\otimes \varrho_2, h_1\otimes h_2)$
    is a Paschke dilation of~$\varphi_1 \otimes \varphi_2$.
\item
In \sref{paschke-corner}, we will
    show that~$(\cceil{p}\scrA, h_{\cceil{p}}, h_p)$
    is a Paschke dilation of the
    ncp-map~$h\colon \scrA \to p \scrA p$
    given by~$h(a) = pap$ for some projection~$p\in \scrA$.
\item
In \sref{paschke-pure}, it's proven
    that a map~$\varphi$
    with Paschke dilation~$(\scrP, \varrho, h)$
    is pure (in the sense of \sref{pure})
    if and only if~$\varrho$ is surjective.
\item
In \sref{dils-filter-basics-exercise}
    the Paschke dilations of~$\varphi$
    and~$c \after \varphi$ for a filter~$c$ are related.
\end{enumerate}
\spacingfix{}
\end{point}
\end{parsec}

\section{Hilbert C$^*$-modules}
\begin{parsec}{1410}%
\begin{point}{10}%
Recall that by GNS
    every state (pu-map)~$\omega\colon \scrA \to \C$
    on a von Neumann algebra~$\scrA$
    splits as~$h \after \varrho$
    for some nmiu-map~$\varrho\colon \scrA \to \scrB(\scrH)$,
        Hilbert space~$\scrH$
        and vector state~$h (T) =  \left<x, T x\right>$
         for some~$x \in \scrH$.
It is easy to see~$(\scrH, \varrho, h)$ is a Paschke dilation of~$\omega$.

Paschke showed that GNS can be generalized by replacing the scalars~$\C$
    with an arbitrary von Neumann algebra~$\scrB$:
    every `state'~ncp-map $\varphi\colon \scrA \to \scrB$
    splits as~$h \after \varrho$
    for some nmiu-map~$\varrho\colon \scrA \to \scrB^a (X)$,
    and `vector state' $h(T) = \left<x, Tx\right>$,
    where~$X$ is a self-dual Hilbert C$^*$-module over~$\scrB$, $x \in X$
     and~$\scrB^a(X)$ are the~$\scrB$-linear bounded operators on~$X$.
It turns out~$(\scrB^a(X), \varrho, h)$
    is a Paschke dilation of~$\varphi$.
Before we will prove this,
    we first develop the theory of (self-dual) Hilbert C$^*$-modules.

We already encountered Hilbert C$^*$-modules
    in the development of the basics of completely positive maps
    in my twin brother's thesis \cite{bram}.
(They appear in \sref{chilb-basic}.)
As Hilbert C$^*$-modules are central to the present topic,
    we will prove some basic results again, which were already
    covered by my twin.
\end{point}
\begin{point}{20}[dils-basicdfns]{Definition}%
    Let~$\scrB$ be a C$^*$-algebra and~$X$
        be a right $\scrB$-module.
        A \Define{$\scrB$-valued inner product}
            \index{inner product!$\scrB$-valued}
            \index{*langle@$\langle\,\cdot\,,\,\cdot\,\rangle$!$\scrB$-valued}
            on~$X$
        is a map~$\left<\,\cdot\,,\,\cdot\,\right>\colon X \times X \to \scrB$
        such that
        \begin{enumerate}
            \item $\left<x, \,\cdot\,\right> \colon X \to  \scrB$
                    is~$\scrB$-linear for all~$x \in X$
                    --- viz.~$\left< x, yb\right> = \left<x, y\right>b$;
            \item $\left<x,y\right>^* = \left<y,x\right>$
                for all~$x,y\in X$ \emph{and}
            \item $\left<x,x\right> \geq 0$ for all~$x \in X$.
        \end{enumerate}
        A $\scrB$-valued inner product is called \Define{definite}\index{inner product!$\scrB$-valued!definite}
        provided~$x = 0$ whenever~$\left<x,x\right>=0$.
        A \Define{pre-Hilbert $\scrB$-module}\index{Hilbert $\scrB$-module!pre-} $X$
        is a right $\scrB$-module $X$
        together with a definite~$\scrB$-valued inner product.

It will turn out~$\Define{\| x \|} = \| \left<x,x\right> \|^{\frac{1}{2}}$\index{$\|\ \|$!Hilbert module}
        defines a norm on a pre-Hilbert~$\scrB$-module $X$.
A \Define{Hilbert $\scrB$-module}\index{Hilbert $\scrB$-module}
    (also called a \Define{Hilbert C$^*$-module over~$\scrB$})
    is a pre-Hilbert~$\scrB$-module
    which is complete in this norm.
\begin{point}{21}[moved-dfn-selfdual]%
A pre-Hilbert~$\scrB$-module is called \Define{self dual}
                \index{Hilbert~$\scrB$-module!self-dual}
    if for each bounded~$\scrB$-linear~$\tau\colon X \to \scrB$,
                there is a~$t \in X$ such that~$\tau(x) = \left<t,x\right>$
                    for all~$x \in X$.
\end{point}
\begin{point}{22}{Beware}%
In his original definition~\cite[\S2.1]{paschke}, Paschke
    chose to make his~$\scrB$-valued inner product~$\scrB$-linear
    on the left instead of on the right.
He also assumes his inner products to be definite.
\end{point}
\end{point}
\begin{point}{30}{Examples}%
    Several examples\cite{paschke} of (pre-)Hilbert C$^*$-modules follow,
    to each of which we will return in the remainder.
\begin{itemize}
\item
Hilbert $\C$-modules are the same as Hilbert spaces.
They are all self dual.

\item
Any C$^*$-algebra~$\scrB$ is a self-dual Hilbert C$^*$-module over itself
    with $\scrB$-valued inner product~$\left<a,b\right> = a^*b$.
By the C$^*$-identity, the Hilbert~$\scrB$-module
    norm on~$\scrB$ coincides with
    the norm of the C$^*$-algebra.

\item
More generally: if~$J \subseteq \scrB$
    is any right ideal of~$\scrB$,
    then~$J$ is a pre-Hilbert~$\scrB$-module
    with~$\left<a,b\right> = a^*b$.
If~$J$ is closed with respect to the norm, then~$J$ is a Hilbert~$\scrB$-module.
It might not be self dual.

\item
Principal right-ideals are self dual:
    the right ideal~$e \scrB$ for a projection~$e \in \scrB$
    is a self-dual Hilbert~$\scrB$-module with~$\left<ea,eb\right> = a^*eb$.
Later we will see that every self-dual
    Hilbert~$\scrB$-module over a von Neumann algebra~$\scrB$
        is built up from such~$e \scrB$.
    
\item
If~$X$ and~$Y$ are (pre-)Hilbert~$\scrB$-modules,
        then their direct sum~$X \oplus Y$ (as right~$\scrB$-modules)
    is a (pre-)Hilbert~$\scrB$-module
    with the familiar inner product~$\left<(x,y), (x',y')\right>
                = \left<x,x'\right>+\left<y,y'\right>$.
If~$X$ and~$Y$ are self dual, then so is~$X \oplus Y$.
\end{itemize}
\end{point}
\end{parsec}
\subsection{The basics of Hilbert~C$^*$-modules}
\begin{parsec}{1420}%
\begin{point}{10}%
    Perhaps the most important stones in the foundation
        of the theory of operator algebras
        are the Cauchy--Schwarz inequality
        $| \langle x,y\rangle|^2 \leq \|x\|^2\|y\|^2$
        and the polarization
        identity~$\langle x,y\rangle = \frac{1}{4} \sum_{k=0}^3 i^k \|i^k x+ y\|^2$,
        see \sref{inner-product-basic}.
    Also in the development of Hilbert C$^*$-modules,
        their generalizations support the bulk of the structure.
\end{point}
\begin{point}{20}[module-innerprod-state]{Definition}%
Let~$X$ be a right~$\scrB$-module with~$\scrB$-valued inner
    product~$\left<\,\cdot\,,\,\cdot\,\right>$ for some
    C$^*$-algebra~$\scrB$.
For any positive functional~$f\colon \scrB \to \C$,
write~$\Define{\left<x,y\right>_f} = f(\left<x,y\right>)$\index{$\langle \,\cdot\,,\,\cdot\,\rangle_f$!$\scrB$-module}
and~$\Define{\|x\|_f} = \smash{ \langle x,x\rangle^{\nicefrac{1}{2}}_f}$\index{$\|\ \|_f$!$\scrB$-module}.
This~$\left<\,\cdot\,,\,\cdot\,\right>_f$
    is a complex-valued inner product on~$X$.
Hence, by~\sref{inner-product-basic},
    we know that~$\|\,\cdot\,\|_f$ is a seminorm
        (i.e.~a not-necessarily-positive-definite norm.)
\begin{point}{21}%
We start with the Cauchy--Schwarz inequality
    for~$\scrB$-valued inner products.
The proof is the same as Paschke gives
        in \cite[prop.~2.3 (ii)]{paschke},
        but we use it to prove a slightly more general result.
\end{point}
\end{point}
\begin{point}{30}[module-CS]{Proposition (Cauchy--Schwarz)}%
Let~$X$ be a right~$\scrB$-module
    with~$\scrB$-valued inner product~$\left<\,\cdot\,,\,\cdot\,\right>$
    for some~C$^*$-algebra~$\scrB$.
Then~$\left<x,y\right>\left<y,x\right> \ \leq\  \|\left<y,y\right>\| \left<x,x\right>$.
\begin{point}{40}{Proof}%
Let~$f\colon \scrB \to \C$ be any state (pu-map).
Since the states on~$\scrB$ are order separating,
see \sref{states-order-separating},
it suffices to show~$f(\left<x,y\right>\left<y,x\right>)
\leq \|\left<y,y\right>\| f(\left<x,x\right>)$.
If~$f(\left<x,y\right>\left<y,x\right>) = 0$
the inequality holds trivially,
so assume~$f(\left<x,y\right>\left<y,x\right>) \neq 0$.
Using Cauchy--Schwarz for~$\left<\,\cdot\,,\,\cdot\,\right>_f$,
    we derive
\begin{align*}
    f(\left<x,y\right>\left<y,x\right>)^2
        &\  =\  \left<x,y\left<y,x\right>\right>_f^2 \\
        &\  \leq\  \|x\|_f^2 \, \|y \left<y,x\right>\|_f^2 \\
    & \ =\  \|x\|_f^2\, f(\left<x,y\right>\left<y,y\right>\left<y,x\right>)\\
    & \ \leq\  \|x\|_f^2 \,
    \|\left<y,y\right>\| \,
    f(\left<x,y\right>\left<y,x\right>),
\end{align*}
    which yields the inequality by dividing by
$f(\left<x,y\right>\left<y,x\right>)$. \qed
\end{point}
\end{point}
\begin{point}{50}[module-seminorm]{Exercise}
Let~$X$ be a right $\scrB$-module
    with~$\scrB$-valued inner product~$\left<\,\cdot\,,\,\cdot\,\right>$
    for some C$^*$-algebra~$\scrB$.
First use Cauchy--Schwarz (\sref{module-CS})
    to show~$\|\left<x,y\right>\| \leq \|x\|\|y\|$.
With this, derive $\|x\| =\|\left<x,x\right>\|^{\frac{1}{2}}$ is a seminorm on~$X$ with
    $\|x \cdot b \| \leq \|x\|\|b\|$
    for all~$b \in \scrB$.
    (See \cite[prop.~2.3]{paschke}.)
\begin{point}{60}%
From this, it follows easily that
that the addition and the module action
are jointly uniformly continuous
    and that the inner product is uniformly continuous in each argument
    separately.
(The inner product is also jointly continuous,
    but not uniformly.)
\end{point}
\end{point}
\begin{point}{70}{Definition}%
Let~$V$ be a right $\scrB$-module for some~C$^*$-algebra~$\scrB$.
A~\Define{$\scrB$-sesquilinear form}\index{sesquilinear form!$\scrB$-valued} on~$V$
    is a map~$B\colon V\times V \to \scrB$
    such that for each~$y\in V$
    the maps~$x \mapsto B(y, x)$
        and~$x \mapsto B(x,y)^*$ are~$\scrB$-linear --- that is:
\begin{equation*}
\beta^* \,\langle x,y\rangle \,b \ =\ \langle \,x \beta\,,\, y b \,\rangle 
    \qquad \text{for all $x,y \in X$ and~$b,\beta \in \scrB$}.
\end{equation*}
\spacingfix{}
\begin{point}{71}{Remark}%
It seems that~$\scrB$-sesquilinear forms do not explicitly appear
    in the established literature,
    even though its associated polarization identity
    is often used implicitly.
\end{point}
\end{point}
\begin{point}{80}{Example}%
Let~$X$ be a pre-Hilbert~$\scrB$-module.
For every~$\scrB$-linear~$T\colon X \to X$
    the map~$\left<(\,\cdot\,), T (\,\cdot\,)\right>$
    is a~$\scrB$-sesquilinear form.
In \sref{hilbmod-sesquilinear-forms}
    we will see that on a self-dual Hilbert~C$^*$-module
    all \emph{bounded}
    $\scrB$-sesquilinear forms
    are of this form.
\end{point}
\begin{point}{90}[hilbmod-polarization]{Exercise}%
Let~$B$ be a $\scrB$-sesquilinear form.
Prove the \Define{polarization identity}\index{polarization identity}:
\begin{equation*}
    B(x,y) \ =\  \frac{1}{4} \sum_{k=0}^3 i^k B(i^k x + y, i^k x + y).
\end{equation*}
\end{point}
\spacingfix{}
\end{parsec}

\begin{parsec}{1430}%
\begin{point}{10}{Definition}%
Let~$X$, $Y$ be pre-Hilbert~$\scrB$-modules
    for some~C$^*$-algebra~$\scrB$.
    A ($\C$-)linear map~$T\colon X \to Y$
    is said to be \Define{adjointable}\index{adjointable}
    if there is a linear map~$S \colon Y \to X$
    such that~$\left<y, Tx\right>_Y = \left<Sy, x\right>_X$
    for all~$x \in X$ and~$y \in Y$.
Adjoints are unique:
    if~$S$ and~$S'$ are both adjoints of~$T$
    an easy computation shows~$\left<(S-S')y,(S-S)'y\right>=0$
    for all~$y \in Y$.
We write \Define{$T^*$} for the adjoint of~$T$ (if it exists)
and~$\Define{\scrB^a(X)}$\index{$\scrB^a(X)$}
    for the subset of (norm-)bounded operators on~$X$
    which are adjointable. (See \cite[\S2]{paschke}.)
\end{point}
\begin{point}{11}%
The second part of the following lemma is well-known,
    but is proven in a new way using
    the equivalence which is inspired
    by~\cite[thm.~2.8 (ii)]{paschke}.
\end{point}
\begin{point}{20}[adjointable-cstar-identity]{Lemma}%
For a linear map~$T \colon X \to Y$
between pre-Hilbert~$\scrB$-modules~$X$,~$Y$ \emph{and}
real number~$B > 0$,
the following are equivalent.
\begin{enumerate}
    \item $\|T x\| \leq B \|x\|$ for all~$x \in X$.
    \item $\|\left<y,Tx\right>\| \leq B \|y\|\|x\|$
            for all~$x \in X$ and~$y \in Y$.
\end{enumerate}
If~$T$ is bounded and adjointable, then 
    ~$\|T^*\| = \|T\|$ and~$\|T^*T\|=\|T\|^2$.
\begin{point}{30}{Proof}%
   Assume~$\|Tx\| \leq B \|x\|$ for all~$x \in X$.
   Then by \sref{module-seminorm}
   we find that for all~$y \in X$ we have~$\|\left<y, Tx\right>\|
            \leq \|y\|\|Tx\|
            \leq B \|y\|\|x\|$, as desired.

For the converse, pick~$x\in X$
    and assume~$\|\left<y,Tx\right>\| \leq B \|y\|\|x\|$
            for all~$y \in Y$.
Then~$\|Tx \|^2 = \|\left<Tx,Tx\right>\|
                \leq B\|Tx\|\|x\|$.
So~$\|Tx\| \leq B\|x\|$
    by dividing~$\|Tx\|$ if~$\|Tx\|\neq0$
    and trivially otherwise.

Now assume~$T$ is adjointable and bounded.
Then~$\|\left<x,T^*y\right>\| = \|\left<y, Tx\right>\| \leq \|T\|\|y\|\|x\|$
    for all~$x \in X$, $y\in Y$
    and so by the previous~$\| T^*\| \leq \|T\|$.
As adjoints are unique and~$T$ is an adjoint of~$T^*$
    we get~$T^{**}=T$
    and so~$\|T\| = \|T^{**}\| \leq \|T^*\| \leq \|T\|$, as desired.

For the final identity, note that for any~$x \in X$ we have
\begin{equation*}
\|Tx\|^2 \ =\  \|\left<Tx,Tx\right>\| 
            \ =\  \| \left<x, T^*Tx\right>\|
            \ \leq\  \|x \| \|T^*Tx \|
            \ \leq\  \|x\|^2 \|T^*T\|
\end{equation*}
and so~$\|T\| \leq \|T^*T\|^{\frac{1}{2}}$.
Hence~$\|T\|^2 \leq \|T^*T\| \leq \|T^*\|\|T\| = \|T\|^2$.
    \qed
\end{point}
\end{point}
\begin{point}{40}[hilbmod-cstar]{Proposition}%
Let $X$ be a Hilbert~$\scrB$-module for some~C$^*$-algebra~$\scrB$.
With composition as multiplication,
    adjoint as involution
    and operator norm,
    the set of adjointable operators
     $\scrB^a(X)$ is a C$^*$-algebra. \cite[\S2]{paschke}
\begin{point}{50}{Proof}
It is easy to see~$T^*+S^*$ is an adjoint of~$T+S$
    for~$T,S \in \scrB^a(X)$
    and so~$T+S$ is adjointable with~$(T+S)^* = T^*+ S^*$.
    Similarly it is easy to see~$\scrB^a(X)$
        is closed under scalar multiplication, composition and involution
        with~$(\lambda T)^* = \overline{\lambda} T^* $,
            $(TS)^* = S^*T^*$ and
            $T^{**} =T$.
Also the identity~$1$ is clearly bounded and self-adjoint.
So~$\scrB^a(X)$ is a unital~$*$-algebra.
By \sref{adjointable-cstar-identity} 
    the C$^*$-identity holds.
It only remains to be shown that~$\scrB^a(X)$ is complete.
Let~$T_1,T_2,\ldots$ be a Cauchy sequence in~$\scrB^a(X)$.
As~$X$ is a Banach space $\scrB(X)$ is complete
    and so~$T_n \to T$ for some~$T \in \scrB(X)$.
We have to show~$T$ is adjointable.
For any~$n,m \in \N$ we have~$\|T_n^* - T_m^*\|
= \|(T_n - T_m)^*\| = \|T_n - T_m\|$
and so~$T^*_1, T^*_2,\ldots$
is also Cauchy
and converges to, say, $S \in \scrB(X)$.
So
\begin{equation*}
\left<S x, y\right> =
\langle\lim_{n} T^*_n x, y\rangle = \lim_n \left<T_n^* x, y\right>
                = \lim_n \left<x, T_n y\right>
                =  \left<x, T y\right>,
\end{equation*}
for any~$x \in X$ and~$y \in Y$.
Thus~$S = T^*$ and~$T \in \scrB^a(X)$ as desired. \qed
\end{point}
\end{point}
\end{parsec}
\begin{parsec}{1440}%
\begin{point}{10}[hilbmod-ordersep]{Proposition}%
Let~$X$ be a Hilbert~$\scrB$-module for some~C$^*$-algebra~$\scrB$.
The vector states on~$\scrB^a(X)$ are order separating ---
that is: $T \geq 0$ if and only if~$\left<x,Tx\right> \geq 0$
for all~$x \in X$.
\begin{point}{20}{Proof}%
(Proof from \cite[lem.~4.1]{lance}.)
Assume~$T \in \scrB^a(X)$ with~$T \geq 0$.
Then~$T = S^*S$ for some~$S \in \scrB^a(X)$
    and so~$\left<x, Tx\right> = \left<Sx,Sx\right> \geq 0$
    for all~$x \in X$.

For the converse, assume~$\left<x, Tx\right> \geq 0$
    for all~$x \in X$.
We claim $T$ is self-adjoint.
By the polarization identity
$\left<x, Ty\right>
    = \frac{1}{4} \sum^3_{k=0} i^k \left< i^k x+y, T(i^k x+y)\right>$,
    see \sref{hilbmod-polarization},
    we have~$\left<x,Ty \right> = \left<Tx, y\right>$ for all~$x,y \in X$
    and so~$T^*=T$.
There are positive~$T_+,T_- \in \scrB^a(X)$
with~$T = T_+ - T_-$ and~$T_+T_- = 0$, see \sref{cstar-pos-neg-part}.
By assumption~$0 \leq \left<T_-x,TT_-x\right>$
    from which~$-\left<x,T_-^3x\right> \geq 0$,
    but~$T_-^3 \geq 0$ so~$T_-^3 = 0$.
    Thus~$T_- = 0$ (by the functional calculus, \sref{functional-calculus})
    and so~$T \geq 0$. \qed
\end{point}
\end{point}
\begin{point}{30}{Lemma}%
Let~$T\colon X \to Y$
    be an adjointable linear map
    between pre-Hilbert $\scrB$-modules~$X,Y$.
    Then~$T$ is~$\scrB$-linear. \cite[\S2]{paschke}
\begin{point}{40}{Proof}%
We have
    $\left<y, (Tx)b\right>
    =\left<T^*y, x\right>b
    =\left<y, T(xb)\right>$
for any~$x \in X$, $y \in Y$ and~$b \in \scrB$.
In particular we
get~$\left<(Tx)b-T(xb), (Tx)b-T(xb)\right>=0$
taking~$y = (Tx)b-T(xb)$
    and so~$T(x)b=T(xb)$. \qed
\end{point}
\end{point}
\begin{point}{50}[blinear-inprod-inequality]{Proposition}%
Let~$X$ and~$Y$ be right~$\scrB$-modules with~$\scrB$-valued
    inner products for some~C$^*$-algebra~$\scrB$.
If~$T \colon X \to Y$ is a bounded~$\scrB$-linear map,
    then we have the
        inequality~$\left<Tx,Tx\right>\leq \|T\|^2 \left<x,x\right>$
    for any~$x \in X$.  (cf.~\cite[rem.~2.9]{paschke}.)
\begin{point}{60}{Proof}
(The proof is a slight variation on the first part
of \cite[thm.~2.8]{paschke}.)
Pick~$\varepsilon > 0$.
By \sref{cstar-positive-1}
    we know $\left<x,x\right> + \varepsilon$ is invertible.
Define~$h = (\left<x,x\right> + \varepsilon)^{-\frac{1}{2}}$.
Clearly
\begin{equation*}
0 \ \leq\  \left<xh,xh\right>
    \ =\  h (\left<x,x\right> + \varepsilon)h - h^2 \varepsilon
    \ =\  1 - h^2 \varepsilon \ \leq\  1.
\end{equation*}
Thus~$\|xh\| \leq 1$
and so~$\|h\left<Tx, Tx\right>h\| = \|T(xh)\|^2 \leq \|T\|^2$.
From this and~$0 \leq h\left<Tx,Tx\right>h$
we get~$h\left<Tx,Tx\right>h \leq \|T\|^2$. Hence, by dividing by~$h$
    on both sides,
\begin{equation*}
\left<Tx,Tx\right> \ \leq \ \|T\| h^{-2}\  =
    \  \|T\|^2(\left<x,x\right> + \varepsilon)
\end{equation*}
for all~$\varepsilon > 0$
and so~$\left<Tx,Tx\right> \leq \|T\|^2 \left<x,x\right>$ as desired. \qed
\end{point}
\end{point}
\end{parsec}

\begin{parsec}{1450}%
\begin{point}{10}[hilbmod-vectstates-cp]{Proposition}%
Let~$\scrB$ be a von Neumann algebra,
    $X$ be a Hilbert~$\scrB$-module
    and~$x \in X$.
Then the map~$h\colon \scrB^a(X )\to \scrB$
defined by~$h(T) = \left<x,Tx\right>$
is completely positive.
\begin{point}{20}{Proof}%
Pick any~$n \in \N$, ~$T_1, \ldots, T_n \in \scrB^a(X)$
    and~$b_1, \ldots, b_n \in \scrB$
    and compute
\begin{align*}
\sum_{i,j} b_i^* h(T_i^*T_j)b_j
&\ =\ \sum_{i,j} b_i^* \left<x,T_i^*T_j x\right> b_j\\
&\ =\ \sum_{i,j} \left<T_i x b_i,T_j x b_j\right>  \\
&\ =\ \bigl\langle\sum_i T_i x b_i, \sum_i T_i x b_i\bigr\rangle \ \geq\  0,
\end{align*}
which shows~$h$ is indeed completely positive. \qed
\end{point}
\end{point}
\end{parsec}

\subsection{Ultranorm uniformity}
\begin{parsec}{1460}%
\begin{point}{10}%
In \sref{prop-complete-into-hilbert-space}
    we saw how to complete a complex vector space with inner product
    to a Hilbert space.
Now we will see how to complete a right $\scrB$-module
    with $\scrB$-valued inner product
    to a \emph{self-dual} Hilbert~$\scrB$-module
    under the assumption~$\scrB$ is a von Neumann algebra.
    We will use a different construction than Paschke used~\cite{paschke}.
Paschke showed that for a pre-Hilbert $\scrB$-module~$X$
    the set of functionals~$X'$
        (i.e.~$\scrB$-module homomorphisms into~$\scrB$)
    turns out to  be a self-dual Hilbert~$\scrB$-module,
    which he uses as the completion of~$X$.
A considerable part of his paper~\cite{paschke}
    is devoted to this construction.
It requires Sakai's characterization of von Neumann algebras,
    which we have not covered.
To avoid developing Sakai's theory,
    we give a different construction.
Instead of embedding~$X$ into a dual space,
    we will stay closer to the similar fact for Hilbert spaces
    and use a topological completion.
A simple metric completion will not do,
    we will need to complete $X$ as a \emph{uniform space}.
Uniformities are structures that sit between topological spaces
    and metric spaces: every metric space is a uniformity
        (see~\sref{dils-uniformity-examples})
        and every uniformity is a topological space.
We will refresh the basics of uniformities ---
for a thorough treatment, see for instance \cite[ch.~9]{willard}.
\end{point}
\begin{point}{20}[dils-dfn-uniformity]{Definition}%
    A \Define{uniform space}\index{uniform space} is a set~$X$
    together with a family of relations
    $\Phi \subseteq \wp (X\times X)$ called \Define{entourages}\index{entourages}
    \index{p@$\wp$, power set}
        satisfying the following conditions.
    \begin{enumerate}
    \item
        The set of entourages~$\Phi$ is a filter.
        That is:
        \begin{inparaenum}
        \item
        if~$\varepsilon,\delta \in \Phi$,
            then~$\varepsilon \cap \delta \in \Phi$
        \emph{and}
        \item
        if~$\varepsilon \subseteq \delta$ and~$\varepsilon \in \Phi$,
            then~$\delta \in \Phi$.
        \end{inparaenum}
    \item
        For every~$\varepsilon \in \Phi$
            and~$x \in X$ we have~$x \mathrel\varepsilon x$.
    \item
        For each~$\varepsilon \in \Phi$,
            there is a~$\delta \in \Phi$
            such that~$\delta^2 \subseteq \varepsilon$.
        So, if~$x \mathrel\delta y$ and~$y \mathrel\delta z$
            then~$x \mathrel\varepsilon z$ for any $x,y,z \in X$.
    \item
        For every~$\varepsilon \in \Phi$,
            there is a~$\delta \in \Phi$
                with~$\delta^{-1}\subseteq \varepsilon$.
        So, if~$x \mathrel\delta y$ then~$y \mathrel\varepsilon x$
            for any~$x,y \in X$.
    \end{enumerate}
    A uniform space is \Define{Hausdorff}\index{uniform space!Hausdorff}
        whenever~$x \mathrel\varepsilon y$
            for all~$\varepsilon \in \Phi$
            implies that~$x=y$.
\begin{point}{30}%
The elements of~$X$ are the points of the uniform space.
The entourages~$\varepsilon \in \Phi$
    are generalized distances:
one reads~$x \mathrel\varepsilon y$ as
    `the point $x$ is at most~$\varepsilon$ far from to~$y$'.
With this in mind, the second axiom states every point is arbitrarily close
    to itself.
The third axiom requires that for every entourage~$\varepsilon$
    there is an entourage which acts like~$\nicefrac{\varepsilon}{2}$.
There might be several~$\nicefrac{\varepsilon}{2}$
    which fit the bill.
Whenever we write~$\nicefrac{\varepsilon}{2}$
    we implicitly pick some entourage that
    satisfies~$(\nicefrac{\varepsilon}{2})^2 \subseteq \varepsilon$.
Also we will use the obvious shorthand~$\nicefrac{\varepsilon}{4}
=   \nicefrac{(\nicefrac{\varepsilon}{2})}{2}$.
\end{point}
\end{point}
\begin{point}{31}{Definition}%
Let~$X$ be a set together with
    a family of relations~$B \subseteq \wp (X \times X)$.
The set~$B$ is said to be a \Define{subbase (for a uniformity on~$X$)}
    if it satisfies axioms 2, 3 and 4
    of \sref{dils-dfn-uniformity}. (Cf.~\cite[dfn.~5.5]{willard}.)
\end{point}
\begin{point}{40}[exc-subbase]{Exercise}%
Let~$X$ be a set together with a subbase~$B$.
Write~$\Phi$ for the filter generated by~$B$
    (that is: $\delta \in \Phi$
        iff~$\varepsilon_1 \cap \ldots \cap \varepsilon_n \subseteq \delta$
            for some~$\varepsilon_1, \ldots, \varepsilon_n \in B$).
Show that~$(X,\Phi)$ is a uniform space.
\end{point}
\begin{point}{50}[dils-uniformity-examples]{Examples}%
Using the definition of a subbase, it is easy
to describe the entourages of some common uniformities.
    \begin{enumerate}
        \item
    Let~$(X,d)$ be a metric space.
Define~$\hat\varepsilon \equiv \{(x,y);\ d(x,y) \leq \varepsilon\}$
for any $\varepsilon > 0$.
The set~$B \equiv \{ \hat\varepsilon; \ \varepsilon > 0\}$
is a subbase and so fixes a uniformity~$\Phi$ on~$X$.
In this sense every metric space is a uniformity.
        \item
    Let~$X$ be a set together with a (infinite)
            family~$(d_{\alpha})_{\alpha\in I}$
            of pseudometrics.\footnote{Like a metric,
                a pseudometric~$d$ is symmetric,
                obeys the triangle inequality
                    and~$d(x,x)=0$,
                    but unlike a metric~$d(x,y)=0$
                    is not required to entail~$x=y$.}
    Define~$E_{\alpha,\varepsilon} = \{ (x,y); \ d_\alpha(x,y)
            \leq \varepsilon\}$.
            Then~$B \equiv \{ E_{\alpha,\varepsilon}; \ \varepsilon > 0, \ 
                    \alpha \in I\}$
                    is a subbase and fixes
                    a uniformity~$\Phi$ on~$X$.
(Every uniform space is of this form, see e.g.~\cite[thm.~39.11]{willard}.)
\item
Let~$V$ be a vector space with a
family of seminorms~$(\|\  \|_\alpha)_{\alpha \in I}$.
As a special case of the previous,
 $V$ has a uniformity fixed by
    pseudometrics~$d_\alpha(x,y) = \|x-y\|_\alpha$.
    (If additionally~$x =0$ provided~$\|x\|_\alpha = 0$ for all~$\alpha$,
    then~$V$ is called a locally convex space.)

\item
Assume~$\scrB$ is a von Neumann algebra.
There are several uniformities on~$\scrB$ --- two of which
    are of particular interest to us.
    The \Define{ultrastrong uniformity}\index{ultrastrong} on~$\scrB$
is fixed by the seminorms~$\|b\|_f = f(b^*b)^{\frac{1}{2}}$
        where~$f\colon \scrB \to \C$ is an np-map.
        The \Define{ultraweak uniformity}\index{ultraweak} on~$\scrB$
    is given by the seminorms~$|f(b)|$
        for np-maps~$f\colon \scrB \to \C$.

    \end{enumerate}
    \spacingfix{}
\begin{point}{60}%
We are ready to define a uniformity for Hilbert $\scrB$-modules.
Later, we will complete~$X$ with respect to this uniformity.
\end{point}
\end{point}
\begin{point}{70}[dils-ultranorm]{Definition}%
Let~$\scrB$ be a von Neumann algebra.
Assume~$X$ is a right $\scrB$-module
    with~$\scrB$-valued inner product~$[\,\cdot\,,\,\cdot\,]$.
    \index{$[\,\cdot\,,\,\cdot\,]$!inner product!$\scrB$-valued}
We call the uniformity on~$X$
    given by seminorms~$\|x\|_f = f([x,x])^{\frac{1}{2}}$
    for np-maps~$f\colon \scrB \to \C$
    the \Define{ultranorm uniformity}\index{ultranorm}.
\begin{point}{80}%
The ultranorm uniformity will play a very similar r\^ole
    to the norm for Hilbert spaces.
If~$\scrB=\C$, then the ultranorm uniformity is
    the same as the uniformity induced by the norm.
If~$X=\scrB$ with~$[a,b]=a^*b$,
    then the ultranorm uniformity coincides with the ultrastrong uniformity.
\end{point}
\begin{point}{90}{Beware}%
    The ultranorm uniformity is (in general) not given by a single norm.
    Furthermore the ultranorm uniformity is weaker than the norm uniformity
        (that is: norm convergence implies ultranorm convergence,
            but not necessarily the other way around),
        even though on von Neumann algebras the
        ultraweak uniformity
        is stronger than the weak uniformity.
\end{point}
\begin{point}{100}{Remarks}%
The ultranorm uniformity has appeared earlier in the literature
    under various different names; for instance as the~$\tau_1$-topology
    in~\cite[thm.~3.5.1]{manuilov2000hilbertc} and as the s-topology
    in~\cite[2.9]{ghez1985w}.
In previous work the ultranorm uniformity only plays a minor r\^ole,
    often amongst other topologies.
In this thesis we put the ultranorm uniformity center-stage
    as the right generalization of the norm
    and give it the `ultranorm' name to match this view.
\end{point}
\end{point}
\end{parsec}

\begin{parsec}{1470}%
\begin{point}{10}[uniformity-basics]{Definition}%
Let~$X$ be a uniform space with entourages~$\Phi$.
It is easy to translate the familiar notions for
    metric spaces to uniform spaces.
\begin{enumerate}
    \item A net~$(x_\alpha)_\alpha$ is said to
        \Define{converge}\index{uniform space!convergence} to~$x$
            if for each~$\varepsilon \in \Phi$
            there is an~$\alpha_0$
            such that for all~$\alpha > \alpha_0$
                we have~$x \mathrel\varepsilon x_\alpha$.
Note that a net can converge to two different points.
Indeed, this is the case if we start of with a pseudometric
    that is not a metric.
\item A net~$(x_\alpha)_\alpha$ is called \Define{Cauchy}\index{Cauchy net}
        if for each~$\varepsilon \in \Phi$
            there is an~$\alpha_0$
            such that for all~$\alpha,\beta > \alpha_0$
            we have~$x_\alpha \mathrel\varepsilon x_\beta$.
            The uniform space $X$ is \Define{complete}\index{uniform space!complete}
    when every Cauchy net converges.
\item
Let~$(Y,\Psi)$ be a second uniform space.
A map~$f\colon X \to Y$
is said to be \Define{uniformly continuous}\index{uniform space!uniformly continuous}
if for each~$\varepsilon \in \Psi$
there is a~$\delta \in \Phi$
such that for all~$x \mathrel\delta y$
we have~$f(x) \mathrel\varepsilon f(y)$.
The map~$f$ is merely \Define{continuous}\index{uniform space!continuous}
if for each~$x$ and~$\varepsilon\in\Psi$,
there is a~$\delta \in \Phi$
such that for all~$x \mathrel{\delta} y$
we know~$f(x) \mathrel\varepsilon f(y)$.
\item
We say two Cauchy nets
$(x_\alpha)_{\alpha\in I}$
and~$(y_\beta)_{\beta\in J}$ are \Define{equivalent}\index{Cauchy net!equivalent},
in symbols: $(x_\alpha)_\alpha \sim (y_\beta)_\beta$,
when for every~$\varepsilon \in \Phi$
there are~$\alpha_0 \in I$ and~$\beta_0 \in J$
such that for all~$\alpha \geq \alpha_0$ and~$\beta \geq \beta_0$
we have~$x_\alpha \mathrel\varepsilon x_\beta$.

\item
    A subset~$D \subseteq X$ is said to be \Define{dense}\index{uniform space!dense}
    if for each~$\varepsilon \in \Phi$ and~$x \in X$,
    there is a~$y \in D$
    with~$x \mathrel\varepsilon y$.
\end{enumerate}
\spacingfix{}
\end{point}
\begin{point}{20}[dils-uniform-spaces-basics]{Exercise}%
    In the same setting as the previous definition.
    \begin{enumerate}
\item
    Show that equivalence of Cauchy nets is an equivalence relation.
    Show that if~$(x_\alpha)_\alpha$ is a subnet of
    a Cauchy net~$(y_\alpha)_\alpha$,
        that then~$(x_\alpha)_\alpha$ is equivalent to~$(y_\alpha)_\alpha$.
\item
    Prove that if~$(x_\alpha)_\alpha$ and~$(y_\alpha)_\alpha$
        are equivalent Cauchy nets and~$x_\alpha \to x$,
        then also~$y_\alpha \to x$.
\item
    Assume~$(x_\alpha)_\alpha$ is a net with~$x_\alpha \to x$
        and~$x_\alpha \to y$.  Prove~$x = y$ whenever~$X$ is Hausdorff.
\item\label{ex-continuous-preserves-lims}
    Show that if~$f\colon X \to Y$ is a continuous map
    between uniform spaces
        and we have~$x_\alpha \to x$ in~$X$,
        then~$f(x_\alpha) \to f(x)$ in~$Y$.
\item
    Assume~$f$ is a uniformly continuous map between uniform spaces.
    Show that $f$ maps Cauchy nets to Cauchy nets
    and furthermore that~$f$ maps equivalent Cauchy nets to equivalent
        Cauchy nets.
\item\label{ex-cauchy-from-dense-subset}
    Suppose~$D \subseteq X$ is a dense subset.
    Show that for each~$x \in X$
    there is a Cauchy net~$(d_\alpha)_{\alpha \in \Phi}$
    in~$D$ with~$d_\alpha \to x$ in~$X$.
\item
    Assume~$f,g\colon X \to Y$ are continuous maps between uniform spaces
        where~$Y$ is Hausdorff.
Conclude from \ref{ex-continuous-preserves-lims} and
    \ref{ex-cauchy-from-dense-subset}
    that~$f$ and~$g$ are equal whenever they agree on a dense subset of~$X$.
\end{enumerate}
\spacingfix{}
\end{point}
\begin{point}{30}[dils-product-uniformity]{Exercise}%
    Let~$(X_i)_{i \in I}$ be a family of sets with
        uniformities~$(\Phi_i)_{i \in I}$.
    For each~$i_0 \in I$
    and~$\varepsilon \in \Phi_{i_0}$,
    define a relation on
    $\Pi_{i \in I} X_i$ by
    $(x_i)_{i \in I} \mathrel{\hat\varepsilon} (y_i)_{i \in I}
    \iff x_{i_0} \mathrel\varepsilon y_{i_0}$.
    Show~$\{ \hat\varepsilon;\ \varepsilon \in \Phi_i, \ i\in I \}$
    is a subbase for $\Pi_{i \in I} X_i$;
    that the projections~$\pi_i \colon \Pi_{i \in I} X_i \to X_i$
    are uniformly continuous with respect to them
    \emph{and} that they make~$\Pi_{i \in I} X_i$
    into the product of~$(X_i)_{i \in I}$
    in the category of uniform spaces with uniformly continuous maps.
\end{point}
\end{parsec}

\begin{parsec}{1480}%
\begin{point}{10}[blinear-bounded-is-ultranorm]{Proposition}%
Let~$X$ and~$Y$ be right $\scrB$-modules
    with~$\scrB$-valued inner products.
A bounded $\scrB$-linear map~$T\colon X \to Y$
    is uniformly ultranorm continuous.
\begin{point}{20}{Proof}%
Let~$f\colon \scrB \to \C$ be an np-map
    and~$\varepsilon > 0$.
Assume~$\|x-y\|_f \leq \frac{\varepsilon}{\|T\|}$.
By \sref{blinear-inprod-inequality}
we have~$\left<T(x-y), T(x-y)\right> \leq \|T\|^2\left<x-y,x-y\right>$
and so
\begin{equation*}
    \|Tx - Ty\|_f \  =\  f(\left<T(x-y),T(x-y)\right>)^{\frac{1}{2}}
\ \leq\  \|T\| \|x-y\|_f \leq \varepsilon.
\end{equation*}
Thus~$T$ is uniformly continuous w.r.t.~the ultranorm uniformity. \qed
\end{point}
\end{point}
\begin{point}{30}[ultranormcontstruct]{Corollary}%
Let~$\scrB$ be a von Neumann algebra
    and~$X$ a right~$\scrB$-module
    with~$\scrB$-valued inner product.
For~$x_0 \in X$ the maps
\begin{align*}
    (x,y) &\mapsto x+y & x &\mapsto \left<x_0, x\right>
    & b &\mapsto x_0 b \\
    X\times X &\to X & X &\to \scrB & \scrB & \to X
\end{align*}
are all uniformly continuous w.r.t.~the ultranorm uniformity.
\end{point}
\begin{point}{40}[ultranormscalar]{Exercise}%
Let~$\scrB$ be a von Neumann algebra and~$X$ a right~$\scrB$-module
    with~$\scrB$-valued inner product.
Show that~$x \mapsto xb$
    is ultranorm continuous for any~$b \in \scrB$.
\end{point}
\begin{point}{50}[innerprod-ultraweak]{Proposition}%
Let~$X$ be a right~$\scrB$-module
    with~$\scrB$-valued inner product~$[\,\cdot\,,\,\cdot\,]$
    \index{$[\,\cdot\,,\,\cdot\,]$!inner product!$\scrB$-valued}
    for some C$^*$-algebra~$\scrB$.
If~$x_\alpha \to x$ and~$y_\alpha \to y$
    in the ultranorm uniformity,
    for some nets~$(x_\alpha)_{\alpha\in I}$, $(y_\alpha)_{\alpha \in I}$,
    then~$[x_\alpha,y_\alpha] \to [x,y]$ ultraweakly.
\begin{point}{60}{Proof}%
Let~$\varepsilon > 0$ and np-map~$f\colon \scrB \to \C$ be given.
Note that~$\|y_\alpha\|_f \to \|y\|_f$;
    indeed by the reverse triangle inequality 
    $|\|y_\alpha\|_f - \|y\|_f| \leq \| y_\alpha - y\|_f \to 0$.
Thus we can find~$\alpha_1$ such that~$\|y_\alpha\|_f \leq \|y\|_f + 1$
    for~$\alpha \geq \alpha_1$.

We can find~$\alpha_2$
    such that for all~$\alpha \geq \alpha_2$
    we have $\|x-x_\alpha\|_f \leq \frac{1}{2}\varepsilon (\|y\|_f + 1)^{-1}$
    and
    $\|y-y_\alpha\|_f \leq \frac{1}{2}\varepsilon (\|x\|_f + 1)^{-1}$.
Pick~$\alpha_0 \geq \alpha_1,\alpha_2$.
Then for~$\alpha \geq \alpha_0$ we have
\begin{align*}
    |f([x_\alpha,y_\alpha] - [x,y])| 
    & \ = \ |[x_\alpha - x, y_\alpha]_f + [x,y_\alpha - y]_f| \\
    & \ \leq \ 
         \|x_\alpha - x\|_f \| y_\alpha\|_f + \|x\|_f\|y_\alpha - y\|_f
                && \text{ by \sref{module-CS}}\\
    & \ \leq \ 
         \|x_\alpha - x\|_f (\|y\|_f+1) + (\|x\|_f+1)\|y_\alpha - y\|_f \\
    & \ \leq \ 
         \varepsilon.
\end{align*}
Thus~$[x_\alpha,y_\alpha] \to [x,y]$ ultraweakly. \qed
\end{point}
\end{point}
\begin{point}{70}[hilbmod-denseordersep]{Corollary}%
Let~$X$ be a Hilbert~$\scrB$-module
    together with an ultranorm-dense subset~$D \subseteq X$.
The vector states from~$D$ are order separating:
for each~$T \in \scrB^a(X)$ we have
$T \geq 0$ if and only if~$\left<x,Tx\right> \geq 0$
    for all~$x \in D$.  Cf.~\sref{hilbmod-ordersep}.
\end{point}
\end{parsec}

\begin{parsec}{1490}%
\begin{point}{10}[dfn-selfdual-basis]{Definition}
    Let~$X$ be pre-Hilbert $\scrB$-module
        for a von Neumann algebra~$\scrB$.
    Let~$E \subseteq X$ be some subset.
    \begin{enumerate}
        \item $E$ is \Define{orthogonal}\index{orthogonal} if~$\left<e,d\right> = 0$
            for~$e,d \in E$ with~$e \neq d$.
        \item $E$ is \Define{orthonormal}\index{orthonormal} if additionally
                    $\left<e,e\right>$ is a non-zero projection
                    for~$e \in E$.
        \item A family $(b_e)_{e \in E}$ in~$\scrB$
                is said to be \Define{$\ell^2$-summable}
                if the partial sums of~$\sum_{e \in E} b_e^*b_e$
                are bounded.
        \item $E$ is an \Define{(orthonormal) basis}\index{basis@(orthonormal) basis} if
                it is orthonormal and additionally
            \begin{enumerate}
                \item for each~$x \in X$ we have
                    \begin{equation*}
                        x \ =\  \sum_{e \in E} e\left<e,x\right>,
                    \end{equation*}
                    where the sum converges in the ultranorm uniformity
                    \emph{and}
                \item  the sum~$\sum_{e \in E} eb_e$
                        converges in the ultranorm uniformity
                    for any~$\ell^2$-summable $(b_e)_{e \in E}$
                        in~$\scrB$.
            \end{enumerate}
    \end{enumerate}
\spacingfix{}
\begin{point}{20}{Notation}%
As with Hilbert spaces, it will sometimes be more convenient
    to view an orthonormal basis as a sequence/family~$(e_i)_{i \in I}$,
    instead of a subset. This should cause no confusion.
\end{point}
\begin{point}{21}{Remark}%
The first assumption `(a)' in the definition of an orthonormal basis~$E$
    requires that every element~$x$ can be reconstructed
    from its coefficients~$(\left<e, x\right>)_{e \in E}$.
This assumption fails, for instance, if there is a non-zero element
    orthogonal to all basis elements.
The second assumption `(b)' requires that behind
    any a reasonable family of coefficients,
    there is a corresponding element.
This second assumption fails if the module is not ultranorm complete,
    see~\sref{dils-selfdual}.
\end{point}
\begin{point}{22}{Beware}%
In the literature, there are several different definitions of an
    orthonormal basis
    for a Hilbert~$\scrB$-module~$X$.
\begin{enumerate}
\item
Some~\cite{landi2012orthogonal,manuilov2000hilbertc}
    define an orthonormal basis as an orthonormal subset~$E \subseteq X$
        such that~$x = \sum_{e \in E} e \left<e, x\right>$
    for all~$x \in X$
    where the sum converges in the \emph{norm}.
This notion of basis is independent of ours:
        there is a basis in our sense
         that is not a basis
        in their sense ($\{ \ketbra{0}{0}, \ketbra{1}{1}, \ldots\}$
        in~$\scrB(\ell^2)$)
        and vice versa (the standard Hilbert module~$H_{\scrB(\ell^2)}$
            defined in~\cite[\S3]{landi2012orthogonal} has a basis
            in their sense, but not a basis in ours
            as it is not self dual, cf.~\sref{dils-selfdual}.)
\item
Others such as~\cite[lemma~8.5.23]{blecher2004operator}
    define an orthonormal basis
        only when~$X$ is a self-dual (\sref{dils-selfdual})
        Hilbert C$^*$-module
            over a von Neumann algebra.
    In this case an orthonormal basis
        is defined
        as an orthonormal subset~$E \subseteq X$
    such that~$x = \sum_{e \in E} e \left<e, x \right>$
        for all~$x \in X$
        where the sum converges in the~$w^*$-topology of~$X$
        (with respect to the unique Banach space predual
    of~$X$ for which the inner product is
    separately~$w^*$-continuous,
        see \cite[lemma 8.5.4]{blecher2004operator}.)
This notion is equivalent to ours, but only defined
    for self-dual modules over von Neumann algebras.
\end{enumerate}
\spacingfix{}
\end{point}
\end{point}
\begin{point}{30}[mod-projelabs]{Exercise}%
Let~$X$ be a pre-Hilbert $\scrB$-module for a C$^*$-algebra~$\scrB$.
Show that we have~$ e \langle e,e \rangle = e$
    for any~$e \in X$
   where~$\langle e, e\rangle$ is a projection.

    (Hint: consider~$\| e(1-\langle e, e \rangle) \|^2$;
                 cf.~\cite[thm.~3.12]{paschke}.)
\end{point}
\begin{point}{40}[mod-parseval]{Exercise}
Let~$X$ be a pre-Hilbert~$\scrB$-module for a von Neumann algebra~$\scrB$
    with an orthonormal basis~$E \subseteq X$.
    Show \Define{Parseval's identity}\index{Parseval's identity} holds, that is:
\begin{equation*}
    \langle x, x \rangle \ =\ 
    \sum_{e \in E} \langle x, e\rangle \langle e, x\rangle,
\end{equation*}
where the sum converges ultraweakly. \cite{paschke}
(Hint: use~\sref{innerprod-ultraweak}.)
\end{point}
\begin{point}{50}[dils-selfdual]{Theorem}%
    Assume~$\scrB$ is a von Neumann algebra.
    For any pre-Hilbert~$\scrB$-module~$X$,
        the following conditions are equivalent.
        \begin{enumerate}
            \item $X$ is self dual.
            \item $X$ is ultranorm complete.
            \item Every norm-bounded ultranorm-Cauchy net in~$X$ converges.
            \item $X$ has an orthonormal basis.
        \end{enumerate}
        \spacingfix{}
\begin{point}{60}{Proof}%
We will prove $1 \Rightarrow 3 \Rightarrow 4 \Rightarrow 2 \Rightarrow 4 \Rightarrow 1$.
The equivalence~$1 \Leftrightarrow 4$
    was already shown by Paschke in \cite[thm.~3.12]{paschke},
    but we give a different proof.
The equivalence~$1 \Leftrightarrow 3$ turns out
    proven by Frank already \cite[thm.~3.2]{frank1990self}
    in a different way.
The final equivalence~$1 \Leftrightarrow 2$ is new, as far as we know.
\begin{point}{70}{1 $\Rightarrow$ 3, self-dual~$X$ are bounded complete}%
    Let~$X$ be a self-dual pre-Hilbert~$\scrB$-module.
    To show norm-bounded ultranorm completeness,
     assume~$(x_\alpha)_\alpha$ is an ultranorm-Cauchy net in~$X$
        with~$\|x_\alpha\| \leq B$ for some~$B \geq 0$.
     For the moment pick an~$y \in X$.
By \sref{module-CS} we
    have~$f(\left<x_\alpha,y\right>\left<y,x_\alpha\right>) \leq \|y\|^2
        f(\left<x_\alpha,x_\alpha\right>)$ for any~np-map~$f\colon \scrB \to \C$
        and so~$\left<y,x_\alpha\right>$ is an ultrastrong-Cauchy net in~$\scrB$.
Von Neumann algebras are ultrastrongly complete, see \sref{vn-complete},
so we can define~$\tau(y) = (\uslim_\alpha \left<y, x_\alpha\right>)^*$,
    where~\Define{$\uslim$}\index{uslim@$\uslim$} denotes
    the ultrastrong limit.
As addition and multiplication by a fixed element
are ultrastrongly continuous (see \sref{mult-uws-cont})
    $\tau$ is~$\scrB$-linear.
We will show~$\tau$ is bounded.
For each~$\alpha$ we have~$\left<x_\alpha, y\right>\left<y, x_\alpha\right>
        \leq \|y\|^2\left<x_\alpha, x_\alpha\right> \leq \|y\|^2B^2$ and so
\begin{equation*}
    \tau(y)\tau(y)^*
    \ \overset{\sref{usconv}}{=} \ 
    \uwlim_\alpha \left<x_\alpha,y\right>\left<y,x_\alpha\right>
    \ \leq \ \|y\|^2B^2,
\end{equation*}
where \Define{$\uwlim$}\index{uwlim@$\uwlim$}
    denotes the ultraweak limit.
Thus~$\|\tau(y)\|^2 = \|\tau(y)\tau(y)^*\| \leq \|y\|^2B^2$.
So~$\tau$ is bounded.
By self-duality, there is a~$t \in X$ such that~$\tau(y) = \left<t,y\right>$
for all~$y \in X$.
And so for each np-map~$f\colon \scrB \to \C$ we have
\begin{equation}\label{selfdual-bramA}
    \|\left<y, t-x_\alpha\right>\|_f \ =\ 
    \|\left<y, t\right>-\left<y,x_\alpha\right>\|_f \ =\ 
    \|\tau(y)^*-\left<y,x_\alpha\right>\|_f\  \rightarrow \  0.
\end{equation}
In~\sref{module-innerprod-state}
    we introduced the inner product~$\left<\,\cdot\,,\,\cdot\,\right>_f$
    defined by~$\left<x,y\right>_f \equiv f(\left<x,y\right>)$.
From Kadison's inequality \sref{omega-norm-basic}
it follows~$|f(a)|^2 \leq f(a^*a)f(1)$
and so
\begin{equation}\label{selfdual-bramB}
    |\left<y,t-x_\alpha\right>_f|^2
    \ \leq \ f( \left<t-x_\alpha,y\right> \left<y,t-x_\alpha\right>) f(1)
    \ =\  \| \left<y,t-x_\alpha\right> \|_f^2 f(1).
\end{equation}
Combining~\eqref{selfdual-bramA} and~\eqref{selfdual-bramB},
    we see~$|\left<y,t-x_\alpha\right>_f| \rightarrow 0$.
Let~$\varepsilon > 0$ be given.
As~$\|x_\alpha\|_f$ is norm bounded
we can find~$B > 0$ such that~$\|x_\alpha\|_f \leq B$
for all~$\alpha$.
As~$(x_\alpha)_\alpha$ is ultranorm Cauchy,
we can pick~$\alpha$ such that~$\|x_\beta - x_\alpha\|_f \leq (\|t\|_f + B)^{-1}
\frac{\varepsilon}{2}$ for all~$\beta \geq \alpha$.
Find~$\beta$ such that~$|\left<t-x_\alpha,t-x_\beta\right>_f|
    \leq \frac{\varepsilon}{2}$.
Putting it all together:
\begin{align*}
    \left<t-x_\alpha,t-x_\alpha\right>_f
    & \ \leq\  |\left<t-x_\alpha,t-x_\beta\right>_f| \,
    +\, |\left<t-x_\alpha,x_\beta-x_\alpha\right>_f| \\
    & \ \leq\  \frac{\varepsilon}{2}
    + \| t- x_\alpha \|_f \|x_\beta - x_\alpha\|_f\\
    & \ \leq\  \frac{\varepsilon}{2}
    + (\| t\|_f + B) \,\|x_\beta - x_\alpha\|_f \ \leq\  \varepsilon.
\end{align*}
Thus~$x_\alpha$ converges ultranorm to~$t$.
So~$X$ is ultranorm bounded complete.
\end{point}
\begin{point}{80}[selfdual-bcompl-then-basis]{3
    $\Rightarrow$ 4, bounded complete~$X$ has an orthonormal basis}
    Assume norm-bounded ultranorm-Cauchy nets converge in~$X$.
By Zorn's lemma, there is a maximal orthonormal subset~$E \subseteq X$.
We will show~$E$ is a basis.

Let~$(b_e)_{e \in E}$ be an~$E$-tuple in~$\scrB$
such that the partial sums of~$\sum_{e \in E} b_e^*b_e$ are bounded.
We want to show~$\sum_{e \in E} e b_e$ converges in the ultranorm
    uniformity.
As the summands of~$\sum_{e \in E}b_e^*b_e$ are positive,
    the partial sums form a directed bounded net
    that converges ultrastrongly
    to its supremum~$\sum_{e \in E} b_e^*b_e$.
By~\sref{mod-projelabs},
    we have~$e \left<e,e\right> = e$ for each~$e \in E$
    and so~$\left<x,e\right>\left<e,e\right> = \left<x,e\right>$
    for every~$x \in X$.
Thus for any finite subset~$S \subseteq E$
\begin{equation*}
    \Bigl<\sum_{e \in S} eb_e, \sum_{e \in S} eb_e\Bigr>
    \ =\  \sum_{e \in S} b_e^* \left<e,e\right>b_e
    \ \leq\  \sum_{e \in S} b_e^*b_e.
\end{equation*}
To show~$\sum_e eb_e$ converges ultranorm,
    pick any np-map~$f\colon \scrB \to \C$.
As~$\sum_{e \in E} b^*_eb_e$ converges ultraweakly,
we have~$\sum_{e \in E} f(b^*_eb_e) < \infty$.
Thus the tails of the series~$\sum_{e \in E - S} f(b_e^*b_e)$ tend to zero, hence for finite~$S,T\subseteq E$ we know
\begin{equation*}
    \bigl\| \sum_{e \in S} eb_e - \sum_{e\in T} eb_e \bigr\|_f^2
        \ =\  \bigl\| \sum_{e \in S\Delta T} eb_e  \bigr\|_f^2
        \ \leq \ \sum_{e \in S\Delta T} f(b_e^*b_e).
\end{equation*}
The latter vanishes as~$S \cap T$ grows.
Thus the partial sums of ~$\sum_{e \in E} eb_e$ are ultranorm Cauchy
and norm-bounded by~$\| \sum_{e \in E} b_e^*b_e \|$,
so~$\sum_{e \in E} eb_e$ converges in the ultranorm uniformity.

Pick~$x \in X$.
We have to show~$\sum_{e \in E} x\left<e,x\right>$
converges in the ultranorm uniformity to~$x$.
For any finite subset~$S \subseteq E$ we have
\begin{equation*}
    0\  \leq\  \Bigl< x - \sum_{e \in S} e\left<e,x\right>,
        x - \sum_{e \in S} e\left<e,x\right> \Bigr>
       \  =\  \left<x,x\right> - \sum_{e \in S} \left<x,e\right>\left<e,x\right>.
\end{equation*}
Rearranging we find \Define{Bessel's inequality}\index{Bessel's inequality}:
\begin{equation*}
\sum_{e \in S} \left<x,e\right>\left<e,x\right>
   \  \leq \ \left<x,x\right>.
\end{equation*}
Hence the $E$-tuple~$(\left<e,x\right>)_e$
is~$\ell^2$-summable
and so~$\sum_{e \in E} e \left<e,x\right>$
converges in the ultranorm uniformity.
Consider~$x' = x - \sum_{e \in E} e\left<e,x\right>$.
We are done if we can show~$x' = 0$.

To see~$x'=0$, we
first show that~$X$ has \Define{polar decomposition}\index{polar decomposition}:
    for each~$y \in X$, there is an~$u\in X$
    with~$\left<u,u\right> = \ceil{\left<y,y\right>}$
        and~$y=u\left<y,y\right>^{\frac{1}{2}}$.
By \sref{approximate-pseudoinverse}
    the positive element~$\left<y,y\right>^{\frac{1}{2}}$
    has an approximate pseudoinverse~$h_1, h_2, \ldots$
    --- that is:
    $h_n \left<y,y\right>^{\frac{1}{2}} =
         \left<y,y\right>^{\frac{1}{2}} h_n
         = \ceil{h_n}$
         and~$\sum_n h_n \left<y,y\right>^{\frac{1}{2}} = \sum_n \ceil{h_n}
                    = \ceil{\left<y,y\right>}$. Note
\begin{equation*}
    \left<yh_n,yh_m\right> \ = \ 
        (h_n \left<y,y\right>^{\frac{1}{2}})
        (\left<y,y\right>^{\frac{1}{2}} h_m) \ =\  \ceil{h_n}\ceil{h_m}
\end{equation*}
and so using that~$\ceil{h_1}, \ceil{h_2}, \ldots$ are pairwise orthogonal,
    we see
\begin{equation*}
    \bigl< \sum_{n=1}^N y h_n, \sum_{n=1}^N y h_n \bigr> 
    \ =\  \sum_{n=1}^N \sum_{m=1}^N \ceil{h_n} \ceil{h_m}
                \ =\  \sum_{n=1}^N \ceil{h_n}
                \ \leq\  \ceil{\left<y,y\right>}.
\end{equation*}
Thus~$( \sum_{n=1}^N y h_n)_N$ forms an ultranorm-Cauchy net
    norm bounded by~$1$.
By assumption this bounded net converges:
    write~$u \equiv \sum_n yh_n$.
(So in the special case~$X = \scrB$,
    we have~$u = y / \langle y,y \rangle^{\frac{1}{2}}$,
    see~\sref{proto-douglas}.)
We derive
\begin{align*}
u \left<y,y\right>^{\frac{1}{2}}
&\ =\  \sum_n yh_n \left<y,y\right>^{\frac{1}{2}} 
    && \text{by dfn.~$u$ and \sref{ultranormscalar}}
\\
&\ =\  \sum_n y \ceil{h_n} 
    && \text{by dfn.~$h_n$}\\
    &\ =\  y \sum_n\ceil{h_n}
&& \text{by \sref{ultranormcontstruct}} \\
&\ =\  y \ceil{\left<y,y\right>} && \text{by dfn.~$h_n$}\\
&\ =\  y,
\end{align*}
where the last equality 
follows from~$y (1- \ceil{\left<y,y\right>})=0$,
    which is justified by
\begin{equation*}
    \bigl<
        y (1- \ceil{\left<y,y\right>}),
        y (1- \ceil{\left<y,y\right>}) \bigr>
        \ = \ 
        (1- \ceil{\left<y,y\right>}) \left<y,y\right>
        (1- \ceil{\left<y,y\right>})
        \ = \  0.
\end{equation*}
This finishes the proof of polar decomposition.

Recall we want to
show~$x' \equiv x - \sum_{e \in E} e\left<e,x\right>= 0$.
Reasoning towards contradiction, assume~$x' \neq 0$.
By polar decomposition~$x' = u\left<x',x'\right>^{\frac{1}{2}}$
for some~$u \in X$ with~$\left<u,u\right> = \ceil{\left<x',x'\right>}$.
Note~$\left<u,u\right> \neq 0$
    for otherwise~$x' = 0$, quod non.
For any~$e_0 \in E$, we have
\begin{equation*}
   \left<e_0, x'\right> \ =\ 
   \bigl<e_0, x-\sum_{e\in E} e\left<e,x\right>\bigr> 
   \ \overset{\sref{ultranormcontstruct}}{=} \ 
   \left<e_0, x\right> - \sum_{e \in E}\left<e_0, e\right>\left<e,x\right> \ =\  0
\end{equation*}
and so~$0 = \left<e,x'\right> = \left<e, u\right>\left<x',x'\right>^{\frac{1}{2}}$ for any~$e \in E$.
Thus
\begin{equation*}
        \left<e, u\right> \ =\ 
        \left<e, u\right> \left<u,u\right> \ = \ 
        \left<e, u\right> \ceil{\left<x',x'\right>}
        \ \overset{\sref{mult-cancellation}}{=}\  0.
\end{equation*}
Thus~$E \cup \{ u \}$ contradicts the maximality of~$E$.
Hence~$x'=0$ and so~$x = \sum_{e \in E} e\left<e,x\right>$ as desired.
\end{point}
\begin{point}{90}{4 $\Rightarrow$ 2, basis implies completeness}%
    Assume~$X$ has an orthonormal basis~$E$.
    Let~$(x_\alpha)_\alpha$ be an ultranorm-Cauchy net.
    Pick~$e \in E$.
    By \sref{ultranormcontstruct}
    $(\left<e, x_\alpha\right>)_\alpha$ is an ultrastrong-Cauchy net.
    As von Neumann algebras are ultrastrongly complete
    (see \sref{vn-complete})
        there is a~$b_e \in \scrB$
            to which~$\left<e,x_\alpha\right>$ converges ultrastrongly.
    We will show~$x_\alpha$ converges in the ultranorm uniformity
        to~$\sum_{e \in E} e b_e$.

We will first show that~$\sum_{e \in E} e b_e$ converges.
As before, it is sufficient to show~$\sum_{e \in E} b_e^*b_e$ is bounded.
For this, it suffices to show~$\sum_{e\in E} f( b_e^*b_e ) $ is bounded
        for each np-map~$f\colon \scrB \to \C$,
    as ultraweakly-bounded nets are norm bounded,
    see \sref{ultraweakly-bounded-implies-bounded}.
For any~$e \in E$
    we have~$\|b_e\|_f = \lim_\alpha \|\left<e,x_\alpha\right>\|_f$
    and so~$f(b_e^*b_e) =\lim_\alpha f(\left<x_\alpha,e\right> \left<e, x_\alpha\right>)$.
    Pick a finite subset~$S \subseteq E$. Then we have
\begin{alignat*}{3}
    f\bigl(\sum_{e \in S} b_e^*b_e\bigr)
    &\ =\  \sum_{e \in S} \lim_\alpha 
     f(\left<x_\alpha,e\right> \left<e, x_\alpha\right>) \\
     &\  = \ \lim_\alpha 
            f\bigl(\sum_{e \in S}\left<x_\alpha,e\right>
                \left<e, x_\alpha\right>\bigr) \\
            &\ \leq\  f(\left<x_\alpha,x_\alpha\right>)
                &\qquad&\text{by Bessel's inequality.}
 \end{alignat*}
As~$x_\alpha$ is ultranorm Cauchy,
the net~$f(\left<x_\alpha,x_\alpha\right>)$
    must be Cauchy and so bounded.
Hence~$\sum_{e \in E} e b_e$ converges in the ultranorm uniformity.

To prove~$x_\alpha$ converges to~$\sum_{e \in E} e b_e$,
    pick any~$\varepsilon > 0$
        and np-map $f\colon \scrB \to \C$.
We want to find~$\alpha_0$
    such that~$\| x_\alpha - \sum_{e \in E} e b_e\|_f \leq \varepsilon$
        for all~$\alpha \geq \alpha_0$.
As~$(x_\alpha)_\alpha$ is ultranorm Cauchy,
we can pick an~$\alpha_0$
such that~$\| x_\alpha - x_\beta \|_f \leq \frac{1}{2\sqrt{2}} \varepsilon$
for all~$\alpha,\beta \geq \alpha_0$.
As~$x_\alpha = \sum_{e \in E} e \left<e,x_\alpha\right>$
we find with Parseval's identity (\sref{mod-parseval}) that
\begin{equation*}
    \bigl\| x_\alpha - \sum_{e \in E} eb_e \bigr\|_f^2
    \ =\  \bigl\| \sum_{e \in E} e(\left<e, x_\alpha\right> - b_e) \bigr\|_f^2
    \ =\  \sum_{e \in E}\| \left<e, x_\alpha\right> - b_e \|_f^2.
\end{equation*}
Take a finite subset~$S \subseteq E$
such that
    $\sum_{e \in E - S}\| \left<e, x_\alpha\right> - b_e \|_f^2 \leq
        \frac{1}{2}\varepsilon^2$.
If we can also show
    $\sum_{e \in  S}\| \left<e, x_\alpha\right> - b_e \|_f^2 \leq
        \frac{1}{2}\varepsilon^2$,
        then we are done.
For any~$\beta \geq \alpha$
    we have using the triangle inequality of~$\|\,\cdot\,\|_f$
        on~$X$
\begin{align*}
    \Bigl(\sum_{e \in S}\| \left<e, x_\alpha\right>
        - b_e \|_f^2\Bigr)^{\frac{1}{2}}
        &\ =\  \bigl\|\sum_{e \in S} e\left<e, x_\alpha\right> - eb_e \bigr\|_f\\
        &\ \leq\  \bigl\|\sum_{e \in S} e\left<e, x_\alpha-x_\beta\right> \bigr\|_f\
        +\bigl\|\sum_{e \in S} e\left<e, x_\beta\right> - eb_e \bigr\|_f\\
        &\ =\  \Bigl( \sum_{e \in S}\|\left<e, x_\alpha-x_\beta\right> \|_f^2 \Bigr)^{\frac{1}{2}}
        +\Bigl(\sum_{e \in S}\| \left<e, x_\beta\right> - b_e \|_f^2
        \Bigr)^{\frac{1}{2}}.
\end{align*}
We want to show the previous is bounded by~$\frac{1}{\sqrt{2}}\varepsilon$.
We bound the two terms separately.
For the first term using Bessel's
inequality and our choice of~$\alpha_0$
we have
$\sum_{e \in S}\|\left<e, x_\alpha-x_\beta\right> \|_f^2 
\leq \|x_\alpha - x_\beta\|^2_f \leq \bigl(\frac{1}{2\sqrt{2}}\varepsilon\bigr)^2.$
As~$\|\left<e,x_\beta\right> - b_e\|_f$
vanishes for each~$e \in S$ and~$S$ is finite,
we can find sufficiently large~$\beta$
such that the right term is also bounded by~$\frac{1}{2\sqrt{2}}\varepsilon$.
Thus~$x_\alpha$ converges in the ultranorm uniformity
to~$\sum_{e \in E} e b_e$, as desired.
\end{point}
\begin{point}{100}{2 $\Rightarrow$ 4, completeness implies basis}%
    Follows from~3 $\Rightarrow$ 4
        as~2 $\Rightarrow$ 3 is trivial.
\end{point}
\begin{point}{110}{4 $\Rightarrow$ 1, basis implies self-duality}%
    Assume~$X$ has an orthonormal basis~$E$.
    Suppose~$\tau\colon X \to \scrB$ is a bounded $\scrB$-linear map.
    We want to show there is some~$t \in X$
    such that~$\tau(x) = \left<t,x\right>$ for all~$x \in X$.
    Pick any~$x \in X$.
    By assumption~$\sum_{e \in E} e  \left<e,x\right>$
        converges in the ultranorm uniformity to~$x$.
As by \sref{blinear-bounded-is-ultranorm} $\tau$ is ultranorm continuous
        we see
        \begin{equation*}
            \tau(x) \ =\  \sum_{e \in E} \tau(e) \left<e,x\right>
                \ =\  \sum_{e \in E} \left<e \tau(e)^*,x\right>.
        \end{equation*}
    Thus if we can show~$\sum_e e \tau(e)^*$ converges in the ultranorm,
        we are done. By assumption this is
            the case if~$(\tau(e)^*)_{e \in E} $ is $\ell^2$-summable
            --- that is: the partial sums of~$\sum_{e \in E} \tau(e)\tau(e)^*$
                must be bounded.
        We will show
        that for finite~$S \subseteq E$
        we have
        $\|\sum_{e \in S} e \tau(e)^* \|\leq \|\tau\|$,
        which is sufficient as then~$ \| \sum_{e \in S}  \tau(e)\tau(e)^* \|
        =\|\sum_{e \in S} e \tau(e)^* \|^2  \leq \|\tau\|^2$.
For the moment pick an arbitrary~$x \in X$.
By Bessel's inequality
we have~$
\| \sum_{e\in S} e \left<e,x\right>\|^2
    = \sum_{e\in S} \|\left<x,e\right>\left<e,x\right>\|
    \leq \|x\|^2$.
Thus
    \begin{equation*}
        \bigl\|\bigl< \sum_{e \in S} e \tau(e)^*,x\bigr>\bigr\|
        \ =\  \bigl\| \tau \bigl(\sum_{e \in S} e\left<e,x\right>\bigr) \bigr\|
        \leq \|\tau\| \|x\|.
    \end{equation*}
Substituting~$\sum_{e \in S} e \tau(e)^*$ for~$x$, we find
\begin{equation*}
    \| \sum_{e\in S}e \tau(e)^*\|^2 \ =\ 
    \bigl\|\bigl< \sum_{e \in S} e \tau(e)^*,\sum_{e \in S} e\tau(e)^*\bigr>\bigr\|
    \ \leq \ \|\tau\| \bigl\| \sum_{e\in S}e \tau(e)^*\bigr\|.
\end{equation*}
If~$ \| \sum_{e\in S}e \tau(e)^*\| = 0$,
then we are done.
Otherwise, divide both sides of the previous equation
by~$ \| \sum_{e\in S}e \tau(e)^*\|$
to find~$ \| \sum_{e\in S}e \tau(e)^*\| \leq \|\tau\|$ as desired. \qed
\end{point}
\end{point}
\end{point}
\end{parsec}

\subsection{Self-dual completion}
\begin{parsec}{1500}%
\begin{point}{10}%
In~\sref{prop-complete-into-hilbert-space}
    we saw that every vector space with inner product
    can be completed into a Hilbert space.
We continue with the generalization to~Hilbert C$^*$-modules.
Paschke proves a similar result~\cite[thm.~3.2]{paschke}
    in a completely different way.
\end{point}
\begin{point}{20}[dils-completion]{Theorem}%
Let~$\scrB$ be a von Neumann algebra
and~$V$ a right~$\scrB$-module with~$\scrB$-valued
    inner product~$[\,\cdot\,,\,\cdot\,]$.
    \index{$[\,\cdot\,,\,\cdot\,]$!inner product!$\scrB$-valued}
There is a self-dual Hilbert $\scrB$-module~$X$
    together with a bounded~$\scrB$-linear~$\eta\colon V \to X$
    such that \begin{inparaenum}
    \item $[v,w] = \left<\eta (v),\eta(w)\right>$
    and \item the image of~$\eta$ is ultranorm dense in~$X$.
    \end{inparaenum}
\begin{point}{30}{Proof}%
Before we start the proof proper,
    we give an overview.
\begin{point}{40}{Overview}
By~\sref{dils-selfdual} and~\sref{innerprod-ultraweak}
    one would expect that
    one can simply extend the inner product of~$V$
    to its ultranorm completion~$\overline{V}$
    by
    \begin{equation}
        [(x_\alpha)_\alpha, (y_\alpha)_\alpha] \ =\  \uwlim_\alpha [x_\alpha,y_\alpha], \label{eq-definition-inprod-on-completion}
    \end{equation}
    which will turn $\overline{V}$ into a self-dual Hilbert~$\scrB$-module.
It will turn out that~$\overline{V}$
    is indeed a self-dual Hilbert $\scrB$-module
    with this inner product,
    but it isn't clear at all  how to show directly that the ultraweak limit
    in \eqref{eq-definition-inprod-on-completion} converges.

We will sketch an indirect construction,
    before delving into the details.
As with a metric completion,
    $\overline{V}$ consists of equivalence classes of ultranorm
        Cauchy nets in~$V$.
Let~$V_0\subseteq \overline{V}$
    denote the equivalence classes of constant Cauchy nets.
$V_0$ is a right $\scrB$-module with~$\scrB$-valued inner product
    using the operations of~$V$.
Let~$V_1\subseteq \overline{V}$ denote
    the limits of norm-bounded ultranorm Cauchy nets over~$V_0$.
Using the norm boundedness we can extend
    the module structure and inner product of~$V_0$ to~$V_1$.
If all norm-bounded ultranorm Cauchy-nets over~$V_1$ would
    converge in~$V_1$,
    then~$V_1$ would be a self-dual Hilbert~$\scrB$-module
    and~$V_1$ would even be ultranorm complete.
Later, when we have established that~$\overline{V}$ is a Hilbert~$\scrB$-module,
    we will find out that in fact~$\overline{V} \subseteq V_1$
    due to~\sref{kaplansky-hilbmod}
    ---
    i.e.~$V_1$ is already ultranorm complete.
However, in a most frustrating state of affairs,
    we do not know how to prove this directly
    without first proving that~$\overline{V}$ is a Hilbert~$\scrB$-module.
Instead, we repeat: let~$V_2 \subseteq \overline{V}$
    denote the limits of norm-bounded ultranorm Cauchy nets over~$V_1$.
We have to go further:
    for any~$n \in \N$
    define~$V_{n+1} = \sigma(V_n)$
    where $\sigma(U)\subseteq \overline{V}$
    denotes the limits of norm-bounded ultranorm Cauchy over~$U$.
\begin{equation*}
   V_0 \ \subseteq \ V_1 \ \subseteq\  V_2\  \subseteq\  V_3 \ \subseteq \ \cdots\ \subseteq \ \bigcup_{n \in \N} V_n.
\end{equation*}
Will we be able to show that~$\bigcup_n V_n$ is finally ultranorm complete?
Not yet.
We have to go even further.
Define~$V_\omega = \sigma(\bigcup_{n \in \N} V_n)$,
$V_{\omega+1} = \sigma(V_\omega)$
and so on:
for any ordinal number~$\alpha > 0$,
    define~$V_\alpha = \sigma(\bigcup_{\beta < \alpha} V_\beta)$.
Now~$V_\alpha$
    is an ascending chain of subsets of~$\overline{V}$.
    Thus~$V_{\alpha_0} = V_{{\alpha_0}+1}$
for some~$\alpha_0 \leq |2^{\overline{V}}|$
(for otherwise~$|2^{\overline{V}}| \leq |\overline{V}|$).
Thus norm-bounded ultranorm Cauchy nets in~$V_{\alpha_0}$
    must converge.
We will see that at each step we can extend the module structure
    and inner product and with those~$V_{\alpha_0}$ is self-dual.
    Thus~$V_{\alpha_0}$ is ultranorm complete
    and actually~$V_{\alpha_0} = \overline{V}$.
\begin{equation*}
   V_0 \subseteq V_1 \subseteq \cdots \subseteq V_\omega \subseteq
   V_{\omega+1} \subseteq \cdots \subseteq V_{2\omega} \subseteq
       \ \cdots\ 
       \ \cdots\  \subseteq
   V_{\alpha_0} = \overline{V}
\end{equation*}
To avoid requiring familiarity with transfinite induction,
    we will phrase the proof using Zorn's lemma instead.
\end{point}
\begin{point}{50}{$\overline{V}$: fast nets}%
As we need some details in its construction,
    we will explicitly define~$\overline{V}$,
    the completion of~$V$ in the ultranorm uniformity
    using Cauchy nets.
There are other ways to construct a completion of a uniform space,
    see for instance~\cite[thm.~39.12]{willard}.

Let~$\Phi$ denote the set of entourages of the ultranorm
    uniformity on~$V$; for its definition see~\sref{dils-ultranorm}.
$\Phi$ is a filter and thus can be used as index set for a net
    using reverse inclusion as order.
    (Thus~$\varepsilon \geq \delta
    \ \Leftrightarrow\ \varepsilon \subseteq \delta$.)
We say a net~$(x_\alpha)_{\alpha \in \Phi}$  indexed by entourages
is \Define{fast}\index{Cauchy net!fast} if for every~$\varepsilon \in \Phi$
        and~$\alpha,\beta \geq \varepsilon$
        we have~$x_\alpha \mathrel{\varepsilon^2} x_\beta$.

Every Cauchy net is equivalent to a fast one, but this is not as
    evident as in the metric case.
Let~$(x_\alpha)_{\alpha \in I}$
    be an arbitrary Cauchy net in~$V$.
By definition we can find~$\alpha_\varepsilon \in I$
    such that for each~$\varepsilon \in \Phi$
    and~$\alpha,\beta \geq \alpha_\varepsilon$
    we have~$x_\alpha \mathrel{\varepsilon} x_\beta$.
Unfortunately~$\alpha_{(\,\cdot\,)}$ need not be order preserving
    and so~$(x_{\alpha_{\varepsilon}})_{\varepsilon\in\Phi}$
    need not be a subnet of~$(x_\alpha)_{\alpha \in I}$.
However, we claim the net~$(x_{\alpha_{\varepsilon}})_{\varepsilon\in\Phi}$
    is Cauchy, fast and equivalent to~$(x_{\alpha})_{\alpha\in I}$.
To show it's Cauchy and fast, pick any~$\varepsilon \in \Phi$
    and~$\zeta,\xi \geq \varepsilon$
    (that is: $\zeta,\xi \subseteq \varepsilon$).
As~$I$ is directed, we can find~$\beta\in I$
    with~$\beta \geq \alpha_\zeta, \alpha_\xi$.
By definition of~$\alpha_{(\,\cdot\,)}$
    we have~$x_{\alpha_\zeta} \mathrel\zeta x_{\beta}
        \mathrel\xi x_{\alpha_\xi}$
        and so~$x_{\alpha_\zeta} \mathrel{\varepsilon^2}
                x_{\alpha_\xi}$, as desired.
To show equivalence, assume~$\varepsilon \in \Phi$ is given.
Assume~$\delta \geq \nicefrac{\varepsilon}{2}$
and~$\beta \geq \alpha_{\nicefrac{\varepsilon}{2}}$.
There is a~$\gamma \in I$ with~$\gamma \geq \alpha_\delta,\beta$.
Then~$x_\beta \mathrel{\nicefrac{\varepsilon}{2}}
x_{\gamma} \mathrel{\delta} x_{\alpha_{\delta}}$
and so as~$\delta\subseteq \nicefrac{\varepsilon}{2}$
we get~$x_\beta \mathrel\varepsilon x_{\alpha_\delta}$.
We have shown equivalence.

A different fact about fast nets will be useful later on:
if fast Cauchy nets~$(x_\alpha)_{\alpha \in \Phi}$
    and~$(y_\alpha)_{\alpha \in \Phi}$ are equivalent,
    then we can find for every~$\varepsilon \in \Phi$
    some~$\beta \in \Phi$
    such that for all~$\gamma \geq \beta$,
    we have~$x_\gamma \mathrel{\varepsilon} y_\gamma$.
\end{point}
\begin{point}{60}{$\overline{V}$: the uniform space~$N$}%
Write~$N$ for the set of fast Cauchy nets over~$V$.
Later we will define~$\overline{V}$ as~$N$ modulo equivalence.
Because of a subtlety with the definition of the uniformity
    on~$\overline{V}$ later,
    it is helpful to consider~$N$ separately.
Let~$\varepsilon \in \Phi$ be given.
For nets~$(x_\alpha)_{\alpha}$
    and~$(y_\alpha)_{\alpha}$ in~$N$,
    define
    \begin{equation*}
        (x_\alpha)_\alpha \mathrel{\Define{\hat\varepsilon}\index{*epsi@$\hat\varepsilon$}}
            (y_\alpha)_\alpha
            \quad\Leftrightarrow\quad
        \exists \beta \in \Phi \, \forall \gamma \geq \beta. \ 
        x_\gamma \mathrel{\varepsilon} y_\gamma.
    \end{equation*}
If~$\varepsilon \subseteq \delta$,
then $\hat\varepsilon \subseteq \hat\delta$
and
$\hat\varepsilon_1 \after \hat\varepsilon_2 
\subseteq\widehat{\varepsilon_1 \after \varepsilon_2}$.
So~$
\widehat{\nicefrac{\varepsilon}{2}} \after
\widehat{\nicefrac{\varepsilon}{2}} \subseteq
\widehat{\nicefrac{\varepsilon}{2} \after
\nicefrac{\varepsilon}{2}} \subseteq \widehat{\varepsilon}$,
which is one requirement for~$\{ \hat\varepsilon; \ \varepsilon \in \Phi\}$
to be a subbase for~$N$.
The others are easy as well. 
Also~$\widehat{\varepsilon_1 \cap \varepsilon_2} = \hat{\varepsilon}_1
    \cap \hat{\varepsilon}_2$
    and so each entourage of~$N$ has some~$\hat\varepsilon$ as subset.
For the following,
    it is helpful to remark now
    that Cauchy nets~$(x_\alpha)_\alpha$ and~$(y_\alpha)_\alpha$ over~$V$
    are equivalent if and only if~$(x_\alpha)_\alpha \mathrel{\hat\varepsilon}
    (y_\alpha)_\alpha$ for all~$\varepsilon \in \Phi$.
\end{point}
\begin{point}{70}{$\overline{V}$: $N$ is complete}%
As every Cauchy net is equivalent to a fast one,
it is sufficient to show convergence of fast Cauchy nets.
Thus let~$((x^\gamma_\alpha)_\alpha)_\gamma$
    be a fast Cauchy net in~$N$.
First we will show~$(x^{\hat\alpha}_\alpha)_\alpha$
    is a Cauchy net.
It might not be fast, so
    formally~$((x^\gamma_\alpha)_\alpha)_\gamma$
    cannot converge to it,
    but it will converge to an equivalent fast Cauchy net.
As~$((x^\gamma_\alpha)_\alpha)_\gamma$
    is fast,
    we know that for all~$\gamma_1,\gamma_2 \geq \hat\varepsilon$
we have~$
(x^{\gamma_1}_\alpha)_\alpha \mathrel{\hat\varepsilon^2}
        (x^{\gamma_2}_\alpha)_\alpha$.
So there is a~$\zeta_{\gamma_1,\gamma_2} \in \Phi$
such that for~$\alpha \geq \zeta_{\gamma_1,\gamma_2}$
    we have~$x^{\gamma_1}_\alpha
        \mathrel{\varepsilon^2}
        x^{\gamma_2}_\alpha$.
Suppose~$\alpha,\beta \geq \varepsilon$.
Pick~$\xi \geq \alpha,\beta,\zeta_{\hat\alpha,\hat\beta}$.
Then we have $ x^{\hat{\alpha}}_\alpha
            \mathrel{\varepsilon^2}
        x^{\hat{\alpha}}_\xi
            \mathrel{\varepsilon^2}
        x^{\hat{\beta}}_\xi
            \mathrel{\varepsilon^2}
            x^{\hat{\beta}}_\beta$ and so~$(x_\alpha^{\hat\alpha})_\alpha$
            is Cauchy.
Let~$(y_\alpha)_\alpha$ be a fast Cauchy net equivalent
to~$(x^{\hat\alpha}_\alpha)_\alpha$.
We want to prove~$(x^\gamma_\alpha)_\alpha \to (y_\alpha)_\alpha$.
As~$(y_\alpha)_\alpha$ is equivalent to~$(x^{\hat\alpha}_\alpha)_\alpha$
    we can find~$\alpha_0$
    such that~$y_\alpha \mathrel{\nicefrac{\varepsilon}{4}} x^{\hat\alpha}_\alpha$
    whenever~$\alpha \geq \alpha_0$.
As~$((x^\gamma_\alpha)_\alpha)_\gamma$
is fast Cauchy,
we know that~$(x^{\hat\alpha}_\beta)_\beta
\mathrel{\widehat{\nicefrac{\varepsilon}{4}}}
(x^\gamma_\beta)_\beta $
if~$\hat{\alpha}, \gamma \geq \widehat{\nicefrac{\varepsilon}{8}}$.
Thus there is some~$\beta_0$ such
that~$x^{\hat\alpha}_\beta
\mathrel{\nicefrac{\varepsilon}{4}}
x^\gamma_\beta $
for~$\beta \geq \beta_0$
and~$\hat{\alpha}, \gamma \geq \widehat{\nicefrac{\varepsilon}{8}}$.
Thus for~$\beta \geq \beta_0 \cap \nicefrac{\varepsilon}{8}$,
$\alpha \geq \alpha_0 \cap \nicefrac{\varepsilon}{8}$
and~$\gamma \geq \widehat{\nicefrac{\varepsilon}{8}}$
we get~$y_\alpha
\mathrel{\nicefrac{\varepsilon}{4}}
x^{\hat{\alpha}}_\alpha
\mathrel{\nicefrac{\varepsilon}{4}}
x^{\hat{\alpha}}_\beta
\mathrel{\nicefrac{\varepsilon}{4}}
x^\gamma_\beta
\mathrel{\nicefrac{\varepsilon}{4}}
x^\gamma_\alpha$.
Thus~$(y_\alpha)_\alpha \mathrel{\hat\varepsilon} (x^\gamma_\alpha)_\alpha$
whenever~$\gamma \geq \widehat{\nicefrac{\varepsilon}{8}}$,
so~$(x^\gamma_\alpha)_\alpha \to (y_\alpha)_\alpha$.
We have shown~$N$ is complete.
\end{point}
\begin{point}{80}[dils-uniformity-completion]{$\overline{V}$: uniformity}%
Let~$\overline{V}$ denote the set of fast Cauchy nets over~$V$
    modulo equivalence.
Write~$\eta\colon V \to \overline{V}$
for the map that sends~$x \in V$
to the equivalence class of the constant Cauchy net~$(x)_{\alpha \in \Phi}$.
For brevity, we may write~$\hat{x} = \eta(x)$.
As announced before, define~$V_0 \equiv \eta(V)$.

What uniformity to put on~$\overline{V}$?
Unfortunately the relations~$\hat\varepsilon$ do not necessarily preserve
equivalence ($\sim$) of Cauchy nets:
    there might be~$(x_\alpha)_\alpha \sim (x'_\alpha)_\alpha$
    and~$(y_\alpha)_\alpha \sim (y'_\alpha)_\alpha$ in~$N$
such that~$(x_\alpha)_\alpha \mathrel{\hat\varepsilon} (y_\alpha)_\alpha$,
but not~$(x'_\alpha)_\alpha \mathrel{\hat\varepsilon} (y'_\alpha)_\alpha$
for some~$\varepsilon \in \Phi$.
Instead we define for any~$\varepsilon \in \Phi$
the following entourages on~$\overline{V}$:
\begin{equation*}
    (x_\alpha)_\alpha \mathrel{\Define{\tilde\varepsilon}\index{*epsiti@$\tilde\varepsilon$}}
    (y_\alpha)_\alpha \quad\Leftrightarrow\quad
    (x'_\alpha)_\alpha \mathrel{\hat\varepsilon}
    (y'_\alpha)_\alpha \quad
    \text{for all
        $(x'_\alpha)_\alpha \sim (x_\alpha)_\alpha$ 
    and $(y'_\alpha)_\alpha \sim (y_\alpha)_\alpha$.}
\end{equation*}
By definition~$\tilde\varepsilon$
respects equivalence and can be considered as a relation on~$\overline{V}$.
It is not hard to verify~$\{ \tilde\varepsilon ;\ \varepsilon \in \Phi\}$
is a subbase for~$\overline{V}$
and~$\widetilde{\varepsilon_1 \cap \varepsilon_2} = \tilde{\varepsilon}_1
\cap \tilde{\varepsilon}_2$.
It is easy to see~$\overline{V}$ is a Hausdorff uniform space.
Write~$\tilde\Phi$ for the generated set of entourages on~$\overline{V}$.
The entourages of~$\overline{V}$ are exactly those relations
    which have some~$\tilde\varepsilon$
    as a subset.

Furthermore, if~$(x_\alpha)_\alpha \mathrel{\hat\varepsilon}
(y_\alpha)_\alpha $,
then~$(x_\alpha)_\alpha \mathrel{\tilde\varepsilon^3} (y_\alpha)_\alpha$.
Thus if a net~$(x^\gamma_\alpha)_\alpha$
converges to~$(y_\alpha)_\alpha$ in~$N$,
then so do their equivalence classes in~$\overline{V}$.
Hence~$\overline{V}$ is complete.
(The inclusion~$\hat\varepsilon \subseteq \tilde\varepsilon^3 \subseteq
        \hat\varepsilon^3$
        implies that the~$\hat\varepsilon$ and~$\tilde\varepsilon$
        generate the same uniformity on~$N$.)
\end{point}
\begin{point}{90}{$\overline{V}$ is a right $\scrB$-module}%
    First we will define an addition on~$\overline{V}$.
    The addition on~$V$ is uniformly continuous
    (as the uniformity on~$V$ is given by seminorms,
    see \sref{dils-product-uniformity} for the uniformity on~$V^2$.)
    Thus for any~$(x_\alpha)_\alpha$
    and~$(y_\alpha)_\alpha$ in~$N$,
    the net~$(x_\alpha+y_\alpha)_\alpha$ is again Cauchy,
        see \sref{dils-uniform-spaces-basics}.
    Also because of the uniform continuity of addition:
    if~$(x_\alpha)_\alpha \sim (x'_\alpha)_\alpha$
    and~$(y_\alpha)_\alpha \sim (y'_\alpha)_\alpha$,
    then~$(x_\alpha+y_\alpha)_\alpha \sim
            (x'_\alpha+y'_\alpha)_\alpha$.
    There is a net in~$N$ equivalent to $(x_\alpha+y_\alpha)_\alpha$.
    This fixes an addition on~$\overline{V}$,
        which turns it into an Abelian group.
        By construction~$\eta(x+y) = \eta(x)+\eta(y)$.

Pick any~$b \in \scrB$. We show the map~$r_b\colon V \to V$
given by~$r_b(x) = xb$ is also uniformly continuous,
which requires us to unfold the definitions further
than with addition.
Assume we are given an entourage in~$V$,
that is: np-maps~$f_1,\ldots,f_n\colon \scrB \to \C$
and~$\varepsilon > 0$.
For any np-map~$f\colon \scrB \to \C$,
    the map~$b*f$ given by~$(b * f)(x) \equiv f(b^* x b)$ is also an np-map,
    see \sref{bstaromega-basic}.
    Clearly~$\|xb\|_f = f([xb,xb])^{\frac{1}{2}}
    = f(b^* [x,x] b)^{\frac{1}{2}} = \|x\|_{b*f}$
    for any~$x \in V$.
Thus if~$\|x-y\|_{b*f_i} \leq \varepsilon$ for~$1 \leq i \leq n$,
then~$\|xb-yb\|_{f_i} = \|x-y\|_{b*f_i} \leq \varepsilon$
as well. Hence~$r_b$ is uniformly continuous.
As before, this uniform continuity
allows us to define a right~$\scrB$-action on~$\overline{V}$
by sending the equivalence class of~$(x_\alpha)_\alpha$
to~$(x_\alpha b)_\alpha$. By definition~$\eta(x)b = \eta(xb)$.
It is straightforward to check this turns~$\overline{V}$
into a right $\scrB$-module.
\end{point}
\begin{point}{100}[dils-extending-norm-to-overlineV]{Extending $\|\,\cdot\,\|_f$ to~$\overline{V}$}%
Let~$(x_\alpha)_\alpha$ be a Cauchy net over~$V$
    and~$f\colon \scrB \to \C$ be any np-map.
From the reverse triangle inequality~$| \|x_\alpha \|_f - \|x_\beta\|_f |
        \leq \|x_\alpha - x_\beta \|_f$,
        it follows~$(\| x_\alpha \|_f)_\alpha$ is Cauchy.
Define~$\|(x_\alpha)_\alpha\|_f = \lim_\alpha \|x_\alpha\|_f$.
Again using the reverse triangle
    inequality~$| \|x_\alpha \|_f - \|x'_\alpha\|_f |
        \leq \|x_\alpha - x'_\alpha \|_f$
        we see $\|(x_\alpha)_\alpha\|_f =
        \|(x'_\alpha)_\alpha\|_f$
        whenever $(x_\alpha)_\alpha \sim (x'_\alpha)_\alpha$.
Thus~$\|\,\cdot\,\|_f$ lifts to~$\overline{V}$.

These extended norms~$\|\,\cdot\,\|_f$ also induce a uniformity
    on~$\overline{V}$.  We will now show it is the same uniformity.
    Let~$((x_\alpha^\gamma)_\alpha)_\gamma$ be a net in~$\overline{V}$.
    It is sufficient to show that~$(x_\alpha^\gamma)_\alpha \to 0$
    in the original uniformity
    if and only if~$\|(x_\alpha^\gamma)_\alpha\|_f \to 0$
        for all np-maps~$f\colon \scrB \to \C$.

Assume~$(x^\gamma_\alpha)_\alpha \to 0$ in~$\overline{V}$.
That is: for each~$f\colon \scrB \to \C$ and~$\varepsilon > 0$
    there is a~$\gamma_0$
    such that for each~$\gamma \geq \gamma_0$
    there is an~$\alpha_0$
    such that~$\|x^\gamma_\alpha\|_f \leq \varepsilon$
    when~$\alpha \geq \alpha_0$.
Let any~$f\colon \scrB \to \C$ and~$\varepsilon > 0$ be given.
Find such~$\gamma_0$ for~$f$ and~$\varepsilon$.
Then~$\|(x^\gamma_\alpha)_\alpha \|_f
    = \lim_\alpha \|x^\gamma_\alpha\| \leq \varepsilon$
    whenever~$\gamma \geq \gamma_0$
    and so~$\|(x^\gamma_\alpha)_\alpha\|_f \to 0$ as desired.

For the converse, assume~$\|(x^\gamma_\alpha)_\alpha\|_f \to 0$
    for all np-maps~$f\colon \scrB \to \C$.
Let~$\varepsilon > 0$ and np-map~$f\colon \scrB \to \C$ be given.
There is some~$\gamma_0$ such that
$\lim_\alpha \|x^\gamma_\alpha \|_f = \|(x^\gamma_\alpha)_\alpha\|_f
\leq \frac{1}{2}\varepsilon$
for~$\gamma\geq\gamma_0$.
Hence~$\|x^\gamma_\alpha\|_f \leq \varepsilon$
    for sufficiently large~$\alpha$,
    which is to say~$(x^\gamma_\alpha)_\alpha \to 0$ in~$\overline{V}$.
\end{point}
\begin{point}{110}[dils-completion-setup]{Induction set-up}%
The majority of the remaining work is to show
    we can define an inner product on~$\overline{V}$.
As sketched earlier,
    we will extend the inner product
    from~$V_0$ step-by-step.
Call a subset~$W \subseteq \overline{V}$
together with a~$\scrB$-valued inner product~$[\,\cdot\,,\,\cdot\,]$ a
    \index{$[\,\cdot\,,\,\cdot\,]$!inner product!$\scrB$-valued}
\Define{compatible extension} if
\begin{enumerate}
\item $V_0 \subseteq W$;
\item $W$ is a $\scrB$-submodule of~$\overline{V}$;
\item the inner product turns~$W$ into a pre-Hilbert $\scrB$-module \emph{and}
\item $\|x\|_f^{\overline{V}} = \|x\|_f^W$
        for all~$x \in W$
        and np-maps~$f\colon \scrB \to \C$,
        where~$\|x\|_f^W = f([x,x]_W)^{\frac{1}{2}}$
        and~$\|\,\cdot\,\|^{\overline{V}}_f$
        is the extension of~$\|\,\cdot\,\|_f$
        on~$V$ to~$\overline{V}$
        as in \sref{dils-extending-norm-to-overlineV}.
\end{enumerate}
As~$\|\,\cdot\,\|^{\overline{V}}_f$ also
    generate the uniformity on~$\overline{V}$,
    the last requirement implies
    that on a compatible extension~$W$
    the ultranorm uniformity and
    the uniformity induced by~$\overline{V}$
    coincide.

For compatible extensions~$W_1,W_2 \subseteq \mathscr{V}$
        we say~$W_1 \leq W_2$
            iff~$W_1 \subseteq W_2$
            and~$[v,w]_{W_1} = [v,w]_{W_2}$
        for all~$v,w\in W_1$.
Later, we will apply Zorn's lemma to the compatible extensions
    ordered in this way.
\end{point}
\begin{point}{120}{Induction base case}%
    The inner product of~$V$ lifts to~$V_0$.
    Indeed, $\eta(x) = 0$ if and only if~$\| x\| = 0$.
    So by Cauchy--Schwarz~$[x,y] = 0$ whenever~$\eta(x)=0$.
    Thus~$[x,y] = [x',y']$
        if~$\eta(x)=\eta(x')$
        and~$\eta(y)=\eta(y')$.
    Define~$[\eta(x),\eta(y)] = [x,y]$.
    It's easy to see~$V_0$ is a pre-Hilbert~$\scrB$-submodule of~$\overline{V}$.
    Finally, for any np-map~$f\colon \scrB\to \C$
    we have $\| \eta(x) \|_f^{\overline{V}}
    =   \| x\|_f^{V} = \| \eta(x) \|_f^{V_0}$, as desired.
\end{point}
\begin{point}{130}{Induction step}%
Assume~$W\subseteq \overline{V}$ with~$[\,\cdot\,,\,\,\cdot\,]$
    is a compatible extension.
Let~$\Define{\sigma(W)}$ denote the limits of norm-bounded Cauchy nets over~$W$.
We will define a~$\scrB$-valued inner product on~$\sigma(W)$
    that turns it into a compatible extension.
In fact~$W \leq \sigma(W)$.
By assumption the induced uniformity on~$W$
    is the ultranorm uniformity
    and so addition and the~$\scrB$-action
    are uniformly continuous.
Hence~$\sigma(W)$ is again a submodule.

To define the inner product,
assume~$x,y \in \sigma(W)$.
There are~$\Phi$-indexed norm-bounded Cauchy
    nets~$(x_\alpha)_\alpha$ and~$(y_\alpha)_\alpha$
    over~$W$
    such that~$x_\alpha \to x$ and~$y_\alpha \to y$.
For any np-map $f\colon \scrB \to \C$ we have
    by Cauchy--Schwarz
\begin{align*}\label{eq-dils-V-inprod-cauchy}
    \bigl|f([x_\alpha, y_\alpha] - [x_\beta,y_\beta])\bigr|
    & \ =\  \bigl| [x_\alpha, y_\alpha - y_\beta]_f
        + [x_\alpha - x_\beta, y_\beta]_f\bigr| \\
    &\  \leq\  \|x_\alpha\|_f \|y_\alpha - y_\beta\|_f \numberthis
        + \|x_\alpha - x_\beta\|_f \| y_\beta \|_f.
\end{align*}
(If not otherwise specified $\|\,\cdot\,\|_f = \|\,\cdot\,\|_f^W$.)
As~$(x_\alpha)_\alpha$ is norm bounded and
$\|x_\alpha\|_f^2 = f([x_\alpha,x_\alpha]) \leq \|f\| \|x_\alpha\|^2$,
we see~$(\|x_\alpha\|_f)_\alpha$ is bounded.
Similarly~$(\|y_\beta\|_f)_\beta$ is bounded.
Thus from \eqref{eq-dils-V-inprod-cauchy}
it follows~$([x_\alpha,y_\alpha])_\alpha$
    is a Cauchy net in the ultraweak uniformity of~$\scrB$,
    see \sref{dils-uniformity-examples}.
As~$\|[x_\alpha,y_\alpha]\| \leq \|x_\alpha\| \|y_\alpha\|$,
see \sref{module-CS}, it is also a norm-bounded net.
By \sref{vn-complete} norm-bounded Cauchy nets in the ultraweak uniformity
converge, so~$\uwlim_\alpha [x_\alpha,y_\alpha]$ exists.

Assume we are given any~$\Phi$-indexed
norm-bounded Cauchy nets~$(x')_\alpha$ and~$(y')_\alpha$ over~$W$
with~$ (x'_\alpha)_\alpha \to x$
and~$ (y'_\alpha)_\alpha \to y$.
Then, like \eqref{eq-dils-V-inprod-cauchy}:
\begin{align*}
    \bigl|f([x_\alpha, y_\alpha] - [x'_\alpha,y'_\alpha])\bigr|
    & \ =\  \bigl| [x_\alpha, y_\alpha - y'_\alpha]_f
        + [x_\alpha - x'_\alpha, y'_\alpha]_f\bigr| \\
    &\  \leq\  \|x_\alpha\|_f \|y_\alpha - y'_\alpha\|_f
        + \|x_\alpha - x'_\alpha\|_f \| y'_\alpha \|_f.
\end{align*}
From this it follows~$ \uwlim_\alpha [x_\alpha,y_\alpha]=
    \uwlim_\alpha [x'_\alpha,y'_\alpha]$.
We are justified to define
\begin{equation*}
    [x,y] \ \equiv\  \uwlim_{\alpha} \,[x_\alpha,y_\alpha].
\end{equation*}
We will show~$[\,\cdot\,,\,\cdot\,]$ is a~$\scrB$-valued inner product.
    \index{$[\,\cdot\,,\,\cdot\,]$!inner product!$\scrB$-valued}
Pick any~$x,y,z \in \sigma(W)$
    and~$\Phi$-indexed
norm-bounded Cauchy nets
$(x_\alpha)_\alpha$, $(y_\alpha)_\alpha$ and $(z_\alpha)_\alpha$ over~$W$
with~$x_\alpha \to x$, $y_\alpha \to y$ and~$z_\alpha \to z$ in~$\overline{V}$.
So~$x_\alpha+ y_\alpha \to x + y.$
With ultraweak continuity of addition
\begin{equation*}
    [x+y,z] = \uwlim_\alpha [x_\alpha +y_\alpha, z_\alpha]
        = \uwlim_\alpha [x_\alpha, z_\alpha] +
                \uwlim_\alpha [y_\alpha, z_\alpha]
                = [x,z] + [y,z].
\end{equation*}
In a similar fashion one proves the other axioms of
    a~$\scrB$-valued inner product:
    $[x,y]=[y,x]^*$ follows from the ultraweak continuity of~$(\,\cdot\,)^*$;
    $\scrB$-homogeneity follows from
        ultraweak continuity of~$c \mapsto cb$ (see \sref{mult-uws-cont})
        \emph{and}~$[x,x]\geq 0$ follows from ultraweak-closedness of the
        positive cone of~$\scrB$.

        To prove definiteness of~$[\,\cdot\,,\,\cdot\,]$, assume~$[x,x]=0$
            for some~$x \in \sigma(W)$.
    Pick some $\Phi$-indexed norm-bounded
        Cauchy net~$(x_\alpha)_\alpha$ over~$W$
        with~$x_\alpha \to x$.
For any np-map~$f\colon \scrB \to \C$
we have~$\lim_\alpha f([x_\alpha,x_\alpha])
         = f(\uwlim_\alpha [x_\alpha,x_\alpha])
         = f([x,x])  =  0$.
Thus~$x_\alpha$ converges in the ultranorm uniformity on~$W$ to~$0$.
Hence~$x=0$.

Thus~$\sigma(W)$ is a pre-Hilbert~$\scrB$-module.
It remains to be shown that~$\|\,\cdot\,\|_f^{\sigma(W)}$
    and $\|\,\cdot\,\|_f^{\overline{V}}$
    agree.
Pick any~$x \in \sigma(W)$ and~$f\colon \scrB \to \C$.
By definition there is a norm-bounded Cauchy net~$x_\alpha$ over~$W$
    for which~$x_\alpha \to x$ in~$\overline{V}$.
Then
\begin{equation*}
    \|x\|^{\sigma(W)}_f
       \ =\ f(\uwlim_\alpha [x_\alpha,x_\alpha])^{\frac{1}{2}}
       \ =\ \lim_\alpha \|x_\alpha\|_f^W
       \ =\ \lim_\alpha \|x_\alpha\|_f^{\overline{V}}
       \ =\ \|x\|_f^{\overline{V}},
\end{equation*}
as desired.
We have shown that~$\sigma(W)$ is a compatible extension.
\end{point}
\begin{point}{140}{Induction limit step}%
    Assume~$\mathscr{V}$
        is a non-empty chain of compatible extensions.
        See \sref{dils-completion-setup} for the order on compatible extensions.
    Write~$W = \bigcup \mathscr{V}$.
    We will turn~$W$ into a compatible extension.
    Clearly~$W$ is a submodule of~$\overline{V}$.
    As~$\mathscr{V}$ is not empty, $V_0 \subseteq W$.
    For any~$v,w \in W$
    there is some~$W' \in \mathscr{V}$ with~$v,w \in W'$.
    Define~$[v,w]_W = [v,w]_{W'}$.
    This turns~$W$ into a pre-Hilbert~$\mathscr{B}$-module
        as all its axioms only involve finitely many elements,
        which are contained in some single compatible extension far enough up
        in the chain.
    Finally, to see the seminorms agree,
        pick any np-map~$f\colon \scrB \to \C$
        and~$x \in W$.
    Then~$x \in W'$ for some~$W' \in \mathscr{V}$.
    Consequently
    \begin{equation*}
        \|x\|^{\overline{V}}_f
        \ =\ \|x\|^{W'}_f
        \ =\ f([x,x]_{W'})^{\frac{1}{2}}
        \ =\ f([x,x]_{W})^{\frac{1}{2}}
        \ =\  \|x\|^{W}_f,
    \end{equation*}
    as desired.  We have shown~$W$ is a compatible extension.
\end{point}
\begin{point}{150}{Self-duality}
    We have seen
    that every non-empty chain of compatible extensions
        has an upper bound.
    Thus by Zorn's lemma there is a maximal compatible
        extension~$W\subseteq \overline{V}$.
    By maximality~$\sigma(W) = W$
        and so in~$W$ norm-bounded ultranorm-Cauchy nets converge.
    Thus by~\sref{dils-selfdual}
        $W$ must be self-dual and ultranorm complete.
    By completeness~$W = \overline{V}$.
    Define~$X = \overline{V}$
        and~$\left<x,y\right> = [x,y]_{W}$.
    Clearly~$[v,w] = \left<\eta(v),\eta(w)\right>$
        and the image of~$\eta$ is norm dense. \qed
\end{point}
\end{point}
\end{point}
\end{parsec}
\begin{parsec}{1510}%
\begin{point}{10}%
The completion~\sref{dils-completion} we just constructed
    has a universal property,
    which can be seen as a generalization of~\cite[prop.~3.6]{paschke}.
\end{point}
\begin{point}{11}[selfdual-completion-univ]{Lemma}%
Assume~$\scrB$ is a von Neumann algebra
    and~$V$ is a right $\scrB$-module
    with~$\scrB$-valued inner product.
Let~$X$ be any self-dual Hilbert~$\scrB$-module
    together with inner product preserving
    $\scrB$-linear map~$\eta\colon V \to X$
    such that the image of~$\eta$ is ultranorm dense in~$X$.
(E.g.~the completion into a self-dual Hilbert~$\scrB$-module
from \sref{dils-completion}.)
Then~$\eta$ has the following universal property.
\begin{quote}
Let~$Y$ be a self-dual Hilbert~$\scrB$-module.
For each bounded~$\scrB$-linear
    map~$T\colon V \to Y$,
    there is a unique bounded~$\scrB$-linear
    map~$\hat{T} \colon X \to Y$
    such that~$\hat{T} \after \eta = T$.
\end{quote}
    \spacingfix{}
\begin{point}{20}{Proof}%
By \sref{blinear-bounded-is-ultranorm} bounded~$\scrB$-linear
    maps are ultranorm uniformly continuous.
As~$\eta(V)$ is ultranorm dense in~$X$,
    the extension~$\hat{T}$ is unique if it exists.
Pick any~$x \in X$.
There is a net~$x_\alpha$ in~$V$ such
    that~$\eta(x_\alpha) \to x$ ultranorm.
So~$\eta(x_\alpha)$ is ultranorm Cauchy.
As~$\eta$ preserves inner product,
    $x_\alpha$ is ultranorm Cauchy as well.
As~$T$ is ultranorm uniformly continuous
    $T(x_\alpha)$ is ultranorm Cauchy.
With the same motions we see~$T(x_\alpha) \sim T(x'_\alpha)$
    when~$x_\alpha \sim x'_\alpha$.
Using this and that~$Y$ is ultranorm complete, see \sref{dils-selfdual},
    we may define~$\hat{T} x = \unlim_\alpha x_\alpha$,
    where~\Define{$\unlim$}\index{unlim@$\unlim$} denotes the limit in the ultranorm uniformity.
Clearly~$\hat{T}\after\eta=T$.
It is straightforward to check
    that $\scrB$-linearity
    of~$\hat{T}$
    follows from~$\scrB$-linearity of~$T$
    and ultranorm continuity of addition and
        the~$\scrB$-module action.

It remains to be shown~$\hat{T}$ is bounded.
By \sref{blinear-inprod-inequality}  we have
\begin{equation*}
\left<T x_\alpha,T x_\alpha\right>
\ \leq\  \| T\|^2 [x_\alpha,x_\alpha]
\ =\  \|T\|^2 \left<\eta(x_\alpha),\eta(x_\alpha)\right>
\end{equation*}
and so using \sref{innerprod-ultraweak} we see that
\begin{equation*}
\langle\hat{T} x, \hat{T} x\rangle \ =\ 
\uwlim_\alpha \left<T x_\alpha, T x_\alpha\right>
            \ \leq \ \|T\|^2 \uwlim_\alpha \left<\eta(x_\alpha), \eta(x_\alpha)\right>
            \ =\  \|T\|^2 \left<x,x\right>,
\end{equation*}
which implies that~$\|\hat{T}\| \leq \|T\|$ as desired.
(In fact~$\|\hat{T}\| = \|T\|$.) \qed
\end{point}
\end{point}
\end{parsec}
\begin{parsec}{1520}%
\begin{point}{10}%
    Before we show that every ncp-map has a Paschke dilation,
        we study~$\scrB^a(X)$ a bit more for self-dual~$X$.
    We start with~$\scrB$-sesquilinear forms.
\end{point}
\begin{point}{20}{Definition}%
Let~$V$ be a normed right $\scrB$-module.
A sesquilinear form~$B\colon V\times V \to \scrB$
is said to be \Define{bounded}\index{sesquilinear form!$\scrB$-valued!bounded}
    if there is an~$r \geq 0$ such that
$\|B(x,y)\| \leq r \|x\|\|y\|$
for all~$x,y \in V$.
\end{point}
\begin{point}{30}{Example}%
Let~$X$ be a pre-Hilbert~$\scrB$-module.
For every~$T \in \scrB^a(X)$
    the map~$\left<(\,\cdot\,), T (\,\cdot\,)\right>$
    is a bounded~$\scrB$-sesquilinear form.
\begin{point}{40}%
    For self-dual~$X$ every bounded~$\scrB$-sesquilinear form
        arises in this way.
\end{point}
\end{point}
\begin{point}{50}[hilbmod-sesquilinear-forms]{Proposition}%
Let~$X$ be a self-dual Hilbert~$\scrB$-module.
For every bounded $\scrB$-sesquilinear form~$B$ on~$X$,
    there is a unique~$T \in \scrB^a(X)$
    with~$B(x,y) = \left<x,Ty\right>$
    for all~$x,y \in X$.
\begin{point}{60}{Proof}
For each~$y \in X$
    the map~$B(\,\cdot\,,y)^*$ is~$\scrB$-linear and bounded.
Thus by self-duality of~$X$
there is a unique~$t_y$
    for each~$y \in X$
    such that~$\left<t_y, x\right> = B(x,y)^*$
    for all~$x \in X$.
It is easy to see~$\left<t_{y+y'},x\right>
                = \left<t_y + t_{y'}, x\right>$
                and~$\left<t_{yb},x\right>
                = \left<t_y b, x\right>$,
    and so by uniqueness~$y \mapsto t_y$ is~$\scrB$-linear.
Define~$Ty = t_y$.
Clearly~$\left<x,Ty\right> = \left<t_y,x\right>^* = B(x,y)$
    and~$T$ must be the unique such~$\scrB$-linear map.
Reasoning in the same way,
    we see there is a unique~$S$
    with~$B(y,x)^* = \langle x, S y \rangle$.
Hence~$\langle x, T y\rangle = B(x,y) 
        = \langle y, Sx\rangle^* = \langle Sx, y \rangle$,
        so~$T$ is adjointable.
It remains to be shown~$T$ is bounded.
There is an~$r \geq 0$ such that for all~$x \in X$:
\begin{equation*}
    \| Tx\|^2 
    \ =\  \|\left<Tx,Tx\right>\|
            \ =\  \| B (Tx, x) \|
            \ \leq\  r \|Tx\|\|x\|.
\end{equation*}
So~$r$ is a bound
    by dividing out~$\| T x \|$
    if~$\| Tx \| \neq 0$
    and trivially otherwise.  \qed
\end{point}
\begin{point}{70}{Remark}%
In \sref{chilb-form-representation} a more general result was shown.
\end{point}
\end{point}
\begin{point}{80}[hilbmod-adjoint-exists]{Exercise}%
Suppose~$T\colon X \to Y$ is a bounded~$\scrB$-linear map between
    Hilbert~$\scrB$-modules.
Show that if~$X$ is self dual,
    then~$T$ is adjointable. \cite[prop.~3.4]{paschke}
\end{point}
\begin{point}{90}[hilmod-fixed-on-V]{Exercise}%
Let~$V$ be a right $\scrB$-module with~$\scrB$-valued inner product
        for some von Neumann algebra~$\scrB$ and
    $\eta\colon V \to X$ be the ultranorm completion of~$V$
    from~\sref{dils-completion}.
Show that
\begin{equation*}
    \{\ \langle \hat{x}, (\,\cdot\,)\, \hat{x}\rangle \colon X \to \scrB \ ;\ x \in V\  \} \qquad \text{(where $\hat{x} \equiv \eta(x)$)}
\end{equation*}
is an order separating set of ncp-maps, see \sref{separating}.
    (That is:~$T \geq 0$ iff~$\langle \hat{x}, T \hat{x}\rangle \geq 0$
    for all~$x \in V$.
    Consequently~$S = T$ iff~$\langle \hat{x}, T \hat{x} \rangle 
    = \langle \hat{x}, S \hat{x} \rangle$ for all~$x \in V$.)
\end{point}
\begin{point}{100}{Theorem}%
Suppose~$X$ is a self-dual Hilbert~$\scrB$-module
    for a von Neumann algebra~$\scrB$.
Then~$\scrB^a(X)$ is a von Neumann algebra.
\begin{point}{110}{Proof}%
The theorem is due to Paschke \cite[prop.~3.10]{paschke};
    we give a new proof, which also appeared as \sref{bah-vn}.
We already know that~$\scrB^a(X)$ is a C$^*$-algebra, see \sref{hilbmod-cstar}.
\begin{point}{120}{bounded order completeness}%
Let~$(T_\alpha)_\alpha$ be a norm-bounded net of self-adjoint elements
    of~$\scrB^a(X)$.  We have to show it has a supremum.
Pick~$x\in X$ and  $r\in \R$, $r\geq0 $
such that~$\|T_\alpha\| \leq r$ for all~$\alpha$.
By \sref{hilbmod-ordersep}~$(\left<x,T_\alpha x\right>)_\alpha$
    is a norm-bounded net of self-adjoint elements of~$\scrB$
    and so it has a supremum to which
    it converges ultrastrongly by~\sref{vna-supremum-uslimit}.
In particular
\begin{equation*}
\frac{1}{4} \sum^3_{k=0} i^k \left< i^kx+y, T_\alpha(i^k x+y)\right>
    \ =\  \left<x, T_\alpha y\right>
\end{equation*}
converges ultrastrongly for all~$x,y \in X$.
Define~$B(x,y) = \uslim_\alpha \left<x,T_\alpha y\right>$.
As addition and multiplication by a fixed element are
    ultrastrongly continuous (see \sref{mult-uws-cont}),
    we see that~$B$ is a~$\scrB$-sesquilinear form.
By~\sref{module-CS} and~\sref{blinear-inprod-inequality}
    we have~$\left<T_\alpha y, x\right>\left<x, T_\alpha y\right>
        \leq r^2 \|x\|^2 \left<y,y\right>$
    for each~$\alpha$, hence
\begin{align*}
    \| B(x,y) \|^2 &\ \ = \ \ 
    \| B(x,y)^* B(x,y) \| \\
    &\ \ \overset{\mathclap{\sref{usconv}}}{=} \ \ \| \uwlim_\alpha \left<T_\alpha y, x\right>\left< x, T_\alpha y\right>\| \\
    &\ \ \leq \ \  r^2 \|x\|^2 \|y\|^2.
\end{align*}
So~$B$ is a bounded~$\scrB$-sesquilinear form.
By \sref{hilbmod-sesquilinear-forms}
    there is a unique~$T \in \scrB^a(X)$
    with~$\left<x,Ty\right> = B(x,y) = \uslim_\alpha \left<x,T_\alpha y\right>$.
Noting~$B(x,y) = \uwlim_\alpha \langle x, T_\alpha y \rangle$,
    it is easy to see~$T$ is self-adjoint.
Using \sref{hilbmod-ordersep} and 
the fact that ultraweak limits respect the order (\sref{vn-positive-basic}),
    we see~$\left<x,T_\alpha x\right> \leq \left<x,T x\right>$
    for all~$x \in X$ and so~$T_\alpha \leq T$.
Now suppose we are given a self-adjoint~$S \in \scrB^a(X)$
    with~$T_\alpha \leq S$ for all~$\alpha$.
With a similar argument we see~$T \leq S$ and so~$S$ is
    the supremum of~$(T_\alpha)_\alpha$.
\end{point}
\begin{point}{130}[hilbmod-vecstates-normal]{separating normal states}%
By~\sref{hilbmod-ordersep},
    the states~$\left<x,(\,\cdot\,)x\right>$
    are separating.
We are done if we can show that~$\left<x,(\,\cdot\,)x\right>$
    is normal, i.e.~preserves suprema of bounded directed sets of
    self-adjoint elements.
To this end, pick~$x \in X$ and let~$(T_\alpha)_\alpha$
    be a net with supremum~$T$.
As before~$(\left<x,T_\alpha x\right>)_\alpha)$
    is a norm bounded net which converges ultrastrongly to its
    supremum.
Now, we just saw
\begin{equation*}
    \left<x,Tx\right> \ =\  \uwlim_\alpha \left<x, T_\alpha x\right>
               \  =\  \sup_\alpha \left<x, T_\alpha x\right>
\end{equation*}
and so indeed~$\left<x, (\,\cdot\,)x\right>$ is normal;
the normal states are separating
and consequently~$\scrB^a(X)$ is a von Neumann algebra. \qed
\end{point}
\end{point}
\end{point}
\end{parsec}

\begin{parsec}{1530}%
\begin{point}{10}[hilbmod-ad-ncp]{Proposition}%
Assume~$T \colon X \to Y$ is an adjointable bounded module map
    between Hilbert~$\scrB$-modules.
    Define $\Define{\ad_T}\colon \scrB^a(Y) \to \scrB^a(X)$\index{$\ad_V$!Hilbert $\scrB$-modules}
    by~$\ad_T (B) = T^*BT$.
The map~$\ad_T$ is completely positive.
If~$X$ and~$Y$ are self-dual,
    then~$\ad_T$ is normal.
\begin{point}{20}{Proof}%
For any~$n \in N$, $B_1,\ldots, B_n \in \scrB^a(X)$
    and~$A_1,\ldots,A_n \in \scrB^a(Y)$
    we have
\begin{align*}
    \sum_{i,j} B_j^* T^* A_j^*A_i T B_i \ = \ 
    \bigl( \sum_i A_i T B_i \bigr)^*
    \bigl( \sum_j A_j T B_j \bigr) \ \geq\  0
\end{align*}
and so~$\ad_T$ is completely positive.
\end{point}
\begin{point}{30}%
Now we show~$\ad_S$ is normal.
By \sref{normal-faithful} and \sref{hilbmod-ordersep},
    it suffices to show that
    for every~$x\in X$, the
    map~$T \mapsto \langle x, \ad_S (T) x\rangle
                = \langle S x, T \, S x\rangle$
    is normal,
    which it indeed is by \sref{hilbmod-vecstates-normal}. \qed
\end{point}
\end{point}
\begin{point}{40}[hilbmod-adj-vector-ncp]{Exercise}%
Let~$\scrA$ be a C$^*$-algebra.
Assume~$n \in \N$ and  $a_1, \ldots, a_n \in \scrA$.
Use \sref{hilbmod-ad-ncp}
to show~$\varphi \colon \scrA \to M_n \scrA$
    given by~$\varphi(d) = (a_i^*da_j)_{ij}$
    is an ncp-map.
\end{point}
\end{parsec}

\section{Paschke dilations}
\begin{parsec}{1540}%
\begin{point}{10}%
We are ready to show that every ncp-map has a Paschke dilation.
    That is: for every ncp-map~$\varphi$
        there is a triple~$(\scrP, \varrho, h)$
        with~$\varphi = h \after \varrho$ via~$\scrP$
            for some ncp-map~$h$, nmiu-map~$\varrho$ and von Neumann algebra~$\scrP$,
            which is minimal in the sense
            that for any such triple~$(\scrP', \varrho', h')$
            there is a unique ncp-map~$\sigma\colon \scrP' \to \scrP$
            with~$h' =h \after \sigma$ and~$\varrho=\sigma\after\varrho'$.
    See \sref{def-paschke}.
\end{point}
\begin{point}{20}[phi-compatible-paschke]{Definition}%
    Let~$\varphi\colon \scrA \to \scrB$ be any ncp-map between
        von Neumann algebras.
    A complex bilinear map~$B \colon \scrA \times \scrB \to X$,
        where~$X$ is a self-dual Hilbert~$\scrB$-module is
        called \Define{$\varphi$-compatible}\index{$\varphi$-compatible}
        if there is an~$r > 0$
        such that for all~$n\in \N$, $a_1, \ldots, a_n \in \scrA$
        and~$b_1, \ldots, b_n \in \scrB$ we have
        \begin{equation}
            \bigl\| \sum_i B(a_i,b_i)\bigr\|^2
             \    \leq \ r \cdot \bigl\| \sum_{i,j} b_i^* \varphi(a_i^*a_j)b_j
                \bigr\| \label{phi-compatible}
        \end{equation}
        and~$B(a,b_1)b_2 = B(a,b_1b_2)$
        for all~$a \in \scrA$ and~$b_1,b_2 \in \scrB$.
\end{point}
\begin{point}{30}[existence-paschke]{Theorem}%
    Let~$\varphi\colon \scrA \to \scrB$ be any ncp-map between
        von Neumann algebras.
\begin{enumerate}
    \item There is a self-dual Hilbert~$\scrB$-module~$\scrA \otimes_\varphi
            \scrB$ and~$\varphi$-compatible bilinear
            \index{*tensorphi@$\otimes_\varphi$}
    \begin{equation*}
        \otimes \colon \scrA \times \scrB \to \scrA \otimes_\varphi \scrB
    \end{equation*}
    such that for every~$\varphi$-compatible bilinear
    map~$B \colon \scrA \times \scrB \to Y$
    there is a unique bounded
    module map~$T\colon \scrA \otimes_\varphi \scrB \to Y$
    such that~$T(a \otimes b) = B(a,b)$ for all~$a \in \scrA$ and~$b \in \scrB$.
\item
    For~$a_0 \in \scrA$,
    there is a unique~$\varrho(a_0) \in \scrB^a(\scrA \otimes_\varphi \scrB)$
            fixed by
            \begin{equation*}
                \varrho(a_0)(a\otimes b) = (a_0 a)\otimes b.
                    \text{\qquad ($a\in \scrA$, $b \in \scrB$)}
            \end{equation*}
    Furthermore~$a \mapsto \varrho(a)$
        yields a
        nmiu-map~$\varrho\colon \scrA \to \scrB^a(\scrA\otimes_\varphi \scrB)$.
\item
    The map
    $h\colon \scrB^a(\scrA\otimes_\varphi\scrB) \to \scrB$
    given by~$h(T) =\left<1 \otimes 1, T(1 \otimes 1)\right>$
    is ncp.
\item
The map~$\varrho \colon \scrA \to \scrB^a(\scrA \otimes_\varphi \scrB)$
    together with~$1\otimes 1 \in \scrA \otimes_\varphi \scrB$
    has the following universal property.
\begin{quote}
Let~$\varrho'\colon \scrA \to \scrB^a(X)$
    be an nmiu-map
    for a self-dual Hilbert~$\scrB$-module~$X$
    together with an element~$e \in X$
    such that~$\varphi  = h' \after \varrho'$
        for~$h'(T) \equiv \left<e,T e\right>$.
Then: there is a unique
        inner product preserving~$\scrB$-linear
        map~$S \colon \scrA\otimes_\varphi\scrB \to X$
        such that~$\ad_S \after \varrho' = \varrho$
        and~$S( 1\otimes 1) = e$.
\end{quote}
\item
$(\scrB^a(\scrA \otimes_\varphi \scrB), \varrho, h)$
        is a Paschke dilation of~$\varphi$, see \sref{def-paschke}.
\end{enumerate}
\spacingfix{}
\begin{point}{40}{Proof}%
(This is a simplified version of the proof
    we published earlier in~\cite{wwpaschke}.
The construction of~$\scrA \otimes_\varphi \scrB$
    is essentially due to Paschke \cite[thm.~5.2]{paschke}:
    the self-dual completion of~$X$ in \cite[thm.~5.2]{paschke}
    is isomorphic to~$\scrA \otimes_\varphi \scrB$.)
\begin{point}{50}{1: $\scrA \otimes_\varphi \scrB$}%
The algebraic tensor product~$\scrA \odot \scrB$
    is a right $\scrB$-module
with the action~$(\sum_i a_i \otimes b_i)\beta = \sum_i a_i \otimes(b_i\beta)$.
On~$\scrA \odot \scrB$, define
\begin{equation*}
    \bigl[\sum_i a_i \otimes b_i, \sum_j \alpha_j \otimes \beta_j\bigr]
    \ \equiv\  \sum_{i,j} b_i^* \varphi(a_i^*\alpha_j)\beta_j.
\end{equation*}
By complete positivity of~$\varphi$
this is a~$\scrB$-valued inner product on~$\scrA \odot \scrB$.
By \sref{dils-completion},
    there is a self-dual Hilbert~$\scrB$-module~$\scrA\otimes_\varphi \scrB$
    and~$\scrB$-linear inner product-preserving
    $\eta\colon \scrA \odot \scrB \to \scrA \otimes_\varphi \scrB$
    with ultranorm-dense range.
Define a bilinear map~$\otimes \colon \scrA \times \scrB \to \scrA \otimes_\varphi \scrB$
    by~$a \otimes b = \eta(a \otimes b)$.
By definition we have
\begin{equation*}
\bigl\| \sum_i a_i \otimes b_i \bigr\|^2
    \ =\  \bigl\|\bigl[ \sum_i a_i\otimes b_i, \sum_j a_j \otimes b_j \bigr]\bigr\|
    \ =\  \bigl\|\sum_{i,j} b_i^* \varphi(a_i^*a_j) b_j\bigr\|
\end{equation*}
and so~$\otimes$ is a~$\varphi$-compatible bilinear map.

Let~$B \colon \scrA \times \scrB \to Y$
    be a~$\varphi$-compatible
    bilinear map to some self-dual Hilbert~$\scrB$-module~$Y$.
We must show that there is a unique bounded module
    map~$T \colon \scrA \otimes_\varphi \scrB \to Y$
    such that~$T(a\otimes b) = B(a,b)$ for all~$a \in \scrA$ and $b \in \scrB$.
By the defining property of the algebraic
    tensor product, there is a unique linear
    map~$T_0\colon \scrA \odot \scrB \to Y$
    such that~$T_0(a\otimes b) = B(a,b)$ for all~$a\in \scrA$ and~$b\in \scrB$.
By definition of~$\varphi$-compatibility
    and the inner product on~$\scrA\odot \scrB$,
    the map~$T_0$ is in fact bounded and~$\scrB$-linear.
By \sref{selfdual-completion-univ}
    this map extends uniquely
    to a bounded module map~$T \colon \scrA \otimes_\varphi \scrB \to Y$
    with~$T(a\otimes b) = B(a,b)$, as desired.
\end{point}
\begin{point}{60}{2: $\varrho \colon \scrA \to \scrB^a(\scrA \otimes_\varphi \scrB)$}%
Let~$a_0 \in \scrA$ be given.
To show the module map~$\varrho(a_0)$ exists,
    it suffices to show that the bilinear
    map~$B\colon \scrA \times \scrB \to \scrA \otimes_\varphi \scrB$
    given by~$B(a,b) = (a_0a)\otimes b$ is $\varphi$-compatible.
Clearly~$B(a,b \beta) = B(a,b)\beta$.
To prove~\eqref{phi-compatible},
    let~$n \in \N$, $a_1, \ldots, a_n \in \scrA$
    and~$b_1, \ldots, b_n \in \scrB$
    be given.
The row vector
    $(a_1\ \cdots \ a_n)$
    is an~$\scrA$-linear map~$a\colon \scrA^{\oplus n} \to \scrA$
    in the usual way:~$a(\alpha_1, \ldots, \alpha_n) = a_1\alpha_1 + \cdots + a_n \alpha_n$.
    Similarly the column vector~$(b_1\  \ldots\ b_n)^\T$
    is a~$\scrB$-linear map~$b\colon \scrB \to \scrB^{\oplus n}$
    with~$b(\beta) = (b_1 \beta, \ldots, b_n \beta)$.
    We compute
\begin{align*}
    \bigl\| \sum_i B(a_i,b_i) \bigr\|^2
    & \ =\  \bigl\| \sum_i (a_0 a_i) \otimes b_i \bigr\|^2 \\
    & \ =\  \bigl\| \sum_{i,j} b_i^* \varphi(a_i^* a_0^* a_0 a_j)b_j\bigr\| \\
    & \ =\   \| b^* (M_n \varphi) (a^*a_0^*a_0a) b\| \\
    &\  \leq\  \|a_0^*a_0\| \|b^* (M_n \varphi) (a^*a) b\| \\
    & \ =\  \|a_0\|^2  \bigl\| \sum_{i,j} b_i^* \varphi(a_i^*a_j)b_j\bigr\|.
\end{align*}
Thus~$B$ is~$\varphi$-compatible
    and so there is a unique~$\scrB$-linear bounded module
    map~$\varrho(a_0)\colon \scrA \otimes_\varphi \scrB
            \to \scrA \otimes_\varphi \scrB$
            with~$\varrho(a_0) (a\otimes b) = (a_0 a) \otimes b$.
Clearly~$\varrho$ is a unital, multiplicative and involution preserving map.
It remains to be shown~$\varrho$ is normal.
By~\sref{normal-faithful} and~\sref{hilmod-fixed-on-V},
    it suffices to show
    that~$ d \ \mapsto \ \langle \hat{x}, \varrho(d) \hat{x}\rangle$
    is normal for every~$x \in \scrA \odot \scrB$.
Find~$n \in \N$,
    row vector~$a \equiv(a_1\ \cdots\ a_n)$
    and column vector~$b \equiv (b_1\ \cdots\ b_n)$
    such that~$x = \sum^n_{i=1} a_i \otimes b_i$.
We compute
\begin{equation*}
    \langle \hat{x}, \varrho(d) \hat{x} \rangle
     \ =\  \sum_{i,j} b_i^* \varphi(a_i^* (d) a_j) b_j \\
         \ =\  b^* (M_n \varphi) (a^* \,d\, a) b.
\end{equation*}
So~$ d \ \mapsto \ \langle \hat{x}, \varrho(d) \hat{x}\rangle$
    is normal by~\sref{hilbmod-ad-ncp} and~\sref{mn-vna},
    whence~$\varrho$ is normal.
\end{point}
\begin{point}{70}{3: $h \colon \scrB^a(\scrA \otimes_\varphi \scrB) \to \scrB$}%
Define~$h(T) \equiv \left<1\otimes 1, T\, 1\otimes 1 \right>$.
It is completely positive by~\sref{hilbmod-vectstates-cp}
    and normal by~\sref{hilbmod-vecstates-normal}.
\end{point}
\begin{point}{80}[paschke-uniqueness]{Uniqueness~$\sigma$}%
Before we continue with 4, we will already prove the uniqueness property
    for point 5.
Note that~$(h \after \varrho)(a) = \left<1\otimes1, a \otimes 1\right>
    = \varphi(a)$ for all~$a \in \scrA$
    and so~$\varphi = h \after \varrho$.
Assume~$\varphi = h' \after \varrho'$
    for some nmiu-map~$\varrho'\colon \scrA \to \scrP'$,
        ncp-map~$h\colon \scrP' \to \scrB$
        and von Neumann algebra~$\scrP'$.
For point 5 we must show there is a
    unique ncp-map~$\sigma \colon \scrP' \to \scrB^a (\scrA \otimes_\varphi \scrB)$
    such that~$h \after \sigma = h'$ and~$\sigma \after \varrho' = \varrho$.
Let~$\sigma_1,\sigma_2\colon \scrP \to \scrB^a(\scrA\otimes_\varphi\scrB)$
    be ncp-maps with~$h \after \sigma_k = h'$ and~$\sigma_k \after \varrho'
        = \varrho$, $k=1,2$.
We must show~$\sigma_1=\sigma_2$.
Let~$c \in \scrP'$ and~$x \in \scrA \odot \scrB$ be given.
By \sref{hilmod-fixed-on-V} it suffices to prove
that~$\left<\hat{x}, \sigma_1(c) \hat{x} \right>= \left<\hat{x}, \sigma_2(c) \hat{x} \right>$.
Find~$n \in \N$, $a_1, \ldots, a_n\in \scrA$
    and~$b_1,\ldots,b_n \in \scrB$
    such that~$x = \sum_i a_i\otimes b_i$.
Note~$a_i \otimes b_i = \varrho(a_i) (1 \otimes 1)b_i$.
    So for any~$k=1,2$, we find
\begin{alignat*}{3}\label{equation-sigma}
    \left<\hat{x}, \sigma_k(c) \hat{x}\right>
    &\ =\ \sum_{i,j} b_i^* h(\, \varrho(a_i^*) \sigma_k(c) \varrho(a_j)\,)b_j\\
    &\ =\ \sum_{i,j} b_i^* h(\sigma_k(\, \varrho'(a_i^*) c \varrho'(a_j)\,))b_j &\qquad&\text{by \sref{dils-univlemma}}\\
    &\ =\ \sum_{i,j} b_i^* h'(\varrho'(a_i^*) c \varrho'(a_j)) b_j.\numberthis
\end{alignat*}
Thus~$\sigma_1=\sigma_2$ as desired.
One can show such~$\sigma$ exists by using~\eqref{equation-sigma}
    as defining formula.
We will use, however, an indirect but shorter approach
    by first considering the `spatial case',
    which was suggested by Michael Skeide.
\end{point}
\begin{point}{90}[paschke-spatial]{4: spatial case}%
Let~$X$ be a self-dual Hilbert~$\scrB$-module
    together with~$e \in X$ and nmiu-map~$\scrA \to \scrB^a(X)$
    such that~$\varphi = h' \after \varrho'$
    with~$h'(T) = \left<e,Te\right>$.
To prove uniqueness,
    assume that for~$k=1,2$,
    we have~$\scrB$-linear
    inner product preserving~$S_k \colon \scrA \otimes_\varphi \scrB \to X$
    with~$\ad_{S_k} \after \varrho' = \varrho$
    and~$S_k \, 1\otimes 1 = e$.
Then~$h \after \ad_{S_k} = h'$
    and so (as~$\ad_{S_k}$ is an ncp-map by \sref{hilbmod-ad-ncp}),
    we know~$\ad_{S_1} = \ad_{S_2}$ by~\sref{paschke-uniqueness}.
Hence there is a~$\lambda\in\C,\lambda\neq 0$
    with~$S_1 = \lambda S_2$, cf.~\cite[lemma 9]{westerbaan2016universal}.
If~$e \neq 0$,
    then from~$e = S_1\, 1\otimes 1 = \lambda S_2 \,1\otimes 1 = \lambda e$
    it follows~$\lambda = 1$ and so~$S_1 = S_2$.
In the other case, if~$e = 0$,
then~$h' = 0$ and so~$\varphi=0$,
    whence~$\scrA \otimes_\varphi \scrB = \{0\}$
    and indeed~$S_1=S_2=0$ trivially.
We continue with existence.
There is a unique linear~$S_0\colon \scrA \odot \scrB \to X$
    fixed by~$S_0(a\otimes b) = \varrho'(a) e b$.
Clearly~$S_0$ is also~$\scrB$-linear.
For any~$x,y \in \scrA\odot \scrB$,
    say~$x = \sum_i a_i\otimes b_i$ and~$y = \sum_j \alpha_j \otimes \beta_j$,
    we have
\begin{align*}
    \left< S_0 x, S_0 y \right>
    & \ =\ \sum_{i,j} \left<\varrho'(a_i)e b_i, \varrho'(\alpha_j)e \beta_j \right> \\
    & \ =\ \sum_{i,j} b_i^*\left<e,  \varrho'(a_i^*\alpha_j)e \right>\beta_j \\
    & \ =\ \sum_{i,j} b_i^* \varphi(a_i^*\alpha_j) \beta_j \\
    & \ =\ [x,y].
\end{align*}
Thus~$S_0$ preserves the inner product. In particular~$S_0$ is bounded
    and so there is a unique bounded~$\scrB$-linear
    $S \colon \scrA \otimes_\varphi \scrB \to X$
    with~$S \after \eta = S_0$.
    For all~$x \in \scrA \odot \scrB$
        we have~$\langle \hat{x}, S^*S \hat{x}\rangle = \langle \hat{x}, \hat{x} \rangle$
    and so~$S^*S=1$
    by~\sref{hilmod-fixed-on-V}.
Thus~$S$ preserves the inner product.
By definition~$S (1\otimes1) = \varrho'(1)e 1 = e$.
Pick any~$a \in \scrA$.
A straightforward computation shows~$S \varrho(a) = \varrho'(a) S$
hence~$S^* \varrho'(a) = \varrho(a) S^*$
and~$S^* \varrho'(a) S =  \varrho(a) S^*S= \varrho(a)$.
We have shown~$\ad_S \after \varrho' = \ad_S$.
\end{point}
\begin{point}{100}{5, $\sigma$ existence}
Assume~$\varphi = h' \after \varrho'$
    for some nmiu-map~$\varrho'\colon \scrA \to \scrP'$,
        ncp-map~$h\colon \scrP' \to \scrB$
        and von Neumann algebra~$\scrP'$.
It remains to be shown there is a ~$\sigma\colon \scrP' \to \scrB^a(\scrA\otimes_\varphi\scrB)$
with~$\sigma\after\varrho'=\varrho$ and~$h \after \sigma = h'$.
To apply the previous point,
    we perform the whole construction for~$h'$ instead of~$\varphi$
    yielding
    $h' = h_{h'} \after \varrho_{h'}$
    with~$\varrho_{h'} \colon \scrP' \to \scrB^a (\scrP' \otimes_{h'} \scrB)$
    and~$h_{h'} \colon \scrB^a(\scrP' \otimes_{h'} \scrB) \to \scrB$.
By~\sref{paschke-spatial}
    there is a unique 
    $\scrB$-linear map~$S\colon \scrA \otimes_{\varphi} \scrB \to
                                \scrP' \otimes_{h'} \scrB$
with~$S^*S = 1$, $S 1\otimes1 = 1\otimes1$ and
$\ad_S \after \varrho_{h'}\after \varrho' = \varrho$.
Define~$\sigma = \ad_S \after \varrho_{h'}$.
This~$\sigma$ fits the bill:
$\sigma \after \varrho' = \ad_S \after \varrho_{h'} \after \varrho'
                                = \varrho$
and~$h \after \sigma = h \after \ad_S \after \varrho_{h'}
= h_{h'} \after \varrho_{h'} = h'$. \qed
\end{point}
\end{point}
\end{point}
\end{parsec}

\begin{parsec}{1550}[ksgns]%
\begin{point}{10}%
We have shown that any ncp-map~$\varphi\colon \scrA \to \scrB$
    admits a dilation using Paschke's generalization of GNS
    to Hilbert~C$^*$-modules.
There is also a generalization of Stinespring's theorem
    to Hilbert~C$^*$-modules
    due to Kasparov \cite{ksgns}:
\end{point}
\begin{point}{20}{Theorem (KSGNS)}%
Let~$\varphi\colon \scrA \to \scrB^a(X)$
    be a cp-map
    for some~C$^*$-algebras~$\scrA$, $\scrB$
    and Hilbert~$\scrB$-module~$X$.
There exists a Hilbert~$\scrB$-module~$Y$,
    miu-map~$\varrho\colon \scrA \to \scrB^a(Y)$
    and adjointable~$\scrB$-linear~$T\colon Y\to X$
    with~$\varphi = \ad_T \after \varrho$.
\end{point}
\begin{point}{30}%
For our purposes (the dilation of an arbitrary ncp-map~$\scrA \to \scrB$),
    Paschke's GNS construction sufficed.
We leave open whether
    KSGNS admits a universal property like \sref{dils-univ-stinespring}
    and whether it is a Paschke dilation.
\end{point}
\end{parsec}
\begin{parsec}{1560}%
\begin{point}{10}%
We have all the tools to characterize
    the ncp-maps~$\varphi$
    for which the Paschke representation is injective.
This is a generalization of our answer\cite{stineinj}
    to the same question for the Stinespring embedding.
\end{point}
\begin{point}{20}[paschke-injective]{Theorem}%
Let~$\varphi\colon \scrA \to \scrB$ be any ncp-map
    with Paschke dilation~$(\scrP, \varrho, h)$.
Then $\ceil{\varrho} = \cceil{\varphi}$.
($\cceil{\varphi}$ is the central carrier of~$\varphi$, see~\sref{cceil-map-def}.)
Thus~$\varrho$ is injective ($\ceil{\varrho}=1$) if and only if $\varphi$
maps no non-zero central projection to zero ($\cceil{\varphi}=1$).
\begin{point}{30}{Proof}%
(This is a simplified version of the proof we published earlier
    in~\cite[thm.~30]{paschke}.)
By~\sref{paschke-unique-up-to-iso},
    it is sufficient to prove the equivalence
    for the dilation constructed in \sref{existence-paschke}.
Let~$p \in \scrA$ be any projection.
We will show~$p \leq \ceil{\varrho}^\perp$ iff~$p \leq \cceil{\varphi}^\perp$.
Assume~$p \leq \ceil{\varrho}$,
    viz.~$\varrho(p) = 0$.
Pick any~$a \in \scrA$ and~$b\in \scrB$.
Then~$\varrho(p)\, a\otimes b = 0$.
By definition of~$\varrho$ and the inner product,
    this is the case iff~$b^* \varphi(a^*p^* pa) b = 0$.
In particular~$\varphi(a^* p a) = 0$.
In other words~$a^* p a \leq \ceil{\varphi}^\perp$.
This happens if and only if~$a \ceil{\varphi} a^* \leq p^\perp$
    by~\sref{ad-contraposed}.
Equivalently~$\ceil{a \ceil{\varphi} a^*} \leq p^\perp$.
As~$a$ was arbitrary,
    we see (using \sref{cceil-fundamental}),
    that~$\cceil{\varphi} = \bigcup_{a \in \scrA} \ceil{a \ceil{\varphi} a^*}
            \leq p^\perp$, as desired.
For the converse, note that
    the two implications we just used are, in fact, equivalences:
        if~$\varphi(a^* p a) =0$,
        then~$b^*\varphi(a^* pa) b = 0$ for any~$b \in \scrB$
        \emph{and} if~$\varrho(p)\, a \otimes b = 0$
        for all~$a \in \scrA$ and~$b\in \scrB$,
        then~$\varrho(p) = 0$, by~\sref{hilmod-fixed-on-V}.
\qed
\end{point}
\end{point}
\end{parsec}

\begin{parsec}{1570}%
\begin{point}{10}%
We will now prove a useful connection between
    the ncp-maps \emph{ncp-below}~$\varphi$ (as defined in a moment)
    and the commutant of the image
    of a Paschke representation (``$\varrho$'') of~$\varphi$.
We will use it to give an alternative proof
    of the fact that pure maps and nmiu-maps
    are extreme among all ncpu-maps.
\end{point}
\begin{point}{20}{Definition}%
For linear  maps~$\varphi, \psi\colon \scrA \to \scrB$,
    we say $\varphi$ is ncp-below~$\psi$
    (in symbols:~$\Define{\varphi \mathrel{\leq_\ncp} \psi}$\index{$\leq_\ncp$})
    whenever~$\psi - \varphi$ is an ncp-map.
    Furthermore, write
    \begin{equation*}
        \Define{[0,\varphi]_\ncp}\index{$[0,\varphi]_\ncp$} =\  
        \{\psi; \ \psi \colon \scrA \to \scrB;
            \ 0 \leq_\ncp\psi \leq_\ncp \varphi \}.
    \end{equation*}
\end{point}
\spacingfix{}
\begin{point}{30}{Definition}%
    Let~$\varphi\colon \scrA \to \scrB$ be any ncp-map with Paschke
        dilation~$(\scrP, \varrho, h)$.
        For~$t \in \varrho(\scrA)^\square$, the commutant of~$\varrho(\scrA)$,
        define $\Define{\varphi_t}\index{$\varphi_t$} = h(t \varrho(a))$.
\end{point}
\begin{point}{31}%
The following is a generalization~\cite[prop.~5.4]{paschke}
    to arbitrary Paschke dilations.
\end{point}
\begin{point}{40}[paschke-correspondence]{Theorem}%
Assume~$\varphi\colon \scrA \to \scrB$ is an ncp-map
    with Paschke dilation~$(\scrP, \varrho, h)$.
The map~$t \mapsto \varphi_t$ is a linear
order isomorphism~$[0,1]_{\varrho(\scrA)^\square} \to [0,\varphi]_\ncp$.
\begin{point}{50}{Proof}
First we show the correspondence holds
    for the Paschke dilation constructed in~\sref{existence-paschke}.
This proof is a slight variation on \cite[prop.~5.4]{paschke}.
Then we show the correspondence carries over to arbitrary Paschke dilations.
\begin{point}{60}{Set-up}%
Let~$(\scrB^a(\scrA \otimes_\varphi \scrB), \varrho, h)$
    denote the Paschke dilation constructed in~\sref{existence-paschke}.
Clearly~$T \mapsto \varphi_T$ is linear
    for~$T \in \varrho(\scrA)^\square$.
Pick~$T \in \varrho(\scrA)^\square$.
As~$\sqrt{T} \in \varrho(\scrA)^\square$
    we have~$\varphi_T(a) = h(\sqrt{T} \varrho(a)\sqrt{T})$.
    Thus~$\varphi_T$ is ncp.
In particular, if~$T \leq S$
    for some~$S \in \varrho(\scrA)^\square$,
    then~$\varphi_{S-T} = \varphi_S - \varphi_T$ is ncp
    and so~$\varphi_T \mathrel{\leq_\ncp} \varphi_S$.
Consequently, if~$T\leq 1$, then $\varphi_T \leq \varphi_1 = \varphi$.
\end{point}
\begin{point}{70}{Order embedding}%
Let~$T \in \varrho(\scrA)^\square$
    be given such that~$\varphi_T$ is ncp.
We must show~$T \geq 0$.
Pick any~$x \in \scrA \odot \scrB$.
By~\sref{hilmod-fixed-on-V} and
    the construction of~$\scrA \otimes_\varphi \scrB$,
    it is sufficient to show~$\langle \hat{x}, T \hat{x}\rangle \geq 0$.
Say~$\hat{x} \equiv \sum_i a_i \otimes b_i$.
Then
\begin{equation*}
    \langle\hat{x}, T \hat{x} \rangle
    \ =\  \sum_{i,j} b_i^* \left< 1\otimes 1, T \varrho(a_i^*a_j) 1 \otimes 1\right> b_j
    \ = \ \sum_{i,j} b_i^* \varphi_T (a_i^*a_j) b_j \ \geq\  0.
\end{equation*}
Thus~$T \mapsto \varphi_T$ is an order embedding
    and in particular an injection.
\end{point}
\begin{point}{80}{Surjectivity}%
Pick any~$\psi \in [0,\varphi]_\ncp$.
We will show there is a~$T \in \varrho(\scrA)^\square$
    with~$\varphi_T = \psi$ and~$0 \leq T \leq 1$.
Let $(\scrB^a(\scrA \otimes_\psi \scrB),
    \varrho_\psi, h_\psi)$ denote the Paschke dilation of~$\psi$
    constructed in~\sref{existence-paschke}.
Define~$B\colon \scrA \times \scrB \to \scrA \otimes_\psi \scrB$
    by~$B(a,b) = a \otimes b$.
This~$B$ is $\varphi$-compatible
--- indeed, for any~$n \in \N$, $a_1, \ldots, a_n \in \scrA$
    and~$b_1, \ldots, b_n \in \scrB$,
    it follows from~$\psi \mathrel{\leq_\ncp} \varphi$ that
\begin{equation*}
    \bigl\| \sum_i B(a_i, b_i) \bigr\|^2
     \ = \  \bigl\| \sum_{i,j}
                b_i^* \psi(a_i^*a_j) b_j \bigr\|
     \ \leq \  \bigl\| \sum_{i,j}
                b_i^* \varphi(a_i^*a_j) b_j \bigr\|.
\end{equation*}
So by \sref{existence-paschke},
    there is a unique bounded module
    map~$W \colon \scrA \otimes_\varphi \scrB \to \scrA \otimes_\psi \scrB$
    fixed by~$W a \otimes b = a \otimes b$.
As~$\scrA \otimes_\psi \scrB$ is self-dual,
    the map~$W$ has an adjoint~$W^*$ by \sref{hilbmod-adjoint-exists}.
For any~$a\in\scrA$, we have~$W \varrho(a) = \varrho_\psi(a) W$
and so~$\varrho(a) W^* = W^* \varrho_\psi(a)$.
Thus~$W^*W \varrho(a) = W^* \varrho_\psi(a) W = \varrho(a) W^*W$.
Apparently~$W^*W \in \varrho(\scrA)^\square$.
Define~$T \equiv W^*W$. Clearly~$0\leq T$.
As~$\psi \mathrel{\leq_\ncp} \varphi$,
    we have for each~$x \equiv \sum_i a_i \otimes b_i$,
    the inequality~$\langle x,x\rangle_\psi \leq \langle x,x\rangle_\varphi$
    and so
\begin{equation*}
    \left<x, Tx\right>_\varphi \ =\  \left<Wx,Wx\right>_\psi
    \ =\  \left<x,x\right>_\psi \ \leq\  \langle x,x \rangle_\varphi.
\end{equation*}
    Thus~$T \leq 1$ by~\sref{hilmod-fixed-on-V}.
For any~$a\in \scrA$ we have
\begin{align*}
    \varphi_T(a) &\ =\  \left<1\otimes 1, W^*W\, a \otimes 1\right>_\varphi \\
    &  \ =\  \left<W\, 1\otimes 1, W\, a\otimes1\right>_\psi \\
    &  \ = \ \left<1\otimes1, a\otimes1\right>_\psi \\
    &  \ = \ \psi(a),
\end{align*}
    as desired.  We have shown the correspondence holds for
    the Paschke dilation~$(\scrB^a(\scrA \otimes_\varphi \scrB), \varrho, h)$
                constructed in~\sref{existence-paschke}.
\end{point}
\begin{point}{90}{Arbitrary dilation}%
Let~$(\scrP', \varrho', h')$ be any Paschke dilation
    of~$\varphi$.
For brevity write~$\scrP \equiv \scrA \otimes_\varphi \scrB$
    and for clarity
    write~$\varphi_T^{\scrP}$
        and~$\varphi_t^{\scrP'}$
    to distinguish between~$\varphi_T$ and~$\varphi_t$
    for~$T \in \varrho(\scrA)^\square$ and~$t \in \varrho'(\scrA)^\square$.
By~\sref{paschke-unique-up-to-iso}
    there is a unique nmiu-isomorphism
    $\vartheta\colon \scrP' \to \scrP$
    such that~$\vartheta \after \varrho' = \varrho $
    and~$h \after \vartheta = h'$.
It is easy to see~$\vartheta$ restricts to a linear order isomorphism
    $[0,1]_{\varrho'(\scrA)^\square} \to [0,1]_{\varrho(\scrA)^\square}$.
For any~$t \in \varrho'(\scrA)^\square$ and~$a\in \scrA$ we have
\begin{equation*}
    \varphi_t^{\scrP'} (a)
        \ = \ h'(t \varrho'(a))
        \ = \ h(\vartheta(t \varrho'(a)))
        \ = \ h(\vartheta(t) \varrho(a))
        \ = \ \varphi_{\vartheta(t)}^{\scrP}(a)
\end{equation*}
and so~$t \mapsto \varphi_t^{\scrP'}$
    is a linear order isomorphism as it
    is the composition of the linear order isomorphisms
    $\vartheta^{-1}$ and~$t \mapsto \varphi_t^\scrP$. \qed
\end{point}
\begin{point}{100}%
We return to~$t \mapsto \varphi_t$
    in~\sref{paschke-corresp-pure},
    where we connect it to extremeness.
\end{point}
\end{point}
\end{point}
\end{parsec}

\section{Kaplansky density theorem for Hilbert C$^*$-modules}
\begin{parsec}{1580}%
\begin{point}{10}%
Before we continue our study of Paschke dilations,
    we need to develop some more theory on (self-dual) Hilbert C$^*$-modules.
First we will prove a generalization of the Kaplansky density
    theorem for Hilbert C$^*$-modules,
    which we will need to compute the tensor product of Paschke dilations.

The original density theorem (\sref{kaplansky})
    is an important foundational result
    in the theory of operator algebras;
    quoting Pedersen:
    ``The density theorem is Kaplansky's great gift to mankind.
    It can be used every day, and twice on
    Sundays.'' \cite[\S2.3.4]{pedersen1979c}.
We used the original theorem already several times;
    notably to show that
        a von Neumann algebra is
        ultrastrongly complete~(\sref{vn-complete}
        via~\sref{vnsac}).
This ultrastrong completeness was used
    to define operator division~(\sref{division})
        (and in turn filters~\sref{filter}),
        to give the spatial tensor product a universal
        property~(\sref{tensor-universal-property})
    and in this thesis to show ultranorm completeness
        of self-dual Hilbert C$^*$-modules~(\sref{dils-selfdual}).

So what does is this density theorem again?
In essence, it is a remedy for the fact
    that strongly converging nets might not be bounded.
The theorem is usually stated as follows:
    the unit ball of the SOT-closure of a self-adjoint algebra
    of operators~$\scrA$
    is contained in the SOT-closure of the unit ball of~$\scrA$.
For our generalization, it is more convenient to consider the following
    variation.
\end{point}
\begin{point}{11}{Kaplansky density theorem}%
Let~$\scrB$ be a von Neumann algebra with ultrastrongly-dense
    C$^*$-subalgebra~$\scrA  \subseteq \scrB$.
Then: for every element~$a \in \scrA$,
    there is a net~$b_\alpha$ in~$\scrB$
    with~$b_\alpha \to a$ ultranorm
    and~$\|b_\alpha\| \leq \|b\|$ for all~$\alpha$.
\end{point}
\begin{point}{12}%
Now we are ready to state and prove our generalization,
    which is inspired by the proof of the regular Kaplansky density theorem
    given in \cite[thm.~1.2.2]{arveson2012invitation}.
Cf.~\sref{kaplansky}.
\end{point}
\begin{point}{20}[kaplansky-hilbmod]{Kaplansky density theorem for Hilbert C$^*$-modules}%
Let~$X$ be a Hilbert~$\scrB$-module for a von Neumann algebra~$\scrB$
    with an ultranorm-dense $\scrA$-submodule $D \subseteq X$,
    where~$\scrA \subseteq \scrB$ is some C$^*$-subalgebra
    with~$\langle y,y \rangle \in \scrA$ for all~$y \in D$.

Then: for every element~$x \in X$,
    there is a net~$x_\alpha$ in~$D$
    with~$x_\alpha \to x$ ultranorm
    and~$\|x_\alpha\| \leq \|x\|$ for all~$\alpha$.
\begin{point}{30}{Proof}%
Let~$x \in X$ be given.
The case~$x=0$ is trivial.
Assume~$x \neq 0$.
Without loss of generality, we may assume~$\| x \| = 1$.
For this~$x$ and any~$y \in X$, define
\begin{align*}
    h(y) &\ \equiv \ y \cdot\frac{2}{1+ \langle y, y \rangle}
    &
    g(x) &\ \equiv \ x \cdot\frac{1}{1+ \sqrt{1- \langle x, x\rangle}}.
\end{align*}
We claim~$h$ is ultranorm continuous,
$h(g(x)) = x$ and $\| h(y) \| \leq 1$ for all~$y \in X$.
From this, the promised result follows.
Indeed, let~$x_\alpha$ be a net in~$D$
    that converges ultranorm to~$g(x)$.
Then~$h(x_\alpha)$ is
    still in~$D$
    and furthermore~$h(x_\alpha) \to h(g(x)) = x$ ultranorm
    and~$\| h(x_\alpha) \| \leq 1 = \| x \|$, as desired.
\begin{point}{40}%
We start with~$\| h(y) \| \leq 1$.
The real map~$h_r\colon \lambda \mapsto 4 \lambda (1+\lambda)^{-2}$
    is bounded by~$1$ and so by the functional calculus
    we see~$\| h(y) \|^2 = \| 4 \langle y,y\rangle (1+ \langle y,y\rangle)^{-2} \|
                = \| h_r(\langle y,y\rangle) \|  \leq 1$.
We continue with the proof of~$h(g(x)) = x$.
Note~$\langle g(x), g(x) \rangle
= (1 + \sqrt{1 - \langle x,x \rangle })^{-2} \langle x,x\rangle$
and so with basic algebra
\begin{alignat*}{2}
    h(g(x)) &\ = \ x \cdot
\frac{1}{ 1 + \sqrt{1 - \langle x,x\rangle }}
\frac{2}{1 + \bigl(1 + \sqrt{1 - \langle x,x \rangle }\bigr)^{-2} \langle x,x\rangle} \\
& \ =\  x \cdot 2 \left( \frac{
    \bigl(1 + \sqrt{1-\langle x, x\rangle}\bigr)^2 + \langle x,x \rangle
}{1 + \sqrt{1 - \langle x,x \rangle}} \right)^{-1} \\
& \ = \ x.
\end{alignat*}
\end{point}
\spacingfix{}
\begin{point}{50}{Ultranorm continuity~$h$}%
The real work is in the proof of the ultranorm continuity of~$h$.
Suppose~$y_\alpha \to y$ ultranorm in~$X$. Define
\begin{align*}
    A_1 & \ \equiv \ \frac{\langle y,y\rangle}{(1 + \langle y,y \rangle)^2} 
            - \frac{1}{1+\langle y_\alpha,y_\alpha \rangle} \frac{\langle y,y \rangle}{1+\langle y,y \rangle} \\
    A_1' & \ \equiv \ 
    \frac{\langle y_\alpha,y_\alpha\rangle}{(1 + \langle y_\alpha,y_\alpha \rangle)^2}  
    - \frac{1}{1+\langle y,y \rangle}
    \frac{\langle y_\alpha,y_\alpha\rangle}{1+\langle y_\alpha,y_\alpha \rangle} \\
    A_2 & \ \equiv \ 
    \frac{1}{1+ \langle y,y \rangle}
       \langle y_\alpha- y, y_\alpha\rangle
     \frac{1}{1+ \langle y_\alpha,y_\alpha \rangle} \\
    A_2' & \ \equiv \  \frac{1}{1+ \langle y_\alpha,y_\alpha \rangle}
       \langle y- y_\alpha, y\rangle
     \frac{1}{1+ \langle y,y \rangle}.
\end{align*}
A straightforward computations shows
\begin{equation}\label{kaplansky-splitting}
    \langle h(y) - h(y_\alpha) ,
         h(y) - h(y_\alpha) \rangle  \ = \ A_1 + A_1' + A_2 + A_2'.
\end{equation}
We will show that each of these terms converges ultraweakly to~$0$
 as~$\alpha \to \infty$.
We start with~$A_1$.
We need two facts.
First, for any positive~$b \in \scrB$, we have~$0 \leq \frac{1}{1+b} \leq 1$
and~$0 \leq \frac{b}{1+b} \leq 1$.
Secondly
\begin{align*}
    \frac{1}{1+\langle y,y\rangle }
        -
    \frac{1}{1+\langle y_\alpha,y_\alpha\rangle }
    & \ = \ 
    \frac{1}{1+\langle y,y\rangle }
    \bigl((1+\langle y_\alpha,y_\alpha \rangle) - (1+\langle y,y \rangle)\bigr) 
    \frac{1}{1+\langle y_\alpha,y_\alpha\rangle }\\
    & \ = \ 
    \frac{1}{1+\langle y,y\rangle }
   \bigl(\langle y_\alpha,y_\alpha \rangle - \langle y,y \rangle\bigr) 
    \frac{1}{1+\langle y_\alpha,y_\alpha\rangle }.
\end{align*}
Putting these facts and the definition of~$A_1$ together, we get
\begin{align*}
     A_1 
    &\ = \ 
    \left(
    \frac{1}{1+\langle y,y\rangle } -
    \frac{1}{1+\langle y_\alpha,y_\alpha\rangle }\right)
            \frac{\langle y,y \rangle}{1+\langle y,y \rangle} \\
    &\ = \ 
    \frac{1}{1+\langle y,y\rangle }
   \bigl(\langle y_\alpha,y_\alpha \rangle - \langle y,y \rangle\bigr) 
    \frac{1}{1+\langle y_\alpha,y_\alpha\rangle }
            \frac{\langle y,y \rangle}{1+\langle y,y \rangle} \\
    &\ = \ 
    \frac{1}{1+\langle y,y\rangle }
   \bigl(\langle y_\alpha -y,y_\alpha \rangle + \langle y,y_\alpha - y \rangle\bigr) 
    \frac{1}{1+\langle y_\alpha,y_\alpha\rangle }
            \frac{\langle y,y \rangle}{1+\langle y,y \rangle}.
\end{align*}
As multiplying with
constants is ultraweakly continuous (see \sref{mult-uws-cont}),
    it is sufficient to show that ultraweakly, as~$\alpha \to \infty$, we have
\begin{align}\label{kaplanskytodo2}
    \langle y_\alpha - y, y_\alpha \rangle \frac{1}{1 + \langle y_\alpha, y_\alpha\rangle} & \ \to\  0 &
    \langle y, y_\alpha  -y\rangle \frac{1}{1 + \langle y_\alpha, y_\alpha\rangle} & \ \to \ 0.
\end{align}
For any np-map~$f\colon \scrB \to \C$, we have
\begin{align*}
   \Bigl| f\Bigl(\left\langle y_\alpha - y, y_\alpha \right\rangle \frac{1}{1+\langle y_\alpha, y_\alpha\rangle}\Bigr)\Bigr|^2
    & \ = \ 
    |\langle y_\alpha - y, y_\alpha (1 + \langle y_\alpha, y_\alpha  \rangle)^{-1}
    \rangle_f|^2 \\
    & \ \leq \ 
    \| y_\alpha - y \|^2_f \cdot \| y_\alpha (1 + \langle y_\alpha, y_\alpha  \rangle)^{-1} \|^2_f \\
    & \ \leq \  \|
    y_\alpha - y \|^2_f \cdot \| y_\alpha (1 + \langle y_\alpha, y_\alpha  \rangle)^{-1} \|^2 \cdot
    \|f\|
    \\
    & \ \leq \ 
    \| y_\alpha - y \|^2_f \cdot \|f\| \ \rightarrow \ 0.
\end{align*}
We have shown the LHS of \eqref{kaplanskytodo2}.
The proof for the RHS is different, but simpler.
So~$A_1$ vanishes ultraweakly.
In a similar way one sees~$A_1' \to 0$ ultraweakly.
The proofs for~$A_2, A_2' \to 0$ are very similar. \qed
\end{point}
\end{point}
\end{point}
\end{parsec}

\section{More on self-dual Hilbert C$^*$-modules}
\begin{parsec}{1590}%
\begin{point}{10}%
We have a second look at
    self-dual Hilbert C$^*$-modules over von Neumann algebras
    and their orthonormal bases.
We start with some good news.
In~$\scrB(\scrH)$ the linear span of~$\ketbra{e_i}{e_j}$
    is ultraweakly dense for any orthonormal basis~$(e_i)_{i \in I}$
    of~$\scrH$.
For self-dual Hilbert~$\scrB$-modules, we have a similar result.
\end{point}
\begin{point}{20}{Definition}%
Let~$X$ be a Hilbert~$\scrB$-module.
For any~$x, y \in X$,
define the bounded operator~$\ketbra{x}{y} \in \scrB^a(X)$
by~$\Define{\ketbra{x}{y}}\index{*ketbra@$\ketbra{x}{y}$ for Hilbert $\scrB$-modules} z = x \left<y,z\right>$.
\begin{point}{21}{Remark}%
The operator~$\ketbra{x}{y}$ is usually denoted by~$\theta_{x,y}$
    in the literature,
    see for instance~\cite[\S2.2]{manuilov2000hilbertc}.
\end{point}
\begin{point}{30}[hilbmodketbrarules]%
The following rules are easy to check.
For~$x,y,v,w \in X$ and~$b \in \scrB$
    we have
\begin{equation*}
\ketbra{x}{y}^* = \ketbra{y}{x} \qquad
\ketbra{xb}{y} = \ketbra{x}{yb^*} \qquad
\ketbra{x}{y} \ketbra{v}{w} = \ketbra{x\left<y,v\right>}{w}.
\end{equation*}
If~$\left<e,e\right>$ is a projection,
    then~$\ketbra{e}{e}$ is a projection
    (using \sref{mod-projelabs}).
For any~$T \in \scrB^a(X)$
    we have~$T \ketbra{x}{y} = \ketbra{Tx}{y}$
    and~$\ketbra{x}{y}T^* = \ketbra{x}{Ty}$.
\end{point}
\end{point}
\begin{point}{40}[ketbra-ultraweakly-dense]{Proposition}%
Let~$X$ be a self-dual Hilbert~$\scrB$-module
    for a von Neumann algebra~$\scrB$.
If~$(e_i)_{i \in I}$ is an orthonormal basis of~$X$,
then the linear span of
    the operators~$\{ \ketbra{e_i b}{e_j}; i,j \in I, \ b \in \scrB\}$
    is ultraweakly dense in~$\scrB^a(X)$.
\begin{point}{50}{Proof}%
We start with some preparation.
\begin{point}{60}[ketbra-dense-pt1]%
For a finite subset~$S \subseteq I$,
    write~$p_S = \sum_{i \in S} \ketbra{e_i}{e_i}$.
The~$(\ketbra{e_i}{e_i})_{i \in I}$ are pairwise orthogonal projections.
Thus~$p_S$ is a projection.
Also $p_S \leq p_{S'}$ when~$S \subseteq S'$.
Pick any~$x \in X$.
As~$p_S x = \sum_{i \in S} e_i \left<e_i, x\right>$, we have
\begin{equation*}
\left<x, p_S x\right> \ =\ 
        \left<p_S x, p_S x\right> \ =\ 
        \sum_{i,j \in S}
            \left<x,e_i\right> \left<e_i,e_j\right> \left<e_j,x\right>
            \ =\ 
        \sum_{i\in S} \left<x,e_i\right>\left<e_i,x\right>,
\end{equation*}
and so by Parseval~$\sup_S \left<x, p_S x\right> = \left<x,x\right>$.
Hence~$\sup_S p_S = 1$ by \sref{states-order-separating}
    and so~$p_S$ converges ultraweakly to~$1$
    by \sref{vna-supremum-uslimit}.
\end{point}
\begin{point}{70}%
Pick any~$T \in \scrB^a(X)$.
For~$i,j \in I$ we have
\begin{equation*}
    \ketbra{e_i}{e_i} T \ketbra{e_j}{e_j}
    \ =\  \ketbra{e_i}{e_i} \ketbra{T e_j}{e_j}
    \ =\ \ketbra{e_i \left<e_i, T e_j\right>}{e_j}.
\end{equation*}
and so~$p_S T p_S$ is in the linear
    span of~$\{ \ketbra{e_i b}{e_j} \}$.
Thus it is sufficient to prove~$p_S T p_S$ converges ultraweakly to~$T$.
\end{point}
\begin{point}{80}[err159IV]%
Pick any np-map~$f\colon \scrB \to \C$.
We want to show~$|f(T-p_S T p_S)| \to 0$. Clearly
\begin{equation*}
    |f(T - p_S T p_S)| \ \leq \  |f(\, (1-p_S) T \,) | +| f(\, p_S T (1-p_S) \,)|.
\end{equation*}
Using Cauchy--Schwarz and \sref{ketbra-dense-pt1} we see the second term
    vanishes
\begin{equation*}
|f(p_S T (1-p_S))|^2 \ \leq \ f(p_S T T^* p_S) f(1-p_S) \ \leq\ 
            \| T \|^2 f(1) f(1-p_S) \ \to \ 0.
\end{equation*}
Similarly the first term vanishes and
    so indeed~$p_S T p_S \to T$ ultraweakly. \qed
\end{point}
\end{point}
\end{point}
\begin{point}{90}[ketbra-ultranorm-continuous]{Proposition}%
Let~$X$ be a self-dual Hilbert~$\scrB$-module for a von Neumann algebra~$\scrB$.
If~$x_\alpha \to x$ ultranorm for a norm-bounded net~$x_\alpha$ in~$X$,
    then~$\ketbra{x_\alpha}{y} \to \ketbra{x}{y}$ ultraweakly
    for any~$y \in Y$.
\begin{point}{100}{Proof}%
We need some preparation.
Define
\begin{equation*}
    \Omega \ \equiv \  \{ f(\langle x, (\,\cdot\,) x \rangle); \ 
    f \colon \scrB \to \C \text{ np-map},\ x \in X\}.
\end{equation*}
From~\sref{hilbmod-ordersep}
    it follows~$\Omega$ is order separating.
We will use this fact twice.
First we will prove~$\| \ketbra{z}{y} \| \leq \| z\| \|y\|$.
For any~$\omega \in \Omega$,
    say~$\omega \equiv f(\langle x,(\,\cdot\,) x \rangle)$, we have
\begin{equation*}
    \omega(\ketbra{y}{y})
     \ = \ f(\langle x,y \rangle\langle y, x \rangle) 
     \ \leq \ f(\langle x, x\rangle) \| y\|^2 
     \ = \ \| \omega \| \| y\|^2
\end{equation*}
and so from \sref{order-separating-norm} it follows~$\| \ketbra{y}{y} \| \leq \|y\|^2$. Whence
\begin{align*}
    \| \ketbra{z}{y} \|^2
        & \ = \ \| \ketbra{y}{z} \ketbra{z}{y} \| \\
        & \ = \ \| \ketbra{\smash{y \langle z,z\rangle^{\frac{1}{2}}}}{
\smash{y \langle z,z\rangle^{\frac{1}{2}}}} \| \\
        & \ \leq \ \| y \langle z,z \rangle^{\frac{1}{2}} \|^2.
\end{align*}
Hence~$\| \ketbra{z}{y} \| \leq \|y \langle z,z\rangle^{\frac{1}{2}} \|
                            \leq \|y\|  \|z\|$, as promised.
\begin{point}{110}%
Now we start the proof proper.
As~$\ketbra{\,\cdot\,}{y}$ is linear
    and~$x_\alpha - x$ is norm-bounded,
    we may assume without loss of generality
    that~$x_\alpha \to 0$ ultranorm.
    We have to prove~$\ketbra{x_\alpha}{y} \to 0$ ultraweakly.
Let~$g\colon \scrB^a(X) \to \C$ be an arbitrary np-map.
Using the order separation of~$\Omega$ a second time
        and the fact that~$\omega(T^* (\,\cdot\,)T) \in \Omega$
        for all~$T \in \scrB^a(X)$ and~$\omega \in \Omega$,
        we see that the linear span of~$\Omega$
    is operator norm dense among  all np-maps~$\scrB^a(X) \to \C$
    by \sref{vn-center-separating-fundamental}.
So there are~$n \in \N$ and~$\omega_1, \ldots, \omega_n \in \Omega$
    with~$\| g - g' \| \leq \varepsilon$,
    where~$g' \equiv \sum^n_{i=1} \omega_i$.
By unfolding definitions and~\sref{innerprod-ultraweak},
    it follows easily that~$\omega_i(\ketbra{\,\cdot\,}{y})$
    is ultranorm continuous.
    So~$g'(\ketbra{\,\cdot\,}{y})$ is ultranorm continuous as well.
Pick any~$\varepsilon > 0$.
By the previous,
    there is an~$\alpha_0$ such that for all~$\alpha \geq \alpha_0$
    we have~$|g'(\ketbra{x_\alpha}{y})| \leq \varepsilon$,
    hence
\begin{align*}
        | g(\ketbra{x_\alpha}{y}) |
            &\ \leq\  | (g-g')(\ketbra{x_\alpha}{y})|  +
                    | g'(\ketbra{x_\alpha}{y}) | \\
            &\ \leq\  \|g - g'\| \cdot \| \ketbra{x_\alpha}{y} \|  + \varepsilon\\
            &\ \leq\  \varepsilon \,(\| x_\alpha\| \cdot \|y\| + 1).
\end{align*}
As~$x_\alpha$ is norm-bounded
    and~$\varepsilon$ was arbitrary,
        we conclude~$g(\ketbra{x_\alpha}{y}) \to 0$.
    Thus~$\ketbra{x_\alpha}{y} \to 0$ ultraweakly, as desired.
    \qed
\end{point}
\end{point}
\end{point}
\end{parsec}

\begin{parsec}{1600}%
\begin{point}{10}%
For a Hilbert space~$\scrH$
    with linear subspace~$V \subseteq \scrH$
    we have~$V \subseteq V^{\perp\perp} = \overline{V}$
    and~$\scrH \cong V^{\perp\perp} \oplus  V^\perp$.
We generalize this to self-dual Hilbert~C$^*$-modules
    over von Neumann algebras.
Let's introduce the cast:
\end{point}
\begin{point}{20}[direct-prod-self-dual-basis]{Exercise}%
Let~$\scrB$ be a von Neumann algebra.
Assume~$X$ and~$Y$ are self-dual Hilbert~$\scrB$-modules
    with orthonormal bases~$E \subseteq X$ and~$F \subseteq Y$.
Write~$\kappa_1\colon X \to X \oplus Y$
    and~$\kappa_2 \colon Y \to X \oplus Y$
    for the maps~$\kappa_1(x) = (x,0)$ and~$\kappa_2(y) = (0, y)$.
Show that~$\kappa_1(E) \cup \kappa_2(F)$
    is an orthonormal basis of~$X \oplus Y$.
    Conclude that~$X \oplus Y$ is self-dual.
\end{point}
\begin{point}{30}{Definition}%
Let~$X$ be a Hilbert~C$^*$-module
    with subset~$V \subseteq X$.
    Write~$\Define{V^\perp}$ for the \Define{orthocomplement}
    \index{orthocomplement!w.r.t.~a Hilbert~$\scrB$-module}
    \index{*perp@$(\ )^\perp$!subset Hilbert~$\scrB$-module}
    of~$V$, defined by
\begin{equation*}
    V^\perp = \{ x; \ x \in X;\ \left<x,v\right>=0\text{ for all }v\in V\}.
\end{equation*}
\end{point}
\spacingfix{}
\begin{point}{40}[hilbmod-projthm]{Proposition}%
Assume~$X$ is a self-dual Hilbert $\scrB$-module
    for a von Neumann algebra~$\scrB$.
Let~$V \subseteq X$ be any subset. Then
\begin{enumerate}
\item
$V^\perp$ is a ultranorm-closed Hilbert~$\scrB$-submodule of~$X$
    (and so is~$V^{\perp\perp})$;
\item
    $V^{\perp\perp}$ is the ultranorm closure of the~$\scrB$-linear
    span of~$V$ \emph{and}
\item
    $V^{\perp\perp} \oplus V^\perp \cong X$ as Hilbert~C$^*$-modules
        via~$(x,y) \mapsto x+y$.
\end{enumerate}
\spacingfix{}
\begin{point}{50}{Proof}%
It is easy to see~$V^\perp$ is a submodule of~$X$
    and~$V \subseteq V^{\perp\perp}$.
To show~$V^\perp$ is ultranorm closed in~$X$,
    assume~$x_\alpha$ is a net in~$V^\perp$
    converging ultranorm to~$x \in X$.
For each~$v \in V$ we have
$ \left<v, \unlim_\alpha x_\alpha\right>
    \overset{\sref{ultranormcontstruct}}{=}
    \uslim_\alpha \left<v, x_\alpha\right>  = 0 $
    and so~$x \in V^\perp$ as well.
Thus~$V^\perp$ is ultranorm complete
    (and so self dual by \sref{dils-selfdual}).
\begin{point}{60}%
Write~$W$ for the ultranorm closure of
    the~$\scrB$-linear span of~$V$ in~$X$.
It follows from \sref{ultranormcontstruct}, that $W$ is a submodule.
Clearly the induced uniformity of~$X$
    on~$W$ is the same as its own ultranorm uniformity:
    thus~$W$ is ultranorm complete
        and~$W \subseteq V^{\perp\perp}$.
By \sref{dils-selfdual} there
    is an orthonormal basis~$(e_i)_{i \in I}$ of~$W$.
Going back to the construction of~$(e_i)_{i \in I}$
    (in \sref{selfdual-bcompl-then-basis}),
    we see that we can extend it to a maximal orthonormal subset of~$X$,
    which must be a basis for the whole of~$X$.
Pick any~$(d_j)_{j \in J} $ in~$X$
    such that~$\{e_i\}_{i \in I} \cup \{d_j\}_{j \in J}$
    is an orthonormal basis of~$X$.
\end{point}
\begin{point}{70}%
By construction~$d_j \in W^\perp \subseteq V^\perp$
    for every~$j \in J$.
For any~$x \in X$
    we have~$x = \sum_i e_i \left<e_i, x\right>
                + \sum_j d_j \left<d_j,x\right>$.
If~$x \in V^{\perp\perp}$
    then~$x = \sum_i e_i \left<e_i, x\right>$
    as~$\left<x, d_j\right>= 0$ for each~$j \in J$.
Thus~$V^{\perp\perp} \subseteq W$, so~$V^{\perp\perp} = W$.
In particular, as~$V^\perp$ is already an ultranorm closed submodule,
    we find~$W^\perp = (V^\perp)^{\perp\perp} = V^\perp$.
\end{point}
\begin{point}{80}%
So~$(e_i)_i$ is a basis for~$W = V^{\perp\perp}$
    and~$(d_j)_j$ is a basis for~$W^\perp = V^\perp$.
For brevity, let~$\vartheta\colon V^{\perp\perp} \oplus V^\perp \to X$
    denote the map~$(y,y') \mapsto y+y'$.
Clearly~$\vartheta$ is bounded~$\scrB$-linear.
As~$\left<y,y'\right> = 0$ for~$(y,y') \in V^{\perp\perp} \oplus V^\perp$
    we easily see~$\vartheta$ is inner product preserving and thus injective.
As~$\{e_i\}_{i \in I} \cup \{d_j\}_{j \in J}$
    is an orthonormal basis of~$X$, $\vartheta$ is surjective.
As~$\vartheta$ is inner product preserving, its inverse is bounded
    --- hence~$\vartheta$ is an isomorphism of~Hilbert C$^*$-modules,
    as promised.
    \qed
\end{point}
\end{point}
\end{point}
\begin{point}{90}[selfdual-orthn-basis]{Exercise}%
Assume~$X$ is a self-dual Hilbert~$\scrB$-module for a von Neumann
    algebra~$\scrB$.
Suppose~$E \subseteq X$ is an orthonormal set. Show
\begin{enumerate}
\item
    $E$ is a basis of~$E^{\perp\perp}$ \emph{and}
\item
    for any~$x \in X$,
        we have
    $x \in E^{\perp\perp}$
    if and only if~$\left<x,x\right> = \sum_{e \in E} \left<x,e\right>
                \left<e,x\right>$.
\end{enumerate}
\spacingfix{}
\end{point}
\begin{point}{100}[selfdual-gramschmidt]{Exercise}%
Let~$X$ be a self-dual Hilbert~$\scrB$-module for some von Neumann
    algebra~$\scrB$.
Show that for any~$x_1, \ldots, x_n \in X$,
    there is a finite orthonormal basis
    of~$\{x_1, \ldots, x_n\}^{\perp\perp}$
    consisting of at most~$n$ elements.
(Use the orthonormalization in the last part of \sref{selfdual-bcompl-then-basis}.)
\end{point}
\end{parsec}

\begin{parsec}{1610}%
\begin{point}{10}[thel2matter]%
Every Hilbert space~$\scrH$ is isomorphic to some~$\ell^2(I)$.
In fact, the cardinality of~$I$ is the only thing that matters,
    in the following sense:
    if~$\ell^2(I) \cong \ell^2(J)$, then~$I$ and~$J$
    have the same cardinality.
What about Hilbert C$^*$-modules?
\end{point}
\begin{point}{20}[hilbmod-el2]{Exercise}%
In this exercise we prove \cite[thm.~3.12]{paschke}.
Let~$\scrB$ be a von Neumann algebra
    and~$(p_i)_{i \in I}$
    a family of projections from~$\scrB$.
    Write~$\Define{\ell^2((p_i)_{i \in I})}$\index{$\ell^2((p_i)_{i \in I})$}
    for the set of~$I$-tuples~$(b_i)_{i \in I}$
    from~$\scrB$ that are~$\ell^2$-summable (\sref{dfn-selfdual-basis})
    and furthermore~$\ceil{b_ib_i^*} \leq p_i$ for every~$i \in I$.
Show that~$\ell^2((p_i)_{i \in I})$
    is a right $\scrB$-module
    with coordinatewise operations.
Prove that for any~$(b_i)_i$ and~$(c_i)_i$
    in~$\ell^2((p_i)_{i \in I})$
    the sum~$\sum_{i \in I} b_i^* c_i =: \left<(b_i)_i, (c_i)_i\right>$
    converges ultraweakly
    and turns~$\ell^2((p_i)_{i\in I})$ into a pre-Hilbert~$\scrB$-module.
    Conclude~$\ell^2((p_i)_{i \in I})$ is self-dual
and, in fact, that for every self-dual Hilbert~$\scrB$-module~$X$
    with any orthonormal basis~$(e_i)_{i \in I}$,
    we have~$X \cong \ell^2(\,(\left<e_i,e_i\right>)_{i \in I}\,)$.
\end{point}

\begin{point}{30}[hilbel-matter]%
    Let~$\scrB$ be a von Neumann algebra
        and~$(p_i)_{i \in I}$, $(q_j)_{j \in J}$
        families of projections from~$\scrB$.
    What can we say if~$\ell^2((p_i)_{i \in I}) 
    \cong \ell^2((q_i)_{j \in J})$?
    ($\ell^2$ as defined in \sref{hilbmod-el2}.)
    Let's start with some counterexamples.
\begin{enumerate}
\item
Do all orthonormal bases have the same cardinality?
$\scrB$ is a self-dual Hilbert~$\scrB$-module over itself.
Clearly, the element~$1$ by itself is an orthonormal basis of~$\scrB$.
However, for every projection~$p \in \scrB$
        the pair~$\{ p,1-p \}$ is also an orthonormal basis.
    Worse still, for~$\scrB \equiv \scrB(\ell^2(\N))$,
        the set
\begin{equation*}
    \{\ketbra{0}{0}, \ \ketbra{1}{1},\  \ketbra{2}{2},\  \ldots \}.
\end{equation*}
is an orthonormal basis.

\item
Perhaps any two orthonormal bases have a common coarsening in some sense?
Take~$\scrB = M_3$.
It might not be obvious at this moment, but both
\begin{equation*}
    \{\,\ketbra{0}{0} + \ketbra{2}{2}, \ketbra{1}{+}\,\}
    \quad \text{and} \quad
    \{\,\ketbra{0}{0}, \ketbra{2}{2}+ \ketbra{1}{+}\,\}
\end{equation*}
are orthonormal bases of~$M_3$ over itself.
(Combine \sref{onb1} with~\sref{onb2}.)
There does not seem to be an obvious way in which these
    two bases have a common coarsening of sorts.

\item
Maybe~$\sum_{i\in I} \left<e_i,e_i\right> = \sum_{j\in J} \left<d_j,d_j\right>$
    for finite orthonormal bases~$(e_i)_{i \in I}$ and~$(d_j)_{j \in J}$?
For commutative~$\scrB$ this indeed holds using Parseval's identity.
However, in general it fails:
        both~$\{\ketbra{0}{0}, \ketbra{1}{1}\}$ and
        $\{\ketbra{0}{0}, \ketbra{1}{0} \}$
        are orthonormal bases of~$M_2$ over itself,
    but~$\ketbra{0}{0}+ \ketbra{1}{1} \neq 2 \ketbra{0}{0}$.
\end{enumerate}
\spacingfix{}
\end{point}
\begin{point}{40}[onb1]{Exercise}%
Let~$(e_i)_{ i \in I}$ be some orthonormal basis of a Hilbert~$\scrB$-module~$X$.
Show that if~$(u_i)_{i \in I}$
    is a family of partial isometries
    from~$\scrB$
    (see \sref{partial-isometry-equivalents})
    such that~$u_iu_i^* = \left<e_i,e_i\right>$
    for all~$i \in I$,
    then~$(e_iu_i)_{i \in I}$
    is an orthonormal basis of~$X$.

Use this and~\sref{hilbmod-el2} to conclude
    that~$\ell^2((p_i)_{i \in I}) \cong \ell^2((q_i)_{i \in I})$
    for any projections~$(p_i)_{i\in I}, (q_i)_{i\in I} $
    from~$\scrB$
    with~$p_i \sim q_i$ for~$i \in I$.
    (Here, $\sim$ denotes the \Define{Murray--von Neumann equivalence}.
    That is: \Define{$ p \sim q$} iff there is a partial isometry~$u$
    with~$u^*u = p$ and~$uu^*=q$.
    Eg.~$\ketbra{0}{0} \sim \ketbra{+}{+}$
        via~$\ketbra{0}{+}$.)
\end{point}

\begin{point}{50}[onb2]{Exercise}%
Let~$(e_i)_{ i \in I}$ be some orthonormal basis of a Hilbert~$\scrB$-module~$X$
    with distinguished~$1,2 \in I$.
Show that if~$\left<e_1,e_1\right> + \left<e_2,e_2\right> \leq 1$,
    then we can make a new orthonormal basis of~$X$ by  removing~$e_1$ and~$e_2$
        and inserting~$e_1 + e_2$.

Conclude~$p \scrB \oplus q \scrB \cong (p+q) \scrB$
    for projections~$p,q \in \scrB$ with $p+q \leq 1$.
\end{point}
\end{parsec}
\begin{parsec}{1620}%
\begin{point}{10}%
Before we can prove our normal form (in \sref{selfdual-normalish-form}),
    we need some more \emph{comparison theory of projections},
    which was touched upon in~\sref{vmleq}.
Recall that for projections~$p,q$ in a von Neumann algebra~$\scrA$,
    we write~$\Define{p \lesssim q}$ if there is a partial isometry~$u$
    with~$u^*u=p$ and~$uu^* \leq q$.\index{*lesssim@$\lesssim$}
In~\sref{mvn-preorders} we saw~$\lesssim$ preorders the projections.
\end{point}
\begin{point}{20}[total-mv-order]{Proposition}%
Let~$\scrB$ be a von Neumann algebra that is a factor (i.e.~with center~$\C 1$).
For any two projections~$p,q \in \scrB$,
    either~$p \lesssim q$ or~$q \lesssim p$ (see \sref{vmleq}).
\begin{point}{30}{Proof}%
(For a more traditional proof, see e.g.~\cite[prop.~6.2.4]{kr}.)
Let~$p,q \in \scrB$ be any projections.
By Zorn's lemma, there is a maximal
    set~$U\subseteq \scrB$ of partial isometries
    with~$\sum_{u\in U} u^*u \leq p$
    and~$\sum_{u\in U} uu^* \leq q$.
Define~$p_0 \equiv p - \sum_{u\in U} u^*u$
    and~$q_0 \equiv q - \sum_{u \in U} uu^*$.
Note that if~$p_0 = 0$, then~$p \lesssim q$
    via~$\sum_{u \in U} u$, c.f.~\sref{summing-partial-isometries}.
Similarly, if~$q_0 = 0$, then~$q \lesssim p$.
So, reasoning towards contradiction,
    suppose~$p_0 \neq 0$ and~$q_0 \neq 0$.
As~$\scrB$ is a factor~$\cceil{q_0} = 1$
    and so~$0 < p_0 = p_0\cceil{q_0} p_0 = p_0 (\bigcup_{a \in \scrB} \ceil{a^* q_0 a }) p_0
                = \bigcup_{a \in \scrB} \ceil{p_0a^*q_0ap_0}$
            by~\sref{cceil-fundamental} and~\sref{ncp-union}.
Pick an~$a \in \scrB$ with~$\ceil{p_0a^*q_0ap_0} \neq 0$.
Then~$u_0 \equiv [q_0ap_0]$, see \sref{polar-decomposition},
    is a non-zero partial isometry
    satisfying~$u_0^*u_0 = \ceil{p_0a^*q_0ap_0} \leq p_0$
    and~$u_0 u_0^* = \ceil{q_0ap_0a^*q_0} \leq q_0$,
    contradicting maximality of~$U$. \qed
\end{point}
\end{point}
\begin{point}{40}[selfdual-normalish-form]{Theorem}%
Let~$\scrB$ be a von Neumann algebra that is a factor (i.e.~with center~$\C 1$).
Suppose~$X$ is a self-dual Hilbert~$\scrB$-module.
Either there is an infinite cardinal~$\kappa$
    such that~$X \cong \ell^2((1)_{\alpha \in \kappa})$,
    with~$\ell^2$ as in~\sref{hilbmod-el2},
    or there is a natural number~$n\in \N$
    and a projection~$p \in \scrB$
    such that~$X \cong \ell^2((1, \ldots, 1, p))$,
    where~$1$ occurs~$n$ times.
\begin{point}{50}{Proof}%
The case~$X = \{0\}$ is covered by~$n=0$ and~$p=0$.
Assume~$X \neq \{0\}$.
\begin{point}{60}[selfdual-normalish-form1]%
As a first step we will prove that
either
\begin{enumerate}
\item $X$ has an orthonormal basis~$E$
    such that~$\left<e,e\right>=1$ for some~$e \in E$ \emph{or}
\item there is a single vector~$e \in X$ such that~$\{e\}$ is an orthonormal basis.
\end{enumerate}
Pick any orthonormal basis~$E_0$ of~$X$.
Let~$P$ denote the poset of subsets~$U \subseteq E_0 \times \scrB$
    satisfying
\begin{enumerate}
    \item
    $u$ is a partial isometry with $uu^* = \left<e,e\right>$
        for every~$(e,u) \in U$ \emph{and}
    \item
    the domains of these~$u$ are orthogonal:  $\sum_{(e,u) \in U} u^*u \leq 1$.
\end{enumerate}
For any non-empty chain~$C \subseteq P$
    it is easy to see~$\bigcup C \in P$
    and so by Zorn's lemma,
    there is a maximal element~$U_0 \in P$.
Define
\begin{equation*}
e_0 \ =\  \sum_{(e,u) \in U_0} eu
    \qquad \qquad E_1 \ = \ E_0 \backslash \{e;\ (e,u) \in U_0\}.
\end{equation*}
The set $\{e_0\} \cup E_1$ is an orthonormal basis of~$X$
    --- indeed, it is clearly orthonormal and
    if $x = \sum_{e \in E_0} e \left<e,x\right>$,
    then
\begin{equation*}
    x \ =\  e_0 \Bigl(\sum_{(e,u) \in U_0} u^*\left<e,x\right>\Bigr)
            + \sum_{e \in E_1} e\left<e,x\right>.
\end{equation*}
If~$E_1 = \emptyset$, then~$\{e_0\}$ is an orthonormal basis of~$X$ and
    we have shown point 2.

For the other case, assume~$E_1 \neq \emptyset$.
Pick~$e_1 \in E_1$.
Write~$p_0 = \left<e_0,e_0\right>$
    and~$p_1 = \left<e_1,e_1\right>$.
Suppose~$p_1 \lesssim 1 - p_0$.
Then~$vv^* = p_1$ and~$v^*v \leq 1- p_0$ for some~$v \in \scrB$.
    Hence~$v^*v + p_0 \leq 1$.
    So~$U_0 \cup \{ (e_1, v) \} \in P$ contradicting maximality of~$U_0$.
Apparently~$p_1 \not\lesssim 1-p_0$.
So by \sref{total-mv-order}, we must have~$1-p_0 \lesssim p_1$.
Let~$v \in \scrB$ be such that~$vv^*=1-p_0$
    and~$v^*v \leq p_1$.
Write~$p_0' = 1- v^*v$.
As~$1-p_0' \leq p_1 =1-(1-p_1)$
we have~$(1-p_0') + (1-p_1) \leq 1$.
Define~$q = 1- ((1-p_0') + (1-p_1))$.
We now have the following splittings of~$1$.
\begin{center}
\begin{tikzpicture}
    \draw
            (0,0) rectangle (1.2,.5)
            (0,0) node[anchor=south west] {$1-p_1$}
            (1.2,0) rectangle (2,.5)
            (1.35,0) node[anchor=south west] {$q$} 
            (2,0) rectangle (2.8,.5)
            (2,0) node[anchor=south west] {$v^*v$} ;
    \draw [decorate,thick,decoration={brace}](2.8,-.05) -- (1.2,-.05) node[midway,anchor=north] {$p_1$};
    \draw [decorate,thick,decoration={brace}](0,.55) -- (2,.55) node[midway,anchor=south] {$p_0'$};
    \draw
            (4,0) rectangle (6,.5)
            (4.65,0) node[anchor=south west] {$p_0$} 
            (6,0) rectangle (6.8,.5)
            (6,0) node[anchor=south west] {$vv^*$} ;
\end{tikzpicture}
\end{center}
Note~$q = p_0'+p_1 - 1$ and~$p_0', p_1 \geq q$.
Define
\begin{equation*}
    D \ = \ \{d_0, d_1\} \qquad \text{where} \quad
    d_0 \ =\   e_1 q \qquad
    d_1 \ = \ e_0 + e_1v^*.
\end{equation*}
From~$p_1q=q$ it follows~$D$ is orthogonal.
Clearly~$\left<d_0,d_0\right>=q$ and
\begin{equation*}
    \left<d_1, d_1\right>
        \ =\ \left<e_0, e_0\right> + vp_1v^*    
        \ =\ p_0 + vv^* \ =\ 1.
\end{equation*}
So~$D$ is an orthonormal set.
Using~$0 = (1-p_0)p_0 = vv^*p_0 $ and~$q+v^*v = p_1$, we see
that for any~$x \equiv e_0 b_0 + e_1 b_1$
with~$b_0 \in p_0\scrB$ and~$b_1 \in p_1 \scrB$, we have
\begin{align*}
    d_0 \langle d_0, x \rangle + d_1 \langle d_1 , x \rangle
    &\ = \ d_0 b_1 + d_1(vb_1 + b_0) \\
    &\ = \ e_1 q b_1 + e_0 vb_1 + e_1v^*vb_1 + e_0 b_0 + e_1v^* b_0 \\
    &\ = \ e_1 q b_1 + e_0 p_0 vv^* vb_1 + e_1v^*vb_1 + e_0 b_0 + e_1v^* vv^* p_0 b_0 \\
    &\ = \ e_1 q b_1 + e_1v^*vb_1 + e_0 b_0  \\
    &\ = \ e_1 (q + v^*v) b_1 + e_0 b_0  \\
    &\ = \ e_1 b_1 + e_0 b_0 \ = \ x.
\end{align*}
Thus~$D$ is an orthonormal basis of~$\{e_0, e_1\}^{\perp\perp}$.
Hence~$E \equiv D \cup (E_1 - \{e_1\})$ is an orthonormal basis of~$X$
with~$d_1 \in E$ such that~$\left<d_1,d_1\right>=1$.
\end{point}
\begin{point}{70}%
For brevity, call~$X$ \Define{1-dim} if there is a one-element orthonormal basis of~$X$.
In~\sref{selfdual-normalish-form1} we saw
    how to create an orthonormal basis~$D$ of a non-1-dim~$X$
    with an~$e_0 \in D$ such that~$\left<e_0,e_0\right>=1$.
If~$\{e_0\}^\perp$ is not 1-dim,
    we can apply \sref{selfdual-normalish-form1}
    on~$\{e_0\}^\perp$ (instead of~$X$)
    to find an orthonormal basis~$D'$ of~$\{e_0\}^\perp$
    with~$e_1 \in D'$ such that~$\left<e_1,e_1\right>=1$.
This procedure can be continued using Zorn's lemma as follows.
Let~$P$ denote the poset of orthogonal subsets~$E \subseteq X$
    satisfying~$\langle e,e\rangle = 1$ for all~$e\in E$.
Clearly~$\emptyset\in P$
and for any non-empty chain~$C \subseteq P$
    we have~$\bigcup C \in P$.
By Zorn's lemma there is a maximal~$E \in P$.
Suppose~$E^\perp \neq \{0\}$ and~$E^\perp$ is not 1-dim.
Then by \sref{selfdual-normalish-form1}
    we can find an orthonormal basis~$E'$ of~$E^\perp$
    together with~$e \in E'$ such that~$\left<e,e\right>=1$.
Now~$E' \cup \{e\} \in P$, contradicting maximality of~$E$.
Apparently~$E^\perp = \{0\}$ or~$E^\perp$ is 1-dim.
If~$E^\perp = \{0\}$
    then we are done
    taking for~$\kappa$ the cardinality of~$E$.
For the other case, assume~$E^\perp = \{e\}^{\perp\perp}$
    for some~$e \in X$ with~$\left<e,e\right>$ a non-zero projection.
If~$E$ is finite, then we are done as well.
So, assume~$E$ is infinite.
Pick a sequence~$e_1, e_2, \ldots \in E$ of distinct elements.
Write~$p = \left<e,e\right>$
    and~$E_0 = \{e_n;\ n\in \N\}$.
Define
\begin{equation*}
    E_1 \ =\  \{e + e_1(1-p), \ e_1p + e_2(1-p),\  e_2p  +e_3(1-p),\  \ldots \}.
\end{equation*}
It is easy to see~$E_1$ is an orthogonal set
    with~$\left<d,d\right>=1$ for all~$d \in E_1$.
Furthermore~$E_1$ is an orthonormal basis for~$E_0 \cup \{e\}$,
    hence~$(E - E_0) \cup E_1$ is the desired orthonormal basis of~$X$. \qed
\end{point}
\end{point}
\end{point}
\begin{point}{80}%
Does this settle the issue raised
        in~\sref{thel2matter}?
      One might hope
    that~$\ell^2((1_{\mathscr{B}})_{\alpha \in \kappa}) \cong
        \ell^2((1_{\mathscr{B}})_{\beta\in\lambda})$
    implies~$\kappa=\lambda$ for all cardinals~$\kappa$,$\lambda$.
This is not the case.
Indeed, if~$\scrB$ is a type III factor
    of operators on a separable Hilbert space,
    then for any non-zero projection~$p$,
    we have~$p \sim 1 -p$ (by combining~\cite[dfn.~6.5.1]{kr}
        and~\cite[cor.~6.3.5]{kr})
    and so~$\scrB \cong p\scrB \oplus (1-p)\scrB \cong 
    \scrB \oplus \scrB$ as Hilbert C$^*$-modules.
\end{point}
\end{parsec}

\begin{parsec}{1630}%
\begin{point}{10}%
Before we continue to the next topic,
    the exterior tensor product of Hilbert C$^*$-modules,
    we will show that the ultranorm completion~\sref{dils-completion}
    is determined by its universal property.
    We will need this fact in our treatment
        of the exterior tensor product.
\end{point}
\begin{point}{20}[selfdual-compl-defining]{Proposition}%
Let~$\scrB$ be a von Neumann algebra and~$V$ a right~$\scrB$-module
    with~$\scrB$-valued inner product~$[\,\cdot\,,\,\cdot\,]$.
    \index{$[\,\cdot\,,\,\cdot\,]$!inner product!$\scrB$-valued}
There is an up-to-isomorphism unique
    self-dual Hilbert~$\scrB$-module~$X$
    together with inner product preserving $\scrB$-linear~$\eta \colon V \to X$
    with the following universal property.
    \begin{quote}
    For every bounded~$\scrB$-linear map~$T\colon V \to Y$
    to some self-dual Hilbert~$\scrB$-module~$Y$,
    there is a unique bounded~$\scrB$-linear map~$\hat{T}\colon X \to Y$
    with~$\hat{T} \after \eta = T$.
    \end{quote}
Moreover, for such~$\eta\colon V \to X$,
    the image of~$V$ under~$\eta$ is ultranorm dense.
\begin{point}{30}{Proof}%
We already know such a~$X$ and~$\eta\colon V \to X$
    exist by~\sref{dils-completion}
    and~\sref{selfdual-completion-univ}.
For this one~$\eta(V)$ is ultranorm dense in~$X$.
Assume there is another self-dual Hilbert~$\scrB$-module~$X_2$
    with inner product preserving~$\scrB$-linear~$\eta_2\colon V \to X_2$
    satisfying the universal property.
By the universal property of~$\eta$ applied to~$\eta_2$,
    there must exist a unique~$U \colon X  \to X_2$
    with~$U \after \eta  = \eta_2$.
Similarly, there is a~$V \colon X_2 \to X$
    with~$V \after \eta_2 = \eta$.
Clearly~$\id_{X} \after \eta = V \after \eta_2 = V \after U \after \eta$
    and so by the uniqueness~$\id = V\after U$.
    Similarly~$U \after V = \id$.
It remains to be shown~$U$ preserves the inner product.
As
\begin{equation*}
    \left<\eta(x),\eta(y)\right>
    \ =\  [x,y]
    \ =\  \left<\eta_2(x), \eta_2(y)\right>
    \ =\  \left<U \eta(x), U \eta(y)\right>
\end{equation*}
for all~$x,y \in X$, we know
$\eta(X)$ is ultranorm dense.
So by ultranorm continuity and bijectivity of~$U$,
    we know~$U(\eta(X)) = \eta_2(X)$ is ultranorm dense in~$X_2$.
Hence~$U$ preserves the inner product by~\sref{innerprod-ultraweak}
\qed
\end{point}
\end{point}
\end{parsec}

\subsection{Exterior tensor product}
\begin{parsec}{1640}%
\begin{point}{10}%
Suppose~$X$ is a Hilbert~$\scrA$-module
    and~$Y$ is a Hilbert~$\scrB$-module
    for some C$^*$-algebras~$\scrA$ and~$\scrB$.
The algebraic tensor product~$X \odot Y$
    is a right $\scrA \odot \scrB$-module
    via~$(x\otimes y)\cdot (a\otimes b) = (xa) \otimes (yb)$
    and has an~$\scrA \odot \scrB$-valued inner
    product fixed by~$[x\otimes y, x'\otimes y']
        = \left<x,x'\right> \otimes \left<y,y'\right>$.
Write~$\scrA \otimes_\sigma \scrB$ for the spatial
    C$^*$-tensor product~\cite[\S11.3]{kr} of~$\scrA$ and~$\scrB$.
The inner product on~$X \odot Y$ can be extended
    to an~$\scrA \otimes_\sigma \scrB$-valued inner product
    on the norm completion of~$X \odot Y$
    on which it is definite~\cite{lance}.
This is a Hilbert C$^*$-module,
    which is called the \emph{exterior tensor product} of~$X$ and~$Y$.
It is quite difficult to perform this construction:
    Lance discusses the difficulties in \cite[ch.~4]{lance}.
But if we assume~$\scrA$ and~$\scrB$ are von Neumann algebras
    \emph{and} that $X$ and~$Y$ are self-dual,
    it is rather easier.
In fact, we will directly construct the ultranorm completion
    of the exterior tensor product,
    which is a self-dual Hilbert $\scrA \otimes \scrB$-module,
    which we call the \emph{self-dual exterior tensor product}.
This self-dual exterior tensor product
    behaves much better than the plain exterior tensor product.
    As a first example, we will see (in \sref{ba-ext-tensor-pres})
    that~$\scrB^a(X) \otimes \scrB^a(Y) \cong \scrB^a(X\otimes Y)$
        for the self-dual exterior tensor product,
        which is false in general for the plain exterior tensor product.
The self-dual exterior tensor product also appears naturally
    when computing the Paschke dilation
    of a tensor product~$\varphi \otimes \psi$, see \sref{paschke-tensor}.
\end{point}
\begin{point}{20}[univprop-ext-tensor]{Theorem}%
Suppose~$\scrA$ and~$\scrB$ are von Neumann algebras.
For any self-dual Hilbert $\scrA$-module $X$
    and self-dual Hilbert~$\scrB$-module $Y$,
    there is an up-to-isomorphism unique
    self-dual Hilbert~$\scrA\otimes \scrB$-module
    $\Define{X \otimes Y}$, called the \Define{(self-dual exterior) tensor product}%
    \index{tensor product!of Hilbert $\scrB$-modules}%
    \index{*tensor@$\otimes$!Hilbert $\scrB$-modules},
    together with a~$\scrA \odot \scrB$-linear injective
    inner product preserving map~$\eta \colon X \odot Y \to X\otimes Y$
    with the following universal property.
    \begin{quote}
    For every self-dual~Hilbert~$\scrA\otimes\scrB$-module
    $Z$ with bounded~$\scrA \odot \scrB$-linear
        $T\colon X \odot Y \to Z$,
    there is a unique bounded~$\scrA \otimes \scrB$-linear
        map~$\hat{T}\colon X \otimes Y \to Z$
    with~$\hat{T} \after \eta= T$,
    \end{quote}
where we norm~$X \odot Y$ with~$\|t \| = \|[t,t]\|^{\frac{1}{2}}_{\scrA\otimes \scrB}$.
Write~$x \hotimes y = \eta(x\otimes y)$.
For any such~$X \otimes Y$ we have the following.
\begin{enumerate}
    \item The image of $X \odot Y$ under~$\eta$ is ultranorm dense in~$X \otimes Y$.
    \item If~$(e_i)_{i \in I}$ is an orthonormal basis of~$X$
                and~$(d_j)_{j \in J}$ of~$Y$,
                then
        \begin{enumerate}
            \item $(e_i \hotimes d_j)_{i,j\in I\times J}$
                is an orthonormal basis of~$X \otimes Y$ \emph{and}
            \item
                the linear span of
                \begin{equation*}
                    \{\, \ketbra{\,(e_i a) \hotimes (d_j b) \,}{\,
                        e_{k}\otimes d_{l} \,}; \ 
                        a \in \scrA, \ 
                        b \in \scrB,\ 
                        i,k \in I,\ 
                        j,l \in J \,
                    \} 
                \end{equation*}
                is ultraweakly dense in~$\scrB^a(X \otimes Y)$.
        \end{enumerate}
                
\end{enumerate}
\spacingfix{}
\begin{point}{30}{Proof}%
Pick orthonormal bases~$(e_i)_{i \in I}$ and~$(d_j)_{j \in J}$ of~$X$
    and~$Y$ respectively.
Write~$p_{ij} \equiv \left<e_i, e_i\right>\otimes \left<d_j,d_j\right>
    \in \scrA \otimes \scrB$.
Define~$X \otimes Y \equiv \ell^2((p_{ij})_{i,j \in I\times J})$,
    with~$\ell^2$ as in~\sref{hilbmod-el2}.
This definition of~$X \otimes Y$ seems to depend
    on the choice of bases~$(e_i)_i$ and~$(d_j)_h$.
    We show it does not (up-to-isomorphism), in \sref{ext-tensor-uniqueness}.
\begin{point}{40}[ext-tensor-dfn-eta]{Definition~$\eta$}%
Pick~$x \in X$ and~$y \in Y$.
Write~$x_i \equiv \left<e_i, x\right>$
and~$y_j \equiv \left<d_j, y\right>$.
Then
\begin{align*}
    \left<x,x\right> \otimes \left<y,y\right>
    &\ = \ \bigl( \sum_i x_i^*x_i\bigr) \otimes\bigl(\sum_jy_j^*y_j\bigr)
            &\quad&\text{by Parseval}  \\
            &\ = \ \sum_{i,j}  x_i^*x_i  \otimes y_j^*y_j \\
            &\ = \ \sum_{i,j}  (x_i \otimes y_j)^* (x_i  \otimes y_j).
\end{align*}
Hence~$(x_i \otimes y_j)_{i,j}$ is~$\ell^2$-summable,
    and so there is a unique~$\scrA \odot\scrB$-linear
    map~$\eta\colon X \odot Y \to X \otimes Y$ fixed
    by~$\eta(x \otimes y) = (x_i \otimes y_j)_{i,j}$.
\end{point}
\begin{point}{50}[ext-tensor-preserves-inner-prod]{$\eta$ preserves inner product}%
Pick additional~$x' \in X$ and~$y' \in Y$.
Writing~$x'_i = \left<e_i,x'\right>$
    and~$y'_j = \left<d_j, y'\right>$,
we have
\begin{equation*}
    \left<\eta(x\otimes y), \eta(x' \otimes y')\right>
    \ = \ \sum_{i,j} (x_i\otimes y_j)^* \otimes (x_i' \otimes y_j')
    \ = \ \sum_{i,j} x_i^*x_i' \otimes y_j^* y_j'
\end{equation*}
and
\begin{equation*}
    \bigl(\sum_i x_i^*x_i'\bigr) \otimes \bigl( \sum_j y_j^* y_j' \bigr)
    \ = \ \left<x,x'\right> \otimes \left<y,y'\right>
    \ = \ \left<x\otimes y,x' \otimes y'\right> .
\end{equation*}
Thus it is sufficient to show
    $\sum_{i,j} x_i^*x_i' \otimes y_j^* y_j'
    = \bigl(\sum_i x_i^*x_i'\bigr) \otimes \bigl( \sum_j y_j^* y_j' \bigr)$,
    which holds as the product np-functionals
    are separating by definition, see \sref{tensor}.
Thus indeed~$\eta$ preserves the inner product.
\end{point}
\begin{point}{60}{$\eta$ is injective}%
Let~$n \in \N$,~$x_1, \ldots, x_n \in X$
    and~$y_1,\ldots, y_n \in Y$
    be given such that~$\eta(\sum_l x_l \otimes y_l) = 0$
By~\sref{selfdual-gramschmidt}
    there are orthonormal~$e'_1,\ldots, e'_m \in X$
    such that~$x_l = \sum_i e'_i \langle e'_i, x_l\rangle$.
Similarly~$y_l = \sum_j d'_j \langle d'_j, y_l\rangle$
for some orthonormal~$d'_1, \ldots, d'_{m'} \in Y$.
As~$\eta$ preserves the inner product,
    the elements~$e_i'\hotimes d_j'$ are again orthonormal
    and so from
\begin{equation*}
    0 \ =\  \sum_l x_l \hotimes y_l
    \ =\  \sum_{i,j} (e'_i \hotimes d'_j) \bigl(\sum_l
         \langle e'_i, x_l\rangle \otimes \langle d'_j, y_l\rangle\bigr)
\end{equation*}
         it follows
         $\sum_l \langle e_i',x_l\rangle \otimes \langle d_j',y_l\rangle = 0$
         for all~$i,j$.
Consequently
\begin{equation*}
    \sum_l x_l \otimes y_l
        \ =\  \sum_{i,j} (e'_i \otimes d'_j) \bigl(\sum_l
        \langle e'_i, x_l \rangle \otimes \langle d'_j, y_l \rangle 
        \bigr) \ = \ 0,
\end{equation*}
which shows~$\eta$ is injective
    (and that the inner product on~$X \odot Y$ is definite).
\end{point}
\begin{point}{70}[ultranorm-dense-tensor-base]{Image $\eta$ ultranorm dense}%
For brevity, write~$E \equiv \{ e_i \hotimes d_j; \ i,j \in I\times J\}$.
By~\sref{tensor} (and \sref{ultraclosed})
    $\scrA \odot \scrB$ is ultrastrongly dense in~$\scrA \otimes \scrB$.
So~$E (\scrA \odot \scrB)$
    is ultranorm dense in~$E (\scrA \otimes \scrB)$,
    see~\sref{ultranormcontstruct}.
In turn the linear span of~$E (\scrA \otimes \scrB)$
    is (by construction) ultranorm dense in~$X \otimes Y$.
Thus the linear span of~$E (\scrA \odot \scrB)$
    is ultranorm dense in~$X \otimes Y$.
Hence the image of~$\eta$,
    which contains the linear span of~$E (\scrA \odot \scrB)$,
    is ultranorm dense in~$X \otimes Y$.
\end{point}
\begin{point}{80}{Universal property}%
Let~$Z$ be a self-dual~Hilbert $\scrA \otimes \scrB$-module
 with some bounded~$\scrA \odot \scrB$-linear $T\colon X \odot Y \to Z$.
If~$X \odot Y$ were an~$\scrA \otimes \scrB$-module
    \emph{and} both~$\eta$ and $T$ were~$\scrA \otimes \scrB$-linear,
    we could simply apply~\sref{selfdual-completion-univ}
    to get the desired~$\hat{T}$.
Instead, we will retrace the steps of its proof
    and modify it for the present situation.
Uniqueness of~$\hat{T}$ is the same:
    $\hat{T}$ is fixed by the ultranorm-dense image of~$\eta$
    as it is bounded~$\scrA\otimes\scrB$-linear
    and so ultranorm continuous by \sref{blinear-bounded-is-ultranorm}.
The argument for existence of~$\hat{T}$ is more subtle.
We will show
\begin{equation}\label{eq:exttensorcontT}
    \left<T t, T t\right> \leq \|T\|^2 [t,t]
        \qquad \text{for all~$t \in X \odot Y$}
\end{equation}
as a replacement for the ultranorm continuity of~$T$.
We cannot apply~\sref{blinear-inprod-inequality} directly:
    as~$\scrA \odot \scrB$ is not a C$^*$-algebra,
    we cannot construct~$(\left<t,t\right>+\varepsilon)^{-\frac{1}{2}}$
    required for its proof.
To prove \eqref{eq:exttensorcontT}, we will work in the norm completions.
Write~$\overline{\scrA \odot \scrB}$
    for the norm-closure of~$\scrA \odot \scrB$ in~$\scrA \otimes \scrB$.
Clearly~$\overline{\scrA \odot \scrB}$
    is a~C$^*$-algebra
    (in fact it is the spatial C$^*$-tensor product of~$\scrA$ and~$\scrB$).
Write~$\overline{X \odot Y}$
    for the norm-closure of the image of~$\eta$ in~$X \otimes Y$.
By norm continuity
    (see \sref{module-seminorm})
    the operations of~$X \otimes Y$
    restrict to turn~$\overline{X \odot Y}$ into a
    Hilbert~$\overline{\scrA \odot \scrB}$-module.
As~$T$ is uniformly continuous,
    there is a unique bounded~$\overline{\scrA\odot\scrB}$-linear
    $\overline{T}\colon \overline{X\odot Y} \to Z$
    with~$\overline{T} \after \eta = T$ and~$\|T\| = \|\overline{T}\|$.
Now we can apply the proof of~\sref{blinear-inprod-inequality}
    to~$\overline{T}$
    to find~$\langle \overline{T}s,\overline{T}s\rangle \leq
    \|\overline{T}\|^2\left<s,s\right>$
    for~$s \in \overline{X \odot Y}$.
Substituting~$s = \eta(t)$,
    we get \eqref{eq:exttensorcontT}.

To define~$\hat{T}$, pick any~$t \in X \otimes Y$.
As the image of~$\eta$ is ultranorm dense,
there is a net~$(\eta(t_\alpha))_\alpha$ with~$\eta(t_\alpha) \to t$ ultranorm.
For any~np-map~$f\colon \scrA\otimes \scrB \to \C$
\begin{align*}
    \|T(t_\alpha) - T(t_\beta) \|_f^2 
        &\ =\ f(\langle T(t_\alpha - t_\beta), T(t_\alpha - t_\beta) \rangle) \\
        & \ \leq \ \|T\|^2 f([t_\alpha-t_\beta, t_\alpha - t_\beta]) \\
        & \ = \ \|T\|^2 f(\left<\eta (t_\alpha)-\eta(t_\beta), \eta(t_\alpha) - \eta(t_\beta)\right>) \ \to\  0.
\end{align*}
So~$(T t_\alpha)_\alpha$ is ultranorm Cauchy in~$Z$.
With a similar argument we see~$(Tt_\alpha)_\alpha $
    is equivalent to~$ (T t'_\alpha)_\alpha$
whenever~$(\eta t_\alpha)_\alpha $ is equivalent
    to~$ (\eta t'_\alpha)_\alpha$.
So we may define
\begin{equation*}
\hat{T} t \ = \ \unlim_\alpha T t_\alpha.
\end{equation*}
Clearly~$\hat{T} \after \eta = T$.
With the same argument as for \sref{selfdual-completion-univ},
    we see~$\hat{T}$ is bounded and additive.
To show~$\hat{T}$ is~$\scrA\otimes \scrB$-linear,
    pick any~$b \in \scrA \otimes \scrB$.
There is a net~$b_\beta$ in~$\scrA \odot \scrB$
    with~$b_\beta \to b$ ultrastrongly.
If $\hat{T}$ is ultranorm continuous, we are done:
\begin{equation*}
    (\hat{T} t) b
    = \unlim_\beta (\unlim_\alpha T t_\alpha) b_\beta\\
     = \unlim_\beta \unlim_\alpha T (t_\alpha b_\beta) \\
     = \unlim_\beta \hat{T} (t b_\beta)\\
     = \hat{T} (tb).
\end{equation*}
To show~$\hat{T}$ is ultranorm continuous,
    assume~$t^\alpha \to 0$ ultranorm in~$X \otimes Y$.
It is sufficient to show~$\hat{T} t^\alpha \to 0$ ultranorm.
Let~$f\colon \scrA\otimes \scrB \to \C$ be any np-map.
Write~$\Phi$ for the set of entourages for ultranorm uniformity
    on~$X \otimes Y$.
As the image of~$\eta$ is ultranorm dense,
there is a~$\Phi$-indexed Cauchy net~$(t^\alpha_E)_{E \in \Phi}$ in
    $\eta (X \odot Y)$
    with~$t^\alpha_E \mathrel{E} t^\alpha$
    for each~$E \in \Phi$.
Thus~$\hat{T} t^\alpha = \unlim_E T t^\alpha_E$.
As~$t^\alpha \to 0$
    there is some~$\alpha_0$
    such that for all~$\alpha \geq \alpha_0$
    we have~$\|t^\alpha \|_f \leq \varepsilon$.
We can find~$E_0 \geq E_{f,\varepsilon}$
such that for all~$E \geq E_0$
we have~$\| \hat{T} t^{\alpha_0} -T t^{\alpha_0}_E \|_f \leq \varepsilon$.
For~$E \geq E_0$ we have~$E \geq E_{f,\varepsilon}$ and so
\begin{equation*}
 \|T t^{\alpha_0}_{E_0}\|_f
 \  \overset{\eqref{eq:exttensorcontT}}{\leq}\  \|T\| \|t^{\alpha_0}_{E_0}\|_f
  \  \leq\ \|T\|
  ( \|t^{\alpha_0}_{E_0} - t^{\alpha_0}\|_f +
  \| t^{\alpha_0} \|_f)
   \ \leq\  2 \|T\| \varepsilon.
\end{equation*}
    Consequently
\begin{equation*}
    \|\hat{T} t^{\alpha_0} \|_f \ \leq  \ 
    \|\hat{T} t^{\alpha_0} - T t^{\alpha_0}_{E_0} \|_f
                + \|T t^{\alpha_0}_{E_0}\|_f
    \ \leq\  \varepsilon + 2\|T\|\varepsilon.
\end{equation*}
So~$\hat{T}t^\alpha \to 0$ ultranorm, as desired.
\end{point}
\begin{point}{90}[ext-tensor-uniqueness]{Uniqueness}%
Assume there is a self-dual Hilbert~$\scrA\otimes\scrB$-module~$X \otimes_2 Y$
    together with a injective
    inner product-preserving~$\scrA\odot\scrB$-linear
    map~$\eta_2 \colon X \odot Y \to X \otimes_2 Y$
    also satisfying the universal property.
With the same reasoning as in \sref{selfdual-compl-defining}
    there is an invertible bounded~$\scrA\otimes\scrB$-linear
    $U\colon X \otimes Y \to X \otimes_2 Y$
    with~$U \after \eta = \eta_2$.
Clearly~$\left<Ux, U y\right> = \left<x,y\right>$
    for~$x,y \in E (\scrA\odot\scrB)$,
    where we defined~$E = \{e_i \hotimes d_j; \ i,j \in I\times J\}$.
As we saw the linear span of~$E (\scrA\odot\scrB)$
    is ultranorm dense in~$X \otimes Y$,
    we see~$U$ must preserve the inner product by \sref{innerprod-ultraweak}.
\end{point}
\begin{point}{100}%
Property 1 from the statement of the Theorem follows immediately
    from the fact that the exterior tensor product is unique
    up to an isomorphism which respects the embeddings.
Assume~$(e'_i)_{i \in I'}$ is an orthonormal basis of~$X$
    and~$(d'_j)_{j \in J'}$ is an orthonormal basis of~$Y$.
We will show $E_2 =\{ e'_i \hotimes d'_j; \ i,j \in I'\times J'\}$
is an orthonormal basis of~$X \hotimes Y$.
Clearly~$E_2$ is orthonormal.
By \sref{selfdual-orthn-basis} and \sref{hilbmod-projthm}
    it is sufficient to show~$e_{i_0} \hotimes d_{j_0} \in E_2^{\perp\perp}$
    for all~$i_0 \in I$ and $j_0 \in J$.
Write~$a_i = \left<e_i',e_{i_0}\right>$
    and~$b_j = \left<d_j',d_{j_0}\right>$.
By \sref{selfdual-orthn-basis} and \sref{ext-tensor-preserves-inner-prod}
    it suffices to show the latter equality in
    \begin{equation*}
        \bigl(\sum_i a_i^*a_i\bigr) \otimes \bigl(\sum_j b_j^*b_j\bigr) \ = \ 
    \langle e_{i_0}\hotimes d_{j_0}, e_{i_0}\hotimes d_{j_0}\rangle
    \ =\  \sum_{i,j} a_i^*a_i \otimes b_j^*b_j,
    \end{equation*}
    which holds as product np-functionals are separating
    by definition, see~\sref{tensor}.
\end{point}
\begin{point}{110}%
    Finally we prove that the linear span of
                \begin{equation*}
                    D \ \equiv \ \{ \ketbra{(e'_i a) \hotimes (d'_j b)}{
                        e_{k}\otimes d_{l} }; \ 
                        a \in \scrA, \ 
                        b \in \scrB,\ 
                        i,k \in I',\ 
                        j,l \in J'
                    \} 
                \end{equation*}
            is ultraweakly dense in~$\scrB^a(X \otimes Y)$.
By \sref{ketbra-ultraweakly-dense},
    the linear span
    of (the set of operators of
    the form)~$\ketbra{(e'_i \hotimes d'_j)t}{e'_k \hotimes d'_l}$
    with~$t \in \scrA \otimes \scrB$
    is ultraweakly dense in~$\scrB^a(X \otimes Y)$.
We will show that~$t \in \scrA \odot \scrB$ suffices.
    By~\sref{tensor}, \sref{ultraclosed} and~\sref{kaplansky},
    each~$t \in \scrA \otimes \scrB$
    is the ultrastrong limit
    of a norm bounded net~$t_\alpha$ in~$\scrA \odot \scrB$.
    So by \sref{ketbra-ultranorm-continuous} (and \sref{ultranormcontstruct}),
we see~$\ketbra{(e'_i \hotimes d'_j)t_\alpha}{e'_k \hotimes d'_l}$
converges ultraweakly to~$\ketbra{(e'_i \hotimes d'_j)t}{e'_k \hotimes d'_l}$.
So indeed, the linear span of~$D$
    is ultraweakly dense in~$\scrB^a(X \otimes Y)$. \qed
\end{point}
\end{point}
\end{point}

\begin{point}{120}{Examples}%
We give a few examples of the self-dual exterior tensor product
    and relate it to other notions of tensor product.
To avoid confusion, we will denote the self-dual exterior
    tensor product by~$\otimes_{\text{ext}}$
    (instead instead of the plain~$\otimes$.)
\begin{enumerate}
\item
Any Hilbert space is a self-dual Hilbert~$\C$ module.
For Hilbert spaces~$\scrH$ and~$\scrK$
    we have~$\scrH \otimes_{\text{Hilb}} \scrK
        \cong \scrH \otimes_{\text{ext}} \scrK$,
    where~$\otimes_{\text{Hilb}}$ denotes the regular tensor product of
    Hilbert spaces.
\item
Every von Neumann algebra is a self-dual Hilbert C$^*$-module
    over itself.
For any von Neumann algebras~$\scrA$ and~$\scrB$
    we have~$\scrA\otimes_{\text{vN}} \scrB
        \cong \scrA \otimes_{\text{ext}} \scrB$
        as Hilbert~$\scrA\otimes_{\text{vN}} \scrB$-modules,
        where~$\scrA\otimes_{\text{vN}}\scrB$ denotes the (spatial)
        tensor product of von Neumann algebras.
\item
Let~$X$ be a self-dual Hilbert~$\scrA$-module
    and~$Y$ be a self-dual Hilbert~$\scrB$-module
        for von Neumann algebras~$\scrA,\scrB$.
The norm closure of~$X\odot Y$
    in~$X \otimes_{\text{ext}} Y$
        is the plain exterior tensor product of~$X$ and~$Y$.
In turn, the ultranorm completion of the exterior tensor product
    of~$X$ and~$Y$ is isomorphic to~$X \otimes_{\text{ext}} Y$,
        the self-dual exterior tensor product.
\item
In~\sref{paschke-tensor} it will be shown that
\begin{equation*}
(\scrA_1 \otimes_{\varphi_1} \scrB_1)
    \otimes_{\text{ext}} (\scrA_2 \otimes_{\varphi_2} \scrB_2)
\ \cong \ 
(\scrA_1 \otimes \scrA_2) \otimes_{\varphi_1\otimes \varphi_2} (\scrB_1 \otimes \scrB_2)
\end{equation*}
        for all~ncp-maps $\varphi_i\colon \scrA_i \to \scrB_i$ ($i=1,2$)
    between von Neumann algebras.
\end{enumerate}
\spacingfix{}
\end{point}
\end{parsec}
\begin{parsec}{1650}%
\begin{point}{10}%
For Hilbert spaces~$\scrH$ and~$\scrK$,
    we have~$\scrB(\scrH \otimes \scrK)
        \cong \scrB(\scrH) \otimes \scrB(\scrK)$.
However, for Hilbert modules~$X$ and~$Y$
    and the regular (not necessarily self-dual)
        exterior tensor product
        it is not true in general
        that~$\scrB^a(X) \otimes \scrB^a(Y) \cong \scrB^a(X \otimes Y)$.
However, if~$X$ and~$Y$ are self-dual Hilbert modules over von Neumann
    algebras, then we do have the familiar formula
    for the self-dual exterior tensor product.
We need some preparation.
\end{point}
\begin{point}{20}{Setting}%
Let~$X$ be a self-dual Hilbert~$\scrA$-module
    and~$Y$ be a self-dual Hilbert~$\scrB$-module
    for von Neumann algebras~$\scrA$ and~$\scrB$.
\end{point}
\begin{point}{30}[dfn-tensor-of-hilbmod-maps]{Proposition}%
For any operators~$S \in \scrB^a(X)$ and~$T \in \scrB^a(Y)$,
there is a unique operator~$\Define{S \otimes T} \in \scrB^a(X\otimes Y)$%
    \index{*tensor@$\otimes$!module maps}%
    \index{tensor product!of module maps}
    fixed by~$(S \otimes T) (x \otimes y) = (S x) \otimes (T y)$
    for any~$x\in X$ and~$y \in Y$.
\begin{point}{40}{Proof}%
We will define~$S \otimes T$ using the universal
    property proven in \sref{univprop-ext-tensor}.
Write~$\Theta \colon X \odot Y \to X \otimes Y$
    for the map~$\Theta(x \otimes y) = (Sx) \otimes (Ty)$.
It is easy to see~$\Theta$ is~$\scrA \odot \scrB$-linear.
The hard part is to show that~$\Theta$ is bounded.
Pick any~$t \in X \odot Y$ --- say~$t \equiv \sum^n_{i=1} x_i \otimes y_i$.
To start, note~$(\left<Sx_i, Sx_j\right>)_{ij}$ is positive in $M_n \scrA$ ---
    indeed, for any~$a_1, \ldots, a_n \in \scrA$
    we have~$\sum_{i,j} a_i^* \left<S x_i, S x_j\right> a_j
            = \left<Sx,Sx\right> \geq 0$
                where~$x \equiv \sum_i x_i a_i$
                and so it is positive
                by \sref{when-a-matrix-over-a-cstar-algebra-is-positive}.
By reasoning in the same way for~$\sqrt{\|S\|^2 - S^*S}$
instead of~$S$,
it follows~$(\left<Sx_i, Sy_i \right>)_{ij} \leq (\|S\|^2 \left<x_i,x_j\right>)_{ij}$.
Similarly~$0 \leq (\left<Ty_i, Ty_j\right>)_{ij} \leq (\|T\|^2 \left<y_i, y_j\right>)_{ij}$ in~$M_n \scrB$.
Write~$s\colon (\scrA \otimes \scrB)^{\oplus n} \to \scrA \otimes \scrB$
    for the~$\scrA \otimes \scrB$-linear map
    given by~$s(t_1, \ldots, t_n) = t_1 + \cdots + t_n$.
    ($s$ is the \emph{row vector} filled with~$1$s).
It is easy to see~$s A s^* 1= \sum_{i,j} A_{ij}$
    for any matrix~$A$ over~$\scrA \otimes \scrB$.
As a final ingredient,
    recall that the
    bilinear map~$M_n(\otimes) \colon M_n \scrA \times M_n \scrB \to
                M_n( \scrA \otimes \scrB)$
        defined in the obvious way
        is positive by \sref{cp-bilinear} and \sref{tensor}.
    Putting everything together:
\begin{alignat*}{2}
    \langle \Theta t, \Theta t \rangle
        & \ = \ \sum_{i,j} \langle 
    (Sx_i) \otimes (Ty_i),
    (Sx_j) \otimes (Ty_j)
            \rangle  \\
        & \ = \ \sum_{i,j} 
    \langle Sx_i, S x_j \rangle
    \otimes \langle Ty_i, T y_j \rangle \\
    & \ = \ 
    s \,
    \left( 
    \langle Sx_i, S x_j \rangle
    \otimes \langle Ty_i, T y_j \rangle
\right)_{ij}  \,
     s^*1 \\
    & \ = \ 
     s \,
M_n(\otimes) \bigl(\,
    \left( 
    \langle Sx_i, S x_j \rangle \right)_{ij},
\left( \langle Ty_i, T y_j \rangle \right)_{ij} \,\bigr) 
\, s^*1 \\
    & \ \leq \ 
s \,
M_n(\otimes) \bigl(\,
    ( 
    \|S\|^2 \langle x_i,  x_j \rangle )_{ij},
(\|T\|^2 \langle y_i,  y_j \rangle )_{ij} \,\bigr) 
\, s^*1 \\
        & \ = \
    \|S\|^2\|T\|^2 \left<t,t\right>.
\end{alignat*}
We have proven that~$\Theta$ is bounded (by~$\|S\|\|T\|$).
Thus by \sref{univprop-ext-tensor},
    there exists a unique map~$S \otimes T \colon X\otimes Y \to X\otimes Y$
    with~$(S \otimes T) (x \otimes y) = \Theta(x \otimes y)
        \equiv (S x) \otimes (Ty)$,
        which is what we desired to prove. \qed
\end{point}
\end{point}
\begin{point}{50}[hilbmod-tensor-ketbra]{Exercise}%
Prove that for the
operation~$\otimes$ defined in \sref{dfn-tensor-of-hilbmod-maps}, we have
\begin{enumerate}
    \item $\ketbra{x_1}{x_2} \otimes \ketbra{y_1}{y_2}
            = \ketbra{x_1\otimes y_1}{x_2\otimes y_2}$,
        \item $1 \otimes 1 = 1$;
    \item $(S \otimes T) (S' \otimes T') = (SS' ) \otimes (TT')$
                \emph{and}
            \item $(S \otimes T)^* = S^* \otimes T^*$.

    (Hint:
    use that by the reasoning of \sref{ultranorm-dense-tensor-base}
    the linear
    span of
    \begin{equation*}
    \{ (e a) \otimes (d b); \ (a,b,e,d) \in \scrA \times \scrB
                    \times E \times D\}
    \end{equation*}
        is ultranorm dense in~$X \otimes Y$
        for any bases~$E$ of~$X$ and~$D$ of~$Y$.)

\end{enumerate}
\spacingfix{}
\end{point}
\begin{point}{60}[ba-ext-tensor-pres]{Theorem}%
There is an
    nmiu-isomorphism~$\vartheta\colon \scrB^a(X) \otimes \scrB^a(Y) \cong \scrB^a(X \otimes Y)$
    fixed by~$\vartheta(S\otimes T) = S \otimes T$.
\begin{point}{70}{Proof}%
We will show that
    the bilinear map~$\Theta \colon 
            \scrB^a(X) \times \scrB^a(Y) \to
    \scrB^a(X \otimes Y)$
    given by~$\Theta(S,T) = S \otimes T$
    is a tensor product in the sense of \sref{tensor}
    using~\sref{tensor-characterization}.
From this it follows
    by~\sref{tensor-uniqueness},
    that there is a
    nmiu-isomorphism
    \begin{equation*}
    \vartheta\colon \scrB^a(X) \otimes \scrB^a(Y) \ \cong\  \scrB^a(X \otimes Y)
    \end{equation*}
    with~$\vartheta(S\otimes T) = \Theta(S , T) = S \otimes T$, as promised.
\begin{point}{80}%
From \sref{hilbmod-tensor-ketbra}, it follows~$\Theta$ is an miu-bilinear map.
By~\sref{univprop-ext-tensor}
    the linear span of operators of the form~$
        \ketbra{e_i a}{e_k} \otimes \ketbra{d_j b }{d_l}
        = \ketbra{(e_i a) \otimes (d_j b)}{e_k \otimes d_l} $
        is ultraweakly dense in~$\scrB^a(X \otimes Y)$
        for an orthonormal basis~$e_i \otimes d_j$ of~$X \otimes Y$.
        So the image of~$\Theta$ generates~$\scrB^a(X\otimes Y)$.
\begin{point}{90}%
    It remains to be shown that there are sufficiently many product functionals
    (see~\sref{tensor} for the definition of product functional
    and their sufficiency).
    Write
    \begin{align*}
        \Omega_X &\ \equiv \ 
                    \{ f(\langle x, (\,\cdot\,) x \rangle);
                    \ f \colon \scrA \to \C \text{ np-map}, x\in X \} \\
        \Omega_Y& \  \equiv\ 
                    \{ g(\langle y, (\,\cdot\,) y \rangle);
                    \ g \colon \scrB \to \C \text{ np-map}, y\in Y \}.
    \end{align*}
From~\sref{hilbmod-ordersep}
    it follows~$\Omega_X$  and~$\Omega_Y$ are order separating.
For~$\sigma \in \Omega_X$ and~$\tau \in \Omega_Y$
 --- say~$\sigma \equiv f(\langle x, (\,\cdot\,) x\rangle)$,
    $\tau \equiv g(\langle y, (\,\cdot\,) y\rangle)$ ---
    define~$\sigma \otimes \tau \equiv
        (f \otimes g) (x \otimes y, (\,\cdot\,) x\otimes y) $
        and~$\Omega \equiv \{\sigma \otimes \tau; \sigma \in \Omega_X,
                    \ \tau \in \Omega_Y\}$.
We will show
    that~$\Omega$ is faithful
    and that~$(\sigma \otimes \tau)(\Theta(S,T)) = \sigma(S) \tau(T)$
    for all~$S \in \scrB^a(X)$ and~$T \in \scrB^a(Y)$,
    which is sufficient to show~$\Theta$ is a tensor product.
Indeed for any~$S \in \scrB^a(X)$ and $T \in \scrB^a(Y)$, we have
\begin{align*}
    (\sigma \otimes \tau)(\Theta(S, T))
        &\ = \ (f \otimes g) (\, \langle x \otimes y, 
                (S \otimes T) (x \otimes y) \rangle\,)\\
        &\ = \ (f \otimes g) (\,\langle x \otimes y, 
                (S x) \otimes (T y)\rangle\,)\\
        &\ = \ (f \otimes g) (\,
            \langle x, Sx\rangle \otimes
            \langle y, Ty\rangle \rangle\,)\\
        &\ = \ 
            f(\langle x, Sx\rangle ) \cdot
            g(\langle y, Ty\rangle \rangle)\\
        &\ = \ 
            \sigma(S) \tau(T).
\end{align*}
\end{point}
\spacingfix{}
\begin{point}{100}%
Finally, to show that~$\Omega$ is faithful,
        assume~$A \in \scrB^a(X \otimes Y)$, $A \geq 0$
        and $(\sigma \otimes \tau)(A) = 0$
        for all~$\sigma \otimes \tau \in \Omega$.
Then for all~$x \in X$ and~$y \in Y$,
    we have~$\langle x\otimes y, A (x\otimes y)\rangle = 0$,
        thus~$\| \sqrt{A} (x \otimes y) \|^2 = 0$,
        hence~$\sqrt{A} (x \otimes y) = 0$.
So~$\sqrt{A}$ vanishes on~$X \odot Y$,
    which is ultranorm dense in~$X \otimes Y$,
    whence~$\sqrt{A} = 0$. So~$A = 0$.
\qed
\end{point}
\end{point}
\end{point}
\end{point}
\end{parsec}

\begin{parsec}{1660}%
\begin{point}{10}%
To compute the tensor product of Paschke dilations,
    we will need to know a bit more about
    the ultranorm continuity of the exterior tensor product.
\end{point}
\begin{point}{20}[ultranorm-continuity-ext-tensor]{Lemma}%
Let~$X$ be a self-dual Hilbert~$\scrA$-module
    and~$Y$ be a self-dual Hilbert~$\scrB$-module
    for von Neumann algebras~$\scrA$ and~$\scrB$.
If~$x_\alpha \to x$ and~$y_\alpha \to y$ ultranorm
    for norm-bounded nets~$x_\alpha$ in~$X$ and~$y_\alpha$ in $Y$,
    then~$x_\alpha \otimes y_\alpha \to x \otimes y$ ultranorm.
\begin{point}{30}{Proof}%
As~$x_\alpha \otimes y_\alpha - x \otimes y
            = (x_\alpha -x )\otimes y_\alpha + x\otimes (y_\alpha - y)$,
            it is sufficient
    to show that both~$(x_\alpha -x) \otimes y_\alpha \to 0$
    and~$x\otimes (y_\alpha - y) \to 0$ ultranorm.
Clearly
\begin{align*}
    \langle (x_\alpha - x) \otimes y_\alpha,
        (x_\alpha - x) \otimes y_\alpha \rangle
    & \ = \ 
    \langle x_\alpha - x , x_\alpha - x \rangle\otimes 
    \langle y_\alpha ,y_\alpha \rangle \\
    & \ = \ 
    (\langle x_\alpha - x , x_\alpha - x \rangle\otimes 1) \cdot (1 \otimes 
    \langle y_\alpha ,y_\alpha \rangle ).
\end{align*}
The maps~$b \mapsto 1 \otimes b$ and~$a \mapsto a \otimes 1$
    are ncp % by \TODO{},
    so both~$1 \otimes \langle y_\alpha, y_\alpha \rangle$
    and~$\langle x_\alpha - x , x_\alpha - x \rangle\otimes 1$
    are bounded nets of positive elements.
    Also~$\langle x_\alpha - x , x_\alpha - x \rangle\otimes 1 \to 0$
        as~$x_\alpha \to x$ ultranorm.
Thus \sref{vanishing-effects} applies and we see
    $(\langle x_\alpha - x , x_\alpha - x \rangle\otimes 1) \cdot (1 \otimes 
    \langle y_\alpha ,y_\alpha \rangle ) \to 0$ ultraweakly,
    hence~$(x_\alpha - x) \otimes y_\alpha \to 0$ ultranorm.
With a similar argument
    one shows~$x\otimes (y_\alpha - y) \to 0$ ultranorm. \qed
\end{point}
\end{point}
\begin{point}{40}[exttensor-dense-subsets]{Lemma}%
Let~$X$ be a self-dual Hilbert~$\scrA$-module
    and~$Y$ be a self-dual Hilbert~$\scrB$-module
    for von Neumann algebras~$\scrA$ and~$\scrB$.
Suppose~$U \subseteq X$ and~$V \subseteq Y$
    are ultranorm-dense submodules.
Then the linear span of
\begin{equation*}
    U \otimes V \ \equiv\  \{ u \otimes v;\  u \in U,\ v \in V\}
\end{equation*}
is ultranorm dense in~$X \otimes Y$.
\begin{point}{50}{Proof}%
It is sufficient to show that every~$x \otimes y \in X \otimes Y$
    is the ultranorm limit of elements from~$U \otimes V$
    as the linear span of elements of the form~$x \otimes y$
    is ultranorm dense in~$X \otimes Y$.
So pick any element of the form~$x \otimes y \in X \otimes Y$.
By the Kaplansky density theorem for Hilbert C$^*$-modules,
 see \sref{kaplansky-hilbmod},
    there are norm-bounded nets~$x_\alpha$ in~$U$
    and~$y_\alpha$ in~$V$
    with~$x_\alpha \to x$ and~$y_\alpha \to y$ ultranorm.
So by \sref{ultranorm-continuity-ext-tensor},
    we see~$x_\alpha \otimes y_\alpha \to x \otimes y$. \qed
\end{point}
\end{point}
\begin{point}{60}[dilationspace-dense-subset]{Lemma}%
Let~$\varphi \colon \scrA \to \scrB$ be an ncp-map between von Neumann algebras
    together with ultrastrongly
    dense~$*$-subalgebras~$\scrA' \subseteq \scrA$
        and~$\scrB' \subseteq \scrB$.
Then: the linear span
    of~$ \scrT  \equiv \{ a \otimes b;\ a \in \scrA', \ b\in \scrB'\}$
is ultranorm dense in~$\scrA \otimes_\varphi \scrB$,
which was defined in  \sref{existence-paschke}.
\begin{point}{70}{Proof}%
It is sufficient to show that any~$a\otimes b \in \scrA \otimes_\varphi \scrB$
    is the ultranorm limit of elements from~$\scrT$.
    Pick any~$a\otimes b\in \scrA \otimes_\varphi \scrB$.
By \sref{dense-subalgebra}
    there are norm-bounded nets~$a_\alpha$ in~$\scrA'$ and
    $b_\alpha$ in~$\scrB'$
    with~$a_\alpha \to a$ and~$b_\alpha \to b$ ultrastrongly.
Say~$\|b_\alpha \|, \|a_\alpha\| \leq B$ for some~$B > 0$.
We will show~$a_\alpha \otimes b_\alpha - a\otimes b \to 0$ ultranorm.
As usual, we split it into two:
\begin{equation*}
a_\alpha \otimes b_\alpha - a \otimes b
\ = \ a_\alpha \otimes(b_\alpha -b )\  +\  (a_\alpha - a) \otimes b.
\end{equation*}
For the first term, note
\begin{align*}
    \langle a_\alpha \otimes (b_\alpha  - b),
    a_\alpha \otimes (b_\alpha  - b) \rangle
    & \ = \ 
    (b_\alpha - b)^* \,\varphi (a_\alpha^* a_\alpha) \,(b_\alpha - b) \\
    & \ \leq \ 
    B^2 \|\varphi\| \,(b_\alpha - b)^*  (b_\alpha - b),
\end{align*}
which converges to~$0$ ultraweakly as~$b_\alpha \to b$ ultrastrongly.
The second
\begin{align*}
    \langle
        (a_\alpha - a) \otimes b
        (a_\alpha - a) \otimes b \rangle
        & \ = \ 
        b^* \,\varphi((a_\alpha - a)^* (a_\alpha - a)) \, b,
\end{align*}
converges ultraweakly to~$0$
as and~$ b^* \varphi(\,\cdot\,) b$ is normal
and~$(a_\alpha - a)^*(a_\alpha - a)$ converges ultraweakly to~$0$
as~$a_\alpha \to a$ ultrastrongly. \qed
\end{point}
\end{point}
\end{parsec}

\begin{parsec}{1670}%
\begin{point}{10}[paschke-tensor]{Theorem}%
Suppose~$\varphi_i \colon \scrA_i \to \scrB_i$
    is an ncp-map between von Neumann algebras
    with Paschke dilation~$(\scrP_i, \varrho_i, h_i)$
    for~$i=1,2$.
Then~$(\scrP_1 \otimes \scrP_2, \varrho_1 \otimes \varrho_2, h_1\otimes h_2)$
    is a Paschke dilation of~$\varphi_1 \otimes \varphi_2$.
    Furthermore
\begin{equation}\label{tensor-of-paschke-dilation-spaces}
(\scrA_1 \otimes_{\varphi_1} \scrB_1)
\otimes (\scrA_2 \otimes_{\varphi_2} \scrB_2)
\ \cong \ 
(\scrA_1 \otimes \scrA_2) \otimes_{\varphi_1\otimes \varphi_2} (\scrB_1 \otimes \scrB_2).
\end{equation}
\spacingfix{}
\begin{point}{20}{Proof}%
We will prove \eqref{tensor-of-paschke-dilation-spaces} first.
We proceed as follows.
To start, we prove that~$
(\scrA_1 \odot \scrB_1) \odot (\scrA_2 \odot \scrB_2)$
    (with the obvious inclusion)
    is ultranorm dense in~$(\scrA_1 \otimes_{\varphi_1} \scrB_1)
\otimes (\scrA_2 \otimes_{\varphi_2} \scrB_2)$.
Then we show that~$
(\scrA_1 \odot \scrA_2) \odot (\scrB_1 \odot \scrB_2) $
is ultranorm dense in~$
(\scrA_1 \otimes \scrA_2) \otimes_{\varphi_1\otimes \varphi_2} (\scrB_1 \otimes \scrB_2)$.
We continue and show that the natural bijection
$U_0\colon (\scrA_1 \odot \scrB_1) \odot (\scrA_2 \odot \scrB_2) \cong
(\scrA_1 \odot \scrA_2) \odot (\scrB_1 \odot \scrB_2) $
preserves the induced inner products and
    so extends uniquely to an
    isomorphism~\eqref{tensor-of-paschke-dilation-spaces}.
Using this isomorphism,
    we prove that the tensor product of the standard Paschke dilations
    is a dilation of the tensor product.
    From this, we finally derive the stated result.
\begin{point}{30}[paschke-iso-tensor-density]{Density}%
By construction, $\scrA_1 \odot \scrB_1$ is an ultranorm-dense
        $\scrB_1$-submodule
    of~$\scrA_1 \otimes_{\varphi_1} \scrB_1$.
Similarly~$\scrA_2 \odot \scrB_2 $ is ultranorm dense
    in~$\scrA_2 \otimes_{\varphi_2} \scrB_2$.
Thus by \sref{exttensor-dense-subsets}
    $(\scrA_1 \odot \scrB_1) \odot(
    \scrA_2 \odot \scrB_2) $ is ultranorm dense in~$(\scrA_1 \otimes_{\varphi_1} \scrB_1)
\otimes (\scrA_2 \otimes_{\varphi_2} \scrB_2)$.
By \sref{tensor}
        (with~\sref{ultraclosed} and \sref{dense-subalgebra}),
    we know~$\scrA_1 \odot \scrA_2$
    is an ultrastrongly-dense subalgebra
    of~$\scrA_1 \otimes \scrA_2$.
Similarly~$\scrB_1 \odot \scrB_2$
    is ultrastrongly dense in~$\scrB_1 \otimes \scrB_2$.
Hence~$(\scrA_1 \odot \scrA_2) \odot 
    (\scrB_1 \odot \scrB_2)$ is ultranorm dense
    in
    $(\scrA_1 \otimes \scrA_2) \otimes_{\varphi_1 \otimes \varphi_2}
    (\scrB_1 \otimes \scrB_2)$
        by~\sref{dilationspace-dense-subset},
\end{point}
\begin{point}{40}{Inner product preservation}%
Clearly
\begin{alignat*}{2}
    &
    \bigl\langle\,
    (a_1 \otimes b_1)\otimes (a_2 \otimes b_2)\,,\,
    (\alpha_1 \otimes \beta_1 )\otimes( \alpha_2 \otimes \beta_2)
    \,\bigr\rangle \\
    & \qquad\ = \ 
    \langle a_1 \otimes b_1 , \alpha_1 \otimes \beta_1 \rangle \,\otimes\,
    \langle a_2 \otimes b_2 , \alpha_2 \otimes \beta_2 \rangle \\
    & \qquad\ = \ 
    (b_1^* \varphi_1( a_1^* \alpha_1) \beta_1)
    \,\otimes \,
    (b_2^* \varphi_2( a_2^* \alpha_2) \beta_2) \\
    & \qquad \ = \ 
    (b_1 \otimes b_2)^* \,(\varphi_1 \otimes \varphi_2)
    \bigl((a_1 \otimes a_2)^* (\alpha_1 \otimes \alpha_2)\bigr)
    \,(\beta_1 \otimes \beta_2) \\
    & \qquad \ = \ 
    \bigl\langle\,
    (a_1\otimes a_2)\otimes (b_1 \otimes b_2) \,,\,
    (\alpha_1\otimes \alpha_2)\otimes (\beta_1 \otimes \beta_2)
    \,\bigr\rangle;
\end{alignat*}
so $U_0 \colon (\scrA_1 \odot \scrB_1) \odot (\scrA_2 \odot \scrB_2) \cong
(\scrA_1 \odot \scrA_2) \odot (\scrB_1 \odot \scrB_2) $
preserves the inner product by linearity.
\end{point}
\begin{point}{50}%
As~$U_0$ preserves the inner product,
    it is ultranorm continuous.
So with the usual reasoning, there is a unique bounded
    $\scrB_1 \odot \scrB_2$-linear extension
\begin{equation*}
    U_1\colon (\scrA_1 \otimes_{\varphi_1} \scrB_1) \odot (\scrA_2 \otimes_{\varphi_2} \scrB_2) \ \to \ 
        (\scrA_1 \otimes \scrA_2) \otimes_{\varphi_1 \otimes \varphi_2} (\scrB_1 \otimes \scrB_2).
\end{equation*}
In turn, by the universal property of the self-dual exterior tensor product,
    there exists a unique
    bounded~$\scrB_1\otimes \scrB_2$-linear extension
\begin{equation*}
    U\colon (\scrA_1 \otimes_{\varphi_1} \scrB_1) \otimes (\scrA_2 \otimes_{\varphi_2} \scrB_2) \ \to\ 
        (\scrA_1 \otimes \scrA_2) \otimes_{\varphi_1 \otimes \varphi_2} (\scrB_1 \otimes \scrB_2).
\end{equation*}
    By the ultranorm density~\sref{paschke-iso-tensor-density}
    and \sref{innerprod-ultraweak}
    we know~$U$ also preserves the inner product.
    Also from~\sref{paschke-iso-tensor-density},
        it follows that~$U$ is surjective.
Thus~$U$ is an isomorphism with~$U^{-1} = U^*$.
    We have shown~\eqref{tensor-of-paschke-dilation-spaces}.
By \sref{ba-ext-tensor-pres}
there is a unique nmiu-isomorphism
\begin{equation*}
    \vartheta\colon
    \scrB^a (\scrA_1 \otimes_{\varphi_1} \scrB_1) \otimes
     \scrB^a(\scrA_2 \otimes_{\varphi_2} \scrB_2) \ \to\ 
    \scrB^a
     ((\scrA_1 \otimes_{\varphi_1} \scrB_1) \otimes
     (\scrA_2 \otimes_{\varphi_2} \scrB_2))
\end{equation*}
    fixed by~$\vartheta(S \otimes T) = S \otimes T$.
Let~$(\scrB^a((\scrA_1 \otimes \scrA_2) \otimes_{\varphi_1 \otimes \varphi_2} 
        (\scrB_1 \otimes \scrB_2)), \varrho, h)$
        and~$(\scrB^a(\scrA_i \otimes_{\varphi_i} \scrB_i), \varrho_i, h_i)$
    (for~$i=1,2$)
    denote the Paschke dilations
    of~$\varphi_1 \otimes \varphi_2$  and~$\varphi_i$, respectively,
    constructed in~\sref{existence-paschke}.
\end{point}%
\begin{point}{60}%
    Our next goal is to show
    \begin{align}\label{tensorpaschkesecondgoal}
    \ad_{U^*} \after \vartheta
                \after (\varrho_1 \otimes \varrho_2) &\ =\  \varrho
                & h \after \ad_{U^*} \after \vartheta 
                       & \ =\  h_1 \otimes h_2,
    \end{align}
    which is sufficient to show that~$(
        \scrB^a(\scrA_1 \otimes_{\varphi_1} \scrB_1)
        \otimes \scrB^a(\scrA_2 \otimes_{\varphi_2} \scrB_2),
        \varrho_1 \otimes \varrho_2,
        h_1 \otimes h_2)$
    is a Paschke dilation of~$\varphi_1 \otimes \varphi_2$.
    We start with the equation on the left.
\begin{align*}
    &\bigl((\ad_{U^*} \after \vartheta \after  (\varrho_1 \otimes \varrho_2)) (\alpha_1 \otimes \alpha_2)\bigr)
    \, (a_1 \otimes a_2) \otimes (b_1 \otimes b_2) \\
    &\qquad \ = \ 
    \bigl((\ad_{U^*} \after \vartheta) (\varrho_1 (\alpha_1)\otimes \varrho_2(\alpha_2)) \bigr)
    \, (a_1 \otimes a_2) \otimes (b_1 \otimes b_2) \\
    &\qquad \ = \ 
    U \,(\varrho_1 (\alpha_1)\otimes \varrho_2(\alpha_2)) \, U^*
    \, (a_1 \otimes a_2) \otimes (b_1 \otimes b_2) \\
    &\qquad \ = \ 
    U \,(\varrho_1 (\alpha_1)\otimes \varrho_2(\alpha_2)) 
    \, (a_1 \otimes b_1) \otimes (a_2 \otimes b_2) \\
    &\qquad \ = \ 
    U \,  ((\alpha_1 a_1) \otimes b_1) \otimes ((\alpha_2 a_2) \otimes b_2) \\
    &\qquad \ = \ 
     ((\alpha_1 \otimes \alpha_2)\cdot (a_1 \otimes a_2)) \otimes (b_1 \otimes b_2) \\
    &\qquad \ = \ 
    \varrho(\alpha_1 \otimes\alpha_2) \,(a_1 \otimes a_2) \otimes (b_1 \otimes b_2).
\end{align*}
As the linear span of
    elements of the form~$(a_1 \otimes a_2) \otimes (b_1 \otimes b_2)$
    is ultranorm dense in~$(\scrA_1 \otimes \scrA_2)
                            \otimes_{\varphi_1 \otimes \varphi_2}
                            (\scrB_1 \otimes \scrB_2) $,
    we see~$(\ad_{U^*} \after \vartheta \after (\varrho_1 \otimes \varrho_2))
            (\alpha_1 \otimes \alpha_2)
            = \varrho(\alpha_1 \otimes \alpha_2)$
            for all~$\alpha_1\in\scrA_1$ and~$\alpha_2\in\scrA_2$.
In turn, as the linear span of elements of the form~$\alpha_1 \otimes \alpha_2$
    is ultrastrongly dense in~$\scrA_1 \otimes \scrA_2$
    and ncp-maps are ultrastrongly continuous,
    we see~$\ad_{U^*} \after \vartheta \after (\varrho_1 \otimes \varrho_2) =
        \varrho$.
        We continue with the equation on the right
        of~\eqref{tensorpaschkesecondgoal}.
\begin{align*}
    (h \after \ad_{U^*} \after \vartheta) (T \otimes S)
        & \ = \ \langle (1 \otimes 1) \otimes (1 \otimes 1),
                        \ U\, (T \otimes S)\, U^*(1 \otimes 1) \otimes (1 \otimes 1)
                        \rangle \\
        & \ = \ \langle U^* \,(1 \otimes 1) \otimes (1 \otimes 1),
                        \  (T \otimes S)\, U^*(1 \otimes 1) \otimes (1 \otimes 1)
                        \rangle \\
        & \ = \ \langle (1 \otimes 1) \otimes (1 \otimes 1),
                        \  (T \otimes S)\, (1 \otimes 1) \otimes (1 \otimes 1)
                        \rangle \\
        & \ = \
        \langle 1 \otimes 1, \,T\, 1\otimes 1\rangle \,\otimes\,
        \langle 1 \otimes 1, \,S\, 1\otimes 1\rangle \\
        & \ = \
        h_1(T) \otimes h_2(S) \\
        & \ = \
        (h_1\otimes h_2)(T \otimes S)
\end{align*}
From ultrastrong density, again, it
    follows~$h \after \ad_{U^*} \after \vartheta = h_1 \otimes h_2$. \qed
\end{point}
% \begin{point}%
% Assume that~$(\scrP_i, \varrho_i', h_i')$
%     is any (other) Paschke dilation of~$\varphi_i$ (for~$i=1,2$).
% Then there are isomorphisms~$\beta_1,\beta_2$
%     with~$\beta_i \after \varrho_i = \varrho_i'$
%     and~$h_i' \after \beta_i = h_i$.
% So~$\varrho'_1 \otimes \varrho'_2
%         = (\beta_1 \otimes \beta_2) \after (\varrho_1 \otimes \varrho_2)$
%         and~$(h_1' \otimes h_2') \after (\beta_1 \otimes \beta_2) = h_1 \otimes h_2$.
% As~$\beta_1 \otimes \beta_2$ is an isomorphism,
%     the previous shows that~$(\scrP_1 \otimes \scrP_2, \varrho_1' \otimes \varrho_2', h_1' \otimes h_2')$ is also a Paschke dilation of~$\varphi_1 \otimes \varphi_2$. \qed
% \end{point}
\end{point}
\end{point}
\end{parsec}

\section{Pure maps}
\begin{parsec}{1680}%
\begin{point}{10}%
Schr\"odinger's equation is invariant under the reversal of time
    and so the isolated quantum mechanical processes (ncp-maps) described by
    it are invertible.
In stark contrast, the ncp-maps corresponding to measurement and discarding
    are rarely invertible.
There are, however, many  processes in between
    which might not be invertible, but are still pure enough to be `reversed'
        in some consistent and useful fashion.
For instance, the ncp-map~$\ad_V$
    has the obvious reversal~$(\ad_V)^\dagger = \ad_{V^*}$
        (which is in general not the inverse.)
    Stinespring's theorem (\sref{stinespring-theorem})
        implies that every ncp-map into~$\scrB(\scrH)$
    splits as a reversible~$\ad_V$ after a (possibly) non-reversible
    nmiu-map~$\varrho$.
Using Paschke's dilation,
    we see that every ncp-map factors
    as an nmiu-map~$\varrho$ before some particular~ncp-map~$h$.
Is this~ncp-map~$h$ reversible in some sense?
In~\sref{paschke-pure} we will see that these~$h$
    are pure (in the sense of \sref{pure})
    and in~\sref{dagger-theorem}
    we will see that these pure maps are (in a sense) reversible.
\begin{point}{20}%
The reversibility of pure maps (in a much more general setting)
    is the main result of the next chapter.
In the remainder of this chapter, we study purity of ncp-maps a bit more.
We consider a seemingly unrelated  alternative notion of purity,
    which is based on Paschke's dilation,
    and show (in \sref{paschke-pure})
    that this alternative notion is in fact equivalent to the original.
We end this chapter by giving a different proof
    of the ncp-extremeness of pure and nmiu-maps.
\end{point}
\begin{point}{30}[dils-pure-discussion]%
Before we continue, let's rule out generalizations
    of more familiar notions of purity
    to arbitrary ncp-maps.
\begin{enumerate}
\item
A state~$\varphi\colon \scrA \to \C$ is called pure if it is an extreme
    point among all states.
It is unreasonable to define an arbitrary ncp-map to be pure
    if it is extreme, as every nmiu-map (including 
    von Neumann measurements and discarding)
    is extreme among the unital ncp-maps.
This follows from~\sref{rigid-ncp-extreme} and~\sref{nmiu-rigid},
    but will be proven again in~\sref{nmiu-ncp-extreme}.
In fact, in \sref{ncp-extreme-comp}
    we show that every ncp-map is the composition
    of two extreme maps.

\item
Inspired by the GNS-correspondence between pure states and
    irreducible representations,
St\o rmer defines a map~$\varphi\colon \scrA \to \scrB(\scrH)$ to be pure
    if the only maps below~$\varphi$ in the completely positive ordering
    are scalar multiples of~$\varphi$.
For every central element~$z$ of~$\scrA$,
    the map~$a \mapsto za$ is below~$\id_\scrA$
    and so the identity on non-factors is not pure in the sense
        of St\o rmer,
        even though the identity map is clearly reversible.
\end{enumerate}
\spacingfix{}
\end{point}
\begin{point}{40}%
Let's recall how we came to the notion of purity in \sref{pure},
    by generalizing~$\ad_V$.
Pick~$V \colon \scrH \to \scrK$.
    By polar decomposition, see~\sref{polar-decomposition}, we know that~$V = U A$
    for~$A \equiv \sqrt{V^*V} \colon \scrH \to \ceil{A} \scrH$
    and some isometry~$U \colon \ceil{A} \scrH \to \scrK$.
Thus we have~$\ad_V = \ad_A \after \ad_U$.
These two maps admit dual universal properties:
    $\ad_U$ is a \emph{corner}
    and~$\ad_A$ is a \emph{filter}.
An ncp-map~$\varphi$ is pure
    if it is the composition of a filter after a corner.
It turns out that~$\varphi$ is pure in this sense
    if and only if the map on the left-hand side of
    its Paschke dilation is surjective.
\end{point}
\end{point}
\end{parsec}

\begin{parsec}{1690}%
\begin{point}{10}%
We recall the definition of corner
        from \cite[dfn.~2]{westerbaan2016universal} and \sref{corner}.
\end{point}
\begin{point}{20}[dils-corner]{Definition}%
An ncp-map~$h \colon \scrA \to \scrB$
is a \Define{corner}\index{corner} for~$a \in [0,1]_\scrA$ if~$h(a)=h(1)$
    and
    \begin{quote}
        for every ncp-map $f\colon \scrA \to \scrC$
            with~$f(1)=f(a)$,
            there is a unique ncp-map~$f'\colon \scrB \to \scrC$
            with~$f = f' \after h$.
    \end{quote}
\spacingfix{}
\begin{point}{30}{Remark}%
In the next chapter we will study maps with a similar universal property
    as corners in a more general setting,
    but with all arrows in the reverse direction.
To help alleviate some confusion caused by the change of direction,
    we call those maps \emph{comprehensions},
    see \sref{dfn-comprehension}.  See also~\cite{effintro}.
\end{point}
\end{point}
\begin{point}{40}[standard-corner-dils]{Example}%
    The map~$h_a \colon \scrA \to \floor{a} \scrA \floor{a}$ given
    by~$b \mapsto \floor{a} b \floor{a}$
    is a corner of~$a \in [0,1]_\scrA$.
    We call~$h_a$ the \Define{standard corner}\index{corner!standard} of~$a$.
    See \sref{dfn-standard-corner-and-filter}, \sref{prop-corner}
        or \cite[prop.~5]{westerbaan2016universal}.
\end{point}
\begin{point}{50}[h-is-corner-for-unital-map]{Lemma}%
If~$(\scrP, \varrho, h)$ is a Paschke dilation of a unital ncp-map,
    then~$h$ is a corner.
\begin{point}{60}{Proof}%
(This is a different proof of our similar result \cite[cor.~27]{wwpaschke}.)
Let~$\varphi \colon \scrA \to \scrB$ be any unital ncp-map.
It is sufficient to prove~$h$
    is a corner for the
    standard dilation~$(\scrP , \varrho, h)$
    constructed in \sref{existence-paschke},
    so~$\scrP \equiv \scrB^a(\scrA \otimes_\varphi \scrB)$
    and~$h(T) = \langle e, T e\rangle$, where~$e \equiv 1\otimes 1$.
    Write~$p \equiv \ketbra{e}{e}$.
From~$\langle e,e\rangle = \varphi(1) = 1$
    it follows~$p$ is a projection.
Define~$\vartheta\colon \scrB \to p \scrP p$
    by~$\vartheta(b) = \ketbra{e \cdot b}{e}$.
    Some straightforward computations (using~\sref{hilbmodketbrarules})
    show~$\vartheta$ is an miu-map
    with inverse~$pTp \mapsto h(pTp)$.
So~$\vartheta$ is an miu-isomorphism and thus also normal.
    Furthermore~$h = \vartheta^{-1}(p (\,\cdot\,) p)$
     and so~$h$ is indeed a corner (of~$p$). \qed
\end{point}
\end{point}
\begin{point}{70}%
We recall the definition of filter, see \sref{filter}. 
\end{point}
\begin{point}{80}[dils-def-filter]{Definition}%
An ncp-map~$c \colon \scrA \to \scrB$
is a \Define{filter}\index{filter} for~$b \in\scrB$, $b \geq 0$ if~$c(1)\leq b$
    and
    \begin{quote}
        for every ncp-map $f\colon \scrC \to \scrB$
            with~$f(1) \leq b$,
            there is a unique ncp-map~$f'\colon \scrC \to \scrA$
            with~$f = c \after f'$.
    \end{quote}
\spacingfix{}
\begin{point}{90}{Remark}%
The direction-reversed counterpart
    of filters are called \emph{quotients},
    see~\sref{dfn-quotient} and~\cite{effintro} .
Filters were introduced by us in \cite[dfn.~2]{westerbaan2016universal}
    under the name \emph{compressions}.
To avoid confusion with \cite{alfsen2012}
    and be closer to e.g.~\cite{wilce2016royal},
    we changed the name from compression to filter.
Filters are named after polarization filters:
    the action of a polarization filter on the polarization
        of a photon is the same as a filter for
        the projection of the space of polarizations
        the polarization filter allows to pass without disturbance.
Combining two polarization filters corresponds
    to composing the corresponding filters.
This composed filter is in general a filter for an effect instead
    of a projection.
This corresponds with the fact that there might not be a polarization
    of the photon anymore,
    with which it can pass both filters undisturbed.
\end{point}
\end{point}
\begin{point}{100}[dils-stand-filter]{Example}%
    The map~$c_b\colon \ceil{b} \scrB\ceil{b} \to \scrB $
        given by~$a \mapsto \sqrt{b} a \sqrt{b}$ is a filter
        of~$b\in \scrB$, $b\geq 0$.
        We call~$c_b$ the \Define{standard filter}\index{filter!standard} of~$b$.
    See \sref{canonical-filter}, \sref{dfn-standard-corner-and-filter}
        or \cite[prop.~6]{westerbaan2016universal}.
\end{point}
\begin{point}{110}[dils-filter-basics-exercise]{Exercise}%
Let~$\varphi\colon \scrA \to \scrB$ be any ncp-map
    between von Neumann algebras.
\begin{enumerate}
\item
    Let~$c \colon \scrB \to \scrC$ be any filter.
    Show that if~$(\scrP, \varrho, h)$ is a Paschke dilation
            of~$\varphi$,
    that then~$(\scrP, \varrho, c \after h)$ is a Paschke dilation
        of~$c \after \varphi$.
\item
    Assume~$c'\colon \scrC' \to \scrB$ is a filter of~$\varphi(1)$.
    Prove that there is a unique unital ncp-map~$\varphi'$
        with~$\varphi = c \after \varphi'$.
    Conclude that if~$(\scrP, \varrho, h)$ is a Paschke dilation
            of~$\varphi'$,
            that then~$(\scrP, \varrho, c' \after h)$
            is a Paschke dilation of~$\varphi$.
\end{enumerate}
(See \sref{dils-abstract-basics} and \sref{quotient-total} for the proofs
        in the more general case.)
\end{point}
\begin{point}{120}[dils-filters-injective]{Exercise}%
Derive (perhaps from~\sref{dils-stand-filter} and~\sref{mult-cancellation}),
    that filters are injective.
\end{point}
\end{parsec}
\begin{parsec}{1700}%
\begin{point}{10}[dils-def-pure]{Definition}%
    Filters, corners and their compositions are called \Define{pure}.
    \index{pure!ncp-map}
\begin{point}{11}{Remark}%
This definition is equivalent to our previous definition
    of purity given in \cite[dfn.~21]{wwpaschke};
    see~\sref{pure-fundamental}.
\end{point}
\end{point}
\begin{point}{20}[dils-examples-pure]{Examples}%
Isomorphisms are pure.  Less trivial examples follow.
\begin{enumerate}
\item
The pure maps~$\scrB(\scrH) \to \scrB(\scrK)$
    are precisely of the form~$\ad_T$
    for some bounded operator~$T\colon \scrK \to \scrH$.
\item
Let~$(\scrP, \varrho, h)$ be a Paschke dilation of some ncp-map~$\varphi$.
In \sref{h-is-corner-for-unital-map} we saw that if~$\varphi$ is unital,
        that then~$h$ is a corner (and thus pure).
In general, $h = c \after h'$ for some filter~$c$ and corner~$h'$
    by \sref{dils-filter-basics-exercise},
        and so~$h$ is pure.
\end{enumerate}
\spacingfix{}
\end{point}
\begin{point}{30}{Remark}%
By definition, the set of pure maps is closed under composition.
In the next chapter we will see (\sref{dagger-theorem}) that
    there is a unique dagger~$\dagger$
    which turns the subcategory of von Neumann
    algebras with pure maps into a~$\dagger$-category
    such that~$f^\dagger$ is contraposed to~$f$;
            $\diamond$-positive maps are~$\dagger$-positive
            and~$\dagger$-positive maps have unique roots.
For this unique~$\dagger$, we have~$(\ad_V)^\dagger = \ad_{V^*}$.
\end{point}
\begin{point}{40}[surjective-nmiu]{Exercise}%
Show that any surjective nmiu-map between von Neumann algebras
        is a corner of a central projection, hence pure
        --- and conversely,
        that any corner of a central projection is a surjective nmiu-map.
\end{point}
\end{parsec}

\begin{parsec}{1710}%
\begin{point}{10}%
We will prove that a map is pure if and only if
    the left-hand side map of its Paschke dilation
    is surjective.  We need the following result.
\end{point}
\begin{point}{20}[paschke-corner]{Theorem}%
Let~$\scrA$ be a von Neumann algebra together with a projection~$p \in \scrA$.
Then a Paschke dilation of the standard corner~$h_p \colon \scrA \to p\scrA p$,
$a \mapsto pap$
is given by~$(\cceil{p}\scrA, h_{\cceil{p}}, h'_p)$,
    where~$\cceil{p}$ is the central  carrier of~$p$,
    see \sref{cceil-fundamental};
    $h_{\cceil{p}}$ is the standard corner for~$\cceil{p}$
        and~$h'_p \colon \cceil{p} \scrA \to p\scrA p$
        is the restriction of~$h_p$.
\begin{point}{30}{Proof}%
(This is a simplified proof of the result we published earlier.
    \cite[thm.~28]{wwpaschke})
Let~$(\scrA \otimes_{h_p} p\scrA p, \varrho, h)$
    denote the Paschke dilation constructed in~\sref{existence-paschke}.
We proceed in three steps:
    first we show~$\scrA p$ is a self-dual Hilbert~$p\scrA p$-module,
    then we prove~$\scrA \otimes_{h_p} p\scrA p \cong \scrA p$
    and finally we show~$\scrB^a(\scrA p) \cong \cceil{p} \scrA$.
\begin{point}{40}%
Note $\scrA p$ is a pre-Hilbert $p \scrA p$-module
    with scalar multiplication~$(\alpha p)\cdot (pap) = \alpha pap$
    and inner product~$\left<\alpha p, ap\right> = p\alpha^*ap$.
Since by the C$^*$-identity the norm on~$\scrA p$
    as a Hilbert module coincides with that as a subset of~$p \scrA p$;
    $\scrA p$ is norm closed in~$\scrA$ \emph{and}~$\scrA$ is norm complete,
    we see~$\scrA p$ is a Hilbert~$p \scrA p$-module.

Recall~$\scrA$ itself is a self-dual Hilbert~$\scrA$-module.
The uniformity induced on~$\scrA p$ by the ultranorm uniformity
    of~$\scrA$ coincides with the ultranorm uniformity of~$\scrA p$ itself.
The module $\scrA p$ is ultranorm closed in~$\scrA$
    as~$x_\alpha \to x$ ultranorm
    implies~$x_\alpha p \to xp$ ultranorm.
As self-duality is equivalent to ultranorm completeness,
    we see~$\scrA p$ is self dual as well.
\end{point}
\begin{point}{50}%
Concerning~$\scrA \otimes_{h_p} p \scrA p \cong \scrA p$,
    define~$B\colon \scrA \times  p \scrA p \to \scrA p$
    by~$B(\alpha, pap) = \alpha pap$.
A straightforward computation
    shows~$B$ is a~$h_p$-compatible complex bilinear map
    (see \sref{phi-compatible-paschke})
    and so by \sref{existence-paschke},  there is a unique bounded module
    map~$U \colon \scrA \otimes_{h_p} p \scrA p \to \scrA p$
    fixed by~$U (\alpha \otimes pap) = B(\alpha,pap) \equiv \alpha pap $.
We will show~$U$ is an isomorphism.
   Clearly~$U$ is surjective.

An easy computation shows~$a \otimes p\alpha p - ap\alpha p \otimes p = 0$
    for all~$a, \alpha \in \scrA$.
So the set~$\{a \otimes p;\ a\in \scrA\}$
    is ultranorm dense in~$\scrA \otimes_{h_p} p \scrA p$.
From this and
    $\left<\alpha \otimes p, a\otimes p\right> = 
        p \alpha^* a p
        =\left< U \alpha \otimes p, U a \otimes p\right> $
        we see~$U$ preserves the inner product.
Hence~$U$ is an isomorphism with~$U^* = U^{-1}$
    given by~$U^* a = a\otimes p$.
\end{point}
\begin{point}{60}%
To show~$\scrB^a(\scrA p) \cong \cceil{p}\scrA$,
we will first find a convenient orthonormal basis of~$\scrA p$.
By \sref{cceil-sum},
    there are partial isometries~$(u_i)_i$
    with~$\sum_i u_iu_i^* = \cceil{p}$
    and~$u_i^*u_i \leq p$.
Write~$E = \{u_i; \ i \in I\}$.
As~$u_i = u_i p$ and~$u_i= u_i u_i^* u_i$
    we see~$E$ is an orthonormal subset of~$\scrA p$.
For any~$ap \in \scrA p$, we have
\begin{equation*}
    ap \ =\  a \cceil{p} p
       \ =\  \cceil{p} a p
       \ = \ \bigl(\sum_i u_i u_i^*\bigr) ap
       \ = \ \sum_i u_i u_i^* ap
       \ = \ \sum_i u_i \left<u_i,ap\right>,
\end{equation*}
where the sums converge ultrastrongly in~$\scrA$.
Thus the last sum converges ultranorm in~$\scrA p$ as well.
Hence~$E$ is an orthonormal basis of~$\scrA p$.

Define~$\varrho_0 \colon \scrA \to \scrB^a(\scrA p)$
    by~$\varrho_0 = \ad_{U^*} \after \varrho$.
Note~$\varrho_0$ is an nmiu-map and~$\varrho_0(\alpha)ap = \alpha ap$.
The map~$\varrho_0$ is surjective:
    indeed for any~$T \in \scrB^a(\scrA p)$
\begin{equation*}
    T (ap) \ = \ \sum_i T( u_i ) u_i^* ap
           \ = \ \bigl(\sum_i T( u_i ) u_i^*\bigr) ap
           \ = \ \varrho_0 \bigl(\sum_i T( u_i ) u_i^*\bigr) ap.
\end{equation*}
Next, we show~$\ceil{\varrho_0} = \cceil{p}$.
It suffices to prove for all~$\alpha \in \scrA$,$\alpha \geq 0$, that
\begin{equation*}
    \alpha ap = 0 \text{ \ for all\  } a\in \scrA
    \quad\iff \quad
    \alpha \cceil{p} = 0.
\end{equation*}
    Clearly,
    $\alpha \cceil{p} = 0$ implies~$\alpha ap = \alpha \cceil{p} ap= 0$.
For the converse, assume that we have~$\alpha ap = 0$ for all~$a\in\scrA$.
Then in particular~$\alpha u_ip = 0$,
    hence~$\alpha \cceil{p} = \alpha \sum_i u_iu_i^*  = \sum_i \alpha u_ipu_i^*=0$.
    Indeed~$\ceil{\varrho_0} = \cceil{p}$.

    By~\sref{surjective-nmiu}
    there is an nmiu-isomorphism~$\varphi\colon
    \cceil{p} \scrA \to \scrB^a(\scrA p)$
    with~$\varrho_0 = \varphi \after h_{\cceil{p}}$.
Hence~$h_{\cceil{p}} = (\varphi^{-1} \after \ad_{U^*}) \after \varrho$.
We compute
\begin{equation*}
h'_p \after h_{\cceil{p}}
    \ = \ h_p
    \ =\  h \after \varrho
    \ =\ h \after \ad_U \after \varrho_0
    \ = \ h \after \ad_U \after \varphi \after h_{\cceil{p}}.
\end{equation*}
Using surjectivity of~$h_{\cceil{p}}$
and rearranging, we find~$h'_p \after (\varphi^{-1} \after \ad_{U^*}) = h$.
Thus indeed~$(\cceil{p} \scrA, h_{\cceil{p}}, h'_p)$
is a Paschke dilation of~$h_p$.  \qed
\end{point}
\end{point}
\end{point}
\begin{point}{70}[paschke-pure]{Theorem}%
    Let~$\varphi\colon \scrA \to \scrB$ be an ncp-map
    between von Neumann algebras with
    Paschke dilation~$(\scrP, \varrho, h)$.
The map~$\varphi$ is pure if and only if~$\varrho$ is surjective.
\begin{point}{80}{Proof}%
Assume~$\varrho$ is surjective.
The kernel of~$\varrho$ is an ultraweakly closed
two-sided ideal and so~$\ker\varrho = z \scrA$ for some central projection~$z$,
see \sref{prop:weakly-closed-ideal}.
The standard corner~$h_{z^\perp} \colon \scrA \to z^\perp \scrA$
    is the quotient-map of~$z \scrA$
    and so by the isomorphism theorem $\varrho$ is a corner as well.
In \sref{dils-examples-pure} we saw~$h$ is pure.
So~$\varphi$ is the composition of pure maps, hence pure itself.

For the converse, assume~$\varphi$ is pure.
For brevity, write~$p \equiv \ceil{\varphi}$.
There is a unique ncp-map~$c$ with~$\varphi = c \after h_p$,
    where~$h_{\ceil{\varphi}}$ is the standard corner of~$p$.
The map~$c$ is a filter,
    see \sref{square-f}.
By \sref{paschke-corner},
    the triple~$(\cceil{p}\scrA, h_{\cceil{p}}, h'_p)$
    is a Paschke dilation of~$h_p$,
    where~$h'_p$ denotes the restriction of~$h_p$ to~$\cceil{p} \scrA$.
So by \sref{dils-filter-basics-exercise}
    we see~$(\cceil{p}\scrA, h_{\cceil{p}}, c\after h'_p)$
    is a Paschke dilation of~$ \varphi$.
Clearly~$h_{\cceil{p}}$ is surjective.
If~$(\scrP', \varrho', h')$ is any other Paschke dilation of~$\varphi$,
    then~$\varrho' = \vartheta \after h_{\cceil{p}}$
    for some isomorphism~$\vartheta$
    and so~$\varrho'$ is surjective as well. \qed
\end{point}
\end{point}
\end{parsec}

\begin{parsec}{1720}[paschke-corresp-pure]%
\begin{point}{10}%
The order correspondence \sref{paschke-correspondence}
    has several corollaries.
\end{point}
\begin{point}{20}{Definition}%
Assume~$\varphi\colon \scrA \to \scrB$ is an ncp-map between von Neumann
    algebras.
We say~$\varphi$ is \Define{ncp-extreme}\index{ncp-extreme}
    if~$\varphi$ is an extreme point among
    the ncp-maps with the same value on~$1$~\cite{stormer1963positive}
    ---
    that is:
    for any~$0 < \lambda < 1$, ncp-maps~$\varphi_1,\varphi_2\colon \scrA \to \scrB$
    with~$\varphi_1(1) = \varphi_2(1) = \varphi(1)$
    we have
    \begin{equation*}
    \lambda \varphi_1  +  (1-\lambda)\varphi_2 \ =\  \varphi
    \quad \implies \quad
    \varphi_1\ =\ \varphi_2\ =\ \varphi.
    \end{equation*}
\end{point}
\spacingfix{}
\begin{point}{30}[ncp-extreme-paschke]{Theorem}%
Assume~$\varphi\colon \scrA \to \scrB$ is an ncp-map
    with Paschke dilation~$(\scrP, \varrho, h)$.
Then the following are equivalent.
\begin{enumerate}
    \item The map $h$ is injective on~$\varrho(\scrA)^\square$,
                the commutant of~$\varrho(\scrA)$.
    \item The map $h$ is injective on~$[0,1]_{\varrho(\scrA)^\square}$.
    \item The map~$\varphi$ is ncp-extreme.
\end{enumerate}
\spacingfix{}
\begin{point}{40}{Proof}%
(Based on \cite[prop.~1.4.6]{arveson}
    and \cite[thm.~5.4]{paschke}.)
We prove~$1 \Rightarrow 2 \Rightarrow 3 \Rightarrow 1$.
\begin{point}{50}{$1 \Rightarrow 2$}%
Trivial.
\end{point}
\begin{point}{60}{$2 \Rightarrow 3$}%
Assume~$h$ is injective on~$[0,1]_{\varrho(\scrA)^\square}$.
Suppose~$\varphi = \lambda \psi + (1 - \lambda) \psi'$
    for some~$0 < \lambda < 1$
    and ncp~$\psi,\psi'\colon \scrA \to \scrB$
    with~$\psi(1) = \psi'(1) = \varphi(1)$.
Note~$0 \leq_\ncp \lambda\psi \leq_\ncp \varphi$
hence~$\lambda\psi = \varphi_t$ for some~$t \in [0,1]_{\varrho(\scrA)^\square}$
by~\sref{paschke-correspondence}.
So
\begin{equation*}
    h(\lambda 1)
    \ =\  \varphi_{\lambda1}(1) \ =\  \lambda \varphi(1) 
    \ =\  \lambda \psi(1)
    \ =\  \varphi_t (1)  \ =\  h(t).
\end{equation*}
By the injectivity of~$h$,
    we find~$t = \lambda 1$,
    thus~$\lambda \psi = \lambda \varphi$
    and so~$\psi = \varphi$.
The proof of~$\psi' = \varphi$ is the same.
Hence~$\varphi$ is ncp-extreme.
\end{point}
\begin{point}{70}{$3 \Rightarrow 1$}%
Suppose~$\varphi$ is ncp-extreme.
Pick any~$a \in \varrho(\scrA)^\square$ with~$h(a) = 0$.
We want to show~$a = 0$.
Clearly~$\Real{h(a)} = h(\Real{a}) = 0$
and~$\Imag{h(a)} = h(\Imag{a}) = 0$,
so it is sufficient to prove~$\Real{a} = \Imag{a}=0$.
Thus we may assume without loss of generality that~$a$ is self adjoint.

By~\sref{dils-examples-pure}~$h$ is pure
    and thus splits as~$h = c \after h_p$
    for some filter~$c \colon p\scrP p \to \scrB$
    and the standard corner~$h_p\colon \scrP \to p\scrP p$
    for~$p = \ceil{h}$.
By \sref{dils-filters-injective} filters are injective
    and so from the assumption~$0 = h(a) = c(pap)$
    it follows~$pap = 0$.
Using \sref{cstar-positive} and some easy algebra,
    we can find~$0 < \mu,\lambda$
    such that~$\frac{1}{4} \leq \mu a + \lambda \leq \frac{3}{4}$.
Define~$t = \mu a + \lambda$.
For the moment, note~$\varphi_t$ is ncp.
As $pap = 0$, we see~$ptp = \lambda p$ and
so~$\frac{1}{4} p \leq  ptp = \lambda p
    \leq \frac{3}{4} p$.
Thus either~$0 < \lambda < 1$ or~$p = 0$.
If~$p = 0$ then~$\varphi=0$ hence~$\scrP = \{0\}$
    and~$a = 0$ as desired.
For the non-trivial case, assume~$0 < \lambda <1$.
Note~$\varphi_t(1) = h(t) = c(ptp) = \lambda c(p1p) = \lambda \varphi(1)$
    and similarly~$\varphi_{1-t}(1) = (1-\lambda) \varphi(1)$.
Define~$\psi_1 = \lambda^{-1} \varphi_t$
and~$\psi_2 = (1-\lambda)^{-1}\varphi_{1-t}$.
By the previous,~$\psi_1$ and~$\psi_2$ are ncp-maps
    with~$\psi_1(1)=\psi_2(1)=\varphi(1)$.
    Also~$\lambda \psi_1 + (1-\lambda) \psi_2 = \varphi$.
As~$\varphi$ is ncp-extreme,
    we must have~$\psi_1 = \psi_2 = \varphi$.
So~$\varphi_t = \lambda \varphi = \varphi_{\lambda 1}$,
    thus~$\lambda = t \equiv \mu a + \lambda$.
    Rearranging:~$\mu a = 0$.
    As~$\mu > 0$, we find~$a=0$, as desired. \qed
\end{point}
\end{point}
\end{point}
\begin{point}{80}[nmiu-ncp-extreme]{Corollary}%
Any nmiu-map~$\varrho\colon \scrA \to \scrB$ is ncp-extreme.
    (For a different proof, see e.g.~\cite[thm.~3.5]{stormer1963positive}.)
\begin{point}{90}{Proof}%
Easy: $(\scrB, \varrho, \id)$ is a Paschke dilation of~$\varrho$
and~$\id$ is injective.  \qed
\end{point}
\end{point}
\begin{point}{100}{Theorem}%
Any pure ncp-map~$\varphi\colon \scrA \to \scrB$ is ncp-extreme.
\begin{point}{110}{Proof}%
By \sref{pure-fundamental} $\varphi = c \after h_p$
    for a filter~$c\colon p\scrA p \to \scrB$
    and projection~$p \equiv \ceil{\varphi}$.
    Combining \sref{dils-filter-basics-exercise} and \sref{paschke-corner}
    we see~$(\cceil{p}\scrA, h_{\cceil{p}}, c\after h_p)$
    is a Paschke dilation of~$\varphi$.
Appealing to~\sref{ncp-extreme-paschke},
    we are done if we can show~$c \after h_p$
    is injective on~$[0,1]_{h_{\cceil{p}}(\scrA)^\square}$.
As~$h_{\cceil{p}}$ is surjective,
    $h_{\cceil{p}}(\scrA)^\square = Z(\cceil{p} \scrA)$.
    Assume~$t \in [0,1]_{Z(\cceil{p} \scrA)}$ with~$c(ptp)=0$
We want to show~$t = 0$.
As~$c$ is injective by~\sref{dils-filters-injective},
        we must have~$ptp = 0$.
So~$p \leq 1- \ceil{t}$
and~$\cceil{p} \leq 1- \ceil{t}$.
Hence~$\cceil{p} \leq \cceil{p} - \ceil{t} \leq \cceil{p}$.
So~$t \leq \ceil{t }=0$, as desired. \qed
\end{point}
\end{point}
\begin{point}{120}[ncp-extreme-comp]{Corollary}%
    Any~ncp-map is the composition of two ncp-extreme maps.
\end{point}
\end{parsec}

% vim: ft=tex.latex

\chapter{Diamond, andthen, dagger}\label{chapter2}

\begin{parsec}{1730}%
\begin{point}{10}%
In the previous chapter and in \cite{bram}
    we have studied categorical properties of von Neumann algebras.
In this chapter we change pace:
    we study categories that resemble the
    category of von Neumann algebras.
The goal of this line of study is to identify
    axioms which uniquely pick out the category of von Neumann algebras.
We do not reach this ambitious goal ---
    instead we will build up to a characterization of
    the existence of a unique $\dagger$-structure
    on the pure maps of a von Neumann algebras-like category.
\end{point}
\begin{point}{20}%
As a basic axiom we will require the categories we consider to be
    effectuses (which will be defined in \sref{dfn-effectus}).
For us the definition of an effectus will play a similar role
    to that of a topological space for a geometer.
In many applications,
    general topological spaces on their own are of little interest:
    the axioms are so weak that there are
    many (for the application) pathological examples.
    These weak axioms, however, are very expressive in the sense
    that they allow for the definition of many useful notions.
    Herein lies the strength of topological spaces ---
    as a stepping stone:
    many important classes of mathematical spaces
    are just plain topological spaces with a few additional axioms.
Similarly, there will be many pathological effectuses.
Their use for us lies in their expressiveness.
\end{point}
\begin{point}{30}%
Before we dive into effectuses,
    we need to learn about effect algebras (and related structures),
    which were introduced by mathematical physicists as
    a generalization of Boolean algebra to study fuzzy predicates
    in quantum mechanics.
For us, effectuses are the
    categorical counterpart of effect algebras.\footnote{Effect
    algebroids~\cite{roumen,roumen2016cohomology} are a different
    categorification with applications in cohomology.}
After this, we continue with a brief, but thorough development
    of the basic theory of effectuses.
There is a lot more to say about effectuses
    (see \cite{effintro}), which has been omitted
    to avoid straining the reader.
    We will indulge, however,
    in one tangent,
    which can be skipped:
    the study of abstract convex sets in
    \S\ref{more-aconvm}.
For the moment, think of an effectus as a generalization
    of the opposite category of von Neumann algebras:
    the objects represent the physical systems
    (or data types, if you like)
    and the maps represent the physical processes
    (or properly typed programs).
\end{point}
\begin{point}{40}%
The first two axioms we add on top of those of an effectus
    are the existence of quotients and comprehension.
It's helpful to discuss the origin of these axioms.
It started with the desire to axiomatize categorically
    the sequential product on a von Neumann algebra~$\scrA$,
    that is the operation~$\andthen{a}{b} \equiv  \sqrt{a} b \sqrt{a}$,
    which represents sequential measurement
    of first~$a$ \emph{andthen} $b$.
The sequential product is a quantum mechanical generalization
    of classical conjunction (logical `and').
In  contrast to the tame~$\wedge$ in a Boolean algebras,
    the sequential product does not obey a lot of insightful
    formulas.
The sequential product as a binary operation
        has received quite some attention \cite{gudder2008characterization,gheondea2004sequential,gudder2001sequential,li2011sequential,gudder2005open,shen2009not,gudder2005uniqueness,jun2009remarks,weihua2009uniqueness,tkadlec2008atomic,jia2010entropy,arias2004almost}.
We approach the sequential product in a rather different way:
    we take a step back and
    consider the map~$\asrt_a(b) \equiv  \andthen{a}{b}$
        for fixed~$a$.  ($\asrt_a$ is named after the \texttt{assert} statement
                    used in many programming languages,
                    see~\sref{asrt-remarks}.)
    This map~$\asrt_a$ factors in two:
\begin{equation*}
    \xymatrix@C+2pc@R-2pc{
        \scrA \ar[r]^\pi & \ceil{a} \scrA \ceil{a} \ar[r]^\xi
        & \scrA \\
        b \ar@{|->}[r] & \ceil{a}b\ceil{a}
            \ar@{|->}[r] & \sqrt{a}b\sqrt{a}.
    }
\end{equation*}
It turned out that both~$\pi$ and~$\xi$ have
    nice and dual (in a sense) defining universal properties
    which can be expressed in an effectus:
    $\pi$ is a comprehension and~$\xi$ is a quotient.\footnote{%
            Beware: an effectus quotient is not the same thing
        as a von Neumann-algebra quotient.}
\end{point}
\begin{point}{50}%
The existence of quotients and comprehension in some effectus~$C$
    is interesting on its own, but with
    two additional axioms (the existence of images and
    preservation of images under orthocomplementation),
    we get a surprisingly firm handle
    on the \emph{possibilistic} side of~$C$,
    by means of the existence
    of a certain functor~$(\ )_\diamond\colon C \to \mathsf{OMLat}$
        to the category of orthomodular lattices.
In the guiding example of von Neumann algebras,
    we have~$f_\diamond = g_\diamond$
        if and only
    there are no post-measurement~$a$ and initial state~$\omega$
    that can determine with certainty
    whether~$f$ or~$g$ has been performed.\footnote{%
    In symbols: $
    f_\diamond = g_\diamond \iff 
\bigl(\forall a,\omega.\  \omega(f(a))  =  0 \ \Leftrightarrow\ \omega(g(a))  =  0\bigr).$}
With this functor we can introduce several
    possibilistic notions of which~$\diamond$-adjoint
    and~$\diamond$-positive are the most important.
Because of the central role of~$\diamond$,
    we call an effectus with these axioms 
    a~$\diamond$-effectus.
\end{point}
\begin{point}{60}%
The axioms of a~$\diamond$-effectus do not force
    any coherence between the quotients and comprehension.
This has been a roadblock to the axiomatization
    of the sequential product for quite a while.
The key insight was the following:
    the map~$\asrt_a$ (i.e.~$b\mapsto \andthen{a}{b}$)
    is the unique~$\diamond$-positive map on~$\scrA$
    with~$\asrt_a(1) = a$.
We turn this theorem into an axiom:
    an~$\&$-effectus is
    a~$\diamond$-effectus
    where there are such unique~$\diamond$-positive maps~$\asrt_a$
    and where an additional polar decomposition axiom holds.
In an~$\&$-effectus, the predicates on an object carry a canonical
    binary operation~$\&$
    which is the intended sequential product in
    the case of von Neumann algebras.
\end{point}
\begin{point}{70}[pure-effectus]%
In an effectus, we call a map pure if it can be written
    as the composition of quotients and comprehensions.
In~\sref{paschke-pure} we saw that in the case of von Neumann algebras,
    a map is pure (in this sense) if and only if the corresponding Paschke
    embedding is surjective.
In particular, the pure maps~$\scrB(\scrH) \to \scrB(\scrK)$
    are exactly of the form~$\ad_v$
        (for bounded operators~$v\colon \scrK \to \scrH$)
    which carry a~$\dagger$-structure,
    namely~$(\ad_v)^\dagger = \ad_{v^*}$.
This begs the question: can we define a~$\dagger$ on all pure maps
    in an~$\&$-effectus?
The main result of this chapter
    is Theorem~\sref{dagger-theorem},
    which gives  necessary and sufficient conditions
    for an~$\&$-effectus
    to have a well-behaved~$\dagger$-structure on its pure maps.
After this apogee,
    we tie up some lose ends
    and discuss how these
    structures compare to those already known in the literature.
\end{point}
\end{parsec}
\section{Effect algebras and related structures}
\begin{parsec}{1740}%
\begin{point}{10}%
Before we turn to effectuses,
    it is convenient to introduce some lesser-known algebraic structures
    of which the effect algebra is the most important.
It will turn out that the set of predicates associated to an object in
    an effectus can be arranged into an effect algebra.
In this way, effect algebras play the same role
    for effectuses as Heyting algebras for toposes.
    First, we will recall the definition of partial commutative monoid (PCM)
    --- as later in \sref{cho-thm} 
    we will see that the partial maps between
    two objects in an effectus can be arranged into a PCM.
\end{point}
\begin{point}{20}[dfn-pcm]{Definition}%
    A \Define{partial commutative monoid} (\Define{PCM})~$M$
        \index{PCM}
        is a set~$M$ together with distinguished element~$0 \in M$
        and a partial binary operation~$\ovee$ such that
        for all~$a,b,c \in M$ 
        --- writing~$a \perp b$ whenever~$a \ovee b$ is defined
        --- we have
\begin{enumerate}
    \item \emph{(partial commutativity)}
        if $a \perp b$, then~$b \perp a$ and~$a \ovee b = b \ovee a$;
    \item \emph{(partial associativity)}
        if $a \perp b$ and~$a \ovee b \perp c$,
        then~$b \perp c$, ~$a \perp b \ovee c$
            and~$(a \ovee b) \ovee c = a \ovee (b \ovee c)$ \emph{and}
    \item \emph{(zero)}
        $0 \perp a$ and~$0 \ovee a = a$.
\end{enumerate}
A map~$f\colon M \to N$
    between PCMs~$M$ and~$N$
    is called an \Define{PCM homomorphism}
    \index{PCM!homomorphism}
    if for all~$a, b \in E$
    with~$a \perp b$,
    we have~$f(a) \perp f(b)$
    and~$f(a) \ovee f(b) = f(a \ovee b)$.
Write~$\Define{\textsf{PCM}}$ for the category
    \index{\textsf{PCM}}
    of PCMs with these homomorphisms. (Cf.~\cite{corefl,Wehrung2017})

For~$a,b$ in a PCM~$M$,
    we say~$\Define{a \leq b}$
    \index{*leq@$\leq$!in PCM/EA}
    iff~$a \ovee c = b$ for some~$c \in M$.
\end{point}
\begin{point}{30}[pcm-preorder]{Exercise}%
Show a PCM is preordered by~$\leq$.
\end{point}
\begin{point}{40}%
The partial commutativity and associativity of a PCM
    ensure that a sum only depends on which elements
    occur (and how often), instead of their order.
Formally:
    if~$(\cdots( x_1 \ovee x_2) \ovee \cdots) \ovee x_n$
    exists some~$x_1, \ldots, x_n$ in some PCM, then so
    does~$(\cdots( x_{\pi(1)} \ovee x_{\pi(2)}) \ovee \cdots) \ovee x_{\pi(n)}$
    for any permutation~$\pi$ of~$\{1,\ldots,n\}$ and
\begin{equation*}
    (\cdots( x_1 \ovee x_2) \ovee \cdots) \ovee x_n
    \ = \ (\cdots( x_{\pi(1)} \ovee x_{\pi(2)}) \ovee \cdots) \ovee x_{\pi(n)}.
\end{equation*}
    Thus parantheses are superfluous and we will often leave them out:
    \begin{equation*}
        \bigovee^n_{i=1} x_i \ \equiv \ \Define{x_1 \ovee \cdots \ovee x_n}
        \ \equiv\ 
    (\cdots( x_1 \ovee x_2) \ovee \cdots) \ovee x_n.
    \end{equation*}
\end{point}
\spacingfix{}
\end{parsec}

\begin{parsec}{1750}%
\begin{point}{10}[dfn-ea]{Definition}%
A PCM~$E$ together with distinguished element~$1 \in E$
    is called an \Define{effect algebra} (\Define{EA}) \cite{ea}
    \index{EA}
    \index{effect algebra}
    provided that
\begin{enumerate}
\item
    \emph{(orthocomplement)}
    for every~$a$
    there is a unique~$\Define{a^\perp}$
        \index{orthocomplement!in an effect algebra}
        \index{*perp@$(\ )^\perp$!element effect algebra}
   with~$a \ovee a^\perp = 1$ \emph{and}
\item
    \emph{(zero--one)}
    if~$a \perp 1$, then~$a = 0$.
\end{enumerate}
A map~$f\colon E \to F$
between effect algebras~$E$ and~$F$
is called an \Define{effect algebra homomorphism}
    \index{effect algebra!homomorphism}
if it is a PCM homomorphism that preserves~$1$;
    concretely:
\begin{enumerate}
    \item \emph{(additive)}
    $a \perp b$ implies $f(a) \perp f(b)$ and~$f(a)\ovee f(b) = f(a\ovee b)$
        \emph{and}
    \item \emph{(unital)}
    $f(1) = 1$.
\end{enumerate}
(It follows that~$f(0)=0$ and~$f(a^\perp) = f(a)^\perp$. See~\sref{exc-eamorphism}.)
Write~$\Define{\textsf{EA}}$
    \index{\textsf{EA}}
    for the category
    of effect algebras with these homomorphisms.
A subset~$D \subseteq E$ is a \Define{sub-effect algebra}
    \index{effect algebra!sub-} of~$E$
    if~$0,1 \in D$ and for any~$a,b \in D$
    with~$a\perp b$, we have~$a \ovee b \in D$
    and~$a^\perp \in D$.
\end{point}
\begin{point}{20}[eaexamples]{Examples}%
There are many examples of effect algebras
    --- we only give a few.
\begin{enumerate}
\item
The unit interval
$[0,1]$ with partial addition is an effect algebra ---
i.e.:~$x \perp y$
        whenever~$x +y \leq 1$ and then $x \ovee y = x + y$,
        $x^\perp = 1-x$.
\item
Generalizing the previous:
if~$G$ is an ordered group
with distinguished element~$1$,
then the order interval~$[0,1]_G \equiv \{x;\ x\in G;\ 0 \leq x\leq 1\}$
is an effect algebra
with~$x \perp y$ whenever~$x +y \leq 1$;
$x \ovee y = x + y$ and~$x^\perp = 1-x$.
\item
In particular,
    if~$\scrA$ is any von Neumann algebra,
    then the set of \Define{effects}~$[0,1]_\scrA$
        \index{effects}
        \index{*01@$[0,1]_\scrA$}
    forms an effect algebra
    with~$a \perp b$ whenever~$a +b \leq 1$;
    $a \ovee b = a + b$ and~$a^\perp = 1-a$.
The `effect' in effect algebra originates from this example.
\item
Any orthomodular lattice~$L$ (defined in \sref{dfn-orthomodular-lattice})
    is an effect algebra
    with the same orthocomplement,
    $x \perp y$ whenever~$x \leq y^\perp$
    and~$x \ovee y = x \vee y$.  (See e.g.~\cite[prop.~27]{basmsc}.)
\item
In particular,
    any Boolean algebra~$L$
    is an effect algebra
    with complement as orthocomplement,
    ~$x \perp y$ whenever~$x \wedge y = 0$ and
    $x \ovee y = x \vee y $.
\item
The one-element Boolean algebra~$1 \equiv \{0=1\}$
    is the final object in \textsf{EA}
    and the two-element Boolean algebra~$2 \equiv \{0,1\}$
    is the initial object in \textsf{EA}.
\end{enumerate}
\spacingfix{}
\end{point}
\begin{point}{30}[ea-product]{Exercise}%
Let~$E$ and~$F$ be effect algebras.
First show that the cartesian product~$E \times F$
    is an effect algebra with componentwise operations.
Show that this is in fact the categorical product of~$E$ and~$F$
    in \textsf{EA}.
(The category \textsf{EA} is in fact complete and cocomplete;
    for this and more categorical properties,
        see \cite{corefl}.)
\end{point}
\begin{point}{40}[ea-redund]{Exercise}%
There is a small redundancy in our definition of effect algebra:
show that the zero axiom ($x \ovee 0 = x$)
follows from the remaining axioms (partial commutativity,
partial associativity, orthocomplement and zero--one.)
\end{point}
\begin{point}{50}[eabasics]{Proposition}%
In any effect algebra~$E$ with~$a,b,c\in E$, we have
\begin{enumerate}
    \item \emph{(involution)}
        $a^{\perp\perp} = a$;
    \item
        $1^\perp= 0$ and~$0^\perp = 1$;
    \item \emph{(positivity)}
        if~$a\ovee b = 0$, then~$a = b= 0$;
    \item \emph{(cancellation)}
        if~$a\ovee c = b\ovee c$, then~$a = b$;
    \item the relation  $\leq$ (from \sref{dfn-pcm}) partially orders~$E$;
    \item $a \leq b$ if and only if $b^\perp \leq a^\perp$;
    \item if~$a \leq b$ and~$b \perp c$, then~$a \perp c$
        and~$a \ovee c \leq b \ovee c$ \emph{and}
    \item $a \perp b$ if and only if~$a \leq b^\perp$.
\end{enumerate}
\spacingfix{}
\begin{point}{60}{Proof}%
(These proofs are well-known,
    see for instance~\cite{dvurecenskij2013new}.)
By partial commutativity and definition
of orthocomplement both~$a^\perp \ovee a = 1$
and~$a^\perp \ovee  a^{\perp\perp} = 1$.
So by uniqueness of orthocomplement,
    we must have~$a= a^{\perp\perp}$, which is point 1.
Clearly~$0 \ovee 1 = 1$,
    so~$1^\perp = 0$ and~$0^\perp = 1$,
    which is point 2.
For point 3, assume~$a \ovee b = 0$.
Then~$a \ovee b \perp 1$
    and so by partial associativity~$b \perp 1$.
    By zero--one, we get~$b = 0$.
    Similarly~$a=0$, which shows point 3.
For point 4, assume~$ a \ovee c = b \ovee c$.
From partial associativity and commutativity, we get
$ ((a \ovee b)^\perp\ovee a) \ovee c = 
    ((a \ovee b)^\perp\ovee a) \ovee b  =  1$
and so by uniqueness of
orthocomplement~$c = ((a\ovee b)^\perp\ovee a)^\perp = b$,
    which is point~4.
By \sref{pcm-preorder},
    we only need to show that~$\leq$ is antisymmetric
    for point~5.
So assume~$a \leq b$ and~$b \leq a$
    for some~$a,b\in E$.
Pick~$c,d \in E$ with~$a \ovee c = b$ and~$b \ovee d = a$.
Then~$a = (a \ovee c) \ovee d = a \ovee (c \ovee d)$.
By cancellation~$c \ovee d = 0$.
So by positivity~$c = d= 0$.
Hence~$a = b$.
For point 6, assume~$a\leq b$.
Pick~$c \in E$ with~$a \ovee c = b$.
Clearly~$a \ovee a^\perp = 1 = b \ovee b^\perp = a \ovee c \ovee b^\perp$,
so by cancellation~$a^\perp = c \ovee b^\perp$,
    which is to say~$b^\perp \leq a^\perp$.
For point 7, assume~$a \leq b$ and~$b \perp c$.
Pick~$d$ with~$a \ovee d = b$.
By partial associativity and commutativity
    we have~$b \ovee c = (a \ovee d) \ovee c = (a \ovee c) \ovee d$,
    so~$a \perp c$ and~$a \ovee c \leq b \ovee c$.
For point 8, first assume~$a \perp b$.
Then~$a \ovee b \ovee (a \ovee b)^\perp = 1 = b \ovee b^\perp$.
So by cancellation~$b^\perp = a \ovee (a \ovee b)^\perp$,
hence~$a \leq b^\perp$.
For the converse, assume~$a \leq b^\perp$.
Then~$a \ovee c = b^\perp$ for some~$c$.
Hence~$a \ovee c \perp b$
    and so  by partial associativity and commutativity,
        we get~$a \perp b$, as desired.
    \qed
\end{point}
\end{point}
\end{parsec}%

\begin{parsec}{1760}%
\begin{point}{10}{Definition}%
Suppose~$E$ is an effect algebra.
Write~$\Define{b \ominus a}$
        \index{*ominus@$\ominus$}
for the (by cancellation) unique element (if it exists)
with~$a \ovee (b \ominus a) = b$.
\end{point}
\begin{point}{20}[exc-dposet]{Exercise*}%
Show that for any effect algebra~$E$, we have
\begin{enumerate}
    \item[(D1)]~$a \ominus b$ is defined if and only if~$b \leq a$;
    \item[(D2)]~$a \ominus b \leq a$ (if defined);
    \item[(D3)]~$a \ominus (a \ominus b) = b$ (if defined) \emph{and}
    \item[(D4)]~if $a \leq b \leq c $,
                then~$c \ominus b \leq c \ominus a$
                and~$(c \ominus a) \ominus (c \ominus b) = b \ominus a$.
\end{enumerate}
\spacingfix{}
\begin{point}{30}%
Let~$E$ be a poset with maximum element~$1$
    and partial binary operation~$\ominus$
    satisfying~(D1)--(D4).
Define~$a \ovee b = c \Leftrightarrow c \ominus b = a$
    and~$a^\perp = 1 \ominus a$.
Show that this turns~$E$ into an effect algebra
    with compatible order and~$\ominus$.
\end{point}
\begin{point}{40}{Remark}%
    Such a structure~$(E,\ominus,1)$ is called a \Define{difference-poset}
    \index{D-poset}
    (\Define{D-poset}) \cite{kopka1994d} and is, as we have just seen,
    an alternative way to axiomatize effect algebras.
\end{point}
\end{point}

\begin{point}{50}[exc-eamorphism]{Exercise}%
Show that for an effect algebra homomorphism~$f\colon E \to F$,
    we have
    \begin{enumerate}
        \item~\emph{(preserves zero)} $f(0) = 0$;
        \item~\emph{(order preserving)} if $a \leq b$, then~$f(a) \leq f(b)$;
        \item~if~$a\ominus b$ is defined,
            then $f(a \ominus b) = f(a) \ominus f(b)$ \emph{and}
        \item consequently~$f(a^\perp) = f(a)^\perp$.
    \end{enumerate}
\spacingfix{}
\end{point}
\end{parsec}

\begin{parsec}{1770}%
\begin{point}{10}%
For real numbers we have~$a + b = \min \{a,b\} + \max \{a, b\}$.
The following is a generalization to effect algebras.
\end{point}
\begin{point}{11}[ea-modularity-prop]{Proposition}%
Suppose~$E$ is an effect algebra.
If the infimum~$a \wedge b$
    exists for some~$a,b \in E$ with~$a \perp b$,
    then their supremum~$a \vee b$ exists as well and
\begin{equation*}
    a \ovee b \ = \ (a \wedge b) \,\ovee\, (a \vee b).
\end{equation*}
\spacingfix{}
\begin{point}{20}[modularity-lemma-proof]{Proof}%
(This result appeared in my master's thesis~\cite[prop.~15]{basmsc}.
    It turns out that it was already proven before
        (in D-poset form) \cite[prop.~1.8.2]{dvurecenskij2013new}.)
The result follows from this lemma:
    if~$(x \ominus c) \wedge (x \ominus d)$ exists,
    then~$c \vee d$ exists as well
    and~$(x \ominus c) \wedge (x \ominus d)
        = x \ominus (c \vee d)$.
Indeed:
\begin{equation*}
    a \wedge b \ = \ 
        ((a \ovee b) \ominus a) \wedge ((a \ovee b) \ominus b) \ = \ 
        (a \ovee b) \ominus (a \vee b)
\end{equation*}
    and so~$(a \wedge b) \ovee (a \vee b) = a \ovee b$, as desired.
\begin{point}{30}%
Now, we prove the lemma.
Assume~$(x \ominus c) \wedge (x \ominus d)$ exists.
Note~$x \ominus c = (x^\perp \ovee c)^\perp$.
As~$z \mapsto z^\perp$ is an order anti-isomorphism,
    we see~$(x^\perp \ovee c) \vee (x^\perp \ovee d)$
        exists and~$(x^\perp \ovee c) \vee (x^\perp \ovee d)
                                    = ((x \ominus c)\wedge (x \ominus d))^\perp$.
Clearly~$ (x^\perp \ovee c) \vee
        (x^\perp \ovee d) \geq x^\perp$.
    Write~$r := ((x^\perp \ovee c) \vee
        (x^\perp \ovee d)) \ominus x^\perp$.
We will show~$r = c \vee d$.
It is easy to see~$r \geq c$ and~$r \geq d$.
Assume~$s$ is any element with~$s \geq c$ and~$s \geq d$.
Note~$x^\perp \ovee s \geq x^\perp \ovee c$
and~$x^\perp \ovee s \geq x^\perp \ovee d$,
    hence~$x^\perp \ovee s \geq (x^\perp \ovee c) \vee (x^\perp \ovee d)$
    and so~$ s \geq  r$.
We have shown~$r = c \vee d$.
Consequently~$ (x \ominus c)\wedge (x \ominus d)
        = ((x^\perp \ovee c) \vee (x^\perp \ovee d))^\perp
        = (r \ovee x^\perp)^\perp = 
            x \ominus (c \vee d) $, as promised. \qed
\end{point}    
\end{point}
\end{point}
\begin{point}{40}[dfn-orthomodular-lattice]{Definition}%
An \Define{ortholattice}~$L$
    \index{ortholattice}
    is a bounded lattice with orthocomplement
    --- that is: it has a minimum~$0$,
    a maximum~$1$
    and a unary operation~$(\ )^\perp$
    satisfying
\begin{multicols}{2}
\begin{enumerate}
\item $a \wedge a^\perp = 0$;
\item $a \vee a^\perp = 1$;
\item $a \leq b  \ \Rightarrow\ b^\perp \leq a^\perp$ and
\item $a^{\perp\perp} = a$.
\end{enumerate}    
\end{multicols}
\noindent An ortholattice is \Define{orthomodular}\index{orthomodular!lattice}
 provided
\begin{equation*}
    a \ \leq \ b \quad \implies \quad a \vee (a^\perp \wedge b) \ =\  b.
\end{equation*}
See, for instance~\cite{kalmbach1983orthomodular,dvurecenskij2013new,birkhoff1936logic}.
\end{point}
\begin{point}{50}{Example}%
The lattice of projections in a von Neumann algebra
    is an orthomodular lattice with
    orthocomplement~$p^\perp \equiv 1 - p$.
\end{point}
\begin{point}{60}[orth-ea-is-orthomodular]{Proposition}%
An effect algebra~$E$
    that is an ortholattice
    is also orthomodular.
\begin{point}{70}{Proof}%
(For a different proof, see \cite[prop.~1.5.8]{dvurecenskij2013new}.)
Assume~$a \leq b$.
We have to show~$a \vee (a^\perp \wedge b) = b$.
Note~$a \wedge (a^\perp \wedge b) \leq a \wedge a^\perp = 0$
and so by \sref{ea-modularity-prop}
we have~$a \vee (a^\perp \wedge b) = a \ovee (a^\perp \wedge b)$.
Thus it is sufficient to prove~$a^\perp \wedge b = b \ominus a$.
Note~$a^\perp = (b \ominus a) \ovee b^\perp$
    and~$b = (b \ominus a) \ovee a$.
Similar to~\sref{modularity-lemma-proof}
    one can show that
\begin{equation*}
    (b \ominus a) \ovee (b^\perp \wedge a)
    \ = \
        ((b \ominus a) \ovee b^\perp) \wedge 
        ((b \ominus a) \ovee a) \ \equiv \ a^\perp \ovee b,
\end{equation*}
but~$b^\perp \wedge a \leq b^\perp \wedge b = 0$
    and so~$a^\perp \ovee b = b\ominus a$, as desired.
\qed
\end{point}
\end{point}
\end{parsec}

\subsection{Effect monoids}
\begin{parsec}{1780}%
\begin{point}{10}%
In \sref{dfn-mandso} we will see that the scalars of an effectus
form an effect monoid:
\end{point}
\begin{point}{20}[dfn-effect-monoid]{Definition}%
An \Define{effect monoid}~$M$ \cite{probdistrconv}
    \index{effect monoid}
    \index{*odot@$\odot$!effect monoid}
    is an effect algebra
    together with a binary operation~$\odot$
    such that for all~$a,b,c,d \in M$, we have
\begin{enumerate}
    \item \emph{(unit)}
    $1 \odot a = a \odot 1 = a$;
\item \emph{(associativity)}
    $(a \odot b) \odot c 
    =a \odot (b \odot c)$ \emph{and}
\item \emph{(distributivity)}
    if~$a \perp b$ and~$c \perp d$,
        then the following sum exists and
        furthermore~$(a \odot c) \ovee (b \odot c) \ovee
            (a \odot d) \ovee (b \odot d) = (a \ovee b) \odot (c \ovee d)$.
\end{enumerate}
(Phrased categorically:
    an effect monoid is a monoid in \textsf{EA}
    with the obvious tensor product
    that relates bimorphisms to morphisms,
    see~\cite{corefl,probdistrconv}.)
An effect monoid~$M$ is said to be \Define{commutative}
    \index{effect monoid!commutative}
    if we have~$a\odot b = b\odot a$ for all~$a,b \in M$.
A map~$f\colon M \to N$ between effect monoids
    is called an \Define{effect monoid homomorphism}
    \index{effect monoid!homomorphism}
    if it is an effect algebra homomorphism
    and furthermore~$f(a \odot b) = f(a) \odot f(b)$
    for all~$a,b \in M$.
\end{point}
\begin{point}{30}[eff-monoid-examples]{Examples}%
Effect monoids are less abundant than effect algebras.
\begin{enumerate}
\item The effect algebra~$[0,1]$
        is a commutative effect monoid with the usual product.
        (This is the only way to turn~$[0,1]$ into an effect monoid
                \cite[prop.~41]{basmsc}.)
\item We saw earlier that every Boolean algebra is an effect algebra
        with~$x \ovee y = x \vee y$ defined iff~$x \wedge y = 0$.
    The Boolean algebra is turned into an effect monoid
        with~$x \odot y \equiv x \wedge y$.
    In fact, every finite (not necessarily commutative) effect monoid
        is of this form \cite[prop.~40]{basmsc} and thus commutative.
\item
    In particular: the two-element Boolean algebra~$2 = {0,1}$
    from~\sref{eaexamples}
    is an effect monoid
        with~$x \odot y = x \wedge y$.
\item There is a non-commutative effect monoid
        based on the lexicographically ordered vector space~$\R^5$,
        see ~\cite[cor.~51]{basmsc}.
\item Let~$M$ be an effect monoid.
    Write~$M^{\mathsf{op}}$
        for the effect monoid~$M$ with the opposite multiplication
            --- that is: $a \odot_M b = b \odot_{M^\mathsf{op}} a$.
\end{enumerate}
\spacingfix{}
\end{point}
\begin{point}{31}[exc-emonzero]{Exercise}%
Show that~$a \odot 0 = a = 0 \odot a$
    for any~$a$ in an effect monoid~$M$.
\begin{point}{40}%
Later, in a tangent,
    we will need the following specific fact about effect monoids.
\end{point}
\end{point}
\begin{point}{50}[emond-lemma-for-conv]{Exercise}%
Assume~$M$ is an effect monoid
with~$a_1, \ldots, a_n, b_1, \ldots, b_n \in M$
such that~$\bigovee_i a_i = 1$
and~$\bigovee_i a_i \odot b_i = 1$.
Prove~$a_i \odot b_i = a_i$ for every~$1 \leq i \leq n$.
\end{point}
\end{parsec}

\subsection{Effect modules}
\begin{parsec}{1790}%
\begin{point}{10}%
We will see that in an effectus,
    the effect monoid of scalars will act on
    every effect algebra of predicates,
    turning them into effect modules (see \sref{dfn-mandso}):
\end{point}
\begin{point}{20}[dfn-effect-module]{Definition}%
Suppose~$M$ is an effect monoid.
An \Define{effect module} $E$ over~$M$ \cite{corefl}
    \index{effect module}
    is an effect algebra together with an operation
            $M \times E \to E$
            denoted by~$(\lambda, a) \mapsto \lambda \cdot a$
    such that for all~$a,b \in E$ and~$\lambda,\mu \in M$, we have
\begin{enumerate}
\item
    $(\lambda \odot \mu) \cdot a = \lambda \cdot (\mu \cdot a)$;
\item
    if~$a \perp b$,
     then~$\lambda \odot a \perp \lambda \odot b$
     and~$(\lambda \odot a) \ovee (\lambda\odot b) = \lambda \odot(a \ovee b)$;
\item
    if~$\lambda \perp \mu$,
     then~$\lambda \odot a \perp \mu \odot a$
     and~$(\lambda \odot a) \ovee (\mu \odot a) = (\lambda \ovee \mu) \odot a$
            \emph{and}
\item
    $1 \odot a = a$.
\end{enumerate}
(Categorically: an effect module over~$M$
    is an~$M$-action.)
An effect algebra homomorphism~$f\colon E \to F$
    between effect modules over~$M$
    is an $M$-\Define{effect module homomorphism}
    \index{effect module!homomorphism}
    provided~$\lambda \cdot f(a) = f(\lambda \cdot a)$
    for all~$\lambda \in M$ and~$a \in E$.
Write~$\Define{\mathsf{EMod}_M}$
    \index{$\mathsf{EMod}_M$}
    for the category of effect modules over~$M$
    with effect module homomorphisms between them.
\end{point}
\begin{point}{30}{Examples}%
There are many effect modules.
\begin{enumerate}
\item
Every effect algebra is an effect module over
    the two-element effect monoid~$2$.
    (In fact $\mathsf{EA} \cong \mathsf{EMod}_{2}$.)
The only effect module up-to-isomorphism over the one-element effect monoid~$1$
    is the one-element effect algebra~$1$ itself.
\item
Effect modules over~$[0,1]$
    are the same thing as
        \Define{convex effect algebras},
        see \cite{gudder1998representation,gudder1999convex}.
        \index{effect algebra!convex}
If~$V$ is an ordered real vector space with order unit~$u$,
        then~$[0,u]$ is an effect module over~$[0,1]$.
In fact, every effect module over~$[0,1]$
    is of this form \cite{gudder1998representation}.
See also \cite[thm.~3]{jacobs2016expectation}
    for the stronger categorical equivalence.
\end{enumerate}
\spacingfix{}
\end{point}
\end{parsec}

\section{Effectuses}
An effectus comes in two guises:
    axiomatizing either a category of total maps
    or a category of partial maps.
We will start off with the total form as it has the simplest axioms.
Later we will prefer to work with the partial form.
\begin{parsec}{1800}%
\begin{point}{10}[dfn-effectus]{Definition}%
A category $C$ is said to be an \Define{effectus in total form}
    \index{effectus!in total form}
    \cite{effintro,newdirections,statesofconvexsets}
    if
\begin{enumerate}
\item $C$ has finite coproducts (hence an initial object~$0$)
        and a final object~$1$;
    \item all diagrams of the following
        \index{kappa@$\kappa_i$, coprojection}
        \index{$[\,\cdot\,,\,\cdot\,]$!cotupling}
        \index{*bang@$"!$, initial/final map}
        form\footnote{%
        We write~$\Define{\kappa_i}$ for coproduct coprojections;
            square brackets~$\Define{[f,g]}$ for coproduct cotupling;
            $\Define{h+k} = [\kappa_1 \after h, \kappa_2 \after k]$
            and~$\Define{!}$ for the unique maps associated to
            either the final object~$1$ or initial object~$0$.} are pullbacks
\begin{equation}\label{pullbacks}
    \vcenter{\xymatrix{
        X+Y \ar[r]^{\id+{!}} \ar[d]_{!+\id} & X+1\ar[d]^{!+\id} \\
    1+Y\ar[r]_{\id+!} & 1+1
}}
    \qquad\qquad
    \vcenter{\xymatrix{
    X \ar[r]^{!} \ar[d]_{\kappa_1} & 1 \ar[d]^{\kappa_1} \\
    X+Y\ar[r]_{{!}+{!}} & 1+1
}}
\end{equation}
\item\label{eff-joint-monicity} and the following two arrows are jointly monic.
    \begin{equation*}
        \xymatrix@C+2pc  {
            1+1+1  \ar@/^/[r]^{[\kappa_1,\kappa_2,\kappa_2]}
                    \ar@/_/[r]_{[\kappa_2,\kappa_1,\kappa_2]} & 1+1
        }
    \end{equation*}
\end{enumerate}
An arrow~$f\colon X \to Y+1$ is called a \Define{partial map} and written
    \index{partial map}
    $f\colon \Define{X \pto Y}$.
    \index{*pto@$\pto$}
\begin{point}{20}%
One with an interest in physics might think of the objects
    of an effectus as physical systems and its arrows as
    the physical operations between them.
    The final object~$1$ is the physical system with a single state.
The coproduct~$X+Y$ is the system that can be prepared as either~$X$ or as~$Y$.
An arrow~$1 \to X$ corresponds to a physical preparation of the system~$X$
    and an arrow~$X \to 1+1$ is a yes--no measurement.
\end{point}
\begin{point}{30}%
In this sense an effectus can be seen as a generalized
        probabilistic theory
        --- not unlike the operational-probabilistic theories (OPT)
        of Chiribella et al \cite{chiribella2010probabilistic}.
There are differences:
for instance, OPTs are equipped with a parallel composition of systems,
     effectuses are not and OPTs always have the unit interval as scalars,
     effectuses might not.
In \cite{tull2016operational} Tull compares effectuses and OPTs.
There are other approaches to categorical probabilistic theories;
 see, for instance~\cite{wilceshortcut,gogioso2017categorical}.
\end{point}
\begin{point}{40}%
Studying programming languages, one would do better thinking of
    the objects of an effectus as data types and
    its arrows as the allowed operations between them (semantics of programs).
The final object~$1$ is the unit data type.
The coproduct~$X + Y$ is the union data type of~$X$ and~$Y$.
An arrow~$1 \to X$ is a value of~$X$
    and an arrow~$X \to 1+1$ is a predicate on~$X$.
\end{point}
\begin{point}{50}[effectus-vn]%
Our main example of an effectus in total form
    is the category~$\op{\vN}$\index{vN@$\vN$}
    of von Neumann algebras with completely positive normal unital
    maps in the opposite direction.
(To see~$\op\vN$ is an effectus in total form,
    adapt the proof of \sref{emod-effectus}.)
The partial maps correspond to
    ncp-maps~$f$ with~$f(1) \leq 1$.
(Equivalently: contractive ncp-maps.)
Some other examples appear later on in
    \sref{emod-effectus}, \sref{exc-dm-effectus},
    \sref{aconvm-is-effectus} and~\sref{effexamplesintro}.
For a comprehensive list of examples, see~\cite{effintro}.
\end{point}
\begin{point}{60}%
Let~$C$ be an effectus in total form.
Given two arrows~$f\colon X \to Y+1$
and~$g \colon Y \to Z+1$ (i.e.~partial maps~$X \pto Y$ and~$Y \pto Z$)
    their composition as partial maps
    is defined as~$g \hafter f \equiv  [g, \kappa_2] \after f$.
Write~$\Define{\Par C}$ for the \Define{category of partial maps},
    \index{ParC@$\Par C$}
    which has the same objects
    as~$C$, but as arrows~$X \to Y$ in~$\Par C$
    we take arrows of the form~$X \to Y+1$ in~$C$,
    which we compose using~$\hafter$
    and with identity on~$X$ in~$\Par C$
    given by~$\kappa_1 \colon X \to X+1$.\footnote{In categorical
            parlance: $\Par C$ is the Kleisli category of
            the monad~$(\ )+1\colon C \to C$.}
The category~$\Par C$ is not an effectus in total form --- instead it
    is an effectus \emph{in partial form}.
\end{point}
\end{point}
\begin{point}{70}[effectus-in-partial-form]{Definition}%
A category~$C$ is called an \Define{effectus in partial form}
    \index{effectus!in partial form}
    \cite{effintro,kentapartial} if
\begin{enumerate}
\item
    $C$ is a finPAC\footnote{Finitarily partially additive category.} \cite{kentapartial} (cf.~\cite{arbib}) --- that is:
    \begin{enumerate}
        \item 
            $C$ has finite coproducts $(+,0)$;
        \item $C$ is $\mathsf{PCM}$-enriched --- that is:
            \begin{enumerate}
            \item
            every homset $C(X,Y)$ has a partial binary operation~$\ovee$
                    and distinguished element~$0 \in C(X,Y)$
                    that turns~$C(X,Y)$ into a partial commutative monoid,
                    see \sref{dfn-pcm};
            \item
            if $f \perp g$ then both
                $(h \after f) \perp (h \after g)$ and
                $(f \after k) \perp (g \after k)$ \emph{and}
            \begin{equation*}
                 h \after (f \ovee g) = 
                (h \after f) \ovee (h \ovee g)\qquad
                (f \ovee g) \after k = 
                (f \after k) \ovee (g \after k)
            \end{equation*}
                for any~$f,g \colon X \to Y$,
                $h\colon Y \to Y'$
                    and~$k \colon X'\to X$ \emph{and}
            \item
                \emph{(zero)}
                $0 \after f = 0$ and $f \after 0 = 0$ for
                    any~$f\colon X \to Y$ \emph{and}
            \end{enumerate}
        \item
            \emph{(compatible sum)}
            for any~$b\colon X \to Y + Y$ we have
            $\pproj_1 \after b \perp \pproj_2 \after b$,
            where~$\Define{\pproj_i}\colon Y + Y \to Y$
            \index{*pproj@$\pproj_i$}
            are \Define{partial projectors}\footnote{%
                Later, in~\sref{coprod-prod},
                    we will use the more slightly more general
                    ${\pproj_1}\colon X_1 + X_2 \to X_1$ and
                    ${\pproj_2}\colon X_1 + X_2 \to X_2$
                    defined by~$\pproj_1 = [\id,0]$
                        and~$\pproj_2 = [0,\id]$.} defined
            \index{partial projectors}
            by~$\pproj_1 \equiv [\id, 0]$
            and~$\pproj_2 \equiv [0, \id]$ \emph{and}
        \item
            \emph{(untying)} if~$f\perp g$,
            then~$\kappa_1\after f \perp \kappa_2 \after g$,
            where~$\kappa_1$ and~$\kappa_2$ are coprojections
            on the same coproduct
            \emph{and}
    \end{enumerate}
    \item it has effects --- that is: there is a distinguished object~$I$
            such that
    \begin{enumerate}
        \item for each object~$X$, the PCM~$C(X,I) \equiv \Define{\Pred X}$
            \index{predX@$\Pred X$}
            is an effect algebra, see \sref{dfn-ea}
            --- in particular~$C(X,I)$
                has a maximum element~$1$;
        \item if~$1 \after f \perp 1 \after g$,
            then~$f \perp g$ \emph{and}
        \item if~$1 \after f = 0$, then~$f = 0$.
    \end{enumerate}
\end{enumerate}
In this setting,
a map~$f\colon X \to Y$ is called \Define{total} if~$1 \after f = 1$.
    \index{total map}
\begin{point}{80}{Remark}%
The untying and zero axioms are redundant: they follows from the others.
We include them, as they are part of the definition of finPAC.
\end{point}
\begin{point}{90}%
At first glance an effectus in partial form seems
    to have a much richer structure than an effectus in total form.
This is not the case ---
    effectuses in total and partial form are two views
    on the same thing \cite{kentapartial,effintro}:
\end{point}
\end{point}
\begin{point}{100}[cho-thm]{Theorem (Cho)}%
Let~$C$ be an effectus in total form
and~$D$ be an effectus in partial form.
\begin{enumerate}
\item
The category~$\Par C$ is an effectus in partial form
        with~$I = 1$.
\item
The total maps of~$D$ form an effectus in total form~$\Tot D$.
\item
    Nothing is lost:
        $\Par (\Tot D) \cong D$ and~$\Tot (\Par C )\cong C$.
\end{enumerate}
\spacingfix{}
\begin{point}{110}%
(The result can be rephrased categorically as
    a 2-equivalence of the 2-category
    of effectuses in partial form and
    the 2-category of effectuses in total form.
    For the details, see~\cite[\S5]{kentapartial}.)
To prove the Theorem
    (in \sref{proof-cho-thm}),
    we need some preparation.
\end{point}
\end{point}
\end{parsec}
\subsection{From partial to total}
\begin{parsec}{1810}%
\begin{point}{10}%
We will first show that the subcategory of
    total maps of an effectus in partial form
    is an effectus in total form.
This proof and especially the demonstration
    that the squares in \eqref{pullbacks} are pullbacks,
    will elucidate the axioms of an effectus in total form
    and will make the proof in the opposite direction more palatable.
\end{point}
\begin{point}{20}[coproj-total]{Lemma}%
In an effectus in partial form,
    coprojections are total.
\begin{point}{30}{Proof}%
    (For the original and different proof see \cite[lem.~4.7(5)]{kentapartial}.)
Let~$\kappa_1 \colon X \to X+Y$ be any coprojection.
We have to show~$1 \after \kappa_1 = 1$.
Note that the first~$1$ denotes the maximum of~$\Pred X+Y$
    and the second the maximum of~$\Pred X$.
    Hence~$1_X = 1_X \after \id_X = (1_X \after [\id_X, 0]) \after \kappa_1
                    \leq 1_{X+Y} \after \kappa_1 \leq 1_X$.
Indeed~$1 \after \kappa_1 = 1$. \qed
\end{point}
\end{point}
\begin{point}{40}[cotupl-pcm]{Proposition}%
In an effectus in partial form,
the cotupling bijection~$(f,g) \mapsto [f,g]$
    is a PCM-isomorphism \cite{effintro }---
    that is:
\begin{enumerate}
\item
    $[f,g] \perp [f',g']$ if and only if
        $f \perp f'$ and~$g \perp g'$;
\item
    if~$[f,g] \perp [f', g']$,
    then~$[f,g] \ovee [f',g'] = [f \ovee f', g\ovee g']$ \emph{and}
\item
    $[0,0] = 0$.
\end{enumerate}
Furthermore~$[1,1]=1$ for maps into~$I$,
so the cotupling map is an effect algebra
isomorphism~$\Pred (X) \times \Pred (Y) \cong \Pred (X+Y) $.
\begin{point}{50}{Proof}%
First we show~$[h,l] = [h,0] \ovee [0,l]$.
By the compatible sum axiom
\begin{equation*}
    [h, 0] \ = \ 
    \pproj_1 \after (h + l)
    \ \perp\  \pproj_2 \after (h + l)
    \ = \ [0, l].
\end{equation*}
By $\mathsf{PCM}$-enrichment
$([h, 0] \ovee [0, l]) \after \kappa_1
         =  ([h, 0] \after \kappa_1) \ovee
        ([0, l] \after \kappa_1)  = h$.
Similarly~$([h, 0] \ovee [0, l]) \after \kappa_2 = l$.
Thus indeed~$[h,l] = [h,0] \ovee [0,l]$.

Assume~$[f,g] \perp [f',g']$.
    By $\mathsf{PCM}$-enrichment
    we have~$f = [f,g] \after \kappa_1 \perp [f',g'] \after \kappa_1 = f'$.
Similarly~$g \perp g'$.
Conversely, assume~$f \perp f'$ and~$g \perp g'$.
Again, by $\mathsf{PCM}$-enrichment~$
[f,0] = f  \after \pproj_1 \perp f' \after \pproj_1 = [f', 0]$
    and~$[f\ovee f', 0 ] = [f,0] \ovee [f', 0]$.
Similarly~$[0,g \ovee g'] = [0,g] \ovee [0,g']$.
Putting it all together:
\begin{align*}
[f, g] \ovee [f', g']
&\ = \ 
 [f, 0] \ovee [0,g] \ovee [f', 0] \ovee [0, g'] \\
 &\ = \ 
 [f\ovee f', 0]\ovee [0, g\ovee g'] \\
 &\ = \ 
 [f\ovee f', g \ovee g'].
\end{align*}
To show cotupling is a PCM-isomorphism,
    it only remains to be shown that~$[0,0]=0$.
As~$0 \after \kappa_1 = 0$
    and~$0 \after \kappa_2 = 0$,
    we indeed have~$0 = [0,0]$.
Similarly by \sref{coproj-total} we have~$1 \after \kappa_1 = 1$
    and~$1 \after \kappa_2 = 1$,
    so~$1 = [1,1]$. \qed
\end{point}
\begin{point}{60}%
The coproduct in an effectus in partial form
is almost a (bi)product:
\end{point}
\end{point}
\begin{point}{70}[coprod-prod]{Proposition}%
In an effectus in partial form, we have a bijective correspondence
\begin{prooftree}
\AxiomC{$h\colon Z \to X+Y$}
\doubleLine
\UnaryInfC{$f\colon Z \to X \quad g\colon Z \to Y \quad 1 \after f \perp 1 \after g$}
\end{prooftree}
as follows \cite[lem.~4.8]{kentapartial}:
\begin{quote}
for every~$f\colon Z \to X$
    and~$g\colon Z \to Y$
    with~$1 \after f \perp 1 \after g$,
    there is a unique map~$\Define{\langle f, g\rangle} \colon Z \to X +Y$
    \index{*langle@$\langle\,\cdot\,,\cdot\,\rangle$!in an effectus}
    such that~$\pproj_1 \after \langle f, g \rangle = f$
    and~$\pproj_2 \after \langle f, g \rangle = g$,
    where~$\pproj_1 = [\id,0]$ and~$\pproj_2=[0,\id]$.
In fact~$\langle f, g\rangle = (\kappa_1 \after f) \ovee (\kappa_2 \after g)$.
\end{quote}
Conversely~$h = \langle \pproj_1 \after h, \pproj_2 \after h \rangle$
    with~$1 \after \pproj_1 \after h \perp 1 \after \pproj_2 \after h$
    for~$h \colon Z \to X+Y$.
\begin{point}{80}{Proof}%
Let~$f\colon Z \to X$ and~$g\colon Z \to Y$ be given
        such that~$1\after f \perp 1\after g$.
By \sref{coproj-total},
we have~$1 \after \kappa_1 \after f = 1 \after f \perp 1 \after g =
        1 \after \kappa_2 \after g$.
Thus~$\kappa_1 \after f \perp \kappa_2 \after g$.
Define~$\langle f, g\rangle = (\kappa_1 \after f) \ovee (\kappa_2 \after g)$.
    By the $\mathsf{PCM}$-enrichedness, we have
\begin{align*}
\pproj_1 \after \langle f, g\rangle 
&\  =\  [\id, 0] \after ((\kappa_1 \after f) \ovee (\kappa_2 \after g)) \\
 &\  =\  ([\id, 0] \after \kappa_1 \after f) \ovee 
    ([\id, 0] \after \kappa_2 \after g) \\
    & \ =\  f \ovee (0 \after g) \ = \ f.
\end{align*}
Similarly~$\pproj_2 \after \langle f, g \rangle = g$.
To show uniqueness, assume~$f = \pproj_1 \after h$
    and~$g = \pproj_2 \after h$ for some
    $h\colon Z \to X+Y$.
Note~$\kappa_1 \after \pproj_1 = [\kappa_1, 0]$
and~$\kappa_2 \after \pproj_2 = [0, \kappa_2]$,
so by~\sref{cotupl-pcm} we have~$(\kappa_1 \after \pproj_1)
    \ovee (\kappa_2 \after \pproj_2) = [\kappa_1, \kappa_2] = \id$,
    hence
\begin{align*}
    \langle f, g\rangle & \ = \ 
    \langle \pproj_1 \after h, \pproj_2 \after h \rangle \\
    & \ = \ (\kappa_1 \after\pproj_1 \after h) 
    \ovee (\kappa_2 \after\pproj_2 \after h)  \\
    & \ = \ ((\kappa_1 \after\pproj_1 )
    \ovee (\kappa_2 \after\pproj_2 ))  \after h \\
    & \ = \   \id \after h \ = \  h,
\end{align*}
which demonstrates uniqueness.

Finally, let~$h\colon Z\to X+Y$ be any map.
We must show that~$h = \langle \pproj_1 \after h, \pproj_2 \after h\rangle$.
Note~$1 \after \pproj_1 = [1,0]$
    and~$1 \after \pproj_2 = [0,1]$.
    So by $\mathsf{PCM}$-enrichment~$1 \after \pproj_1 \after h \perp 1 \after \pproj_2 \after h$ and so~$\pproj_1 \after h \perp \pproj_2 \after h$.
By the previous
$\langle \pproj_1 \after h, \pproj_2 \after h\rangle
= (\kappa_1 \after \pproj_1 \after h) \ovee
 (\kappa_2 \after \pproj_2 \after h) 
 = h $ as desired.\qed
\end{point}
\end{point}

\begin{point}{90}[eff-prod-rules]{Exercise}%
Show that in an effectus in partial form, we have
\begin{enumerate}
    \item~$[a,b] \after \langle f, g \rangle = (a \after f) \ovee (b \after g)$;
    \item~$1 \after \langle f, g\rangle = (1 \after f) \ovee(1 \after g)$;
    \item~$(k+l) \after \langle f, g\rangle
        = \langle k \after f, l \after g \rangle$ \emph{and}
    \item~$\langle f, g\rangle \after k = \langle f \after k,
                                g \after k \rangle$
\end{enumerate}
        assuming~$1 \after f \perp 1 \after g$. \cite{effintro}
\begin{point}{100}{Remark}%
In an effectus in partial form,
    there is a straightforward generalization
    of the bijective correspondence~\sref{coprod-prod} to
\begin{prooftree}
\AxiomC{$h\colon Z \to  X_1 + \cdots + X_n$}
\doubleLine
\UnaryInfC{$f_1\colon Z \to X_1 \ \cdots \  f_n\colon Z \to X_n \quad \bigperp_n  1 \after f_n$}
\end{prooftree}
with~$f_i = \pproj_i \after h$
    and~$h = \langle f_1, \ldots, f_n\rangle \equiv
    (\kappa_1 \after f_1 )\ovee \cdots \ovee (\kappa_n \after f_n)$,
    where~$\pproj_i$ is defined in the obvious way.
The expected generalizations  of the rules in \sref{eff-prod-rules} also hold.
\end{point}
\end{point}

\begin{point}{110}[eff-partial-to-total]{Theorem}%
Let~$C$ be an effectus in partial form.
The category of total maps of~$C$ is an effectus in total form.
    \cite[thm.~4.10]{kentapartial}
\begin{point}{120}{Proof}%
The total maps indeed form a subcategory: $1 \after \id = 1$
    and $1 \after f \after g = 1 \after g = 1$ for composable total~$f,g$.
In \sref{coproj-total} we saw coprojections are total.
By \sref{cotupl-pcm}
    we have~$1 \after [f,g] = [1 \after f, 1 \after g] = [1,1] = 1$
    for total~$f\colon X \to Z$ and~$g\colon Y \to Z$,
    so~$\Tot C$ has binary coproducts.
The unique map~$!\colon 0 \to X$ must be total
    as~$1 \after ! = 1$ is the unique map~$0 \to I$,
    so~$\Tot C$ has initial object~$0$, hence all finite coproducts.
\end{point}
\begin{point}{130}[one-m-is-id]%
To show~$I$ is the final object of~$\Tot C$,
    we need~$\id_I = 1$.
As~$C(I,I)$ is an effect algebra~$1 = \id \ovee \id^\perp$
    for some~$\id^\perp$.
    So by $\mathsf{PCM}$-enrichment~$1 = 1 \after 1 = (1 \after \id) \ovee (1 \after\id^\perp) = 
1 \ovee (1 \after \id^\perp)$.
By the zero--one axiom~$1 \after \id^\perp = 0$.
So~$\id^\perp = 0$ and indeed~$\id = 1$.
To show~$I$ is final in~$\Tot C$, pick any object~$X$ in~$\Tot C$.
We claim~$1 \colon X \to I$ is the unique total map.
Indeed, by the previous $1 = \id_I \after 1 = 1 \after 1$, so~$1$ is total
and if~$h\colon X \to I$ is total, then~$1=1 \after h = \id_I \after h = h$.
\end{point}
\begin{point}{140}%
To show the square on the left of \eqref{pullbacks} is a pullback
    in~$\Tot C$,
let~$f\colon Z \to X+I$ 
and~$g\colon Z \to I+Y$ be total maps
    with~$(\id+1) \after g = (1+\id) \after f$.
By~\sref{coprod-prod}, $f = \langle \alpha, a \rangle$
    and~$g = \langle b, \beta \rangle$
    for some maps~$\alpha\colon Z \to X$, $\beta\colon Z \to Y$
        and~$a,b\colon Z \to I$.
\begin{equation*}
\xymatrix@C+1pc@R-1pc{ 
Z \ar@/^1pc/[rrd]^{\langle\alpha,a \rangle \equiv f}
    \ar@/_1pc/[rdd]_{g \equiv \langle b, \beta\rangle}
    \ar@{.>}[rd]|{\langle \alpha, \beta\rangle}
    \\
    &  X+Y \ar[r]^{\id+1} \ar[d]_{1+\id} & X+I\ar[d]^{1+\id} \\
    &I+Y\ar[r]_{\id+1} & I+I
    }
\end{equation*}
By~\sref{eff-prod-rules},
    we have~$1 = 1 \after f = 1 \after \langle \alpha, a\rangle
                = (1 \after \alpha) \ovee a$
        and so~$a^\perp = 1\after \alpha$.
Similarly~$b^\perp = 1 \after \beta$.
Again, using~\sref{eff-prod-rules},
    we see that
\begin{equation*}
    \langle b, 1 \after \beta \rangle \ =\ 
    (\id + 1) \after \langle b, \beta \rangle \ =\ 
    (\id + 1) \after g \ =\ 
    (1 + \id) \after f
               \  = \ \langle 1 \after \alpha, a\rangle,
\end{equation*}
so~$1 \after \beta = a = (1 \after \alpha)^\perp$,
hence~$1 \after \beta \perp 1 \after \alpha$,
so~$\langle \alpha,\beta\rangle \colon Z \to X+Y$ exists
    and is total as~$1 \after \langle \alpha,\beta\rangle = 
        (1 \after \alpha) \ovee (1 \after \beta) = 1$.
We compute~$(\id + 1) \after \langle\alpha, \beta\rangle = \langle \alpha,
    1 \after \beta\rangle = \langle \alpha, a\rangle = f$.
    Similarly~$(1 + \id) \after \langle \alpha, \beta\rangle = g$.
    Assume~$h\colon Z \to X+Y$ is any map with~$(1 +\id) \after h
        = g$ and~$(\id+1) \after h = f$.
Say~$h = \langle h_1, h_2 \rangle$.
Then~$\langle \alpha, a\rangle = f= (\id + 1)\after h = \langle h_1, 1 \after h_2 \rangle$, so~$\alpha = h_1$. Similarly~$\beta=h_2$.
Thus~$h = \langle \alpha, \beta\rangle$,
    which shows our square is indeed a pullback.
\end{point}
\begin{point}{150}%
To show the square on the right of \eqref{pullbacks} is a pullback
    in~$\Tot C$,
assume (using \sref{coprod-prod})
$\langle \alpha, \beta \rangle\colon Z \to X+Y$
is some total map
with~$(1+1) \after \langle\alpha,\beta\rangle = \kappa_1 \after 1$.
    \begin{equation*}
\xymatrix@C+1pc@R-1pc{ 
    Z \ar@{.>}@/^1pc/[rrd]^{1}
    \ar@/_1pc/[rdd]_{\langle \alpha,\beta\rangle}
    \ar@{.>}[rd]|{\alpha}\\
    &X \ar[r]^{1} \ar[d]_{\kappa_1} & I \ar[d]^{\kappa_1} \\
    &X+Y\ar[r]_{1+1} & I+I
    }
    \end{equation*}
With \sref{eff-prod-rules},
we see~$\langle 1 \after \alpha, 1 \after \beta \rangle
        = (1 + 1) \after \langle \alpha, \beta \rangle
        = \kappa_1 \after 1
        = \langle 1, 0 \rangle$.
So~$\alpha$ is total and~$\beta=0$.
Hence~$\langle \alpha, \beta\rangle = \langle \alpha, 0\rangle
        = \kappa_1 \after \alpha$ as desired.
\end{point}
\begin{point}{160}%
    Finally, to show~$m_1\equiv [\kappa_1,\kappa_2,\kappa_2],
        m_2 \equiv [\kappa_2,\kappa_1,\kappa_2]\colon
    I+I+I \to I+I$ are jointly monic,
    let~$f_1 \equiv \langle a_1,b_1,c_1\rangle,f_2\equiv \langle a_2,b_2,c_2\rangle\colon X\to I+I+I$ be any total maps
    with~$m_1 \after f_1 = m_1 \after f_2$
    and~$m_2 \after f_1 = m_2 \after f_2$.
Then
\begin{equation*}
    a_1 \ = \ [\id, 0, 0] \after f_1
        \ = \ \pproj_1 \after m_1 \after f_1
        \ = \ \pproj_1 \after m_1 \after f_2
        \ = \ a_2.
\end{equation*}
and similarly from the equality involving~$m_2$, we get~$b_1 = b_2$.
As~$f_1$ is total, we have~$1 = 1 \after f_1=  a_1 \ovee b_1 \ovee c_1$,
    so~$c_1 = (a_1 \ovee b_1)^\perp$.
With the same reasoning~$c_2 = (a_2 \ovee b_2)^\perp$.
Thus~$c_1 = c_2$ and so~$f_1 = f_2$, as desired. \qed
\end{point}
\end{point}
\end{parsec}

\subsection{From total to partial}
\begin{parsec}{1820}%
\begin{point}{10}%
Let~$C$ be an effectus in total form.
In this section we will show that~$\Par C$ is an effectus in partial form.
Before we get to work, it is helpful to discuss the axioms
    of an effectus in total form, now we have some experience
    with an effectus in partial form.
\begin{enumerate}
\item
In \sref{coprod-prod} we saw that the coproduct in an effectus
    (in partial form) is almost a biproduct.
This structure is hidden (for the most part)
    in the left pullback square of
    \eqref{pullbacks}, which
    allows the formation of~$\langle \alpha, \beta \rangle$ given
    partial maps~$\alpha,\beta$.
\item
The right pullback square of \eqref{pullbacks}
    is used to extract a total map in~$C$
    from a partial map~$f$ in~$\Par C$ that is total (i.e.~$1 \hafter f = 1$).
        See~\sref{pardp}.
\item
The joint-monicity
    of~$[\kappa_1,\kappa_2,\kappa_2]$
    and~$[\kappa_2,\kappa_1,\kappa_2]$
    will imply the joint monicity of~$\pproj_1$ and $\pproj_2$,
    which is required for uniqueness of the partial sum of partial maps.
        See~\sref{pproj-joint-monicity}.
\end{enumerate}
To prove the theorem, we need to study pullbacks:
    first in any category,
    then in an effectus~$C$ in total form
    and finally in~$\Par C$.
Our proof is different from the original proof of Cho
    by using several more general facts about pullbacks,
    which we will prove next.
\end{point}
\end{parsec}

\begin{parsec}{1830}%
\begin{point}{10}%
    We start with two classic facts about pullbacks.
\end{point}
\begin{point}{20}[exc-jointly-monic-pullback]{Exercise}%
Show that if we have any pullback square --- say
\begin{equation*}
    \vcenter{\vbox{\xymatrix@R-1pc{
        {P\ } \pullback \ar[r]^{m_1} \ar[d]_{m_2}
        & B \ar[d]^f
                \\ {A\ } \ar[r]_{g}
    & X}}} \text{,}
\end{equation*}
then~$m_1$ and~$m_2$ are jointly monic.
\end{point}
\begin{point}{30}[pullback-lemma]{Exercise}%
Prove the \Define{pullback lemma} --- that is:
    \index{pullback lemma}
    if we have a commuting diagram
\begin{equation*}
    \vcenter{\vbox{\xymatrix@R-1pc{
                A \ar[r]^f \ar[d]_k
                & B \ar[r]^g \ar[d]_l
                & C \ar[d]^m
                \\ X \ar[r]_{f'}
                & Y \ar[r]_{g'}
                & Z
    }}} \text{,}
\end{equation*}
then we have the following two implications.
\begin{enumerate}
\item
If the left and right inner squares are pullbacks,
    then so is the outer square.
\item
If the outer square is a pullback
    and~$l$ and~$g$ are jointly monic,
    then the left inner square is a pullback.
\end{enumerate}
\spacingfix{}
\begin{point}{40}{Remark}%
In the literature (e.g.~\cite[III.5 exc.~8]{maclane}),
    one often finds a weaker second implication,
    which assumes that the right square is a pullback.
\end{point}
\end{point}
\end{parsec}

\begin{parsec}{1840}%
\begin{point}{10}%
It is well-known that monos are stable under pullbacks ---
that is:
\begin{equation*}
    \text{if} \quad
    \vcenter{\vbox{\xymatrix@R-1pc{
                {P\ } \pullback \ar[r]^n \ar[d]_g
        & X \ar[d]^f
                \\ {A\ } \ar@{>->}[r]_{m}
    & B}}} ,
    \quad \text{then} \quad
    \vcenter{\vbox{\xymatrix@R-1pc{
                {P\ } \pullback \ar@{>->}[r]^n \ar[d]_g
        & X \ar[d]^f
                \\ {A\ } \ar@{>->}[r]_{m}
    & B}}},
\end{equation*}
    where we used~$\xymatrix@C-1pc{\ar@{>->}[r]&}$ to denote a monic map.
We will need an analogous result for jointly monic maps.
\end{point}
\begin{point}{20}[joint-monicity-stable]{Lemma}%
If the pairs~$(m_1,m_2)$, $(n_1,g_1)$, $(n_2,g_2)$ and~$(h_1,h_2)$
in the following commuting diagram
are jointly monic (for instance: if they span pullback squares),
then~$(n_1 \after h_1, n_2 \after h_2)$ is jointly monic as well.
\begin{equation*}
    \xymatrix@C+1pc{
        P \ar[r]^{h_1} \ar[d]_{h_2}
        &P_1\ar[r]^{n_1} \ar[d]^{g_1}
        &X_1 \ar[d]^{f_1}
        \\ P_2\ar[r]_{g_2} \ar[d]_{n_2}
        & A \ar[r]_{m_1} \ar[d]^{m_2}
        & B_1 
        \\ X_2 \ar[r]_{f_2}
        & B_2
    }
\end{equation*}
\spacingfix{}
\begin{point}{30}{Proof}%
To show joint monicity of~$n_1 \after h_1$ and~$n_2 \after h_2$,
assume~$\alpha_1,\alpha_2  \colon Z \to P$
    are maps with~$
        n_1 \after h_1 \after \alpha_1
        = n_1 \after h_1 \after \alpha_2$
    and~$n_2 \after h_2 \after \alpha_1
        = n_2 \after h_2 \after \alpha_2$.
To start
\begin{equation*}
    m_1 \after g_1 \after h_1 \after \alpha_2
    \ =\  f_1 \after n_1 \after h_1 \after \alpha_2 
    \ =\  f_1 \after n_1 \after h_1 \after \alpha_1 
            \ =\  m_1 \after g_1 \after h_1 \after \alpha_1.
\end{equation*}
Reasoning on the other side of the diagram, we find
\begin{align*}
    m_2 \after g_1 \after h_1 \after \alpha_2
    &\ =\  m_2 \after g_2 \after h_2 \after \alpha_2 \\ 
    &\ =\  f_2 \after n_2 \after h_2 \after \alpha_2 \\ 
    &\ =\  f_2 \after n_2 \after h_2 \after \alpha_1 \\ 
    &\ =\  m_2 \after g_2 \after h_2 \after \alpha_1 \\ 
    &\ =\  m_2 \after g_1 \after h_1 \after \alpha_1.
\end{align*}
So by joint monicity of~$m_1$ and~$m_2$,
    we conclude~$g_1 \after h_1 \after \alpha_2
                = g_1 \after h_1 \after \alpha_1$.
So by the joint monicity of~$n_1$ and~$g_1$
    (and by assumption~$n_1 \after h_1 \after \alpha_2
    = n_1 \after h_1 \after \alpha_1$),
    we conclude~$h_1 \after \alpha_1 = h_1 \after \alpha_2$.
Reasoning in the same way mirrored over the diagonal,
    we find~$h_2 \after \alpha_1 = h_2 \after \alpha_2$.
Thus, by joint monicity of~$h_1$ and~$h_2$, we conclude that~$\alpha_1 = \alpha_2$, as desired. \qed
\end{point}
\end{point}
\end{parsec}

\begin{parsec}{1850}%
\begin{point}{10}[tot-pullbacks]{Proposition}%
In an effectus in total form,
    any square of one of the two following forms is a pullback.
\begin{equation*}
    \xymatrix{
        X+A \pullback \ar[r]^{\id+f} \ar[d]_{g + \id}
        & X+B \ar[d]^{g+\id}
        \\ Y+A \ar[r]_{\id+f}
        & Y+B
    } \qquad
    \xymatrix{
        X \pullback \ar[d]_{\kappa_1} \ar[r]^{f}
        & Y \ar[d]^{\kappa_1}
        \\
        X+A \ar[r]_{f+g}
        & Y+B
    }
\end{equation*}
\spacingfix{}
\begin{point}{20}{Proof}%
To start with the left square, consider the following commuting diagram.
\begin{equation*}
    \xymatrix{
        X+A \ar@{}[rd]|{(1)} \ar[r]^{\id+f} \ar[d]_{g + \id}
        & X+B \ar@{}[rd]|{(2)} \ar[d]^{g+\id} \ar[r]^{\id+!}
        & X+1 \ar[d]^{g+\id}
        \\ Y+A \ar@{}[rd]|{(3)} \ar[r]_{\id+f} \ar[d]_{!+\id}
        & Y+B \ar@{}[rd]|{(4)} \ar[r]_{\id+!} \ar[d]^{!+\id}
        & Y+1 \ar[d]^{!+\id}
        \\ 1+A \ar[r]_{\id+f}
        & 1+B \ar[r]_{\id+!}
        & 1+1
    }
\end{equation*}
We want to show (1) is a pullback.
The inner square~(4) and right rectangle~(2,4) are pullbacks by axiom.
Thus the inner square~(2) is also a pullback
    by the pullback lemma, see \sref{pullback-lemma}.
With the same reasoning, we see square (3) is a pullback.
Thus by the pullback lemma, it is sufficient to show that
    left rectangle (1,3) is a pullback.
The left rectangle is indeed a pullback
as both the outer square (1,2,3,4)
    and the right rectangle (2,4) are pullbacks.

For the right square, we consider the following diagram.
\begin{equation*}
    \xymatrix{
         X \ar[r]^f \ar[d]_{\kappa_1}
            %\ar@{}[rd]|{(1)}
        & Y \ar[r]^{!} \ar[d]^{\kappa_1}
            %\ar@{}[rd]|{(2)}
        &1 \ar[d]^{\kappa_1}
        \\ X+A \ar[r]_{f+g}
        & Y+B \ar[r]_{!+!}
        & 1+1
    }
\end{equation*}
The inner right and outer square are pullbacks by axiom.
So the inner left square is also a pullback by the pullback lemma, as desired.
    \qed
\end{point}
\end{point}
\end{parsec}

\begin{parsec}{1860}%
\begin{point}{10}{Definition}%
Assume~$C$ is an effectus in total form.
For a map~$f\colon X \to Y$,
    write~$\Define{\hat{f}} = \kappa_1 \after f \colon X \to Y+1$.
    \index{*hat@$\hat{f}$!in an effectus}
    (This is the Kleisli embedding~$C \to \Par C$.) \cite{kentapartial}
\end{point}
\begin{point}{11}
    (Recall that we use~$f\colon X \pto Y$
        to denote a map~$f\colon X \to Y+1$,
            see~\sref{dfn-effectus}.)
\end{point}
\begin{point}{20}[par-c-coprod]{Exercise}%
Assume~$C$ is an effectus in total form.
Show that if
\begin{equation*}
\kappa_1 \colon X \to X+Y \leftarrow Y \colon \kappa_2
\end{equation*}
    is a coproduct in~$C$, then~
\begin{equation*}
    \hat\kappa_1 \colon X \pto X+Y \pfrom Y \colon \hat\kappa_2
\end{equation*}
    is a coproduct in~$\Par C$.
Show~$0$ is also the initial object of~$\Par C$.
\begin{point}{30}{Beware}%
While cotupling in~$C$ and~$\Par C$ coincide:
$[f,g]_{C} = [f,g]_{\Par C}$
for any (properly typed) partial maps~$f,g$
and $\widehat{f+_C g} = \hat{f} +_{\Par C} \hat{g}$,
    however, in general
    \begin{equation*}
    [\kappa_1 \after f, \kappa_2 \after g] \ =\  f+_C g \ \neq\  f +_{\Par C} g
        \ =\  [[\kappa_1,\kappa_3] \after f,
            [\kappa_2, \kappa_3] \after g].
    \end{equation*}
\end{point}
\end{point}
\spacingfix{}
\begin{point}{40}[par-pullbacks]{Proposition}%
Let~$C$ be an effectus in total form.
Squares of the following form are pullbacks in~$\Par C$.
\begin{equation*}
    \xymatrix{
        X+A \pullback \ar@^{>}[r]^{\id+\hat{f}} \ar@^{>}[d]_{\hat{g} + \id}
        & X+B \ar@^{>}[d]^{\hat{g}+\id}
        \\ Y+A \ar@^{>}[r]_{\id+\hat{f}}
        & Y+B
    } \qquad
    \xymatrix{
        A+X \pullback \ar@^{>}[r]^{\hat{f}+\id} \ar@^{>}[d]_{\pproj_1}
        & B+X \ar@^{>}[d]^{\pproj_1}
        \\ A \ar@^{>}[r]_{\hat{f}}
        & B
    }
\end{equation*}
\spacingfix{}
\begin{point}{50}{Proof}%
(A different proof for the pullback on the right
    is given in \cite[lem.~4.1(5)]{kentapartial}.)
It is sufficient to show that the following squares
    are pullbacks in~$C$.
\begin{equation*}
    \xymatrix{
        X+(A+1) \ar[r]^{\id + (f+\id)} \ar[d]_{g+ \id}
        & X+(B+1) \ar[d]^{g+\id}
        \\ Y+(A+1) \ar[r]_{\id+(f +\id)}
        & Y+(B+1)
    } \qquad
    \xymatrix{
        A+(X+1) \ar[r]^{f+\id} \ar[d]_{\id+!}
        & B+(X+1)  \ar[d]^{\id+!}
        \\ A+1 \ar[r]_{f+\id}
        & B+1
    }
\end{equation*}
These are indeed pullbacks by~\sref{tot-pullbacks}. \qed
\end{point}
\end{point}
\begin{point}{60}[zero-and-one-parc]{Definition}%
Assume~$C$ is an effectus in total form.
For arbitrary objects~$X,Y$ in~$C$,
    write~$\Define{0} \equiv \kappa_2 \after !\colon X \to Y+1$
    and~$\Define{1} \equiv \kappa_1 \after ! \equiv \hat! \colon X \to 1+1$.
    \cite{kentapartial}
\end{point}
\begin{point}{70}[toteff-zero]{Exercise}%
    Show~$0$ is a zero object in~$\Par C$
        with unique zero map
        as in \sref{zero-and-one-parc}.
\end{point}
\begin{point}{80}[pardp]{Proposition}%
Assume~$C$ is an effectus in total form
    and~$f$ is a map in~$\Par C$.
\begin{enumerate}
\item
If~$1 \hafter f = 1$,
    then~$f = \hat{g}$ for a unique~$g$ in~$C$.
\item
If~$1 \hafter f = 0$, then~$f = 0$.
\end{enumerate}
\spacingfix{}
\begin{point}{90}{Proof}%
    (This proof is essentially
        the same as \cite[lem.~4.7(1) \& prop.~4.6]{kentapartial}.)
Both points follow from the right pullback square of~\eqref{pullbacks}
 as follows.
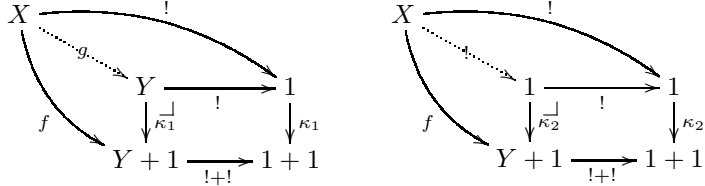
\begin{equation*}
    \xymatrix@R-.8pc{
    X \ar@{.>}[rd]|g
    \ar@/^1pc/[rrd]^{!}
        \ar@/_1pc/[rdd]_f
        \\& Y \pullback
        \ar[r]_{!}
        \ar[d]^{\kappa_1}
    & 1
        \ar[d]^{\kappa_1}
    \\& Y+1
        \ar[r]_{!+!}
&1+1
}
\qquad
    \xymatrix@R-.8pc{
        X \ar@{.>}[rd]|{!}
    \ar@/^1pc/[rrd]^{!}
        \ar@/_1pc/[rdd]_f
        \\& 1 \pullback
        \ar[r]_{!}
        \ar[d]^{\kappa_2}
    & 1
        \ar[d]^{\kappa_2}
    \\& Y+1
        \ar[r]_{!+!}
&1+1
}
\end{equation*}
If~$1 \hafter f = 1$, then~$(!+!)\after f = \kappa_1 \after !$
    and so there is a unique~$g$ with~$f = \kappa_1 \after g$
    as shown above on the left.
For the other point, if~$1 \hafter f = 0$,
    then~$(!+!)\after f = \kappa_2 \after !$
    and so~$f = \kappa_2 \after !$ as shown above on the right. \qed
\end{point}
\end{point}
\begin{point}{100}[pproj-joint-monicity]{Proposition}%
    Assume~$C$ is an effectus in total form.
    The partial projectors~$\pproj_1 \equiv [\id,0]$
        and~$\pproj_2 \equiv [0, \id]$
        are jointly monic in~$\Par C$.
        \cite[lem.~4.1(4)]{kentapartial}
\begin{point}{110}{Proof}%
    (This is a new proof.)
Consider the following diagram in~$\Par C$.
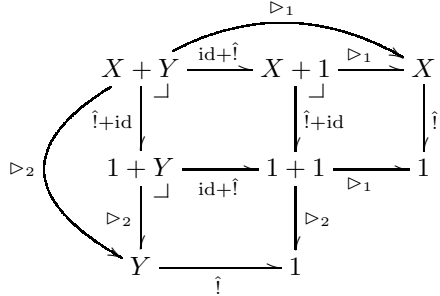
\begin{equation*}
\xymatrix {
    X+Y \pullback
    \ar@/^1.5pc/[rr]^{\pproj_1}
    \ar@/_3pc/[dd]_{\pproj_2}
        \ar@^{>}[r]^{\id+\hat!}
        \ar@^{>}[d]_{\hat! + \id}
& X+1 \pullback
        \ar@^{>}[r]^{\pproj_1}
        \ar@^{>}[d]^{\hat! + \id}
& X
        \ar@^{>}[d]^{\hat!}
\\ 1+Y \pullback
        \ar@^{>}[r]_{\id+\hat!}
        \ar@^{>}[d]_{\pproj_2}
& 1+1
        \ar@^{>}[r]_{\pproj_1}
        \ar@^{>}[d]^{\pproj_2}
& 1
\\ Y
        \ar@^{>}[r]_{\hat!}
& 1
}
\end{equation*}
This diagram commutes and
    each of the inner squares is a pullback by~\sref{par-pullbacks}.
The maps~$\pproj_1,\pproj_2\colon 1 + 1 \pto 1$
    are jointly monic,
    which is a reformulation
    of the joint monicity axiom of an effectus in total form.
Thus by \sref{joint-monicity-stable}
    the outer~$\pproj_1$ and~$\pproj_2$ are jointly monic. \qed
\end{point}
\end{point}
\end{parsec}

\begin{parsec}{1870}%
\begin{point}{10}[eff-total-to-partial]{Theorem}%
If~$C$ is an effectus in total form,
    then~$\Par C$ is an effectus in partial form. \cite[thm.~4.2]{kentapartial}
\begin{point}{20}{Proof}%
In \sref{par-c-coprod}, we saw that~$\Par C$ has finite coproducts.
\begin{point}{30}{$\mathsf{PCM}$-enrichment, I}%
Assume~$f,g \colon X \pto Y$.
We define that~$f \perp g$
    iff there is a \Define{bound}~$b\colon X \pto Y+Y$
    \index{bound}
    with~$\pproj_1 \hafter b = f$
    and~$\pproj_2 \hafter b = g$.
    (By \sref{pproj-joint-monicity}
    this~$b$ is unique if it exists.)
In this case, define~$f \ovee g = \nabla \hafter b$,
where $\Define{\nabla} \equiv [\id,\id]$.
    \index{*nabla@$\nabla$}
We will show~$\ovee$ with~$0$
    as defined in \sref{zero-and-one-parc}
    forms a PCM.
Partial commutativity is obvious.
The map~$\kappa_1 \hafter f\colon X \pto Y+Y$
    is a bound for~$f \perp 0$ and
    so~$f \ovee 0 = \nabla \hafter \kappa_1 \hafter f = f$,
    which shows~$0$
    is indeed a zero for~$\ovee$.
Only partial associativity remains.
Assume~$f \perp g$ via bound~$b$
and~$f \ovee g \perp h$ via bound~$c$.
Note~$\nabla \hafter b = f \ovee g = \pproj_1 \hafter c$
and thus by the right pullback square of \sref{par-pullbacks} (see diagram below),
    there is a unique~$d\colon X \pto (Y+Y)+Y$
    with~$\pproj_1 \hafter d = b$ and~$(\nabla +\id)\hafter d = c$.
\begin{equation*}
    \xymatrix@R-.8pc{
    X \ar@^{.>}[rd]|d
        \ar@/^1pc/[rrd]^c
        \ar@/_1pc/[rdd]_b
        \\& (Y+Y)+Y \pullback
        \ar@^{>}[r]_{\nabla+\id}
        \ar@^{>}[d]^{\pproj_1}
    & Y+Y
        \ar@^{>}[d]^{\pproj_1}
    \\& Y+Y
        \ar@^{>}[r]_{\nabla}
&Y
}
\end{equation*}
The map $(\pproj_2 + \id) \hafter d$
    is a bound for~$g \perp h$,
    indeed:
\begin{align*}
    \pproj_1 \hafter (\pproj_2 + \id) \hafter d
    & \ =\  \pproj_2 \hafter \pproj_1 \hafter d 
    & \pproj_2 \hafter (\pproj_2 + \id) \hafter d
    & \ =\  \pproj_2 \hafter d
    \\
    &\ =\  \pproj_2 \hafter b
    && \ =\  \pproj_2 \hafter (\nabla + \id) \hafter d
    \\
    & \ =\  g
    && \ = \ \pproj_2 \hafter c\\
    &&& \ = \ h.
\end{align*}
We need one more:
    $[\id,\kappa_2]\hafter d$ is a bound for~$f \perp g \ovee h$, indeed
\begin{align*}
    \pproj_1 \hafter [\id,\kappa_2] \hafter d
        & \ = \ \pproj_1 \hafter \pproj_1 \hafter d
    & \pproj_2 \hafter [\id,\kappa_2]\hafter d & \ =\ \nabla \hafter (\pproj_2 + \id) \hafter d \\
        & \ = \ \pproj_1 \hafter b
        && \ =\  g \ovee h \\
        & \ = \ f.
\end{align*}
Finally,
we compute
\begin{equation*}
f \ovee (g \ovee h)
\ = \ \nabla \hafter [\id, \kappa_2] \hafter d
\ = \ \nabla \hafter (\nabla+ \id) \hafter d
\ = \ \nabla \hafter c
\ = \ (f \ovee g) \ovee h,
\end{equation*}
which shows the partial associativity.
Homsets of~$\Par C$ are indeed PCMs.
\end{point}
\begin{point}{40}{$\mathsf{PCM}$-enrichment, II}%
Assume~$f \perp g$ with bound~$b$.
It is easy to see~$b \hafter k$
    is a bound for~$f \hafter k \perp g \hafter k$
    and consequently~$(f \hafter k) \ovee (g \hafter k)
                = \nabla \hafter b \hafter k = 
                (f \ovee g) \hafter k$.
For the other side:
$(h + h) \hafter b$ is a bound for~$h \hafter f \perp h \hafter g$,
    indeed
\begin{align*}
    \pproj_1 \hafter (h+h) \hafter b
    &\ = \ 
    h\hafter \pproj_1 \hafter b \ =\  h \hafter f \\
    \pproj_2 \hafter (h+h) \hafter b
    &\ = \ 
    h\hafter \pproj_2 \hafter b\  =\  h \hafter g.
\end{align*}
And so~$(h \hafter f) \ovee (h \hafter g)
            = \nabla \hafter (h+h) \hafter b
            = h \hafter \nabla \hafter b = h \hafter (f \ovee g)$.
\end{point}
Finally, the zero axiom holds by~\sref{toteff-zero}.
\begin{point}{50}{FinPAC}%
We just saw~$\Par C$ is $\mathsf{PCM}$-enriched.
We already know~$\Par C$ has coproducts.
The compatible sum axiom holds by definition of~$\perp$.
To show~$\Par C$ is a finPAC,
    only the untying axiom remains to be proven.
If~$f \perp g$ with bound~$b$,
then~$(\kappa_1 + \kappa_2) \hafter b$ is a bound for~$\kappa_1 \hafter f
    \perp \kappa_2 \hafter g$,
    which proves the untying axiom.
\end{point}
\begin{point}{60}{Effect algebra of predicates}%
Pick any predicate~$p \colon X \to 1+1$.
We define~$p^\perp = [\kappa_2,\kappa_1]\after p$.
($p^\perp$ is~$p$ with swapped outcomes.)
We compute
\begin{align*}
    \pproj_1 \hafter \hat{p}
    &\ = \ [\id,\kappa_2]\after\kappa_1 \after p
        \  = \ p \\
        \pproj_2 \hafter \hat{p}
    &\ = \ [[\kappa_2, \kappa_1],\kappa_2]\after\kappa_1 \after p
        \  = \ [\kappa_2,\kappa_1]\after p \\
    \nabla \hafter \hat{p}
    &\ = \ [[\kappa_1,\kappa_1],\kappa_2]\after\kappa_1 \after p
        \  = \ \kappa_1 \after ! \ = \ 1.
\end{align*}
So~$\hat{p}$ is a bound for~$p \perp p^\perp$
and~$p \ovee p^\perp = 1$.
To show~$p^\perp$ is the unique orthocomplement,
    assume~$p \ovee q = 1$ for some~$q$ via a bound~$b\colon
        X \pto 1+1$.
Note~$1 \hafter b =\nabla \hafter b = 1$
    and so by \sref{pardp}
    we know~$b = \hat{c}$ for some~$c\colon X \to 1+1$.
    And so~$p = \pproj_1 \hafter b = [\id, \kappa_2]\after \kappa_1 \after c
                    = c$. Hence
\begin{equation*}
q \ =\  \pproj_2 \hafter b \ =\  [[\kappa_2, \kappa_1],\kappa_2] \after
    \kappa_1\after c \ =\  [\kappa_2,\kappa_1]\after p \ =\  p^\perp.
\end{equation*}
To show~$\Pred X$ is an effect algebra,
it only remains to be proven that the zero--one axiom holds.
So assume~$1 \perp p$ via some bound~$b$. We must show~$p=0$.
By the following instance of the right pullback square of~\eqref{pullbacks},
    we see~$b = \kappa_1 \after \kappa_1 \after !$.
\begin{equation*}
    \xymatrix@R-.8pc{
        X \ar@{.>}[rd]|{!}
    \ar@/^1pc/[rrd]^{!}
        \ar@/_1pc/[rdd]_b
        \\& 1 \pullback
        \ar@{=}[r]
        \ar[d]^{\kappa_1\after\kappa_1}
    & 1
        \ar[d]^{\kappa_1}
        \\& (1+1)+1
        \ar[r]_{[\id,\kappa_2]}
&1+1
}
\end{equation*}
Hence~$p = \pproj_2 \hafter b = [[\kappa_2,\kappa_1],\kappa_1]\after
                        \kappa_1 \after \kappa_1 \after ! = \kappa_2 \after ! = 0$.
\end{point}
\begin{point}{70}{Final axioms}%
In \sref{pardp} we proved that~$f = 0$ whenever~$1 \hafter f = 0$.
It only remains to be shown that~$f \perp g$
    provided~$1 \hafter f \perp 1 \hafter g$.
Let~$b$ be a bound for~$1 \hafter f \perp 1 \hafter g$.
We apply the right pullback square of \sref{par-pullbacks}
    twice  in succession as follows
    (using~$\pproj_2 \hafter c
                = \pproj_2 \hafter (\hat{!} + \id) \hafter c
                = \pproj_2 \hafter b = 1\hafter g$ for the right one.)
\begin{equation*}
    \xymatrix@R-.8pc{
        X \ar@^{.>}[rd]|{c}
    \ar@/^1pc/[rrd]^{b}
        \ar@/_1pc/[rdd]_f
        \\& Y+1 \pullback
        \ar@^{>}[r]^{\hat! + \id}
        \ar@^{>}[d]^{\pproj_1}
    & 1+1
        \ar@^{>}[d]^{\pproj_1}
        \\& Y
        \ar@^{>}[r]_{\hat!}
&1
} \qquad
    \xymatrix@R-.8pc{
        X \ar@^{.>}[rd]|{d}
    \ar@/^1pc/[rrd]^{c}
        \ar@/_1pc/[rdd]_g
        \\& Y+Y \pullback
        \ar@^{>}[r]^{\id + \hat!}
        \ar@^{>}[d]^{\pproj_2}
    & Y+1
        \ar@^{>}[d]^{\pproj_2}
        \\& Y
        \ar@^{>}[r]_{\hat!}
&1
}
\end{equation*}
Clearly~$\pproj_2 \hafter d = g$
    and~$
    \pproj_1 \hafter d
     =  \pproj_1 \hafter (\id + \hat! )\hafter d
     =  \pproj_1 \hafter c  =  f$, so~$d$ is a bound for~$f \perp g$,
     as desired. \qed
\end{point}
\end{point}
\end{point}
\end{parsec}
\begin{parsec}{1880}%
\begin{point}{10}%
We are ready to show the equivalence of effectuses in partial and total form.
\begin{point}{20}[proof-cho-thm]{Proof of \sref{cho-thm}}%
Let~$C$ and~$D$ be effectuses in total and respectively partial form.
The first point, that~$\Par C$ is an effectus in partial form,
    is shown in \sref{eff-total-to-partial}.
The second point, that~$\Tot D$ is an effectus in total form,
    is shown in \sref{eff-partial-to-total}.
It remains to be shown that nothing is lost.
\begin{point}{30}{$\Par \Tot D \cong D$}%
Let~$D$ be an effectus in partial form.
A map~$f\colon X \to Y$ in~$\Par \Tot D$
    is by definition a map~$f\colon X \to Y+1$
        in~$D$ with~$1 \after f= 1$.
Let~$P \colon \Par \Tot D \to D$
    denote the identity-on-objects map~$f \mapsto \pproj_1 \after f$.
This is a functor: clearly~$P\hat\id = P\kappa_1 = \id $ and
    if~$g\colon Y \to Z$ in~$\Par\Tot D$, then
\begin{align*}
    \pproj_1 \after (g \hafter f) & \ = \ 
    \pproj_1 \after [g, \kappa_2] \after f \\
    & \ = \ [\pproj_1 \after g, \pproj_1 \after \kappa_2] \after f \\
    & \ = \ [\pproj_1 \after g,  0] \after 
                \langle \pproj_1 \after f, \pproj_2 \after f\rangle \\
                & \ = \ \pproj_1 \after g  \after \pproj_1 \after f.
\end{align*}
The functor~$P$ is an isomorphism with inverse~$P' f \equiv \langle f, (1 \after f)^\perp\rangle$.
Indeed, we have
    $P P' f = \pproj_1 \after \langle f, (1 \after f)^\perp \rangle = f $
and
\begin{equation*}
    P' P f  \ =\  \langle \pproj_1 \after f, (1 \after \pproj_1 \after f)^\perp\rangle \ =\  (\kappa_1 \after \pproj_1 \after f) \ovee (\kappa_2 \after \pproj_1 \after f)^\perp\ =\  f.
\end{equation*}
\spacingfix{}
\begin{point}{40}{$\Tot \Par C \cong C$}%
Assume~$C$ is an effectus in total form.
Recall~$1 \hafter \hat{g} = 1$
and~$\smash{\widehat{f \after g}} = \hat{f} \hafter \hat{g}$,
    so~$Q\colon C \to \Tot \Par C$
    given by~$Q g = \hat{g}$ is an identity-on-objects functor.
It is an isomorphism by the first part of \sref{pardp}:
    for every~$f$ in~$\Tot \Par C$,
    there is a unique~$g$ in~$C$
    with~$f = \hat{g}$. \qed
\end{point}
\end{point}
\end{point}
\end{point}
\end{parsec}

\begin{parsec}{1890}%
\begin{point}{10}[distinction-part-tot-eff]{Exercise}%
In this exercise we will
    distinguish effectuses in total and partial form.
\begin{enumerate}
\item
To start, show that the initial object of an effectus in total form is strict
    --- that is: show that any map into~$0$ is an isomorphism.
    (Hint: use the right pullback square of~\sref{tot-pullbacks}
        with~$X=Y=0$.)
\item
Conclude that if~$C$ is both an effectus in total and partial form,
then every object in~$C$ is isomorphic to~$0$.
\end{enumerate}
\spacingfix{}
\end{point}
\begin{point}{20}[eff-convention]{Convention}%
For the remainder of this text, we will work with effectuses in partial form
    and write~$1$ instead of~$I$.
To spare ink, a phrase like ``let $C$ be an \Define{effectus}''
    \index{effectus}
should be read as~``let $C$ be an effectus in partial form''
and simply~``$C$ is an effectus'' means either ``$C$ is an effectus in total form''
    or ``$C$ is an effectus in partial form''.
(In non-trivial cases, this
    is unambiguous, as seen in \sref{distinction-part-tot-eff}.)
As before, when we took the effectus in partial form as base category,
    we will denote partial maps by~$f\colon X \to Y$
    (instead of~$f\colon X \pto Y$)
        and their composition as~$f \after g$
        (instead of~$f \hafter g$).
\end{point}
\end{parsec}

\begin{parsec}{1891}%
\begin{point}{10}[effexamplesintro]{Examples}%
Our main example of an effectus (in total form)
    is the category~$\op\vN$ of von Neumann algebras with
    ncpu-maps between them in the opposite direction,
        see~\sref{effectus-vn}.
Its full subcategory~$\op\CvN$ of commutative von Neumann algebras
    is an effectus as well. \index{CvN@$\CvN$}
\begin{point}{20}%
In~\cite{effintro} several other examples of effectuses are
    discussed in detail:
\begin{enumerate}
\item
    The category~$\op\OUS$ of order unit spaces
        (i.e.~ordered real vector spaces with distinguished order unit)
        with positive unit-preserving linear maps in the opposite
        direction is an effectus (in total form). \index{OUS@$\OUS$}
\item
    A related example is the category~$\op\OUG$ of order unit groups
        (i.e.~ordered abelian groups with distinguished order unit)
        with positive
        unit-preserving homomorphisms in the opposite direction,
        which is also an effectus. \index{OUG@$\OUG$}
\item
    Every extensive category~\cite{carboni1993introduction}
        with final object is an effectus in total form.
    Examples of extensive categories with final object
            include
    \begin{enumerate}
        \item $\SET$, the category
            of sets with maps between them; \index{Set@$\SET$}
        \item $\op\CRng$, the category
            of commutative unit rings with unit-preserving
            homomorphisms between them
            in the opposite direction; \index{CRng@$\CRng$}
        \item $\bCH$, the category of
            compact Hausdorff spaces with continuous maps
            between them; \index{CH@$\bCH$} \emph{and}
    \end{enumerate}
\end{enumerate}
\end{point}
\spacingfix{}
\begin{point}{30}%
    Together with van de Wetering, we recently investigated~\cite{eja}
    the category~$\op\EJA$ of Euclidean Jordan algebras
    with positive unit-preserving linear maps
    in the opposite direction, which is also an effectus (in total form).
\end{point}
\end{point}
\end{parsec}

\subsection{Predicates, states and scalars}
\begin{parsec}{1900}%
\begin{point}{10}%
We turn to the internal structures of an effectus.~\cite{effintro}
\end{point}
\begin{point}{20}[dfn-mandso]{Definition}%
Let~$C$ be an effectus (in partial form).
\begin{enumerate}
\item
As alluded to in the definition of an effectus in partial form,
a \Define{predicate} on~$X$ is a map~$X \to 1$.
        \index{predicate}
The set of predicates~$\Pred X$ on~$X$ is (by definition of an effectus
    in partial form) an effect algebra.
\item
A \Define{scalar} is a predicate on~$1$; that is, a map~$1 \to 1$.
        \index{scalar}
We write~$\Define{\Scal C} \equiv \Define{M} \equiv \Pred 1$ for the set of scalars.
        \index{$M$}
        \index{ScalC@$\Scal C$}
We multiply two scalars~$\lambda,\mu \colon 1 \to 1$
    simply by composing~$\lambda \odot \mu \equiv \lambda \after \mu$.
As~$1_M = \id_M$ (see \sref{one-m-is-id}) and~$C$ is $\mathsf{PCM}$-enriched,
    the set of scalars~$M$ with~$\odot$
    is an effect monoid,
    see \sref{dfn-effect-monoid}.
\item
A \Define{real effectus} is an effectus
        \index{effectus!real}
where the set of scalars~$M$ is isomorphic (as effect monoid) to~$[0,1]$.
\item
For a scalar~$\lambda\colon 1\to 1$ and a predicate~$p\colon X \to 1$,
we write~$\Define{\lambda \cdot p} \equiv \lambda \after p$.
        \index{*lambdacdot@$\lambda \cdot p$}
Again due to $\mathsf{PCM}$-enrichment and~$1_M = \id_M$,
this scalar multiplication
turns~$\Pred X$ into an~$M$-effect module, see \sref{dfn-effect-module}.
\item
For~$f\colon X \to Y$ in~$C$,
        define~$\Define{\Pred(f)} \colon \Pred Y \to \Pred X$
        \index{Pred@$\Pred(f)$}
    by~$\Pred(f)(p) = p \after f$.
It is easy to see
$\Pred(f)$ is an~$M$-effect module homomorphism
if~$f$ is total
and that, in fact,~$\Pred\colon \Tot C \to \mathsf{EMod}^{\mathsf{op}}_M$
is a functor.
\item
A \Define{substate} on~$X$ is a map~$\omega\colon 1 \to X$.
A \Define{state} is a total substate.
        \index{state}
        \index{state!sub-}
We denote the set of states on~$X$ by~$\Define{\Stat X}$.
        \index{StatX@$\Stat X$}
\item
An effectus~$C$ has \Define{separating predicates}
    \index{separating!predicates}
        if~$\Pred X$ (as a set of maps~$X \to 1$)
        is jointly monic for every~$X$ in~$C$.
    Similarly, an effectus has \Define{separating states}
        \index{separating!states}
        if~$\Stat X$ (as a set of maps~$1 \to X$)
        is jointly epic for every~$X$ in~$C$.
\end{enumerate}
\spacingfix{}
\end{point}
\begin{point}{30}{Examples}%
    The category~$\op\vN$ (see \sref{effectus-vn})
    is a real effectus ($M \cong [0,1]$)
    with separating states and predicates.~\cite{effintro}
The predicates on a von Neumann algebra~$\scrA$
    correspond to the set of effects~$[0,1]_\scrA$
    and~$\Stat \scrA$ can be identified with the convex set of normal states.
\begin{point}{40}%
We briefly cover the examples from~\cite{effintro}
    mentioned in~\sref{effexamplesintro}.
\begin{enumerate}
\item
The category~$\op\OUS$ (see~\sref{effexamplesintro})
    is a real effectus with separating predicates,
    but without seperating states.
    (The states on a single order unit space are separating
        if and only if the order unit space is archimedean.)
A predicate on an order unit space~$X$
    corresponds  to a point~$x \in X$ with~$0 \leq x \leq 1$,
    where~$1$ is the distinguished order unit.
The states are exactly what are called states for order unit spaces
    in the literature.
\item
The category~$\op\OUG$  (see~\sref{effexamplesintro})
    has the two-element effect monoid~$2$ as scalars
        and seperating predicates.
The predicates of an order unit group~$G$ correspond
    to the elements~$x \in G$ with~$0 \leq x \leq 1$,
        where~$1$ is the distinguished order unit.
States on~$G$ correspond to unit-preserving positive
    homomorphisms~$G \to \Z$, which are not separating.
\item
Any extensive category with final object
    (such as~$\SET$, $\op\CRng$ and $\bCH$)
    has as scalars the two-element effect monoid~$2$.
\begin{enumerate}
\item
In~$\SET$ the predicates on a set~$X$ correspond to subsets~$U \subseteq X$
        and are thus separating. States on~$X$ correspond t
        elements~$x \in X$, which are also separating.
\item
In~$\op\CRng$ the predicates on~$R$ correspond to idempotents of~$R$,
    which are not separating.
        States correspond to unit-preserving
            homomorphisms~$R \to \Z$, which are not separating either.
\item
In~$\bCH$ the predicates on~$X$ correspond to clopen subsets~$U \subseteq X$,
    which are not separating.
States correspond on~$X$ correspond to its points~$x \in X$,
    which are separating.
\end{enumerate}
\end{enumerate}
\end{point}
\spacingfix{}
\begin{point}{50}%
The category~$\op\EJA$ is a real effectus with separating states
    and predicates. The predicates on an Euclidean Jordan algebra (EJA)~$E$
    correspond to its effects; i.e.~the elements~$0 \leq a \leq 1$.
    The states are exactly what are usually considered
    states for EJAs.
\end{point}
\end{point}
\end{parsec}
\begin{parsec}{1910}%
\begin{point}{10}%
Does every effect monoid~$M$ occur as the effect monoid of scalars of some
    effectus?  Does every $M$-effect module occur as
        effect module of predicates?
We will show they do, in the following strong sense.
\end{point}
\begin{point}{20}[emod-effectus]{Theorem}%
    Let~$M$ be any effect monoid.
    The category~$\EMod^{\mathsf{op}}_M$ is an effectus in total form
        with scalars~$M$ and separating predicates.~\cite{effintro}
In fact: every effectus~$C$ in total form
    with scalars~$M$ and separating predicates
    is equivalent to a subcategory of~$\EMod^{\mathsf{op}}_M$.
\begin{point}{30}{Proof}%
To show~$\EMod^{\mathsf{op}}_M$ is an effectus,
    we have to show the dual  axioms for~$\EMod_M$.
\begin{point}{40}{Finite products}%
For every~$M$-effect module~$E$,
    there is a unique $M$-effect module map~$!\colon M \to E$,
    which is given by~$\lambda \mapsto \lambda \cdot 1$.
So~$M$ is the initial object of~$\EMod_M$.
The one-element $M$-effect module~$\{0=1\}$ is the final object.

If~$E$ and $F$ are~$M$-effect modules,
    then the set of pairs~$E \times F$
    with componentwise operations
    is again an~$M$-effect module.
The effect module~$E \times F$
    with obvious projections~$E \xleftarrow{\pi_1} E \times F \xrightarrow{\pi_2} F$
    is the categorical product of~$E$ and~$F$ in~$\EMod_M$
    (as is the case with effect algebras, see \sref{ea-product}).
\end{point}
\begin{point}{50}{Pushout diagrams}%
To show the pushout diagrams corresponding to \eqref{pullbacks} hold,
assume we are given $M$-effect modules~$E,F,G$
    with module maps~$\alpha,\beta,\delta$
    that make the  outer squares of the following diagrams commute.
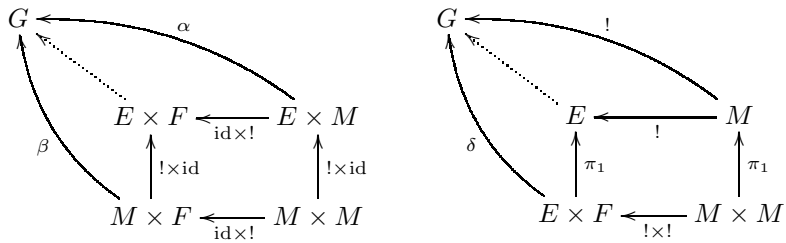
\begin{equation*}
    \xymatrix{
    G\\
    & E \times F \ar@{.>}[lu]
    & E \times M \ar[l]^{\id\times !}
    \ar@/_1pc/[llu]_{\alpha}
    \\& M \times F \ar[u]_{! \times \id}
    \ar@/^1pc/[luu]^{\beta}
    & M \times M \ar[l]^{\id\times!} \ar[u]_{! \times \id}
    } \qquad
    \xymatrix{
    G\\
    & E \ar@{.>}[lu]
    & M \ar[l]^{!}
    \ar@/_1pc/[llu]_{!}
    \\& E \times F \ar[u]_{\pi_1}
    \ar@/^1pc/[luu]^{\delta}
    & M \times M \ar[l]^{!\times!} \ar[u]_{\pi_1}
    }
\end{equation*}
    We have to show that there are unique dashed arrows (as shown)
    that make these diagrams commute.
We start with the left diagram.
By assumption~$\alpha \after (! \times \id) = \beta \after (\id \times !)$,
    so in particular~$\alpha(1,0) = \beta(1,0)$.
    For any~$(x,y) \in E \times F$,
        we have~$\alpha(x,0) \leq \alpha(1,0) = \beta(1,0)
            \perp \beta(0,1) \geq \beta(0,y)$,
            so~$\alpha(x,0) \perp \beta(0,y)$.
Define~$f\colon E\times F \to G$
    by~$f(x,y) = \alpha(x,0) \ovee \beta(0,y)$.
    It is easy to see~$f$ is (partially) additive
    and~$f(1,1) = \alpha(1,0) \ovee \beta(0,1) = \beta(1,1)=1$,
    so~$f$ is a module map.
We compute
\begin{equation*}
    f(x, \lambda \cdot 1)
\ = \ \alpha(x,0) \ovee \beta(0,\lambda \cdot 1)
\ = \ \alpha(x,0) \ovee \lambda \cdot \alpha(0, 1)
\ = \ \alpha(x,\lambda)
\end{equation*}
and so~$f \after (\id \times !) = \alpha$.
Similarly~$\beta = f \after (! \times \id)$.
It is easy to see~$f$
is the unique map with~$f \after (\id \times !) = \alpha$
and~$\beta = f \after (! \times \id)$
and so the left square of \eqref{pullbacks}
    is a pullback in~$\EMod_M^{\mathsf{op}}$.

Concerning the right diagram:
    as~$\delta \after (!\times!) = ! \after \pi_1$,
    we have~$\delta(\lambda \cdot 1, \mu \cdot 1) = \lambda \cdot 1$
    and so~$\delta(0, y) \leq \delta(0,1) = 0$.
    Hence~$\delta(x,y) = \delta(x,0)$.
Define~$g\colon E \to G$ by~$g(x) = \delta(x,0)$.
Clearly~$g$ is additive, $g(1) = \delta(1,0) = \delta(1,1) = 1$
    and~$\delta = g \after \pi_1$ by definition.
Obviously~$g$ is the unique such map.
Thus the right square of~\eqref{pullbacks}
    is a pullback in~$\EMod_M^{\mathsf{op}}$.
\end{point}
\begin{point}{60}{Joint epicity}%
To show~$\EMod^{\mathsf{op}}_M$ is an effectus,
    it only remains to be shown
    that~$
    \langle \pi_1, \pi_2, \pi_2\rangle,
    \langle \pi_2, \pi_1, \pi_2\rangle\colon M \times M \to M \times M \times M
    $ are jointly epic.
So assume~$f,g\colon M \times M \times M \to E$
    are two $M$-effect module maps
    with~$
        f \after \langle \pi_1,\pi_2,\pi_2\rangle =
        g \after \langle \pi_1,\pi_2,\pi_2\rangle$
        and~$
        f \after \langle \pi_2,\pi_1,\pi_2\rangle =
        g \after \langle \pi_2,\pi_1,\pi_2\rangle$.
From the first equality, it follows that~$f(1,0,0) = g(1,0,0)$.
The other one implies~$f(0,1,0) = g(0,1,0)$.
Thus~$f(0,0,1) = f(1,1,0)^\perp = g(1,1,0)^\perp = g(0,0,1)$.
As for~$(\lambda_1,\lambda_2,\lambda_3) \in M^3$,
we have~$(\lambda_1,\lambda_2,\lambda_3)
            = \lambda_1 \cdot(1,0,0)
                \ovee \lambda_2\cdot (0,1,0)
                \ovee \lambda_3\cdot (0,0,1)$
    we get~$f=g$.
\end{point}
\begin{point}{70}{Representation}%
Let~$C$ be an effectus in total form.
Clearly $\Pred f = \Pred g $ (for~$f,g\colon X \to Y$ in~$C$)
if and only if~$p \after f = (\Pred f)(p)= (\Pred g)(p) = p \after g$
for every~$p \in \Pred X$.
Thus the functor~$\Pred\colon C \to \EMod^{\mathsf{op}}_M$
    is faithful if and only if~$C$ has separating predicates.
So if~$C$ has separating predicates,
    $C$ is equivalent to the subcategory~$\Pred C$
    of~$\EMod_M^{\mathsf{op}}$.
\end{point}
    \qed
\end{point}
\end{point}
\begin{point}{80}[exc-rng-eff]{Exercise}%
    Show that the category~$\op\Rng$ of unit rings \index{Rng@$\Rng$}
    with unit-preserving homomorphisms in the opposite direction
    is an effectus in total form.
\begin{enumerate}
    \item 
        Show that the predicates on a ring~$R$
            correspond to its idempotents;
            that~$p \perp q$ iff~$pq = qp = 0$,
            $p^\perp = 1-p$
            and~$p \ovee q = p+q$.
        (So~$2$ is its effect monoid of scalars.)
        Conclude~$\op\Rng$ does not have separating predicates.
    \item
        Show that there is no unit-preserving
            ring-homomorphism~$\Z_2 \to \Z$.
            Conclude~$\op\Rng$ does not have separating states.
\end{enumerate}
\spacingfix{}
\end{point}
\end{parsec}
\begin{parsec}{1920}%
\begin{point}{10}%
We studied the structure of predicates in an effectus.  What about the states?
They will turn out to form an abstract~$M$-convex set.
Before we define these, we introduce a generalized distribution monad.
\end{point}
\begin{point}{20}{Definition}%
Let~$M$ be an effect monoid
    and~$X$ be any set.
A \Define{formal $M$-convex combination} over~$X$
    \index{convexcombination@(formal) $M$-convex combination}
    is a function~$p\colon X \to M$
    with finite support
    such that~$\bigovee_{x \in X} p(x) = 1$.
We use~$\Define{\lambda_1 \ket{x_1} \ovee \cdots \ovee \lambda_n \ket{x_n}}$
    as a shorthand for the formal $M$-convex combination~$p\colon X \to M$
    given by~$p(x) = \bigovee_{x_i = x} \lambda_i$.
    (So~$p(x_i) = \lambda_i$ if the~$x_i$ are distinct.)
We write~$\Define{\mathcal{D}_M} X$
    \index{$\mathcal{D}_M$}
    for the set of all formal~$M$-convex combinations over~$X$.~\cite{probdistrconv}
\end{point}
\begin{point}{30}[exc-dm-effectus]{Exercise*}%
    In this exercise you study~$\mathcal{D}_M$ as a monad.
    See \cite{probdistrconv,basmsc}.
\begin{enumerate}
\item
For~$f\colon X \to Y$,
    define~$\mathcal{D}_M f \colon \mathcal{D}_M X \to \mathcal{D}_M Y$
    by
    \begin{equation*}
    (\mathcal{D}_M f)(p)(y) \ =\  \bigovee_{x; f(x)=y} p(x).
    \end{equation*}
    (That is: $(\mathcal{D}_M f) (\lambda_1 \ket{x_1} \ovee \cdots
            \ovee \lambda_n \ket{x_n})
            = \lambda_1 \ket{f(x_1)} \ovee \cdots
            \ovee \lambda_n \ket{f(x_n)} $.)
Show that this turns~$\mathcal{D}_M$ into a functor~$\mathsf{Set} \to \mathsf{Set}$.
\item
Define~$\eta \colon X \to \mathcal{D}_M X$
    by~$\eta(x) = 1\ket{x}$
    and~$\mu\colon \mathcal{D}_M \mathcal{D}_M X \to \mathcal{D}_M X$ by
    \begin{equation*}
        \mu(\Phi)(x) \ = \ 
        \bigovee_{\varphi}
        \Phi(\varphi) \odot \varphi(x).
    \end{equation*}
(That is: 
        $\mu (\bigovee_i \sigma_i \ket{\bigovee_j \lambda_{ij} \ket{x_{ij}}})
            = \bigovee_{i,j} (\sigma_i \odot \lambda_{ij})\ket{x_{ij}}$.)\\
Show that~$(\mathcal{D}_M, \eta,\mu)$
    is a monad.
\item
    Show that~$\Kl \mathcal{D}_M$,
    the Kleisli category of~$\mathcal{D}_M$
    (see e.g.~\cite[\S 2.6]{basmsc}),
    is an effectus
    with~$M$ as scalars.
\end{enumerate}
\spacingfix{}
\end{point}
\begin{point}{40}{Definition}%
Let~$X$ be a set together
    with a map~$h\colon \mathcal{D}_M X \to X$.
We say~$(X,h)$ is an \Define{abstract $M$-convex set}
    \index{convex@(abstract) $M$-convex set}
    provided~$h(\ket{x})=x$
    and
\begin{equation}\label{convex-mult}
    h \Bigl( \bigovee_{i,j} \sigma_i \odot \lambda_{ij} \ket{x_{ij}} \Bigr)
    \ = \ 
    h \Bigl( \bigovee_i \sigma_i \, \ket{
        h\bigl(\, \bigovee_j \lambda_{ij} \ket{x_{ij}}\, \bigr)
    } \, \Bigr)
\end{equation}
for every~$\bigovee_i \sigma_i \ket{\bigovee_j \lambda_{ij} \ket{x_{ij}}}$
    in~$\mathcal{D}_M \mathcal{D}_M X$.
(Equivalently: $h \after \mu = h \after \mathcal{D}_M h$
and~$h \after \eta = \id$.)
    We call~$h$ the convex structure on~$X$.
A map~$f\colon X \to Y$
    between abstract~$M$-convex sets
    is \Define{$M$-affine} if
    \index{affine@$M$-affine}
\begin{equation*}
    f \bigl(h_X \bigl(\,\bigovee_i \lambda_i \ket{x_i}\, \bigr)\bigr)
     \ = \ 
     h_Y\bigl(\,\bigovee_i \lambda_i \ket{f(x_i)}\,\bigr)
\end{equation*}
for every~$\bigovee_i \lambda_i \ket{x_i}$ in~$\mathcal{D}_M X$.
(Equivalently: $h_X \after \mathcal{D}_M f = f \after h_Y$.)
Write~$\Define{\AConvM}$ for the category
    \index{AConvM@$\AConvM$}
    of abstract~$M$-convex sets with
    $M$-affine maps between them.
(Equivalently: $\AConvM$
    is the Eilenberg--Moore category
    of the monad~$\mathcal{D}_M$.)
We call an abstract $M$-convex set~$(X,h)$ \Define{cancellative}
    \index{cancellative!$M$-convex set}
provided
\begin{equation*}
    h(\lambda \ket{y} \ovee \lambda^\perp \ket{x_1})\  = \
    h(\lambda \ket{y} \ovee \lambda^\perp \ket{x_2})
    \quad\text{implies}\quad~x_1=x_2
\end{equation*}
    for any~$x_1,x_2,y \in X$ and~$\lambda \in M, \lambda \neq 1$.
\begin{point}{41}{Remark}%
    For an abstract~$[0,1]$-convex set~$(X,h)$,
    binary convex combinations are sufficient:
    for instance, for any~$\lambda \ovee \mu \neq 1$
        and~$x,y,z \in X$ we have
\begin{equation*}
    h \bigl(\lambda \ket{x} \ovee \mu \ket{y} \ovee
        (\lambda \ovee \mu)^\perp \ket{z}\bigr)
        \ = \ 
    h\Bigl(\lambda \ket{x}
        \ovee \lambda^\perp \ket{
            h\bigl( \nicefrac{\mu}{\lambda^\perp} \ket{y}
        \ovee \nicefrac{(\lambda \ovee \mu)^\perp}{\lambda^\perp}
        \ket{z}}\bigr)\Bigr).
\end{equation*}
This reduction to binary convex combinations uses
    the division on~$[0,1]$.
Some effect monoids have a suitable division as well,
        but many do not
    as we will see later on
        in~\sref{dfn-effect-divisoid} and \sref{basic-divisoid-equiv}.
It seems plausible that
    it is impossible to write an arbitrary
    three-element~$M$-convex combination
    over~$\mathcal{D}_M \{x,y,z\}$
    using only binary convex combinations
    in the case that~$M$ is the effect monoid
        on the unit interval of~$C[0,1]$.
I do not have a proof to  back this claim.
\end{point}
\end{point}
\begin{point}{50}{Examples}%
We give some examples of abstract convex sets.
\begin{enumerate}
\item
Every convex subset~$X \subseteq V$
    of some real vector space~$V$
    is a cancellative abstract~$[0,1]$-convex set
    with~$h(\lambda_1 \ket{x_1} \ovee \cdots \ovee \lambda_n \ket{x_n})
            = \lambda_1 x_1 + \cdots + \lambda_n x_n$.
\item
Not every abstract~$[0,1]$-convex set is
    a convex subset of a real vector space.
Consider~$T \subseteq [0,1]^2$
    defined by~$T \equiv \{ (x,y);\ x+y \leq 1\}$.
Say~$(x,y) \sim (x',y')$
    whenever~$x = x'$ and either $x\neq 0$ or~$y=y'$.
This equivalence relation~$\sim$
is a congruence (we will cover these in \sref{aconv-cong})
        for the convex structure on~$T$
    (inherited by~$\R^2$)
    and so~$T /_\sim$ is an abstract~$[0,1]$-convex set.
That is: we start off with a filled triangle
$T \equiv \begin{tikzpicture}[scale=0.2]
    \draw (0,0) -- (0,1) ; 
    \draw (1,0) -- (0,0) ; 
    \draw (1,0) -- (0,1) ; 
\end{tikzpicture}$,
identify all vertical lines
except for the y-axis
and are left with
$T/_\sim \equiv
\begin{tikzpicture}[scale=0.2]
    \draw (0,0) -- (0,1) ; 
    \draw (1,0) -- (0,0) ; 
\end{tikzpicture}$.
The abstract~$[0,1]$-convex set~$T/_\sim$
    is not cancellative:
\begin{equation*}
    h\Bigl( 
    \frac{1}{2}\ket{\begin{tikzpicture}[scale=0.2]
    \filldraw( 1,0) circle (3pt);
    \draw (0,0) -- (0,1) ; 
    \draw (1,0) -- (0,0) ; 
\end{tikzpicture}}
\ovee
    \frac{1}{2}\ket{\begin{tikzpicture}[scale=0.2]
    \filldraw( 0,1) circle (3pt);
    \draw (0,0) -- (0,1) ; 
    \draw (1,0) -- (0,0) ; 
\end{tikzpicture}} \Bigr)
\ = \
    \begin{tikzpicture}[scale=0.2]
    \filldraw( .5,0) circle (3pt);
    \draw (0,0) -- (0,1) ; 
    \draw (1,0) -- (0,0) ; 
\end{tikzpicture}
\ = \ 
h \Bigl(\frac{1}{2}\ket{\begin{tikzpicture}[scale=0.2]
    \filldraw( 1,0) circle (3pt);
    \draw (0,0) -- (0,1) ; 
    \draw (1,0) -- (0,0) ; 
\end{tikzpicture}}
\ovee
    \frac{1}{2}\ket{
        \begin{tikzpicture}[scale=0.2]
    \filldraw( 0,.5) circle (3pt);
    \draw (0,0) -- (0,1) ; 
    \draw (1,0) -- (0,0) ; 
\end{tikzpicture}
}\Bigr), \quad\text{but}\quad
        \begin{tikzpicture}[scale=0.2]
    \filldraw( 0,1) circle (3pt);
    \draw (0,0) -- (0,1) ; 
    \draw (1,0) -- (0,0) ; 
\end{tikzpicture} \ \neq \ 
        \begin{tikzpicture}[scale=0.2]
    \filldraw( 0,.5) circle (3pt);
    \draw (0,0) -- (0,1) ; 
    \draw (1,0) -- (0,0) ; 
\end{tikzpicture}.
\end{equation*}
\item
Semilattices are exactly~abstract $2$-convex sets.
    ($2$ is the two-element Effect Monoid,
        see \sref{eff-monoid-examples}.)
Furthermore, every semilattice~$(L, \vee)$
    is an abstract~$[0,1]$-convex set
    with
    \begin{equation*}
    h\bigl(\lambda_1 \ket{x_1} \ovee \cdots \ovee \lambda_n \ket{x_n}\bigr)
        \ =\  \bigvee_{i; \lambda_i \neq 0} x_i.
    \end{equation*}
    A semilattice~$L$ is cancellative as abstract $[0,1]$-convex set
    if and only if~$x=y$ for all~$x,y\in L$.
    See~\cite{neumann1970quasivariety}.
\item
Every cancellative~$[0,1]$-convex set
    is isomorphic to a convex subset of a real vector space.
    See e.g.~\cite[thm.~8]{statesofconvexsets}.
\end{enumerate}
\spacingfix{}
\begin{point}{60}{Remarks}%
The study of
    convex subsets (of real vector spaces)
    as algebraic structures goes back a long time,
        see e.g.~\cite{stone1949postulates,neumann1970quasivariety,flood1981semiconvex,tz2009convex,gudder1979general}.
The more pathological abstract~$[0,1]$-convex sets have been studied before
    under different names:
    they are semiconvex sets in \cite{flood1981semiconvex,swirszcz1975monadic} and
    convex spaces in \cite{tz2009convex}.
The description of convex sets using Eilenberg--Moore algebras is probably
    due to \'Swirszcz \cite{swirszcz1975monadic}.
For the even more exotic abstract~$M$-convex
    sets, see~\cite{effintro} or
        \cite{probdistrconv}.
\end{point}
\end{point}
\begin{point}{70}{Proposition}%
Let~$C$ be an effectus with scalars~$M$.
For any object~$X$, the set of states
$\Stat X$ is an abstract $M^{\mathsf{op}}$-convex set
with~$h\colon \mathcal{D}_{M^{\mathsf{op}}} \Stat X \to \Stat X$
    defined by
    $h\bigl( \lambda_1 \ket{\varphi_1} \ovee \cdots \ovee \lambda_n \ket{\varphi_n} \bigr)
     =  [\varphi_1, \ldots, \varphi_n] \after \langle
                    \lambda_1, \ldots, \lambda_n\rangle$.
Furthermore,
    for any total map~$f\colon X \to Y$,
    the map~$ (\Stat f)(\varphi) \ = \ f \after \varphi$
    is an~$M^{\mathsf{op}}$-affine map~$\Stat X \to \Stat Y$ yielding a
    functor~$\Stat \colon \Tot C \to \mathsf{AConv}_{M^{\mathsf{op}}}$.~\cite{effintro}
\begin{point}{80}{Proof}%
Clearly~$h(\ket{\varphi}) = \varphi$.
To show the other axiom~\eqref{convex-mult},
assume we have any~$\bigovee^n_{i=1} \sigma_i \ket{\bigovee^{m_i}_{j=1} \lambda_{ij} \ket{\varphi_{ij}}}$
in~$\mathcal{D}_{M^\mathsf{op}} \mathcal{D}_{M^{\mathsf{op}}} \Stat X$.
Write~$\lambda_i \equiv \langle\lambda_{i1}, \ldots, \lambda_{im_i}\rangle$
and $\varphi_i \equiv [\varphi_{i1}, \ldots, \varphi_{i{m_i}}]$.
Then
\begin{align*}
    h \Bigl( \bigovee_i \sigma_i \, \ket{
        h\bigl(\, \bigovee_j \lambda_{ij} \ket{\varphi_{ij}}\, \bigr)
    } \, \Bigr)
    &\ = \ 
    [\varphi_1 \after \lambda_1, \ldots, 
    \varphi_n \after \lambda_n] \after \langle
    \sigma_1, \ldots, \sigma_n
    \rangle
    \\ & \ = \ 
    \bigovee_{i} \varphi_i \after \lambda_i \after \sigma_i
    \\ & \ = \ 
    \bigovee_{i} \bigl(\bigovee_j \varphi_{ij} \after \lambda_{ij}\bigr) \after \sigma_i
    \\ & \ = \ 
    \bigovee_{i,j} \varphi_{ij} \after (\lambda_{ij} \after \sigma_i)
    \\ & \ = \ 
    \bigovee_{i,j} \varphi_{ij} \after (\sigma_i \odot_{M^{\mathsf{op}}} \lambda_{ij})
    \\ & \ = \ 
    h \Bigl( \bigovee_{i,j} \sigma_i \odot_{M^{\mathsf{op}}} \lambda_{ij} \ket{\varphi_{ij}} \Bigr).
\end{align*}
So~$\Stat X$ is indeed an abstract~$M^{\mathsf{op}}$-convex set.
Let~$f\colon X \to Y$ be any total map.
To show~$\Stat f$ is affine,
it is sufficient to show that
for every~$\bigovee_i \lambda_i \ket{\varphi_i}$
 in~$\mathcal{D}_{M^{\mathsf{op}}}$,
 we have
        $f \after \bigovee_i  \varphi_i \after \lambda_i 
         = 
        \bigovee_i f \after  \varphi_i \after \lambda_i$,
        which is clearly true.
    Functoriality of~$\Stat$ is obvious. \qed
\end{point}
\end{point}
\end{parsec}

\subsection{Effectus of abstract $M$-convex sets}\label{more-aconvm}
\begin{parsec}{1930}%
\begin{point}{10}%
For our next project,
    we indulge ourselves in a tangent:
    we investigate whether the category~$\AConvM$ is an effectus.
We can show it is, if~$M$ has a certain partial division.
Unfortunately it will remain unclear if~$\AConvM$ is an effectus in general.
The case~$M=[0,1]$ is much easier and has already been dealt with
        in \cite{statesofconvexsets}. 
The first step is to study the coproduct in~$\AConvM$.
To construct the coproduct, we will need to know
    about quotients of abstract~$M$-convex sets.
\end{point}
\begin{point}{20}[aconv-cong]{Exercise}%
In this exercise you construct the quotient of an
    abstract~$M$-convex set by a congruence.
Assume~$(X,h)$ is an abstract~$M$-convex set
    and~${\sim} \subseteq X^2$ is an equivalence relation.
Write~$X/_\sim$ for the set of equivalence classes
    of~$\sim$
    and~$q\colon X \to X/_\sim$
    for~$q(x) = [x]_\sim$, where~$[x]_\sim$ is the equivalence class of~$x$.
For~$\varphi,\psi \in \mathcal{D}_M X$,
    we write~$\varphi \sim \psi$
    provided~$(\mathcal{D}_M q)(\varphi) = (\mathcal{D}_M q)(\psi)$.
We say~$\sim$ is a \Define{congruence} for~$(X,h)$
    \index{congruence!$M$-convex set}
if~$\varphi \sim \psi$ implies~$h(\varphi) \sim h(\psi)$.
\begin{enumerate}
\item
Prove that~$q$, $\mathcal{D}_M q$ and~$\mathcal{D}_M \mathcal{D}_M q$
    are all surjective.
\item
Show that the following are equivalent.
\begin{enumerate}
\item
$\sim$ is a congruence.
\item
There is a map~$h_\sim \colon \mathcal{D}_M X/_\sim \to X/_\sim$
    such that~$h_\sim \after \mathcal{D}_M q = q \after h$.
\end{enumerate}
By the previous point $h_\sim$ is unique, if it exists.
\item
Assume~$\sim$ is a congruence.
Use the previous two points to show
\begin{equation*}
h_\sim \after \mathcal{D}_M h_\sim \after \mathcal{D}_M \mathcal{D}_M q
            \ =\  h_\sim \after \mu \after \mathcal{D}_M \mathcal{D}_M q.
\end{equation*}
Conclude that~$(X/_\sim,h_\sim)$ is an abstract~$M$-convex set
    and that the qoutient map~$q\colon X \to X/_\sim$ is~$M$-affine.
\end{enumerate}
\spacingfix{}
\end{point}
\begin{point}{30}[affine-kernel-cong]{Exercise}%
Assume~$f\colon X \to Y$ is an affine map between abstract
    $M$-convex sets.
Show that the kernel~$\{(x,y);\ f(x) = f(y)\}$ of~$f$
    is a congruence.
\end{point}
\begin{point}{40}[least-conv-cong]{Exercise}%
Assume~$(X,h)$ is an abstract~$M$-convex set and~$R \subseteq X^2$
    is some relation.
It is easy to see there is a least congruence containing~$R$.
We need to know a bit more than mere existence.
Write~$R^*$ for the reflexive symmetric transitive closure of~$R$
    (i.e.~the least equivalence relation containing~$R$).

For~$\psi_1,\psi_2 \in \mathcal{D}_M X$,
    we write~$\psi_1 \approx \psi_2$
    if there is a \emph{derivation}: a tuple~$\varphi_1, \ldots, \varphi_n \in \mathcal{D}_M X$
(say~$\varphi_i \equiv \bigovee^{n_i}_{j = 1}  \lambda_{ij}\ket{x_{ij}}$)
such that~$\varphi_1 = \psi_1$, $\varphi_n = \psi_2$ 
    and for each~$1 \leq i<n$ \emph{one} of the following two conditions holds.
\begin{enumerate}
\item
    $n_i = n_{i+1}$
    and for all~$1 \leq j \leq n_i$ we have
    $x_{ij} \mathrel{R^*} x_{(i+1)j}$
    and~$\lambda_{ij}  = \lambda_{(i+1)j}$.
\item
    $h (\varphi_i) = h(\varphi_{i+1})$.
\end{enumerate}
We will show that~$\sim$
    given by~$x \sim y \iff \eta (x) \approx \eta(y)$
    is the least congruence containing~$R$:
\begin{enumerate}
\item
To start, show that~$\eta ( h(\psi)) \approx \psi $
    and so~$\varphi \approx \psi$ implies~$h(\varphi) \sim h(\psi)$.
\item
Show that if~$\varphi \approx \psi$,
    then
    \begin{equation*}
    \mu \Bigl(\lambda_0 \ket{\psi} \ovee \bigovee_{j=1}^n \lambda_j \ket{\chi_j}\Bigr)
    \ \approx  \ 
    \mu \Bigl(\lambda_0 \ket{\varphi} \ovee \bigovee_{j=1}^n \lambda_j \ket{\chi_j}\Bigr).
    \end{equation*}
    for any~$\bigovee_j \lambda_j=1$  in~$M$ and~$\chi_1,\ldots, \chi_n \in \mathcal{D}_M X$.
\item
    Use the previous point to show that~$\varphi \sim \psi$
        implies~$\varphi \approx \psi$.
Conclude~$\sim$ is the smallest congruence containing~$R$.
\end{enumerate}
\spacingfix{}
\end{point}
\begin{point}{50}[aconv-coprod]{Proposition}%
Assume~$(X,h_X)$ and~$(Y,h_Y)$
    are abstract~$M$-convex sets.
    We will construct their coproduct.
Note that~$\mu$ turns~$\mathcal{D}_M (X+Y)$
    into an abstract~$M$-convex set.
Let~$\sim$ denote the least congruence on~$\mathcal{D}_M (X+Y)$
with
\begin{equation}
    (\mathcal{D}_M \kappa_1) (\varphi)
    \ \sim \ \eta(\kappa_1 ( h_X (\varphi)))
    \quad \text{and} \quad
    (\mathcal{D}_M \kappa_2) (\psi)
    \ \sim \ \eta(\kappa_2 ( h_Y (\psi))) \label{congruence-coprod-conv}
\end{equation}
    for all~$\varphi \in \mathcal{D}_M X$
    and~$\psi \in \mathcal{D}_M Y$.
(Alternatively: $\sim$ is the least congruence
    that makes~$\eta \after \kappa_1$ and~$\eta\after\kappa_2$ affine.)
Write~$q\colon \mathcal{D}_M (X+Y) \to C$
    for the quotient of~$\mathcal{D}_M(X+Y)$ by~$\sim$, see~\sref{aconv-cong}.
    Denote the convex structure on~$C$ by~$h$.
Define~$c_1 \colon X \to C$ and~$c_2 \colon Y \to C$
    by~$c_1 \equiv q \after \eta \after \kappa_1$
    and~$c_2 \equiv q \after \eta \after \kappa_2$.

Then: $c_1 \colon X \to C \leftarrow Y \colon c_2$
is a coproduct of~$X$~and~$Y$ in~$\AConvM$.
\begin{point}{60}{Proof}%
We start with the affinity of the coprojections:
\begin{alignat*}{2}
    c_1 \after h_X & \ = \ 
    q \after \eta \after \kappa_1 \after h_X \\
    &\ =\  q \after \mathcal{D}_M \kappa_1
        &\qquad& \text{by dfn.~of~$\sim$}
    \\
    & \ =\  q \after \mu \after \mathcal{D}_M \eta \after \mathcal{D}_M \kappa_1. \\
    & \ =\  \mu \after \mathcal{D}_Mq \after \mathcal{D}_M \eta \after \mathcal{D}_M \kappa_1. \\
    & \ =\  \mu \after \mathcal{D}_M c_1,
\end{alignat*}
so~$c_1$ is affine. With similar reasoning, we see~$c_2$ is affine.
\begin{point}{70}%
Assume~$f\colon X \to Z$ and~$g\colon Y \to Z$
    are two affine maps to some abstract~$M$-convex set~$(Z,h_Z)$.
We want to show~$h_Z \after \mathcal{D}_M [f,g]$
    lifts to~$C$.
With an easy calculation,
    one sees~$h_Z \after \mathcal{D_M} [f,g]$ is affine.
We compute
\begin{align*}
    h_Z \after \mathcal{D}_M [f,g] \after \mathcal{D}_M \kappa_1
     & \ = \ 
    h_Z \after \mathcal{D}_M f \\
     & \ = \ f \after h_X \\
     & \ = \ h_Z \after \eta \after [f,g] \after \kappa_1 \after h_X \\
     & \ = \ h_Z \after \mathcal{D}_M [f,g] \after \eta \after \kappa_1 \after h_X.
\end{align*}
So~$
(\,(\mathcal{D}_M \kappa_1) (\varphi), \,
\ \eta(\kappa_1 ( h_X (\varphi)))\,)$
is in the kernel (see \sref{affine-kernel-cong})
of~$h_Z \after \mathcal{D}_M [f,g]$
for~$\varphi \in \mathcal{D}_M$.
Similarly~$(\,(\mathcal{D}_M \kappa_2) (\psi), \,
\ \eta(\kappa_2 ( h_Y (\psi)))\,)$
is also in the kernel of~$h_Z \after \mathcal{D}_M [f,g]$
for all~$\psi \in \mathcal{D}_M$.

As kernels of affine maps are congruences
    and~$\sim$ is the least congruence that contains these pairs,
    we conclude~$\sim$ is a subset of the kernel
    of~$h_Z \after \mathcal{D}_M [f,g]$.
Thus there is a unique affine~$k\colon C \to Z$
    fixed by~$k \after q = h_Z \after \mathcal{D}_M [f,g]$.
Now
\begin{align*}
    k \after c_1 
    &\ = \ k \after q \after \eta \after \kappa_1 \\
    &\ = \ h_Z \after \mathcal{D}_M [f,g]
                            \after \eta\after \kappa_1 \\
    &\ = \ h_Z \after \eta \after  [f,g]
                            \after \kappa_1 \\
    &\ = \  f
\end{align*}
and similarly~$k \after c_2 = g$.
\end{point}
\begin{point}{80}%
Only uniqueness of~$k$ remains.
So assume~$k'\colon C \to Z$ is an affine map
    such that~$k' \after c_1 = f$
    and~$k' \after c_2 = g$.
We turn the wheel:
\begin{align*}
    k' \after q
    & \ = \ k' \after q \after \mu \after \mathcal{D}_M \eta \\
    & \ = \ k' \after h \after \mathcal{D}_M (q \after \eta) \\
    & \ = \ h_Z \after \mathcal{D}_M (k' \after q \after \eta) \\
    & \ = \ h_Z \after \mathcal{D}_M 
    [ k' \after q \after \eta \after \kappa_1,
    k' \after q \after \eta \after \kappa_2]  \\
    & \ = \ h_Z \after \mathcal{D}_M 
    [ k' \after c_1,
    k' \after c_2]  \\
    & \ = \ h_Z \after \mathcal{D}_M [ f, g] \\
    & \ = \ k \after q.
\end{align*}
So~$k = k'$ by surjectivity of~$q$. \qed
\end{point}
\end{point}
\begin{point}{90}[elements-coprod-conv]{Remark}%
The existence of coproducts of abstract~$M$-convex sets already
    follows from the more general theorem of Linton \cite{linton}.
    (Cf.~\cite{barr}.)
However, for the results in the remainder of this section we need to know
    more about the coproduct than Linton's construction provides.

By our construction, we know that
    every~$z \in X+Y$
    is of the
    form
\begin{equation*}
    z \ =\  h 
    \Bigl( \bigovee^n_{i=1} \lambda_i \ket{c_1(x_i)}
    \ovee \bigovee^m_{i=1} \sigma_i \ket{c_2(y_i)} \Bigr)
\end{equation*}
for some~$\lambda_i,\sigma_i \in M $, 
$x_i \in X$ and~$y_i \in Y$.
In general these are not unique.
Indeed
using \sref{least-conv-cong}, we see that
\begin{equation*}
     h 
     \Bigl( \bigovee^n_{i=1} \lambda_i \ket{c_1(x_i)}
     \ovee \bigovee^m_{i=1} \sigma_i \ket{c_2(y_i)} \Bigr) 
     \ = \ 
     h 
     \Bigl( \bigovee^{n'}_{i=1} \lambda'_i \ket{c_1(x'_i)}
     \ovee \bigovee^{m'}_{i=1} \sigma'_i \ket{c_2(y'_i)} \Bigr)
\end{equation*}
iff there is a \emph{derivation}~$\Phi_1, \ldots, \Phi_l \in \mathcal{D}^2_M(X+Y)$,
say $\Phi_i \equiv \bigovee_j \zeta_{ij} \ket{\varphi_{ij}}$,
with
\begin{equation*}
\Phi_1 \ =\  
\ket {\bigovee^n_{i=1} \lambda_i \ket{\kappa_1(x_i)}
\ovee \bigovee^m_{i=1} \sigma_i \ket{\kappa_2(y_i)}}
, \quad
 \Phi_l \ =\  
 \ket {\bigovee^{n'}_{i=1} \lambda'_i \ket{\kappa_1(x_i)}
 \ovee \bigovee^{m'}_{i=1} \sigma'_i \ket{\kappa_2(y_i)}}
\end{equation*}
and for each~$1 \leq i < l$ either
\begin{enumerate}
    \item $\mu(\Phi_i) = \mu(\Phi_{i+1})$
        or
    \item
        for each~$j$ we have~$\zeta_{ij} = \zeta_{(i+1)j}$
        and one of the following.
        \begin{enumerate}
        \item
    $\varphi_{ij} = \varphi_{(i+1)j}$.
\item
    $\varphi_{ij} = (\mathcal{D}_M\kappa_1) (\chi)$
and~$\varphi_{(i+1)j} = (\mathcal{D}_M\kappa_1) (\chi')$
        for some~$\chi,\chi' \in \mathcal{D}_M X$
        with~$h_X( \chi) = h_X(\chi')$.
\item
    $\varphi_{ij} = (\mathcal{D}_M\kappa_2) (\chi)$
and~$\varphi_{(i+1)j} = (\mathcal{D}_M\kappa_2) (\chi')$
        for some~$\chi,\chi' \in \mathcal{D}_M Y$
        with~$h_Y( \chi) = h_Y(\chi')$.
        \end{enumerate}
\end{enumerate}
\end{point}
\end{point}
\spacingfix{}
\begin{point}{100}[n-times-one-aconvm]{Exercise}%
Show that the one-element set is an abstract~$M$-convex set
    and, in fact, the final object of~$\AConvM$.
Moreover, show that in~$\AConvM$
\begin{equation*}
    n \cdot 1 \ \equiv \ \underbrace{1 + \cdots + 1}_{\text{$n$ times}}
    \ \cong \ \mathcal{D}_M \{1, \ldots, n\}.
\end{equation*}
(You might want to use that~$\mathcal{D}_M$, as a left-adjoint,
    preserves coproducts.)
\end{point}
\end{parsec}

\begin{parsec}{1940}%
\begin{point}{10}[aconvalmosteffectus]{Proposition}%
The category~$\AConvM$ obeys all the axioms of an effectus in total
form except perhaps for the left pullback square of~\eqref{pullbacks}.
\begin{point}{20}{Proof}%
In \sref{aconv-coprod} we proved~$\AConvM$ has binary coproducts
    and in \sref{n-times-one-aconvm}
    that the one-element set~$1$ is its final object.
As~$\mathcal{D}_M \emptyset = \emptyset$,
    the empty set is trivially also an abstract~$M$-convex set
    and in fact the initial object of~$\AConvM$.
\end{point}
\begin{point}{30}%
We continue with joint monicity
    of~$[\kappa_1,\kappa_2,\kappa_2],
    [\kappa_2,\kappa_1,\kappa_2]\colon 1+1+1 \to 1+1$.
We can represent the convex set~$1+1+1$
    as the set of triples~$(a,b,c) \in M^3$
    with~$a\ovee b \ovee c = 1$
    with the obvious convex structure,
    cf.~\sref{n-times-one-aconvm}.
Similarly, we identify~$1+1$
    with pairs~$(a,a^\perp) \in M^2$.
    The maps at hand are given by
\begin{equation}\label{conv-jointly-monic}
    (a,b,c) \ \mapsto\  (a , b\ovee c)\quad \text{and} \quad
    (a,b,c) \ \mapsto\  (b, a\ovee c).
\end{equation}
Assume~$
(a, b \ovee c) = 
(a', b' \ovee c') $
and
$(b, a \ovee c) = 
(b', a' \ovee c') $
for~$a,a',b,b',c,c' \in M$ with~$a \ovee b \ovee c= 1$
and~$a' \ovee b' \ovee c'= 1$.
Then clearly~$a=a'$,~$b=b'$
    and so~$c = (a\ovee b)^\perp = (a' \ovee b')^\perp = c'$.
Thus the maps~\eqref{conv-jointly-monic} are jointly injective,
    hence jointly monic.
\end{point}
\begin{point}{40}%
At last, we are ready to prove that the square on the right of \eqref{pullbacks}
    is a pullback diagram in~$\AConvM$.
So assume~$\alpha \colon Z \to X+Y$
    is a map in~$\AConvM$
    with~$(!+!) \after \alpha = \kappa_1 \after !$.
\begin{equation*}
\xymatrix{
    Z \ar@/^1pc/[rrd]^{!}
        \ar@/_1pc/[rdd]_\alpha
        \ar@{.>}[rd]|\gamma\\
        &X \ar[r]^{!} \ar[d]^{\kappa_1}
        & 1\ar[d]^{\kappa_1} \\
        &X+Y\ar[r]_{!+!} & 1+1
    }
\end{equation*}
We have to show~$\alpha = \kappa_1 \after \gamma$,
    for a unique~$\gamma\colon Z \to X+Y$ in~$\AConvM$.
For the moment, assume~$\kappa_1$ is injective.
Then~$\gamma$, if it exists, is unique.
Pick~$z \in Z$.
We know
\begin{equation*}
    \alpha(z) \ =\  h 
    \Bigl( \bigovee^n_{i=1} \lambda_i \ket{\kappa_1(x_i)}
    \ovee \bigovee^m_{i=1} \sigma_i \ket{\kappa_2(y_i)} \Bigr)
\end{equation*}
for some~$\lambda_i,\sigma_i \in M $, 
$x_i \in X$ and~$y_i \in Y$,
where~$h$ is the convex structure on~$X+Y$.
As~$(!+!) \after \alpha = ! \after \kappa_1$,
    we must have~$\bigovee_i \sigma_i = 0$
    and so~$\alpha(z) = \kappa_1(x_z)$
    for some~$x_z \in X$.
This~$x_z$ is unique by supposed injectivity of~$\kappa_1$.
    Define~$\gamma(z) = x_z$.
By definition~$\kappa_1 \after \gamma = \alpha$.
To see~$\gamma$ is affine, pick any~$\varphi \in \mathcal{D}_M (Z)$.
We have
$
\kappa_1 ( \gamma ( h_Z (\varphi )))
    = \alpha(h_Z (\varphi))
    = h(\mathcal{D}_M\alpha(\varphi))
    = h(\mathcal{D}_M\kappa_1 ( \mathcal{D}_M \gamma(\varphi)))
    = \kappa_1  (h_X( \mathcal{D}_M \gamma(\varphi)))
$ and
so indeed~$\gamma(h_Z(\varphi)) = h_X(\mathcal{D}_M \gamma (\varphi))$
    --- i.e.~$\gamma$ is affine.
\end{point}
\begin{point}{50}%
Finally, we will show that~$\kappa_1$ is indeed injective.
Assume~$\kappa_1(x_0) = \kappa_1(x'_0)$
for some~$x_0,x'_0 \in X$.
    Let~$\Phi_1, \ldots, \Phi_l \in \mathcal{D}_M^2{(X+Y)}$
    be a derivation
    of~$\kappa_1(x_0) = \kappa_1(x'_0)$
    as in \sref{elements-coprod-conv}.
The remainder of the proof is, in essence, a simple induction
    over~$1 \leq i < l$, which has been complicated
    by a technicality.
Define
\begin{equation*}
    \mathrm{IH}(i) \ \equiv \ \text{``}
    \bigovee_y \mu(\Phi_i)(\kappa_2 (y)) = 0
        \text{ and }
        h_X\Bigl(\bigovee_{x} \mu(\Phi_i)(\kappa_1(x)) \ket{x}\Bigr)
        =    x_0 \text{''}.
\end{equation*}
Clearly~$\mathrm{IH}(1)$ holds
    and if~$\mathrm{IH}(l)$ should hold,
    then~$x_0 = x_0'$ as desired.
So, to prove the inductive step, assume~$\mathrm{IH}(i)$ holds.
There are two possible derivation steps, see \sref{elements-coprod-conv}.
In the first case (i.e.~$\mu(\Phi_i) = \mu(\Phi_{i+1})$)
    it is obvious~$\mathrm{IH}(i+1)$ holds.
So we continue with the other case,
where (among others) we have~$\zeta_{ij} = \zeta_{(i+1)j}$,
    where~$\Phi_i \equiv \bigovee_j \zeta_{ij} \ket{\varphi_{ij}}$.
Without loss of generality, we may assume~$\zeta_{ij} \neq 0$.
For each~$j$ there are three possibilities, see \sref{elements-coprod-conv}.
The third does not occur:
if~$\varphi_{ij_0} = (\mathcal{D}_M \kappa_2) (\chi)$
    for some~$\chi \in \mathcal{D}_M Y$ and~$j_0$,
    then we reach a contradiction:
\begin{equation*}
    0\ = \ 
    \bigovee_y \mu(\Phi_i)(\kappa_2(y))
    \ = \ \bigovee_{y,j} \zeta_{ij} \odot \varphi_{ij}(\kappa_2(y))
    \ \geq \ \bigovee_y \zeta_{ij_0} \odot \chi(y)\  \neq\  0.
\end{equation*}
We are left with two possibilities ---
    to distinguish those,
    write~$j \in J$ if~$\varphi_{ij} = \varphi_{(i+1)j}$
    and~$j \notin J$
    if~$\varphi_{ij} = (\mathcal{D}_M \kappa_1)(\chi)$
    and~$\varphi_{(i+1)j} = (\mathcal{D}_M \kappa_1)(\chi')$
    for some~$\chi,\chi' \in \mathcal{D}_M X$
    with~$h_X (\chi) = h_X(\chi')$.

    We will show the first part of~$\mathrm{IH}(i+1)$.
Pick any~$y \in Y$.
If~$j \in J$, then~$\zeta_{(i+1)j}\odot \varphi_{(i+1)j}(\kappa_2(y))
= \zeta_{ij} \odot \varphi_{ij}(\kappa_2(y)) = 0$.
In the other case, if~$j \notin J$,
then clearly~$\zeta_{(i+1)j}\odot \varphi_{(i+1)j}(\kappa_2(y)) = 0$
    as~$\varphi_{(i+1)j} = (\mathcal{D}_M \kappa_1)(\chi')$.
So, we have~$\bigovee_y \mu(\Phi_{i+1})(\kappa_2(y))
= \bigovee_{y,j} \zeta_{(i+1)j} \odot \varphi_{(i+1)j} (\kappa_2(y)) = 0 $.

Write~$r_{ij} = (\bigovee_x \varphi_{ij}(\kappa_1(x)))^\perp$.
The remainder of the proof is complicated
    by the fact that in general~$r_{ij} \neq 0$.
    We do have~$\bigovee_j \zeta_{ij} \odot r^\perp_{ij} = 1$
        and~$\bigovee_j \zeta_{ij} = 1$
        so by \sref{emond-lemma-for-conv}
        we see~$\zeta_{ij} \odot r^\perp_{ij} = \zeta_{ij}$,
        which  forces~$\zeta_{ij} \odot r_{ij} = 0$.
We compute
\begin{align*}
    h_X\Bigl(\bigovee_x \mu(\Phi_i) (\kappa_1(x)) \ket{x}\Bigr) 
    &\ = \ 
    h_X \Bigl( \bigovee_{x,j} \zeta_{ij} \odot \varphi_{ij}
        (\kappa_1(x)) \ket{x} \Bigr) \\
        &\ = \ 
    h_X \Bigl(
    \bigovee_{j} \zeta_{ij} \odot r_{ij} \ket{x_0} \ovee \bigovee_x
    \zeta_{ij} \odot \varphi_{ij}
        (\kappa_1(x)) \ket{x} \Bigr) \\
        & \ = \ h_X \Bigl(
    \mu\Bigl( \bigovee_j \zeta_{ij} \ket{
        r_{ij} \ket{x_0} \ovee \bigovee_x \varphi_{ij} (\kappa_1 (x)) \ket{x}
    }
    \Bigr)
    \Bigr) \\
        & \ = \ h_X \Bigl(
    \bigovee_j \zeta_{ij} \ket{
        h_X\Bigl(
        r_{ij} \ket{x_0} \ovee \bigovee_x \varphi_{ij} (\kappa_1 (x)) \ket{x}
    \Bigr)
}
    \Bigr).
\end{align*}
So, if we show that for each~$j$,
    we have
\begin{multline}\label{eqkappainj}
        h_X\Bigl(
        r_{ij} \ket{x_0} \, \ovee \,\bigovee_x \varphi_{ij} (\kappa_1 (x)) \ket{x}
    \Bigr) \\
    \  = \  
        h_X\Bigl(
        r_{(i+1)j} \ket{x_0} \,\ovee \,\bigovee_x \varphi_{(i+1)j} (\kappa_1 (x)) \ket{x}
    \Bigr),
\end{multline}
then
\begin{equation*}
    x_0 \ = \ h_X\Bigl(\bigovee_x \mu(\Phi_i) (\kappa_1(x)) \ket{x}\Bigr) 
    \ =\  h_X\Bigl(\bigovee_x \mu(\Phi_{i+1}) (\kappa_1(x)) \ket{x}\Bigr),
\end{equation*}
as desired.
If~$j \in J$,
then~$\varphi_{ij} = \varphi_{(i+1)j}$
and so~$r_{ij} = r_{(i+1)j}$
and obviously~\eqref{eqkappainj}.
So assume~$j \notin J$.
Then clearly~$r_{ij} = 0 = r_{(i+1)j}$
as~$\varphi_{ij} = (\mathcal{D}_M \kappa_1)(\chi))$
    and~$\varphi_{(i+1)j} = (\mathcal{D}_M \kappa_1) (\chi')$.
Now~\eqref{eqkappainj} follows from~$h_X(\chi) = h_X(\chi')$. \qed
\end{point}
\end{point}
\end{parsec}

\begin{parsec}{1950}%
\begin{point}{10}%
It is unclear whether the left square of~\eqref{pullbacks}
    is a pullback in~$\AConvM$.
We will see that it is, if we assume that~$M$ has
        a partial division in the following sense.
\end{point}
\begin{point}{20}[dfn-effect-divisoid]{Definition}%
An \Define{effect divisoid}
    \index{effect divisoid}
    is an effect monoid~$M$
    with partial binary operation~$(a,b) \mapsto \rfrac{a}{b}$
    that is defined iff~$a \leq b$ and
        satisfies the following axioms.
\begin{enumerate}
\item~$\rfrac{a}{b}$
        is the unique element of~$M$
        with~$\rfrac{a}{b} \leq \rfrac{b}{b}$
        and~$b\odot \rfrac{a}{b} \ =\  a$
\item~$a \ \leq \ \rfrac{a}{a}$
\item~$\rfrac{(\rfrac{a}{a})}{(\rfrac{a}{a})} \ =\  \rfrac{a}{a}$
\end{enumerate}
\spacingfix{}
\begin{point}{21}{Beware}%
The first axiom might be read as ``$\rfrac{a}{b}$ is defined to be the
    unique element such that (...)'' in which case the definition
    is cyclical and thus problematic.
    What is meant instead is that~$\rfrac{a}{b}$ is
        part of the structure and satisfies
\begin{itemize}
\item $\rfrac{a}{b} \leq \rfrac{b}{b}$
        and~$b\odot \rfrac{a}{b} \ =\  a$ and
\item
        if~$c \leq \rfrac{b}{b}$
        and~$b\odot c \ =\  a$,
        then~$c = \rfrac{a}{b}$
\end{itemize}
for all~$a,b,c \in M$ with~$a\leq b$.
In particular, it is not ruled out
    that there is an effect monoid
    which is an effect divisoid in two different ways.
\end{point}
\begin{point}{30}{Remark}%
The element~$\rfrac{a}{a}$ behaves like a support-projection for~$a$.
\end{point}
\end{point}
\begin{point}{40}[exc-divisoid-basics]{Exercise}%
Show that for an effect divisoid~$M$, the following holds.
\begin{enumerate}
\item
    $\rfrac{0}{0}=0$, $\rfrac{1}{1}=1$, $\rfrac{a}{1}=a$, $\rfrac{a}{a} \odot \rfrac{a}{a} = \rfrac{a}{a}$ and
    $\rfrac{a \odot b}{a} = \rfrac{a}{a }\odot b$.
\item
    For any~$a \leq b \leq c$,
    we have~$\rfrac{b}{c} \odot \rfrac{a}{b} = \rfrac{a}{c}$.
\end{enumerate}
\spacingfix{}
\end{point}
\begin{point}{50}{Examples}%
Almost all effect monoids we encountered
    are effect divisoids.
\begin{enumerate}
\item
The effect monoid~$[0,1]$
    is an effect divisoid with~$\rfrac{a}{b} = \frac{a}{b}$
        if~$b \neq 0$ and~$\rfrac{0}{0}=0$.
With the same partial division, we see the two-element effect monoid~$2$
    is also an effect divisoid.
\item
More interestingly,
    if~$M_1$ and~$M_2$
    are effect divisoids,
    then~$M_1 \times M_2$
    is an effect divisoid
    with~$\rfrac{(a_1,a_2)}{(b_1,b_2)} =
    (\rfrac{a_1}{b_1}, \rfrac{a_2}{b_2})$.
In particular~$[0,1]^n$ is an effect divisoid
    with component-wise partial division.
\item
Later, in \sref{andthen-effect-divisoid},
    we will see that if~$C$ is an~$\&$-effectus,
then~$(\Scal C)^{\mathsf{op}}$ (the scalars with multiplication in the
    opposite direction)
    is an effect divisoid.
\item
If~$M$ is a division effect monoid as in \cite[dfn.~6.3]{kentapartial},
then~$M^{\textrm{op}}$ is an effect divisoid in the obvious way.
\item
The effect monoid on the unit interval of~$C[0,1]$
    is \emph{not} an effect divisoid,
    but the effect monoid on the unit interval of~$L^\infty[0,1]$
    is, see \sref{basic-divisoid-equiv}.
\end{enumerate}
\spacingfix{}
\end{point}
\begin{point}{60}[basic-divisoid-equiv]{Exercise*}%
Let~$X$ be a compact Hausdorff space.
In this exercise you will show
    that the unit interval of~$C(X)$
    is an effect divisoid
    if and only if~$X$
    is \Define{basically disconnected}
    \index{basically disconnected}
    \cite[1H]{gillman2013rings}
    --- that is: if~$\overline{\supp f}$ is open
    for every~$f \in C(X)$.
Equivalently:
$C(X)$ is~$\sigma$-Dedekind complete
    \cite[3N.5]{gillman2013rings},
    which is in turn equivalent to bounded~$\omega$-completeness
    (i.e.~$C(X)$ has suprema, resp.~infima, of
            bounded ascending, resp.~descending, sequences).
\begin{enumerate}
\item
    Suppose the unit interval of~$C(X)$ is an effect divisoid.
        Let~$f \in C(X)$ with~$0 \leq f \leq 1$.
    If~$\overline{\supp f} = X$, then you're done.
        Otherwise pick any~$y \notin \overline{\supp f}$.
    Show, using Urysohn's lemma,
    that there is a~$g \in C(X), 0\leq g \leq \rfrac{f}{f}$
    with~$g(y) = 0$ and $g(x) = 1$ for all~$x \in \overline{\supp f}$.
    Show that this implies~$(\rfrac{f}{f})(y) = 0$
        and so~$\rfrac{f}{f}$ is the characteristic function
    of~$\overline{\supp f}$.  Conclude~$X$ is basically disconnected.
\item
Assume~$X$ is basically disconnected.
Let~$f,g\in C(X)$ with~$0 \leq f \leq g \leq 1$.
For~$n > 0$,
show that~$U_n \equiv \{ x;\ g(x) > \frac{1}{n} \}$
has open closure and so
\begin{equation*}
    h_n \ \equiv \  \begin{cases}
        \frac{f(x)}{g(x)} & x \in \overline{U_n}\\
        0 & \text{otherwise}.
    \end{cases}
\end{equation*}
is continuous.
Note~$h_1 \leq h_2 \leq \ldots \leq 1$
    and define~$\rfrac{f}{g} = \sup_n h_n$.
Show~$\rfrac{f}{f}$ is the characteristic function
    of~$\overline{\supp f}$ and
    $\rfrac{f}{g} \leq \rfrac{g}{g}$.
Prove that this partial division turns the unit interval of~$C(X)$
    into an effect divisoid.
\end{enumerate}
\spacingfix{}

\end{point}
\begin{point}{70}{Proposition}%
If $a \perp b$ and~$a \ovee b \leq c$ in an effect divisoid,
then~$\rfrac{a \ovee b}{c } = \rfrac{a}{c} \ovee \rfrac{b}{c}$.
\begin{point}{80}{Proof}%
We will first show that if~$c \odot a \perp c \odot b$,
    then~$\rfrac{c}{c} \odot a \perp \rfrac{c}{c}\odot b$.
To start,
    note~$c \odot b^\perp \leq c \odot b^\perp \ovee c^\perp
        = (c \odot b)^\perp \leq c \odot a$
        and so
\begin{equation*}
     \rfrac{c \odot a}{c} 
     \ = \ 
            \rfrac{c \odot b^\perp}{c} \odot
            \rfrac{c \odot a}{c \odot b^\perp}
     \ = \ 
            \rfrac{c}{c} \odot b^\perp \odot
            \rfrac{c \odot a}{c \odot b^\perp}
    \ \leq \ \rfrac{c}{c} \odot b^\perp.
\end{equation*}
Hence~$(\rfrac{c}{c}\odot b)^\perp
            \geq \rfrac{c}{c} \odot b^\perp
            \geq \rfrac{c \odot a}{c}
            = \rfrac{c}{c} \odot a $.
    Indeed:~$\rfrac{c}{c}\odot a \perp \rfrac{c}{c} \odot b$.
\begin{point}{90}%
Assume~$a\perp b$ and~$a\ovee b \leq c$. 
Clearly~$c \odot \rfrac{a}{c}  = a \perp b = c \odot \rfrac{b}{c}$.
By our initial lemma, we get
    $\rfrac{a}{c } =
    \rfrac{c}{c} \odot \rfrac{a}{c} \perp
    \rfrac{c}{c} \odot \rfrac{b}{c} =
    \rfrac{b}{c }
$.
Clearly~$c \odot (\rfrac{a}{c} \ovee \rfrac{b}{c}) = a \ovee b$, so
\begin{equation*}
    \rfrac{a \ovee b}{c} \ = \  
    \rfrac{c \odot (\rfrac{a}{c} \ovee \rfrac{b}{c})}{c} \ = \ 
\rfrac{c}{c} \odot (\rfrac{a}{c} \ovee \rfrac{b}{c}) \ = \ 
\rfrac{a}{c} \ovee \rfrac{b}{c},
\end{equation*}
as desired. \qed
\end{point}
\end{point}
\end{point}
\end{parsec}

\begin{parsec}{1960}%
\begin{point}{10}%
The following is a generalization
of~\cite[prop.~C.3]{kentapartial}.
See also~\cite[prop.~15]{statesofconvexsets}.
\end{point}
\begin{point}{20}[aconvm-is-effectus]{Theorem}%
If~$M$ is an effect divisoid,
    then~$\AConvM$ is an effectus.
\begin{point}{30}{Proof}%
We only have to show that the square on the left
    of \eqref{pullbacks} is a pullback in~$\AConvM$
    --- the other axioms were proven in \sref{aconvalmosteffectus}.
\begin{equation*}
\xymatrix{
    Z \ar@/^1pc/[rrd]^\beta
        \ar@/_1pc/[rdd]_\alpha
        \ar@{.>}[rd]|\gamma\\
        &X+Y \ar[r]^{\id+!} \ar[d]^{!+\id} & X+1\ar[d]^{!+\id} \\
        &1+Y\ar[r]_{\id+!} & 1+1
    }
\end{equation*}
To this end, assume~$\alpha\colon Z \to X+1$
and~$\beta\colon Z \to Y+1$
are maps in~$\AConvM$
with~$(! + \id) \after \beta = (\id+!)\after \alpha$.
We have to show there is a unique~$\gamma\colon Z \to X+Y$ in~$\AConvM$
with~$(\id+!)\after \gamma =\beta$
and~$(!+\id) \after \gamma =\alpha$.
Pick any~$z \in Z$.
Write~$\cdot$ for the unique element of~$1$.
By construction of the coproduct, see \sref{aconv-coprod},
we have
\begin{align*}
    \alpha(z) &\ =\  h_1\Bigl(\lambda_0 \ket{\kappa_1(\cdot)}
    \ovee \bigovee^n_{i=1} \lambda_i \ket{\kappa_2(y_i)}\Bigr) \\
        \beta(z) &\ =\  h_2 \Bigl(\sigma_0 \ket{\kappa_2(\cdot)}
        \ovee \bigovee^m_{i=1} \sigma_i \ket{\kappa_1(x_i)}\Bigr)
\end{align*}
for some~$\lambda_i, \sigma_i \in M$,
    $x_i \in X$ and $y_i \in Y$,
    where~$h_1$ and~$h_2$ are the convex structures
    on~$1+Y$ and~$X+1$, respectively.
Unfolding definitions
    (and identifying~$1+1 \cong M$ via~$\kappa_1 (\cdot) = 1$),
    it is easy to
    see that by assumption
\begin{equation}\label{aconv-epb-ass}
    \bigovee^m_{i=1} \sigma_i \ = \ 
    (!+\id)(\beta(z))  \ = \ 
    (\id + !)(\alpha(z)) \ =\  \lambda_0
    \ = \ \Bigl(\bigovee^n_{i=1} \lambda_i\Bigr)^\perp.
\end{equation}
This suggests the following definition
\begin{equation}\label{aconvmexpuldefgamma}
    \gamma(z) \ = \  h_3 \Bigl(
    \bigovee^n_{i=1} \lambda_i \ket{\kappa_2(y_i)}
        \ovee \bigovee^m_{i=1} \sigma_i \ket{\kappa_1(x_i)}\Bigr),
\end{equation}
where~$h_3$ is the convex structure on~$X+Y$.
To show this a proper definition,
let~$n',m'\in \N$, $\lambda'_i, \sigma'_i \in M$,
$x'_i \in X$ and $y'_i \in Y$ be any elements such that
\begin{align*}
    \alpha(z) &\ =\ 
    h_1\Bigl(\lambda'_0 \ket{\kappa_1(\cdot)}
    \ovee \bigovee^{n'}_{i=1} \lambda'_i \ket{\kappa_2(y'_i)}\Bigr) \\
        \beta(z) &\ =\  h_2 \Bigl(\sigma'_0 \ket{\kappa_2(\cdot)}
    \ovee \bigovee^{m'}_{i=1} \sigma'_i \ket{\kappa_1(x'_i)}\Bigr).
\end{align*}
\spacingfix{}
\begin{point}{40}%
We have to show
\begin{multline}\label{convexpulgoal}
    h_3 \Bigl(
    \bigovee^{n'}_{i=1} \lambda'_i \ket{\kappa_2(y'_i)}
    \,\ovee\, \bigovee^{m'}_{i=1} \sigma'_i \ket{\kappa_1(x'_i)}\Bigr) \\
    \ = \ 
    h_3 \Bigl(
    \bigovee^n_{i=1} \lambda_i \ket{\kappa_2(y_i)}
        \,\ovee\, \bigovee^m_{i=1} \sigma_i \ket{\kappa_1(x_i)}\Bigr).
\end{multline}
Without loss of generality,
    we may assume~$n,n',m,m' \geq 1$
    (as we allow, for instance,~$\lambda_1=0$).
Due to \eqref{aconv-epb-ass},
 it is clear~$\lambda_0 = \lambda_0'$ and~$\sigma_0=\sigma_0'$.
By \sref{elements-coprod-conv}
    we know there
    are derivations~$\Phi_1, \ldots, \Phi_l \in \mathcal{D}_M^2 (1+Y)$
    and~$\Psi_1, \ldots, \Psi_k \in \mathcal{D}_M^2 (X+1)$,
    say~$\Phi_i \equiv \bigovee_j \zeta_{ij} \ket{\varphi_{ij}}$,
    $\Psi_i \equiv \bigovee_j \xi_{ij} \ket{\psi_{ij}}$
    with
\begin{align*}
    \Phi_1 &\ = \ \ket{\,
    \lambda_0 \ket{\kappa_1(\cdot)}
\ovee \bigovee^{n}_{i=1} \lambda_i \ket{\kappa_2(y_i)} \,} &
    \Phi_l &\ = \ \ket{\,
    \lambda_0 \ket{\kappa_1(\cdot)}
\ovee \bigovee^{n'}_{i=1} \lambda'_i \ket{\kappa_2(y'_i)} \,} \\
    \Psi_1 &\ = \ \ket{\,
    \sigma_0 \ket{\kappa_2(\cdot)}
\ovee \bigovee^{m}_{i=1} \sigma_i \ket{\kappa_1(x_i)} \,} &
    \Psi_k &\ = \ \ket{\,
    \sigma_0 \ket{\kappa_2(\cdot)}
\ovee \bigovee^{m'}_{i=1} \sigma'_i \ket{\kappa_1(x'_i)} \,}.
\end{align*}
We are going to combine the~$\Phi_i$ and~$\Psi_i$
into a single derivation of~\eqref{convexpulgoal}
    by first `applying' the $\Phi_i$ and then
    the~$\Psi_i$.
Define~$\Omega_1, \ldots, \Omega_{l+k} \in \mathcal{D}_M^2 (X+Y)$ by
\begin{align*}
    \Omega_i  &\ = \ \bigovee_j\zeta_{ij} \ket{\omega_{ij}}
        &&\text{for $1 \leq i \leq l$}
    \\
    \Omega_{i+l}  &\ = \ \bigovee_j \xi_{ij} \ket{\omega_{(i+l)j}}
        &&\text{for $1 \leq i \leq k$,}
\end{align*}
where~$\omega_{ij} \in \mathcal{D}_M(X+Y)$
are given by
\begin{align*}
    \omega_{ij} (\kappa_2 (y)) &\ =\ \varphi_{ij}(\kappa_2(y))  \\
    \omega_{ij} (\kappa_1(x)) &\ = \ 
    \varphi_{ij}(\kappa_1(\cdot))\odot
    ( \rfrac{\psi_{11}(\kappa_1(x))}{\lambda_0}
    \ovee r(\kappa_1(x)))
    \\
    \omega_{(i+l)j} (\kappa_1(x)) &\ =\ \psi_{ij}(\kappa_1(x)) \\
    \omega_{(i+l)j} (\kappa_2 (y)) &\ =\ 
    \psi_{ij}(\kappa_2(\cdot)) \odot (
    \rfrac{\varphi_{l1}(\kappa_2(y))}{\sigma_0}
    \ovee r(\kappa_2(y))) \\
    r(z) & \ = \ 
    \begin{cases}
        (\rfrac{\lambda_0}{\lambda_0})^\perp & z = \kappa_1(x_1) \\
        (\rfrac{\sigma_0}{\sigma_0})^\perp & z = \kappa_2(y_1) \\
        0 & \text{otherwise},
    \end{cases}
\end{align*}
    where on the first two lines~$1 \leq i \leq l$
    and on the second two~$1 \leq i \leq k$.
Before we continue, we check whether the~$\omega_{ij}$
    are distributions, as claimed.
For~$1 \leq i \leq l$
    this amounts to $\bigovee_x \omega_{ij}(\kappa_1(x))
= \varphi_{ij}(\kappa_1(\cdot))$.
Indeed
\begin{align*}
\bigovee_x
    \omega_{ij}(\kappa_1(x))
    & \  =\  \varphi_{ij}(\kappa_1(\cdot)) \odot r(\kappa_1(x_1))
            \,\ovee\, \bigovee_x 
        \varphi_{ij}(\kappa_1(\cdot)) \odot
        \rfrac{\psi_{11}(\kappa_1(x))}{\lambda_0} \\
        & \  =\  \varphi_{ij}(\kappa_1(\cdot)) \odot (\rfrac{\lambda_0}{\lambda_0})^\perp
            \,\ovee \,
        \varphi_{ij}(\kappa_1(\cdot)) \odot
        \rfrac{\bigovee_x \psi_{11}(\kappa_1(x))}{\lambda_0} \\
        & \  =\  \varphi_{ij}(\kappa_1(\cdot)) \odot (\rfrac{\lambda_0}{\lambda_0})^\perp
            \,\ovee \,
        \varphi_{ij}(\kappa_1(\cdot)) \odot
        \rfrac{\lambda_0}{\lambda_0} \\
        & \  =\  \varphi_{ij}(\kappa_1(\cdot)).
\end{align*}
With an analogous argument, one checks~$\omega_{(i+l)j}$
    is a distribution for~$1 \leq i \leq k$.

To show~$\Omega_i$ is
    a derivation (in the sense of \sref{elements-coprod-conv}),
    first pick any~$1 \leq i < l$.
We distinguish between the two allowable derivation steps.
If~$ \mu(\Phi_i) = \mu(\Phi_{i+1}) $,
    then~$\mu(\Omega_i) = \mu(\Omega_{i+1})$
    by a straightforward computation
    and so~$\Omega_i$ has a valid step.
In the other case, we have~$\zeta_{ij} = \zeta_{(i+1)j}$ and
    for each~$j$ we have three possibilities.
In the first case, if~$\varphi_{ij} = \varphi_{(i+1)j}$,
    then clearly~$\omega_{ij} = \omega_{(i+1)j}$ as well.
    For the second case, assume~$\varphi_{ij} = (\mathcal{D}_M \kappa_2) (\chi)$
    and~$\varphi_{(i+1)j} = (\mathcal{D}_M \kappa_2) (\chi')$
    for some~$\chi,\chi' \in \mathcal{D}_M Y$
    with~$h_Y(\chi) = h_Y(\chi')$.
    Clearly~$\varphi_{ij}(\kappa_1(\cdot)) = 0$
    and so~$\omega_{ij} = (\mathcal{D}_M \kappa_2)(\chi)$.
    Similarly~$\omega_{(i+1)j} = (\mathcal{D}_M \kappa_2)(\chi')$.
    The third case is trivial.
    So~$\Omega_i$ also makes a valid step of the second kind.

With the same kind of argument, we cover
    $\Omega_{i+l}$ for~$1 \leq i < k$.
The only step left to check is the one
    from~$\Omega_l$ to~$\Omega_{l+1}$.
Note~$\zeta_{l1}=1$,~$\varphi_{l1}(\kappa_1(\cdot)) = \lambda_0$
and~$\lambda_0 \odot  (\rfrac{\lambda_0}{\lambda_0})^\perp = 0$,
so with some easy manipulation, we find
\begin{equation*}
    \Omega_l \ = \  \ket{
        \bigovee^{n'}_{i=1} \lambda'_i \ket{\kappa_2(y'_i)}
\, \ovee\,  \bigovee^m_{i=1} \sigma_i \ket{\kappa_1(x_i)} }.
\end{equation*}
and in a similar fashion
\begin{align*}
    \Omega_{l+1} & \ = \  \ket{
        \bigovee^{n'}_{i=1} \lambda'_i \ket{\kappa_2(y'_i)}
    \, \ovee\,  \bigovee^m_{i=1} \sigma_i \ket{\kappa_1(x_i)} } \\
    \Omega_{1} & \ = \  \ket{
        \bigovee^{n}_{i=1} \lambda_i \ket{\kappa_2(y_i)}
\,\ovee\, \bigovee^m_{i=1} \sigma_i \ket{\kappa_1(x_i)} } \\
    \Omega_{l+k} & \ = \  \ket{
        \bigovee^{n'}_{i=1} \lambda'_i \ket{\kappa_2(y'_i)}
    \,\ovee\, \bigovee^{m'}_{i=1} \sigma'_i \ket{\kappa_1(x'_i)} }.
\end{align*}
So~$\Omega_i$ is a valid derivation of \eqref{convexpulgoal}.
\end{point}
\begin{point}{50}%
We are in the home stretch now.
Define~$\gamma(z)$ as in \eqref{aconvmexpuldefgamma}.
It is easy to see~$\gamma$
    is the unique map with~$(\id+!) \after \gamma = \beta$
    and~$(!+\id)\after \gamma = \alpha$.
It remains to be shown~$\gamma$ is affine.
Write
\begin{equation*}
q_1\colon \mathcal{D}_M(X+Y) \to X+Y \quad
q_2\colon \mathcal{D}_M(1+Y) \to 1+Y \quad
q_3\colon \mathcal{D}_M(X+1) \to X+1
\end{equation*}
for the quotient maps used in the construction of the coproducts.
We worked hard to show
that
for all~$z \in Z$, $\varphi \in \mathcal{D}_M Y$,
$\psi \in \mathcal{D}_M X$, 
and~$\chi \in \mathcal{D}_M (X+Y)$
with~$\chi(\kappa_1 (y)) = \varphi(\kappa_1(y))$
and~$\chi(\kappa_2 (x)) = \psi(\kappa_2(x))$,
we have
\begin{equation*}
    \left.\begin{array}{ll}
    \alpha(z)\ =\ q_1(\varphi) \\
    \beta(z)\ =\ q_2(\psi)
    \end{array}\right] \quad \implies \quad
    \gamma(z) = q_3(\chi).
\end{equation*}
Assume~$z = h_Z(\bigovee_i \lambda_i \ket{z_i})$.
We want to show~$\gamma(z) = h_3(\bigovee_i \lambda_i \ket{\gamma(z_i)})$.
To this end, pick~$\varphi_i$, $\psi_i$, $\chi_i$
with~$\alpha(z_i) = q_1(\varphi_i)$,
$\beta(z_i) = q_2(\psi_i)$,
$\chi_i(\kappa_1(y)) = \varphi_i(\kappa_1(y))$
and $\chi_i(\kappa_2(x)) = \psi_i(\kappa_2(x))$.
So $\gamma(z_i) = q_3(\chi_i)$.
Define~$\varphi,\psi,\chi$ as follows.
\begin{align*}
\varphi &\ \equiv\  \mu\bigl(\bigovee_i \lambda_i \ket{\varphi_i}\bigr)&
\psi &\ \equiv\  \mu\bigl(\bigovee_i \lambda_i \ket{\psi_i}\bigr) &
\chi &\ \equiv\  \mu\bigl(\bigovee_i \lambda_i \ket{\chi_i}\bigr)
\end{align*}
Note
$q_1(\varphi)
    = h_1(\bigovee_i \lambda_i \ket{q_1(\varphi_i)})
    = h_1(\bigovee_i \lambda_i \ket{\alpha(z_i)})
    = \alpha(\bigovee_i \lambda_i \ket{z_i})
    = \alpha(z)$.
Similarly~$q_2(\psi) = \beta(z)$.
For any~$y \in Y$,
    we have
    $
    \chi(\kappa_1(y))
    = \bigovee_i \lambda_i \odot \chi_i(\kappa_1(y))
    = \bigovee_i \lambda_i \odot \varphi_i(\kappa_1(y))
    = \varphi(\kappa_1(y))
    $.
Similarly~$\chi(\kappa_2(x)) = \psi(\kappa_2(x))$ for any~$x \in X$.
So~$\gamma(z) = q_3(\chi)$ and consequently
\begin{equation*}
    \gamma(z) \ = \ q_3(\chi)
    \ = \ h_3\bigl(\bigovee_i \lambda_i\ket{q_3(\chi_i)} \bigr)
    \ = \ h_3\bigl(\bigovee_i \lambda_i\ket{\gamma(z_i)} \bigr),
\end{equation*}
as desired. \qed
\end{point}
\end{point}
\end{point}
\end{parsec}

\section{Effectuses with quotients}
\begin{parsec}{1970}%
\begin{point}{10}%
    We return to the main program
        of this chapter: the attempted axiomatisation
        of the effectus~$\op\vN$.
\end{point}
\begin{point}{20}[dfn-quotient]{Definition}%
Let~$C$ be an effectus.
We say~$C$ is an \Define{effectus with quotients}
    \index{effectus!with quotients}
    \index{*bra@$X/_p$}
    if for each predicate~$p \colon X \to 1$,
    there exists a map~$\xi_p \colon X \to X/_p$
    with~$1 \after \xi_p \leq p^\perp$
    satisfying the following universal property.
\begin{quote}
    For any map~$f\colon X \to Y$
        with~$1 \after f \leq p^\perp$,
        there is a unique map~$f' \colon X/_p \to Y$
        such that~$f' \after \xi_p = f$.
\end{quote}
Any  map with this universal property
    is called a \Define{quotient} for~$p$.
    \index{quotient}
\begin{point}{21}{Remarks}%
The universal property of quotients
    is essentially the same as
    the universal property of (contractive) \emph{compressions},
    which we proposed in~\cite{westerbaan2016universal}.
Effectuses with quotients also appear in~\cite{effintro}.
\end{point}
\begin{point}{30}[quot-not]{Notation}%
    Unless otherwise specified~\Define{$\xi_p$}\index{*xi@$\xi_p$}
    will denote some quotient for~$p$.
\end{point}
\end{point}
\begin{point}{40}{Examples}%
In~$\op{\vN}$ quotients are exactly the same thing
    as contractive filters,
    see~\sref{dils-def-filter}, \sref{filter}
    and~\sref{canonical-filter}.
An example of a quotient map for~$1-a \in \scrA$ is given
    by~$\xi\colon \ceil{a}\scrA \ceil{a} \to \scrA$
    with~$\xi(b) = \sqrt{a} b \sqrt{a}$.
\begin{point}{41}%
The effectus~$\op\EJA$ has quotients,
    which are in essence very similar to those of~$\op\vN$,
    when ignoring the technical complications Jordan algebras
    impose. \cite[prop.~25]{eja}
\end{point}
\begin{point}{42}[opous-quotients]%
The effectus~$\op\OUS$ has quotients,
    which are rather different from the previous two examples.
An example of a  quotient for a predicate~$v$ on an order unit space~$V$
    is given by the inclusion of the order ideal generated by~$1-v$
\begin{equation*}
    \langle 1-v\rangle \ \equiv \ \{
        a; \ a\in V; \ \exists n.\, -n (1-v) \leq a \leq n (1-v) \}
    \ \subseteq\  V,
\end{equation*}
where the order ideal is considered as an order unit space
    with order-unit~$1-v$.
    The effectus~$\op\OUG$ has very similar quotients. See~\cite{effintro}.
\end{point}
\end{point}
\begin{point}{50}[quotient-basics]{Exercise}%
Verify the following basic properties of quotients.
    (See~\cite{effintro}.)
\begin{enumerate}
    \item If~$\xi\colon X \to Y$ is a quotient for~$p$
                and~$\vartheta\colon Y \to Z$ is an isomorphism,
                then~$\vartheta \after \xi$ is a quotient for~$p$
                as well.
    \item Conversely, if~$\xi_1$ and~$\xi_2$
            are both quotients for~$p$,
            then there is a unique isomorphism~$\vartheta$
            with~$\xi_1 = \vartheta \after \xi_2$.
    \item Isomorphisms are quotients (for 0).
    \item Maps into~$0$ are quotients (for 1).
    \item If~$\xi$ is a quotient for~$p$, then~$1\after \xi = p^\perp$
                (Hint: apply the universal property to $p^\perp$).
    \item Quotients are epic.
\end{enumerate}
\spacingfix{}
\begin{point}{60}%
The following proposition is easy to prove,
    but shows an important property of quotients:
    any map~$f$ factors as a total map after a quotient
    for~$(1\after f)^\perp$.
\end{point}
\end{point}
\begin{point}{70}[quotient-total]{Proposition}%
Assume~$\xi_{p^\perp} \colon X \to X/_{p^\perp}$ is a quotient for~$p^\perp$.
For any~$f\colon X \to Z$
    with~$1 \after f = p$,
    there is a unique \emph{total} $g\colon X/_{p^\perp} \to Z$
    with~$f = g \after \xi_{p^\perp}$.
\begin{point}{80}{Proof}%
(This is a simplified version of our previous
    proof in~\cite{effintro}.)
By definition of quotient, there is a unique~$g\colon X/_{p^\perp} \to Z$
    with~$g \after \xi_{p^\perp} = f$.
Note~$1 \after g \after \xi_{p^\perp} = 1 \after f = p = 1 \after \xi_{p^\perp}$.
Thus, as~$\xi_{p^\perp}$ is epi (\sref{quotient-basics}),
    we conclude~$1 \after g = 1$.
That is: $g$ is total. \qed
\end{point}
\end{point}
\begin{point}{90}[quotients-composition]{Proposition}%
In an effectus with quotients,
    quotients are closed under composition.
\begin{point}{100}{Proof}%
(This is essentially the same proof as we gave in~\cite{effintro}.)
Assume~$\xi_1\colon X \to Y$ is a quotient for~$p^\perp$
    and~$\xi_2\colon Y \to Z$ is a quotient for~$q^\perp$.
We will prove~$\xi_2 \after \xi_1$
    is a quotient for~$(q \after \xi_1)^\perp$.
As our effectus has quotients,
    we can pick a quotient~$\xi \colon X \to X/_{(q \after \xi_1)^\perp}$
        of~$(q \after \xi_1)^\perp$.
First, some preparation.
Note~$1 \after \xi = q \after \xi_1 \leq 1 \after \xi_1 = p$.
Thus by the universal property of~$\xi_1$,
there is a unique map~$h_1\colon Y \to X/_{(q \after \xi_1)^\perp}$
        with~$h_1 \after \xi_1 = \xi$.
As~$1 \after h_1 \after \xi_1 = 1 \after \xi = q \after \xi_1$
    and~$\xi_1$ is epi, we see~$1 \after h_1 = q$.
Thus by \sref{quotient-total},
there is a unique total map~$h_2 \colon Z \to X/_{(q \after \xi_1)^\perp}$
    with~$h_2 \after \xi_2 = h_1$.
Let~$g\colon X/_{(q \after \xi_1)^\perp} \to Z$
    be the unique map such that~$g \after \xi = \xi_2 \after \xi_1$.
We are in the following situation.
\begin{equation*}
    \xymatrix@C+2pc{
        X  \ar[r]^{\xi_1} \ar@/_/[rrd]_{\xi}
        & Y \ar[r]^{\xi_2} \ar@{.>}[rd]_{h_1}
        & Z \ar@{.>}@/^/[d]^{h_2} \\
        && X/_{(q \after \xi_1)^\perp} \ar@{.>}@/^/[u]^{g}
    }
\end{equation*}
We claim~$g$ and~$h_2$ are each other's inverse.
Indeed: from $g \after h_2 \after \xi_2 \after \xi_1= g \after \xi = \xi_2 \after \xi_1$
    we get~$g \after h_2 = \id$
    and from~$h_2 \after g \after \xi = h_2 \after \xi_2 \after \xi_1 = \xi$
    we find~$h_2 \after g = \id$.
Thus~$\xi_2 \after \xi_1 = g \after \xi$ for an isomorphism~$g$,
which shows~$\xi_2 \after \xi_1$ is a quotient, see \sref{quotient-basics}. \qed
\end{point}
\end{point}
\begin{point}{110}[quot-fact-system]{Exercise}
By \sref{quotient-total} we know that each
    map~$f$ in an effectus with quotients factors as~$t \after \xi$
    for some total map~$t$ and quotient~$\xi$.
Show that this forms an orthogonal factorization system (cf.~\cite{korostenski1993factorization}) ---
    that is: prove that if~$t' \after \xi' = t \after \xi$
    for some quotient~$\xi'$ and total map~$t'$,
    then there is a unique isomorphism~$\vartheta$
    with~$\xi' = \vartheta \after \xi$
    and~$t = t' \after \vartheta$.
\end{point}
\end{parsec}

\begin{parsec}{1980}%
\begin{point}{10}%
In \sref{quot-fact-system} we saw
    that every  map~$f$
    factors uniquely via a quotient for~$(1 \after f)^\perp$.
In this sense quotients are the universal way to restrict to the
    support on the domain --- so why do we call them quotients
    instead of, say, initial-support-maps?
We answer this question in~\sref{bartequivquotrem},
    but need some preparation first.
\end{point}
\begin{point}{20}[dfn-eff-grothendieck]{Definition}%
For an effectus~$C$,
    write~$\Define{\int \Pred_\square}$
    \index{$\int \Pred_\square$}
    for the category with
\begin{enumerate}
\item as objects pairs~$(X,p)$,
        where~$X$ is an object of~$C$
        and~$p \in \Pred X$ and
\item
    an morphism between~$(X,p)$ and~$(Y,q)$
    corresponds to a partial map~$f\colon X \to Y$ in~$C$
       for which we have~$p \leq (q^\perp \after f)^\perp$.
\end{enumerate}
Write~$U\colon \int \Pred_\square \to C$
for the forgetful functor --- that is:~$U(X,p) = X$
and~$Uf = f$ for arrows~$f$.
The functor~$U$ has left- and right-adjoint
    $0 \dashv U \dashv 1$
    with
    $0,1 \colon C \to \int \Pred_\square$
    given by~$0X = (X,0)$, $1X = (X,1)$
    and~$0f=1f=f$.
\end{point}
\begin{point}{30}[exc-quot-adjoint]{Exercise*}%
Show that for an effectus~$C$, the following are equivalent.
\begin{enumerate}
\item $C$ has quotients.
\item The functor~$0\colon C \to \int \Pred_\square$,
        has a left-adjoint~$Q \colon \int\Pred_\square \to C$.
\end{enumerate}
\spacingfix{}
\end{point}
\begin{point}{40}[bartequivquotrem]%
If the existence of a left adjoint~$Q$ to~$0$
    gives quotients,
    what does the existence of a right adjoint~$K$ to~$1$ amount to?
    It will turn out (\sref{compr-grothendieck})
    this is equivalent to the existence of
        so-called \emph{comprehensions},
        the topic of the next section.
With both quotients and comprehension, we have a chain of four adjunctions:
\begin{equation*}
   \vcenter{\xymatrix{
\int \Pred_\square \ar[d]_{\dashv\;}^{\;\dashv}
   \ar@/_6ex/[d]^{\;\dashv}_{Q} 
   \ar@/^6ex/[d]_{\dashv\;}^{K} \\
   C\ar@/^3ex/[u]^(0.4){0\!}\ar@/_3ex/[u]_(0.4){\!1}
}} 
\end{equation*}
There are numerous other examples
    of categories which also
    have such a chain of four adjunctions
    between it and a natural category of object-with-predicates.
    See \cite{cho2015quotient} for several examples.
    In most of these examples, the predicates correspond to subspaces
        and the left-most functor to quotienting by this subspace.
    This is the reason for the name `effectus with quotients'.
\begin{point}{50}%
Originally, we considered the universal properties of
    quotients and comprehensions separately and under different names.
It was Jacobs who recognized that both universal properties
    appear together in a chain of adjunctions.
\end{point}
\end{point}
\end{parsec}

\section{Effectuses with comprehension and images}
\begin{parsec}{1990}%
\begin{point}{10}%
We continue with the formal introduction
        of our second axiom:
        the existence of~\emph{comprehension}.
\end{point}
\begin{point}{20}[dfn-comprehension]{Definition}%
Let~$C$ be an effectus.
We say~$C$ \Define{has comprehension}
    \index{effectus!with comprehension}%
    \index{*acc@$\cmpr{X}{p}$}%
    if for each predicate~$p \colon X \to 1$,
    there exists a map~$\pi_p \colon \cmpr{X}{p} \to X$
    with~$p \after \pi_p = 1 \after \pi_p$
    satisfying the following universal property.
\begin{quote}
    For any map~$g\colon Z \to X$
        with~$p \after g = 1 \after g$,
        there is a unique map~$g' \colon Z \to \cmpr{X}{p}$
        such that~$\pi_p \after g' = g$.
\end{quote}
Any  map with this universal property
    is called a \Define{comprehension} for~$p$.
    \index{comprehension}
\begin{point}{30}{Beware}%
    We do not assume comprehensions are total
    (in contrast to~\cite{effintro}.)
    In~\sref{compr-total} we will see that
    in an effectus with quotients, comprehensions must be total.
\end{point}
\begin{point}{31}{Remarks}%
The universal property of comprehensions is essentially the same
    as the universal property of (contractive) \emph{corners},
    which we proposed in~\cite{westerbaan2016universal}.
Effectuses with compression also appear in~\cite{effintro}.
In partial form, comprehensions are
    essentially the same
    as categorical kernels, see~\sref{compr-is-kernel} later on.
\end{point}
\begin{point}{40}[compr-not]{Notation}%
    Unless otherwise specified~\Define{$\pi_p$}\index{*pi@$\pi_p$}
        will denote some comprehension of~$p$.
\end{point}
\end{point}
\begin{point}{50}{Examples}%
In~$\op{\vN}$ comprehensions are exactly the same thing
    as corners, see \sref{dils-corner} or \sref{corner}.
An example of a comprehension for~$a \in \scrA$
    is given by~$\pi\colon \scrA \to \floor{a}\scrA\floor{a}$
    with~$\pi(b) = \floor{a}b\floor{a}$
    as proven in \sref{prop-corner}.
\begin{point}{51}%
The effectus~$\op\EJA$ also has comprehension,
    which is very similar to those of~$\op\vN$.~\cite[prop.~24]{eja}

The comprehension maps of~$\op\OUS$ are again quite different
    from those of~$\op\vN$ and~$\op\EJA$.
Confusingly, a comprehension for a predicate~$v$
    of an order unit space~$V$
    is given by the vector-space-quotient
    map~$q\colon V \to V/_{\langle1-v\rangle}$,
    where~$\langle 1-v \rangle$ is the order ideal generated
    by~$1-v$ defined in~\sref{opous-quotients}.
The effectus~$\op\OUG$ has similar comprehension.
    See~\cite{effintro}.

Any extensive category with final object has comprehension,
    which includes~$\SET$, $\bCH$ and~$\CRng$,
    see~\cite{effintro}.
\end{point}
\end{point}
\begin{point}{60}[compr-grothendieck]{Exercise*}%
    Show that an effectus~$C$ has comprehension if and only
        if the functor~$1 \colon C \to \int \Pred_\square$
        from \sref{dfn-eff-grothendieck}
        has a right adjoint~$K$.
\end{point}
\begin{point}{70}[compr-basics]{Exercise}%
Show the following basic properties of comprehensions.
    (See~\cite{effintro}.)
\begin{enumerate}
    \item If~$\pi\colon X \to Y$ is a comprehension for~$p$
                and~$\vartheta\colon Z \to X$ is an isomorphism,
                then~$\pi \after \vartheta$ is a comprehension for~$p$
                as well.
    \item Conversely, if~$\pi_1$ and~$\pi_2$
            are both comprehensions for~$p$,
            then there is a unique isomorphism~$\vartheta$
            with~$\pi_1 = \pi_2 \after \vartheta$.
    \item Isomorphisms are comprehensions (for 1).
    \item Zero maps are comprehensions (for 0).
    \item Comprehensions are monic.
    \item $p^\perp \after \pi = 0$ if~$\pi$ is a comprehension for~$p$.
\end{enumerate}
\end{point}
\end{parsec}
\spacingfix{}
\begin{parsec}{2000}%
\begin{point}{10}%
Comprehensions and categorical kernels are very similar.
\end{point}
\begin{point}{20}{Definition}%
Let~$C$ be a category with a zero object.
    A \Define{(categorical) kernel}~of an arrow~$f$
    \index{kernel!category theoretic}
    (in symbols: $\Define{\ker f}$)
    \index{kerf@$\ker f$}
    is an equalizer of~$f$ with the parallel zero map.
(That is: $f \after (\ker f) = 0$
    and for every~$g$ with~$f \after g = 0$,
    there exists a unique~$g'$ with~$(\ker f) \after g' = g$.)
We will use \Define{$\cok f$}
    \index{$\cok f$}
    to denote a cokernel of~$f$
    --- that is: a kernel of~$f$ in~$\op{C}$.
\end{point}
\begin{point}{30}[effectus-kernels]{Proposition}%
An effectus with comprehension has all kernels.
    The kernel of a map~$f$ is given by
        a comprehension~$\pi_{(1 \after f)^\perp}$
        for~$(1 \after f)^\perp$.
\begin{point}{40}{Proof}%
Clearly~$1 \after f \after \pi_{(1 \after f)^\perp} = 0$
    and so~$f \after \pi_{(1 \after f)^\perp} = 0$.
Assume~$g$ is any map with~$f \after g = 0$.
Then~$1 \after f \after g = 0$
    and so~$(1 \after f)^\perp \after g = 1 \after g$.
Thus there exists a unique~$g'$
    with~$\pi_{(1 \after f)^\perp} \after g' = g$. \qed
\end{point}
\end{point}
\begin{point}{50}[compr-is-kernel]{Exercise}%
Show that in an effectus,
    a map~$f$ is a comprehension for~$p$
    if and only if it is a kernel of~$p^\perp$.
\end{point}
\end{parsec}

\begin{parsec}{2010}%
\begin{point}{10}%
We are ready to define purity in effectuses.
\end{point}
\begin{point}{20}{Definition}%
In an effectus a map~$f$ is called \Define{pure}
    \index{pure!in an effectus}
    if~$f = \pi \after \xi$
    for some comprehension~$\pi$
    and quotient~$\xi$.
\end{point}
\begin{point}{30}{Example}%
Due to \sref{pure-fundamental}
    pure maps in~$\op\vN$
    are precisely the pure maps
    as defined in~\sref{dils-def-pure}
    and~\sref{pure}.

The pure maps~$\scrB(\scrH) \to \scrB(\scrK)$
    are exactly the maps of the form~$\ad_T$
    where~$T$ is a contractive map~$\scrK \to \scrH$.
\end{point}
\begin{point}{40}{Remark}%
In the previous chapter
    we discussed (in \sref{dils-pure-discussion})
    why some more familiar notions of purity for states,
    like extremality,
    do not generalize properly to arbitrary maps.
Recently, three other definitions of purity
        have been proposed in different contexts.
        \cite{cunningham2017purity,chiribella2014distinguishability,selby2017leaks}
All three are inspired by the essential uniqueness
    of purification (see \sref{ess-uniq-pur})
    and require a monoidal structure to state.
The identity map might not be pure
    in the sense of~\cite{chiribella2014distinguishability},
    which is too restrictive for our taste.
Hazarding a guess,
    it seems likely
    that~\cite{cunningham2017purity} considers maps pure which we don't
     and that \cite{selby2017leaks} does not contain all our pure maps.
The exact relation between
 our definition of pure
    and~\cite{cunningham2017purity,selby2017leaks}
    is left open.
\end{point}
\begin{point}{50}%
We do not have a good handle on the structure of pure maps in
    an arbitrary effectus.  We will need a few additional assumptions.
\end{point}
\end{parsec}

\begin{parsec}{2020}%
\begin{point}{10}{Definition}%
Let~$C$ be an effectus.
\begin{enumerate}
\item
    We say~$C$ has \Define{images}\index{effectus!with images}
    if for each map~$f\colon X \to Y$,
    there is a least predicate~$\Define{\IM f}$ on~$Y$
        \index{im@$\IM$, $\IMperp$}
        with the property~$(\IM f) \after f = 1 \after f$ ---
        that is, 
    for every predicate~$p$ on~$Y$
    with~$p \after f = 1 \after f$,
    we must have~$\IM f \leq p$.
For brevity, write~$\Define{\IMperp f} \equiv (\IM f)^\perp$. 
Note~$\IMperp f$ is the greatest predicate
with the property~$(\IMperp f) \after f = 0$.
\item
We say a map~$f\colon X \to Y$ is \Define{faithful}
    \index{faithful}
    if~$\IM f = 1$.
That is: $f$ is faithful if and only if
    $p \after f = 0$ implies~$p = 0$
    for every predicate~$p$.
\end{enumerate}
\spacingfix{}
\begin{point}{20}{Notation}%
The expression~$\IM f \after g$
    is read as~$\IM (f \after g)$.
\end{point}
\begin{point}{30}{Remarks}%
The predicate~$\IM f$ will play for comprehension the same role
    as~$1 \after f$ plays for quotients.
Similarly faithful is the analogue of total.
Images in effectuses are also studied in \cite{effintro}.
\end{point}
\end{point}
\begin{point}{40}{Examples}%
In $\op{\vN}$ the image of a map~$f\colon \scrA \to \scrB$
    is given by the least projection~$e \in \scrA$
    with the property~$f(1-e) = 0$.
The map~$f$ is faithful if and only
    if~$f(a^*a)=0$ implies~$a^*a = 0$
    for all~$a\in \scrA$.
\begin{point}{41}%
The effectus~$\op\EJA$ also has all images,
    which are defined similarly, see~\cite[prop.~29]{eja}.
In contrast, the effectuses~$\op\OUG$ and~$\op\OUS$
    do not have all images, see~\cite{effintro}.
\end{point}
\end{point}
\begin{point}{50}[im-ineq]{Exercise}%
Show~$\IM f \after g \leq \IM f$.
Conclude~$\IM f \after \alpha = \IM f$
    for any iso~$\alpha$.
\end{point}
\begin{point}{60}[exc-quot-faithful]{Exercise}%
Show that in an effectus, any quotient is faithful.
\end{point}
\begin{point}{70}%
In contrast to quotients,
    it is not clear at all whether (without additional assumptions)
    comprehensions are total;
    they are part of a factorization system
    or are closed under composition.
\end{point}
\begin{point}{80}[compr-total]{Lemma}%
In an effectus with quotients, comprehensions are total.
\begin{point}{90}{Proof}%
Assume~$\pi$ is some comprehension for~$p$.
Let~$\xi$ be a quotient for~$(1 \after \pi)^\perp$.
There is a total~$\pi_t$ with~$\pi = \pi_t \after \xi$.
As~$ 1 \after \pi_t \after \xi
        = 1 \after \pi
        = p \after \pi
        = p \after \pi_t \after \xi$
        and~$\xi$ is epi,
        we see~$1 \after \pi_t = p \after \pi_t$.
Thus, as~$\pi$ is a comprehension for~$p$,
    there exists an~$f$ with~$\pi_t = \pi \after f$.
Now~$\pi = \pi_t \after \xi = \pi \after f \after \xi$
    and so~$\id = f \after \xi$, since~$\pi$ is mono.
Hence~$1 = 1 \after \id =1 \after f \after \xi \leq 1 \after \xi = 1 \after \pi$
    and so~$\pi$ is total. \qed
\end{point}
\end{point}
\end{parsec}
\subsection{Sharp predicates}
\begin{parsec}{2030}%
\begin{point}{10}{Definition}%
Let~$C$ be an effectus with comprehension and images.
\begin{enumerate}
\item
    We say a predicate~$p$ is (image) \Define{sharp} \index{sharp}
        if~$p = \IM f$ for some map~$f$.
        Write~$\Define{\SPred X}$ for the set of all sharp predicates on~$X$.
        \index{SPred@$\SPred X$}
\item
    We define~$\Define{\floor{p}} = \IM \pi_p$, \index{*floor@$\floor{\ }$, floor}
    where~$\pi_p$ is some comprehension for~$p$.
(Comprehensions for the same predicate have the same image
    by \sref{im-ineq}.)
Also write~$\Define{\ceil{p}} = \floor{p^\perp}^\perp$.
    \index{*floorc@$\ceil{\ }$, ceiling}
\end{enumerate}
\spacingfix{}
\begin{point}{20}{Beware}%
In \cite{effintro} $p$ is called sharp whenever~$p \wedge p^\perp=0$.
In general this is weaker than image-sharpness, which we use in this text.
In~\sref{floor-basics}
    we will see~$\floor{p}$ is sharp.
It is unclear whether~$\ceil{p}$ is sharp (without additional
    assumptions).
\end{point}
\end{point}
\begin{point}{30}{Example}%
In $\op{\vN}$ the sharp predicates are the projections.
For a predicate~$a \in \scrA$,
    the sharp predicate~$\ceil{a}$
    is the least projection above~$a$.
\end{point}
\begin{point}{40}[floor-basics]{Lemma}%
In an effectus with comprehension and images, we have
\begin{multicols}{2}
\begin{enumerate}
\item
    $\floor{p} \leq p$
\item
    $\pi_p = \pi_{\floor{p}} \after \alpha$
        for some iso~$\alpha$;
\item
    $\floor{\floor{p}} = \floor{p}$;
\item
    $p \leq q$ $\implies$ $\floor{p} \leq \floor{q}$;
\item
    $\ceil{p} \after f \leq \ceil{p \after f}$ \emph{and}
\item
    $\ceil{p} \after f =0$ iff~$p \after f = 0$,
\end{enumerate}
\end{multicols}
\noindent for any map~$f\colon X \to Y$,
    predicates~$p,q$ on~$X$ and
$\pi_p$ and~$\pi_{\floor{p}}$
    comprehensions for~$p$ and~$\floor{p}$ respectively.
\begin{point}{50}{Proof}%
We will prove the statements in listed order.
\begin{point}{60}{Ad 1}%
    Let~$\pi$ be a comprehension for~$p$.
    By definition~$p \after \pi = 1 \after \pi$.
    Thus~$\floor{p} = \IM \pi \leq p$, as desired.
\end{point}
\begin{point}{70}{Ad 2}%
It is sufficient to show~$\pi_p$ is a comprehension for~$\floor{p}$.
First, note~$\floor{p}\after \pi_p = (\IM \pi_p) \after \pi_p = 1 \after \pi_p$.
To show the universal property,
    assume~$g\colon Z \to X$
    is some map with~$\floor{p} \after g = 1 \after g$.
Then~$1 \after g = \floor{p} \after g \leq p \after g \leq 1 \after g$
    and so~$1 \after g = p \after g$.
As~$\pi_p$ is a comprehension for~$p$,
    there is a unique~$g'$ with~$\pi_p \after  g' = g$
    and so~$\pi_p$ is indeed a comprehension for~$\floor{p}$ as well.
\end{point}
\begin{point}{80}{Ad 3}%
Follows from the previous point and \sref{im-ineq}.
\end{point}
\begin{point}{90}{Ad 4}%
Pick a comprehension~$\pi_p$ for~$p$ and~$\pi_q$ for~$q$.
Note~$1 \after \pi_p = p \after \pi_p
                \leq q \after \pi_p \leq 1 \after \pi_p $
so~$q \after \pi_p = 1 \after \pi_p$
and thus~$\pi_p = \pi_q \after f$ for some~$f$.
By~\sref{im-ineq}
    we see~$\floor{p} = \IM \pi_p = \IM \pi_q \after f \leq \IM \pi_q = \floor{q}$.
\end{point}
\begin{point}{100}{Ad 5}%
Clearly~$p \after f \after \pi_{(p \after f)^\perp} = 0$.
(Recall our convention \sref{compr-not} for $\pi$'s.)
Thus there is some~$h$
with~$f \after \pi_{(p \after f)^\perp} = \pi_{p^\perp}\after h$.
By point 2 there is some (isomorphism)~$\alpha$
    with~$\pi_{p^\perp} =\pi_{\floor{p^\perp}} \after \alpha$.
    We compute
\begin{equation*}
    \ceil{p} \after f \after \pi_{(p \after f)^\perp}
        = \ceil{p} \after \pi_{p^\perp}\after h
        = \ceil{p} \after \pi_{\floor{p^\perp}}\after \alpha\after h
    = \ceil{p} \after \pi_{\ceil{p}^\perp}\after \alpha\after h
    = 0.
\end{equation*}
Thus~$\ceil{p} \after f \leq \IMperp \pi_{(p \after f)^\perp}
= \lfloor(p \after f)^\perp\rfloor^\perp = \ceil{p \after f}$, as promised.
\end{point}
\begin{point}{110}{Ad 6}%
Assume~$\ceil{p}\after f = 0$.
From~$p \leq \ceil{p}$
it follows~$p \after f \leq \ceil{p} \after f = 0$.
Thus~$p \after f = 0$.
For the converse, assume~$p \after f = 0$.
Then~$\ceil{p \after f} = \ceil{0} = (\IM \id)^\perp = 0$
    as~$\id$ is a comprehension for~$1$.
    Thus~$\ceil{p} \after f \leq \ceil{p \after f} = 0$. \qed
\end{point}
\end{point}
\end{point}
\begin{point}{120}[img-of-compr]{Exercise}%
Let~$C$ be an effectus with comprehension and images.
Show that~$p$ is sharp if and only if~$\floor{p}= p$.
Conclude~$\IM \pi_s = s$ for sharp~$s$.
\end{point}
\begin{point}{130}[ceiling-within-ceiling]{Exercise}%
Show that in an effectus with comprehension and images
    we have the equality
    $\ceil{\ceil{p}\after f} = \ceil{p \after f}$
    for any map~$f\colon X \to Y$ and predicate~$p$ on~$X$.
\end{point}
\begin{point}{140}[img-tupling]{Exercise}%
Show that in an effectus with
images~$\IM \left<f,g \right> = [\IM f, \IM g]$.
Conclude that a predicate~$[p,q]$ is sharp
    if and only if~$p$ and~$q$ are sharp.
\end{point}
\end{parsec}

\begin{parsec}{2040}%
\begin{point}{10}[compr-is-full]{Lemma}%
In an effectus with comprehension and images
    we have
\begin{equation*}
    s \leq t
    \quad \iff \quad
    \pi_s = \pi_t \after h
    \quad \text{for some $h$},
\end{equation*}
for all sharp predicates~$s,t$ on the same object.
\begin{point}{20}{Proof}%
(This simpler proof than ours in~\cite{effintro}.)
Assume~$\pi_s = \pi_t \after h$.
Then
\begin{equation*}
    s \ \overset{\sref{img-of-compr}}{=}\  \floor{s} \ =\  \IM \pi_s \  = \ \IM \pi_t \after h
    \ \overset{\sref{im-ineq}}{\leq} \ \IM \pi_t = \floor{t} 
            \ \overset{\sref{img-of-compr}}{=} \ t.
\end{equation*}
as desired.
Conversely assume~$s \leq t$.
Then~$t^\perp \after \pi_s \leq s^\perp \after \pi_s = 0$
    and so~$\pi_s = \pi_t \after h$
    for some~$h$ by the universal property of~$\pi_t$. \qed
\end{point}
\end{point}
\begin{point}{30}{Lemma}%
In an effectus with images,
    we have~$\IM [f,g] = (\IM f) \vee (\IM g)$.
\begin{point}{40}{Proof}%
We get~$\IM[f,g] \geq \IM f$ and~$\IM[f,g] \geq \IM  g$
    from
\begin{align*}
    [1 \after f, 1 \after g] 
        &\ =\  1 \after [f,g]  \\ 
        &\ =\ (\IM [f,g]) \after [f,g] \\
        & \ = \  [(\IM [f,g]) \after f, (\IM [f,g]) \after g].
\end{align*}
To show~$\IM [f,g]$ is the least upper-bound,
assume~$p \geq \IM f$ and~$p \geq \IM g$.
Then
\begin{equation*}
    1 \after [f,g] \ \geq\ p \after [f,g] \ =\  [p \after f, p \after g]
                    \ \geq\  [(\IM f) \after f, (\IM g)\after g]
                    \ = \ 1 \after [f,g],
\end{equation*}
    hence~$1 \after [f,g] = p \after [f,g]$
    and so~$p \geq \IM [f,g]$, as desired. \qed
\end{point}
\begin{point}{50}[lattice-compr]{Corollary}%
In an effectus with comprehension and images,
we have for any sharp~$s,t$
    a supremum (among all predicates)
    given by~$s \vee t = \IM[\pi_s, \pi_t]$.
    In particular:~$s \vee t$ is sharp.
\end{point}
\end{point}
\end{parsec}

\begin{parsec}{2050}%
\begin{point}{10}%
In \sref{compr-is-kernel} we saw that comprehensions
    are precisely kernels of predicates.
What about cokernels and quotients?
\end{point}
\begin{point}{20}[effectus-cokernels]{Proposition}%
An effectus with quotients and images has all cokernels.
A cokernel of a map~$f$
    is given by a quotient~$\xi_{\IM f}$ of~$\IM f$.
\begin{point}{30}{Proof}%
As~$0 = (\IMperp f) \after f = 1 \after \xi_{\IM f} \after f$,
        we have~$\xi_{\IM f} \after f = 0$.
    Assume~$g$ is a map with~$g \after f = 0$.
Then~$1 \after g \leq \IMperp f$
    and so~$g = g' \after \xi_{\IM f}$
    for a unique~$g'$.
Thus indeed, $\xi_{\IM f}$ is a cokernel of~$f$. \qed
\end{point}
\end{point}
\begin{point}{40}[exc-cokernels]{Exercise}%
Show that in an effectus with comprehension and images,
    a map~$f$ is a quotient of sharp~$s$
    if and only if~$f$ is a cokernel
    of a comprehension~$\pi_s$ of~$s$.
\end{point}
\end{parsec}

\section{\texorpdfstring{$\diamond$-effectuses}{%
                    diamond-effectuses}}
\begin{parsec}{2060}%
\begin{point}{10}%
In quantum physics it is not uncommon to restrict one's attention
    to sharp predicates (i.e.~projections).
Does this restriction hurt the expressivity?
It does not: every (normal) state on a von Neumann algebra
        is determined by its values on projections.
In fact, on~$\scrB(\scrH)$ with~$\dim \scrH \geq 3$,
    Gleason's famous theorem states that every
    measure on the projections extends uniquely to a state.
We take the idea of restricting oneself to projections
    one step further: we also want to restrict to projections
    for our outcomes.
It is rare for quantum processes to send projections to projections
    (see \sref{exa-sharp-vn}),
    so instead we consider the least projection above the outcome.
The simple idea of taking the restriction to sharp predicates
    seriously leads to a host of interesting new notions.
We will give these right off the bat and study their relevance later on.
\end{point}
\begin{point}{20}[diamond-basics]{Definition}%
    A~\Define{$\diamond$-effectus} (``diamond effectus'')
    \index{effectus!$\diamond$-}
    \index{*dia@$\diamond$, diamond!-effectus}
    is an effectus with quotients, comprehension and images
    such that~$s^\perp$ is sharp for every sharp predicate~$s$.
In a~$\diamond$-effectus,
    define for~$f\colon X \to Y$
    the following restrictions to sharp predicates.
    \begin{equation*}
        \xymatrix{
            \SPred X  \ar@/^/[r]^{f_\diamond}
            & \SPred Y \ar@/^/[l]^{f^\diamond}}
            \quad \text{by} \quad
            \Define{f^\diamond(s)} = \ceil{s \after f}
            \index{*par1@$(\ )_\diamond$, $(\ )^\diamond$}
            \quad
            \text{and}
            \quad
            \Define{f_\diamond(s)} = \IM f \after \pi_s
    \end{equation*}
    \begin{enumerate}
        \item
    We say maps~$f \colon X \leftrightarrows Y\colon g$
        are~\Define{$\diamond$-adjoint}
        if~$f^\diamond = g_\diamond$.
    \index{*dia@$\diamond$, diamond!-adjoint}
    \item
    An endomap~$f\colon X \to X$ is \Define{$\diamond$-self-adjoint}
        if~$f$ is $\diamond$-adjoint to itself.
    \index{*dia@$\diamond$, diamond!-self-adjoint}
    \item
    Two maps~$f,g\colon X \to Y$
        are~\Define{$\diamond$-equivalent}
        if~$f^\diamond = g^\diamond$
        (or equivalently whenever~$f_\diamond = g_\diamond$,
                see~\sref{diamond-equiv-equiv}.)
            \index{*dia@$\diamond$, diamond!-equivalent}
    \item
        A pure endomap~$f$ is~\Define{$\diamond$-positive}
            if $f = g\after g$ for some~$\diamond$-self-adjoint~$g$.
            \index{*dia@$\diamond$, diamond!-positive}
    \end{enumerate}
For brevity, write $\Define{f^\BOX(s)} = f^\diamond(s^\perp)^\perp$
    and~$\Define{f_\BOX(s)} = f_\diamond(s^\perp)^\perp$
    \index{*par2@$(\ )_\BOX$, $(\ )^\BOX$}
\end{point}
\begin{point}{30}{Examples}%
The categories~$\op\vN$, $\op\CvN$, $\op\EJA$
        and $\SET$ are all~$\diamond$-effectuses.
\end{point}
\end{parsec}
\begin{parsec}{2070}%
\begin{point}{10}%
Let's investigate the basic properties of~$(\ )^\diamond$, $(\ )_\diamond$
    and $(\ )^\BOX$.
\end{point}
\begin{point}{20}[exc-diam-order-pres]{Exercise}%
Show that in a~$\diamond$-effectus
    both~$f^\diamond$ and~$f^\BOX$ are order preserving maps.
\end{point}
\begin{point}{30}[diamond-adjunction]{Proposition}%
For~$f\colon X \to Y$ in a~$\diamond$-effectus we have
\begin{equation}
    f^\diamond(s) \ \leq\  t^\perp
    \quad \iff
    \quad f_\diamond(t) \ \leq\  s^\perp \label{diamond-main-lemma}
\end{equation}
for all sharp~$s,t$.
In other words: $f_\diamond$ is the left order-adjoint of~$f^\BOX$.
\begin{point}{40}{Proof}%
To start, let's prove the order-adjunction reformulation
\begin{equation*}
    f_\diamond(s) \leq t
    \ \overset{\eqref{diamond-main-lemma}}{\iff} \ 
    f^\diamond(t^\perp) \equiv f^\BOX(t)^\perp \leq s^\perp
            \ \iff \
    s \leq f^\BOX(t).
\end{equation*}
To prove \eqref{diamond-main-lemma},
first assume~$f^\diamond(s) \leq t^\perp$.
Then~$s \after f \leq \ceil{s \after f} = f^\diamond(s) \leq t^\perp
= \IMperp \pi_t$, where the last equality is
due to \sref{img-of-compr}.
Thus~$s \after f \after \pi_t = 0$
    which is to say~$s \leq \IMperp f \after \pi_t$,
    so~$f_\diamond(t) = \IM f \after \pi_t \leq s^\perp$.

For the converse, assume
    $f_\diamond(t) \leq s^\perp$.
Then as before (but in the other direction)
    we find~$s \after f \after \pi_t = 0$
    and so~$s \after f \leq t^\perp$.
    Hence~$f^\diamond(s) =  \ceil{s \after f} \leq \lceil t^\perp\rceil
                = \floor{t}^\perp = t^\perp$, as desired. \qed
\end{point}
\end{point}
\begin{point}{50}[order-adj-basics]{Exercise}%
Use the fact that there is an
order adjunction between~$f_\diamond$ and $f^\BOX$
to show that in a~$\diamond$-effectus
\begin{multicols}{2}
\begin{enumerate}
    \item $f_\diamond$ is order preserving;
    \item $f_\diamond$ preserves suprema;
    \item $f^\BOX$ preserves infima;
    \item $f^\diamond$ preserves suprema;
    \item $f_\diamond \after f^\BOX \after f_\diamond = f_\diamond$ \emph{and}
    \item $f^\BOX \after f_\diamond \after f^\BOX = f^\BOX$.
\end{enumerate}
\end{multicols}
\end{point}
\spacingfix{}
\begin{point}{60}[diamond-functor]{Lemma}%
    In a~$\diamond$-effectus
        $(\ )^\diamond$,
        $(\ )^\BOX$ and
        $(\ )_\diamond$
        are functorial --- that is
\begin{multicols}{3}
\begin{enumerate}
\item
$(\id)^\diamond = \id$,
\item
$(f\after g)^\diamond
            =   g^\diamond \after f^\diamond$,
\item
$(\id)^\BOX = \id$,
\item
$(f\after g)^\BOX
            =   g^\BOX \after f^\BOX$,
\item
$(\id)_\diamond= \id$ and
\item
            $(f \after g)_\diamond = f_\diamond \after g_\diamond$.
\end{enumerate}
\end{multicols}
\spacingfix{}
\begin{point}{70}{Proof}%
We get~$(\id)^\diamond = \id$
directly from~\sref{img-of-compr}.
For 2 we only need a single line:
\begin{equation*}
(f\after g)^\diamond(s)
    \ =\  \ceil{s \after f \after g}
    \ \overset{\sref{ceiling-within-ceiling}}{=}\ 
    \ceil{\ceil{s \after f} \after g}
    \ =\  g^\diamond(f^\diamond(s)).
\end{equation*}
Note that we used that ceilings are sharp.
Point 3 and 4 follow easily from~1 and 2 respectively.
The identity~$(\id)_\diamond = \id$ is again~\sref{img-of-compr}.
We claim~$f_\diamond \after g_\diamond$
    is left order-adjoint to~$(f\after g)^\BOX$,
     indeed
\begin{equation*}
    f_\diamond (g_\diamond(s)) \leq t
        \ \iff\  g_\diamond(s) \leq f^\BOX (t)
        \ \iff\  s \leq g^\BOX (f^\BOX (t)) = (f \after g)^\BOX(t).
\end{equation*}
Thus by uniqueness of order adjoints, we find
$(f\after g)_\diamond = f_\diamond \after g_\diamond$. \qed
\end{point}
\end{point}
\begin{point}{71}[diamond-equiv-equiv]{Exercise}%
Derive from~\sref{diamond-adjunction}
    that~$f^\diamond = g^\diamond$
    if and only if~$f_\diamond = g_\diamond$.
\end{point}
\end{parsec}
\begin{parsec}{2080}%
\begin{point}{10}[image-sharp-is-order-sharp]{Lemma}%
In a $\diamond$-effectus,
sharp predicates are \Define{order sharp} ---
    \index{sharp!order-}
    that is:
    for any predicate~$p$
    and sharp predicate~$s$
    with~$p \leq s$ and~$p \leq s^\perp$,
    we must have~$p = 0$.
\begin{point}{20}{Proof}%
Note~$\ceil{p} \leq \ceil{s} = s$
    and~$\ceil{p} \leq \lceil s^\perp\rceil = s^\perp$.
So by~\sref{compr-is-full},
there is an~$h$ with~$\pi_{\ceil{p}} = \pi_s \after h$.
We compute
\begin{equation*}
    1 \after \pi_{\ceil{p}} \ =\  \ceil{p} \after \pi_{\ceil{p}}  \ = \
    \ceil{p} \after \pi_s \after h \ \leq \ 
        s^\perp \after \pi_s \after h \ = \ 0.
\end{equation*}
Thus~$\pi_{\ceil{p}} = 0$
and so~$p \leq \ceil{p} = \IM\pi_{\ceil{p}} = 0$, as desired.\qed
\end{point}
\end{point}
\begin{point}{30}[diamond-oml]{Proposition (Cho)}%
In a~$\diamond$-effectus,
the poset~$\SPred X$ of sharp predicates on any object~$X$,
is a sub-effect algebra of~$\Pred X$
and an orthomodular lattice.
\begin{point}{40}{Proof}%
Let~$C$ be a~$\diamond$-effectus with some object~$X$.
We will show~$\SPred X$ is a sub-effect algebra of~$\Pred X$,
    which is additionally an ortholattice.
    By \sref{orth-ea-is-orthomodular} this is sufficient to show
    $\SPred X$ is orthomodular.
\begin{point}{50}{Ortholattice}%
In \sref{lattice-compr}
    we already saw that sharp~$s,t$ have a sharp supremum
    $s\vee t = \IM [\pi_s,\pi_t]$ 
    in~$\Pred X$.
So~$s\vee t$ is also the supremum of~$s$ and~$t$ in~$\SPred X$.
As~$(\ )^\perp$ is an order anti-automorphism of both~$\Pred X$ and~$\SPred X$,
    we know~$(s^\perp \vee t^\perp)^\perp$
    is the infimum of~$s$ and~$t$ in~$\Pred X$ and~$\SPred X$.
As any sharp predicate~$s$ is order sharp by \sref{image-sharp-is-order-sharp},
    we find~$s \wedge s^\perp = 0$ (and consequently~$s^\perp \vee s = 1$).
    We have shown~$\SPred X$ is an ortholattice.
\end{point}
\begin{point}{60}{Sub-EA}%
Clearly~$0,1$ are sharp and~$s^\perp$ is sharp for sharp~$s$
        by definition of~$\diamond$-effectus.
To prove~$\SPred X$ is a sub-effect algebra of~$\Pred X$,
    it only remains to be shown~$s \ovee t$ is sharp
        for sharp and summable~$s,t$.
So, assume~$s,t$ are sharp and~$s \perp t$.
Note~$s \wedge t \leq s \wedge s^\perp =0$ as~$\SPred X$ is an ortholattice
    and so by \sref{ea-modularity-prop}
    we find~$s \ovee t = s \vee t$, which is indeed sharp. \qed
\end{point}
\end{point}
\begin{point}{70}{Corollary}%
In an $\diamond$-effectus~$C$
the assignment~$X\mapsto \SPred X$, ~$f \mapsto (f_\diamond,f^\BOX)$
    yields a functor from~$C$
    to~$\mathsf{OMLatGal}$,
    the category of orthomodular lattices
    with Galois connection between them,
    as defined in~\cite{jacobs2009orthomodular}.
\end{point}
\begin{point}{80}%
    There is a rather different
    formula for the infima of sharp predicates,
    which will be useful later on.
\end{point}
\end{point}
\begin{point}{90}[spred-infimum]{Lemma}%
In a~$\diamond$-effectus, we have
$s \wedge t = (\pi_s)_\diamond (\pi_s^\BOX (t))$
    for sharp predicates~$s,t$.
\begin{point}{100}{Proof}%
It is easy to see
$(\pi_s)_\diamond (\pi_s^\BOX (t))$ is a lower bound:
indeed
$(\pi_s)_\diamond (\pi_s^\BOX (t)) \leq t$
and~$(\pi_s)_\diamond (\pi_s^\BOX (t)) \leq (\pi_s)_\diamond(1)
    = \IM \pi_s = s$.
We have to show $(\pi_s)_\diamond (\pi_s^\BOX (t))$ is the greatest lower bound.
Let~$r$ be any sharp predicate with~$r \leq s$ and~$r \leq t$.
By~\sref{compr-is-full}
    we have~$\pi_r = \pi_s \after h$ for some~$h$.
Thus
\begin{equation*}
    (\pi_r)_\diamond
        \ = \  (\pi_s)_\diamond \after h_\diamond
        \ =\  (\pi_s)_\diamond \after \pi_s^\BOX  
            \after (\pi_s)_\diamond \after h_\diamond
        \ =\  (\pi_s)_\diamond \after \pi_s^\BOX \after (\pi_r)_\diamond.
\end{equation*}
Hence~$r = (\pi_r)_\diamond(1)  
=((\pi_s)_\diamond \after \pi_s^\BOX \after (\pi_r)_\diamond )(1)
=  (\pi_s)_\diamond ( \pi_s^\BOX (r))
\leq  (\pi_s)_\diamond ( \pi_s^\BOX (t))$. \qed
\end{point}
\begin{point}{110}{Remark}%
In \cite{effintro} there appears a similar result due to Jacobs,
    which is based on subtly different assumptions.
\end{point}
\end{point}
\begin{point}{120}[spred-sup]{Exercise}%
Show that in a~$\diamond$-effectus,
we have
        $(\xi^\BOX \after \xi_\diamond )(t) = s \vee t$
    for sharp~$s,t$ and  quotient~$\xi$ of~$s$.
    (Hint: mimic \sref{spred-infimum}.)
\end{point}
\end{parsec}

\begin{parsec}{2090}%
\begin{point}{10}%
We turn to~$\diamond$-adjointness.
\end{point}
\begin{point}{20}[exc-diamond-adj]{Exercise}%
Show the following basic properties of~$\diamond$-adjointness
\begin{enumerate}
    \item
$f^\diamond = g_\diamond$
    ($f$ is $\diamond$-adjoint to~$g$)
    if and only if~$f_\diamond = g^\diamond$.
    \item
If~$f$ and~$g$ are~$\diamond$-adjoint,
    then~$\IM f = \ceil{1 \after g}$.
\end{enumerate}
\end{point}
\spacingfix{}
\begin{point}{30}[diamond-squares]{Exercise}%
Show in order:
\begin{enumerate}
\item
If~$f$ is~$\diamond$-self-adjoint,
    then~$f \after f$ is~$\diamond$-self-adjoint.
\item
If~$f$ is~$\diamond$-positive,
    then~$f$ is~$\diamond$-self-adjoint.
\item
If~$f$ is~$\diamond$-positive and~$f\after f$ is pure,
    then~$f \after f$ is~$\diamond$-positive.
\end{enumerate}
\end{point}
\spacingfix{}
\begin{point}{40}[iso-diamond-adjoint]{Lemma}%
Let~$\alpha$ be an isomorphism in a~$\diamond$-effectus.
Then
\begin{enumerate}
\item
    $s\after\alpha$ is sharp for sharp predicates~$s$ \emph{and}
\item
    $\alpha^\diamond(s) = s \after \alpha$
    and~$\alpha_\diamond(s) = s\after \alpha^{-1}$
    (so~$\alpha$ and~$\alpha^{-1}$ are~$\diamond$-adjoint).
\end{enumerate}
\spacingfix{}
\begin{point}{50}{Proof}%
Let~$s$ be a sharp predicate. Then~$s = \IM \pi_s$.
Note~$\IM \alpha^{-1}\after \pi_s = s \after \alpha$
    --- indeed, $s \after \alpha \after \alpha^{-1} \after \pi_s=1$
    and when~$p \after \alpha^{-1} \after \pi_s=1$,
    we must have~$p \after \alpha^{-1} \geq s$,
    which gives~$p \geq s \after \alpha$ as desired.
    So~$s \after \alpha$ is indeed sharp.

So~$\alpha^\diamond(s) = \ceil{s \after \alpha} = s\after \alpha$
    and~$\alpha_\diamond(s) = \IM \alpha \after \pi_s = s \after \alpha^{-1}$
    as promised. \qed
\end{point}
\end{point}
\end{parsec}

\begin{parsec}{2100}%
\begin{point}{10}[sharp-map]{Definition}%
A map~$f$ in a~$\diamond$-effectus is a~\Define{sharp map}
    \index{sharp!map}
    provided~$s \after f$ is sharp for all sharp predicates~$s$.
\end{point}
\begin{point}{20}[sharp-ceil]{Exercise}%
Show that the following are equivalent.
\begin{enumerate}
    \item $f$ is a sharp map.
    \item $\ceil{p \after f} = \ceil{p}\after f$
            for every predicate~$p$.
\end{enumerate}
\end{point}
\spacingfix{}
\begin{point}{30}[exa-sharp-vn]{Example}%
In~$\op{\vN}$ the sharp maps are exactly the~mni-maps
(i.e.~the normal $*$-homomorphisms).
See \sref{sharp-multiplicative}.
\end{point}
    
\end{parsec}

\section{$\&$-effectuses}
\begin{parsec}{2110}%
\begin{point}{10}%
In a~$\diamond$-effectus
    quotient and comprehension are not tied together by its axioms.
    With two additional axioms, we will see quotient and comprehension
    become tightly interwoven.
\end{point}
\begin{point}{20}{Definition}%
An~\Define{$\&$-effectus} (``andthen effectus'')
    \index{effectus!$\&$-}
is a $\diamond$-effectus
such that
\begin{enumerate}
\item
    for each object~$X$
    and each predicate~$p$ on~$X$,
        there is a unique $\diamond$-positive
        map $\Define{\asrt_p}\colon X \to X$
        \index{asrt@$\asrt_p$, assert}
        (``assert p'')
        with~$1 \after \asrt_p = p$ (see \sref{diamond-basics})
        \emph{and}
\item
    for every quotient~$\xi \colon Y \to Z$
    and comprehension~$\pi \colon X \to Y$
    the composite~$\xi \after \pi$ is pure.
\end{enumerate}
In an $\&$-effectus,
we define~$\Define{\andthen{p}{q}} \equiv q \after \asrt_p$,
    \index{$\&$, andthen}
    pronounced ``$p$ andthen~$q$''.
For brevity, we will write~$\Define{p^2} \equiv \andthen{p}{p}$.
\begin{point}{30}[asrt-remarks]{Remarks}%
The $\asrt_p$ maps are named after the \texttt{assert} statement
    found in many programming languages.
The \texttt{assert} statement checks whether
    the provided Boolean expression (i.e.~predicate)
    is true (at the time of execution) --- and if it is, it will continue the program
    without further action; if it isn't, it will halt execution
    immediately.
Our~$\asrt_p$ maps can be thought of in the same way.
Alternatively, $\asrt_p$ can be viewed as a filter
    that blocks states for which~$p^\perp$ holds.
With this intuition,
    the predicate~$\andthen{p}{q}$
    corresponds to whether~$p$ and~$q$ are true,
    by first checking~$p$ and then checking~$q$.

The~$\asrt_p$ maps can be used to construct more complicated measurements.
For instance, let~$p, q, r$ be three predicates (measurement outcomes)
    on an object~$X$
    with~$p \ovee q \ovee r = 1$.
Then~$\left<\asrt_p,\asrt_q,\asrt_r\right>\colon X \to X + X + X $
    models measuring whether~$p,q,r$ hold
    without discarding the state of the (possibly affected) system
    after measurement
    (in contrast to~$\left<p,q,r\right>\colon X \to 1+1+1$).
If we discard the measurement outcome,
    we find the map~$\asrt_p \ovee \asrt_q \ovee \asrt_r \colon X \to X$,
    which is the side-effect of measuring $p,q,r$.
\begin{point}{31}
The second axiom is closely related to polar decomposition:
    in the case of~$\op\vN$ (flipping direction of arrows now),
    the unique isomorphism~$\varphi$
    such that~$h \after c = c' \after \varphi \after h'$
    for the suitable standard filter~$c'$ and corner~$h'$
    is given by~$\ad_u$
    where~$u$ is the isomtry from the
    polar decomposition~$u \sqrt{\ceil{h}c(1)\ceil{h}}
        \equiv \sqrt{c(1)} \ceil{h}$
    restricted to a unitary.
\end{point}
\end{point}
\end{point}
\begin{point}{40}[vn-is-andthen-eff]{Examples}%
The category~$\op{\vN}$ is an~$\&$-effectuses
    with~$\asrt_a\colon b \mapsto \sqrt{a}b\sqrt{a}$.
The first axiom is proven in \sref{positive-map-uniqueness}
    and the second in \sref{pure-fundamental}.
The full subcategory~$\op\CvN$ of commutative von Neumann algebras
    is a~$\&$-effectus as well.

The only other known example of an~$\&$-effectus
    is the category~$\op\EJA$ of Euclidean Jordan algebras
    with positive unital maps in the opposite direction,
    see~\cite{eja}.
\end{point}
\begin{point}{50}[sharp-prop]{Proposition}%
For a predicate~$p$ in an $\&$-effectus the following are equivalent.
\begin{enumerate}
\item $p$ is sharp,
\item $\andthen{p}{p}=p$ and
\item $\asrt_p \after \asrt_p = \asrt_p$.
\end{enumerate}
\spacingfix{}
\begin{point}{60}{Proof}%
First we prove that~$p$ is sharp if and only if~$\andthen{p}{p}=p.$
So, assume~$p$ is sharp.
As $\diamond$-positive maps are~$\diamond$-self-adjoint
    we have~$\IM \asrt_p = \ceil{1 \after \asrt_p} = \ceil{p} = p$.
Thus~$\andthen{p}{p} = p \after \asrt_p = 1 \after \asrt_p = p$.
For the converse, assume~$\andthen{p}{p}=p$.
From~$\IM \asrt_p = \ceil{1 \after \asrt_p} = \ceil{p}$,
we get~$\ceil{p} \after \asrt_p = \andthen{p}{\ceil{p}} = 1 \after \asrt_p = p$.
By assumption~$\andthen{p}{p}=p$.
So~$\andthen{p}{(\ceil{p} \ominus p)} = 0$.
Hence~$\ceil{p}\ominus p \leq \IMperp \asrt_p = \ceil{p}^\perp$.
However $\ceil{p}\ominus p \leq \ceil{p}$.
Thus by~\sref{image-sharp-is-order-sharp}
    get~$\ceil{p} \ominus p = 0$. So~$p$ is indeed sharp.

Clearly, if~$\asrt_p \after \asrt_p = \asrt_p$,
then~$p = 1 \after \asrt_p = 1\after\asrt_p\after\asrt_p = \andthen{p}{p}$.
It only remains to be shown~$\asrt_p \after \asrt_p = \asrt_p$
    whenever~$p$ is sharp.
So assume~$p$ is sharp.
By definition~$\asrt_p$ is pure: $\asrt_p = \pi \after \xi$
    for some quotient~$\xi$ and comprehension~$\pi$.
By assumption~$\xi \after \pi$ is pure as well,
    so there is a quotient~$\xi'$ and comprehension~$\pi'$
    with~$\xi \after \pi = \pi' \after \xi'$.
Note~$1 \after \xi = 1 \after \asrt_p = p$
and~$1 \after \asrt_p \after \asrt_p = \andthen{p}{p} = p = 1 \after \xi$,
hence
\begin{equation*}
 1\after \xi  
    \ =\  1 \after \asrt_p \after \asrt_p 
    \ =\  1 \after \pi \after \xi \after \pi \after \xi 
    \ =\  1 \after \pi \after \pi' \after \xi' \after \xi 
    \ =\  1 \after \xi' \after \xi,
\end{equation*}
so~$1 \after \xi' = 1$. Thus~$\xi'$ is an iso.
Also~$
    (\IM \pi') \after \xi \after \pi
    = (\IM \pi') \after \pi' \after \xi'
    = 1 \after \xi' = 1 $
    from which it follows~$p = \IM \pi \leq  (\IM \pi') \after \xi 
            \leq 1 \after \xi = p$.
Thus~$(\IM \pi') \after \xi = p = 1\after \xi$.
Hence~$\IM \pi' = 1$ and so~$\pi'$ is an isomorphism.
Now we know~$\pi \after \pi'$
    is a comprehension and~$\xi' \after \xi$ is a quotient,
    we see~$\asrt_p\after\asrt_p$ is pure.
As~$\asrt_p$ is $\diamond$-self-adjoint,
    we see~$\asrt_p\after\asrt_p$ is $\diamond$-positive.
By uniqueness of positive maps~$\asrt_p \after \asrt_p = \asrt_p$,
    as desired. \qed
\end{point}    
\end{point}
\begin{point}{70}[prop-corr-zeta-pi]{Proposition}%
Let~$C$ be an~$\&$-effectus and~$s$
     be any sharp predicate.
There exist comprehension~$\pi_s$ of~$s$
        and quotient~$\zeta_s$ of~$s^\perp$ such that
\begin{equation*}
    \zeta_s \after \pi_s \ = \ \id \quad \text{and} \quad
        \pi_s \after \zeta_s \ = \ \asrt_s.
\end{equation*}
In fact, for every comprehension~$\pi$ of~$s$,
    there is a quotient~$\xi$ of~$s^\perp$
    with~$\xi \after \pi = \id$ and~$\pi \after \xi = \asrt_s$
    \emph{and} conversely for every quotient~$\xi$ of~$s^\perp$
    there exists a comprehension~$\pi$ of~$s$
    with~$\xi \after \pi = \id$ and~$\pi \after \xi = \asrt_s$.
\begin{point}{80}{Proof}%
Let~$s$ be any sharp predicate.
By definition~$\diamond$-positive maps are pure and so
    there is a quotient~$\xi$ and comprehension~$\pi$
    with~$\pi \after \xi = \asrt_s$.
    So
\begin{equation*}
   \pi \after \xi \ =\  \asrt_s 
    \ \stackrel{\mathclap{\sref{sharp-prop}}}{=}\  \asrt_s\after\asrt_s \ =\  \pi \after \xi \after \pi \after \xi.
\end{equation*}
Thus~$\xi \after \pi = \id$.
We compute~$s = 1 \after \asrt_s = 1 \after \xi$
and so
\begin{alignat*}{2}
    \IM \pi &\ = \ 
    \IM \pi \after \xi \after \pi &\qquad& \text{as $\xi \after \pi=\id$}\\
                  &\ = \ \IM \asrt_s \after \pi &&\text{by dfn.~$\xi$ and~$\pi$}\\
                  &\ = \ (\asrt_s)_{\diamond}(\IM \pi)&& \text{by dfn.~$(\ )_\diamond$} \\
                  & \ = \ (\asrt_s)^{\diamond}(\IM \pi)&&
        \text{by $\diamond$-s.a.~$\asrt_s$}\\
        &\ = \ \ceil{(\IM \pi) \after \pi \after \xi} && \text{by dfn.~$(\ )^\diamond$} \\
&\ = \ \ceil{1 \after \xi} \\
    & \ = \ s.
\end{alignat*}
Thus~$\pi$ is a comprehension of~$s$
and~$\xi$ is a quotient of~$s^\perp$.
We have proven the first part.

For the second part,
let~$\pi'$ be any comprehension of~$s$.
Then~$\pi' = \pi \after \alpha$ for some iso~$\alpha$.
Define~$\xi' = \alpha^{-1} \after \xi$.
It is easy to see~$\pi' \after \xi' = \asrt_s$
and~$\xi' \after \pi' = \id$.
The other statement is proven in a similar way. \qed
\end{point}
\begin{point}{90}[zeta-s-convention]{Notation}%
In an~$\&$-effectus
    together with chosen comprehension~$\pi_s$ for~$s$,
    we will write~$\zeta_s$
    for the unique quotient for~$s^\perp$
    satisfying~$\zeta_s \after \pi_s = \id$
    and~$\pi_s \after \zeta_s = \asrt_s$.
We call this~$\zeta_s$
    \index{quotient!corresponding}
    \index{comprehension!corresponding}
    the \Define{corresponding quotient} of~$\pi_s$
    and vice versa~$\pi_s$ the \Define{corresponding comprehension}
    of~$\zeta_s$. \index{*zeta@$\zeta_s$}\index{*pi@$\pi_p$}
\end{point}
\begin{point}{100}{Warning}%
    $\xi_{s^\perp} = \zeta_s$!
\end{point}
\end{point}
\begin{point}{110}[upm-closed]{Proposition}%
In an $\&$-effectus
    both comprehensions and pure maps are closed under composition.
\begin{point}{120}{Proof}%
We will first prove that comprehensions are closed under composition.
Assume~$\pi_1 \colon X \to Y$
    and~$\pi_2 \colon Y \to Z$
    are comprehensions with~$s = \IM \pi_1$
    and~$t = \IM \pi_2$.
We will show~$\pi_2 \after \pi_1$
    is a comprehension for~$\IM \pi_2 \after \pi_1$.
To this end, let~$f\colon V\to Z$ be any map
with~$(\IM \pi_2 \after \pi_1)^\perp \after f = 0$.
As~$\IM \pi_2 \after \pi_1 \leq \IM \pi_2 = t$
    we get~$t^\perp \after f = 0$,
    so~$f = \pi_2 \after g_2$ for a unique~$g_2\colon V \to Y$.
Let~$\zeta_2$ be a quotient for~$t^\perp$
such that~$\zeta_2 \after \pi_2 = \id$,
which exists by \sref{prop-corr-zeta-pi}.
Then~$s \after \zeta_2 \after \pi_2 \after \pi_1
            = s \after \pi_1 = 1$
            so~$s \after \zeta_2 \geq \IM \pi_2 \after \pi_1$.
Consequently
\begin{equation*}
    s \after g_2 \ =\  s \after \zeta_2 \after \pi_2 \after g_2
        \ =\  s \after \zeta_2 \after f
        \ \geq\   (\IM \pi_2 \after \pi_1) \after f \ =\  1 \after f
                \ =\  1 \after g_2.
\end{equation*}
Hence there exists a unique~$g_1\colon V \to X$
    with~$\pi_1 \after g_1 = g_2$
    and so~$\pi_2 \after \pi_1 \after g_1 = f$.
By monicity of~$\pi_2 \after \pi_1$,
    this~$g_1$ is unique and so~$\pi_2 \after \pi_1$
    is indeed a comprehension.
\begin{point}{130}%
Now we will prove that pure maps are closed under composition.
Assume~$g\colon X \to Y$ and~$f\colon Y \to Z$ are pure maps.
Say~$f = \pi_1 \after \xi_1$
    and~$g = \pi_2 \after \xi_2$
    for some comprehensions~$\pi_1$, $\pi_2$
        and quotients~$\xi_1$, $\xi_2$.
By definition of $\&$-effectus
    there is a comprehension~$\pi'$ and quotient~$\xi'$
    such that~$\xi_1 \after \pi_2 = \pi' \after \xi'$.
By the previous point~$\pi_1 \after \pi'$
    is a comprehension
    and~$\xi' \after \xi_2$ is a quotient by \sref{quotients-composition}.
Thus
$f \after g 
=\pi_1 \after \pi' \after \xi' \after \xi_2$
is pure. \qed
\end{point}
\end{point}
\end{point}
\begin{point}{140}[andthen-square-rule]{Exercise}%
Show that in an~$\&$-effectus,
    we have~$\asrt_p \after \asrt_p = \asrt_{\andthen{p}{p}}$.
\end{point}
\begin{point}{150}[asrt-absorp-rule]{Exercise}%
    Show that in an~$\&$-effectus,
    we have
    \begin{align*}
        \IM f \leq s \quad&\iff\quad  \asrt_s \after f = f \\
        1 \after f \leq t \quad&\iff\quad  f \after \asrt_t = f
    \end{align*}
for any sharp predicates~$s$ and~$t$.
\end{point}
\begin{point}{160}{Definition}%
Let~$C$ be an~$\&$-effectus.
In \sref{upm-closed} we saw pure maps are closed under composition.
Write~$\Define{\Pure C}$ for the subcategory
    \index{PureC@$\Pure C$} of pure maps.
\end{point}
\end{parsec}

\begin{parsec}{2120}%
\begin{point}{10}[zeta-asrt-quot]{Lemma}%
    In an~$\&$-effectus,
    $\zeta_{\ceil{p}} \after \asrt_p$
    is a quotient for~$p^\perp$.
\begin{point}{20}{Proof}%
Recall~$\IM \asrt_p = \ceil{1 \after \asrt_p} = \ceil{p}$ and so
\begin{equation*}
\IM \zeta_{\ceil{p}} \after \asrt_p
\ =\  (\zeta_{\ceil{p}})_\diamond( \ceil{p})
\ =\  \IM \zeta_{\ceil{p}} \after \pi_{\ceil{p}}
\ =\  \IM \id \ =\  1,
\end{equation*}
where~$\pi_{\ceil{p}}$ is the comprehension corresponding to~$\zeta_{\ceil{p}}$.
By~\sref{upm-closed}~$\zeta_{\ceil{p}} \after \asrt_p = \pi \after \xi$
    for some comprehension~$\pi$ and quotient~$\xi$.
Putting it together:
$1 = \IM \zeta_{\ceil{p}} \after \asrt_p = \IM \pi \after \xi \leq \IM \pi$,
so~$\IM \pi = 1$ and thus~$\pi$ is an isomorphism,
hence~$\zeta_{\ceil{p}}\after \asrt_p$
is a quotient for the orthocomplement of~$\ceil{p} \after \asrt_p = p$. \qed
\end{point}
\end{point}
\begin{point}{30}[standard-form-map]{Proposition}%
Every map~$f$
in an~$\&$-effectus
factors as
\begin{equation}\label{eq-standard-form-map}
        f \ =\ \pi_{\IM f} \after g \after \zeta_{\ceil{1 \after f}} \after \asrt_{1 \after f} 
    \end{equation}
    for a unique total and faithful map~$g$.
Furthermore the following holds.
\begin{enumerate}
\item
If~$f$ is pure, then~$g$ is an isomorphism.
\item
If~$f$ is pure and~$\IM f = 1$, then~$f$ is a quotient.
\item
If~$f$ is pure and~$1 \after f = 1$, then~$f$ is a comprehension.
\end{enumerate}
\spacingfix{}
\begin{point}{40}{Proof}%
By the universal property of~$\pi_{\IM f}$,
    there is a unique~$g'$ with~$f = \pi_{\IM f} \after g'$.
To show~$g'$ is faithful, assume~$p\after g' = 0$
    for some predicate~$p$.
Then~$0 = p \after g' = p \after \zeta_{\IM f} \after \pi_{\IM f} \after g'
                = p \after \zeta_{\IM f} \after f$
    and so~$p \after \zeta_{\IM f} \leq \IMperp f$,
    hence~$p = p \after \zeta_{\IM f} \after \pi_{\IM f}
                \leq (\IMperp f )\after \pi_{\IM f} = 0$, which
                shows~$g'$ is indeed faithful.

Note~$1 \after g' = 1 \after \pi_{\IM f} \after g' = 1 \after f$.
By~\sref{quotient-total} and~\sref{zeta-asrt-quot}
    there is a unique total~$g$ with~$g' = g \after \zeta_{\ceil{1 \after f}} \after \asrt_{1 \after f}$.
Clearly~\eqref{eq-standard-form-map} holds
and~$g$ is the unique map
    for which \eqref{eq-standard-form-map} holds
    as  comprehensions are mono and quotients are epi.
    As~$1 = \IM g' = \IM g \after \zeta_{\ceil{1 \after f}}
            \after \asrt_{1 \after f} \leq \IM g$,
            we see~$\IM g = 1$ and so~$g$ is faithful.

To prove point 1, assume~$f$ is pure.
That is: $f = \pi \after \xi$ for some comprehension~$\pi$
        and quotient~$\xi$.
As~$1 \after \pi=1$ and~$\IM \xi = 1$,
    we see~$\IM \pi = \IM f$ and~$1 \after \xi = 1\after f$.
Thus~$\pi = \pi_{\IM f} \after \alpha$
    and~$\xi = \beta \after \zeta_{\ceil{1 \after f}}
                \after \asrt_{1 \after f}$
        for some isomorphisms~$\alpha$ and~$\beta$.
Thus~$f = \pi_{\IM f} \after \alpha \after \beta
                \after \zeta_{\ceil{1 \after f}}
                \after \asrt_{1 \after f}$
        and so by uniqueness of~$g$,
        we see~$g = \alpha \after \beta$ is an isomorphism.

To prove point 2, additionally assume~$\IM f= 1$.
Then~$\pi_{\IM f}$ is an iso and so
using the previous, we see~$f$ is indeed a quotient.
Point 3 is just as easy. \qed
\end{point}
\end{point}
\end{parsec}

\begin{parsec}{2130}%
\begin{point}{10}[andthen-effect-divisoid]{Proposition}%
If~$C$ is an~$\&$-effectus,
    then~$(\Scal C)^{\mathsf{op}}$ is an effect divisoid,
    see~\sref{dfn-effect-divisoid}.
\begin{point}{20}{Proof}%
Let~$\lambda,\mu$ be scalars with~$\lambda \leq \mu$.
Recall that the scalar~$1=\id$
    and so~$\mu = 1 \after \asrt_\mu = \asrt_\mu$.
By~\sref{zeta-asrt-quot},
    there is a unique~$\lambda'$
    with~$\lambda' \after \zeta_{\ceil{\mu}} \after \mu = \lambda$.
Define~$\rfrac{\lambda}{\mu} = \lambda'  \after \zeta_{\ceil{\mu}}$.

For the moment, assume~$\lambda = \mu$.
Then~$\lambda' \after \zeta_{\ceil{\mu}} \after \mu
            = \mu = 1 \after \zeta_{\ceil{\mu}} \after \mu$
            as~$\IM \mu = \IM \asrt_{\mu} = \ceil{\mu}$.
        By epicity of~$\zeta_{\ceil{\mu}} \after \mu$,
            we see~$\lambda' = 1$, hence~$\rfrac{\mu}{\mu}=\ceil{\mu}$.
As we have~$\mu \leq \ceil{\mu}$ and~$\ceil{\ceil{\mu}} = \ceil{\mu}$,
    we see axioms 2 and 3 hold.

We return to the general case to prove axiom 1.
Clearly~$\rfrac{\lambda}{\mu} = \lambda' \after \zeta_{\ceil{\mu}}
            \leq 1 \after \zeta_{\ceil{\mu}} = \ceil{\mu} 
            = \rfrac{\mu}{\mu}$
    and~$\mu \odot_{\mathsf{op}} \rfrac{\lambda}{\mu}
= \lambda'  \after \zeta_{\ceil{\mu}} \after \mu = \lambda$
as required.
Assume~$\sigma$ is an arbitrary
    scalar with~$\mu \odot_{\mathsf{op}} \sigma = \lambda$
    and~$\sigma \leq \rfrac{\mu}{\mu} \equiv \ceil{\mu}$.
Then~$\sigma = \sigma' \after \zeta_{\ceil{\mu}}$
    for a unique~$\sigma'$.
As~$
\lambda'  \after \zeta_{\ceil{\mu}} \after \mu = \lambda
= \sigma' \after \zeta_{\ceil{\mu}} \after \mu$
we must have~$\lambda' = \sigma'$, whence~$\rfrac{\lambda}{\mu}=\sigma$. \qed
\end{point}
\end{point}
\begin{point}{30}[perp-sharp-is-orth]{Lemma}%
In an~$\&$-effectus, if~$s \perp t$
    for sharp~$s,t$, then~$\andthen{s}{t} = 0 = \andthen{t}{s}$.
\begin{point}{40}{Proof}%
Write~$r \equiv (s \ovee t)^\perp$.
    As~$1 = s \ovee t \ovee r$ and~$\andthen{s}{s}=s$, we have
\begin{equation*}
    s \ = \ 
    \andthen{s}{1} \ = \ 
    \andthen{s}{s} \ovee \andthen{s}{t} \ovee \andthen{s}{r} \ = \ 
    s \ovee \andthen{s}{t} \ovee \andthen{s}{r}
\end{equation*}
and so~$\andthen{s}{t} \leq \andthen{s}{t} \ovee \andthen{s}{r} = 0$. \qed
\end{point}
\end{point}
\begin{point}{50}[simple-andthen-absorption]{Exercise}%
Show that in an~$\&$-effectus, we have
\begin{equation*}
p \ \leq\  s \quad \iff \quad \andthen{s}{p} \ =\  p
\end{equation*}
for sharp~$s$ and any predicate~$p$.
\end{point}
\begin{point}{60}[exc-prod-sharp-maps]{Exercise}%
Show that in a~$\&$-effectus,
    a map~$\left<f,g\right>$
    is sharp if and only if both~$f$ and~$g$ are sharp.
    (Hint: use \sref{diamond-oml} and \sref{img-tupling}.)
\end{point}
\end{parsec}

\section{\texorpdfstring{$\dagger$-effectuses}{%
                    dagger-effectuses}}
\begin{parsec}{2140}%
\begin{point}{10}[dagger-effectus]{Definition}%
A \Define{$\dagger$-category} (``dagger category'')
    \index{*dagg@$\dagger$, dagger!-category}
is a category~$C$ together with an involutive identity-on-objects
    functor~$(\ )^\dagger \colon C \to \op{C}$
    ---
    that is, for all objects~$X$ and maps~$f,g$ in~$C$, we have
\begin{multicols}{2}
\begin{enumerate}
    \item $(f \after g)^\dagger = g^\dagger \after f^\dagger$
    \item $\id^\dagger = \id$
    \item $f^{\dagger\dagger} = f$ \emph{and}
    \item $X^\dagger = X$.
\end{enumerate}
\end{multicols}
\noindent
Cf.~\cite{burgin1970categories,selinger,heunenphd,karvonen}.
In any~$\dagger$-category we may define the following.
\begin{enumerate}
\item
    An endomap~$f$ is called \Define{$\dagger$-self-adjoint}
    \index{*dagg@$\dagger$, dagger!-self-adjoint}
    if~$f^\dagger = f$.
\item
    An endomap~$f$ is \Define{$\dagger$-positive}
    \index{*dagg@$\dagger$, dagger!-positive}
    if~$f = g^\dagger \after g$ for some other map~$g$.
\item
    An isomorphism~$\alpha$ is called~\Define{$\dagger$-unitary}
        \index{*dagg@$\dagger$, dagger!-unitary}
        whenever~$\alpha^{-1} = \alpha^\dagger$.
\end{enumerate}
\end{point}
\spacingfix{}
\begin{point}{20}{Example}%
The category~$\mathsf{Hilb}$
    of Hilbert spaces with bounded linear maps
    is a~$\dagger$-category
    with the familiar adjoint as~$\dagger$.
\end{point}
\end{parsec}
\begin{parsec}{2150}%
\begin{point}{10}{Definition}%
We call an~$\&$-effectus~$C$
    a~\Define{$\dagger$-effectus} (``dagger effectus'') provided
    \index{effectus!$\dagger$-}
    \index{*dagg@$\dagger$, dagger!-effectus}
\begin{enumerate}
\item
    $\Pure C$ is a  $\dagger$-category
     satisfying~$\asrt_p^\dagger = \asrt_p$
        and~$f$ is $\diamond$-adjoint to~$f^\dagger$;
\item
    for every~$\dagger$-positive~$f$,
        there is a unique~$\dagger$-positive~$g$
        with~$g \after g = f$ \emph{and}
\item
    $\diamond$-positive maps are $\dagger$-positive.
\end{enumerate}
\end{point}
\spacingfix{}
\begin{point}{20}{Examples}%
In \sref{vn-is-dagger-category}
we will see that the category~$\op\vN$ is a~$\dagger$-effectus.
Recently we have shown with van de Wetering
    that the category~$\op\EJA$
    of Euclidean Jordan algebras
    with positive unital maps in the opposite direction,
    is a~$\dagger$-effectus as well. \cite{eja}
\end{point}
\begin{point}{30}[dagger-theorem]{Theorem}%
    An $\&$-effectus
        is a~$\dagger$-effectus if and only if
\begin{enumerate}
\item
for every predicate~$p$, there is a unique predicate~$q$
    with~$\andthen{q}{q} = p$;
\item
    $\asrt^2_{\andthen{p}{q}}
        = \asrt_p \after \asrt^2_q \after \asrt_p$
     for all predicates~$p,q$ \emph{and}
\item
    a quotient for a sharp predicate (e.g.~$\zeta_s$)
    is a sharp map, see \sref{sharp-map}.
\end{enumerate}
\spacingfix{}
\begin{point}{40}{Proof}%
Necessity will be proven in
    \sref{dagger-thm-necessity}
        and sufficiency in \sref{dagger-thm-sufficiency}. \qed
\end{point}
\begin{point}{50}%
Especially the sufficiency requires quite some preparation.
For convenience, call~$C$ a \Define{$\dagger'$-effectus}
    \index{effectus!$\dagger'$-}
    if~$C$ is an~$\&$-effectus
    satisfying axioms 1, 2 and 3 from the Theorem above.
\end{point}
\begin{point}{60}[vn-is-dagger-category]{Corollary}%
The category~$\op\vN$ of von Neumann algebras with
    ncpu-maps in the opposite direction,
    which is an~$\&$-effectus by~\sref{vn-is-andthen-eff},
    is also a~$\dagger$-effectus.
\end{point}
\begin{point}{61}{Remarks}%
    The dagger on the pure maps of~$\vN$ (that exists by the previous
    corollary) is fixed by the following two rules.
\begin{enumerate}
\item
    For any (nmiu-)isomorphism~$\vartheta$
        we have~$\vartheta^\dagger = \vartheta^{-1}$
        as~$\vartheta$ is~$\diamond$-adjoint to~$\vartheta^{-1}$
         and pure maps are rigid, see~\sref{pure-is-rigid}.
\item
    The standard filter~$c\colon \scrA \to \ceil{b} \scrA \ceil{b}$
        for~$b \in \scrA$
        (ie.~$c(a) = \sqrt{b}a\sqrt{b}$, see~\sref{dils-def-filter})
        has as dagger~$c^\dagger\colon \ceil{b} \scrA  \ceil{b} \to \scrA$
        given by~$c^\dagger(a) = \sqrt{b} a \sqrt{b}$.
        (This follows from~$\asrt_p^\dagger = \asrt_p$
            and~\sref{dagger-of-zeta}.)
\end{enumerate}
More concretely,
    for any pure ncp-map~$\varphi\colon \scrA \to \scrB$
    we can find
    a unique isomorphism~$\vartheta \colon \ceil{\varphi} \scrA \ceil{\varphi}
    \to \ceil{\varphi(1)} \scrB \ceil{\varphi(1)}$
    such that
\begin{equation*}
    \varphi(a) \ =\  \sqrt{\varphi(1)} \,\vartheta
    \bigl(\ceil{\varphi} a \ceil{\varphi}\bigr) \,\sqrt{\varphi(1)} \qquad \text{(for all $a\in \scrA$)}
\end{equation*}
    (this follows from~\sref{pure-fundamental} and~\sref{square-f})
    and so by the previous two rules we get
\begin{equation*}
    \varphi^\dagger(b) \ = \  \ceil{\varphi} \,\vartheta^{-1} 
    \Bigl( \sqrt{\varphi(1)} \, b \, \sqrt{\varphi(1)} \Bigr) \, \ceil{\varphi}
        \qquad \text{(for all $b \in \scrB$)}.
\end{equation*}
\spacingfix{}
\begin{point}{72}%
In special cases, we can give a simpler definition of the dagger.
For instance, if the
    pure map is of the form~$\varphi\colon \scrB(\scrH) \to \scrB(\scrK)$,
    then~$\varphi = \ad_T$ for some operator~$T \colon \scrK \to \scrH$
    and in that case~$\varphi^\dagger = \ad_{T^*}$.

Another example is a pure map~$\varphi\colon \scrA \to \scrB$
    between finite-dimensional von Neumann algebras.
The Hilbert--Schmidt inner product~$[a, b] \equiv \TR a^*b$
    turns a finite-dimensional von Neumann algebra into a Hilbert space.
The adjoint of~$\varphi$ with respect to this inner product
    equals~$\varphi^\dagger$.
This can be seen by considering the special cases
    of nmiu-isomorphisms and standard filters
    (noting nmiu-isomorphisms are trace-preserving).
\end{point}
\end{point}
\begin{point}{70}{Remarks}%
    In an~$\&$-effectus
    with separating predicates
    (see \sref{dfn-mandso}),
    the second axiom of a~$\dagger'$-effectus is equivalent
    to
\begin{equation*}
    \andthen{(\andthen{p}{q})^2}{r}
    \ =\  \andthen{p}{(\andthen{q^2}{(\andthen{p}{r})})}
\end{equation*}
    for all predicates~$p,q,r$.
This is essentially the \emph{fundamental formula
    of quadratic Jordan algebras} \cite[\S4.2]{mccrimmon2006taste}
    (with~$U_x y \equiv \andthen{x^2}{y}$).

The second axiom of a~$\dagger$-effectus
    has been considered before in~$\dagger$-categories
    by Selinger \cite{selinger2008idempotents}
    who calls it the \emph{unique square root axiom}.
\end{point}
\end{point}
\end{parsec}

\begin{parsec}{2160}%
\begin{point}{10}[diamond-is-dagger-positive]{Lemma}%
In a~$\dagger$-effectus,
    a map is~$\dagger$-positive if and only if it is~$\diamond$-positive.
\begin{point}{20}{Proof}%
Assume~$f$ is~$\dagger$-positive.
By assumption 2, there is a (unique) $\dagger$-positive
    $g$ with~$f = g \after g$.
By~$\dagger$-positivity of~$g$,
    there is an~$h$ with~$g = h^\dagger \after h$.
Using the fact~$h$ is~$\diamond$-adjoint
    to~$h^\dagger$, we see~$g$ is~$\diamond$-self-adjoint:
\begin{equation*}
    g^\diamond
    \ =\  (h^\dagger \after h)^\diamond
    \ =\  h^\diamond \after (h^\dagger)^\diamond 
    \ =\  (h^\dagger)_\diamond \after h_\diamond 
    \ =\  (h^\dagger\after h)_\diamond
    \ =\  g_\diamond.
\end{equation*}
Thus~$f$ is the square of the~$\diamond$-self-adjoint map $g$,
    hence~$f$ is~$\diamond$-positive. \qed
\end{point}
\end{point}

\begin{point}{30}[dagger-eff-square-root]{Lemma}%
Predicates in a~$\dagger$-effectus have a unique square root:
    for every predicate~$p$,
    there is a unique predicate~$q$
    with~$\andthen{q}{q}=p$.
\begin{point}{40}{Proof}%
Let~$p\colon X \to 1$ be any predicate.
By assumption 3,
    the~$\diamond$-positive map~$\asrt_p$
        is also~$\dagger$-positive.
So by assumption 2,
    there is (a unique)~$\dagger$-positive map~$f$
    with~$f \after f = \asrt_p$.
Define~$q = 1 \after f$.
By \sref{diamond-is-dagger-positive}
    $f$ is~$\diamond$-positive
    and so by uniqueness of~$\diamond$-positive maps,
    we get~$f = \asrt_{1\after f} \equiv \asrt_q$.
We compute
\begin{equation*}
    \andthen{q}{q}
        \ =\    
        q \after \asrt_q \ \equiv \ 
        1 \after f \after \asrt_{1 \after f} \ = \ 
        1 \after f \after f \ = \ 1 \after \asrt_p \ = \ p,
\end{equation*}
which shows~$p$ has as square root~$q$.

To show uniqueness, assume~$p = \andthen{r}{r}$ for some
    predicate~$r$.
Note
\begin{equation*}
    \asrt_p
        \ = \ \asrt_{\andthen{r}{r}}
        \ \overset{\smash{\sref{andthen-square-rule}}}{=} \ \asrt_r \after \asrt_r.
\end{equation*}
As~$\asrt_r$ is~$\diamond$-positive,
    it is also~$\dagger$-positive by the third axiom.
    So by the second axiom~$\asrt_r = \asrt_q$.
Thus~$r = q$, which shows uniqueness of the square root.
\qed
\end{point}
\end{point}
\begin{point}{50}[asrt-iso]{Proposition}%
In an~$\&$-effectus with square roots
    (e.g.~$\dagger$- or~$\dagger'$-effectus) we have
\begin{equation*}
    \asrt_p \after \alpha
        \ =\  \alpha \after \asrt_{p\after \alpha}
\end{equation*}
for every isomorphism~$\alpha$ and predicate~$p$.
\begin{point}{60}{Proof}%
There is some~$q$ with~$\andthen{q}{q}=p$.
The map~$\alpha^{-1} \after \asrt_q \after \alpha$
    is~$\diamond$-self-adjoint:
\begin{alignat*}{2}
    (\alpha^{-1} \after \asrt_q \after \alpha)_\diamond
    &\ = \ \alpha^{-1}_\diamond \after (\asrt_q)_\diamond \after \alpha_\diamond \\
    &\ \overset{\mathclap{\sref{iso-diamond-adjoint}}}{=} \ \alpha^\diamond \after (\asrt_q)^\diamond \after (\alpha^{-1})^\diamond\\
    &\ = \ (\alpha^{-1} \after \asrt_q \after \alpha)^\diamond.
\end{alignat*}
Thus~$
\alpha^{-1} \after \asrt_q \after \alpha \after
\alpha^{-1} \after \asrt_q \after \alpha
= \alpha^{-1} \after \asrt_p \after \alpha
$ is~$\diamond$-positive.  By uniqueness of~$\diamond$-positive maps,
we get~$\alpha^{-1} \after \asrt_p \after \alpha
    = \asrt_{1 \after \alpha^{-1}\after\asrt_p \after \alpha}
    = \asrt_{p \after \alpha} $.
Postcomposing~$\alpha$, we find~$\asrt_p\after\alpha = 
        \alpha \after \asrt_{p \after \alpha}$, as desired. \qed
\end{point}
\end{point}
\begin{point}{70}[dagger-of-zeta]{Proposition}%
In a $\dagger$-effectus:~$\zeta_s^\dagger = \pi_s$
(recall convention \sref{zeta-s-convention} for~$\pi_s$ and~$\zeta_s$).
\begin{point}{80}{Proof}%
As~$\pi_s$ is~$\diamond$-adjoint to~$\pi_s^\dagger$,
we have
\begin{equation*}
\IM \pi_s^\dagger
 \ = \  (\pi^\dagger_s)_\diamond(1) \ =\ 
 (\pi_s)^\diamond(1) \ = \  \ceil{1 \after \pi_s}
 \ =\  1
\end{equation*}
and similarly~$\ceil{1 \after \pi_s^\dagger} = \IM \pi_s = s$.
As~$1 \after \pi_s^\dagger \leq \ceil{1 \after \pi_s^\dagger} = s$,
    there is some pure~$h$ with~$\pi_s^\dagger = h \after \zeta_s$.
By \sref{standard-form-map},
    there is also some pure and faithful~$g$
    with~$\zeta_s^\dagger = \pi_s \after g$.
Using~$\zeta_s \after \pi_s = \id$ twice we find
\begin{equation}\label{eq-dagger-zeta-pi-conn}
    \id \ =\  \id^\dagger
    \ = \ \pi_s^\dagger \after \zeta_s^\dagger
    \ = \ h \after \zeta_s \after \pi_s \after g 
    \ = \ h \after g.
\end{equation}
Now~$1= 1 \after h \after g \leq 1 \after g$ and so~$g$ is total.
Clearly~$\zeta_s^\dagger\after\zeta_s$ is~$\dagger$-positive,
    hence~$\diamond$-positive and so by uniqueness of~$\diamond$-positive maps:
\begin{equation*}
\zeta_s^\dagger \after \zeta_s
    \ = \ \asrt_{1 \after \zeta_s^\dagger \after \zeta_s}
    \ =\  \asrt_{1 \after \pi_s \after g\after \zeta_s}
    \ =\  \asrt_s
    \ =\  \pi_s \after \zeta_s.
\end{equation*}
So by epicity of~$\zeta_s$, we find~$\zeta_s^\dagger = \pi_s$, as desired.\qed
\end{point}
\end{point}
\begin{point}{90}[dagger-of-iso]{Corollary}%
    In a~$\dagger$-effectus,
        $\pi_s$ is~$\diamond$-adjoint to~$\zeta_s$.
    Also~$\alpha^\dagger = \alpha^{-1}$
        for any iso~$\alpha$.
\end{point}
\begin{point}{100}[zeta-through-asrt]{Exercise}%
Show that in an~$\&$-effectus
        where every predicate has a square root
        and where~$\pi_s$ is~$\diamond$-adjoint to~$\zeta_s$
        (e.g.~a~$\dagger$-effectus)
        we have
    $\asrt_p \after \zeta_s \ = \ \zeta_s \after \asrt_{p \after \zeta_s}$.
    (Hint: mimic the proof of \sref{asrt-iso}.)
\end{point}
\begin{point}{110}[dagger-thm-necessity]{Theorem}%
    A~$\dagger$-effectus is a~$\dagger'$-effectus.
\begin{point}{120}{Proof}%
Assume~$C$ is a~$\dagger$-effectus.
Axiom 1 is already proven in \sref{dagger-eff-square-root}.
\begin{point}{130}[pqqp-from-dagger]{Ax.~2}%
Let~$p,q$ be predicates.
To start, note~$1 \after \asrt_q\after  \asrt_p = \andthen{p}{q}$ and
\begin{equation*}
    \IM \asrt_q \after \asrt_p \ =\ 
        (\asrt_q)_\diamond(p) \ = \ 
        (\asrt_q)^\diamond(p) \ = \ 
        \ceil{\andthen{q}{p}}.
\end{equation*}
So using \sref{standard-form-map}
    we know
\begin{equation*}
    \asrt_q \after \asrt_p \ =\ 
    \pi_{\ceil{\andthen{q}{p}}}
    \after \alpha \after \zeta_{\ceil{\andthen{p}{q}}} \after \asrt_{\andthen{p}{q}}
\end{equation*}
for some iso~$\alpha$.
Applying the dagger to both sides
we get, using \sref{dagger-of-zeta} and \sref{dagger-of-iso}:
\begin{align*}
    \asrt_p \after \asrt_q &\ =\ 
    (\asrt_q \after \asrt_p)^\dagger \\
    & \ =\  \asrt_{\andthen{p}{q}}
    \after \zeta^\dagger_{\ceil{\andthen{p}{q}}}
    \after \alpha^\dagger \after
    \pi^\dagger_{\ceil{\andthen{q}{p}}} \\
    & \ = \ 
    \asrt_{\andthen{p}{q}}
    \after \pi_{\ceil{\andthen{p}{q}}}
    \after \alpha^{-1} \after
    \zeta_{\ceil{\andthen{q}{p}}}.
\end{align*}
Combining both:
\begin{align*}
    & \asrt_p \after \asrt_q^2 \after \asrt_p \\
    &\qquad = \ \asrt_{\andthen{p}{q}}
    \after \pi_{\ceil{\andthen{p}{q}}}
    \after \alpha^{-1} \after
    \zeta_{\ceil{\andthen{q}{p}}}
    \after \pi_{\ceil{\andthen{q}{p}}}
    \after \alpha \after \zeta_{\ceil{\andthen{p}{q}}} \after \asrt_{\andthen{p}{q}} \\
    &\qquad = \ \asrt_{\andthen{p}{q}}
    \after \pi_{\ceil{\andthen{p}{q}}}
    \after \zeta_{\ceil{\andthen{p}{q}}} \after \asrt_{\andthen{p}{q}} \\
    &\qquad = \ \asrt_{\andthen{p}{q}}
    \after \asrt_{\ceil{\andthen{p}{q}}}
    \after \asrt_{\andthen{p}{q}} \\
    &\qquad \overset{\smash{\mathclap{\sref{asrt-absorp-rule}}}}{=} \ 
    \asrt_{\andthen{p}{q}}^2,
\end{align*}
as desired.
\end{point}
\begin{point}{140}{Ax.~3}%
Pick any sharp predicates~$s,t$.
We want to show~$t \after \zeta_s$ is sharp.
To this end, we will show~$t \after \zeta_s$
    is the image of~$\pi_s \after \pi_t$.
Clearly~$t \after \zeta_s \after \pi_s \after \pi_t = 1$.
Let~$p$ be any sharp predicate with~$p \after \pi_s \after \pi_t = 1$.
Then~$p \after \pi_s \geq \IM \pi_t = t$.
So~$  p^\perp \after \pi_s \leq t^\perp$
and hence~$\ceil{p^\perp \after \pi_s} \leq t^\perp$.
As~$\pi_s$ is~$\diamond$-adjoint to~$\zeta_s$ by \sref{dagger-of-iso},
we get~$t \after \zeta_s \leq \ceil{t \after \zeta_s} \leq p$,
    which shows~$t \after \zeta_s$
    is the image of~$\pi_s \after \pi_t$
    and consequently sharp. \qed
\end{point}
\end{point}
\end{point}

\end{parsec}

\begin{parsec}{2170}%
\begin{point}{10}%
Let~$f$ be a pure map in a~$\dagger'$-effectus.
We will work towards the definition of~$f^\dagger$.
By \sref{standard-form-map} we
know there is an iso~$\alpha$ with
\begin{equation*}
f  \ =\   \pi_{\IM f} \after \alpha \after \zeta_{\ceil{1\after f}}
            \after \asrt_{1 \after f}.
\end{equation*}
In a~$\dagger$-effectus
    we have~$\asrt_p^\dagger = \asrt_p$,
    $\zeta_s^\dagger = \pi_s$,
    $\alpha^\dagger = \alpha^{-1}$
    and~$\pi_s^\dagger = \zeta_s$ 
    (for corresponding~$\zeta_s$ and~$\pi_s$)
    ---
    so we are forced to define~
\begin{equation}\label{dagger-definition}
    f^\dagger
        \ =\  \asrt_{1 \after f} \after
    \pi_{\ceil{1\after f}} \after
    \alpha^{-1} \after
    \zeta_{\IM f},
\end{equation}
    where~$\zeta_{\IM f}$
    is the unique corresponding
        quotient of~$\pi_{\IM f}$
        and~$\pi_{\ceil{1 \after f}}$
        the unique corresponding comprehension of
        and~$\zeta_{\ceil{1 \after f}}$,
        see \sref{zeta-s-convention}.
Before we declare~\eqref{dagger-definition} a definition,
 we have to check whether it is independent
    of choice of~$\pi$ (and corresponding $\zeta$).
    So suppose~$ f  =  \pi' \after \alpha' \after \zeta'
        \after \asrt_{1 \after f}$
        for some iso~$\alpha'$, comprehension~$\pi'$ of~$\IM f$
        and quotient~$\zeta'$ of~$\smash{\ceil{1 \after f}^\perp}$.
There are isos~$\beta$ and~$\gamma$
    such that~$\pi' = \pi_{\IM f} \after \beta$
    and~$\zeta' = \gamma \after \zeta_{\ceil{1 \after f}}$.
We will take a moment to relate~$\alpha$ and~$\alpha'$:
as~$\pi_{\IM f}$ is mono
    and~$\zeta_{1 \after f}\after\asrt_{1 \after f}$
    is epic,
    we have~$\beta\after\alpha'\after\gamma = \alpha$
    and so~$(\alpha')^{-1}= \gamma \after\alpha^{-1}\after\beta$.
To continue,
it is easy to see~$\beta^{-1} \after \zeta_{\IM f}$
    is the unique corresponding quotient to~$\pi'$
    and~$\pi_{\ceil{1 \after f}}\after \gamma^{-1}$
    is the unique corresponding comprehension to~$\zeta'$.
So with this choice of quotient and comprehension,
    we are forced to define
\begin{align*}
    f^\dagger 
    &\ = \ \asrt_{1 \after f}\after \pi_{\ceil{1 \after f}}
    \after \gamma^{-1} \after {\alpha'}^{-1}
                \after \beta^{-1}\after \zeta_{\IM f} \\
    &\ = \ \asrt_{1 \after f}\after \pi_{\ceil{1 \after f}}
                \after \gamma^{-1} \after \gamma \after \alpha^{-1} \after \beta 
                \after \beta^{-1}\after \zeta_{\IM f} \\
    &\ = \ \asrt_{1 \after f}\after \pi_{\ceil{1 \after f}}
                \after \alpha^{-1} \after \zeta_{\IM f},
\end{align*}
which is indeed consistent with~\eqref{dagger-definition}.
So we are justified to declare:
\end{point}
\begin{point}{20}[dagger-definition2]{Definition}%
In a~$\dagger'$-effectus, for a pure map~$f$, define
    \index{*dagg@$\dagger$, dagger!in a $\dagger'$-effectus}
\begin{equation*}
    \Define{f^\dagger}
        \ =\  \asrt_{1 \after f} \after
    \pi_{\ceil{1\after f}} \after
    \alpha^{-1} \after
    \zeta_{\IM f},
\end{equation*}
where~$\alpha$ is the unique iso such that
$f  =  \pi_{\IM f} \after \alpha \after \zeta_{\ceil{1\after f}}
            \after \asrt_{1 \after f}$.
\end{point}
\begin{point}{30}[dagger-prime-basics]{Exercise}%
Show that in a~$\dagger'$-effectus,
    we have
\begin{align*}
    \asrt_p^\dagger &= \asrt_p &
    \pi_s^\dagger &= \zeta_s &
    \zeta_s^\dagger &= \pi_s &
    \alpha^\dagger &= \alpha^{-1}
\end{align*}
for a quotient~$\zeta_s$ corresponding to~$\pi_s$
and iso~$\alpha$.
\end{point}
\end{parsec}

\begin{parsec}{2180}%
\begin{point}{10}%
To compute~$f^{\dagger\dagger}$,
    we need to put~$f^\dagger$
    in the standard form of~\sref{standard-form-map}.
To do this, we need to pull~$\asrt_p$ from one side to the other.
In this section we will work towards a general result for this.
\end{point}
\begin{point}{20}[quotcompr-diamond-adjoint]{Lemma}%
In a~$\dagger'$-effectus,
$\pi_s$ is~$\diamond$-adjoint to~$\zeta_s$.
\begin{point}{30}{Proof}%
    Pick any sharp~$t,u$.
We have to
    show~$t \after \zeta_s \leq u^\perp$
    if and only if~$u \after \pi_s \leq t^\perp$.
So assume~$t \after \zeta_s \leq u^\perp$.
Then~$t =  t\after \zeta_s \after \pi_s  \leq  u^\perp \after \pi_s
                 =  (u \after \pi_s)^\perp$
    so~$u \after \pi_s \leq t^\perp$.
For the converse, assume~$u \after \pi_s \leq t^\perp$.
Then~$u \after \asrt_s = u \after \pi_s \after \zeta_s \leq t^\perp \after \zeta_s$.
From this,
    the
    $\diamond$-self-adjointness of~$\asrt_s$
    and the fact that~$t^\perp \after \zeta_s$ is sharp,
    we find~$(t^\perp \after \zeta_s)^\perp \after \pi_s \after \zeta_s \leq
        u^\perp$
    and so
\begin{alignat*}{2}
    t \after \zeta_s
    &\ =\  (t^\perp \after \id)^\perp \after \zeta_s \\
    &\ =\  (t^\perp \after \zeta_s \after \pi_s)^\perp \after \zeta_s \\
    &\ =\  (t^\perp \after \zeta_s)^\perp \after \pi_s \after \zeta_s \\
    &\ \leq \ u^\perp,
\end{alignat*}
as desired. \qed
\end{point}
\end{point}
\begin{point}{40}[dfn-pristine]{Definition}%
    In an~$\&$-effectus, a map~$f$ is \Define{pristine}\index{pristine}
        if it is pure and~$1\after f$ is sharp.
\begin{point}{50}{Remarks}%
    In general, pristine maps are not closed under composition:
    in~$\op\vN$, the maps~$\asrt_{\ketbra{0}{0}}$
            and~$\asrt_{\ketbra{+}{+}}$ on~$M_2$
            are both pristine,
            but their composite is not.
However, with a non-standard composition
    of pristine maps
    (namely~$f \cdot g \equiv \asrt_{\floor{1 \after f \after g}}
    \after f \after g$),
    the pristine maps are closed under composition and even
    form a~$\dagger$-category \cite{effintro}.
Essentially the same construction
    has been found independently by~\cite{hines2010structure}
    for partial isometries between Hilbert spaces.
\end{point}
\end{point}
\begin{point}{60}[standard-form-pristine]{Exercise}%
Show that in an~$\&$-effectus, every pristine map~$h$ is of the form
        \begin{equation*}
            h \ =\ \pi_{\IM h} \after \alpha \after \zeta_{1 \after h}
        \end{equation*}
        for some iso~$\alpha$.
\end{point}
\begin{point}{70}[pristine-asrt]{Proposition}%
In a~$\dagger'$-effectus,
    with pristine map~$h$, we have
\begin{equation*}
    \asrt_p \after h
        \ =\  h \after \asrt_{p \after h}
\end{equation*}
for any predicate~$p$ with~$p \leq \IM h$.
\begin{point}{80}{Proof}%
For brevity,
write~$t = 1 \after h$
    and~$s = \IM h$.
By~\sref{standard-form-pristine}
    there is some iso~$\alpha$
    such that~$h = \pi_s \after \alpha \after \zeta_t$.
    Note~$\IM \asrt_p = \ceil{p } \leq \ceil{s} = s$
    and so~$\asrt_p = \asrt_{s} \after \asrt_p$
    by the first rule of \sref{asrt-absorp-rule}.
By the second rule of \sref{asrt-absorp-rule}
    and~$1 \after p = p \leq s$,
    we see~$p = p \after \asrt_{s}$.
Thus
\begin{alignat*}{2}
    \asrt_p \after \pi_{s}
    &\ =\  \asrt_{s} \after \asrt_{p \after \asrt_{s}} \after \pi_{s} \\
    &\ =\  \pi_{s} \after \zeta_{s} \after
    \asrt_{p \after \pi_{s} \after \zeta_{s}} \after \pi_{s} \\
    &\ =\  \pi_{s}  \after
    \asrt_{p \after \pi_{s}} \after \zeta_{s} \after \pi_{s} 
    &\qquad& \text{by \sref{zeta-through-asrt}
                    and \sref{quotcompr-diamond-adjoint} }\\
    &\ =\  \pi_{s}  \after
    \asrt_{p \after \pi_{s}}.
\end{alignat*}
Putting everything together
\begin{alignat*}{2}
   \asrt_p \after h
   & \ =\  \asrt_p \after \pi_s \after \alpha \after \zeta_t \\
   & \ =\  \pi_s \after \asrt_{p\after\pi_s} \after \alpha \after \zeta_t  \\
   & \ =\  \pi_s \after \alpha \after \asrt_{p\after\pi_s \after \alpha} \after \zeta_t 
   &\qquad&\text{by \sref{asrt-iso}}\\
   & \ =\  \pi_s \after \alpha \after \zeta_t \after \asrt_{p\after\pi_s \after \alpha \after \zeta_t} 
    &\qquad& \text{by \sref{zeta-through-asrt}
                    and \sref{quotcompr-diamond-adjoint} }\\
   & \ = \ h \after \asrt_{p \after h},
\end{alignat*}
as desired.\qed
\end{point}
\end{point}
\begin{point}{90}[asrt-pristine-reverse]{Exercise}%
Working in a~$\dagger'$-effectus, show in order
\begin{enumerate}
    \item if~$h \equiv \pi_{\IM h} \after \alpha \after \zeta_{1 \after h}$
            is some pristine map,
            then~$h^\dagger = \pi_{1 \after h} \after \alpha^{-1} \after \zeta_{\IM h}$;
    \item $h^{\dagger\dagger} = h$ for any pristine~$h$;
    \item $h^\dagger \after h = \asrt_{1 \after h}$ for any pristine map~$h$;
    \item $p \after h^\dagger \leq \IM h$ for any predicate~$p$
        and pristine map~$h$ \emph{and}
    \item
        if~$p \leq 1 \after h$,
        then~$\asrt_{p\after h^\dagger} \after h = h \after \asrt_p$
        for any pristine map~$h$.
\end{enumerate}
\end{point}
\spacingfix{}
\begin{point}{100}[prist-asrt-decomp]{Proposition}%
In a~$\dagger'$-effectus,
        for every pure map~$f$,
    there exists a unique pristine map~$h$
    with~$1 \after h = \ceil{1 \after f}$
    and~$f = h \after \asrt_{1 \after f}$.
    Furthermore~$f^\dagger = \asrt_{1 \after f}\after h^\dagger$.
\begin{point}{110}{Proof}%
By~\sref{standard-form-map}
 we know~$f = \pi_{\IM f} \after \alpha \after \zeta_{\ceil{1 \after f}}
        \after \asrt_{1 \after f}$
        for some iso~$\alpha$.
Define~$h = \pi_{\IM f} \after \alpha \after \zeta_{\ceil{1 \after f}}$.
Clearly~$1 \after h = \ceil{1 \after f}$
    and $f = h \after \asrt_{1 \after f}$.
    Also
    \begin{alignat*}{2}
        f^\dagger &\ \equiv\ 
        \asrt_{1 \after f}
            \after \pi_{\ceil{1 \after f}}
            \after \alpha^{-1}
            \after \zeta_{\IM f} \\
        &\ = \ 
        \asrt_{1 \after f}
        \after \asrt_{\ceil{1 \after f}}
            \after \pi_{\ceil{1 \after f}}
            \after \alpha^{-1} 
            \after \zeta_{\IM f} 
            &\qquad&\text{by \sref{asrt-absorp-rule}}
            \\
        &\ \equiv \ 
        \asrt_{1 \after f}
        \after h^\dagger.
\end{alignat*}
Only uniqueness remains.
Assume~$f = h' \after \asrt_{1 \after f}$
    for some pristine
    map~$h'$ with~$\ceil{1 \after h'} = \ceil{1 \after f}$.
By \sref{standard-form-map} and \sref{asrt-absorp-rule},
    there is an iso~$\alpha'$
    with~$h' = \pi_{\IM h} \after \alpha' \after \zeta_{\ceil{1 \after f}}$.
As~$\zeta_{\ceil{1 \after f}} \after \asrt_{1 \after f}$
    is a quotient (by \sref{zeta-asrt-quot}),
    quotients are faithful and~$f = h' \after \asrt_{1 \after f}$,
    we see~$\IM h' = \IM f$
    and so~$\alpha = \alpha'$ by \sref{standard-form-map}.\qed
\end{point}
\end{point}
\begin{point}{120}[dagger-idempotent]{Proposition}%
    In a~$\dagger'$-effectus, we have~$f^{\dagger\dagger}=f$
        for any pure map~$f$.
\begin{point}{130}{Proof}%
By~\sref{prist-asrt-decomp}
    we have~$f = h \after \asrt_{1\after f}$
    and~$f^\dagger = \asrt_{1 \after f} \after h^\dagger$
    for some pristine~$h$
    with~$1 \after h = \ceil{1 \after f}$.
Clearly~$1 \after f \leq \ceil{1 \after f} = 1 \after h = \IM h^\dagger$,
    so by \sref{pristine-asrt}
    we get~$f^\dagger = \asrt_{1 \after f} \after h^\dagger
              = h^\dagger \after \asrt_{1 \after f \after h^\dagger}$.
Consequently
\begin{equation*}
    f^{\dagger\dagger}
    \ =\  \asrt_{1 \after f\after h^\dagger} \after h^{\dagger\dagger}
    \ \overset{\sref{asrt-pristine-reverse}}{=}\  \asrt_{1 \after f\after h^\dagger} \after h
    \ \overset{\sref{asrt-pristine-reverse}}{=}\  h\after \asrt_{1 \after f}
    \ = \ f,
\end{equation*}
as desired. \qed
\end{point}
\end{point}
\end{parsec}

\begin{parsec}{2190}%
\begin{point}{10}%
Now we tackle the most tedious part
    of \sref{dagger-theorem}:
    we will show~$(f \after g)^\dagger = g^\dagger \after f^\dagger$
    in a~$\dagger'$-effectus.
To avoid too much repetition, let us fix the setting.
\end{point}%
\begin{point}{20}[dagger-setting]{Setting}%
Let~$f,g$ be two composable pure maps in a~$\dagger'$-effectus.
For brevity, write~$p = 1 \after f$, $q = 1 \after g$,
    $s = \IM f$ and $t = \IM g$.
Let~$\varphi$ and~$\psi$
be the unique isomorphisms (see \sref{standard-form-map}) such that
\begin{equation*}
    f \ =\  \pi_s \after \varphi \after \zeta_{\ceil{p}} \after \asrt_p 
    \qquad \text{and} \qquad
    g \ =\  \pi_t \after \psi \after \zeta_{\ceil{q}} \after \asrt_q.
\end{equation*}
Define~$
    h = \pi_s \after \varphi \after \zeta_{\ceil{p}}$
    and~$k = \pi_t \after \psi \after \zeta_{\ceil{q}}$.
To compute~$(f\after g)^\dagger$,
    we have to put~$f\after g$
    in the standard form
    of~\sref{standard-form-map}.
We will do this step-by-step,
    first we put~$\zeta_{\ceil{p}} \after \asrt_p \after \pi_t$
        in standard form.
\begin{align*}
1 \after \zeta_{\ceil{p}} \after \asrt_p \after \pi_t
&\ =\  \ceil{p} \after \asrt_p \after \pi_t
    \ \overset{\sref{asrt-absorp-rule}}{=}\  p \after \pi_t\\
    \IM \zeta_{\ceil{p}}  \after \asrt_p \after \pi_t &\ = \ 
(\zeta_{\ceil{p}}  \after \asrt_p)_\diamond(t)
\ \overset{\sref{quotcompr-diamond-adjoint}}{=} \ \ceil{t \after \asrt_p \after \pi_{\ceil{p}}}.
\end{align*}
Thus there is a unique isomorphism~$\chi$ with
\begin{equation}\label{dagger-iso-chi}
    \zeta_{\ceil{p}} \after \asrt_p \after \pi_t
    \ = \ \pi_{\lceil t \after \asrt_p \after \pi_{\ceil{p}}\rceil} \after \chi \after \zeta_{\ceil{p \after \pi_t}} \after \asrt_{p \after \pi_t}.
\end{equation}
Next, we consider~$\asrt_{p \after k}
                        \after \asrt_q$, clearly
\begin{equation*}
    1 \after \asrt_{p \after k}
\after \asrt_q 
\ =\ p \after k
\after \asrt_q \ = \ p \after g.
\end{equation*}
Concerning the image, first
    note~$p \after k = p \after \pi_t \after \psi \after \zeta_{\ceil{q}}
    \leq 1 \after \zeta_{\ceil{q}} = \ceil{q}$
        and so we must have~$\IM \asrt_{p \after k} \leq \ceil{q}$, which implies
\begin{equation*}
\IM \asrt_{p \after k}
\after \asrt_q
\ = \ (\asrt_{p \after k})_\diamond(q)
\ =\  \ceil{\ceil{q} \after \asrt_{p \after k}}
\ =\  \ceil{p \after k}.
\end{equation*}
So there is a unique isomorphism~$\omega$ such that
\begin{equation}\label{dagger-iso-omega}
\asrt_{p \after k} \after \asrt_q
    \ =\   \pi_{\lceil p \after k \rceil}
    \after \omega \after \zeta_{\ceil{p \after g}} \after \asrt_{p \after g}.
\end{equation}
Next, we consider~$
\zeta_{\ceil{p \after \pi_t}} \after \psi \after
                \zeta_{\ceil{q}}$.
Note~$\psi \after \zeta_{\ceil{q}}$
    is a quotient for a sharp predicate,
    hence sharp and so~$\ceil{p \after \pi_t} \after \psi \after \zeta_{\ceil{q}}
                =\lceil p \after \pi_t \after \psi
                \after \zeta_{\ceil{q}}\rceil = \lceil p \after k \rceil$
                by \sref{sharp-ceil}.
As quotients are closed under composition
(\sref{quotients-composition}),
    there is an iso~$\beta$
    with
    \begin{equation}\label{dagger-iso-beta}
    \zeta_{\ceil{p \after \pi_t}} \after \psi \after \zeta_{\ceil{q}}
    \ = \ \beta \after \zeta_{\ceil{p \after k}}.
\end{equation}
Finally, we deal with~$
\pi_s \after \varphi \after 
\pi_{\lceil t \after \asrt_p \after \pi_{\ceil{p}}\rceil}
$.
By \sref{upm-closed} this is again a comprehension.
We compute
\begin{align*}
    \IM
\pi_s \after \varphi \after 
\pi_{\lceil t \after \asrt_p \after \pi_{\ceil{p}}\rceil}
& \ = \ 
(\pi_s \after \varphi \after 
\pi_{\lceil t \after \asrt_p \after \pi_{\ceil{p}}\rceil})_\diamond(1)
\\
& \ \overset{\smash{\mathclap{\sref{quotcompr-diamond-adjoint}}}}{=}\ 
\zeta_s^\diamond ( \varphi_\diamond (
\lceil t \after \asrt_p \after \pi_{\ceil{p}}\rceil)) \\
& \ =\ 
\lceil t \after \asrt_p \after \pi_{\ceil{p}} \rceil
\after \varphi^{-1} \after \zeta_s \\
& \ \overset{\smash{\mathclap{\sref{sharp-ceil}}}}{=}\ 
\lceil t \after \asrt_p \after \pi_{\ceil{p}}
\after \varphi^{-1} \after \zeta_s \rceil \\
& \ =\ 
\lceil t \after f^\dagger\rceil
\end{align*}
and so there must be a unique iso~$\alpha$ with
\begin{equation}\label{dagger-iso-alpha}
    \pi_s \after \varphi \after \pi_{\lceil t \after \asrt_p \after
    \pi_{\ceil{p}} \rceil} = \pi_{\ceil{t \after f^\dagger}} \after \alpha.
\end{equation}
\end{point}
\spacingfix{}
\begin{point}{30}{Lemma}%
In setting \sref{dagger-setting}, we have
    $f \after g  =  \pi_{\ceil{t \after f^\dagger}}
        \after \alpha \after \chi \after \beta \after \omega
        \after \zeta_{\ceil{p \after g}} \after
        \asrt_{p \after g}$.
\begin{point}{40}{Proof}%
It's a long, but easy verification,
    either with a diagram
\begin{equation*}
    \xymatrix@C+3pc {
        \bullet \ar[rr]^g
        \ar[rd]|{\asrt_q}
        \ar[ddd]_{\rotatebox{90}{$\scriptstyle\omega\after\zeta_{\ceil{p \after g}} \after \asrt_{p \after g}$}}
        && \bullet \ar[rr]^f
            \ar[rd]|{\zeta_{\ceil{p}}\after \asrt_p}
        && \bullet
            \\ \ar@{}[rd]|{\text{\eqref{dagger-iso-omega}}}
            & \bullet
                \ar[r]^{\psi \after \zeta_{\ceil{q}}}
                \ar[d]|{\asrt_{p \after k}}
                \ar@{}[rd]|{\text{\sref{pristine-asrt}}}
 & \bullet
            \ar[u]^{\pi_t}
            \ar[d]|{\asrt_{p \after \pi_t}}
            & \bullet \ar[ru]_{\pi_s\after\varphi}
            \ar@{}[rd]|{\text{\eqref{dagger-iso-alpha}}}
            \\& \bullet \ar[r]^{\psi \after \zeta_{\ceil{q}}}
                        \ar[d]|{\zeta_{\ceil{p \after k}}}
                \ar@{}[rd]|{\text{\eqref{dagger-iso-beta}}}
            &\bullet \ar[d]|{\zeta_{\ceil{p \after \pi_t}}}
                        \ar@{}[ru]|{\text{\eqref{dagger-iso-chi}}}
                        &&
            \\ \bullet \ar@{=}[r]
            \ar[ru]^{\pi_{\ceil{p \after k}}}
            &\bullet \ar[r]_\beta
            &\bullet \ar[rr]_\chi
            &&\bullet \ar[uuu]_{\rotatebox{90}{$\scriptstyle\pi_{\ceil{t \after f^\dagger}} \after \alpha$}}
                        \ar[luu]|{\pi_{\lceil t \after \asrt_p \after \pi_{\ceil{p}} \rceil}}
        }
\end{equation*}
or with equational reasoning
\begin{align*}
   f \after g
    & \ = \ \pi_s \after \varphi \after \zeta_{\ceil{p}}
        \after \asrt_p \after \pi_t \after \psi \after \zeta_{\ceil{q}}
        \after \asrt_q \\
    & \ \overset{\smash{\mathclap{\eqref{dagger-iso-chi}}}}{=} \ 
        \pi_s \after \varphi \after
        \pi_{\lceil t \after \asrt_p \after \pi_{\ceil{p}} \rceil}
        \after \chi \after \zeta_{\ceil{p \after \pi_t}}
        \after \asrt_{p \after \pi_t}
    \after \psi \after 
        \zeta_{\ceil{q}}
        \after \asrt_q \\
    & \ \overset{\smash{\mathclap{\sref{pristine-asrt}}}}{=} \ 
        \pi_s \after \varphi \after
        \pi_{\lceil t \after \asrt_p \after \pi_{\ceil{p}} \rceil}
        \after \chi \after \zeta_{\ceil{p \after \pi_t}}
        \after \psi \after \zeta_{\ceil{q}}
        \after \asrt_{p \after k}
        \after \asrt_q \\
    & \ \overset{\smash{\mathclap{\eqref{dagger-iso-beta}}}}{=} \ 
        \pi_s \after \varphi \after
        \pi_{\lceil t \after \asrt_p \after \pi_{\ceil{p}} \rceil}
        \after \chi
        \after \beta \after \zeta_{\ceil{p \after k}}
        \after \asrt_{p \after k}
        \after \asrt_q \\
    & \ \overset{\smash{\mathclap{\eqref{dagger-iso-alpha}}}}{=} \ 
        \pi_{\ceil{t \after f^\dagger}}\after\alpha
        \after \chi
        \after \beta \after \zeta_{\ceil{p \after k}}
        \after \asrt_{p \after k}
        \after \asrt_q \\
    & \ \overset{\smash{\mathclap{\eqref{dagger-iso-omega}}}}{=} \ 
        \pi_{\ceil{t \after f^\dagger}}\after\alpha
        \after \chi
        \after \beta \after \zeta_{\ceil{p \after k}}
        \after \pi_{\ceil{p \after k}} \after \omega
            \after \zeta_{\ceil{p \after g}}
            \after \asrt_{p \after g} \\
    & \ = \ 
        \pi_{\ceil{t \after f^\dagger}}\after\alpha
        \after \chi
        \after \beta \after
        \omega
            \after \zeta_{\ceil{p \after g}}
            \after \asrt_{p \after g},
\end{align*}
whichever the Reader might prefer. \qed
\end{point}
\end{point}
\begin{point}{50}[dagger-of-fg]{Corollary}%
    $(f \after g)^\dagger = \asrt_{p \after g}
                \after \pi_{\ceil{p \after g}}
                \after \omega^{-1}
                \after \beta^{-1}
                \after \chi^{-1}
                \after \alpha^{-1}
                \after \zeta_{\ceil{t \after f^\dagger}} $ .
\begin{point}{60}%
To show~$(f \after g)^\dagger = g^\dagger \after f^\dagger$,
    it is sufficient to proof that
    the `daggered' version
    of each of the subdiagrams
        \eqref{dagger-iso-alpha},
        \eqref{dagger-iso-beta},
        \eqref{dagger-iso-chi},
        \eqref{dagger-iso-omega} and
        \sref{pristine-asrt} of the above diagram holds.
We start with the simple ones.
\end{point}
\end{point}
\begin{point}{70}[dagger-iso-beta2]{Lemma}%
In setting \sref{dagger-setting},
    the daggered version of
        \eqref{dagger-iso-beta} holds --- that is:
    \begin{equation*}
        \pi_{\ceil{p \after k}} \after \beta^{-1}
            \ =\  \pi_{\ceil{q}} \after \psi^{-1} \after \pi_{\ceil{p \after \pi_t}}.
    \end{equation*}
\spacingfix{}
\begin{point}{80}{Proof}%
The map~$
\pi_{\ceil{q}} \after \psi^{-1} \after \pi_{\ceil{p \after \pi_t}} \after \beta$
is a comprehension for~$\ceil{p \after k}$ --- indeed
\begin{align*}
    (\pi_{\ceil{q}} \after \psi^{-1} \after \pi_{\ceil{p \after \pi_t}}\after\beta)_\diamond(1)
    &\ = \ 
    (\zeta^\diamond_{\ceil{q}} \after \psi^\diamond) (\ceil{p \after \pi_t})\\
    &\ = \ \ceil{p \after \pi_t} \after \zeta_{\ceil{q}} \after \psi \\
    &\ = \ \ceil{p \after \pi_t \after \zeta_{\ceil{q}} \after \psi} \\
    &\ = \ \ceil{p \after k}.
\end{align*}
Furthermore
\begin{align*}
    &\zeta_{\ceil{p \after k}} \after
    \pi_{\ceil{q}} \after \psi^{-1} \after \pi_{\ceil{p \after \pi_t}}
    \after \beta \\
    & \qquad \ \overset{\mathclap{\eqref{dagger-iso-beta}}}{=}\  
    \beta^{-1}\after \zeta_{\ceil{p \after \pi_t}}
    \after \psi \after \zeta_{\ceil{q}}
    \after \pi_{\ceil{q}} \after \psi^{-1} \after \pi_{\ceil{p \after \pi_t}}
    \after \beta \\
    & \qquad \ = \ \id.
\end{align*}
So~$
\pi_{\ceil{q}} \after \psi^{-1} \after \pi_{\ceil{p \after \pi_t}}
\after \beta
$ is the unique comprehension corresponding to~$\zeta_{\ceil{p \after k}}$ --- that is:~$
\pi_{\ceil{q}} \after \psi^{-1} \after \pi_{\ceil{p \after \pi_t}}
\after \beta \ = \ \pi_{\ceil{p \after k}} $,
as desired. \qed
\end{point}
\end{point}
\begin{point}{90}[dagger-iso-alpha2]{Exercise}%
Show that in the setting \sref{dagger-setting},
    the daggered version of
        \eqref{dagger-iso-alpha} holds --- i.e.
    \begin{equation*}
        \alpha^{-1} \after \zeta_{\ceil{t \after f^\dagger}}
        \ =\   \zeta_{\lceil t \after \asrt_p \after \pi_{\ceil{p}} \rceil}
        \after \varphi^{-1} \after \zeta_s.
    \end{equation*}
(Hint: mimic the proof of~\sref{dagger-iso-beta2}.)
\end{point}
\begin{point}{100}[dagger-iso-zeta2]{Exercise}%
Show that in the setting \sref{dagger-setting},
we have
\begin{equation*}
    \pi_{\ceil{q}} \after \psi^{-1} \after \asrt_{p \after \pi_t}
    \ = \ \asrt_{p \after k} \after \pi_{\ceil{q}} \after \psi^{-1}.
\end{equation*}
(This is the daggered version of the subdiagram
marked~\sref{pristine-asrt}.)
\end{point}
\begin{point}{110}[dagger-iso-mu]{Proposition}%
If in a~$\dagger'$-effectus,
$\nu$ is the unique iso (cf.~\sref{pqqp-from-dagger}) such that
\begin{equation*}
    \asrt_a \after \asrt_b\  =\  \pi_{\ceil{\andthen{a}{b}}} 
    \after \nu \after \zeta_{\ceil{\andthen{b}{a}}} \after \asrt_{\andthen{b}{a}},
\end{equation*}
then~$\asrt_b \after \asrt_a\  =\  \asrt_{\andthen{b}{a}}
        \after \pi_{\ceil{\andthen{b}{a}}} 
        \after \nu^{-1} \after \zeta_{\ceil{\andthen{a}{b}}}$.
\begin{point}{120}{Proof}%
Let~$\mu$ be the unique iso with~$
\asrt_b \after \asrt_a\  =\  
        \pi_{\ceil{\andthen{b}{a}}} 
        \after \mu \after \zeta_{\ceil{\andthen{a}{b}}}
                \after \asrt_{{\andthen{a}{b}}}$.
We will see~$\mu = \nu^{-1}$.
For brevity, write
\begin{align*}
    \overline{\andthen{a}{b}} &\ = \ (\andthen{a}{b}) \after \pi_{\ceil{\andthen{a}{b}}} &
\overline{\andthen{b}{a}} &\ = \ (\andthen{b}{a}) \after \pi_{\ceil{\andthen{b}{a}}}
\end{align*}
By \sref{asrt-pristine-reverse} and~\sref{asrt-absorp-rule}, we have
\begin{equation*}
    \pi_{\ceil{\andthen{a}{b}}} \after
\asrt_{\overline{\andthen{a}{b}}} \after
 \zeta_{\ceil{\andthen{a}{b}}}
    \ = \ 
    \pi_{\ceil{\andthen{a}{b}}} \after
 \zeta_{\ceil{\andthen{a}{b}}} \after
\asrt_{\andthen{a}{b}}
    \ = \ 
\asrt_{\andthen{a}{b}}.
\end{equation*}
Now the second axiom of a~$\dagger'$-effectus comes into play
\begin{align*}
    &\pi_{\ceil{\andthen{a}{b}}} \after
    \asrt^2_{\overline{\andthen{a}{b}}} \after \zeta_{\ceil{\andthen{a}{b}}} \\
    &\qquad \ = \ 
    \pi_{\ceil{\andthen{a}{b}}} \after
\asrt_{\overline{\andthen{a}{b}}} \after
 \zeta_{\ceil{\andthen{a}{b}}} \after
\pi_{\ceil{\andthen{a}{b}}} \after
\asrt_{\overline{\andthen{a}{b}}}
\after \zeta_{\ceil{\andthen{a}{b}}} \\
    &\qquad \ = \ 
\asrt_{\andthen{a}{b}}^2
\\
    &\qquad \ = \ 
\asrt_{a} \after
\asrt^2_{b}\after
\asrt_{a}
\\
    &\qquad \ = \ 
\pi_{\ceil{\andthen{a}{b}}} \after
    \nu \after
    \asrt_{\overline{\andthen{b}{a}}} \after
    \zeta_{\ceil{\andthen{b}{a}}} \after
\pi_{\ceil{\andthen{b}{a}}} \after
    \mu \after
    \asrt_{\overline{\andthen{a}{b}}} \after
    \zeta_{\ceil{\andthen{a}{b}}} 
\\
    &\qquad \ = \ 
\pi_{\ceil{\andthen{a}{b}}} \after
    \nu \after
    \asrt_{\overline{\andthen{b}{a}}} \after
    \mu \after
    \asrt_{\overline{\andthen{a}{b}}} \after
    \zeta_{\ceil{\andthen{a}{b}}}.
\end{align*}
Thus as quotients are epis and comprehensions are monos:
\begin{equation} \label{dagger-second-axiom-intermediate}
\asrt^2_{\overline{\andthen{a}{b}}}
         \ = \ \nu \after \asrt_{\overline{\andthen{b}{a}}}
            \after \mu \after \asrt_{\overline{\andthen{a}{b}}}.
\end{equation}
We want to show~$\asrt_{\overline{\andthen{a}{b}}}$ is an epi.
First note
\begin{equation*}
    \ceil{\overline{\andthen{a}{b}}}\after\zeta_{\ceil{\andthen{a}{b}}}
    \ \overset{\sref{sharp-ceil}}{=} \ 
    \ceil{\overline{\andthen{a}{b}}\after\zeta_{\ceil{\andthen{a}{b}}}}
    \ \overset{\sref{asrt-absorp-rule}}{=} \ 
    \ceil{\andthen{a}{b}} \ = \ 1 \after \zeta_{\ceil{\andthen{a}{b}}}
\end{equation*}
and so~$\IM \asrt_{\overline{\andthen{a}{b}}}
= \ceil{\overline{\andthen{a}{b}}} = 1 $,
    which tells us~$\asrt_{\overline{\andthen{a}{b}}}$
    is a quotient and therefore an epi.
    So, from \eqref{dagger-second-axiom-intermediate} we get
    $ \asrt_{\overline{\andthen{a}{b}}}
     =  \nu \after \asrt_{\overline{\andthen{b}{a}}} \after \mu$ and so
\begin{equation}\label{dagger-seqprod-inversion}
    \overline{\andthen{a}{b}}
    \ = \ 1 \after \asrt_{\overline{\andthen{a}{b}}}
    \ = \ 1 \after \nu \after \asrt_{\overline{\andthen{b}{a}}} \after \mu
    \ = \ \overline{\andthen{b}{a}} \after \mu.
\end{equation}
Using this equation again in the previous, we find
\begin{equation*}
    \nu \after \mu \after \asrt_{\overline{\andthen{a}{b}}}
            \ = \ 
    \nu \after \mu \after \asrt_{\overline{\andthen{b}{a}}\after \mu}
            \ = \ 
    \nu \after \asrt_{\overline{\andthen{b}{a}}} \after \mu
            \ = \ \asrt_{\overline{\andthen{a}{b}}}.
\end{equation*}
Thus~$\nu \after \mu = \id$ and so~$\mu = \nu^{-1}$.
Write~$l = \pi_{\ceil{\andthen{b}{a}}} \after \nu^{-1}
                        \after \zeta_{\ceil{\andthen{a}{b}}}$. Then
\begin{alignat*}{2}
    (\andthen{a}{b}) \after l^\dagger & \ = \ 
    (\andthen{a}{b})
        \after \pi_{\ceil{\andthen{a}{b}}} 
        \after \nu
        \after \zeta_{\ceil{\andthen{b}{a}}}  
        &\qquad& \text{by \sref{asrt-pristine-reverse}}\\
        & \ = \ 
        \overline{\andthen{a}{b}}
        \after \nu
        \after \zeta_{\ceil{\andthen{b}{a}}}  \\
        & \ = \ 
        \overline{\andthen{b}{a}}
        \after \zeta_{\ceil{\andthen{b}{a}}}  
        &&\text{by \eqref{dagger-seqprod-inversion} and~$\mu = \nu^{-1}$}\\
        & \ = \ 
        \andthen{b}{a}.
\end{alignat*}
And so, keeping in mind~$\andthen{a}{b} \leq \ceil{\andthen{a}{b}}
    =  1 \after l$,
we have
\begin{alignat*}{2}
    \asrt_b \after \asrt_a &\ = \ 
        l \after \asrt_{\andthen{a}{b}} \\
        & \ = \ 
        \asrt_{(\andthen{a}{b}) \after l^\dagger} \after l 
    &\qquad&\text{by \sref{asrt-pristine-reverse}} \\
        &\ = \ 
        \asrt_{\andthen{b}{a}} \after l \\
        &\ = \ 
        \asrt_{\andthen{b}{a}}
        \after \pi_{\ceil{\andthen{b}{a}}} 
        \after \nu^{-1} \after \zeta_{\ceil{\andthen{a}{b}}},
\end{alignat*}
as promised. \qed
\end{point}
\end{point}
\begin{point}{130}[dagger-iso-omega2]{Corollary}
In setting \sref{dagger-setting},
    the daggered version of
        \eqref{dagger-iso-omega} holds --- that is:
    \begin{equation*}
\asrt_q \after
\asrt_{p \after k}
    \ =\ 
    \asrt_{p \after g} \after
    \pi_{\ceil{p \after g}} \after
    \omega^{-1} \after
    \zeta_{\lceil p \after k \rceil}.
    \end{equation*}
\end{point}
\spacingfix{}
\begin{point}{140}[dagger-iso-chi2]{Lemma}%
In setting \sref{dagger-setting},
    the daggered version of
        \eqref{dagger-iso-chi} holds --- that is:
    \begin{equation*}
    \zeta_t \after
    \asrt_p \after
    \pi_{\ceil{p}}
    \ = \ 
    \asrt_{p \after \pi_t} \after
    \pi_{\ceil{p \after \pi_t}} \after
    \chi^{-1} \after
    \zeta_{\lceil t \after \asrt_p \after \pi_{\ceil{p}}\rceil}.
    \end{equation*}
\spacingfix{}
\begin{point}{150}{Proof}%
The heavy lifting has been done in \sref{dagger-iso-mu} already. To start, note
    \begin{align*}
\asrt_p \after \asrt_t
& \ \overset{\mathclap{\smash{\sref{asrt-absorp-rule}}}}{=}
    \ \asrt_{\ceil{p}} \after \asrt_p  \after \asrt_t \\
    & \ = \ \pi_{\ceil{p}} \after \zeta_{\ceil{p}} \after \asrt_p \after
                    \pi_t \after \zeta_t \\
                    & \ \overset{\mathclap{\smash{\eqref{dagger-iso-chi}}}}{=} \ 
        \pi_{\ceil{p}} \after 
        \pi_{\lceil t \after
        \asrt_p \after \pi_{\ceil{p}}\rceil} \after
        \chi \after \zeta_{\ceil{p \after \pi_t}} \after
        \asrt_{p \after \pi_t} \after
        \zeta_t \\
                    & \ \overset{\mathclap{\smash{\sref{pristine-asrt}}}}{=} \ 
        \pi_{\ceil{p}} \after 
        \pi_{\lceil t \after
        \asrt_p \after \pi_{\ceil{p}}\rceil} \after
        \chi \after \zeta_{\ceil{p \after \pi_t}} \after
        \zeta_t \after
        \asrt_{p \after \pi_t\after \zeta_t} \\
                    & \  = \ 
        \pi_{\ceil{p}} \after 
        \pi_{\lceil t \after
        \asrt_p \after \pi_{\ceil{p}}\rceil} \after
        \chi \after
        \zeta_{\ceil{p \after \pi_t}} \after
        \zeta_t \after
        \asrt_{\andthen{t}{p}} \\
                    & \  = \ 
        \pi_{\ceil{\andthen{p}{t}}} \after \alpha_2 \after
        \chi \after \beta_2 \after
        \zeta_{\ceil{\andthen{t}{p}}} \after
        \asrt_{\andthen{t}{p}},
    \end{align*}
where~$\alpha_2$ and~$\beta_2$ are the unique isomorphisms such that
\begin{align*}
        \pi_{\ceil{p}} \after 
        \pi_{\lceil t \after
        \asrt_p \after \pi_{\ceil{p}}\rceil} &\ = \ 
            \pi_{\ceil{\andthen{p}{t}}} \after \alpha_2 &
        \zeta_{\ceil{p \after \pi_t}} \after
        \zeta_t
        &\ = \ 
        \beta_2 \after \zeta_{\ceil{\andthen{t}{p}}}.
\end{align*}
With the same reasoning as in \sref{dagger-iso-alpha2}
        and~\sref{dagger-iso-beta2}, we see
\begin{align*}
        \zeta_{\lceil t \after
        \asrt_p \after \pi_{\ceil{p}}\rceil} \after
        \zeta_{\ceil{p}}
        &\ = \ 
        \alpha_2^{-1} \after
            \zeta_{\ceil{\andthen{p}{t}}}
            &
        \pi_t \after
        \pi_{\ceil{p \after \pi_t}}
        &\ = \ 
        \pi_{\ceil{\andthen{t}{p}}}
        \after \beta_2^{-1}.
\end{align*}
Now we can apply \sref{dagger-iso-mu}:
\begin{align*}
    \zeta_t \after \asrt_p \after \pi_{\ceil{p}} 
    &\ \overset{\smash{\mathclap{\sref{asrt-absorp-rule}}}}{=} \ 
    \zeta_t \after \asrt_t \after
    \asrt_p \after
    \pi_{\ceil{p}}
    \\
    &\ \overset{\smash{\mathclap{\sref{dagger-iso-mu}}}}{=} \ 
    \zeta_t \after \asrt_{\andthen{t}{p}} \after
    \pi_{\ceil{\andthen{t}{p}}} \after
    \beta_2^{-1} \after
    \chi^{-1} \after
    \alpha_2^{-1} \after
    \zeta_{\ceil{\andthen{p}{t}}} \after
    \pi_{\ceil{p}}
    \\
    &\ = \ 
    \zeta_t \after \asrt_{\andthen{t}{p}} \after
        \pi_t \after
        \pi_{\ceil{p \after \pi_t}} \after
    \chi^{-1} \after
        \zeta_{\lceil t \after
        \asrt_p \after \pi_{\ceil{p}}\rceil} \after
        \zeta_{\ceil{p}} \after
    \pi_{\ceil{p}}
    \\
    &\ = \ 
    \zeta_t \after \asrt_{p \after \pi_t \after \zeta_t} \after
        \pi_t \after
        \pi_{\ceil{p \after \pi_t}} \after
    \chi^{-1} \after
        \zeta_{\lceil t \after
        \asrt_p \after \pi_{\ceil{p}}\rceil}
    \\
    &\ = \ 
    \asrt_{p \after \pi_t} \after 
        \pi_{\ceil{p \after \pi_t}} \after
    \chi^{-1} \after
        \zeta_{\lceil t \after
        \asrt_p \after \pi_{\ceil{p}}\rceil},
\end{align*}
as desired. \qed
\end{point}
\end{point}
\begin{point}{160}[dagger-is-functor]{Proposition}%
In a~$\dagger'$-effectus~$(f \after g)^\dagger = g^\dagger \after f^\dagger$
    holds.
\begin{point}{170}{Proof}%
We work in setting~\sref{dagger-setting}.
The equality~$(f\after g)^\dagger = g^\dagger \after f^\dagger$
follows from~\sref{dagger-of-fg} and the commutativity
of the following diagram
\begin{equation*}
    \xymatrix@C+3pc {
        \bullet \ar@{<-}[rr]^{g^\dagger}
        \ar@{<-}[rd]|{\asrt_q}
        \ar@{<-}[ddd]_{\rotatebox{90}{$\scriptstyle
            \asrt_{p \after g}
            \pi_{\ceil{p \after g}} \after
            \omega^{-1}
        $}}
        && \bullet \ar@{<-}[rr]^{f^\dagger}
            \ar@{<-}[rd]|{ \asrt_p \after \zeta_{\ceil{p}} }
        && \bullet
            \\ \ar@{}[rd]|{\text{\sref{dagger-iso-omega2}}}
            & \bullet
            \ar@{<-}[r]^{ \pi_{\ceil{q}} \after \psi^{-1} }
            \ar@{<-}[d]|{\asrt_{p \after k}}
                \ar@{}[rd]|{\text{\sref{dagger-iso-zeta2}}}
 & \bullet
                \ar@{<-}[u]^{\zeta_t}
                \ar@{<-}[d]|{\asrt_{p \after \pi_t}}
                & \bullet \ar@{<-}[ru]_{\varphi^{-1} \after \zeta_s}
            \ar@{}[rd]|{\text{\sref{dagger-iso-alpha2}}}
            \\& \bullet \ar@{<-}[r]^{\pi_{\ceil{q}} \after \psi^{-1}}
            \ar@{<-}[d]|{\pi_{\ceil{p \after k}}}
                \ar@{}[rd]|{\text{\sref{dagger-iso-beta2}}}
                &\bullet \ar@{<-}[d]|{\pi_{\ceil{p \after \pi_t}}}
                        \ar@{}[ru]|{\text{\sref{dagger-iso-chi2}}}
                        &&
            \\ \bullet \ar@{=}[r]
            \ar@{<-}[ru]^{\zeta_{\ceil{p \after k}}}
            &\bullet \ar@{<-}[r]_{\beta^{-1}}
            &\bullet \ar@{<-}[rr]_{\chi^{-1}}
            &&\bullet \ar@{<-}[uuu]_{\rotatebox{90}{$\scriptstyle
    \alpha^{-1}\after \zeta_{\ceil{t \after f^\dagger}}$}}
    \ar@{<-}[luu]|{ \zeta_{\lceil t \after \asrt_p \after \pi_{\ceil{p}} \rceil}} }
\end{equation*}
or alternatively by
\begin{align*}
    g^\dagger \after f^\dagger
        & \ = \ 
            \asrt_q \after
            \pi_{\ceil{q}} \after
            \psi^{-1} \after
            \zeta_t \after
            \asrt_p \after
            \pi_{\ceil{p}} \after
            \varphi^{-1} \after
            \zeta_s
        \\
        & \ \overset{\smash{\mathclap{\sref{dagger-iso-chi2}}}}{=} \ 
            \asrt_q \after
            \pi_{\ceil{q}} \after
            \psi^{-1} \after
    \asrt_{p \after \pi_t} \after
    \pi_{\ceil{p \after \pi_t}} \after
    \chi^{-1} \after
    \zeta_{\lceil t \after \asrt_p \after \pi_{\ceil{p}}\rceil} \after
            \varphi^{-1} \after
            \zeta_s
        \\
        & \ \overset{\smash{\mathclap{\sref{dagger-iso-zeta2}}}}{=} \ 
            \asrt_q \after
    \asrt_{p \after k} \after
    \pi_{\ceil{q}} \after
    \psi^{-1} \after
    \pi_{\ceil{p \after \pi_t}} \after
    \chi^{-1} \after
    \zeta_{\lceil t \after \asrt_p \after \pi_{\ceil{p}}\rceil} \after
            \varphi^{-1} \after
            \zeta_s
        \\
        & \ \overset{\smash{\mathclap{\sref{dagger-iso-beta2}}}}{=} \ 
            \asrt_q \after
    \asrt_{p \after k} \after
        \pi_{\ceil{p \after k}} \after 
        \beta^{-1} \after
    \chi^{-1} \after
    \zeta_{\lceil t \after \asrt_p \after \pi_{\ceil{p}}\rceil} \after
            \varphi^{-1} \after
            \zeta_s
        \\
        & \ \overset{\smash{\mathclap{\sref{dagger-iso-alpha2}}}}{=} \ 
            \asrt_q \after
    \asrt_{p \after k} \after
        \pi_{\ceil{p \after k}} \after 
        \beta^{-1} \after
    \chi^{-1} \after
        \alpha^{-1} \after
        \zeta_{\ceil{t \after f^\dagger}}
        \\
        & \ \overset{\smash{\mathclap{\sref{dagger-iso-omega2}}}}{=} \ 
    \asrt_{p \after g} \after
    \pi_{\ceil{p \after g}} \after
    \omega^{-1} \after
    \zeta_{\lceil p \after k \rceil} \after
        \pi_{\ceil{p \after k}} \after 
        \beta^{-1} \after
    \chi^{-1} \after
        \alpha^{-1} \after
        \zeta_{\ceil{t \after f^\dagger}}
        \\
        & \ = \ 
    \asrt_{p \after g} \after
    \pi_{\ceil{p \after g}} \after
    \omega^{-1} \after
        \beta^{-1} \after
    \chi^{-1} \after
        \alpha^{-1} \after
        \zeta_{\ceil{t \after f^\dagger}}
        \\
        & \ \overset{\smash{\mathclap{\sref{dagger-of-fg}}}}{=} \ 
        (f \after g)^\dagger,
\end{align*}
whichever the Reader might prefer. \qed
\end{point}
\end{point}
\end{parsec}

\begin{parsec}{2200}%
\begin{point}{10}%
We are ready to finish the proof of \sref{dagger-theorem}.
\end{point}
\begin{point}{20}[dagger-thm-sufficiency]{Theorem}%
    A~$\dagger'$-effectus is a~$\dagger$-effectus
        with~$\dagger$
            as defined in \sref{dagger-definition2}.
\begin{point}{30}{Proof}%
Let~$C$ be a~$\dagger'$-effectus.
\begin{point}{40}{Ax.~1}%
By~\sref{dagger-is-functor}, \sref{dagger-idempotent}
    and~\sref{dagger-prime-basics}
    the~$\dagger$ defined in~\sref{dagger-definition2}
    turns~$\Pure C$ into a~$\dagger$-category.
Also by~\sref{dagger-prime-basics},
    we have~$\asrt_p^\dagger = \asrt_p$
    for any predicate~$p$.
Pick any pure~$f$.
We have to show~$f$ is~$\diamond$-adjoint to~$f^\dagger$.
By~\sref{standard-form-map}
    we have~$f =
    \pi_{\IM f} \after \varphi \after \zeta_{\ceil{1 \after f}}
        \after \asrt_{1 \after f}$
        for some iso~$\alpha$.
We compute
\begin{align*}
   f_\diamond 
   & \ = \ 
   (\pi_{\IM f})_\diamond \after \varphi_\diamond \after (\zeta_{\ceil{1 \after f}})_\diamond
   \after (\asrt_{1 \after f})^\diamond \\
   & \ \overset{\mathclap{\smash{\sref{quotcompr-diamond-adjoint}}}}{=} \ 
   (\zeta_{\IM f})^\diamond \after \varphi_\diamond \after (\pi_{\ceil{1 \after f}})^\diamond
   \after (\asrt_{1 \after f})^\diamond \\
   & \ \overset{\mathclap{\smash{\sref{iso-diamond-adjoint}}}}{=} \ 
   (\zeta_{\IM f})^\diamond \after (\varphi^{-1})^\diamond \after (\pi_{\ceil{1 \after f}})^\diamond
   \after (\asrt_{1 \after f})^\diamond \\
   & \ = \ 
   (\asrt_{1 \after f} \after
   \pi_{\ceil{1 \after f}} \after
   {\varphi^{-1}} \after
   \zeta_{\IM f})^\diamond
   \\
   & \ = \ 
   (f^\dagger)^\diamond,
\end{align*}
so~$f$ is indeed~$\diamond$-adjoint to~$f^\dagger$.
\end{point}
\begin{point}{50}{Ax.~2}%
Let~$f$ be a~$\dagger$-positive map.
That is: $f = h^\dagger \after h$ for some pure map~$h$.
By~\sref{standard-form-map}
    we have~$h =
    \pi_{\IM h} \after \alpha \after \zeta_{\ceil{1 \after h}}
                    \after \asrt_{1 \after h}
                    $
                    for some iso~$\alpha$.
            We compute
\begin{align*}
    f &\ = \ h^\dagger \after h \\
    &\ = \ 
    \asrt_{1 \after h} \after \pi_{\ceil{1 \after  h}} \after \alpha^{-1} \after \zeta_{\IM h} \after
    \pi_{\IM h} \after \alpha \after \zeta_{\ceil{1 \after h}}
                    \after \asrt_{1 \after h}
                    \\
    &\ = \ 
    \asrt_{1 \after h} \after \asrt_{\ceil{1 \after  h}}
                    \after \asrt_{1 \after h}
                    \\
                    &\ \overset{\smash{\mathclap{\sref{asrt-absorp-rule}}}}{=} \ 
    \asrt_{1 \after h} \after \asrt_{1 \after h}.
\end{align*}
Form this it follows that $\dagger$-positive maps are~$\diamond$-positive.
Let~$q$ be the predicate such that~$\andthen{q}{q} = 1 \after h$.
Then
\begin{equation*}
\asrt_{q}^\dagger \after \asrt_{q}  
    \ =\  \asrt_q\after \asrt_q
    \ \overset{\sref{andthen-square-rule}}{=}\  
    \asrt_{\andthen{q}{q}}
    \  =\  \asrt_{1 \after h}.
\end{equation*}
Thus~$\asrt_{1 \after h}$ is a~$\dagger$-positive
    with~$\asrt_{1 \after h} \after \asrt_{1 \after h} = f$.
We have to show~$\asrt_{1 \after h}$
    is the unique map with this property.
Let~$g$ be any~$\dagger$-positive map with~$g \after g = f$.
Recall both~$g$ and~$f$ are~$\diamond$-positive, hence
\begin{equation*}
    \asrt_{1 \after f} \ =\  f\ =\ g \after g \ = \ \asrt_{1 \after g}
            \after \asrt_{1 \after g}
            \ \overset{\smash{\sref{andthen-square-rule}}}{=}\  
            \asrt_{\andthen{(1 \after g)}{(1 \after g)}}.
\end{equation*}
So~$\andthen{(1 \after g)}{(1\after g)} = 1\after f =
            \andthen{(1 \after h)}{(1 \after h)}$.
Hence~$1\after g = 1 \after h$ by uniqueness of the square root.
So~$g = \asrt_{1 \after g} = \asrt_{1 \after h} = h$, as desired.
\end{point}
\begin{point}{60}{Ax.~3}%
Let~$\asrt_p$ be any~$\diamond$-positive map.
Write~$q$ for the unique predicate with~$\andthen{q}{q}=p$.
Then
\begin{equation*}
    \asrt_p
     \ = \ \asrt_{\andthen{q}{q}}
            \ \overset{\smash{\sref{andthen-square-rule}}}{=}\ 
     \asrt_q \after \asrt_q  \ =\ 
     \asrt_q^\dagger \after \asrt_q ,
\end{equation*}
so~$\asrt_p$ is~$\dagger$-positive, as desired. \qed
\end{point}
\end{point}
\end{point}
\end{parsec}

\subsection{Dilations}
\begin{parsec}{2210}%
\begin{point}{10}%
Being the main topic of the first part of this thesis,
    we cannot finish our discussion of~$\dagger$-effectuses
    without mentioning dilations.
The presented abstract theory of dilations in a~$\dagger$-effectus
    is preliminary and as of yet rather unimpressive.
\end{point}
\begin{point}{20}[dfn-eff-dilations]{Definition}%
Let~$f\colon X \to Y$ be any map in an effectus.
    A \Define{dilation} of~$f$ \index{dilation!in an effectus}
    is a triple~$(P, \varrho, h)$
    of a sharp total map~$\varrho \colon P \to Y$
    and a pure map~$h \colon X \to P$
    such that~$\varrho \after h = f$
        and the following universal property holds.
\begin{quote}
    For every  triple~$(P', \varrho', h')$
    with~$\varrho'\colon P' \to Y$ total sharp,
    $h'\colon X \to P$ arbitrary and~$f = \varrho'\after h'$,
    there is a unique~$\sigma \colon P \to P'$
    with~$\sigma \after h = h'$
    and~$\varrho' \after \sigma = \varrho$.
\end{quote}
    We say an effectus \Define{has dilations} if \index{effectus!with dilations}
    there exists a dilation for every map.
\end{point}
\begin{point}{30}{Example}%
    The effectus $\op\vN$~(\sref{effectus-vn}) has dilations,
    as shown in \sref{existence-paschke},
    but the full subcategory~$\op\CvN$
    of commutative von Neumann algebras
    does not as shown in the next exercise.
\end{point}
\begin{point}{31}[exc-cvn-no-dilations]{Exercise}%
In this exercise you will show that~$\op\CvN$
    does not have dilations.
First, show that a corner between commutative von Neumann algebras
    is multiplicative.
Conclude that if~$\op\CvN$ were to have dilations,
    then any ncpu-map would be multiplicative, quod non.
\end{point}
\begin{point}{40}[dils-abstract-basics]{Proposition}%
Let~$C$ be an effectus with dilations.
\begin{enumerate}
\item
    If~$(P_1,\varrho_1,h_1)$
    and~$(P_2,\varrho_2,h_2)$
        are both dilations of~$f$,
    then there is a unique isomorphism~$\alpha\colon P_1 \to P_2$
        with~$\alpha \after h_1 = h_2$ and~$\varrho_2 \after \alpha = \varrho_1$.
\item
    Conversely,
    if~$(P,\varrho,h)$
        is a dilation of~$f$ and~$\alpha\colon P \to P'$ is some isomorphism,
        then~$(P', \varrho \after \alpha^{-1}, \alpha \after h)$
        is also a dilation of~$f$.

\item
    $(X, \varrho, \id)$ is the dilation
        of a sharp total map~$\varrho$.

\item
    If~$(P, \varrho, h)$ is some dilation,
        then~$(P, \id, h)$ is the dilation of~$h$.

\item
    For any quotient~$\xi \colon X \to Q$
            and map~$f\colon Q \to Y$
            with dilation~$(P, \varrho, h)$,
            the triple~$(P, \varrho, h \after \xi)$
            is a dilation of~$f \after \xi$.

\item
    Conversely, if~$(P, \varrho, h)$
            is a dilation of~$f \after \xi$
            for a map~$f\colon Q \to Y$
            and quotient~$\xi \colon X \to Q$,
        then~$(P, \varrho, h'')$ is a dilation of~$f$,
        where~$h''\colon Q \to P$ is the unique map with~$h'' \after \xi = h$.

\item
    If for~$i=1,2$ the triple~$(P_i, \varrho_i, h_i)$
            is a dilation of~$f_i \colon X_i \to Y$,
    then~$(P_1 + P_2, [\varrho_1,\varrho_2], h_1 + h_2)$
        is a dilation of~$[f_1, f_2] \colon X_1 + X_2 \to Y$.
\end{enumerate}
\spacingfix{}
\begin{point}{50}{Proof}%
Point 1 has been proven in \sref{paschke-unique-up-to-iso}.
Point 2 is easily verified.
To prove point 3, assume~$\varrho$ is a sharp total map.
Assume~$\varrho = \varrho' \after h'$ for some~$h'$ and sharp total~$\varrho'$.
Then~$h'$ is clearly the unique map
    satisfying~$h' = h' \after \id$ (and~$\varrho = \varrho' \after h'$),
    which demonstrates point 3.

For point 4, assume~$(P, \varrho, h)$ is some dilation.
Assume~$h = \varrho' \after h'$ for some~$h'$ and sharp total~$\varrho'$.
There is a unique~$\sigma$ with~$h' = \sigma \after h$
    and~$\varrho\after \varrho' \after \sigma = \varrho$.
Now~$\varrho' \after \sigma \after h = \varrho' \after h' = h$.
    So~$\varrho' \after \sigma$
        satisfies the properties of the unique mediating map~$\id\colon P \to P$,
        hence~$\varrho' \after \sigma = \id$.
Only uniqueness of~$\sigma$ remains,
    so assume~$\sigma'$ is a ncp-map with~$h'  = \sigma' \after h$
    and~$\varrho' \after \sigma' = \id$.
Then~$\varrho \after \varrho' \after \sigma' = \varrho$
    and so~$\sigma = \sigma'$.

To demonstrate point 5, assume~$\varrho' \after h' = f \after \xi$
    for some~$h'$ and total sharp~$\varrho'$.
    Then~$1 \after h' = 1 \after \varrho' \after h'
                    = 1 \after f \after \xi \leq 1 \after \xi$,
                    so there is a unique~$h''$
                    with~$h'' \after \xi = h'$.
By epicity of~$\xi$,
    we get~$\varrho'\after h'' = f$
    from~$\varrho' \after h'' \after \xi = \varrho' \after h' = f \after \xi$.
Thus there is a unique~$\sigma$
    with~$\sigma \after h = h''$ and~$\varrho' \after \sigma = \varrho$
    and so~$\sigma \after h \after \xi = h'' \after \xi = h'$.
Only uniqueness of~$\sigma$ remains,
    so suppose~$\sigma'$ is a map with~$\sigma'\after h \after \xi = h'$
            and~$\varrho' \after \sigma = \varrho$.
    Then~$\sigma' \after h = h''$ by epicity of~$\xi$,
    hence~$\sigma'=\sigma$ as desired.

For point 6, assume~$f  = \varrho' \after h'$
    for some~$h'$ and sharp total~$\varrho'$.
    Then~$f \after \xi = \varrho' \after h' \after \xi $
    and so there is a unique~$\sigma$
    with~$\varrho \after \sigma = \varrho'$
    and~$\sigma \after h' \after \xi = h \equiv h'' \after \xi$.
    Thus~$\sigma \after h' = h''$ by epicity of~$\xi$.
To show uniqueness,
    assume~$\sigma'$ is such that~$\sigma' \after h' = h''$
    and~$\varrho \after \sigma' = \varrho'$.
    Then~$\sigma' \after h' \after \xi = h'' \after xi = h$
    and indeed~$\sigma' = \sigma$.

To show point 7, assume~$[f_1,f_2] = \varrho'\after h'$
    for some~$h'$ and total sharp~$\varrho'$.
For any~$i=1,2$
    we have~$\varrho' \after h' \after \kappa_i = f_i$,
    so there is a unique~$\sigma_i$
    with~$\sigma_i \after h_i = h' \after \kappa_i$
    and~$\varrho_i = \varrho' \after \sigma_i$.
Define~$\sigma = [\sigma_1, \sigma_2]$.
Clearly~$\sigma \after (h_1 + h_2)
    = [\sigma_1 \after h_1,\sigma_2 \after h_2]
    = [ h' \after \kappa_1,  h'\after \kappa_2] = h'$
    and~$\varrho' \after \sigma
            = [\varrho' \after \sigma_1, \varrho' \after \sigma_2]
            = [\varrho_1, \varrho_2]$.
Assume~$\sigma'$ is any ncp-map with~$\sigma' \after (h_1 + h_2) = h'$
    and~$\varrho' \after \sigma' = [\varrho_1, \varrho_2]$.
For~$i=1,2$
    we have~$\varrho' \after \sigma' \after \kappa_i = \varrho_i$
    and~$h' \after \kappa_i = \sigma' \after (h_1 + h_2) \after \kappa_i
                =  \sigma' \after h_i$,
                so~$\sigma' \after \kappa_i = \sigma_i$.
    Hence~$\sigma' = [\sigma' \after \kappa_1, 
                \sigma' \after \kappa_2] = [\sigma_1,\sigma_2] = \sigma$.
                \qed
\end{point}
\end{point}
\end{parsec}

\begin{parsec}{2220}%
\begin{point}{10}%
With some additional assumptions,
    the existence of dilations
    forces the existence of some
    familiar quantum gates.
Assume~$C$ is a~$\dagger$-effectus with dilation
    with a scalar~$\lambda$ such that~$\lambda^\perp = \lambda$.
    (For instance: $\lambda = \frac{1}{2}$ in~$\op\vN$.)
Let~$(P, \langle s,s^\perp \rangle, h)$ be a dilation
    of the map~$\left<h,h^\perp\right>\colon 1 \to 1+1$.
We think of~$P$ as an abstract qubit
    --- indeed in $\op\vN$ any
    dilation of~$\left<\frac{1}{2},\frac{1}{2}\right>$
    is of the form~$(M_2, \langle e,e^\perp \rangle, h)$
    where~$h$ is the vector state for some~$v \in \C^2$
    and~$e$ is a projector on~$w\C$,
    with~$\left<v,w\right> = \frac{1}{\sqrt{2}}$
    (for instance~$v = \ket{0}$ and~$w = \ket{+}$).
\begin{point}{20}%
By assumption~$\langle s,s^\perp \rangle \after h = \langle\lambda, \lambda\rangle$,
    so~$s \after h = s^\perp \after h = \lambda$.
and~$\langle s^\perp,s \rangle \after h = \langle\lambda, \lambda\rangle$.
Hence by the definition of dilation,
    there is a unique map~$X \colon P \to P$
    with~$X \after h = h$
    and~$\langle s^\perp, s\rangle\after X = \langle s, s^\perp\rangle $,
    whence~$s \after X = s^\perp \after X$
    and~$s^\perp \after X = s \after X$.
Clearly~$X \after X \after h = h$
    and~$\langle s, s^\perp \rangle  \after X \after X 
            = \langle s, s^\perp \rangle$,
    so by uniqueness of the mediating map~$\id \colon P \to P$,
    we must have~$X \after X = \id$.
Thus the isomorphism~$X$ behaves like the~$X$-gate
    --- indeed, in~$\op\vN$, we have~$X = \ad_X$
        if~$e = \ketbra{0}{0}$.
\end{point}
\begin{point}{30}%
For the next gate,
    assume~$s$ is pure and~$\IM s = 1$.
Note~$s^\dagger$ is total.
As~$\asrt_\lambda = \lambda$,
    we have~$h^\dagger \after s^\dagger = \lambda^\dagger = \lambda$
    and so~$(h^\dagger)^\perp \after s^\dagger = \lambda^\perp = \lambda$.
Combined:~$\langle h^\dagger, (h^\dagger)^\perp \rangle \after s^\dagger
            = \langle \lambda,\lambda \rangle$.
Thus using the universal property of the dilation,
    there is a unique map~$H\colon P \to P$
    with~$H \after h = s^\dagger$
    and~$\langle h^\dagger, (h^\dagger)^\perp \rangle \after H =
            \langle s, s^\perp \rangle$.
So~$h^\dagger \after H = s$.
In general, we do not know a lot more about~$H$.
Assume~$H$ is an comprehension (as it is in~$\op\vN$ in this situation).
    Then~$s^\dagger = H^\dagger \after h$ and
        $s = h^\dagger \after H^\dagger$,
        so $s^\perp = (h^\dagger)^\perp \after H^\dagger$,
        whence~$\langle s,s^\perp \rangle = \langle
                h^\dagger, (h^\dagger)^\perp \rangle \after H^\dagger$.
Thus, by uniqueness of~$H$, we must have~$H^\dagger = H$.
Thus both~$H^\dagger$ and $H$ are total,
    which forces~$H\after H^\dagger = H^\dagger \after H = \asrt_{1} = \id$,
    by uniqueness of~$\diamond$-positive maps.
Hence~$H \after H = \id$,
    $s \after H = h^\dagger$
    and~$h^\dagger \after H = s$.
    So~$H$ is similar to the Hadamard-gate
        --- indeed, in~$\op\vN$, we have~$H = \ad_H$ if~$e = \ketbra{0}{0}$.
\end{point}
\begin{point}{40}%
Do our additional assumptions (purity of~$s$, $\IM s = 1$, etc.)
follow from more general principles?  Can we define the other gates, like
    the CNOT? Do these gates interact in the usual way?
    Unfortunately, we do not have an answer to any of these questions.
\end{point}
\end{point}
\end{parsec}

\begin{parsec}{2230}%
\begin{point}{10}%
In~\sref{paschke-correspondence}
    we saw that for any ncp-map~$\varphi$
    with Paschke dilation~$(\scrP, \varrho, h)$,
       there  is a linear order
       isomorphism~$[0,1]_{\varrho(\scrA)^\square} \cong [0,\varphi]_\ncp$.
To finish the discussion of dilations,
    we mention how to state this correspondence
        in an~$\&$-effectus.
\end{point}
\begin{point}{20}[sefp]{Definition}%
In an~$\&$-effectus,
for any predicate~$p$,
    define~$\Define{\sef_p} \equiv \asrt_p \ovee \asrt_{p^\perp}$.
    \index{sefp@$\sef_p$}
This map models the side-effect of measuring~$p$.
    For any map~$f\colon X \to Y$, define
    \begin{equation*}
        \Define{\Inv f} \ \equiv \ \{p; \ p \in \Pred X;\ f \after \sef_p = f\}.
        \index{Inv@$\Inv f$}
    \end{equation*}
This is the set of predicates whose measurement does not disturb~$f$.
\end{point}
\begin{point}{30}{Lemma}%
    In~$\op\vN$, for any nmiu-map~$\varrho\colon \scrA \to \scrB$,
        we have~$\Inv \varrho = [0,1]_{\varrho(\scrA)^\square}$ ---
        that is:~$\Inv \varrho$
is the set of effects of the commutator of the image of~$\varrho$.
\begin{point}{40}{Proof}%
Assume~$a \in [0,1]_{f(\scrA)^\square}$.
Simply unfolding definitions, we find~$\sef_a(\varrho(b)) = \sqrt{a} \varrho(b) \sqrt{a}
                        + \sqrt{1-a} \varrho(b) \sqrt{1-a}
                        = \varrho(b) $ for any~$b \in \scrA$, so~$a \in \Inv \varrho$.
For the converse, assume~$a \in \Inv \varrho$.
    Recall~$(\scrA, \varrho, \id)$ is a Paschke dilation of~$\varrho$.
    As
\begin{equation*}
    \varrho  \ =\ \sef_a \after \varrho
            \ = \  (\asrt_a \after \varrho)\ovee (\asrt_{1-a} \after \varrho )
            \ \geq_{\ncp} \ \asrt_{a} \after \varrho,
\end{equation*}
    we know by \sref{paschke-correspondence}
    that there is a~$t \in [0,1]_{\varrho(\scrA)^\square}$
    with~$\asrt_a(\varrho(b)) \equiv \sqrt{a} \varrho(b) \sqrt{a}
                =t  \varrho(b)$ for all~$b \in \scrA$,
        but then clearly~$a = \sqrt{a} \varrho(1) \sqrt{a} = t \varrho(1) = t$. \qed
\end{point}
\end{point}
\begin{point}{50}{Definition}%
Let~$C$ be an $\&$-effectus.
For a map~$f\colon X \to Y$,
    write
\begin{equation*}
    \Define{\downarrow\!f} \ \equiv \ 
    \{ g; \ g\colon X \to Y; \ g \leq f \}.
\end{equation*}
    We say a dilation~$(P, \varrho, h)$ of~$f$
        has \Define{the order correspondence} if
        \index{effectus!with dilations!and order correspondence}
        there is an order isomorphism~$\Theta\colon \downarrow\!f \to \Inv \varrho$
        such that for every~$g \in \downarrow\!f$,
        we have
\begin{equation*}
g \ =\  \varrho \after \asrt_{\Theta(g)} \after h.
\end{equation*}
\end{point}
\spacingfix{}
\begin{point}{60}{Example}%
By \sref{paschke-correspondence}
    every dilation in~$\op\vN$
    has the order correspondence.
\end{point}
\end{parsec}

\section{Comparisons}
\subsection{Dagger kernel categories}
\begin{parsec}{2240}%
\begin{point}{10}%
To finish this chapter, we relate the structures in a~$\dagger$-effectus
    with existing structures in the literature.
Just like we try to axiomatize~$\op\vN$,
    Heunen in the third chapter of his Ph.D-thesis attempted
    to axiomatize~$\mathsf{Hilb}$,
    the category of Hilbert spaces with bounded operators.
Heunen came very close:
    he proves that a~$\dagger$-category obeying
    certain additional axioms must embed into~$\mathsf{Hilb}$.
    \cite[3.7.18]{heunenphd}
We will discover that some of his axioms for~$\mathsf{Hilb}$
    hold in~$\Pure(C)$ for a~$\dagger$-effectus~$C$
    and others do not in general.
(We do not cover all his axioms.)
\end{point}
\begin{point}{20}{Definition}%
Let~$C$ be a~$\dagger$-category.
We say that an arrow~$f$
    is \Define{$\dagger$-mono}
    \index{*dagg@$\dagger$, dagger!-mono}
    iff~$f^\dagger \after f = \id$
    and dually~$f$ is \Define{$\dagger$-epi}
    \index{*dagg@$\dagger$, dagger!-epi}
    iff~$f\after f^\dagger = \id$.
An arrow~$f$ is a \Define{$\dagger$-partial isometry}
    \index{*dagg@$\dagger$, dagger!-partial isometry}
    if~$f = m \after e$
    for some~$\dagger$-mono~$m$
    and~$\dagger$-epi~$e$.

A~\Define{$\dagger$-kernel} of~$f$
    \index{kernel!$\dagger-$}
    \index{*dagg@$\dagger$, dagger!-kernel}
is an  equalizer of~$f$ with~$0$, which is~$\dagger$-mono.
A~\Define{$\dagger$-kernel category} is a~$\dagger$-category
    \index{*dagg@$\dagger$, dagger!-kernel!category}
    with zero object and where each arrow has a~$\dagger$-kernel
\cite[3.2.20]{heunenphd}.
\end{point}
\begin{point}{30}{Proposition}%
In~$\Pure(C)$ for some~$\dagger$-effectus~$C$,
    the following holds.
\begin{enumerate}
\item
    $f$ is $\dagger$-mono iff $f$ is a comprehension. \\
    Dually:
    $f$ is $\dagger$-epi iff $f$ is a quotient of a sharp predicate.
\item
    The $\dagger$-partial isometries are exactly the pristine maps, see \sref{dfn-pristine}.
\end{enumerate}
Furthermore: $\Pure C$ is a $\dagger$-kernel category:
    the~$\dagger$-kernel of~$f$ is given by~$\pi_{(1 \after f)^\perp}$.
\begin{point}{40}{Proof}%
Let~$f$ be any~$\dagger$-mono
pure map, say~$f \equiv \pi_s \after \alpha \after \zeta_{\ceil{p}}
        \after \asrt_p$
        for some iso~$\alpha$, $s \equiv \IM f$ and~$p \equiv 1 \after f$.
As we saw before (using \sref{asrt-absorp-rule}),
    we see~$\id = f^\dagger \after f = \asrt_{p^2}$.
Thus~$\andthen{p}{p} = p^2 =  1 = \andthen{1}{1}$
    and so by uniqueness of the square root
    we get~$p = 1$.
Hence~$f$ is a comprehension.

Dualizing, we see~$\dagger$-epis are exactly quotients of sharp predicates.
Now it is clear that the~$\dagger$-partial isometries
are exactly the pristine maps.
\begin{point}{50}%
By the previous~$\pi_{(1 \after f)^\perp}$ is~$\dagger$-mono.
To see it is the~$\dagger$-kernel of~$f$,
we have to show it is the equalizer of~$f$ with~$0$.
Clearly~$(1 \after f )^\perp \after \pi_{(1 \after f)^\perp} = 1$
and so~$1 \after f \after \pi_{(1 \after f)^\perp} = 0$,
whence~$f \after \pi_{(1 \after f)^\perp} = 0$.
Assume~$f \after g = 0$ for some pure map~$g$.
Then~$1 \after f \after g = 0$,
    so~$1 \after f \leq \IMperp g$,
    hence~$(1 \after f)^\perp \geq \IM g$,
    so~$g = \pi_{(1 \after f)^\perp} \after g'$
    for a unique~$g'$. \qed
\end{point}
\end{point}
\end{point}
\begin{point}{60}[exc-purec-no-biproduct]{Exercise*}%
Let~$C$ be a~$\dagger$-effectus.
In this exercise, you will show that
    in general $\Pure C$
    does not have finite coproducts.
In particular, $\Pure C$ does not have~$\dagger$-biproducts,
    see~\cite[3.2.15]{heunenphd}.
Consider~$C \equiv \op\vN$.
Reasoning towards contradiction,
    assume~$(\Pure \op\vN)^{\mathsf{op}}$
    has a product~$ \C \xleftarrow{\pi_1} \scrA \xrightarrow{\pi_2} \C$
    of~$\C $ and~$\C$.
First show that for any non-zero pure map~$f\colon \scrA \to \C$,
    there exists a Hilbert space~$\scrH$, an element~$x \in \scrH$,
    a von Neumann algebra~$\scrC$ and
    isomorphism~$\varphi\colon \scrB(\scrH)\oplus \scrC \to \scrA$
    with~$f (varphi(T,c)) = \langle x, T x\rangle $.
Use this on~$\pi_1$, $\pi_2$ and~$\langle\id_\C,\id_\C\rangle $
    to show there is an
    isomorphism~$\varphi\colon  \scrB(\scrH) \oplus \scrC \cong \scrA$
    with~$\pi_1 \after \varphi = \left<x,\,(\cdot\,)x\right>$
    and~$\pi_2 \after \varphi = \left<y,(\,\cdot\,)y\right>$
    for some von Neumann algebra~$\scrC$,
    Hilbert space~$\scrH$ and $x,y \in \scrH$.
Press on: show~$\scrC$ is trivial.
    Then prove~$\dim \scrH \leq 2$.
    Almost there: show~$\dim \scrH \geq 2$.
    Derive a contradiction.
\end{point}
\begin{point}{70}[exc-purec-equal]{Exercise*}%
In this exercise you will show that~$\Pure(\op\vN)$ does not have
    all coequalizers.  It is helpful to first consider
    two concrete coequalizers in~$\mathsf{Hilb}$ and~$\op\vN$.
(These coequalizers capture unordered pairs of qubits.)

Write~$\sigma\colon \C^2 \otimes \C^2 \to \C^2 \otimes \C^2$
    for the unitary
    fixed by~$\sigma \ket{ij} = \ket{ji}$.
The equalizer of~$\sigma$ with~$\id$ exists in~$\mathsf{Hilb}$
    and is given by~$e_{\mathscr{S}}$, the inclusion
    of the subspace~$\mathscr{S}$
    spanned by~$\{\ket{00}, \ket{11}, \ket{01}+\ket{10}\}$
    into~$\C^2\otimes \C^2$.
    ($\mathscr{S}$ is called the symmetric tensor of~$\C^2$ with itself.)
    The coequalizer of~$\sigma$ with~$\id$ is simply~$e_{\mathscr{S}}^\dagger$.
    There is also an equalizer of~$\ad_{\sigma}\colon M_4 \to M_4$ with~$\id$
    in~$\vN$, but curiously enough, it is given
    by~$e \colon M_3 \oplus \C \to M_4$,
    where~$e(a,\lambda) =
            \ad_{e^\dagger_{\mathscr{S}}}(a) +
            \ad_{e^\dagger_{\mathscr{A}}}(\lambda)$
     and~$\scrA$ is the subspace spanned by~$\ket{01} - \ket{10}$.
     For the proof and similar results, see \cite{bags}.

Assume~$\xi \colon \scrC \to M_4$
    is some filter (i.e.~quotient in the opposite category)
    with~$\ad_\sigma \after \xi = \xi$.
Show that for each~$b \in \scrC$,
        we have~$\xi(b) \sigma = \sigma \xi(b)$.
    From this, derive (perhaps in the same way as \cite{bags}),
    that~$\xi(b) = p_\scrA \xi(b) p_\scrA+ p_\scrS \xi(b) p_\scrS$
        for all~$b\in \scrC$,
        where~$p_\scrA = e_\scrA e_\scrA^\dagger$
        and~$p_\scrS = e_\scrS e_\scrS^\dagger$.
Prove that either~$\xi(b) \leq p_\scrA$ for all~$b\in \scrC$
    or~$\xi(b) \leq p_\scrS$ for all~$b \in \scrC$.
Conclude~$\Pure(\op\vN)$ does not have all coequalizers.
\end{point}
\begin{point}{80}%
To summarize:
    the category~$\Pure(C)$
    of a~$\dagger$-effectus~$C$
    is a~$\dagger$-category,
    has~$\dagger$-kernels
    and each~$\dagger$-mono is a~$\dagger$-kernel,
    but in general it does not have~$\dagger$-biproducts or
        $\dagger$-equalizers,
        which shows it is not, in general,
        a pre-Hilbert category as in \cite[3.7.1]{heunenphd}.
\end{point}
\begin{point}{81}%
Recently Tull axiomatized~\cite{tull}
    the $\dagger$-category~$\mathsf{Hilb}_\sim$
        of Hilbert spaces with bounded operators between
        them modulo global phase (i.e.~$T=S$ in~$\mathsf{Hilb}_\sim$
        iff~$T = \lambda S$ for some~$\lambda \in \C, |\lambda|=1$.)
Write~$\mathsf{Hilb}^c_\sim$
    for the restruction of~$\mathsf{Hilb}_\sim$
    to contractive operators.
This category~$\mathsf{Hilb}^c_{\sim}$
    is a full~$\dagger$-subcategory of~$\Pure(\op\vN)$
    via the functor
\begin{equation*}
    \scrH \ \mapsto \ \scrB(\scrH),  \quad
    \bigl(T\colon \scrH \to \scrK\bigr)
        \ \mapsto \ \bigl(\ad_T \colon \scrB(\scrK) \to \scrB(\scrH)\bigr).
\end{equation*}
The category~$\mathsf{Hilb}^c_\sim$
    (and in fact~$\mathsf{Hilb}_\sim$)
    do not have all (co)products.
Tull shows that~$\mathsf{Hilb}_\sim$ does have
    \emph{phased biproducts}. \cite{tull,tull2018quotient}
This raises the question: do the phased biproducts
    of~$\mathsf{Hilb}_\sim$ extend to~$\Pure(\op\vN)$.
If this is the case, then a
    phased product of~$\C$ with itself in~$\Pure(\op\vN)$
    is given by: $\ketbra{0}{0}, \ketbra{1}{1}\colon \C  \to M_2$.
Unfortunately this leads to a contradiction.
It hasn't been determined yet whether~$\Pure(\op\vN)$
    has phased products at all.
\end{point}
\end{parsec}

\subsection{Sequential product}
\begin{parsec}{2250}%
\begin{point}{10}%
The sequential product --- the
    operation~$\andthen{a}{b} \equiv \sqrt{a}b\sqrt{a}$
    on the effects of a von Neumann algebra
    --- has been studied before.
In \cite{gudder2001sequential,gheondea2004sequential}
    Gudder and coauthors establish some basic algebraic properties of~$\&$
    on~$\scrB(\scrH)$.
Most of them are surprisingly hard to prove ---
    for instance, the proof that~$\andthen{a}{b}=\andthen{b}{a}$
    implies~$ab = ba$, requires
    the Fuglede--Rosenblum--Putnam Theorem.
\begin{point}{20}%
Then, in \cite{gudder2008characterization},
    Gudder and Latr\'emoli\`ere (G\&L)
    characterize~$\&$ on~$\scrB(\scrH)$ as the unique operation satisfying
    (for effects~$a,b$ and density matrix~$\varrho$)
\begin{multicols}{2}
\begin{enumerate}
    \item $\TR [ (\andthen{a}{\varrho}) b] = \TR [\varrho
                    (\andthen{a}{b})]$
    \item $\andthen{a}{1} = \andthen{1}{a} = a$
    \item $ \andthen{a}{(\andthen{a}{b})} =
            \andthen{(\andthen{a}{a})}{b} =
            \andthen{a^2}{b}$
        \item $a \mapsto \andthen{a}{b}$ is strongly continuous
\end{enumerate}    
\end{multicols}
\noindent
In \sref{positive-map-uniqueness} we characterized~$\asrt_a \colon b \mapsto \andthen{a}{b}$
    as the unique~$\diamond$-positive map with~$\andthen{a}{1} = a$.
How do these characterizations compare?
G\&L's first axiom plays a very similar role
        as the~$\diamond$-self adjointness of~$\diamond$-positive maps for us.
Their third axiom is somewhat related to~$\diamond$-positivity.
Their fourth axiom seems unrelated to our characterization
    and conversely the purity of our~$\diamond$-positive maps
    has no counterpart in their axioms.
\end{point}
\begin{point}{30}%
A few years earlier, Gudder and Greechie started
    \cite{gudder2001sequential}
    the abstract study of some of the algebraic properties of the
    sequential product on arbitrary effect algebras,
    which has been picked up by several other authors
    \cite{li2011sequential,gudder2005open,shen2009not,gudder2005uniqueness,jun2009remarks,weihua2009uniqueness,tkadlec2008atomic,jia2010entropy,arias2004almost,van2018three,weteringseqprod}.
\end{point}
\end{point}
\begin{point}{40}{Definition}%
A \Define{sequential effect algebra} (SEA)
    \index{SEA}
    \index{effect algebra!sequential-}
    is an effect algebra~$E$ together
    with a binary operation~$\&$
    satisfying
\begin{enumerate}
    \item[(S1)] if~$a \perp b$,
        then~$\andthen{c}{a} \perp \andthen{c}{b}$
        and~$(\andthen{c}{a}) \ovee (\andthen{c}{b})
                = \andthen{c}{(a \ovee b)}$;
    \item[(S2)] $\andthen{1}{a} =a$;
    \item[(S3)] if~$\andthen{a}{b} = 0$,
            then~$\andthen{a}{b} = \andthen{b}{a}$;
    \item[(S4)]
        if~$\andthen{a}{b} = \andthen{b}{a}$,
        then~$\andthen{a}{b^\perp} = \andthen{b^\perp}{a}$
        and~$\andthen{(\andthen{a}{b})}{c}
            = \andthen{a}{(\andthen{b}{c})} $ \emph{and}
    \item[(S5)]
        if~$\andthen{c}{a} = \andthen{a}{c}$,
        $\andthen{c}{b} = \andthen{b}{c}$ and~$a \perp b$,\\
        then~$\andthen{c}{(\andthen{a}{b})}
                = \andthen{(\andthen{a}{b})}{c}$
                and~$\andthen{c}{(a \ovee b)} = \andthen{(a \ovee b)}{c}$.
\end{enumerate}
\end{point}
\spacingfix{}
\begin{point}{50}{Examples}%
The effect algebra~$[0,1]_\scrA$
    of effects on a von Neumann algebra~$\scrA$
    is a sequential effect algebra with~$\andthen{a}{b} = \sqrt{a}b\sqrt{a}$.
Any commutative effect monoid is a sequential effect algebra
    with~$\andthen{a}{b} = a \odot b$.
\end{point}
\begin{point}{60}{Proposition}%
In a~$\dagger$-effectus,
    the set of predicates~$\Pred X$ on any object~$X$
    with~$\andthen{p}{q} = q \after \asrt_p$
    satisfies axioms (S1), (S2) and (S3) of a SEA.
\begin{point}{70}{Proof}%
The proofs of (S1) and (S2) are obvious.
To show~(S3), assume~$\andthen{p}{q} = 0$
    for some predicates~$p,q$.
    Then~$1 \after \asrt_q \after \asrt_p  = \andthen{p}{q}= 0$
and so~$\asrt_q \after \asrt_p = 0 = \asrt_0$.
Applying the dagger, we find
$0 = \asrt_0^\dagger = \asrt_p^\dagger \after \asrt_q^\dagger
                = \asrt_p \after \asrt_q$.
    Thus~$\andthen{q}{p} = 0 = \andthen{p}{q}$, as desired. \qed
\end{point}
\begin{point}{80}{Remark}%
It's unclear whether the set of predicates in a~$\dagger$-effectus
    forms a sequential effect algebra.
A seemingly related open problem
    is whether~$\andthen{p}{q} = \andthen{q}{p}$
    implies that~$\asrt_p \after \asrt_q = \asrt_q \after \asrt_p$.
\end{point}
\end{point}
\end{parsec}

\subsection{Homological categories}
\begin{parsec}{2260}%
\begin{point}{10}%
Next we compare our effectuses
    to Grandis' homological categories \cite{grandis},
    which generalize abelian categories
    used in homological algebra.
First, we need a lemma.
\end{point}
\begin{point}{20}[homology-lemma]{Lemma}%
Suppose~$s,t$ are any sharp predicates on the same object in a~$\dagger$-effectus.
    If~$s^\perp \leq t$,
    then~$\andthen{s}{t}$ is sharp.
\begin{point}{30}{Proof}%
By \sref{diamond-oml},
it is sufficient to show
    $\andthen{s}{t^\perp}$ is sharp
as~$(\andthen{s}{t})^\perp = \andthen{s}{t^\perp} \ovee s^\perp$.
Note~$s^\perp \perp t^\perp$,
so by \sref{perp-sharp-is-orth},
    we find~$\andthen{t^\perp}{s^\perp} = 0$,
    hence~$\andthen{t^\perp}{s} = t^\perp$ is sharp.
By the same reasoning as in the final
    part of the proof of \sref{dagger-iso-mu},
    there exists a (unique) pristine map~$l$
    with~$\asrt_s \after \asrt_{t^\perp}
        = l\after \asrt_{\andthen{t^\perp}{s}}$,
        $l \after l^\dagger = \asrt_{\ceil{\andthen{t^\perp}{s}}}$
        and $l^\dagger \after l = \asrt_{\ceil{\andthen{s}{t^\perp}}}$.
It follows~$l \after \asrt_{\andthen{t^\perp}{s}} \after l^\dagger
    = \asrt_{\andthen{s}{t^\perp}}$ and so
\begin{align*}
    \asrt_{\andthen{s}{t^\perp}}^2
    & \ = \ l \after \asrt_{\andthen{t^\perp}{s}} \after
            l^\dagger \after l \after 
            \asrt_{\andthen{t^\perp}{s}} \after l^\dagger \\
    & \ = \ 
    l \after \asrt_{\andthen{t^\perp}{s}} \after
            \asrt_{\ceil{\andthen{s}{t^\perp}}} \after
            \asrt_{\andthen{t^\perp}{s}} \after l^\dagger \\
    & \ \overset{\mathclap{\sref{asrt-absorp-rule}}}{=} \ 
    l \after \asrt_{\andthen{t^\perp}{s}} \after
            \asrt_{\andthen{t^\perp}{s}} \after l^\dagger \\
    & \ \overset{\mathclap{\sref{sharp-prop}}}{=} \ 
    l \after \asrt_{\andthen{t^\perp}{s}} \after
            l^\dagger \\
    & \ = \ 
     \asrt_{\andthen{s}{t^\perp}},
\end{align*}
    which shows~$\andthen{s}{t^\perp}$ is sharp. \qed
\end{point}
\end{point}
\begin{point}{40}{Definition}%
A category~$C$ is a \Define{pointed semiexact category} \cite[\S1.1]{grandis}
    \index{semiexact category!pointed}
    if~$C$ has a zero object and all kernels and cokernels.
In a pointed semiexact category:
\begin{enumerate}
\item
We put a partial order on kernels in the usual way:
we write~$\Define{n \leq m}$
        \index{*leq@$\leq$!among kernels}
for kernels~$n\colon A \to X$ and~$m \colon B \to X$,
    if there is an~$f\colon A \to B$
    with~$n = m \after f$.
If both~$n \leq m$ and~$m \leq n$,
    then we write~$n \approx m$.
We denote the poset of kernels on~$A$ modulo~$\approx$
        by~$\Define{\Nsb A}$. (\textbf{N}ormal \textbf{s}ubo\textbf{b}jects,
        \index{Nsb@$\Nsb A$}
            \cite[\S1.5]{grandis}.)
\item
For every map~$f$
    there is a unique map~$g$
        with~$f = (\ker \cok f ) \after g \after (\cok \ker f)$.
    If this~$g$ is an iso, then~$f$ is called \Define{exact}.
        \index{exact}
\end{enumerate}
A pointed semiexact category~$C$
    is a \Define{pointed homological category} if
    \index{homological category!pointed}
\begin{enumerate}
    \item kernels are closed under composition;
    \item cokernels are closed under composition \emph{and}
    \item \emph{(homology axiom)}
        for any kernel~$m\colon M \to A$
        and cokernel~$q\colon A \to Q$
        with~$\ker q \leq m$,
        the composition~$q \after m$ is exact.
\end{enumerate}
\end{point}
\spacingfix{}
\begin{point}{50}{Theorem}%
Any~$\dagger$-effectus (in partial form) is a (pointed) homological category:
\begin{enumerate}
    \item 
        a map is a kernel iff it is a comprehension;
    \item
        a map is a cokernel iff it is a quotient of a sharp predicate \emph{and}
    \item
        a map is exact iff it is pristine.
\end{enumerate}
\spacingfix{}
\begin{point}{60}{Proof}%
Let~$C$ be a~$\dagger$-effectus.
By~\sref{effectus-kernels} and \sref{effectus-cokernels},
    we know~$C$ has all kernels and cokernels:
    a kernel of~$f$ is exactly a comprehension of~$(1 \after f)^\perp$
    and a cokernel of~$f$ is exactly a quotient of~$\IM f$.
By~\sref{upm-closed},
    comprehensions are closed under composition,
    and so kernels are closed under composition as well.
As for the cokernels,
    assume~$\zeta_1, \zeta_2$ are two composable cokernels.
Both~$\zeta_1$ and~$\zeta_2$ are quotients of a sharp predicate.
By \sref{quotients-composition} the composition  $\zeta_1 \after \zeta_2$
    is a quotient as well.
The map~$\zeta_2$ is sharp by \sref{quotients-composition}
    and so the predicate~$1 \after \zeta_1 \after \zeta_2$ is sharp.
Hence~$\zeta_1 \after \zeta_2$ is a cokernel too.
Unfolding definitions, it is easy to see exact maps correspond precisely
    to pristine maps.
\begin{point}{70}{Homology axiom}%
We have to show~$q \after m$
    is exact for any cokernel~$q$ and kernel~$m$ with~$\ker q \leq m$.
As~$q \after m$ is pure,
    it is sufficient to show~$1 \after q \after m$ is sharp.
By \sref{compr-is-full},
the assumption~$\ker q \leq m$
    is equivalent to~$(1 \after q)^\perp \equiv \IM \ker q \leq \IM m$
    and so~$\IMperp m \leq 1 \after q$.
Our lemma \sref{homology-lemma}
    shows~$\andthen{(\IM m)}{(1 \after q)}
    \equiv 1 \after q \after m \after m^\dagger$ is sharp
    and so
\begin{align*}
    \ceil{1 \after q \after m}
    &\ = \ \ceil{1 \after q \after m} \after m^\dagger \after m 
        &\quad&\text{as comprehensions are~$\dagger$-mono}\\
    &\ = \ \ceil{1 \after q \after m \after m^\dagger} \after m 
        &&\text{as $m^\dagger$ is sharp} \\
    &\ = \ 1 \after q  \after  m \after m^\dagger \after m  \\
    &\ = \ 1 \after q  \after  m,
\end{align*}
which show~$q \after m$ is sharp, as desired. \qed
\end{point}
\end{point}
\end{point}
\end{parsec}
\begin{parsec}{2270}%
\begin{point}{10}%
In any homological category, a generalization
    of the famous \emph{Snake Lemma} holds.
Before we can discuss it, we need some more homological category theory.
\end{point}
\begin{point}{20}{Definition}%
Let a (pointed) semiexact category be given.
\begin{enumerate}
\item
    We say the diagram~$A \xrightarrow{f} B \xrightarrow{g} C$
        is \Define{exact at} $B$,
        \index{exact!at}
        if~$\ker \cok (f) \approx \ker (g)$.
We say~$A_1 \xrightarrow{f_1} A_2 \to \cdots
            \to A_{n}$
            is a \Define{(long) exact sequence}
            \index{exact!sequence}
            if it is exact at~$A_2, \ldots, A_{n-1}$.
\item
The poset~$\Nsb A$ of kernels modulo~$\approx$
        is a bounded lattice with minimum~$0$ and maximum~$1\equiv \id$.
        \cite[\S1.5]{grandis}
\item
For any~$f\colon A \to B$,
    define~$f_*\colon \Nsb A \leftrightarrows \Nsb B \colon f^*$
        by~$\Define{f_*} (k) = \ker \cok (f \after k)$
        and~$\Define{f^*}(k) = \ker ((\cok k) \after f)$.
        \index{*par3@$(\ )_*$, $(\ )^*$}
\item
We say~$f$ is \Define{left-modular at} $k$
        \index{modular}
    if~$f^* (f_* (k)) = k \vee f^*(0)$
        and $f$ is \Define{right-modular at} $k$
        if~$f_*(f^* (k)) = k \wedge f_*(1)$.
We call~$f$ simply \Define{left-modular} (resp.~right-modular)
    if it is left-modular (resp.~left-modular) at every~$k$.
We say~$f$ is \Define{modular}
    if it is both left- and right-modular.
\end{enumerate}
\end{point}
\spacingfix{}
\begin{point}{30}{Example}%
In a~$\dagger$-effectus, the previous concepts are related
to our familiar notions as follows.
\begin{enumerate}
\item
    A diagram~$A \xrightarrow{f} B \xrightarrow{g} C$
        is exact at $B$
        if and only if~$\IMperp f = \ceil{1 \after f}$.
\item
    The lattice~$\Nsb A$ is isomorphic to the lattice~$\SPred A$
        via~$k \mapsto \IM k$.
\item
    For any~$f$, we
        have~$\IM f_*(k) = f_\diamond (\IM k)$
        and~$\IM f^*(k) = f^\BOX (\IM k)$.
\item
    A map~$f$ is left-modular at~$k$
            iff~$f^\BOX (f_\diamond ( \IM k)) = 
                (\IM k) \vee \ceil{1 \after f}$
        and right-modular at~$k$
            iff~$f_\diamond (f^\BOX ( \IM k)) = 
                (\IM k) \wedge \IM f$.
\end{enumerate}
\spacingfix{}
\begin{point}{40}{Notation}%
As~$\Nsb A$ and~$\SPred A$ are isomorphic,
    we will abbreviate~``$f$ is left-modular at~$\pi_s$''
    by~``$f$ is left-modular at~$s$''
    and similarly for right-modularity.
\end{point}
\end{point}
\begin{point}{50}[diamondboxlemma]{Lemma}%
In a~$\dagger$-effectus,
    we have~$\pi^\BOX \after \pi_\diamond = \id$
    and~$\zeta_\diamond \after \zeta^\BOX = \id$
    for any comprehension~$\pi$
    and sharp quotient~$\zeta$.
\begin{point}{60}{Proof}%
We start with~$\pi$.  The trick is to use~$\pi^\dagger$ is sharp:
    $\pi^\BOX \after \pi_\diamond (s) =
    \pi^\BOX \after (\pi^\dagger)^\diamond (s) =
    \lceil \lceil s \after \pi^\dagger \rceil^\perp \after \pi\rceil^\perp =
    \lceil(s \after \pi^\dagger)^\perp \after \pi\rceil^\perp =
    \lceil(s \after \pi^\dagger \after \pi)^\perp\rceil^\perp =
    \lceil s ^\perp\rceil^\perp = s $.
As for the other equality:
    $\zeta_\diamond \after \zeta^\BOX (s)
        = (\zeta^\dagger)^\diamond \after \zeta^\BOX (s)
        = \lceil\lceil s^\perp \after  \zeta \rceil^\perp \after \zeta^\dagger\rceil
        = \lceil(\lceil s^\perp \after  \zeta \rceil \after \zeta^\dagger)^\perp\rceil
        = \lceil(\lceil s^\perp\rceil \after \zeta \after \zeta^\dagger)^\perp\rceil
        = \lceil \lceil s^\perp\rceil^\perp \rceil = s $. \qed
\end{point}
\end{point}
\end{parsec}
\begin{parsec}{2280}%
\begin{point}{10}%
The following is a translation of (a part of) Grandis' Snake
    Lemma~\cite[\S3.4]{grandis}
    into~$\dagger$-effectus jargon.
We include a proof as it highlights how
    the reasoning in a homological category is subtly different
    from the reasoning we already saw in a~$\dagger$-effectus:
    maps~$l^\BOX$ appear, where we would have expected~$l^\diamond$.
\end{point}
\begin{point}{20}{Snake Lemma}%
Suppose we have a commuting
    diagram in~$\dagger$-effectus as follows.
\begin{equation*}
    \xymatrix{
        &A \ar[d]^{a} \ar[r]^f
    &B \ar[d]^{b} \ar[r]^g
    &C \ar[d]^{c} \ar[r]^0
    &0
    \\0 \ar[r]^{0}
    &A' \ar[r]^{h}
    &B' \ar[r]^{k}
    &C'
}
\end{equation*}
Furthermore, assume
\begin{enumerate}
    \item the diagram has exact rows;
    \item  $b$ is left-modular over~$\IM f$;
    \item $b$ is right-modular over~$\IM h$;
    \item $f$ is right-modular over~$\ceil{1 \after b}^\perp$ \emph{and}
    \item $k$ is left-modular over~$\IM b$.
\end{enumerate}
These additional assumptions are equivalent to the following conditions.
    \begin{multicols}{2}
    \begin{enumerate}
    \item
        $b^\BOX \after b_\diamond (\IM f) = \ceil{1 \after b}^\perp \vee \IM f$,
    \item
        $b_\diamond \after b^\BOX (\IM h)
                = (\IM h) \wedge \IM b$,
    \item
        $k^\BOX \after k_\diamond (\IM b) = (\IM h) \vee \IM b$,
    \item
        %$f_\diamond \after f^\BOX (\ceil{1 \after b}^\perp)
        $f_\diamond \after f^\BOX \after b^\BOX(0)
                = \ceil{1 \after b}^\perp \wedge \IM f$,
    \item
        $\IMperp f = \ceil{1\after g}$,
    \item
        $\IMperp h = \ceil{1\after k}$,
    \item
        $g$ is a sharp quotient \emph{and}
    \item
        $h$ is a comprehension.
    \end{enumerate}
\end{multicols}\noindent
\emph{Completing the diagram},
write~$a_\pi, b_\pi, c_\pi$ for kernels of~$a$, $b$ and~$c$, respectively;
$a_\zeta$, $b_\zeta$, $c_\zeta$ for the cokernels
    \emph{and}~$\overline{f},\overline{g},\overline{h}, \overline{k}$
    for the induced maps between the kernels and cokernels.
It is easy to see
\begin{align*}
   \overline{f} &= b_\pi^\dagger \after f \after a_\pi
     &   \overline{g} &= c_\pi^\dagger \after g \after b_\pi
     &   \overline{h} &= b_\zeta \after h \after a_\zeta^\dagger
     &   \overline{k} &= c_\zeta \after k \after b_\zeta^\dagger.
\end{align*}
Then: there is a map~$d \colon \ker c \to \cok a$
        that turns\footnote{%
            Beware: we use~$\ker a$ (and~$\cok$) inconsistently ---
         earlier in the text~$\ker f$ refers
            to a kernel map --- here it
            refers to a kernel object instead.}
\begin{equation}\label{thesnake}
    \xymatrix{
        \ker a \ar[r]^{\overline{f}}
        &\ker b \ar[r]^{\overline{g}}
        &\ker c \ar[r]^{d}
        &\cok a \ar[r]^{\overline{h}}
        &\cok b \ar[r]^{\overline{k}}
        &\cok c 
    }
\end{equation}
into a (long) exact sequence. This reveals the snake:
\begin{equation*}
\begin{tikzpicture}[
        every edge/.append style={font=\scriptsize}
        ]
    \matrix[matrix of math nodes,
            column sep={45pt,between origins},
            row sep={45pt,between origins},
            nodes={asymmetrical rectangle}] (s) {
        & |[name=ka]| \ker a
        & |[name=kb]| \ker b
        & |[name=kc]| \ker c \\
        & |[name=A]| A
        & |[name=B]| B
        & |[name=C]| C 
        & |[name=01]| 0\\
         |[name=02]| 0
        & |[name=A']| A'
        & |[name=B']| B'
        & |[name=C']| C' \\
        & |[name=ca]| \cok a
        & |[name=cb]| \cok b
        & |[name=cc]| \cok c \\
    };
\draw[->]
    (ka) edge node[auto] {$a_\pi$} (A)
    (kb) edge node[auto] {$b_\pi$} (B)
    (kc) edge node[auto] {$c_\pi$} (C)
    (A') edge node[auto] {$a_\zeta$} (ca)
    (B') edge node[auto] {$b_\zeta$} (cb)
    (C') edge node[auto] {$c_\zeta$} (cc)
    (C) edge node[auto] {$0$} (01)
    (02) edge node[auto] {$0$} (A')
    (A) edge node[auto] {$f$} (B)
    (B) edge node[auto] {$g$} (C)
    (A') edge node[auto] {$h$} (B')
    (B') edge node[auto] {$k$} (C')
    (A) edge node[auto] {$a$} (A')
    (B) edge node[auto] {$b$} (B')
    (C) edge node[auto] {$c$} (C')
;
\draw[->,darkgreen] 
    (ka) edge node[auto,text=black] {$\overline{f}$} (kb)
    (kb) edge node[auto,text=black] {$\overline{g}$} (kc)
    (ca) edge node[auto,below,text=black] {$\overline{h}$} (cb)
    (cb) edge node[auto,below,text=black] {$\overline{k}$} (cc)
;
\draw[->,darkgreen,rounded corners]
    (kc) 
        -| node[auto,text=black,pos=.7,font=\scriptsize]{$d$}
            ($(01.east)+(.5,0)$)
        |- ($(B)!.35!(B')$)
        -| ($(02.west)+(-.5,0)$)
        |- (ca);
;
\end{tikzpicture}
\end{equation*}
\spacingfix{}
\begin{point}{30}{Proof}%
Before we start, we give an overview of the proof.
We will first show
    there are comprehensions~$m,h'$
        and sharp quotients~$g',v$
        such that the left and right faces of the following cube commute.
\begin{equation*}
\xymatrix{
    & B \ar[rr]^b \ar@{->>}[dd]^g
    && B' \ar@{->>}[dd]^v
    \\ Z \ar@{^{(}->}[ru]^m \ar@{->>}[dd]_{g'}
    && A' \ar@{^{(}->}[ru]^h \ar@{->>}[dd]_{a_\zeta}
\\& C
&& Q
    \\ \ker c \ar@{^{(}->}[ru]_{c_\pi} \ar[rr]_{d}
    && \cok a \ar@{^{(}->}[ru]_{h'}
}
\end{equation*}
(These two faces are, what Grandis calls, subquotients.)
Then we will prove using exactness of the rows
    that there is a unique lifting of~$b$
    to~$b'\colon Z \to A'$ along the comprehensions~$m,h$.
In turn, $d$ is defined as the unique lifting of~$b'$
    along the sharp quotients~$h', a_\zeta$.
($d$ is the map regularly induced by~$b$ along the two subquotients.)
Before demonstrating the exactness of the snake,
    we need some non-trivial identities for~$d_\diamond$ and~$d^\BOX$.
(These also follow from Grandis' calculus of direct and inverse images
    along an induced morphism.)
\begin{point}{40}%
Define~$v \equiv \xi_{h_\diamond(\IM a)}$. 
As~$(v \after h)^\diamond(1) = h^\diamond ( h_\diamond(\IM a)) = \IM a
        = 1\after a_\zeta$,
        there must be a total pure map~$h'$
        such that~$v \after h = h' \after a_\zeta$.
Note~$\ceil{1 \after v}^\perp = h_\diamond(\IM a) \leq h_\diamond(1) = \IM h$.
So by the homology axiom, we know~$v \after h$ is pristine
    and so~$h'$ must be a comprehension.
Moving to the other side,
    define~$m \equiv \pi_{g^\BOX(\ceil{1\after c}^\perp)}$.
Reasoning in a similar way,
    we see there is a unique sharp quotient~$g'$
    with~$g \after m = c_\pi \after g'$.
\end{point}
\begin{point}{50}%
To show~$b$ lifts along~$m$ and~$h$,
    it suffices (by the universal property of comprehensions)
    to show~$(b \after m)_\diamond(1)
        \equiv \IM b\after m \leq \IM h
                    \equiv h_\diamond(1)$.
\begin{alignat*}{2}
    (b_\diamond \after m_\diamond) (1)
    &\ =\  (b_\diamond \after g^\BOX \after c^\BOX )(0) &\qquad& \text{by dfn.~$m$}\\
    &\ = \ (b_\diamond \after b^\BOX \after k^\BOX )(0) \\
    &\ = \ (b_\diamond \after b^\BOX \after h_\diamond )(1) && \text{by exactness at~$B'$} \\
    &\ \leq \  h_\diamond (1) && \text{as $b_\diamond \dashv b^\BOX$}.
\end{alignat*}
So there is a unique~$b'\colon Z \to A'$ with~$h \after b' = b \after m$.
Before we continue,
    we note~$b' = h^\dagger \after b \after m$
    (as ~$h^\dagger \after h= \id$)
    and so surprisingly:
\begin{equation}\label{snakeceilbprime}
    m^\BOX \after b^\BOX \after h_\BOX \ = \ 
    b'^\BOX \ \overset{\mathclap{\sref{diamondboxlemma}}}{=} \   b'^\BOX \after h^\BOX \after h_\diamond
            \ =\    m^\BOX \after b^\BOX \after h_\diamond.
\end{equation}
Next, to show~$b'$ lifts along~$g'$ and~$a_\zeta$,
    it is sufficient to prove~$(b'^\BOX \after a_\zeta^\BOX)(0)
        = \ceil{1\after b' \after a_\zeta}^\perp
                \geq \ceil{1 \after g'}^\perp
                = g'^\BOX (0)$.
\begin{alignat*}{2}
    (b'^\BOX \after a_\zeta^\BOX)(0)
        & \ = \ (b'^\BOX \after a_\diamond)(1)
                &\qquad& \text{by dfn.~$a_\zeta$} \\
        & \ = \ (m^\BOX \after b^\BOX \after h_\diamond \after a_\diamond)(1)
                && \text{by \eqref{snakeceilbprime}} \\
        & \ = \ (m^\BOX \after b^\BOX \after b_\diamond \after f_\diamond)(1) \\
        & \ \geq \ (m^\BOX \after f_\diamond )(1) 
                && \text{as $b_\diamond \dashv b^\BOX$}\\
        & \ =\ (m^\BOX \after g_\BOX )(0) 
                && \text{by exactness at $B$}\\
        & \ = \ (g'^\BOX \after (c_\pi)_\BOX )(0) 
                && \text{by dfn.~$g'$}\\
        & \ \geq \ g'^\BOX (0).
\end{alignat*}
Thus there is a unique~$d\colon \ker c \to \cok a$
    with~$d \after g' = a_\zeta \after b'$.
More concretely, using~$g' \after g'^\dagger = \id$,
    we see
    $ d  =  a_\zeta \after h^\dagger \after b \after m \after 
                g'^\dagger$.
\end{point}
\begin{point}{60}%
Our next goal is to show
\begin{align*}
    d_\diamond &\ =\ 
            (a_\zeta)_\diamond \after h^\BOX \after
                b_\diamond \after m_\diamond \after g'^\BOX
            \ = \ h'^\BOX \after v_\diamond \after
                b_\diamond \after m_\diamond \after g'^\BOX 
                \numberthis\label{snakedidents}
                \\
            &\ = \ (a_\zeta)_\diamond \after h^\BOX \after
                b_\diamond\after g^\BOX \after (c_\pi)_\diamond
            \ = \ h'^\BOX \after v_\diamond \after
                b_\diamond\after g^\BOX \after (c_\pi)_\diamond.
\end{align*}
As a first step,
    note~$\IM m = (g^\BOX \after c_\diamond )(1) \geq
        (g^\BOX\after (c_\pi)_\diamond)(s)$
    and so
\begin{alignat*}{2}
    (g^\BOX\after (c_\pi)_\diamond)(s)
       & \ =\ 
    (g^\BOX\after (c_\pi)_\diamond(s)) \wedge  (\IM m) \\
    &\ = \ (m_\diamond \after m^\BOX \after 
    g^\BOX\after (c_\pi)_\diamond)(s)
    &\qquad&\text{by \sref{spred-infimum}} \\
    &\ = \ (m_\diamond \after g'^\BOX \after
    c_\pi^\BOX\after (c_\pi)_\diamond)(s) &&\text{by dfn.~$g'$}\\
    &\ = \ (m_\diamond \after g'^\BOX)(s)
    &&\text{by \sref{diamondboxlemma}}.
\end{alignat*}
and so we only have to show the first two equalities of \eqref{snakedidents}.
The first is easy:
\begin{alignat*}{2}
    d_\diamond
    &\ = \ d_\diamond \after g'_\diamond \after g'^\BOX 
    &\qquad&\text{by \sref{diamondboxlemma}}\\
    & \ = \ (a_\zeta)_\diamond \after b'_\diamond \after g'^\BOX 
                &&\text{by dfn.~$d$}\\
    & \ = \ (a_\zeta)_\diamond \after
                    h^\BOX \after h_\diamond \after 
    b'_\diamond \after g'^\BOX
    &&\text{by \sref{diamondboxlemma}}\\
    &\ = \ (a_\zeta)_\diamond \after
                    h^\BOX \after b_\diamond \after
    m_\diamond \after g'^\BOX &&\text{by dfn.~$b'$}.
\end{alignat*}
The second is trickier.  We need some preparation.
    Recall~$v^\BOX(0) \leq h_\diamond(1)$
    and so
\begin{equation}\label{snakehdiamondone}
    h_\diamond(1) \ = \ 
    v^\BOX(0)  \vee h_\diamond(1) \ \overset{\mathclap{\sref{spred-sup}}}{=} \ 
    (v^\BOX \after v_\diamond \after h_\diamond)(1) \ = \ 
    (v^\BOX \after h'_\diamond)(1).
\end{equation}
As~$h_\diamond \after a^\BOX_\zeta$ is injective,
    the last equality follows from
\begin{alignat*}{2}
    &(h_\diamond \after a_\zeta^\BOX \after (a_\zeta)_\diamond \after
                    h^\BOX \after b_\diamond \after
    m_\diamond \after g'^\BOX )(s)\\
    &\qquad \ =\ 
    (h_\diamond ( h^\BOX \after b_\diamond \after
    m_\diamond \after g'^\BOX  )(s) \vee a_\zeta^\BOX(0)) 
        &\qquad&\text{by \sref{spred-sup}}\\
    &\qquad \ =\ 
    (h_\diamond \after  h^\BOX \after b_\diamond \after
    m_\diamond \after g'^\BOX )(s) \vee (h_\diamond \after a_\zeta^\BOX)(0)
        &\qquad&\text{as~$h_\diamond \dashv h^\BOX$}\\
    &\qquad \ =\ 
    (h_\diamond \after  h^\BOX \after b_\diamond \after
    m_\diamond \after g'^\BOX )(s) \vee v^\BOX(0)
        &\qquad&\text{by dfns.~$v$,$a_\zeta$}\\
    &\qquad \ =\ 
    (\,( b_\diamond \after
    m_\diamond \after g'^\BOX ) (s) \wedge h_\diamond(1)\,) \vee v^\BOX(0)
        &\qquad&\text{by \sref{spred-infimum}}\\
    &\qquad \ =\ 
    (b_\diamond \after
    m_\diamond \after g'^\BOX  )(s)  \vee v^\BOX(0)
        &\qquad&\text{as $b_\diamond \after m_\diamond (1) \leq h_\diamond(1)$}\\
    &\qquad \ =\ 
    (\, (b_\diamond \after
    m_\diamond \after g'^\BOX ) (s)  \vee v^\BOX(0)\,
    ) \wedge h_\diamond(1)
        &\qquad&\text{as $v^\BOX(0) \leq h_\diamond(1)$}\\
    &\qquad \ =\ 
    (v^\BOX \after v_\diamond \after  b_\diamond \after
    m_\diamond \after g'^\BOX  )(s)  \wedge h_\diamond(1)
        &\qquad&\text{by \sref{spred-sup}} \\
    &\qquad \ =\ 
    v^\BOX \after v_\diamond \after  b_\diamond \after
    m_\diamond \after g'^\BOX  (s)  \wedge v^\BOX \after h'_\diamond(1)
        &\qquad&\text{by \eqref{snakehdiamondone} } \\
    &\qquad \ =\ 
    v^\BOX (\,(v_\diamond \after  b_\diamond \after
    m_\diamond \after g'^\BOX ) (s)  \wedge  h'_\diamond(1)\,)
        &\qquad&\text{as~$v^\BOX \vdash v_\diamond$} \\
    &\qquad \ =\ 
    (v^\BOX \after h'_\diamond
    \after h'^\BOX \after  v_\diamond \after  b_\diamond \after
    m_\diamond \after g'^\BOX  )(s)
        &\qquad&\text{by \sref{spred-infimum}}\\
    &\qquad \ =\ 
    (h_\diamond \after a_\zeta^\BOX \after
    h'^\BOX \after  v_\diamond \after  b_\diamond \after
    m_\diamond \after g'^\BOX )(s),
\end{alignat*}
where~$
    h_\diamond \after a_\zeta^\BOX \after =
    v^\BOX \after h'_\diamond$
    is proven in the same way as $g^\BOX \after (c_\pi)_\diamond =
                                    m_\diamond \after g'^\BOX$.
\end{point}
\begin{point}{70}%
We are ready to show exactness of \eqref{thesnake}.
We will first derive exactness in~$\cok a$
    from the right-modularity of~$b$ over~$\IM h$.
\begin{alignat*}{2}%
    d_\diamond(1) &
    \ = \  ((a_\zeta)_\diamond \after h^\BOX \after b_\diamond
            \after g^\BOX \after (c_\pi)_\diamond)(1)
            &\qquad&\text{by \eqref{snakedidents}} \\
    & \ = \ ((a_\zeta)_\diamond \after h^\BOX \after 
                b_\diamond \after g^\BOX \after c^\BOX)(0) 
                &&\text{by dfn.~$c_\pi$}\\
    & \ = \ ((a_\zeta)_\diamond \after h^\BOX \after 
    b_\diamond \after b^\BOX \after k^\BOX)(0) \\
    & \ = \ ((a_\zeta)_\diamond \after h^\BOX \after 
    b_\diamond \after b^\BOX \after h_\diamond)(1) 
    &\qquad&\text{by exactness in~$B'$}\\
    & \ = \ ((a_\zeta)_\diamond \after h^\BOX) (
    b_\diamond(1) \wedge h_\diamond(1))
        &&\text{by right-modularity}\\
    & \ = \ ((a_\zeta)_\diamond \after h^\BOX \after 
    h_\diamond \after h^\BOX \after b_\diamond)(1)
    &&\text{by \sref{spred-infimum}} \\
    & \ = \ ((a_\zeta)_\diamond \after h^\BOX
    \after b_\diamond)(1)
    &&\text{as $h^\BOX \vdash h_\diamond$}\\
    & \ = \ ((a_\zeta)_\diamond \after h^\BOX
    \after b_\zeta^\BOX)(0)
    &&\text{by dfn.~$b_\zeta$} \\
    & \ = \ ((a_\zeta)_\diamond \after a_\zeta^\BOX
    \after \smash{\overline{h}^\BOX} )(0)
    &&\text{by dfn.~$\overline{h}$}\\
    & \ = \ \smash{\overline{h}^\BOX} (0)
    &&\text{by \sref{diamondboxlemma}.}
\end{alignat*}
    By a dual argument
    one derives exactness of~$\eqref{thesnake}$
        in~$\ker c$ from the left-modularity
        of~$b$ over~$\IM f$.
\end{point}
\begin{point}{80}%
    To show exactness of
    $\eqref{thesnake}$ in~$\cok b$,
    we will use left-modularity of~$k$ in~$\IM b$:
\begin{alignat*}{2}
    \overline{h}_\diamond(1) 
        & \ = \ (\overline{h}_\diamond \after (a_\zeta)_\diamond
                \after a_\zeta^\BOX)(1)
                &\qquad&\text{by \sref{diamondboxlemma}} \\
        & \ = \ (\overline{h}_\diamond \after (a_\zeta)_\diamond)(1)
                &&\text{as $a_\zeta^\BOX \dashv (a_\zeta)_\diamond$}\\
        & \ = \ ((b_\zeta)_\diamond \after h_\diamond)(1)
                &&\text{by dfn.~$\overline{h}$}\\
        & \ = \ ((b_\zeta)_\diamond \after
                b^\BOX_\zeta \after (b_\zeta)_\diamond \after
        h_\diamond)(1)
                &&\text{as~$(b_\zeta)_\diamond \dashv b^\BOX_\zeta$}\\
        & \ = \ (b_\zeta)_\diamond (
                    h_\diamond(1) \vee b_\zeta^\BOX(0))
                &&\text{by \sref{spred-sup}}\\
        & \ = \ (b_\zeta)_\diamond (
                    h_\diamond(1) \vee b_\diamond(1))
                &&\text{by dfn.~$b_\zeta$}\\
        & \ = \ (b_\zeta)_\diamond (
                    k^\BOX(0) \vee b_\diamond(1))
                &&\text{exactness in~$B'$}\\
        & \ = \ ((b_\zeta)_\diamond \after
                    k^\BOX \after k_\diamond \after b_\diamond)(1)
                &&\text{by left-modularity}\\
        & \ = \ ((b_\zeta)_\diamond \after
                    k^\BOX \after c_\diamond \after g_\diamond)(1) \\
        & \ = \ ((b_\zeta)_\diamond \after
                    k^\BOX \after c_\diamond )(1) 
                    && \text{as quotients are faithful}\\
        & \ = \ ((b_\zeta)_\diamond \after
                    k^\BOX \after c_\zeta^\BOX)  (0) 
                    && \text{by dfn.~$c_\zeta$}\\
        & \ = \ ((b_\zeta)_\diamond \after
                    b_\zeta^\BOX \after \smash{\overline{k}^\BOX} ) (0) 
                    && \text{by dfn.~$\overline{k}$}\\
        & \ = \ \smash{\overline{k}^\BOX}  (0) 
                    && \text{by \sref{diamondboxlemma}}.
\end{alignat*}
    Finally, the exactness of~\eqref{thesnake}
        in~$\ker b$ is shown with a dual argument
        using the right-modularity of~$f$ over~$\ceil{1\after b}^\perp$. \qed
\end{point}
\end{point}
\begin{point}{90}{Remark}%
As~$\op\vN$ is a~$\dagger$-effectus,
    Grandis' Snake Lemma also holds for von Neumann algebras.
This is somewhat puzzling for the following reason.
To study spaces,
    one associates in algebraic geometry
    to each space
    an object
    (or even a diagram of objects) of
    an abelian category representing some algebraic invariant.
The key point is that it is easier to compute and understand
    these algebraic invariants than it is to work with the spaces
    themselves.
One normally thinks of a category of spaces as being quite different
    from a category of algebraic invariants, hence it's surprising
    that~$\op\vN$, which is a rather complicated
    category of (non-commutative) spaces,
    behaves in this way, similar to an abelian category.
\end{point}
\end{point}
\end{parsec}

% vim: ft=tex.latex

\chapter{Outlook}
\begin{parsec}{2290}%
\begin{point}{10}%
For my closing remarks, I would like to speculate
    on possible applications
    and suggest possible
        directions for future research.

We started this thesis by asking the question
    whether the incredibly useful Stinespring dilation theorem
    extends to any ncp-map between von Neumann algebras.
Now that we have seen it does, the next question is obvious:
    do the applications that make the Stinespring dilation so useful
    also extend to the Paschke dilation?
I want to draw attention to one application in particular:
    proving security bounds on quantum protocols.
Here one uses a continuous version
    of the Stinespring dilation theorem \cite{werner2}.
Is there also such a continuous version for Paschke dilation?
\begin{point}{20}%
Next, we have seen that self-dual Hilbert~C$^*$-modules
    over von Neumann algebras
    are very well-behaved compared
    to arbitrary Hilbert C$^*$-modules.
    Frankly, I'm surprised they haven't been studied
        more extensively before:
        the way results about Hilbert spaces
        generalize elegantly to self-dual Hilbert~C$^*$-modules
        over von Neumann algebras seems hard to ignore.
Only one immediate open question remains here:
    does the normal form of~\sref{selfdual-normalish-form}
    extend in some way for to arbitrary von Neumann algebras?
\end{point}
\begin{point}{30}%
As announced, in the second part of this thesis, we did not reach
    our goal of axiomatizing the category of von Neumann algebras
    categorically.
In our attempt, we did find several new concepts
    such as~$\diamond$-adjointness, $\diamond$-positivity,
        purity defined with quotient and comprehension and
        the~$\dagger$ on these pure maps.
Building on this work,
    van de Wetering recently announced~\cite{wetering,weteringeffthe}
    a reconstruction of finite-dimensional quantum theory.
    To discuss it, we need a definition: call an effectus \emph{operational}
    \index{effectus!operational}
        if its scalars are isomorphic to the real interval;
        both the states and the predicates are order separating;
        all predicate effect modules are embeddable in
        finite-dimensional order unit spaces
        \emph{and} each state space is a closed subset
        of the base norm space of all
        unital positive functionals on the corresponding
        order unit space of predicates.
It follows from \cite{wetering}
    that the state space of any operational~$\&$-effectus
    is a spectral convex set in the sense of Alfsen and Shultz
    and that any operational~$\dagger$-effectus
    is equivalent to a subcategory of the category of Euclidean Jordan Algebras
    with positive maps between them in the opposite direction.
Are there infinite-dimensional generalizations of these results?
Are the state spaces of a real~$\&$-effectuses perhaps
    spectral convex sets?
Is any real~$\dagger$-effectus equivalent to a subcategory of, say,
    the category of JBW-algebras with positive maps between them
        in the opposite direction?
\end{point}
\end{point}
\end{parsec}

% vim: ft=tex.latex

\backmatter
\fancyfoot[CE]{}
\fancyfoot[CO]{}
\fancypagestyle{plain}{
    \fancyfoot[CE]{}
    \fancyfoot[CO]{}
}

\printindex

\begingroup
\renewcommand\chapter[2]{\backmattertitle{Bibliography}}
\bibliography{main}{}
\endgroup

\bibliographystyle{plain}

\oldchapter{Lekensamenvatting}
    
De volgende twee technische vraagstukken
    in de wiskundige theorie van quantum computers
    liggen ten grondslag aan dit proefschrift.
\begin{enumerate}
\item 
    Breidt de stelling van Stinespring uit naar
        alle von Neumann algebra's?
\item
    Is de categorie van von Neumann algebra's te axiomatiseren?
\end{enumerate}
Het eerste vraagstuk wordt beantwoord: ja, de stelling breidt uit.
Het tweede vraagstuk blijft helaas voor een groot deel onbeantwoord.
De belangrijkste bijdragen van dit proefschrift aan de wetenschap
    zijn niet per se deze directe resultaten,
    maar de bijvangst die het onderzoek erna heeft opgeleverd:
    er worden verscheidene nieuwe stellingen bewezen
        en begrippen ingevoerd
    die toepassing hebben buiten deze vraagstukken
    (zoals bijvoorbeeld een uitbreiding van de stelling
        van Kaplansky en het begrip~$\diamond$-geadjungeerdheid).

Maar waar gaan deze twee originele vraagstukken eigenlijk over?
Wat is die stelling van Stinespring?
Wat zijn von Neumann algebra's?
En wat is een quantum computer \"uberhaupt?

\oldsection*{Quantum computers}
Een quantum computer gebruikt
       bepaalde quantum mechanische eigenaardigheden
    om sommige berekeningen veel effici\"enter uit
    te voeren dan een traditionele computer.
Een zorgelijk voorbeeld is
    dat quantum computer
    erg goed is in het ontbinden van getallen in priemfactoren.
Het grootste deel van de moderne cryptografie
        die ons beveiligt is gestoeld op de aanname
        dat het vinden van zulke ontbindingen lastig is
        en deze wordt dan ook makkelijk door
        een voldoende grote quantum computer gebroken.
Gelukkig zijn de huidige quantum computers nog erg klein
    en hebben wij de tijd om onze cryptografie aan te passen.
Aan de positieve kant zijn
    quantum computer erg geschikt voor het simuleren
    van het gedrag van moleculen,
    waaronder het lokale effect van potenti\"ele medicijnen.

In functie is een quantum computer te vergelijken met
        een grafische kaart van een computer:
    voor de meeste berekeningen is een normale CPU het best,
            maar voor sommige berekeningen
        is een grafische kaart veel effici\"enter.
Dit komt doordat een normale CPU gemaakt is om zo snel mogelijk
    losse rekenstappen op afzonderlijke getallen
    achter elkaar uit te voeren,
    terwijl een grafische kaart rekenstappen
    uitvoert op miljoenen getallen tegelijkertijd.
Niet elke berekening heeft baat bij dat grote parallelisme.
Voor diegene die dat wel hebben (zoals het
    weergeven van een 3D-wereld in een computerspel)
    is het lastig om het originele algoritme aan te passen
    om optimaal gebruik te maken van de grafische kaart.
    Programmeren voor een grafische kaart vraagt om een
    heel andere \emph{mindset}.
Ook de theoretische analyse van algoritmes voor grafische kaarten
    is van een andere aard dan die voor een traditionele computer.

De situatie bij een quantum computer is vergelijkbaar, maar extremer:
    de effici\"entiewinst tussen een quantum computer en PC
    is groter dan die tussen een CPU en grafische kaart.
Helaas is het vinden en analyseren van algoritmes voor quantum computers
    ook vele malen ingewikkelder.
Er is ook geen methode of zelfs vuistregel bekend
    om te bepalen of een berekening effici\"enter uit te voeren
    is op een quantum computer of niet.

\oldsection*{Von Neumann algebra's}

Het gedrag van een algoritme voor een traditionele computer,
    zoals diegene die het kwadraat van een telgetal berekent,
    wordt vaak beschreven met wiskundige functie,
    zoals in dit geval de functie~$f\colon \N \to \N$
        met~$f(n) = n^2$.
De notatie~$f\colon \N \to \N$
    betekent dat de functie~$f$ als invoer
    een telgetal~$\N \equiv \{0, 1, 2, \ldots \}$
    neemt en ook een telgetal als uitvoer geeft.
Een echt programma draait op een computer
    dat maar een beperkte hoeveelheid geheugen heeft:
        er is dus een grootste telgetal dat in de computer past.
    Dat we het algoritme beschrijven
    met willekeurig grote telgetallen is een idealisering van de situatie.
Het zou echter onhandig zijn om alle generiek toepasbare algoritme
    te beschrijven als wiskundige functies tussen
    eindige verzamelingen~(zoals~$\{0, 1, 2, 3\}$)
    in plaats van tussen oneindige verzamelingen zoals~$\N$.

Dit is helaas precies de situatie voor de meeste beschrijvingen van
    algoritmes voor quantum computers.
Het probleem is dat het analoog van~$\N$
    en de wiskundige functies daartussen
    vele malen complexer van aard zijn
    dan de wiskundige functies die het gedrag
    beschrijven voor het eindige geval.
Ook al dit eindige geval op zichzelf is
    stukken ingewikkelder dan de simpele functies
    tussen eindige verzamelingen voor de normale algoritmes.

\emph{Von Neumann algebra's} zijn de wiskundige structuren
    die gebruikt worden om oneindige quantumdatatypes
    te beschrijven zoals het analoog van~$\N$.
Zogenoemde~\emph{ncp-functies} tussen von Neumann algebras
    zijn het analoog van de normale wiskundige functies
    tussen verzamelingen.
Het zijn deze von Neumann algebra's met bijbehorende
    ncp-functies, die in dit proefschrift onderzocht worden
    om ze beter te begrijpen en eenvoudiger toe te kunnen passen.

\oldsection*{De stelling van Stinespring}
Het is lastig om grip te krijgen
    op een willekeurige ncp-functie
    tussen von Neumann algebra's.
De stelling van Stinespring helpt hierbij:
    deze zegt dat elke ncp-functie~$f\colon \scrA \to \scrB_I$
    tussen von Neumann algebra's~$\scrA$ en~$\scrB_I$
        (waar~$\scrB_I$ van een speciale soort is)
        in feite de samenstelling is
        van twee veel simpelere ncp-functies:
        eerst~$f_1\colon \scrA \to \scrP$
            en dan~$f_2\colon \scrP \to \scrB_I$.
    Deze bijzondere opsplitsing wordt
        ook wel een \emph{dilatie} genoemd
            en wordt gebruikt als cruciale stap in het
    bewijs van vele stellingen, niet alleen over ncp-functies zelf,
    maar ook over bijvoorneeld quantum protocollen en quantum informatie.
Het probleem is dat de stelling van Stinespring alleen toepasbaar
    is als~$\scrB_I$ van de speciale soort \emph{type I} is.
Daarmee zijn de stellingen die ermee worden bewezen vaak ook maar toepasbaar
    op zulke speciale von Neumann algebra's.
In het eerste hoofdstuk van het proefschrift
    wordt een generalisatie van de stelling van Stinespring
    bewezen, die toepasbaar is op alle ncp-functies.

\oldsection*{Axiomatisatie}
Waarom zijn juist von Neumann algebra's en ncp-functies
    geschikt voor oneindige quantum datatypes?
In feite is het een wiskundige gok:
    von Neumann algebra's en ncp-functies
    lijken de eenvoudigste structuren die
    de voorbeelden die we kennen goed beschrijven.
Eenvoud hier is relatief: de ingewikkelde definitie
    van een von Neumann algebra
    staat ver van de intu\"itie van een
    informaticus en natuurkundige.
Het valt ook te betwijfelen of elke von Neumann
    algebra een realistisch quantum mechanisch systeem beschrijft,
    laat staan een quantum datatype.

Een mogelijke manier om deze kwestie te beslechten is de volgende:
    wij zoeken een lijst van kenmerkende
    (en voor onze toepassing relevante) eigenschappen waaraan
        von Neumann algebra's samen met ncp-functies aan voldoen.
Als er bewezen kan worden dat alleen von Neumann algebra's
    met ncp-functies aan deze lijst van eigenschappen kan voldoen
    (en geen ander soort wiskundige structuur)
    dan hebben wij met deze lijst van eigenschappen
        (nu verheven tot axioma's)
    een nieuw inzicht in de aard van von Neumann algebra's.

Een dergelijke aanpak is eerder al succesvol gebruikt
    om de regels van quantum theorie voor eindige systemen beter te begrijpen.
Een voorbeeld hiervan zijn de axioma's van Chiribella et al,
    die quantum theorie karakteriseren door
    de wijze waarop informatie verwerkt wordt.

In het tweede hoofdstuk van dit proefschrift
    wordt gezocht naar een dergelijke axiomatisatie.
        Een van de kernaxioma's is het bestaan
            van een soort \emph{omkering}, een zogenoemde~\emph{dagger},
            op de speciale ncp-functies die aan de rechterkant
            voorkomen van een dilatie,
            wat de twee hoofdstuken met elkaar verbindt.
Helaas wordt een volledige axiomatisatie niet bereikt.
Het verhaal eindigt hier niet:
    mijn college van de Wetering
        heeft recent bewezen dat met een paar toevoegingen,
        de axioma's uit dit proefschrift
        zogeheten Euclidische Jordan algebra's (EJA's) karakteriseren.
Hiermee is de kous nog niet af:
    EJA's zijn niet geschikt voor oneindige quantum datatypes.
Dit resultaat suggereert wel de mogelijkheid
    dat JBW-algebra's,
        welke een kruising zijn van von Neumann en Jordan algebra's,
    uiteindelijk de juiste algebra's zijn voor onze toepassing.

Hoewel ons tweede vraagstuk niet beantwoord is, heeft de zoektocht
    naar geschikte axioma's een aantal nieuwe begrippen aan het licht
    gebracht (zoals~$\diamond$-geadjun\-geerdheid, hoeken en filters), die een
    nieuw inzicht verschaffen in von Neumann algebra's
        en daarmee quantum computers in het algemeen.

\oldchapter{About the author}

Bas Westerbaan, born in 1988,
    enrolled to the Radboud Universiteit
     as a physics student in 2007.
He received his bachelor's degree in
    mathematics with a thesis on computability theory
    supervised by~dr.~Wim Veldman.
He continued his study of logic and foundational computer science
    during his master's eduction and developed a keen interest for,
     what he thought to be, an unrelated mathematical field:
        functional analysis.
These interests however, are conveniently combined in the present thesis
    which continues the research of his master's thesis,
    which was supervised by prof.~Bart Jacobs as well.
During his undergraduate studies,
    Bas also served as secretary of the board of the 150 member student club
        Karpe Noktem and he was elected to
        the student board of the faculty of science.
His master's degree in mathematics was awarded cum laude in 2013.

While perfoming the research as a \emph{promovendus}
    at the Digital Security group that led to this thesis,
    Bas coauthered several security audits
        including those of the \emph{DigiD-} and \emph{BerichtenBox}-app;
        the Dutch digital identity and mail system. 
In the year between the submission and defense
    of his thesis, Bas was employed as a post-doctoral researcher
    where he assisted with the development of a cryptographic system
    for polymorphic pseudonymisation and encryption.
During that same year he also taught a quarter-year course on
    computer networking
        and was kindly invited by dr.~Heunen
        to speak about his research at the University of Edinburgh.

% vim: ft=tex.latex

\oldchapter{Solutions to exercises}\label{sols}
\begin{solution}{physics-stinespring}%
To prove the first statement,
let~$\varphi \colon \scrB(\scrH) \to \scrB(\scrK)$ be any ncp-map.
If~$\varphi=0$ then~$\scrK'=0$ and~$V=0$ does the job,
        so for the other case, assume~$\varphi \neq 0$.
    By~\sref{stinespring-theorem}
        there is a Hilbert space~$\scrH'$,
        a bounded operator~$W \colon \scrH \to \scrH'$ and an
        nmiu-map~$\varrho \colon \scrB(\scrH) \to \scrB(\scrH')$
        such that~$\varphi = \ad_W \after \varrho$.
Clearly~$\varrho \neq 0$.
    By~\sref{nmiu-between-type-I}
        there is a Hilbert space~$\scrK'$
        and a unitary~$U \colon \scrH' \to \scrH \otimes \scrK'$
        with~$\varrho(A) = U^* (A \otimes 1) U$
         for all~$A \in \scrB(\scrH)$.
Define~$V \equiv UW$.
    Then~$\varphi(A) = W^* \varrho(A) W = W^*U^* (A \otimes 1) UW
        =  V^* (A \otimes 1 ) V$, as desired.

Before we can continue with the second statement,
    we need to understand the relationship
    between quantum channels and ncpu-maps.
This relationship is best understood with
    predual characterization of von Neumann algebras
    due to Sakai~\cite{sakai}, which we have been avoiding.
The characterization is as thus:
    a C$^*$-algebra~$\scrA$ is a von Neumann algebra
    if and only if it is isomorphic to the dual of a Banach space.
Then this Banach space is unique up-to-isomoprhism
    as it must be isomorphic to the space of normal functionals on~$\scrA$
    (denoted by~$\scrA_*$)
    and is appropriately called the \emph{predual}~of~$\scrA$.
Any normal linear map~$\varphi\colon \scrA \to \scrB$
    between von Neumann algebras~$\scrA$ and~$\scrB$
    yields a linear map~$\varphi_*\colon \scrB_* \to \scrA_*$
    via~$\varphi_*(\omega) = \omega \after \varphi$.
In the other direction, any linear map~$\varphi_* \colon \scrB_* \to \scrA_*$
    gives rise to a normal linear map~$\varphi\colon \scrA \to \scrB$
    by defining~$\varphi(a)(\omega) = \varphi_*(\omega)(a)$
    where we identified~$\scrA \equiv (\scrA_*)^*$.
Clearly~$\varphi$ is positive precisely
    if~$\varphi_*$ maps positive functionals to positive functionals.

A normal state~$\omega\colon \scrB(\scrH) \to \C$
    is precisely of the form~$\omega(a) = \TR[\rho a]$ for some density
        matrix~$\rho$ over~$\scrH$.
Thus the predual of~$\scrB(\scrH)$ can be identified
    with the set of trace-class operators over~$\scrH$.
Let~$\varphi_*$ be a linear map from the density operators
    on~$\scrH$ to those on~$\scrK$.
The map~$\varphi_*$ is completely positive in its usual sense
    if the corresponding map~$\varphi$ is completely positive.
Furthermore~$\varphi$ is unital if and only if~$\varphi_*$
    is trace-preserving.

To prove the second statement,
    let~$\Phi$ be any quantum channel
    mapping density matrices over~$\scrH$ to those of~$\scrH$ again.
(Note: in the printed version of the thesis
    the exercise incorrectly
    assumes~$\Phi$ to map density matrices over~$\scrH$
    to those over some other Hilbert space~$\scrK$.)
It follows from the previous, that there is a unique
        ncpu-map~$\varphi\colon \scrB(\scrH) \to \scrB(\scrH)$ with
\begin{equation*}
    \TR[ \Phi(\rho) A] \ = \ \TR[\rho \varphi(A)]
        \quad\text{for any density matrix~$\rho\in \scrB(\scrH)$}.
\end{equation*}
See also~\cite{tomamichel} for a more direct approach.
By the first part if the exercise,
    we know that there is a Hilbert space~$\scrK'$
    and a bounded operator~$V \colon \scrH \to \scrH \otimes \scrK'$
    with~$\varphi (A) = V^* (A \otimes 1) V$.
Tracing back the definition of~$V$, we see that~$V$
    is an isometry because~$\varphi$ is unital.
Pick any orthonormal bases~$E$ and~$F$ of~$\scrH$ and~$\scrK'$
    respectively.
We may assume, without loss of generality,
    that~$\scrK'$ is not zero-dimensional
    by setting~$\scrK' =\C$ and~$V = 0$
    in the case that~$\Phi = 0$.
Pick any~$f_0 \in F$
    and any unitary~$U\colon \scrH \otimes \scrK' \to \scrH \otimes \scrK'$
    with~$U^* x \otimes f_0 = V x$,
    which exists as~$V$ is an isometry.
Now we compute
\begin{align*}
    \TR [\varphi(A) \rho ]
    & \ = \ \TR[\rho V^* (A \otimes 1) V] \\
    & \ = \ \sum_{e \in E} \langle V \rho e, (A \otimes 1) V e \rangle \\
    & \ = \ \sum_{e \in E}
    \bigl\langle U^* (\rho \otimes 1) \,e\otimes f_0 ,
    \ (A \otimes 1) U^* \,e\otimes f_0 \bigr\rangle. \\
\intertext{Inserting~$\ketbra{f_0}{f_0}$
    in the previous,
    we may sum over all~$f \in F$ and get}
     \TR[\varphi(A) \rho] &\ = \ \sum_{\substack{e \in E \\ f\in F}}
    \bigl\langle U^* (\rho \otimes \ketbra{f_0}{f_0}) \,e\otimes f ,
    \ (A \otimes 1) U^* \,e\otimes f \bigr\rangle \\
    & \ = \ \sum_{\substack{e \in E \\ f\in F}}
    \bigl\langle
    (U^* (e\otimes f)) , \ 
    U^* (\rho \otimes \ketbra{f_0}{f_0})U \,
    (A \otimes 1) \,(U (e\otimes f)) \bigr\rangle \\
    & \ = \ 
    \TR \bigl[
    U^* (\rho \otimes \ketbra{f_0}{f_0})U \,
    (A \otimes 1)  \bigr] \\
& \ = \ 
\TR \bigl[ A   \TR\nolimits_{\scrK'}[ U^* (\rho \otimes \ketbra{f_0}{f_0} ) U] \bigr].
\end{align*}
This show that indeed~$\Phi(\rho)
    = \TR_{\scrK'}[ U^* (\rho \otimes \ketbra{v_0}{v_0})U]$
    as desired with~$v_0 \equiv f_0$.
\end{solution}
\begin{solution}{kraus-exercise}%
Let~$\varphi\colon \scrB(\scrH) \to \scrB(\scrK)$
    be any ncp-map.
By~\sref{physics-stinespring}
    there is a Hilbert space~$\scrK'$
    and a bounded operator~$V\colon \scrK \to \scrH \otimes \scrK'$
    with~$\varphi(A) = V^* (A \otimes 1) V$.
Let~$E$ be any orthonormal basis of~$\scrK'$.
    Then~$1 = \sum_{e\in E} \ketbra{e}{e}$
        where the sum converges ultraweakly
        and so by ultraweak continuity of~$\ad_V$ (\sref{ad-normal})
        and~$B \mapsto A\otimes B$ (\sref{tensor-simple-facts}), we see
    \begin{equation}\label{kraus-exc-eq1}
        \varphi(A) \ =\  V^* \Bigl(A \otimes \sum_{e\in E} \ketbra{e}{e}\Bigr) V
        ) \ =\  \sum_{e \in E} V^* (A \otimes \ketbra{e}{e}) V.
    \end{equation}
For~$e\in E$, define~$P_e \colon \scrH \otimes \scrK' \to \scrH$
    by~$P_e \equiv 1\otimes \bra{e}$,
    i.e.~$P_e(x \otimes y) = x \langle e, y\rangle$.
Define~$V_e \equiv P_e V$.
    Note that~$P_e^*AP_e = A \otimes\ketbra{e}{e}$ and so
\begin{alignat*}{2}
    \varphi(A)
    &\ = \  \sum_{e \in E} V^* (A \otimes \ketbra{e}{e} ) V &\qquad&
    \text{by \eqref{kraus-exc-eq1}} \\
    &\ = \  \sum_{e \in E} V^* P_e^*A P_e V \\
    &\ = \  \sum_{e \in E} V_e^* A V_e,
\end{alignat*}
as desired.
From the special case~$A=1$, we see
    that~$\sum_{e \in E} V_e^*V_e = \varphi(1)$
    and so the partial sums of~$\sum_{e \in E} V_e^* V_e$ are bounded.

For the final part, assume~$\scrH$ and~$\scrK$ are finite dimensional.
Recall that the standard Stinespring dilation space (say~$\scrK''$)
    for~$\varphi$ is constructed using a completion
    and quotient of~$\scrB(\scrH)\odot \scrK$.
As~$\scrB(\scrH)\odot \scrK$ is finite dimensional
    it is already complete.
    Hence~$\scrK''$ has dimension at most~$(\dim\scrH)^2( \dim\scrK)$.
By construction~$\scrH \otimes \scrK' \cong \scrK''$,
    hence~$\dim \scrK' \leq (\dim \scrH )(\dim \scrK)$.
Recall~$E$ is a basis of~$\scrK'$
    and so there are indeed at most~$(\dim \scrH )(\dim \scrK)$
        Kraus operators.
\end{solution}
\begin{solution}{exc-chris-univ-prop}%
We will show that~$U\colon \mathsf{Rep} \to \mathsf{Rep}_{\mathrm{cp}}$
    has a left adjoint by demonstrating the universal mapping property.
    Let~$\varphi\colon \scrA \to \scrB(\scrH)$ be any object of~$\mathsf{Rep}_{\mathrm{cp}}$.
    Pick any minimal Stinespring dilation~$(\scrK, \varrho, V)$ of~$\varphi$.
The map~$\varrho\colon \scrA\to \scrB(\scrK)$ is an object
        of~$\mathsf{Rep}$.
Clearly~$\ad_V \after \varrho \after \id = \varphi$
    and so~$\eta_\varphi\equiv (\id,V)\colon \varphi \to U\varrho$
        is a morphism in~$\mathsf{Rep}_{\mathrm{cp}}$.
We will show that for each~$f\colon \varphi \to U\varrho'$
        in~$\mathsf{Rep}_{\mathrm{cp}}$,
    there is a unique~$f'\colon \varrho \to \varrho'$
    in~$\mathsf{Rep}$ with~$Uf' \after \eta_\varphi = f$.
This is sufficient to show that~$U$
    has a left adjoint.

    So let~$f\colon \varphi\to U\varrho'$ be any morphism
        in~$\mathsf{Rep}_{\mathrm{cp}}$.
Say~$\varrho'\colon \scrA' \to \scrB(\scrK')$.
    Then~$f \equiv (m', V')$ consists of
    a nmiu-map~$m'\colon \scrA \to \scrA'$
    and bounded operator~$V' \colon \scrH \to \scrK'$
    with~$\ad_{V'} \after \varrho' \after m' = \varphi$.
By~\sref{dils-univ-stinespring}
    there is a unique bounded operator~$S\colon \scrK \to \scrK'$
        with~$SV = V'$ and~$\varrho = \ad_S \after \varrho' \after m'$.
This turns~$f' \equiv(m',S)$ into a
    morphism~$\varrho \to \varrho'$ in~$\mathsf{Rep}$.
Furthermore~$Uf' \after \eta_\varphi
                = (m' \after \id, SV) = (m',V') = f$.
To show uniqueness, assume
        there is some~$f'' \colon \varrho \to \varrho'$
        in~$\mathsf{Rep}$
        with~$Uf'' \after \eta_\varphi = f$.
Say~$f'' = (m'',S'')$.
    Then~$(m',V') = f = Uf' \after \eta_\varphi = (m'', S''V)$.
So~$m''=m'$ and~$V' = S''V$.
The fact that~$f''$ is a morphism in~$\mathsf{Rep}$
    is
    equivalent to~$\ad_{S''} \after \varrho' \after m'' = \varrho$.
Thus~$\ad_{S''} \after \varrho' \after m' = \varrho$.
By uniqueness of~$S$, we get~$S'' = S$.
    Hence~$f''=(m'',S'') = (m',S) = f'$, as desired.
\end{solution}
\begin{solution}{ess-uniq-pur}%
Let~$\varphi\colon \scrB(\scrH) \to \scrB(\scrK)$ be any ncp-map.
As in the description of the exercise,
    let~$\scrK$ be a Hilbert space
    and~$V,W\colon \scrK \to \scrH \otimes \scrK'$
    be bounded operators
    with~$V^* (a \otimes 1) V = \varphi(a) = W^* (a\otimes 1) W$.
Write~$\scrV$ for the closed linear span
    of~$\{(a \otimes 1) V x; \ a \in \scrB(\scrH),\ x \in \scrK\}$
        in~$\scrH\otimes \scrK'$
and similarly~$\scrW$ for that
    of~$\{(a \otimes 1) W x; \ a \in \scrB(\scrH),\ x \in \scrK\}$.
Note that for any~$n\in \N$, ~$x_1,\ldots, x_n \in \scrK$
    and~$a_1, \ldots, a_n \in \scrB(\scrH)$ we have
\begin{align*}
    \bigl\| \sum_i (a_i\otimes1) V x_i \bigr\|^2
    &\ = \ 
     \sum_{i,j} \langle x_i,\, V^* ((a_i^*a_j) \otimes 1) V x_j\rangle \\
    &\ = \ 
     \sum_{i,j} \langle x_i,\, W^* ((a_i^*a_j) \otimes 1) W x_j\rangle \\
     &\ = \ 
    \bigl\| \sum_i (a_i\otimes1) W x_i \bigr\|^2.
\end{align*}
Thus there is a unique unitary~$U_0\colon \scrW \to \scrV$
    fixed by~$U_0 (a \otimes 1) W x = U_0 (a \otimes 1) V x$.
We see~$U_0 W = V$ by setting~$a=1$.
Furthermore
    \begin{equation*}
        (\alpha \otimes 1) U_0 (a \otimes 1) W x
        \ = \ ((\alpha a)  \otimes 1) V x
        \ = \ U_0 (\alpha  \otimes 1 )(a  \otimes 1) W x
    \end{equation*}
    for any~$\alpha,a \in \scrB(\scrH)$ and~$x \in \scrK$,
    hence~$(\alpha \otimes 1) U_0 = U_0 (\alpha \otimes 1)$.

For any~$a \in \scrB(\scrH)$,
    the operator~$a \otimes 1 \in \scrB(\scrH \otimes \scrK')$
    restricts to~$\scrB(\scrW)$.
Pick an orthonormal basis~$E$ of~$\scrH$
    and some~$e_0 \in E$.
    Note that~$(\ketbra{e_0}{e_0} \otimes 1) \scrW = e_0 \otimes \scrW'$
    for some closed subspace~$\scrW' \subseteq \scrK'$.
In fact, for any~$e \in E$
    we have~$e \otimes \scrW'
    = (\ketbra{e}{e_0} \otimes 1) (e_0 \otimes \scrW')
    = (\ketbra{e}{e_0}\otimes 1)  (1 \otimes \ketbra{e_0}{e_0}) \scrW
    = (\ketbra{e}{e} T \otimes 1)  \scrW = (\ketbra{e}{e} \otimes 1) \scrW$,
    where~$T$ is the unitary on~$\scrH$ that only swaps~$e$ and~$e_0$.
Hence~$\scrW = \scrH \otimes \scrW'$.
Similarly~$\scrV = \scrH \otimes \scrV'$
    for some closed subspace~$\scrV' \subseteq \scrK'$.

For any non-zero~$w \in \scrW'$ and unit-vector~$x \in \scrH$,
    we have~$U_0 (x \otimes w)
        = U_0 (\ketbra{x}{x} \otimes 1)( x\otimes w)
        = (\ketbra{x}{x} \otimes 1) U_0 (x\otimes w)$
        so~$U_0 (x \otimes w) = x \otimes y$ for some~$y \in \scrV'$.
Clearly~$\| w \| = \| x \otimes w\|=\| U_0 (x\otimes w) \|
        = \|x \otimes y\| = \|y\|$,
        so there is a unique unitary~$U_1\colon \scrW' \to \scrV'$
        with~$U_0 (x \otimes w) = x \otimes U_1 w$.
        It follows~$U_0 = 1 \otimes U_1$.

As~$\scrV$ and~$\scrW$ are isomorphic, they have the same dimension
    and so do~$\scrV^\perp$ and~$\scrW$.
Consequently, there is an unitary~$U\colon \scrK' \to \scrK'$
    extending~$U_1$.
We have~$V = (1\otimes U) W
    =   (1 \otimes U_1) W
     = U_0 W = V$ as desired.
\end{solution}
\begin{solution}{paschke-basics}%
We cover the points in order.
\begin{enumerate}
\item
Let~$\varrho \colon \scrA \to \scrB$ be a mniu-map.
Assume there are nmiu~$\varrho'\colon \scrA \to \scrP'$
    and ncp~$h'\colon \scrP' \to \scrB$
        with~$h' \after \varrho' = \varrho$.
We have to show there is a unique map~$\sigma\colon \scrP \to \scrB$
    with~$\id \after \sigma= h'$ and~$h' \after \varrho' = \varrho$.
        Clearly~$\sigma\equiv h' $ fits the bill.

\item
Let~$(\scrP,\varrho,h)$ be any Paschke dilation.
        We will show that~$(\scrP, \id, h)$ is a Paschke dilation of~$h$.
To this end, let~$\varrho'\colon \scrA \to \scrP'$
        be any nmiu-map and~$h'\colon \scrP' \to \scrB$
        be any ncp-map with~$h' \after \varrho' = h$.
Consider~$\varrho' \after \varrho$ and~$h'$.
    By the universal property of the original dilation,
            there is a unique ncp-map~$\sigma \colon \scrP' \to \scrP$
                with~$\sigma\after\varrho'\after\varrho = \varrho$
                and~$h \after\sigma = h'$.
Furthermore~$\id\colon \scrP \to \scrP$
    is the unique ncp-map
        with~$\id \after \varrho = \varrho$
        and~$h \after \id = h$.
Now~$(\sigma \after \varrho') \after \varrho = \varrho$
    and~$h \after (\sigma \after \varrho') = h' \after \varrho' = h$,
        so~$\sigma \after \varrho' = \id$.
We are are halfway demonstrating that~$\sigma$ is also the mediating map
        for our dilation of~$h$.
It remains to be shown that~$\sigma$ is the unique ncp-map
    with~$\sigma \after \varrho' = \id$
        and~$h \after \sigma = h'$.
So assume there is a ncp-map~$\sigma'\colon \scrP' \to \scrP$
    with~$h \after \sigma' = h'$ and~$\sigma' \after \varrho' = \id$.
    Clearly~$\sigma' \after\varrho' \after\varrho =\varrho$
        and so by uniqueness of~$\sigma$ as the mediating map
        for the orignal dilation,
        we see~$\sigma' = \sigma$, as desired.
    \item

    Let~$
        \left(\begin{smallmatrix}\varphi_1\\\varphi_2 \end{smallmatrix} \right)
        \colon \scrA \to \scrB_1 \oplus \scrB_2$
        be any ncp-map.
    Pick Paschke dilations~$(\scrP_i, \varrho_i, h_i)$
            of~$\varphi_i$ for~$i=1,2$.
        We will show that~$(\scrP_1 \oplus \scrP_2, 
        \left(\begin{smallmatrix}\varrho_1\\\varrho_2 \end{smallmatrix} \right),
            h_1 \oplus h_2 )$
            is a Paschke dilation of~$
    \left(\begin{smallmatrix}\varphi_1\\\varphi_2 \end{smallmatrix} \right)$.
Clearly~$h_1 \oplus h_2 \after 
        \left(\begin{smallmatrix}\varrho_1\\\varrho_2 \end{smallmatrix} \right)=
        (\begin{smallmatrix}h_1 \after\varrho_1\\h_2 \after \varrho_2 \end{smallmatrix} )
            = 
    \left(\begin{smallmatrix}\varphi_1\\\varphi_2 \end{smallmatrix} \right) $.
Let~$\varrho'\colon \scrA \to \scrP'$ be any nmiu-map
    and~$
    \left(\begin{smallmatrix}h_1\\ h_2\end{smallmatrix} \right) 
        \colon \scrP' \to \scrB_1 \oplus \scrB_2$
        any ncp-map with~$
    \left(\begin{smallmatrix}h_1\\ h_2\end{smallmatrix} \right)  \after
        \varrho' = 
    \left(\begin{smallmatrix}\varphi_1\\ \varphi_2\end{smallmatrix} \right) $.
Note~$h_i \after \varrho' = \varphi_i$
    (for~$i=1,2$)
    and so there is a unique~$\sigma_i\colon \scrP' \to \scrP_i$
    with~$\sigma_i \after \varrho' = \varrho_i$
    and~$h_i \after \sigma_i = h_i'$.
We will show that~$
    \left(\begin{smallmatrix}\sigma_1\\ \sigma_2\end{smallmatrix} \right) 
        \colon \scrP' \to \scrP_1\oplus \scrP_2$
        is the unique mediating map.
Clearly~$
    \left(\begin{smallmatrix}\sigma_1\\ \sigma_2\end{smallmatrix} \right) 
        \after \varrho' = 
    \bigl(\begin{smallmatrix}\sigma_1 \after \varrho'\\ \sigma_2 \after \varrho' \end{smallmatrix} \bigr) =
    \left(\begin{smallmatrix}\varrho_1\\ \varrho_2 \end{smallmatrix} \right) $
and~$(h_1 \oplus h_2) \after
\bigl( \begin{smallmatrix} \sigma_1\\ \sigma_2 \end{smallmatrix} \bigr)
    =
\bigl( \begin{smallmatrix}
h'_1 \after \sigma_1\\
h'_2 \after \sigma_2
\end{smallmatrix} \bigr) =
\bigl( \begin{smallmatrix}
h_1 \\
h_2
\end{smallmatrix} \bigr)$.
To show uniqueness of~$
( \begin{smallmatrix}
\sigma_1\\
\sigma_2
\end{smallmatrix})$,
assume
there is ncp-map$
    \bigl(\begin{smallmatrix}\sigma'_1\\ \sigma'_2\end{smallmatrix} \bigr) 
        \colon \scrP' \to \scrP_1\oplus \scrP_2$
    such that~$ \bigl(\begin{smallmatrix}\sigma'_1\\ \sigma'_2\end{smallmatrix} \bigr) 
        \after \varrho'
        = 
    \bigl(\begin{smallmatrix}\varrho_1\\ \varrho_2 \end{smallmatrix} \bigr) $
and~$(h_1 \oplus h_2) \after
\bigl( \begin{smallmatrix} \sigma'_1\\ \sigma'_2 \end{smallmatrix} \bigr)
    =
\bigl( \begin{smallmatrix}
h_1 \\
h_2
\end{smallmatrix} \bigr)$.
Then~$h_i \after \sigma_i' = h_i$ and~$\sigma_i' \after \varrho' = \varrho_i$
    for~$i=1,2$ and so~$\sigma_i=\sigma_i'$ by the uniqueness
    of the seperate~$\sigma_i$.
    Thus indeed~$
( \begin{smallmatrix}
\sigma_1\\
\sigma_2
\end{smallmatrix}) =
\bigl( \begin{smallmatrix}
\sigma'_1\\
\sigma'_2
\end{smallmatrix}\bigr)$.

\item
Let~$\varphi\colon \scrA \to \scrB$ be any ncp-map
    with Paschke dilation~$(\scrP, \varrho, h)$.
Assume~$\lambda \in \R, \lambda > 0$.
We will show~$(\scrP, \varrho, \lambda h)$
    is a Paschke dilation of~$\lambda \varphi$.
Clearly~$\lambda h \after \varrho = \lambda \varphi$.
To this end, assume~$\varrho'\colon \scrA \to \scrP'$ is a nmiu-map
    and~$h' \colon \scrP' \to \scrB$ is an ncp-map
    with~$h' \after \varrho' = \lambda \varphi$.
Then~$\lambda^{-1} h' \after \varrho' = \varphi $.
Thus there is a unique~$\sigma\colon \scrP' \to \scrP$
    with~$\sigma \after \varrho' = \varrho$
    and~$h \after \sigma = \lambda^{-1} h'$.
Clearly~$\lambda h \after \sigma = h'$ and so~$\sigma$
    also serves as the unique mediating map
    for the dilation of~$\lambda \varphi$.
\end{enumerate}
\end{solution}
\spacingfix
\begin{solution}{module-seminorm}%
Let~$X$ be a right~$\scrB$-module
with~$\scrB$-valued inner product~$\langle \,\cdot\,,\,\cdot\,\rangle$
    for some C$^*$-algebra~$\scrB$.
Using the C$^*$-identity, ~\sref{module-CS} and the definition
    of~$\|\,\cdot\,\|$ on~$X$, we
    get~$ \| \langle x, y\rangle\|^2
        = \| \langle x, y\rangle^* \langle x, y\rangle\|
        = \| \langle y,x \rangle\langle x,y\rangle \|
        \leq \bigl\| \|\langle x,x\rangle\| \langle y, y \rangle \bigr\|
        = \| x\|^2 \|y\|^2$,
        so indeed~$\|\langle x, y \rangle \| \leq \|x\|\|y\|$.

    Now we will show~$\|x\| \equiv \| \langle x,x\rangle\|^{\frac{1}{2}}$
    is a seminorm on~$X$.
    Clearly~$\|x\| \geq 0$ for any~$x \in X$.
For any~$\lambda \in \C$ and~$x \in X$
    we have~$\langle \lambda x, \lambda x \rangle
        = \overline\lambda \langle x, x \rangle \lambda
        = | \lambda |^2 \langle x,x\rangle$,
    hence~$
        \|\lambda x\|=
       \|\lambda^2 \langle x, x \rangle \|^{\frac{1}{2}} =
        |\lambda|  \|x\| $.
Next, for any~$x,y \in X$ we have
\begin{alignat*}{2}
    \| x + y \|^2 & \ = \  \| \langle x + y, x+y\rangle \| \\
    & \ \leq \ \| \langle x,x\rangle\| 
        + \| \langle y,y\rangle\|
        + \| \langle x, y \rangle \|
        + \| \langle x, y \rangle^* \| \\
    & \ = \ \|x\|^2 + \|y\|^2 + 2\|\langle x,y\rangle\| \\
    & \ \leq \ \|x\|^2 + \|y\|^2 + 2\|x\|\|y\| \\
    & \ =\ (\|x\| + \|y\|)^2
\end{alignat*}
and thus~$\|\,\cdot\,\|$ is indeed a seminorm.

Finally, we will show~$\|x \cdot b\| \leq \|x\|\|b\|$
    for any~$b\in \scrB$ and~$x \in X$.
As~$\langle x,x\rangle$ is positive,
    we have~$\langle x,x \rangle  \leq \|\langle x,x\rangle \|
            = \|x\|^2$.
Also recall~$a \mapsto b^* a b$ is clearly positive.
Thus~$\|x\cdot b\|^2 \equiv \|\langle xb, xb\rangle \|
        = \| b^* \langle x, x \rangle b\|
        \leq  \bigl\|
        \|x\|^2
        b^*b
        \bigr\|
        = \|b\|^2 \|x\|^2$ as desired.
\end{solution}
\begin{solution}{hilbmod-polarization}%
Let~$B$ be any~$\scrB$-sesquilinear form on a pre-Hilbert~$\scrB$-module~$X$
    for some C$^*$-algebra~$\scrB$.
    Distributing~$B$ in~$\sum_{k=0}^3 i^k B(i^k x+y, i^k x+y)$
    we get 16 terms
    consisting of four of each~$B(x,x)$, $B(y,y)$, $B(x,y)$ and $B(y,x)$
    with the following coefficients.
    \begin{center}
        \begin{tabular}{c|rrrr}
        & $B(x,x)$ & $B(y,y)$ & $B(x,y)$ & $B(y,x)$\\ \hline
        $k\,=\,0$ & $1$ & $1$ & $1$ & $1$ \\
        $k\,=\,1$ & $i\cdot (-i)\cdot i\,=\,i$ & $i$ & $i \cdot (-i)\,=\,1$ & $i^2 \,=\, -1$ \\
        $k\,=\,2$ & $(-1)^3\,=\,-1 $&$ -1 $&$ (-1)^2 \,=\,1$&$ (-1)^2\,=\,1$ \\
            $k\,=\,3$ & $(-i) \cdot i \cdot (-i)\,=\,-i $&$ -i $&$ (-i)\cdot i\,=\, 1$&$ (-i)^2\,=\,-1 $\\
    \end{tabular}
\end{center}
    Note that the coefficients in
    every column sum to~$0$, except for the coefficients for~$B(x,y)$ which
    sum to~$4$.
    Hence~$ \sum_{k=0}^3 i^k B(i^k x+y, i^k x+y) = 4B(x,y)$.
\end{solution}
\begin{solution}{exc-subbase}%
    Let~$X$ be a set together with a subbase~$B$.
    Write~$\Phi$ for the filter generated by~$B$.
    Note~$B \subseteq \Phi$.
    We will show~$(X,\Phi)$ is a uniform space, by proving
        its axioms in order.
\begin{enumerate}
    \item
        By construction~$\Phi$ is a filter.
    \item
        Pick any~$\varepsilon \in \Phi$.
        We have to show~$\Delta\equiv\{(x,x); \ x\in X\} \subseteq \varepsilon$.
        By definition of~$\Phi$, there
        are~$\varepsilon_1, \ldots, \varepsilon_n \in B$
        with~$\varepsilon_1 \cap \ldots \cap \varepsilon_n \subseteq \varepsilon$.
    By definition of a base, we have~$\Delta \subseteq \varepsilon_i$
        for each~$1 \leq i \leq n$
        and so~$\Delta \subseteq \varepsilon_1 \cap \ldots
        \cap \varepsilon_n \subseteq\varepsilon$ as well.
    \item
        Pick any~$\varepsilon \in \Phi$.
    By definition of~$\Phi$, there
        are~$\varepsilon_1, \ldots, \varepsilon_n \in B$
        with~$\varepsilon_1 \cap \ldots \cap
        \varepsilon_n \subseteq \varepsilon$.
    As~$B$ is a base, there are~$\delta_1, \ldots, \delta_n$
     with~$\delta_i \after \delta_i \subseteq \varepsilon_i$
        for~$1 \leq i \leq n$.
    Define~$\delta = \delta_1 \cap \ldots \cap \delta_n$.
        Then~$\delta \after \delta
            \subseteq \bigcap_{i,j} \delta_i \after \delta_j
            \subseteq \bigcap_i \delta_i \subseteq 
            \bigcap_i \varepsilon_i \subseteq \varepsilon$.
\item
        Pick any~$\varepsilon \in \Phi$.
    By definition of~$\Phi$, there
        are~$\varepsilon_1, \ldots, \varepsilon_n \in B$
        with~$\varepsilon_1 \cap \ldots \cap
        \varepsilon_n \subseteq \varepsilon$.
    As~$B$ is a base, there are~$\delta_1, \ldots, \delta_n$
        with~$\delta_i^{-1} \subseteq \varepsilon_i$.
    Define~$\delta = \delta_1 \cap \ldots \cap \delta_n$.
        Then~$\delta^{-1} =
            \bigcap_{i} \delta_i^{-1}  \subseteq
            \bigcap_i \varepsilon_i \subseteq \varepsilon$.
\end{enumerate}
    Thus indeed, $(X, \Phi)$ is a uniform space.
\end{solution}
\begin{solution}{dils-uniform-spaces-basics}%
We prove the subexercises in order.
\begin{enumerate}
\item
First we will show that equivalence of Cauchy nets is an equivalence
    relation.
As for every entourage~$\varepsilon$ we have~$x \mathrel\varepsilon x$,
    we see that every Cauchy net is equivalent to itself.
Assume~$(x_\alpha)_\alpha \sim (y_\beta)_\beta$.
        We will show~$(y_\beta)_\beta \sim (y_\alpha)_\alpha$.
Let~$\varepsilon$ be some entourage.
    There is some entourage~$\delta$ with~$\delta^{-1}\subseteq \varepsilon$.
        By assumption there are~$\alpha_0$ and~$\beta_0$
        such that~$x_\alpha \mathrel{\delta} y_\beta$
        for all~$\alpha \geq \alpha_0$ and~$\beta \geq \beta_0$.
    But then~$y_\beta \mathrel{\varepsilon} x_\alpha$
        for~$\alpha \geq \alpha_0$ and~$\beta \geq\beta_0$.
        Hence~$(y_\beta)_\beta \sim (x_\alpha)_\alpha$.
    To prove transitivity,
assume we are given Cauchy nets~$(x_\alpha)_\alpha \sim (y_\beta)_\beta 
            \sim (z_\gamma)_\gamma$.
Let~$\varepsilon$ be some entourage.
There is an entourage~$\delta$ with~$\delta \after \delta \subseteq \varepsilon$.
There are~$\alpha_0, \beta_0, \gamma_0$
    such that~$x_\alpha \mathrel\delta y_\beta $
            and~$y_\beta \mathrel\delta z_\gamma $
            for~$\alpha \geq \alpha_0$, $\beta \geq \beta_0$
            and~$\gamma \geq \gamma_0$.
        Hence~$x_\alpha \mathrel\varepsilon z_\gamma$ for such~$\alpha$ and~$\gamma$, which shows~$(x_\alpha)_\alpha \sim (z_\gamma)_\gamma$.

Next, assume that~$(x_\alpha)_\alpha$ is a subnet
        of a Cauchy net~$(y_\alpha)_\alpha$.
    To show~$(x_\alpha)_\alpha \sim (y_\alpha)_\alpha$,
            assume~$\varepsilon$ is some entourage.
    By the definition of Cauchy net,
            there is some~$\alpha_0$
            such that~$x_\alpha \mathrel{\varepsilon} x_\beta$
            for all~$\alpha,\beta \geq \alpha_0$.
    In particular~$x_\alpha \mathrel{\varepsilon} y_\beta$
        for all~$\alpha, \beta \geq \alpha_0$
        which shows~$(x_\alpha)_\alpha \sim (x_\beta)_\beta$.

\item
Assume~$(x_\alpha)_\alpha \sim (y_\beta)_\beta$
    and~$x_\alpha \to x$.
To prove $y_\alpha \to x$, pick any entourage~$\varepsilon$.
Pick~$\delta$ such that~$\delta^2 \subseteq \varepsilon$.
There are~$\alpha_0$ and~$\beta_0$
    such that~$y_\beta \mathrel\delta x_\alpha$
        and~$x_\alpha \mathrel\delta x$
        for all~$\alpha \geq \alpha_0 $ and~$\beta \geq \beta_0$.
    Then~$y_\beta \mathrel\varepsilon x$ for all~$\beta \geq \beta_0$,
        whence~$y_\beta \to x$.

\item
    Assume~$(x_\alpha)_\alpha \to x$ and~$(x_\alpha)_\alpha \to y$
            in some Hausdorff uniform space.
    Let~$\varepsilon$ be any entourage.
        Pick~$\delta$ and~$\delta'$ with~$\delta^2 \subseteq \varepsilon$
        and~$\delta'\subseteq \delta^{-1}$.
    There is an~$\alpha_0$ such that~$x_\alpha \mathrel\delta x$
        and~$x_\alpha \mathrel{\delta'} y$
        for all~$\alpha \geq \alpha_0$.
    Thus~$x \mathrel\varepsilon y$.
    As our space is Hausdorff the previous implies~$x=y$.
\item
    Let~$f\colon X\to Y$ be a continuous map between uniform spaces.
    Assume~$x_\alpha \to x$ in~$X$.
    Let~$\varepsilon$ be any entourage of~$Y$.
By continuity there is a~$\delta$
    such that~$x \mathrel\delta y$ implies~$f(x) \mathrel{\varepsilon} f(y)$
    for any~$y \in Y$.
There is an~$\alpha_0$ such that~$x \mathrel{\delta} x_\alpha$
    for all~$\alpha \geq \alpha_0$.
For those~$\alpha$ we also have~$f(x) \mathrel{\varepsilon} f(x_\alpha)$,
    which shows~$f(x_\alpha) \to f(x)$.
\item
    Let~$f\colon X\to Y$ be a uniformly continuous map
        between uniform spaces.
Let~$(x_\alpha)_\alpha$ and~$(y_\beta)_\beta$
    be nets of~$X$
    such that for each entourage~$\varepsilon$ of~$X$
        there are~$\alpha_0,\beta_0$ with~$x_\alpha \mathrel\varepsilon y_\beta$
        for all~$\alpha\geq\alpha_0$ and~$\beta \geq \beta_0$.
The map~$f$ preserves this relation between the nets
    in the following way.
Let~$\varepsilon$ be any entourage of~$Y$
By uniform continuity there is a~$\delta$
    such that~$x \mathrel\delta y$ implies~$f(x) \mathrel\varepsilon f(y)$.
There are~$\alpha_0$ and~$\beta_0$
    with~$x_\alpha \mathrel\delta y_\beta$ for all~$\alpha \geq\alpha_0$
    and~$\beta \geq\beta_0$.
    Hence~$f(x_\alpha) \mathrel\varepsilon f(y_\beta)$
        for such~$\alpha,\beta$.

From the previous it follows that~$f$
    preserves Cauchy nets (by setting~$x_\alpha=y_\alpha$)
    and that it preserves equivalence between Cauchy nets.
\item
Let~$D \subseteq X$ be a dense subset of a uniform space~$X$.
Let~$x \in X$ be any point.
Pick for every~$\varepsilon \in \Phi$
    an element~$d_\varepsilon \in D$
    with~$x \mathrel\varepsilon d_\varepsilon$.
Clearly~$(d_\varepsilon)_{\varepsilon\in\Phi}$ is a net with inverse
inclusion. We have~$d_\varepsilon \to x$
    as~$d_\delta \mathrel\varepsilon x$
    whenever~$\delta \subseteq \varepsilon$.
\item
Assume~$f,g\colon X \to Y$ are continuous maps between uniform
    spaces that agree on a dense subset~$D \subseteq X$.
Let~$x\in X$ be any point.
Pick a net~$x_\alpha$ from~$D$ with~$x_\alpha \to x$.
Then~$f(x) = f(\lim_\alpha x_\alpha) = \lim_\alpha f(x_\alpha)
    =  \lim_\alpha g(x_\alpha) = g(\lim_\alpha x_\alpha) = g(x)$.
    Hence~$f=g$.
\end{enumerate}
\end{solution}
\spacingfix
\begin{solution}{dils-product-uniformity}%
Write~$B \equiv \{ \hat\varepsilon; \ \varepsilon \in \Phi_i;\ i \in I\}$.
First we show~$B$ is a subbase,
    i.e.~that is satisfies conditions 2, 3 and 4
    of~\sref{dils-dfn-uniformity}.
Let~$\hat\varepsilon$ be an arbitrary element of~$B$
    and~$i \in I$ denote the index element with~$\varepsilon \in \Phi_i$.
    Assume~$(x_i)_{i \in I} \in \Pi_i X_i$.
Clearly~$x_i \mathrel\varepsilon x_i$
    and so~$(x_i)_i \mathrel{\hat\varepsilon} (x_i)_i$.
    Thus~$B$ satisfies condition 2 of~\sref{dils-dfn-uniformity}.
Pick a~$\delta \in \Phi_i$ with~$\delta^2 \subseteq \varepsilon$.
    Then~${\hat\delta}^2 =\widehat{\delta^2} \subseteq \hat{\varepsilon}$
    and so~$B$ satisfies condition 3 of~\sref{dils-dfn-uniformity}.
Now let~$\delta$ denote an entourage of~$X_i$
    with~$\delta^{-1} \subseteq \varepsilon$.
    Then~$\hat{\delta}^{-1} = \widehat{\delta^{-1}}
        \subseteq \hat{\varepsilon}$
        and so~$B$ also satisfies condition 4 of~\sref{dils-dfn-uniformity}.

Next we show that the projectors~$\pi_i \colon \prod_i X_i \to X_i$
    are uniformly continuous.
    Assume~$i_0 \in I$ and~$\varepsilon$ is an entourage of~$X_{i_0}$.
Define~$\delta \equiv \hat\varepsilon$.
    Let~$(x_i)_i$ and~$(y_i)_i$
        from~$\prod_i X_i$ be given
        with~$(x_i)_i \mathrel\delta (y_i)_i$.
    Then~$x_{i_0} \mathrel\varepsilon y_{i_0}$.
    Thus~$\pi_{i_0}$ is indeed uniformly continuous.

To show~$(\pi_i)_{i}$ is a categorical product,
    assume we are given a uniform space~$Y$ together
        with uniformly continuous maps~$f_i \colon Y \to X_i$
            for each~$i \in I$.
We have to show that there is a unique uniformly continuous 
    map~$f\colon Y \to \prod_i X_i$
        with~$\pi_i \after f = f_i$ for all~$i \in I$.
    Define~$f$ by~$(f(y))_i \equiv f_i (y)$.
    Clearly~$\pi_i \after f = f_i$.

    To show~$f$ is uniformly continuous,
        pick any entourage~$\varepsilon$ of~$\prod_i X_i$.
        By definition, there are~$i_1, \ldots, i_n$
            and~$\varepsilon_1, \ldots, \varepsilon_n$
            with~$\varepsilon_j \in \Phi_{i_j}$
            and~$\bigcap_j \widehat{\varepsilon_j} \subseteq \varepsilon$.
    For each~$1 \leq j \leq n$
        pick an entourage~$\delta_j$ of~$Y$
        such that~$x \mathrel\delta_j y$
        implies~$f_{i_j}(x) \mathrel{\varepsilon_j} f_{i_j}(y)$.
Define~$\delta \equiv \bigcap_j \delta_j$
    Assume~$x \mathrel\delta y$.
    Then for each~$1 \leq j \leq n$
        we have~$f_{i_j}(x) \mathrel{\varepsilon_j} f_{i_j}(y)$
            and so~$f(x) \mathrel{\widehat{\varepsilon_j}} f(y)$,
            hence~$f(x) \mathrel\varepsilon f(y)$.
        Thus~$f$ is uniformly continuous.

To show uniqueness of~$f$, assume there is a uniformly continuous
    map~$f'\colon Y \to \prod_i X_i$ with~$\pi_i \after f' = f_i$.
Then~$(f'(y))_i = f_i(y) = (f(y))_i$ for every~$y \in Y$
    and so~$f' = f$.
\end{solution}
\begin{solution}{ultranormscalar}%
    Let~$\scrB$ be a von Neumann algebra and~$X$
        a right $\scrB$-module with~$\scrB$-valued inner product.
    We will show that~$x \mapsto xb$ is ultranorm continuous
        for any~$b \in \scrB$.
    It is sufficient to show it is ultranorm continuous at~$0$,
        so assume~$x_\alpha \to 0$ ultranorm for some net~$x_\alpha$ in~$X$.
Let~$f\colon \scrB \to \C$ be any np-map.
    Then~$\|x_\alpha b\|_f^2  
        = f([x_\alpha b, x_\alpha b])
        = f(b^* [x_\alpha , x_\alpha ]b)
        \leq \|b^*b\| f([x_\alpha, x_\alpha])
        = \|b\|^2 \|x_\alpha \|_f^2 \to 0$
        thus~$x_\alpha b \to 0$ ultranorm as well.
\end{solution}
\begin{solution}{mod-projelabs}%
For brevity, write~$p \equiv \langle e,e\rangle$.
Note~$\| e p - e\|^2
    = \| e (1-p) \|^2
    = \langle e (1 - p), e (1-p) \rangle
    = (1-p) \langle e,e\rangle (1-p)
    = (1-p) p (1-p) = 0$.
Thus~$ep - e= 0$ and so~$ep = e$ as desired.
\end{solution}
\begin{solution}{mod-parseval}%
Let~$X$ be a pre-Hilbert~$\scrB$-module for a von Neumann algebra~$\scrB$
    with orthonormal basis~$E \subseteq X$.
Assume~$x \in X$.
By definition of orthonormal basis,
    we know~$x = \sum_{e \in E} e \langle e,x\rangle$
    where the sum converges ultranorm.
That is: $\sum_{e \in S} e \langle e,x\rangle \to x$
    as~$S$ ranges over the finite subsets of~$E$.
    By~\sref{innerprod-ultraweak}
        we
        see~$
        \bigl\langle
        \sum_{e \in S} e \langle e,x\rangle,
        \sum_{e \in S} e \langle e,x\rangle \bigr\rangle
        \to \langle x, x\rangle $ ultraweakly.
Thus~$\sum_{e \in S} \langle x,e\rangle\langle e, x\rangle
    = \sum_{e,d \in S} \langle x,e\rangle\langle e,d\rangle\langle d, x\rangle
    =
        \bigl\langle
        \sum_{e \in S} e \langle e,x\rangle,
        \sum_{e \in S} e \langle e,x\rangle \bigr\rangle
        \to \langle x, x\rangle $ ultraweakly, as desired.
\end{solution}
\begin{solution}{hilbmod-adjoint-exists}%
Let~$T\colon X\to Y$ be a bounded~$\scrB$-linear map between Hilbert~$\scrB$-modules.
Assume~$X$ is self dual.
For any~$y \in Y$, the map~$x \mapsto \langle y, Tx\rangle$
    is~$\scrB$-linear and bounded
    and so by self-duality of~$X$,
    there is a~$t_y \in X$
    with~$\langle t_y, x \rangle = \langle y, Tx\rangle$
    for all~$x \in X$.
For any~$z,y \in Y$ and~$x \in X$,
        we have~$\langle t_z + t_y, x\rangle
            = \langle t_z, x\rangle + \langle t_y, x\rangle
            = \langle z, Tx\rangle + \langle y, Tx \rangle
            = \langle z + y, Tx \rangle
            = \langle t_{z+y}, x\rangle$.
        Thus~$t_z + t_y = t_{z+y}$
For any~$\lambda \in \C$, $y \in Y$ and~$x \in X$
    we have~$\langle \lambda t_y, x \rangle
        =  \langle \lambda y, T x\rangle
        = \langle t_{\lambda y}, x \rangle$
        and so~$\lambda t_y = t_{\lambda y}$.
Hence~$T^* y \equiv t_y$ defines a linear map from~$Y$ to~$X$
    with~$\langle y, Tx\rangle = \langle t_y, x\rangle
    \equiv \langle T^*y, x \rangle$ for all~$x\in X$ and~$y \in Y$.
So~$T^*$ is the adjoint of~$T$.
\end{solution}
\begin{solution}{hilmod-fixed-on-V}%
    Let~$V$ be a right~$\scrB$-module with~$\scrB$-valued inner product
        for some von Neumann algebra~$\scrB$.
    Write~$\eta\colon V \to X$ for the ultranorm completion of~$V$
        from~\sref{dils-completion}.
    Let~$T \in \scrB^a(X)$ be given
        with~$\langle \hat{x}, T \hat{x} \rangle \geq 0$
        for all~$x \in V$.
    We have to show~$T \geq 0$.
    Let~$x \in X$ be an arbitrary vector.
    As all vector states on~$\scrB^a(X)$ are order separating
        by~\sref{hilbmod-denseordersep},
        it is sufficient to show~$\langle x, Tx \rangle \geq 0$.
As the image of~$V$ under~$\eta$ is ultranorm dense in~$X$,
    we can find find a net~$x_\alpha$ with~$\widehat{x_\alpha} \to x$.
Then by~\sref{innerprod-ultraweak}
    we get~$\langle x, Tx\rangle = \uwlim_\alpha \langle \widehat{x_\alpha},
            T\widehat{x_\alpha}\rangle \geq 0$.
            So indeed~$T \geq 0$, as desired.
\end{solution}
\begin{solution}{hilbmod-adj-vector-ncp}%
   Let~$\scrA$ be a C$^*$-algebra
        with~$a_1, \ldots, a_n \in \scrA$.
    Define~$\varphi \colon \scrA \to M_n \scrA$
        by~$\varphi(d) \equiv (a_i^* d a_j)_{ij}$.
        We have to show~$\varphi$ is an ncp-map.
Recall~$\scrA$ is a self-dual Hilbert~$\scrA$-module
    with~$\langle a,b\rangle \equiv a^*b$
    and so its~$n$-fold direct product~$\scrA^n$
    is also self dual (see also \sref{direct-prod-self-dual-basis}.)
Define~$T\colon \scrA^n \to \scrA$
    by~$T((b_i)_i) \equiv \sum_i b_i a_i$
    (i.e.~$T$ is the row-vector~$(a_i)_i$.)
Clearly~$T$ is~$\scrA$-linear.
    It is also bounded: $\| T (b_i)_i \|^2 = \sum_i \|b_i a_i\|^2
        \leq \sum_i A \|b_i\|^2 = A \| (b_i)_i\|^2 $,
            where~$A \equiv \max_i \|a_i\|^2$.
    It's easy that~$T^*(b) \equiv (a_i^* b)_i$ is the adjoint of~$T$.
We may identify~$\scrB^a(\scrA^n) = M_n$
    and then~$\ad_T(d)\,  ((b_i)_i) = T^* d T (b_i)_i
                =  T^* d \sum_i a_i b_i
                = (\sum_i (a_j^* d a_i) b_i)_j 
                = \varphi(d) \, ((b_i)_i)$.
Thus~$\ad_T = \varphi$.
    By~\sref{hilbmod-ad-ncp} the map~$\ad_T$ and thus~$\varphi$ is as well.
\end{solution}
\begin{solution}{direct-prod-self-dual-basis}%
    Assume~$X$ and~$Y$ are self-dual Hilbert~$\scrB$-modules over a
        von Neumann algebra~$\scrB$ with orthonormal bases~$E \subseteq X$
        and~$F \subseteq Y$.
Write~$G \equiv \{(e,0); \ e \in E\} \cup \{ (0,d); \ d \in D \}$.
Clearly~$G$ is orthonormal.
To show~$G$ is an orthonormal basis,
    two conditions remain.
For the first, let~$(x,y) \in X \oplus Y$ be given.
As~$E$ and~$F$ are orthonormal bases
    we know~$x = \sum_{e\in E} e \langle e,x\rangle$
    and~$y = \sum_{f \in F} f \langle f,x\rangle$,
    where the sums converge ultranorm.
    The inclusions~$x \mapsto (x,0)$ and~$y \mapsto (0,y)$
        are bounded~$\scrB$-linear and thus ultranorm continuous,
        hence
\begin{align*}
        \sum_{g\in G} g \langle g, (x,y)\rangle
        &\ =\  \Bigl(\sum_{e\in E}
                        (e,0) \langle (e,0), (x,y) \rangle \Bigr)  \ +\ 
                        \Bigl( \sum_{f \in F} 
                        (0,f) \langle (0,f), (x,y) \rangle \Bigr) \\
        &\ =\  \Bigl(\sum_{e\in E}
                        (e,0) \langle e, x \rangle \Bigr) \ + \ 
                        \Bigl( \sum_{f \in F} 
                        (0,f) \langle f, y \rangle \Bigr) \\
        &\ =\  \Bigl(\sum_{e \in E} e \langle e,x \rangle ,
                \sum_{f \in F} f \langle f, y \rangle  \Bigr)\\
                    & \ =\ (x,y),
\end{align*}
which proves the first condition.
    For the second condition, let~$(b_g)_{g\in G}$
    be some~$\ell^2$-summable family from~$\scrB$.
The subfamilies~$(b_{e_0})_{e \in E}$
    and~$(b_{0,f})_{f \in F}$
    are~$\ell^2$-summable as well.
Hence the
    sums~$\sum_{e \in E} e b_{(e,0)}$
    and~$\sum_{f \in F} f b_{(0,f)}$ converge ultranorm.
Thus~$\sum_{e \in E} (e,0) b_{(e,0)}
    + \sum_{f \in F} (0,f) b_{(0,f)}
    =\sum_{g \in G} g b_g$ converges ultranorm as well.
We have shown~$G$ is an orthonormal basis of~$X\oplus Y$.
Consequently~$X \oplus Y$ is self dual
    by~\sref{dils-selfdual}.
\end{solution}
\begin{solution}{selfdual-orthn-basis}%
    Assume~$X$ is a self-dual Hilbert~$\scrB$-module for a von Neumann
    algebra~$\scrB$. Suppose~$E \subseteq X$ is an orthonormal set.
We will show~$E$ is an orthonormal basis of~$E^{\perp\perp}$.
Clearly~$E$ is orthonormal.
    Because of this, and the fact that~$E^{\perp\perp}$ is ultranorm
        closed (by~\sref{hilbmod-projthm}) 
        we know~$\sum_e b b_e$ converges ultranorm in~$E^{\perp\perp}$
        for any~$\ell^2$-family $(b_e)_e$.
Assume~$x \in E^{\perp\perp}$.
To show~$E$ is an orthonromal base,
    it only remains to be shown that~$x = \sum_e e \langle e, x\rangle$.
Define~$x' \equiv x - \sum_e e \langle e, x\rangle$.
By~\sref{hilbmod-projthm} we know~$E^{\perp\perp}$ is ultranorm closed
    and so~$x' \in E^{\perp\perp}$.
For any~$e_0\in E$
    we also have~$\langle e_0, x'\rangle = \langle e_0,x \rangle
            - \sum_{e \in E} \langle e_0, e\rangle\langle e, x\rangle
            = 0 $ and so~$x' \in E^{\perp}$.
    Hence~$\langle x',x'\rangle = 0$, so indeed~$x = \sum_e e\langle e, x\rangle$
        and~$E$ is an orthonormal basis of~$E^{\perp\perp}$.

For the second part, assume~$x \in X$.
    By Parseval's identity (see \sref{mod-parseval})
        we have~$\langle x, x\rangle = \sum_e \langle x,e \rangle\langle e,x\rangle$
        for any~$x \in E^{\perp\perp}$.
To prove the converse,
    assume~$\langle x, x\rangle = \sum_e \langle x,e \rangle\langle e,x\rangle$.
By~\sref{hilbmod-projthm} we know~$x = x' + x''$
    for~$x' \in E^{\perp\perp}$  and~$x'' \in E^\perp$.
Note~$\langle e, x' \rangle = \langle e, x\rangle$ for any~$e \in E$
    and so by Parseval's identity for~$E^{\perp\perp}$
    we see~$\langle x',x'\rangle = \sum_e \langle x, e\rangle\langle e, x\rangle
        = \langle x, x\rangle$.
Now, using~$\langle x'',x'\rangle = 0$
    we see~$\langle x,x\rangle = \langle x'', x''\rangle + \langle x', x'\rangle
        = \langle x'', x''\rangle + \langle x, x\rangle$
        and so~$\langle x'', x''\rangle=0$, whence~$x \in E^{\perp\perp}$.
\end{solution}
\begin{solution}{selfdual-gramschmidt}%
    Let~$X$ be a self-dual Hilbert~$\scrB$-module
        for some von Neumann algebra~$\scrB$.
    Assume~$x_1, \ldots, x_n \in X$.
We will show that there is a finite orthonormal set~$E$
    of~$n$ or fewer elements such that~$E$ is a
    basis of~$\{x_1, \ldots, x_n\}^{\perp\perp}$.
We do this by induction over~$n$.
For~$n=0$, the set~$E = \emptyset$ suffices.
For the induction step,
    assume~$n > 0$ and~$E'$ is an
    orthonormal basis of~$\{x_1, \ldots, x_{n-1}\}^{\perp\perp}$.
    Write~$x' \equiv x_n - \sum_{e\in E'} e\langle e,x_n\rangle$. 
If~$x'=0$, then~$E \equiv E'$ suffices.
For the other case, assume~$x' \neq 0$.
    By polar decomposition (see the end of~\sref{selfdual-bcompl-then-basis}),
    there is an~$u \in X$
    with~$x' = u\langle x',x'\rangle^{\frac{1}{2}}$
    and~$\langle u, u \rangle = \ceil{\langle x',x'\rangle}$.
Define~$E \equiv E' \cup \{ u \}$.
Clearly~$E$ is an othonormal set of~$n$ or fewer elements.
    By the induction assumption
    and~\sref{selfdual-orthn-basis},
    we know~$x_i \in E'^{\perp\perp} \subseteq E^{\perp\perp}$ for~$i \leq 1 \neq n-1$.
For any~$d \in E^\perp$
    we have~$\langle d, x'\rangle =
    \langle d, u \rangle \langle x',x'\rangle^{\frac{1}{2}} = 0$
    and so~$x' \in E^{\perp\perp}$.
    Clearly~$\sum_{e \in E'}e \langle e,x_n \rangle \in E'^{\perp\perp}
        \subseteq E^{\perp\perp}$
        and so~$x_n = x' + \sum_{e \in E'} e \langle e,x_n\rangle \in
            E^{\perp\perp}$.
    Together with the previous,
    we see~$\{x_1, \ldots, x_n\}^{\perp\perp} \subseteq E^{\perp\perp\perp\perp}
        = E^{\perp\perp} \subseteq \{x_1, \ldots, x_n\}^{\perp\perp}$,
        thus~$E^{\perp\perp} = \{ x_1, \ldots, x_n\}^{\perp\perp}$.
    By~\sref{selfdual-orthn-basis}
        we know~$E$ is an orthonormal basis of~$E^{\perp\perp}$,
        which completes the proof by induction.
\end{solution}
\begin{solution}{hilbmod-el2}%
For brevity write~$\ell^2 \equiv \ell^2((p_i)_{i \in I})$.
We will first proof that~$\ell^2$ is a right~$\scrB$-module.
Assume~$(a_i)_i, (b_i)_i \in \ell^2$.
    We want to show~$(a_i+b_i)_i \in \ell^2$.
    First, we to show~$\sum_i (a_i + b_i)^*(a_i+b_i)$ is bounded.
Pick~$A,B\in \R^+$ with~$\sum_i a_i^*a_i \leq A$
    and~$\sum_i b_i^* b_i \leq B$ which exist
    as~$(a_i)_i$ and~$(b_i)_i$ are~$\ell^2$.
Let~$f$ be any normal state on~$\scrB$
    and~$S \subseteq I$ some finite subset.
Then
\begin{align*}
    f \Bigl(\sum_{i \in S} (a_i + b_i)^*(a_i + b_i) \Bigr)
        & \ =\ \sum_{i \in S} \| a_i + b_i\|_f^2  \\
        & \ \leq\ \sum_{i \in S} \| a_i\|_f^2  +\| b_i\|_f^2 
                                + 2\|a_i\|_f \|b_i\|_f \\
        & \ \leq\ A+B+ 2\sum_{i \in S}  \|a_i\|_f \|b_i\|_f.
\end{align*}
    By Cauchy--Schwarz~$\sum_{i \in S} \|a_i\|_f \|b_i\|_f
        \leq \bigl(\sum_{i \in S} \|a_i\|^2_f\bigr)^{\frac{1}{2}}
            \bigl(\sum_{i \in S} \|b_i\|^2_f \bigr)^{\frac{1}{2}}
            \leq (AB)^{\frac{1}{2}}$.
    As normal states are order separating,
        we see that we have a bounded and thus
        converging sum of positive elements~$\sum_i (a_i+b_i)^*(a_i+b_i) \leq 
                    A+B+(AB)^\frac{1}{2}$.
Suppose~$i \in I$.
    It remains to be shown~$\ceil{(a_i+b_i)(a_i+b_i)^*} \leq p_i$.
Recall from~\sref{ceill-basic}
    that~$\ceil{xx^*} \leq p_i $ if and only if~$p_i x = x$.
Clearly~$p_i (a_i +b_i) =  p_i a_i + p_i b_i = a_i+b_i$
    as~$(a_i)_i, (b_i)_i \in \ell^2$ and so indeed~$(a_i+b_i)_i \in \ell^2$.

Suppose~$b \in \scrB$
    and~$(a_i)_i \in \ell^2$.
    Then~$b^* \bigl( \sum_i a_i^* a_i \bigr) b
            = \sum_i (a_i b)^* a_i b$ as~$x \mapsto b^*xb$ is normal
            by~\sref{ad-normal}
            and so~$(a_ib)_i$ is~$\ell^2$.
    Furthermore~$p_i a_i b = a_i b$ for any~$i \in I$,
     so~$(a_i b)_i \in \ell^2$.
     We have shown~$\ell^2$ is a right~$\scrB$-module
     with coordinatewise operations.

    Next we will show~$\langle (a_i)_i , (b_i)_i \rangle \equiv
        \sum_i a_i^* b_i$ defines a~$\scrB$-valued inner product
        on~$\ell^2$.
    First we have to show the sum converges.
    Suppose~$f$ is any normal state on~$\scrB$.
Pick~$A,B\in \R^+$
    with~$\sum_i a_i^*a_i \leq A$
    and~$\sum_i b_i^*b_i \leq B$.
For any finite subset~$S \subseteq I$
    we have
\begin{equation*}
    \Bigl|f \Bigl( \sum_{ i \in S} a_i^*b_i \Bigr) \Bigr|
    \ \leq \  \sum_{i \in S} | [a_i,b_i]_f | \ \leq \ 
    \sum_{i \in S} \|a_i\|_f \|b_i\|_f
    \ \leq\ (AB)^{\frac{1}{2}},
\end{equation*}
where we used Cauchy--Schwarz for~$\scrB$-valued inner products
    in the second inequality and classic
    Cauchy--Schwarz in the final inequality.
We have shown that~$\sum_{i \in S} a_i^*b_i$
    is a norm-bounded net in~$S$.
    We claim it is ultraweakly Cauchy as well.
For now, pick any finite sets~$S,T \subseteq I$.
Assume~$f$ is any normal state on~$\scrB$.
We want to show that the following quantity vanishes
    for sufficiently large~$S \cap T$.
\begin{equation}\label{e1612eq2}
\Bigl| f \Bigl( \sum_{i \in S} a_i^*b_i
                - \sum_{i \in T} a_i^*b_i \Bigr)\Bigr|
                \ \leq \ 
                \Bigl| \sum_{i \in S - T} f(a_i^*b_i) \Bigr|
                    \ +\  \Bigl| \sum_{i \in T -S} f(a_i^*b_i) \Bigr|.
\end{equation}
Note that~$[ (a_i)_i, (b_i)_i ] \equiv \sum_{i\in S-T} f(a_i^*b_i)$
    is an inner product and so
\begin{equation}\label{e1612eq1}
        \Bigl| \sum_{i \in S - T} f(a_i^*b_i) \Bigr|
        \ \leq \ 
         \Bigl( \sum_{i \in S - T} f(a_i^*a_i) \Bigr)^{\frac{1}{2}}
         \Bigl( \sum_{i \in S - T} f(b_i^*b_i) \Bigr)^{\frac{1}{2}}.
\end{equation}
The sum~$\sum_i f(a_i^*a_i)$ converges
    and so~$\sum_{i \in S - T} f(a_i^*a_i)$
    can be made arbitrarily small
    by picking sufficiently large~$S\cap T$.
And so (with a similar argument for the other factor),
    we see that~\eqref{e1612eq1} vanishes,
    which is the left term of~\eqref{e1612eq2}.
The argument for the other term of~\eqref{e1612eq2} is the same.
Thus~$\sum_{i \in S} a_i^*b_i$ is ultraweakly Cauchy
    and converges by~\sref{bh-bounded-uw-complete}.
From the fact that~$a \mapsto a^*$,
    $a \mapsto ab$ and~$(a,b) \mapsto a+b$
    are all ultraweakly continuous,
    it follows readily that~$\langle (a_i)_i , \rangle (b_i)_i\rangle \equiv
        \sum_i a_i^*b_i$ is an inner product on~$\ell^2$.
    We claim this inner product is definite.
Assume~$0 = \langle (a_i)_i, (a_i)_i\rangle \equiv \sum_i a_i^*a_i$.
Then~$a_i^*a_i = 0$~for each~$i \in I$ and so~$a_i = 0$,
    which is to say~$(a_i)_i = 0$. Indeed our inner product is definite.
    Thus~$\ell^2$ is a pre-Hilbert~$\scrB$-module.

Write~$E \equiv \{ \delta_i; \ i \in I\}$,
    where~$(\delta_i)_j = 0$ for~$i \neq j$
    and~$(\delta_i)_i = p_i$.
Clearly~$E$ is an orthonormal set.
We claim it's an orthonormal basis of~$\ell^2$.
Assume~$(a_i)_i \in \ell^2$.
It is easy to see~$\sum_{i \in I} \delta_i a_i = (a_i)_i$
    and so~$\sum_{\delta_i \in E} \delta_i \langle \delta_i, (a_i)_i\rangle
        = \sum_{i \in I} \delta_i a_i = (a_i)_i$.
    It remains to be shown that~$\sum_{\delta_i \in E} \delta_i b_{\delta_i}$
        converges ultranorm for any~$\ell^2$ summable~$(b_i)_{\delta_i \in E}$,
        which indeed it does to~$(b_{\delta_i})_{i \in I}$
         as we already saw.
Thus~$\ell^2$ is self dual.

For the final part of the exercise, assume~$X$ is a self-dual Hilbert~$\scrB$-module
    over some von Neumanna algebra~$\scrB$ with orthonormal basis~$E \subseteq X$.
Define~$\vartheta\colon X \to \ell^2((\langle e,e\rangle)_{e \in E})$
    by~$\vartheta(x) \equiv (\langle e, x\rangle)_e$.
Clearly~$\vartheta$ is~$\scrB$-linear.
It also preserves the inner product:
$\langle x, y\rangle = \langle x, \sum_e e \langle e, y\rangle \rangle
    = \sum_e \langle x, e\langle e, y\rangle\rangle
    = \sum_e \langle x, e\rangle\langle e, y\rangle
    = \langle \vartheta(x), \vartheta(y)\rangle$.
    which entails it's injective.
To show~$\vartheta$ is surjective,
    let~$(x_e)_e \in \ell^2((\langle e,e\rangle)_{e \in E})$
    be given.
    The family~$(x_e)_e$ is~$\ell^2$
        so~$\sum_e e b_e$ converges ultranorm.
        Clearly~$\vartheta (\sum_e e b_e)
                = (\langle e, \sum_e e b_e\rangle )_e
                = (b_e)_e$,
    so~$\vartheta$ is indeed surjective.
    It follows~$\vartheta$ is an isomorphism~$X \cong
    \ell^2((\langle e,e\rangle)_{e \in E})$.
\end{solution}
\begin{solution}{onb1}%
    Let~$X$ be a self-dual Hilbert~$\scrB$-module~$X$
        for some von Neumann algebra~$\scrB$
        with orthonormal basis~$(e_i)_{i\in I}$.
Assume~$(u_i)_{i \in I}$ is a family partial isometries
    with~$u_iu_i^* = \langle e_i, e_i \rangle$.
    We will show~$(e_iu_i)_{i \in I}$ is an orthonormal basis
    of~$X$.
For brevity, write~$d_i \equiv e_iu_i$.
To start, note~$\langle d_i, d_j\rangle
            = u_i^* \langle e_i, e_j\rangle u_j$,
            which is zero if~$i \neq j$.
If~$i = j$, then~$\langle d_i, d_i\rangle = u_i^* \langle e_i, e_i\rangle u_i
            = u_i^*u_i $ which is a projection.
    So~$(d_i)_i$ is orthonormal.
Hence~$\sum_i d_i b_i$ converges for any~$\ell^2$-family~$(d_i)_i$.
It remains to be shown that~$x = \sum_i d_i \langle d_i, x\rangle$.
Note~$\sum_i d_i \langle d_i, x\rangle
    =  \sum_i e_i \langle e_i u_i u_i^*, x \rangle
    =  \sum_i e_i \langle e_i \langle e_i, e_i\rangle, x \rangle
    =  \sum_i e_i \langle e_i, x \rangle = x$,
    so indeed~$(d_i)_i$ is an orthonormal basis of~$X$.

For the next part, assume~$(p_i)_{i \in I}$ and~$(q_i)_{i \in I}$
    are projections with~$p_i \sim q_i$.
Let~$p_i$ denote the partial projection with~$u_i^* u_i = q_i$
    and~$u_i u_i^* = p_i$.
Consider~$\ell^2((p_i)_{i \in I})$.
Define~$((\delta_i)_i) \in \ell^2$
    by~$ (\delta_i)_j = 0$ if~$i \neq j$
        and~$(\delta_i)_i = p_i$.
    In~\sref{hilbmod-el2} we saw~$(\delta_i)_i$ is an
        orthonormal basis of~$\ell^2((p_i)_{i \in I})$.
    Note~$u_iu_i^* = p_i = \langle \delta_i, \delta_i\rangle$.
    By the previous~$\delta_i u_i$ is another orthonormal basis of~$
        \ell^2((p_i)_{i \in I})$.
    By the second part of~\sref{hilbmod-el2}
        we see~$
        \ell^2((p_i)_{i \in I}) \cong
        \ell^2((\langle \delta_i u_i, \delta_i u_i\rangle)_{i \in I} ) =
        \ell^2((u_i^* p_i u_i  )_{i \in I} ) =
        \ell^2((q_i )_{i \in I} )$,
        as promised.
\end{solution}
\begin{solution}{onb2}%
Let~$X$ be a self-dual Hilbert~$\scrB$-module.
Assume~$E \subseteq X$ and~$e_1,e_2 \in X$
    such that~$E \cup \{ e_1, e_2\} $ is an orthonormal basis
    and~$\langle e_1, e_1 \rangle + \langle e_2, e_2\rangle \leq 1$
We will show that~$E' \equiv E \cup \{e_1 + e_2\}$
        is an orthonormal basis as well.
For brevity, write~$p_1 \equiv \langle e_1,e_1\rangle$
    and~$p_2 \equiv \langle e_2, e_2\rangle$.
Clearly~$E$ itself is an orthonormal set.
For any~$e \in E$ we have~$\langle e, e_1+e_2 \rangle =
    \langle e, e_1\rangle + \langle e,e_2\rangle = 0$
    and so~$E'$ is an orthogonal set.
By assumption~$p_1$ and~$p_2$ are projections
    with~$p_1 + p_2 \leq 1$
    and so by~\sref{orthogonal-tuple-of-projections}
    they are orthogonal and thus in particular~$p_1 + p_2$
        is again a projection.
Hence~$\langle e_1 + e_2, e_1 + e_2 \rangle
        = p_1 + p_2$ is a projection.
    Thus~$E'$ is an orthonormal set.

Let~$x \in X$ be given.
By~\sref{mod-projelabs} we have~$e_2 = e_2 p_2$
    and so~$e_2 \langle e_1 + e_2, x\rangle
                = e_2 \langle (e_1 + e_2)p_2, x\rangle
                = e_2 \langle e_2, x \rangle$.
Similarly~$e_1 \langle e_1 + e_2, x\rangle
                = e_1 \langle e_1, x\rangle$.
Thus~$ x = e_1 \langle e_1, x\rangle + e_2 \langle e_2, x\rangle 
            + \sum_{e \in E} e \langle e, x\rangle
        = (e_1 + e_2) \langle e_1 + e_2, x\rangle
            + \sum_{e \in E} e \langle e, x\rangle$,
            which shows the first condition on an orthonormal
            set to be an orthonormal basis.
The second (and final) condition
    holds automatically as~$E'$ is an orthonormal set
        and~$X$ is ultranorm complete.  Thus~$E'$ is indeed
        an orthonormal basis.

The the last part of the exercise, assume~$p,q \in \scrB$
    are projections with~$p+q \leq 1$.
    Clearly~$\{p+q\}$ is an orthonormal basis of~$(p+q)\scrB$
        and so is~$\{p,q\}$ by the previous.
    Hence by~\sref{hilbmod-el2}
        we see~$(p+q) \scrB \cong \ell^2(\{p, q\})
            = p\scrB \oplus q\scrB$, as promised.
\end{solution}
\begin{solution}{hilbmod-tensor-ketbra}%
Assume~$X$ is a self-dual Hilbert~$\scrA$-module
    and~$Y$ is a self-dual Hilbert~$\scrB$-module
        for von Neumann algebras~$\scrA$ and~$\scrB$.
\begin{enumerate}
\item
    Assume~$x_1,x_2,x \in X$ and~$y_1,y_2,y \in Y$.
        Then
    \begin{align*}
        \ketbra{x_1}{x_2} \otimes \ketbra{y_1}{y_2} 
            \, x \otimes y 
        &\ =\  (\ketbra{x_1}{x_2} x) \otimes (\ketbra{y_1}{y_2} y) \\
        &\ =\ (x_1 \langle x_2, x\rangle) \otimes (y_1 \langle y_2, y\rangle) \\
        &\ =\ (x_1 \otimes y_1)  (\langle x_2, x\rangle \otimes \langle y_2, y\rangle) \\
        &\ =\ (x_1 \otimes y_1)  \langle x_2 \otimes y_2, x \otimes y\rangle \\
        &\ =\ \ketbra{x_1 \otimes y_1}{x_2 \otimes y_2} \, x\otimes y.
    \end{align*}
This is sufficient to show
    that~$\ketbra{x_1}{x_2} \otimes \ketbra{y_1}{y_2} 
        = \ketbra{x_1 \otimes y_1}{x_2 \otimes y_2}$,
    either by appealing to the defining universal property of~$X \otimes Y$
    or by the related property
    that the~$\scrA \odot \scrB$-linear
    span of~$\{x \otimes y; \ x \in X, y\in Y\}$
    is ultranorm dense in~$X \otimes Y$.
        (In fact the~$\C$-linear span is already ultranorm dense.)
\item
    Clearly, for any~$x \in X$ and~$y \in Y$
            we have~$(1 \otimes 1) \, (x \otimes y)
                    = 1 \, x \otimes y$
        and so by same token (as the conclusion of the previous point)
        we see~$1 = 1\otimes 1$.
\item
    Assume~$x \in X$, $y \in Y$, $S,S' \in \scrB^a(X)$ and~$T,T' \in \scrB^a(Y)$.
        Then
\begin{equation*}
    (S \otimes T) (S' \otimes T') \, x \otimes y
         \ = \ (SS' x) \otimes (TT' y) \\
         \ = \ (SS')  \otimes (TT') \, x \otimes y,
\end{equation*}
        which is sufficient to show~$(SS') \otimes (TT') = (S\otimes T) \otimes (S' \otimes T')$ (see the conclusion of the first point.)
\item
    Assume~$S \in \scrB^a(X)$, $T \in \scrB^a(Y)$, $x_1,x_2\in X$
        and~$y_1,y_2\in Y$.
    Then
\begin{align*}
    \langle (S \otimes T)^* \, x_1 \otimes y_1, \ x_2 \otimes y_2\rangle    
& \ = \  \langle x_1 \otimes y_1, \ (S \otimes T) \, x_2 \otimes y_2\rangle \\
    & \ = \  \langle x_1 \otimes y_1, \ (Sx_2) \otimes (Ty_2) \rangle \\
    & \ = \ \langle x_1, Sx_2\rangle \otimes \langle y_1 , T y_2\rangle\\
    & \ = \ \langle S^* x_1, x_2\rangle \otimes \langle T^* y_1 ,  y_2\rangle\\
    & \ = \  \langle (S^* \otimes T^*)\, x_1 \otimes y_1, \ x_2 \otimes y_2\rangle.
\end{align*}
This is sufficient (see the conclusion of the first point)
        to show that the vector functionals
        for~$(S \otimes T)^* x_1\otimes y_1$
        and~$(S^*\otimes T^*)x_1 \otimes y_1$ agree.
        Hence~$(S \otimes T)^* x_1\otimes y_1= (S^*\otimes T^*)x_1 \otimes y_1$.
        In turn this is sufficient to show~$(S \otimes T)^* = S^* \otimes T^*$,
        as desired.
\end{enumerate}
\end{solution}
\begin{solution}{dils-filter-basics-exercise}%
    Assume~$\varphi\colon \scrA \to \scrB$ is an ncp-map between
        von Neumann algebras.
\begin{enumerate}
\item
    Assume~$c\colon \scrB \to \scrC$ is a filter
            and that~$(\scrP, \varrho, h)$ is a Paschke dilation of~$\varphi$.
    We will show that~$(\scrP, \varrho, c \after h)$ is a Paschke
        dilation of~$c \after \varphi$.
    To this end, assume~$\varrho\colon \scrA \to \scrP'$ is some nmiu-map
            and~$h'\colon \scrP' \to \scrC$ is any ncp-map
            with~$h' \after \varrho' = c \after \varphi$.
    Assume~$\varphi(1) \leq 1$.
        (We will reduce the general to this one later on.)
        Note~$h'(1) = h'(\varrho'(1)) = c(\varphi(1)) \leq c(1)$
            and so by the universal property of~$c$,
            there is a unique ncp-map map~$h''\colon \scrP' \to \scrB$
                with~$c \after h'' = h'$.
    Then~$c \after \varphi = h' \after \varrho' =
             c \after h'' \after \varrho'$.
    And so~$\varphi = h'' \after \varrho'$
        as~$c$ is injective, see \sref{dils-filters-injective}.
    There is a unique ncp-map~$\sigma\colon \scrP' \to \scrP$
            with~$h \after \sigma = h''$
                and~$\sigma \after \varrho' =\varrho$.
Then~$h' = c \after h'' = c \after h \after \sigma$.
Assume~$\sigma'\colon \scrP'\to \scrP$
    is any ncp-map with~$\sigma' \after \varrho = \varrho$
            and~$c \after h \after \sigma' = h'$.
Then~$c \after h \after \sigma' = h' = c \after h''$
    and so~$h \after \sigma' = h''$.
    Hence~$\sigma=\sigma'$ by the uniqueness of~$\sigma$.
Thus~$(\scrP, \varrho, c\after h)$ is a Paschke dilation of~$c\after\varphi$.

    Now, in the general case, assume~$\varphi(1) \nleq 1$.
    Then~$\varphi(1) \neq 0$.
        Define~$\varphi' \equiv \|\varphi(1)\|^{-1} \varphi$.
        Then~$\varphi'(1) = \|\varphi(1)\|^{-1} \varphi(1) \leq 1$.
        By~\sref{paschke-basics}
         we know~$(\scrP, \varrho, \|\varphi(1)\|^{-1} h)$
            is a Paschke dilation of~$\varphi'$.
    Now we are in our previous case
        and so~$(\scrP, \varrho, \|\varphi(1)\|^{-1} c \after h)$
            is a Paschke dilation of~$c \after \varphi'$.
    Applying~\sref{paschke-basics} once again,
            we conclude~$(\scrP, \varrho,c \after h)$
        is a Paschke dilation of~$c \after \varphi$ as desired.
\item
Assume~$c' \colon \scrC' \to \scrB$ is a filter of~$\varphi(1)$.
By the universal property of~$c'$,
    there is a unique ncp-map~$\varphi'\colon \scrA \to \scrC'$
            with~$c'\after \varphi' = \varphi$.
Then~$c'(\varphi'(1)) = \varphi(1) = c'(1)$
        and so by injectivity of~$c'$ (see~\sref{dils-filters-injective})
            we get~$\varphi'(1) = 1$ ---
            that is to say: $\varphi'$ is unital.

For the final part of this exercise, assume~$(\scrP, \varrho, h)$
    is a Paschke dilation of~$\varphi'$.
By the previous point,
    we know~$(\scrP, \varrho, c' \after h)$
        is a Paschke dilation of~$c' \after \varphi' = \varphi$,
        as desired.
\end{enumerate}
\end{solution}
\begin{solution}{dils-filters-injective}%
Let~$c\colon \scrC \to \scrB$ be some filter for~$b \in \scrB$.
    By~\sref{dils-stand-filter}
        we know~$c_b\colon \ceil{b} \scrB \ceil{b} \to \scrB$
            given by~$a \mapsto \sqrt{b} a \sqrt{b}$
                is also a filter for~$b$.
    By the universal property of~$c$,
        there is a unique ncp-map~$\vartheta_1$
            with~$c \after \vartheta_1 = c_b$.
    In the other direction, by the universal property of~$c_b$,
        there is a unique ncp-map~$\vartheta_2$
            with~$c_b \after \vartheta_2 = c$.
    Hence~$c \after (\vartheta_1 \after \vartheta_2)
            = c_b \after \vartheta_2 = c \after \id$.
Clearly~$c(1) \leq c(1)$
    and so by the universal property of~$c$
    the map~$\id$ is the unique ncp-map
    with~$c \after \id = c$, hence~$\vartheta_1 \after \vartheta_2 = \id$.
Reasoning similarly on the other side,
    we see~$\vartheta_2 \after \vartheta_1 = \id$
    and so~$\vartheta_2$ is an isomorphism.
Recall~$c = c_b \after \vartheta_2$
    and so it is sufficient to show~$c_b$ is injective.
    To this end, let~$a_1,a_2 \in \ceil{b} \scrB \ceil{b}$ be given
    with~$c_b(a_1) = c_b(a_2)$;
    i.e.~$\sqrt{b}a_1\sqrt{b} = \sqrt{b}a_2\sqrt{b}$.
    Then~$a_1 = a_2$ by~\sref{mult-cancellation}
        and so~$c_b$ is indeed injective.
\end{solution}
\begin{solution}{surjective-nmiu}%
Assume~$\varrho \colon \scrA \to \scrB$
    is a surjective nmiu-map between von Neumann algebras.
We already saw in~\sref{kernel-ultraweak-twosided-ideal-dils}
    that the kernel of an nmiu-map is an ultraweakly-closed
    two-sided ideal and hence a principal ideal of a central projection
        by~\sref{prop:weakly-closed-ideal}.
    Write~$z$ for the central projection of~$\scrA$
        with~$\ker \varrho = (1-z) \scrA$;
        $h_z \colon \scrA \to z\scrA$ for the map~$h_z(a) = za$
        and~$\varrho'\colon z\scrA \to \scrB$
        for the restriction of~$\varrho$ to~$z \scrA$.
    Clearly~$\varrho' \after h_z = \varrho$ and so~$\varrho'$ is surjective.
    Furthermore~$\varrho'$ is injective,
        because for any~$a \in z\scrA$ with~$\varrho(a)=0$
        we have~$\varrho(a)=0$ and so~$a \in \ker \varrho = (1-z)\scrA$,
        hence~$a = 0$.
It is easy to see that~$\varrho'$ as a bijective miu-map
    has an miu inverse and so it is an miu-isomorphism.
    Consequently, it is an~nmiu-isomorphism (as
        miu-maps are order-preserving)
        and in particular an ncp-isomorphism.
The map~$h_z$ is a standard corner for~$z$ (see~\sref{standard-corner-dils})
    and so~$\varrho = \varrho' \after h_z$ is also a corner for~$z$
    as~$\varrho'$ is an ncp-isomoprhism.

Conversely, assume~$h\colon \scrA \to \scrC$
        is a corner for a central projection~$z \in \scrA$.
The map~$h_z\colon \scrA \to z \scrA$ given by~$a \mapsto za$
    is another corner for~$z$, see~\sref{standard-corner-dils}.
By the universal property of~$h$,
    there is a unique ncp-map~$\vartheta_1\colon z\scrA \to \scrC$
            with~$h \after \vartheta_1 =\vartheta_2$.
On the other side, by the universal property of~$h_z$,
    there is a unique ncp-map~$\vartheta_2\colon \scrC\to z\scrA$
    with~$h_z \after \vartheta_2 = \vartheta_1$.
Again, by the universal property of~$h$,
    the map~$\id\colon \scrC \to \scrC$
    is the unique ncp-map~$\scrC \to \scrC$
        with~$h = h \after \id$.
    Note~$h = h_z \after \vartheta_2 = h \after (\vartheta_1 \after \vartheta_2)$ and so~$\id = \vartheta_1 \after \vartheta_2$.
    Similarly~$\vartheta_2 \after \vartheta_1 = \id$.
    Thus~$\vartheta_1$ is an ncp-isomorphism
        and consequently an nmiu-isomoprhism by~\sref{iso}.
Clearly~$h_z$ is a surjective nmiu-map
    and thus~$h = h_z \after \vartheta_2$
    is a surjective nmiu-map as well.
\end{solution}
\begin{solution}{pcm-preorder}%
Reflexivity is easy:~$a \ovee 0 = 0 \ovee a = a$
    by the partial commutativity and zero axiom, so~$a \leq a$.
For transitivity, assume~$a \leq b$ and~$b \leq c$.
Pick~$d,e \in M$
    with~$a \ovee d = b$ and~$b \ovee e = c$.
Then by partial
associativity~$c = b \ovee e = (a \ovee d) \ovee e = a \ovee (d \ovee e)$
so~$a \leq c$.
\end{solution}
\begin{solution}{ea-product}%
Assume~$E$ and~$F$ are effect algebras.
We define~$a \ovee b$ for~$a,b \in E\times F$
    whenever both~$a_1 \perp b_1$ and~$a_2 \perp b_2$
    and in that case~$a\ovee b \equiv (a_1 \ovee b_1, a_2 \ovee b_2)$.
Furthermore, define~$1 \equiv (1,1)$ and~$0 \equiv (0,0)$.
We claim that these operations turn~$E \times F$ into an effect algebra.

We start with partial commutativity.
Assume~$a,b \in E\times F$ with~$a \perp b$.
Then~$a_1 \perp b_1$ and~$a_2 \perp b_2$,
    hence~$b_1 \perp a_1$, $b_2 \perp a_2$
    and~$a \ovee b = (a_1 \ovee b_1, a_2 \ovee b_2)
                    = (b_1 \ovee a_1, b_2 \ovee a_2)
                    = b \ovee a$.

We continue with partial associativity.
Assume~$a,b,c \in E\times F$ with~$a\perp b$ and~$a \ovee b \perp c$.
Then~$a_i \perp b_i$, $a_i \ovee b_i \perp c_i$ for~$i=1,2$
    and so~$b_i \perp c_i$ and~$(a_i \ovee b_i) \ovee c_i = a_i \ovee (
    b_i \ovee c_i)$. Consequently~$(a \ovee b) \ovee c = a \ovee (b\ovee c)$.

    Next the zero axiom: for any~$a \in E\times F$
        we have~$a_1 \perp 0$,~$a_2 \perp 0$, $0 \ovee a_1 =a_1$
        and~$0 \ovee a_2 = a_2$.
            Thus~$0 \equiv (0,0) \perp (a_1,a_2)$ and~$0 \ovee a = a$.

    To prove the orthocomplement law, assume~$a \in E \times F$.
        Then for~$a \equiv (a_1^\perp, a_2^\perp)$
            we have~$a \perp a^\perp$ and~$a \ovee a^\perp = (1,1) \equiv 1$.
    Furthermore, if~$a \vee b = 1$ for some~$b \in E \times F$,
        then~$b_1 = a_1^\perp$ and~$b_2 = a_2^\perp$, hence~$b = a^\perp$.

Only the zero--one law remains.  Asumme~$a\perp 1$ for some~$a \in E \times F$.
Then~$a_1 \perp 1$ and~$a_2 \perp 1$, hence~$a_1=0$ and~$a_2=0$.
Thus~$a = 0$, as desired.
We have shown~$E \times F$ is an effect algebra with componentwise operations.

Write~$\pi_1\colon E\times F \to E$ and~$\pi_2\colon E\times F \to F$
        for the maps given by~$\pi_1(a_1,a_2) = a_1$
            and~$\pi_2(a_1,a_2)=a_2$.
    Clearly~$\pi_1$ and~$\pi_2$ are additive and unital, hence
    effect algebra homomorphisms.
We will show that~$E\times F$ with projections~$\pi_1$ and~$\pi_2$
    forms a categorical product of~$E$ and~$F$ in~$\mathsf{EA}$.
To this end, assume~$f_1\colon G \to E$ and~$f_2\colon G \to F$
    are effect algebra homomorphisms for some effect algebra~$G$.
Define~$g\colon G \to E\times F$ by~$g(a) = (f_1(a), f_2(a))$.
Clearly~$g$ is the unique map with~$\pi_i \after g = f_i$ for~$i=1,2$.
It is also an effect algebra homomorphism and
    so~$E \times F$ is indeed the product of~$E$ and~$F$.
\end{solution}
\begin{solution}{ea-redund}%
Let~$E$ be an effect algebra with any~$a \in E$.
We will first show that~$1^\perp = 0$ and~$a^{\perp\perp} = a$.
By the orthocomplement law~$1 \ovee 1^\perp = 1$,
    so~$1^\perp \perp 1$,
    hence~$1^\perp=0$ by the zero--one law.
By the orthocomplement law~$a \ovee a^\perp = 1$.
Thus~$a^\perp \ovee a = a \ovee a^\perp = 1$ by partial commutativity.
    Hence~$a = a^{\perp\perp}$ by uniqueness of the orthocomplement.
Now, to prove the zero law,
    note~$a^\perp \ovee a^{\perp\perp}
        = 1 = (a^\perp\ovee a^{\perp\perp})
    \ovee (a^\perp \ovee a^{\perp\perp})^\perp
    = (a^\perp \ovee a^{\perp\perp}) \ovee 0 
        = a^\perp \ovee (a^{\perp\perp} \ovee 0)$
    by the orthocomplement law and partial associativity.
Thus by uniqueness of the orthocomplement~$a^{\perp\perp} = a^{\perp\perp} \ovee 0$,
    hence~$0 \ovee a = a$ by partial commutativity and~$a^{\perp\perp} = a$.
\end{solution}
\begin{solution}{exc-dposet}%
We start with~(D1).  Assume~$b \leq a$. Then~$b \ovee c= a$ for some~$c$.
    By definition~$c = a \ominus b $.
Conversely, assume~$a \ominus b$ is defined.
    Then~$b \ovee (b \ominus a) = a$, hence~$b \leq a$.
To prove~(D2), assume~$a \ominus b$ is defined.
    Then~$a \ominus b \leq b \ovee (a \ominus b) = a$ as desired.
For~(D3) assume~$a \ominus b$ and~$a \ominus (a \ominus b)$
    are defined.
Note~$(a \ominus b) \ovee b = b \ovee (a \ominus b) = a$
    and so~$b = a \ominus (a \ominus b)$.

(D4) is next.  Assume~$a\leq b \leq c$.
    Then~$(c \ominus b) \ovee (b\ominus a)  \ovee a
        = c  =  (c \ominus a) \ovee a$
            and so by cancellation~$
               c \ominus b \leq (c\ominus b) \ovee (b \ominus a)
               = c \ominus a$.
    Note~$ ((c \ominus a) \ominus (c \ominus b))  \ovee (c \ominus b)
        \ovee a = c = (b \ominus a) \ovee a \ovee (c \ominus b)$
    and so by cancellation~$ (c\ominus a) \ominus (c \ominus b) = b\ominus a$.

For the final part of the exercise, assume~$E$ is a poset
    with maximum element~$1$ and a partial operation~$\ominus$
    satisfying~(D1), (D2), (D3) and (D4).
We will show~$E$ is an effect algebra with~$a \ovee b = c \Leftrightarrow
            a = c \ominus b$ and~$0 = 1\ominus 1$.

First we have to show that if~$c_1 \ominus b = c_2 \ominus b$
    for some~$b \leq c_1, c_2$ in~$E$,
    then~$c_1 = c_2$
    so that~$\ovee$ is at most single-valued.
By~(D4)
    we have~$(1 \ominus b) \ominus (1 \ominus c_1)
        = c_1\ominus b$
        and so~$(1 \ominus b) \ominus a =
        (1 \ominus b) \ominus (c_1 \ominus b)
            =  (1\ominus b) \ominus ((1 \ominus b) \ominus (1 \ominus c_1))
            = 1 \ominus c_1$ by (D3).
Similarly~$1 \ominus c_2 = (1 \ominus b) \ominus a$,
        so with another application of (D3) we
        see~$c_1 = 1 \ominus(1\ominus c_1) = 1 \ominus (1 \ominus c_2)
            = c_2$.

The first axiom we will prove is partial commutativity.
    Assume~$a,b \in E$ with~$a \perp b$.
    By definition, there is some~$c\geq b$ with~$c \ominus b = a$.
    By~(D3) we have~$c \ominus a = c \ominus (c \ominus b) = b$
        and so~$b \perp a$
        with~$b \ovee a = c = a\ovee b$ as desired.
Note that we have also shown that~$a \leq a\ovee b$
    and~$(a \ovee b) \ominus a = b$
    for~$a \perp b$.

Next, to prove partial associativity, assume~$a,b,c\in E$
    with~$a \perp b$ and~$a \ovee b \perp c$.
Then~$a \leq a \ovee b \leq (a \ovee b) \ovee c$ and thus
\begin{align*}
    b & \ =\ (a \ovee b) \ominus a \\
    & \ \overset{\mathclap{\mathrm{(D4)}}}{=} \ 
    \bigl( ((a \ovee b) \ovee c) \ominus a \bigr) \,
    \ominus \, \bigl( ((a\ovee b) \ovee c) \ominus (a \ovee b) \bigr) \\
    & \ = \ 
    \bigl( ((a \ovee b) \ovee c) \ominus a \bigr) \,
    \ominus \, c,
\end{align*}
which shows~$b \perp c$
    with~$b \ovee c=  ((a \ovee b) \ovee c) \ominus a$.
From that it immediately follows
    that~$b \ovee c \perp a$ with~$(b\ovee c)\ovee a = (a\ovee b) \ovee c$
        and so~$a \perp b \ovee c$
            with~$(a \ovee b) \ovee c = a \ovee (b \ovee c)$
                by partial commutativity.

We continue with the zero axiom.
Assume~$a \in E$.
Clearly~$a \leq 1 \leq 1$ and
    so~$(1 \ominus a) \ominus (1 \ominus 1) = 1 \ominus a$ by~(D4),
    which shows~$1 \ominus a \perp  1\ominus 1 \equiv 0$
    with~$(1 \ominus a) \ovee 0 = 1\ominus a$.
    The zero law follows by partial commutativity.

Now we prove the orthocomplement law.
Assume~$a \in E$.
By~(D3) we have~$1 \ominus( 1\ominus a) = a$,
    which demonstrates~$a \perp 1\ominus a \equiv a^\perp$
    with~$a \ovee a^\perp = 1$.
Assume~$b \in E$ with~$a \ovee b = 1$.
By definition, we have~$a = 1 \ominus b$ and so by~(D3)
    we see~$b = 1 \ominus (1 \ominus b) = 1 \ominus a \equiv a^\perp$.

To show the last axiom (zero--one), assume~$a \in E$ with~$a\perp 1$.
Then~$1 \leq a\ovee 1 \leq 1$ and so~$a \ovee 1 = 1$.
Hence~$a = (a\ovee 1) \ominus 1 = 1\ominus 1 \equiv 0$.
We have shown~$E$ is an effect algebra.

The conclude this exercise, we show that the original
    order and~$\ominus$ coincide with
    those defined for~$E$ as an effect algebra.
For clarity, write~$\leq_{\mathrm{EA}}$ and~$\ominus_{\mathrm{EA}}$
    for the latter.

To start, assume~$a \leq b$ for some~$a,b \in E$ (as a D-poset).
Then~$b = a \ovee (b \ominus a)$
    as trivially~$b \ominus a = b \ominus a$.
Hence~$a \leq_{\mathrm{EA}} b$.
To prove the converse, assume~$a \leq_{\mathrm{EA}} b$.
Then~$a \ovee c = b$ for some~$c \in E$.
Hence~$a \leq a \ovee c \leq b$.  Indeed~$\leq = \leq_{\mathrm{EA}}$.

Assume~$a,b\in E$.
The element~$a \ominus_{\textrm{EA}} b$ is defined
    if and only if~$b \leq a$,
    which is precisely when~$ a \ominus b$ is defined.
        Furthermore~$b \ovee (a \ominus_{\textrm{EA}} b) = a$
        and so~$
        a \ominus_{\textrm{EA}} b =
        (b \ovee (a \ominus_{\textrm{EA}} b)) \ominus b =
        a \ominus b $.
\end{solution}
\begin{solution}{exc-eamorphism}%
   Assume~$f\colon E \to F$ is an homomorphism between effect algebras.
    \begin{enumerate}
        \item
            We have~$0\perp 1$ and so~$f(0) \perp f(1) = 1$,
                hence~$f(0) = 0$ by the zero--one law.
        \item
            Assume~$a \leq b$ in~$E$.
            Then~$a \ovee c = b$ for some~$c \in E$.
            Hence~$f(b) = f(a \ovee c) = f(a) \ovee f(c)$,
            which shows~$f(a) \leq f(b)$.
        \item
            Assume~$a \ominus b$ is defined (i.e.~$a \geq b$).
            Recall~$(a \ominus b) \ovee b = a$
                and so~$f(a \ominus b) \ovee f(b) = f(a)$.
                Hence~$f(a \ominus b) = f(a) \ominus f(b)$.
        \item
            Assume~$a \in E$.
            In~\sref{exc-dposet} we saw~$1 \ominus a = a^\perp$
            and so by the previous point~$f(a^\perp) = f(1\ominus a) =
                f(1) \ominus f(a) = 1 \ominus f(a) = f(a)^\perp$.
    \end{enumerate}
\end{solution}
\spacingfix{}
\begin{solution}{exc-emonzero}%
    Let~$M$ be an effect monoid with any~$a \in M$.
    Then~$a\odot 1 = a \odot (0 \ovee 1) = (a\odot 0) \ovee (a \odot 1)$.
    Hence by cancellation~$a \odot 0 = 0$.
    Similarly~$0 \odot a = 0$.
\end{solution}
\begin{solution}{emond-lemma-for-conv}%
   Assume~$M$ is an effect monoid with~$a_1,\ldots,a_n,b_1, \ldots, b_n \in M$
        such that~$\bigovee_i a_i = 1$ and~$\bigovee_i a_i \odot b_i = 1$.
Note that for any~$1 \leq i \leq n$
    we have~$a_i \odot b_i \leq a_i$
        as~$(a_i \odot b_i) \ovee a_i \odot b_i^\perp = a_i \odot 1 = a_i$.
Thus
\begin{equation*}
    a_i^\perp
    \ = \ \bigovee_{j \neq i} a_j
    \ \geq \ \bigovee_{j \neq i} a_j \odot b_j
    \ = \ (a_i \odot b_i)^\perp
    \ \geq a_i^\perp
\end{equation*}
    and so~$a_i^\perp = (a_i \odot b_i)^\perp$,
        hence~$a_i = a_i \odot b_i$.
\end{solution}
\begin{solution}{eff-prod-rules}%
   Let~$C$ be an effectus in partial form.
Assume~$f\colon X \to A$ and~$g\colon X \to B$
    are any arrows in~$C$ with~$1 \after f \perp 1 \after g$.
\begin{enumerate}
\item Assume~$[a,b]\colon A+B \to Y$ is any map.
Then~$1 \after a \after f \leq 1 \after f \perp 1 \after g
        \geq 1 \after b \after g$
        and so~$1 \after a \after f \perp 1 \after b \after g$.
    Hence~$a \after f \perp b \after g$
    and so~$[a,b] \after \langle f, g \rangle
        = [a,b] \after ((\kappa_1 \after f) \ovee (\kappa_2 \after g))
        = ([a,b] \after (\kappa_1 \after f)) \ovee ([a,b] \after (\kappa_2 \after g)) = (a \after f) \ovee (b \after g)$.
\item By the previous point and~\sref{cotupl-pcm},
        we have~$1 \after \langle f,g\rangle
            = [1,1] \after \langle f,g\rangle
            = (1 \after f) \ovee (1 \after g)$.
\item Assume~$k\colon A \to A'$ and~$l \colon B \to B'$ are two maps in~$C$.
    Then by the first point of this exercise,
            we have~$(k+l) \after \langle f, g\rangle
                \equiv [\kappa_1 \after k, \kappa_2 \after l] \after
                    \langle f, g \rangle 
                = (\kappa_1 \after k \after f) \ovee
                (\kappa_2 \after l \after g) \equiv \langle 
                k \after f, l \after g \rangle $.
\item
Assume~$k\colon X' \to X$ is some map in~$C$.
Then by PCM-enrichment of~$C$, we
    get~$\langle f, g \rangle \after k 
        \equiv ((\kappa_1 \after f) \ovee (\kappa_2 \after g)) \after k
        = (\kappa_1 \after f\after k) \ovee (\kappa_2 \after g \after k)
        \equiv \langle f \after k, g\after k\rangle$.
\end{enumerate}
\end{solution}
\spacingfix{}
\begin{solution}{exc-jointly-monic-pullback}%
    Assume~$h_1,h_2 \colon Q \to P$
        are two arrows with~$m_1 \after h_1 = m_1 \after h_2$
        and~$m_2 \after h_1 = m_2 \after h_2$.
    Note~$(f \after m_1) \after h_1= g \after m_2 \after h_1
        = (g \after m_2) \after h_2$
        and so by the universal property of a pullback
        there is a unique map~$h\colon Q \to P$
        with~$m_1 \after h = m_1 \after h_1$
            and~$m_2 \after h = m_2 \after h_2$
            and so~$h_1 = h_2 = h$.
\end{solution}
\begin{solution}{pullback-lemma}%
    For the first point,
   assume the left and right inner squares are pullbacks.
To prove the outer square is a pullback,
    assume~$\alpha \colon S \to C$ and~$\beta\colon S \to X$
    are any maps with~$g' \after f' \after \beta = m \after \alpha$.
As the right inner square is a pullback,
    there is a unique map~$\gamma\colon S \to B$
        with~$g \after \gamma = \alpha$
            and~$l \after \gamma = f' \after \beta$.
By this last equality and the fact that the left inner square is
    a pullback, there is a unique map~$\delta \colon S \to A$
    with~$f \after \delta = \gamma$
        and~$k \after \delta = \beta$.
    We will prove that this~$\delta$ is also the mediating map
        for the outer square.
Clearly~$g \after f \after \delta g \after \gamma = \alpha$.
Assume~$\delta'\colon S \to A$
    is some other map with~$g \after f \after \delta' = \alpha$
    and~$k \after \delta' \after \beta$.
Note that~$g \after f \after \delta' = \alpha = g \after \gamma$
    and~$l \after f \after \delta' = f' \after k \after \delta'
        = f' \after \beta = l \after \gamma$.
Thus~$f \after \delta' = \gamma$ as~$g$ and~$l$ are jointly monic by~\sref{exc-jointly-monic-pullback}.
Consequently~$\delta=\delta'$ and so the outer square is indeed a pullback.

For the second point, assume the outer square is a pullback
    and that~$l$ and~$g$ are jointly monic.
We will prove that the left inner square is a pullback.
To this end, assume~$\alpha\colon S \to B$ and~$\beta\colon S \to X$
    are maps with~$l \after \alpha = f' \after \beta$.
Note~$m \after g \after \alpha = g' \after l \after \alpha = g' \after f' \after \beta$
    and thus by the fact that the outer square is a pullback,
        there is a unique~$\gamma\colon S \to A$
        with~$g \after f \after \gamma = g \after \alpha$
            and~$k \after \gamma = \beta$.
    This will also be the unique mediating map
        demonstrating that the left inner square is a pullback.
We have~$l \after f \after \gamma = f' \after k \after \gamma
        = f' \after \beta = l \after \alpha$
            and so by the joint monicity of~$l$ and~$g$,
            we conclude~$f \after \gamma = \alpha$.
    Assume there is some map~$\gamma' \colon S \to A$
        with~$f \after \gamma' = \alpha$ and~$k \after \gamma' = \beta$.
Then clearly~$g \after f \after \gamma' = g \after \alpha$
    and so by the uniqueness of~$\gamma$
    as a mediating map for the outer square,
    we find~$\gamma' = \gamma$.
    Thus the left inner square is a pullback as well.
\end{solution}
\begin{solution}{par-c-coprod}%
    Let~$f\colon X \pto Z$ and~$g\colon Y \pto Z$
        be two partial maps in~$C$;
        that is~$f\colon X \to Z+1$ and~$g \colon Y \to Z+1$.
    Then the cotuple~$[f,g]\colon X+Y \to Z+1$
        (with respect to the coproduct that we assumed to exist)
    is also a partial map~$X+Y \pto Z$.
Note~$[f,g] \hafter \hat\kappa_1
    = [[f,g],\kappa_2]\after\kappa_1\after\kappa_1 = f$
and~$[f,g] \hafter \hat\kappa_2
    = [[f,g],\kappa_2]\after\kappa_1\after\kappa_2 = g$.
Suppose~$h \colon X+Y \pto Z$ is a partial map
    with~$h \hafter \hat\kappa_1 = f$
        and~$h \hafter \hat\kappa_2 = g$.
Then~$h \after \kappa_1 = h \hafter \hat\kappa_1 = f$
and~$h \after \kappa_2 = h \hafter \hat\kappa_2 = g$,
    so~$h = [f,g]$.
    This shows~$X+Y$ is also a coproduct in~$\Par C$ with
        coprojectors~$\hat{\kappa_1}$ and~$\hat{\kappa_2}$.

To show~$0$ is an initial object of~$\Par C$,
    assume~$X$ is some object of~$\Par C$, i.e.~of~$C$.
    As~$0$ is initial in~$C$, there
    is a unique map~$!\colon 0 \to X+1$ in~$C$
    and so there is a unique map~$!\colon 0 \pto X$ in~$\Par C$.
\end{solution}
\begin{solution}{toteff-zero}%
We already know that the initial object~$0$ of~$C$
    is also initial in~$\Par C$.
    We will show that it is also final in~$\Par C$.
To this end, assume~$X$ is some object in~$\Par C$, i.e.~of~$C$.
There is a unique map~$!\colon X\to 1$ in~$C$
    and so there is a unique map~$!\colon X \pto 0$ in~$\Par C$.
    Thus~$0$ is final in~$\Par C$.

For any objects~$X$ and~$Y$ from~$\Par C$,
    the zero map~$0\colon X \to Y$ in~$\Par C$
    is by definition equal to~$i \hafter f$,
    where~$f \colon X \pto 0$ is the unique final map
        and~$i \colon 0 \to \pto Y$ is the unique initial map.
Unfolding definitions, we see~$0 = i \hafter f = [i, \kappa_2]\after f
        = \kappa_2\after !$ as desired.
\end{solution}
\begin{solution}{distinction-part-tot-eff}%
For the first point, assume~$C$ is an effectus in total form
    together with a map~$f\colon A \to 0$.
    We will show that~$f$ is an isomorpism.
    By initiality of~$0$, we know~$f \after !_A = !_0 = \id$.
By~\sref{tot-pullbacks}, the following square is a pullback.
\begin{equation*}
    \xymatrix{
        0 \pullback \ar[d]_{!} \ar@{=}[r]
        & 0 \ar@{=}[d]
        \\
        A \ar[r]_{f}
        & 0
    }
\end{equation*}
As~$\id_0 \after f = f \after \id_A$,
    there is a unique map~$h\colon A \to 0$
        with~$\id \after h = f$ and~$! \after h = \id$.
So~$h = f$ and~$! \after f = ! \after h = \id$.
Thus~$f$ is indeed an isomorphism with inverse~$!$.

For the second point, assume~$C$ is an effectus in total and partial form.
Let~$X$ and~$Y$ be any two objects.
As~$0$ is final in an effectus in total form,
    there are maps~$!_X \colon X \to 0$ and~$!_Y \colon Y \to 0$.
As~$0$ is a strict initial object in an effectus in partial form,
    both~$!_X$ and~$!_Y$ are isomorphisms.
Hence~$!_Y \after !_X^{-1}$ is an isomorphism between~$X$ and~$Y$.
\end{solution}
\begin{solution}{exc-rng-eff}%
   To show~$\op\Rng$ is an effectus in total form,
        we shown the dual axioms for~$\Rng$.
    Clearly, for any two unit rings~$R$ and~$S$,
        their cartesian product~$R \times S$ is a categorical product
        with projectors~$\pi_1(r,s) = r$ and~$\pi_2(r,s)=s$.
    The integer ring~$\Z$ is the initial object of~$\Rng$.
The zero ring (the unique unit ring with a single element~$0=1$)
    is the final object of~$\Rng$.

To show the pushout diagrams corresponding to~\eqref{pullbacks} hold,
    assume we are given unit rings~$R,S,T$ with
    unit-preserving homomorphisms~$\alpha, \beta,\delta$
    that make the outer squares of the following diagrams commute.
\begin{equation*}
    \xymatrix{
    T\\
    & R \times S \ar@{.>}[lu]
    & R \times \Z \ar[l]^{\id\times !}
    \ar@/_1pc/[llu]_{\alpha}
    \\& \Z \times S \ar[u]_{! \times \id}
    \ar@/^1pc/[luu]^{\beta}
    & \Z \times \Z \ar[l]^{\id\times!} \ar[u]_{! \times \id}
    } \qquad
    \xymatrix{
    T\\
    & R \ar@{.>}[lu]
    & \Z \ar[l]^{!}
    \ar@/_1pc/[llu]_{!}
    \\& R \times S \ar[u]_{\pi_1}
    \ar@/^1pc/[luu]^{\delta}
    & \Z \times \Z \ar[l]^{!\times!} \ar[u]_{\pi_1}
    }
\end{equation*}
    We have to show that there are unique dashed arrows (as shown)
        that make these diagrams commute.
    We start with the left diagram.
    Define~$f\colon R\times S \to T$
    by~$f(r,s) \equiv \alpha(r,0) + \beta(0,s)$.
Clearly~$f$ is additive and $f(0,0) = \alpha(0,0) + \beta(0,0) = 0$.
    By assumption~$\alpha(n,m) = \beta(n,m)$ for any~$n,m \in \Z$
        and so~$\alpha(1,0) = \beta(1,0)$ in particular,
        hence~$f(1,1) = \alpha(1,0) + \beta(0,1) = \beta(1,1) = 1$.
It remains to be shown that~$f$ is multiplicative.
    First note that~$\alpha(1,0) \beta(0,1)
    = \beta(1,0)\beta(0,1) = \beta(0,0) = 0$.
    Thus~$\alpha(r,0)\beta(0,s)
            = \alpha(r,0)\alpha(1,0)\beta(0,1)\beta(0,s) = 0$
        for any~$r \in R$ and~$s \in S$
            and similarly~$\beta(0,s)\alpha(r,0) = 0$.
Hence~$f(r,s)f(r',s') = \alpha(r,0)\alpha(r',0) + \beta(0,s)\beta(0,s')
        = \alpha(rr',0) + \beta(0,ss') = f(rr', ss')$.
        We have shown~$f$ is a unit-preserving ring homomorphism.
Clearly~$f(r, m) = \alpha(r, 0) + \beta(0, m)
            = \alpha(r,0) + \alpha(0,m) = \alpha(r,m)$
                for any~$r \in R$ and~$m \in \Z$
                hence~$\alpha = f \after (\id \times !)$.
    Similarly~$\beta = f \after (! \after \id)$.
Assume~$f'\colon R\times S \to T$ is any unit-preserving
    ring homomorphism with~$\alpha = f' \after (\id \times !)$
        and~$\beta = f' \after (! \after \id)$.
Then clearly~$f'(r,0) = \alpha(r,0)$ and~$f'(0,s) = \beta(0,s)$,
    so~$f'(r,s) = f'(r,0) + f'(0,s) = \alpha(r,0) + \beta(0,s) = f(r,s)$.
This shows the left square above is indeed a pushout.

We continue with the diagram on the right.
By assumption~$\delta(n,m) = n$ for any~$n,m \in \Z$
    hence~$\delta(0,s) = \delta(0,s)\delta(0,1) = 0$ for any~$s\in S$.
Thus~$\delta(r,s) = \delta(r,0) + \delta(0,s) = \delta(r,0)$
    for any~$r \in R$ and~$s \in S$.
Define~$g\colon R \to T$ by~$g(r) = \delta(r,0)$.
Clearly~$g$ is additive, multiplicative and~$g(0,0) = 0$.
Furthermore~$g(1) = \delta(1,0) = \delta(1,1) = 1$,
    so it is a unit-preserving ring homomorphism.
It is easy to see that~$g \after \pi_1 = \delta$ and that~$g$
    is the unique such unit-preserving ring homomorphism.
We have shown the right diagram is a pushout too.

To show~$\op\Rng$ is an effectus in total form,
    it only remains to be shown that~$
        \langle \pi_1, \pi_2,\pi_2\rangle,
        \langle \pi_2, \pi_1,\pi_2\rangle
        \colon \Z\times\Z
            \to \Z\times\Z\times\Z$ are jointly epic.
    So assume~$f,g\colon \Z\times\Z\times\Z \to R$
        are two unit-preserving ring homomorphisms
        for which we have~$f \after \langle \pi_1, \pi_2,\pi_2\rangle  =
        g \after \langle \pi_1, \pi_2,\pi_2\rangle$
        and~$f \after \langle \pi_2, \pi_1,\pi_2\rangle =
        g \after \langle \pi_2, \pi_1,\pi_2\rangle$.
    By the first equality~$f(k,0,0) = g(k,0,0)$
        and by the second~$f(0,l,0) = g(0,l,0)$
            for any~$k,l \in \Z$.
In particular~$f(0,0,1) = 1-f(1,0,0) - f(0,1,0)
    = 1 - g(1,0,0) - g(0,1,0) = g(0,0,1)$
        and so~$f(0,0,m) = m f(0,0,1) = m g(0,0,1) = g(0,0,m)$
            for any~$m \in \Z$.
    Putting it all together: $f(k,l,m) 
            = f(k,0,0) = f(0,l,m) + f(0,0,m)
            = g(k,0,0) = g(0,l,m) + g(0,0,m) = g(k,l,m)$.
    So~$f=g$ and so we have shown the required joint epicity.
We have shown~$\op\Rng$ is indeed an effectus in total form.

We continue with the two additional points.
Let~$p \colon \Z\times\Z \to R$ be any predicate on~$R$.
Then~$p(1,0) \in R$ with~$p(1,0)^2 = p(1,0^2) = p(1,0)$ and so~$p(1,0)$
    is an idempotent.
    Furthermore~$p(0,1) = p(1,1) - p(1,0) = 1 - p(1,0)$,
        so~$p$ is fixed by the idempotent~$p(1,0)$.
Let~$e \in R$ be any idempotent.
It is easy to see that~$p(n,m) \equiv ne + m(1-e)$ is a unit-preserving
    ring homomorphism with~$p(1,0) = e$.  Thus the predicates
        on~$R$ correspond to its idempotents.
    This also shows that the set of scalars of~$\op\Rng$
        is the two-element effect monoid~$2$
        as~$\Z$ has exactly two idempotents: $0$ and~$1$.

Assume~$p \perp q$ for two predicates~$p,q \colon \Z\times\Z \to R$
    on~$R$.
Then, by definition, $p \perp q$ if
        there is a~$b\colon \Z\times\Z\times\Z\to R$
    with~$b(n,m,m) = p(n,m)$
    and~$b(m,n,m) = q(n,m)$ for all~$n,m \in \Z$
In particular~$p(1,0)q(1,0) = b(1,0,0)b(0,1,0) = 0$,
    thus the idempotents corresponding to~$p$ and~$q$ are orthogonal.
Conversely, if~$p(1,0)$ and~$q(1,0)$ are orthogonal,
    then~$b(k,l,m) = k p(1,0) + l q(1,0) + m (1-p(1,0) - q(1,0))$
        defines a unit-presrving ring homomorphism
        that shows~$p \perp q$ with~$(p \ovee q)(1,0)
            = b(1,1,0) = p(1,0) + q(1,0)$.
From this it is also clear that~$p^\perp(1,0) = 1 - p(1,0)$.

There are many examples to show unit-preserving ring homomorphisms
    need not be fixed on their value on idempotents
        (and thus that~$\op\Rng$ does not have separting predicates.)
For instance, let~$R$ denote the unit ring of continuous real-valued
    functions on the unit interval~$[0,1]$.
This ring has only two idempotents: the functions that are constant~$0$
    and~$1$, which are also the zero and unit element (respectively)
    of the ring.
For every~$x \in [0,1]$
    the map~$\delta_x \colon R \to \R$ given by~$\delta_x(f) =f(x)$
    is a unit-preserving ring homomorphism.
Clearly~$\delta_x(0) = 0 = \delta_y(0)$ and~$\delta_x(1) = 1 = \delta_y(1)$
    for any~$x,y \in [0,1]$ with~$x\neq y$, but not~$\delta_x = \delta_y$.

Finally, we treat the second and last point of the exercise.
Reasoning towards contradiction,
    let~$f\colon \Z_2 \to \Z$ be any unit-preserving ring homomorphism.
Then~$0 = f(0) = f(1+1) = f(1) + f(1) = 1 + 1 = 2$, which is absurd.
Thus there is no such homomorphism.
However, there are two unit-preserving ring
    homomorphisms~$\Z \times \Z \to \Z_2$
    corresponding to the idempotents~$0$ and~$1$ in~$\Z_2$.
Hence the states cannot be separating (for otherwise
    there could at most be a single unit-preserving
    ring homomorphism~$\Z\times\Z\to\Z_2$.)
\end{solution}
\begin{solution}{exc-dm-effectus}%
We will first show that~$\mathcal{D}_M$ is a functor.
To start, clearly
\begin{equation*}
    \mathcal{D}_M(\id)(p)(x) 
        \ \equiv\  \bigovee_{y; \id(y)=x} p(y)
        \ =\ p(x)
\end{equation*}
        and so~$\mathcal{D}_M(\id) = \id$.
Let~$f\colon X \to Y$ and~$g\colon Y \to Z$ be given.
    Then
    \begin{align*}
        \mathcal{D}_M (g \after f) (p) (z)
        & \ \equiv\ \bigovee_{x;\ g(f(x))=z} p(x)\\
        & \ =\ \bigovee_{y;\ g(y) = z}
               \bigovee_{x;\ f(x) = y} p(x) \\
        & \ =\ \bigovee_{y;\ g(y) = z}
                \bigl((\mathcal{D}_M f)(p)\bigr)(y) \\
        & \ =\ \bigl( \, (\mathcal{D}_M g) \bigl(
                (\mathcal{D}_M f)(p)\bigr) \, \bigr)(z)
    \end{align*}
    and so~$\mathcal{D}_M (g \after f)
        = (\mathcal{D}_M g)\after (\mathcal{D}_M f)$.

Next, we show that~$\eta$ is a natural transformation.
Let~$f\colon X \to Y$ be any map.
Then for any~$x_0 \in X$, we have
\begin{align*}
    ((\mathcal{D}_M f) \after \eta_X) (x_0)
        & \ = \ \bigovee_{x;\ f(x)=y} (\eta_X(x_0))(x) \\
        & \ = \ \bigovee_{x;\  \substack{f(x)=y \\ x = x_0}} 1 \\
        & \ = \ (\eta_Y (f(x_0)))
\end{align*}
and so~$\mathcal{D}_M f \after \eta_X = \eta_Y \after f$,
    i.e.~$\eta$ is a natural transformation~$\id \Rightarrow \mathcal{D}_M$.

To show~$\mu$ is a natural transformation,
    assume~$f\colon X \to Y$ is some map, $y \in Y$ and~$\Phi \in
        \mathcal{D}_M\mathcal{D}_M X$.
Then
\begin{align*}
    (\mathcal{D}_M f) ( \mu_X(\Phi) ) (y)
    & \ = \ \bigovee_{x;\ f(x)=y} \mu_X(\Phi)(x) \\
    & \ = \ \bigovee_{x;\ f(x)=y}
            \bigovee_{p} \Phi(p) \odot p(x) \\
    & \ = \ \bigovee_{p}
                \bigovee_{x;\ f(x)=y}
            \Phi(p) \odot p(x) \\
    & \ = \ \bigovee_{p}
            \Phi(p) \odot 
                \bigovee_{x;\ f(x)=y} p(x)  \\
    & \ = \ \bigovee_{p}
            \Phi(p) \odot 
                (\mathcal{D}_M f)(p)(y) \\
    & \ = \ \bigovee_{q}
                \bigovee_{p;\ (\mathcal{D}_M f)(p) = q}
            \Phi(p) \odot q(y) \\
    & \ = \ \bigovee_{q}
               \Bigl( \bigovee_{p;\ (\mathcal{D}_M f)(p) = q}
            \Phi(p)\Bigr) \odot q(y) \\
    & \ = \ \bigovee_{q}
                    (\mathcal{D}_M \mathcal{D}_M f) (\Phi)(q) \odot q(y) \\
    & \ = \ \mu_Y (( \mathcal{D}_M \mathcal{D}_M f)(\Phi)) (y)
\end{align*}
and so~$(\mathcal{D}_M f) \after \mu_X
= \mu_Y \after (\mathcal{D}_M \mathcal{D}_M f)$,
 hence~$\mu$ is a natural transformation~$\mathcal{D}_M \mathcal{D}_M
    \Rightarrow \mathcal{D}_M$.
    We continue with the monad laws. Assume~$X$ is a set, $x \in X$
        and~$\aleph \in \mathcal{D}_M \mathcal{D}_M \mathcal{D}_M X$.
        Then
\begin{align*}
\mu_X ( \mu_{\mathcal{D}_M X} ( \aleph))(x)
    & \ = \ \bigovee_p (\mu_{\mathcal{D}_M X}) (\aleph) (p) \odot p(x) \\
    & \ = \ \bigovee_p \bigovee_\Phi \aleph(\Phi) \odot \Phi(p) \odot p(x) \\
    & \ = \ \bigovee_\Phi  \aleph(\Phi) \odot \bigovee_p \Phi(p) \odot p(x) \\
    & \ = \ \bigovee_\Phi  \aleph(\Phi) \odot \mu_X(\Phi)(x) \\
    & \ = \ \bigovee_p \bigovee_{\Phi; \ \mu_X(\Phi)=p}
        \aleph(\Phi) \odot p(x) \\
    & \ = \ \bigovee_p  (\mathcal{D}_M \mu_X) (\aleph) (p)\odot p(x) \\
    & \ = \ \mu_X ((\mathcal{D}_M \mu_X)(\aleph)) (x)
\end{align*}
    and so~$\mu \after \mu_{\mathcal{D}_M} = \mu \after (\mathcal{D}_M )$.
For any~$x \in X$ and $p \in \mathcal{D}_M X$, we have
\begin{equation*}
    \mu_X (\eta_{\mathcal{D}_M X}(p))(x)
    \ = \ \bigovee_q \eta_{\mathcal{D}_M X}(p)(q) \odot q(x) \ = \ p(x)
\end{equation*}
and so~$\mu \after \eta_{\mathcal{D}_M} = \id$.
For the final monad law, assume~$x \in X$ and~$p \in \mathcal{D}_M X$.
Then
\begin{align*}
    \mu_X( (\mathcal{D}_M \eta_X)(p))(x)
    & \ = \ \bigovee_q (\mathcal{D}_M \eta_X) (p)(q) \odot q(x) \\
    & \ = \ \bigovee_q \bigovee_{y; \ \eta_X(y)=q}
            p(y) \odot q(x) \\
    & \ = \ \bigovee_y p(y) \odot \eta_X (y)(x) \\
    & \ = \ p(x)
\end{align*}
and so~$\mu \after \mathcal{D}_M \eta = \id$.
We have shown~$(\mathcal{D}_M, \eta, \mu)$ is a monad.

Our next project is to show that~$\Kl \mathcal{D}_M$ is an effectus
    in total form.
It is easy to see that the coproduct~$X + Y$ from~$\mathsf{Set}$
    is also a coproduct in~$\mathcal{D}_M$
        with coprojectors~$\eta\after\kappa_i$,
            where~$\kappa_i$ are the coprojectors for~$X+Y$ in~$\mathsf{Set}$.
(In fact, for any category~$C$ and monad~$T$,
    the inclusion functor~$K\colon C \to \Kl T$
        given by~$K X = X$ and~$K f = \eta \after f$
        is a left adjoint and so~$K$ preserves colimits.)
The empty set is also the initial object of~$\Kl \mathcal{D}_M$.
As~$\mathcal{D}_M 1 \cong 1$,
    the category~$\Kl \mathcal{D}_M$
        has as final object the one-element set~$1$.

Now we will show that a square such as that on the
    left of \eqref{pullbacks} is a pullback.
So assume~$\alpha \colon Z \to \mathcal{D}_M (X+1)$
    and~$\beta\colon Z \to \mathcal{D}_M (1+Y)$
        are maps with~$(\hat !+\hat\id) \hafter \alpha = (\hat\id+\hat!) \hafter \beta$,
        where~$\hafter$ denotes the composition in~$\Kl \mathcal{D}_M$
        and~$\hat f \equiv \eta \after f$ the Kleisli embedding.
Assume~$z \in Z$.
By assumption~$ \alpha(z) (\kappa_2*)
                = ((\hat! + \hat\id) \hafter \alpha) (z) (\kappa_2 *)
                = ((\hat\id + \hat!) \hafter \beta) (z) (\kappa_2 *)
                = \bigovee_y \beta(z)(\kappa_2 y) $,
                    where~$*$ is the unique element of~$1$.
We want to define~$\delta\colon Z \to \mathcal{D}_M (X+Y)$
    by~$\delta(z)(\kappa_1 x) = \alpha(z)(\kappa_1 x)$
        and~$\delta(z)(\kappa_2 y) = \beta(z)(\kappa_2 y)$,
            but first need to check the image of~$\delta(z)$ sums to~$1$.
Indeed~$\bigovee_w \delta(z)(w)
    = \bigovee_y \beta(z)(\kappa_2 y)
    \ovee \bigovee_x \alpha(z)(\kappa_1 x)
    = \alpha(z)(\kappa_2*)
    \ovee \bigovee_x \alpha(z)(\kappa_1 x)
    = 1$.
Clearly~$\alpha = (\hat\id+\hat!) \hafter \delta $
    and~$\beta = (\hat!+\hat\id)\hafter \delta$.
Suppose~$\delta'\colon Z \to \mathcal{D}_M (X+Y)$
    is some map with~$\alpha = (\hat\id+\hat!) \hafter \delta'$
        and~$\beta = (\hat!+\hat\id) \hafter \delta'$.
    Then~$\delta'(\kappa_1 x)
        = ((\hat\id + \hat!) \hafter \delta')(\kappa_1 x)
        = \alpha(\kappa_1 x)$ for any~$x \in X$
        and similarly~$\delta'(\kappa_2 y)
        = \beta(\kappa_2 y)$ for any~$y \in Y$, so~$\delta' = \delta$.
We have shown that the square on the left of~\eqref{pullbacks}
    is indeed a pullback in~$\Kl \mathcal{D}_M$.

To show that the square on the right of~\eqref{pullbacks}
    is also a pullback, assume~$\alpha\colon Z \to \mathcal{D}_M 1$
        and~$\beta\colon Z \to \mathcal{D}_M (X+Y)$
            are maps with~$(\hat!+\hat!) \hafter \beta = \hat\kappa_1 \hafter \alpha$.
    Assume~$z \in Z$ and~$y \in Y$.
    Then~$\beta(z)(\kappa_2 y) 
                \leq ((\hat!+\hat!) \hafter \beta)(z)(\kappa_2 *)
                = (\hat\kappa_1 \hafter \alpha)(z)(\kappa_2 *)
                = 0 $.
This allows us to define~$\delta \colon Z \to \mathcal{D}_M X$
        by~$\delta(z)(x) = \beta(\kappa_1 x)$
        as~$\bigovee_x \delta(z)(x)
                = \bigovee_x \beta(\kappa_1 x)
                = \bigovee_x \beta(\kappa_1 x)
                 \ovee \bigovee_y \beta(\kappa_2 y) = 1$.
    Clearly~$\hat\kappa_1 \hafter \delta = \beta$
        and~$\hat! \hafter \delta = \alpha$.
Suppose~$\delta' \colon Z \to \mathcal{D}_M X$ is some map
    with~$\hat\kappa_1 \hafter \delta' = \beta$
        and~$\hat! \hafter \delta' = \alpha$.
    Then~$\delta'(z)(x) = (\hat\kappa_1 \hafter \delta')(\kappa_1 x)
        = \beta(z)(\kappa_1 x) = \delta(z)(x)$
        and so~$\delta = \delta'$.
We have shown that the square on the right of~\eqref{pullbacks}
    is a pullback in~$\Kl \mathcal{D}_M$.

It only remains to be shown that the
    maps~$
    [\hat\kappa_1, \hat\kappa_2, \hat\kappa_2],
    [\hat\kappa_2, \hat\kappa_1, \hat\kappa_2] \colon
        1+1+1 \to \mathcal{D}_M (1+1) $
        are jointly monic in~$\Kl \mathcal{D}_M$.
To this end,
assume~$f_1,f_2 \colon Z \to \mathcal{D}_M (1+1+1)$
are given with~$
    [\hat\kappa_1, \hat\kappa_2, \hat\kappa_2] \hafter f_1 =
    [\hat\kappa_1, \hat\kappa_2, \hat\kappa_2] \hafter f_2$
    and~$[\hat\kappa_2, \hat\kappa_1, \hat\kappa_2] \hafter f_1 =
    [\hat\kappa_2, \hat\kappa_1, \hat\kappa_2] \hafter f_2$.
Assume~$z \in Z$.
Clearly
\begin{align*}
f_1 (z)(\kappa_1 *)
            &\ =\  ([\hat \kappa_1, \hat\kappa_2, \hat\kappa_2] \hafter  f_1)
                (z)(\kappa_1 *)\\
            &\ =\  ([\hat \kappa_1, \hat\kappa_2, \hat\kappa_2] \hafter  f_2)
                (z)(\kappa_1 *)\\
            &\ =\  f_2(z)(\kappa_1 *)
\end{align*}
    and
\begin{align*}
    f_1 (z)(\kappa_2 *)
           &\ =\ ([\hat \kappa_2, \hat\kappa_1, \hat\kappa_2] \hafter  f_1)
                (z)(\kappa_1 *)\\
           &\ =\ ([\hat \kappa_2, \hat\kappa_1, \hat\kappa_2] \hafter  f_2)
                (z)(\kappa_1 *)\\
           &\ =\ f_2(z)(\kappa_2 *).
\end{align*}
So
\begin{align*}
f_1(z)(\kappa_3 *)
       &\ =\ \bigl(f_1(z)(\kappa_1 *) \ovee 
        f_1(z)(\kappa_2 *)\bigr)^\perp\\
       &\ =\ \bigl(f_2(z)(\kappa_1 *) \ovee 
        f_2(z)(\kappa_2 *)\bigr)^\perp\\
            &\ =\  f_2(z)(\kappa_3 *),
\end{align*}
hence~$f_1 = f_2$.
We have shown that~$\Kl \mathcal{D}_M$ is an effectus.

A scalar in~$\Kl \mathcal{D}_M$
    corresponds to a map~$\lambda\colon 1 \to \mathcal{D}_M (1+1)$,
        which corresponds in turn to an element of~$M$.
Unfolding definitions, it is straight-forward to see
    that the effect monoid structure defined on the scalars
        is the same as that of~$M$.
\end{solution}
\begin{solution}{aconv-cong}%
We start with the surjectivity of~$q$, $\mathcal{D}_M q$
    and~$\mathcal{D}_M \mathcal{D}_M q$.
By definition~$q$ is clearly surjective.
Pick a section~$r\colon X/_\sim \to X$ of~$q$;
    i.e.~$q \after r = \id$.
Then~$\mathcal{D}_M r$
    and $\mathcal{D}_M \mathcal{D}_M r$
    are sections of~$\mathcal{D}_M q$ and~$\mathcal{D}_M \mathcal{D}_M q$
    respectively, hence~$\mathcal{D}_M q$ and~$\mathcal{D}_M \mathcal{D}_M q$
        are surjective.

To prove the second point,
    assume~$\sim$ is a congruence.
Assume~$\varphi, \psi \in \mathcal{D}_M X$ with~$
    (\mathcal{D}_M q)(\varphi) =
    (\mathcal{D}_M q)(\psi) $.
    By definition~$\varphi \sim \psi$ and~$h(\varphi) \sim h(\psi)$
        by assumption that~$\sim$ is a congruence.
    Thus~$q \after h (\varphi) = q \after h (\psi)$.
Thus together with the surjectivity of~$\mathcal{D}_M q$,
    there is a unique map~$h_\sim \colon \mathcal{D}_M X/_\sim \to
        X/_\sim$ fixed by~$h_\sim \after (\mathcal{D}_M q) = q \after h$.
To prove the converse, assume there is a
    map~$h_\sim \colon \mathcal{D}_M X/_\sim \to
        X/_\sim$ with~$h_\sim \after (\mathcal{D}_M q) = q \after h$
            and assume~$\varphi \sim \psi$
            for some~$\varphi,\psi \in \mathcal{D}_M X$.
    Then~$q(h(\varphi))
    = h_\sim (\mathcal{D}_M q(\varphi))
    = h_\sim (\mathcal{D}_M q(\psi))
    = q(h(\psi))$ and so~$h(\varphi) \sim h(\psi)$,
    which shows that~$\sim$ is a congruence.

For the third point, assume~$\sim$ is a congruence.
Then
\begin{alignat*}{3}
    h_\sim \after
    \mathcal{D}_M h_\sim  \after
    \mathcal{D}_M \mathcal{D}_M h_\sim 
&\ =\ h_\sim \after
    \mathcal{D}_M (h_\sim \after \mathcal{D}_M q)\\
&\ =\ h_\sim \after
    \mathcal{D}_M (q \after h) &\quad & \text{by the second point}\\
    &\ =\ q \after h \after \mathcal{D}_M h && \text{idem}\\
    &\ =\ q \after h \after \mu && \text{as $(X,h)$ is a.conv.}\\
    &\ =\ h_\sim \after \mathcal{D}_M q \after \mu && \text{by the second point} \\
    &\ =\ h_\sim \after \mu \after
    \mathcal{D}_M \mathcal{D}_M q
    && \text{by naturality~$\mu$.}
\end{alignat*}
    Hence by surjectivity of~$\mathcal{D}_M \mathcal{D}_M q$,
        we find~$h_\sim \after \mathcal{D}_M h_\sim
            = h_\sim \after \mu$.
    Using naturality of~$\eta$,
        the second point
            and the fact that~$(X,h)$ is an abstract convex set
            (in that order),
            we find~$h_\sim \after \eta \after q
        = h_\sim \after \mathcal{D}_M q \after \eta
        = q \after h \after \eta
        = q $ and so~$h_\sim \after \eta = \id$ by surjectivity of~$q$.
We have shown~$(X/_\sim, h_\sim)$ is an abstract~$M$-convex set.
    The equality~$q \after h = h_\sim \after \mathcal{D}_M q$
        from the second point shows that~$q$ is
        an~$M$-affine map.
\end{solution}
\begin{solution}{affine-kernel-cong}%
To show~$\sim$ is a congruence,
    assume~$\varphi \sim \psi$ for~$\varphi,\psi \in \mathcal{D}_M X$.
Write~$q \colon X \to X/_\sim$
    for the quotient map for~$\sim$
        and~$f_\sim \colon X/_\sim \to X$ for
        the unique map with~$f_\sim \after q =f$.
    Note that~$
    \mathcal{D}_M q (\varphi)
    = \mathcal{D}_M q (\psi)$ and so~$f(h(\varphi)) = 
    h (\mathcal{D}_M f(\varphi)) =
    h (\mathcal{D}_M f_\sim (\mathcal{D}_M q (\varphi))) =
    h (\mathcal{D}_M f_\sim (\mathcal{D}_M q (\psi))) =
    f(h(\psi))$
    and so~$h(\varphi) \sim h(\psi)$, which shows~$\sim$
    is a congruence.
\end{solution}
\begin{solution}{least-conv-cong}%
Concerning the first point:
clearly~$h (\eta(h(\psi))) = h(\psi)$
    and so~$(\, \eta(h(\psi)), \, \psi\,)$
    is a derivation of~$\eta(h(\psi)) \approx \psi$.
Assume~$\varphi \approx \psi$.
    Then~$\eta(h(\varphi)) \approx \varphi \approx \psi \approx \eta(h(\psi))$
        and so~$h(\varphi) \sim h(\psi)$.

We continue with the second point.
    Suppose we are given~$\chi_1, \ldots, \chi_n, \varphi,\psi \in \mathcal{D}_M X$
    and~$\lambda_0, \ldots, \lambda_n \in M$
        with~$\bigovee_i \lambda_i = 1$.
To cover the first base case in the definition
    of a derivation, assume~$h(\varphi) = h(\psi)$.
Then
\begin{align*}
    (h \after \mu) \Bigl(\lambda_0 \ket{\psi} \ovee \bigovee_{j=1}^n \lambda_j \ket{\chi_j}\Bigr)
    & \ = \ 
    (h \after \mathcal{D}_M h) \Bigl(\lambda_0 \ket{\psi} \ovee \bigovee_{j=1}^n \lambda_j \ket{\chi_j}\Bigr) \\
    & \ = \ 
    h  \Bigl(\lambda_0 \ket{h(\psi)} \ovee \bigovee_{j=1}^n \lambda_j \ket{h(\chi_j)}\Bigr) \\
    & \ = \ 
    h  \Bigl(\lambda_0 \ket{h(\varphi)} \ovee \bigovee_{j=1}^n \lambda_j \ket{h(\chi_j)}\Bigr) \\
    & \ = \ 
    (h \after \mu) \Bigl(\lambda_0 \ket{\varphi} \ovee \bigovee_{j=1}^n \lambda_j \ket{\chi_j}\Bigr)
\end{align*}
and so
    $\mu \bigl(\lambda_0 \ket{\psi} \ovee \bigovee_{j=1}^n \lambda_j \ket{\chi_j}\bigr)
     \approx  
    \mu \bigl(\lambda_0 \ket{\varphi} \ovee \bigovee_{j=1}^n \lambda_j \ket{\chi_j}\bigr)$.
Next, for the other base case in the definition of a derivation,
    assume~$\varphi \equiv \bigovee_i^m \mu_i \ket{x_i}$
    and~$\psi \equiv \bigovee_i^m \mu_i \ket{y_i}$
    with~$x_i \mathrel{R^*} y_i$ for~$1 \leq i \leq m$.
Then
\begin{align*}
    \mu \Bigl(\lambda_0 \ket{\psi} \ovee \bigovee_{j=1}^n \lambda_j \ket{\chi_j}\Bigr)
    & \ = \  
    \bigovee_{i=1}^m \lambda_0 \odot \mu_i \ket{y_i}
    \ovee \bigovee_{j=1}^n \bigovee_x \lambda_j \odot \chi_j(x) \ket{x}\\
    \mu \Bigl(\lambda_0 \ket{\varphi} \ovee \bigovee_{j=1}^n \lambda_j \ket{\chi_j}\Bigr)
    & \ = \  
    \bigovee_{i=1}^m \lambda_0 \odot \mu_i \ket{x_i}
    \ovee \bigovee_{j=1}^n \bigovee_x \lambda_j \odot \chi_j(x) \ket{x}.
\end{align*}
This shows that~$\mu \bigl(\lambda_0 \ket{\psi} \ovee \bigovee_{j=1}^n \lambda_j \ket{\chi_j}\bigr)
     \approx  
    \mu \bigl(\lambda_0 \ket{\varphi} \ovee \bigovee_{j=1}^n \lambda_j \ket{\chi_j}\bigr)$ by the second base case of a derivation.
Thus by induction~$\mu \bigl(\lambda_0 \ket{\psi} \ovee \bigovee_{j=1}^n \lambda_j \ket{\chi_j}\bigr)
     \approx  
    \mu \bigl(\lambda_0 \ket{\varphi} \ovee \bigovee_{j=1}^n \lambda_j \ket{\chi_j}\bigr)$ whenever merely~$\varphi \approx \psi$.

For the third point, we need some preparation.
Pick representatives~$R \subseteq X$ of~$\sim$:
    for every~$x \in X$ there is a unique~$r_x \in R$ with~$r_x\sim x$.
    By definition~$\eta(r_x) \approx \eta(x)$.
    Let~$\varphi \equiv \bigovee^n_{i=1} \lambda_i\ket{x_i}  \in \mathcal{D}_M X$ given.
By the previous point~
\begin{align*}
    \varphi & \ = \ 
    \mu \Bigl( \lambda_1 \ket{\eta(x_1)}  \ovee \bigovee^n_{i=2} \lambda_i\ket{\eta(x_i)}  \Bigr) \\
    & \ \approx \ 
    \mu \Bigl( \lambda_1 \ket{\eta(r_{x_1})}  \ovee \bigovee^n_{i=2} \lambda_i\ket{\eta(x_i)}  \Bigr) \\
    & \ = \ 
    \mu \Bigl( \lambda_2 \ket{\eta(x_2)}  \ovee \lambda_1 \ket{r_{x_1}}  \ovee  \bigovee^n_{i=3} \lambda_i\ket{\eta(x_i)}  \Bigr) \\
    & \ \approx \ 
    \mu \Bigl( \lambda_2 \ket{\eta(r_{x_2})}  \ovee \lambda_1 \ket{r_{x_1}}  \ovee  \bigovee^n_{i=3} \lambda_i\ket{\eta(x_i)}  \Bigr) \\
    & \qquad\qquad  \vdots\\
    & \ \approx \ \bigovee_{i=1}^n \lambda_i \ket{r_{x_i}} \\
    & \ \approx \ \bigovee_{r \in R} \Bigl(\bigovee_{x \sim r} \varphi(x) \Bigr)\ket{r}.
\end{align*}
Now, assume~$\varphi \sim \psi$.
Then~$\bigovee_{x \sim x_0} \varphi(x)
    = \bigovee_{x \sim x_0} \psi(x)$ for any~$x_0 \in X$
        and so
\begin{equation*}
\varphi \ \approx  \ 
                \bigovee_{r \in R} \Bigl( \bigovee_{x \sim r} \varphi(x) \Bigr) \ket{r}
              \ = \  \bigovee_{r \in R} \Bigl( \bigovee_{x \sim r} \psi(x) \Bigr) \ket{r}
            \ \approx \  \psi.
\end{equation*}

Clearly~$\sim$ is an equivalence relation on~$X$.
To show it is a congruence, assume~$\varphi \sim \psi$.
By the previous, we know~$\varphi \approx \psi$
    and so~$h(\varphi) \sim h(\psi)$ by the first point,
        which shows~$\sim$ is indeed a congruence.
Furthermore~$R$ is contained in~$\sim$;
indeed, if~$x \mathrel{R} y$,
    then~$(\eta(x), \eta(y))$ is a derivation of~$\eta(x) \approx \eta(y)$
    and so~$x \sim y$.

To show~$\sim$ is the smallest congruence containing~$R$,
    assume~$S \subseteq X^2$ is some congruence of~$X$
        with~$R \subseteq S$.
It is sufficient to show that~$\varphi \approx \psi$
    implies~$h(\varphi) \mathrel{S} h(\psi)$.
Indeed if this is the case and~$x \sim y$,
    then by definition~$\eta(x) \approx \eta(y)$
    and so~$x = h(\eta(x)) \mathrel{S} h(\eta(y)) = y$.

We will prove that~$\varphi \approx \psi$
    implies~$h(\varphi) \mathrel{S} h(\psi)$
    by induction over the definition of~$\approx$.
For the first base case in the definition of~$\approx$,
    assume~$\varphi \equiv \bigovee^n_{i=1} \lambda_i \ket{x_i}$
and~$\psi \equiv \bigovee^n_{i=1} \lambda_i \ket{y_i}$
for some~$(x_1,y_1), \ldots, (x_n,y_n) \in R^*$.
Then~$x_i \mathrel{S} y_i$
    and so~$\varphi \mathrel{S} \psi$,
    hence~$h(\varphi) \mathrel{S} h(\psi)$.
For the other base case, suppose~$h(\varphi) = h(\psi)$
    for some~$\varphi,\psi \in \mathcal{D}_M X$.
Then clearly~$h(\varphi) \mathrel{S} h(\psi)$.
Thus, by induction~$h(\varphi) \mathrel{S} h(\psi)$
    whenever~$\varphi \approx \psi$. 
\end{solution}
\begin{solution}{n-times-one-aconvm}%
As the one element set is final in~$\mathsf{Set}$,
    there is a unique map~$!\colon \mathcal{D}_M 1 \to 1$.
        This turns~$\mathcal{D}_M 1$ into an abstract~$M$-convex set:
        ${!} \after \mathcal{D}_M {!} = {!} = {!} \after \mu$
        and~${!} \after \eta = {!} = \id$.
Note that there is a single formal~$M$-convex combination over~$1$
    and so~$\mathcal{D}_M 1 \cong 1$.
As~$\AConvM$ is the Eilenberg--Moore category
    of the monad~$\mathcal{D}_M$,
    it follows that the functor~$F\colon \mathsf{Set} \to \AConvM$
    given by~$F(X) = (\mathcal{D}_M X, \mu)$
    and~$Ff = \mathcal{D}_M f$ is a left adjoint
    and so preserves coproducts.
    Hence~$\mathcal{D}_M \{1,\ldots,n\}
                \cong F(n\cdot 1)
                \cong n \cdot F(1)
                \cong n \cdot 1$.
\end{solution}
\begin{solution}{exc-divisoid-basics}%
    Clearly~$0 \leq \rfrac{0}{0}$
        and~$0\odot 0 = 0$.
        Thus, by uniqueness of~$\rfrac{0}{0}$,
            we have~$\rfrac{0}{0}=0$.
For any~$a \in M$,
    we have~$\rfrac{a}{1} = 1 \odot \rfrac{a}{1} = a$
        and so in particular~$\rfrac{1}{1}=1$.
Next, note~$\rfrac{a}{a} \odot \rfrac{a}{a} \leq \rfrac{a}{a}$
    and~$a \odot \rfrac{a}{a} \odot \rfrac{a}{a} = a \odot \rfrac{a}{a} = a$
        and so by uniqueness of~$\rfrac{a}{a}$,
            we have~$\rfrac{a}{a}\odot \rfrac{a}{a} = \rfrac{a}{a}$.
For the final equation of the first point, assume~$b \in M$.
    As~$a \odot b \leq a$ we know~$\rfrac{a\odot b}{a}$ is defined.
Note~$\rfrac{a}{a} \odot b \leq \rfrac{a}{a}$
    and~$a \odot \rfrac{a}{a} \odot b = a \odot b$.
    So by uniqueness of~$ \rfrac{a\odot b}{a}$
    we see~$ \rfrac{a\odot b}{a}
    = \rfrac{a}{a} \odot b$.

For the second point, suppose~$a,b,c \in M$ with~$a \leq b \leq c$.
    Note~$\rfrac{b}{c} \odot \rfrac{a}{b} \leq \rfrac{b}{c} \leq \rfrac{c}{c}$
        and~$c \odot \rfrac{b}{c} \odot \rfrac{a}{b}
            = b\odot \rfrac{a}{b} = a $.
    So by uniqueness of~$\rfrac{a}{c}$,
        we get~$\rfrac{b}{c} \odot \rfrac{a}{b} = \rfrac{a}{c}$, as desired.
\end{solution}
\begin{solution}{basic-divisoid-equiv}%
Let~$X$ be a compact Hausdorff space.
Assume that the unit interval of~$C(X)$
    is an effect divisoid.
To show~$X$ is basically disconnected,
    assume~$f \in C(X)$.
We have to show that~$\overline{\supp f}$ is open.
Note~$\supp | f | = \supp  f$ and so, without loss of generality,
    we may assume that~$f \geq 0$.
By compactness~$f$ is bounded.  Pick~$B > 0$ with~$B \geq f$.
    Then~$\supp B^{-1} f = \supp f$ and~$B^{-1} f \leq 1$,
    so we may also assume without loss of generality, that~$0 \leq f \leq 1$.
    If~$\overline{\supp f} = X$, then we are done,
    so assume~$\overline{\supp f} \neq X$.
    Pick any~$y \notin \overline{\supp f}$.
By Urysohn's lemma, there is a~$g \in C(X)$
    with~$g(y) = 0$ and~$g(x) = 1$ for every~$x \in \overline{\supp f}$.
    Define~$h \equiv  \rfrac{f}{f} \wedge (0 \vee g)$.
Clearly~$0 \leq h \leq \rfrac{f}{f} \leq 1$
    and~$h(y) \leq (0 \vee g)(y) = 0$.
    Note that~$\rfrac{f}{f}\rfrac{f}{f} \equiv  \rfrac{f}{f} \odot \rfrac{f}{f} = \rfrac{f}{f}$
    and so~$\rfrac{f}{f}$ is zero--one valued: a characteristic function.
Pick any~$x \in \supp f$.
    Then~$0 < f(x) \leq \rfrac{f}{f}(x) \in \{0,1\}$,
        hence~$\rfrac{f}{f}(x)=1$.
    By continuity~$\rfrac{f}{f}(x)=1$
        for every~$x \in \overline{\supp f}$
            and so~$h(x) = \rfrac{f}{f}(x) \wedge g(x) = 1 \wedge 1 = 1$.
Hence~$(f \odot h) (x) = f(x)h(x) = f(x)$ for~$x \in \overline{\supp f}$
    and~$(f \odot h)(x) = f(x)h(x) = 0 = f(x)$ whenever~$x \notin \overline{\supp f}$.
Thus~$f \odot h = f$.
    By uniqueness of~$\rfrac{f}{f}$, we have~$\rfrac{f}{f} = h$.
    Consequently~$\rfrac{f}{f}(y) = h(y) = 0$.
    Recall that~$y$ was an arbitrary element~$y \notin \overline{\supp f}$
        and thus~$\rfrac{f}{f}$ is the characteristic function of~$\overline{\supp f}$,
            which is thus open. Hence~$X$ is basically disconnected.

To prove the converse, assume~$X$ is basically disconnected.
Let~$f,g \in C(X)$ be given with~$0 \leq f \leq g \leq 1$.
    Define~$U_n \equiv \{x; \ g(x) > \frac{1}{n}\}$.
Note~$U_n = \supp ((g - \frac{1}{n}) \vee 0)$
    and so~$\overline{U_n}$ is open as~$X$ is basically disconnected.
Define
\begin{equation*}
    h_n \ \equiv \  \begin{cases}
        \frac{f(x)}{g(x)} & x \in \overline{U_n}\\
        0 & \text{otherwise}.
    \end{cases}
\end{equation*}
To show~$h_n$ is continuous, assume
    there is a net~$(x_\alpha)_\alpha$ with~$x_\alpha \to x$
        for some~$x \in X$.
    Suppose~$x \in \overline{U_n}$.
    As~$\overline{U_n}$ is open,
        we know~$x_\alpha \in U_n$
        for sufficiently large~$\alpha$
        and so~$h_n(x_\alpha) \to h_n(x)$
            as~$\frac{f(x)}{g(x)}$
            is continuous for~$x \in \overline{U_n}$
                as then~$g(x) \geq \frac{1}{n} > 0$.
    For the other case, suppose~$x \notin\overline{ U_n}$.
The set~$X - \overline{U_n}$ is open
    and so~$x_\alpha \notin \overline{U_n}$ for sufficiently large~$\alpha$
        and then~$h_n(x_\alpha)=0 = h_n(x)$.
    Thus~$h_n$ is continuous.

Clearly~$0\leq h_1 \leq h_2 \leq \ldots \leq 1$
 and so we may define~$\rfrac{f}{g} \equiv \sup_n h_n$
 as~$C(X)$ is~$\omega$-complete by basic disconnectedness of~$X$.
 Note~$0 \leq \rfrac{f}{g} \leq 1$.

Suppose~$f=g$ and~$x \in \supp f$.
Then~$h_n(x) = 1$ for all~$n > f(x)^{-1}$
    and so~$\rfrac{f}{f}(x) = 1$.
Thus~$\rfrac{f}{f}(x) = 1$ for all~$x \in \overline{\supp f}$.
Write~$\chi$ for the characteristic function of~$\overline{\supp f}$.
We just saw~$\chi \leq \rfrac{f}{f}(x)$.
As~$\overline{U_n} \subseteq \overline{\supp f}$ and~$h_n \leq 1$,
    we have~$h_n \leq \chi$ for all~$n$ and
    so~$\rfrac{f}{f} \leq \chi$. Thus~$\rfrac{f}{f} = \chi$,
    the characteristic function of~$\overline{\supp f}$.
In particular~$f \leq \rfrac{f}{f}$
    and~$\rfrac{f}{f} \odot \rfrac{f}{f}
        \equiv \rfrac{f}{f}\rfrac{f}{f} = \rfrac{f}{f}$.

Let~$f,g \in C(X)$ be given with~$0 \leq f\leq g \leq 1$
    and~$h_1 \leq h_2\leq  \ldots$ as above.
We will show~$\rfrac{f}{g} \leq \rfrac{g}{g}$.
Assume~$x \in X$.
Suppose~$x \in \supp g$.
Then~$h_n(x) = \frac{f(x)}{g(x)} \leq 1 = \rfrac{g}{g}(x)$
    for all~$n > g(x)^{-1}$
    and so~$\rfrac{f}{g}(x) \leq \rfrac{g}{g}(x)$.
In particular~$\rfrac{f}{g}(x) \leq \rfrac{g}{g}(x)$
    for all~$x \in \overline{\supp g}$.
For the other case, assume~$x \notin \overline{\supp g}$.
Then~$x \notin \overline{U_n}$ and
    so~$h_n(x) = 0$ for all~$n \in \N$,
    whence~$\rfrac{f}{g}(x) = 0 \leq \rfrac{g}{g}(x)$.
Thus indeed~$\rfrac{f}{g} \leq \rfrac{g}{g}$.

In a similar way, it is easy to see
    that~$\rfrac{f}{g}(x) \leq \frac{f(x)}{g(x)}$ for~$x \in \supp g$.
    To show equality, we define
\begin{equation*}
    k_n \ \equiv \  \begin{cases}
        \frac{f(x)}{g(x)} & x \in \overline{U_n}\\
        1 & \text{otherwise}.
    \end{cases}
\end{equation*}
With similar reasoning as for~$h_n$,  we see that these~$k_n$ are continuous.
Furthermore~$h_m \leq k_n$ for all~$m$
    and so~$\rfrac{f}{g} \leq k_n$ for any~$n$.
Pick any~$x \in \supp g$.
Then~$\rfrac{f}{g} (x) \leq k_n(x) = \frac{f(x)}{g(x)}$
    for~$n > g(x)^{-1}$
    and so~$\rfrac{f}{g}(x) = \frac{f(x)}{g(x)}$
        for any~$x \in \supp g$.
Thus~$(g \odot \rfrac{f}{g} )(x)= f(x) $ for~$x \in \supp g$
    and so~$(g\odot \rfrac{f}{g})(x) = f(x)$ for~$x \in \overline{\supp g}$
    by continuity.
For the other case, assume~$x \notin \overline{\supp g}$.
Then~$g \odot \rfrac{f}{g} (x) \leq \rfrac{f}{g} (x) \leq \rfrac{g}{g}(x) = 0
    = g(0) \geq f(0)$.
We have shown~$g\odot \rfrac{f}{g} = f$.

Only uniquness of~$\rfrac{f}{g}$ remains.
So assume~$h \in C(X)$
    with~$0 \leq h \leq \rfrac{g}{g}$
    and~$g \odot  h = f$.
Assume~$x \in \supp g$.
Then~$g(x)h(x) = f(x)$ and so~$h(x) = \frac{f(x)}{g(x)} = \rfrac{f}{g}(x)$.
By continuity~$h(x) = \rfrac{f}{g} (x)$ for all~$x \in \overline{\supp g}$.
For the other case, assume~$x \notin \overline{\supp g}$.
As both~$h, \rfrac{f}{g} \leq \rfrac{g}{g}$,
    we have~$h(x) = 0 = \rfrac{f}{g}(x)$.
We have shown~$h = \rfrac{f}{g}$
    and thus the unit interval of~$C(X)$
    is an effect divisoid.
\end{solution}
\begin{solution}{quotient-basics}%
Assume~$C$ is an effectus.
\begin{enumerate}
\item
Assume~$\xi \colon X \to Y$ is a quotient for~$p$
    and~$\vartheta \colon Y \to Z$ is an isomorphism.
Note~$1 \after \vartheta \after \xi = 1 \after \xi \leq p^\perp$.
To prove~$\vartheta \after \xi$ is a quotient for~$p$,
    assume~$f\colon X \to Y'$ is some map
        with~$1 \after f \leq p^\perp$.
    As~$\xi$ is a quotient, there exists a unique map~$f'\colon Y \to Y'$
        with~$f' \after \xi = f$.
        Clearly~$f = f' \after \xi =
            (f' \after \vartheta^{-1}) \after (\vartheta \after \xi)$.
Assume~$h\colon Z \to Y'$ is any map with~$f = h\after (\vartheta \after \xi)$.
Then~$h \after \vartheta = f'$ by uniqueness of~$f'$
        and so~$h = f' \after \vartheta^{-1}$.
    We have shown that~$\vartheta \after \xi$ is a quotient of~$p$ as well.
\item
    Assume~$\xi_1 \colon X \to Y_1$
            and~$\xi_2 \colon X \to Y_2$
        are quotients for~$p$.
    By the universal property of~$\xi_1$,
            there is a unique map~$\vartheta_1 \colon Y_1 \to Y_2$
            with~$\vartheta_1 \after \xi_1 = \xi_2$.
    Similarly, there is a unique map~$\vartheta_2 \colon Y_2 \to Y_1$
            with~$\vartheta_2 \after \xi_2 = \xi_1$.
Note that~$\xi_1 = \vartheta_2 \after \xi_2 = \vartheta_2 \after \vartheta_1 \after \xi_1$.
By the universal property of~$\xi_1$ again,
        the map~$\id\colon Y_1 \to Y_1$
            is the unique map with~$\id \after \xi_1 = \xi_1$
                and so~$\vartheta_2 \after \vartheta_1 = \id$.
                Similarly~$\vartheta_1 \after \vartheta_2 = \id$.
        Thus~$\vartheta_2$ is an isomorphism
            with~$\vartheta_2 \after \xi_2 = \xi_1$.
            It is the unique such isomorphism by the universal property of~$\xi_2$.
\item
Assume~$f\colon X \to Y$ is some map.
        Clearly~$1 \after f \leq 1 = 0^\perp$.
Furthermore~$f$ is itself the unique map~$f'$ with~$f' \after \id = f$.
Thus~$\id$ is a quotient for~$0$.
        Consequently any isomorphism is a quotient for~$0$ by the first point.
\item
Assume~$X$ is any object.
As~$0$ is a zero object,
    there is a unique map~$0_X \colon X \to 0$.
Clearly~$1 \after 0 = 0 = 1^\perp$.
Assume~$f\colon X \to Y$ is any map with~$1 \after f \leq 1^\perp= 0$.
Then~$f = 0$ and so~$f = 0 = 0 \after 0_X$.
If~$f = h \after 0_X$ for some~$h\colon 0 \to Y$,
    then~$h=0$ and so any map into zero is indeed a quotient for~$1$.
\item
Assume~$\xi \colon X \to Y$ is a quotient for~$p$.
Clearly~$1 \after p^\perp = p^\perp \leq p^\perp$
    and so there is a unique map~$f\colon Y \colon 1$
        with~$f \after \xi = p^\perp$.
Then~$p^\perp = 1\after p^\perp=  1 \after f \after \xi \leq 1 \after \xi = p^\perp$
    and so~$1 \after \xi = p^\perp$.
\item
Assume~$f_1 \after \xi = f_2 \after \xi$
    for some quotient~$\xi$ of~$p$.
        Then~$1 \after f_1 \after \xi \leq 1 \after \xi = p^\perp$.
    Hence there is a unique~$f$ with~$f \after \xi = f_1 \after \xi$.
    Both~$f_1$ and~$f_2$ fit the bill, hence~$f_1=f=f_2$.
    Thus~$\xi$ is epic.
\end{enumerate}
\end{solution}
\spacingfix{}
\begin{solution}{quot-fact-system}%
Let~$C$ be an effectus with quotients.
Assume~$t' \after \xi' = t \after \xi$
    for some quotient~$\xi,\xi'$ and total maps~$t,t'$.
Then~$1 \after \xi = 1 \after t \after \xi = 1 \after t' \after \xi'
        = 1 \after \xi'$.
    Thus by~\sref{quotient-basics} there is a unique isomorphism~$\vartheta$
        with~$\xi' = \vartheta \after \xi$.
Note~$t \after \xi = t' \after \xi' = t' \after \vartheta \after \xi$
    and so by epicity of~$\xi$
    (see~\sref{quotient-basics}),
    we see~$t = t' \after \vartheta$, as desired.
\end{solution}
\begin{solution}{exc-quot-adjoint}%
Assume~$C$ is an effectus with quotients.
We will show that~$0$ has a left adjoint by demonstrating the universal mapping property.
Pick for every object~$(X,p)$ in~$\int \Pred_\square$
    a quotient~$\xi_p \colon X \to X/_p$ of~$p$.
    Note~$(0^\perp \after \xi_p)^\perp = p$
        and so~$\xi_p$ is a map~$(X,p) \to (X/_p, 0)$ in~$\int \Pred_\square$.
We will use these maps~$\xi_p$ as the components of  the unit of the adjunction.
    Let~$f \colon (X,p) \to (Y,0)$ be some map in~$\int \Pred_\square$.
By definition of map~$1 \after f \leq p^\perp$
    and so by the universal property of~$\xi_p$,
    there exists a unique map~$f'\colon X/_p \to Y $
        with~$f = f' \after \xi_p$ in~$C$.
        Clearly also~$(0 f) \after \xi_p = f$ in~$\int\Pred_\square$.
Let~$g\colon X/_p \to Y$ be some map with~$(0 g) \after \xi_p = f$
    as well.
    Then~$g \after \xi_p = f$ and so~$g' = f'$.
This shows that~$0$ has a left adjoint.

To prove the converse, assume~$C$ is an effectus
    where~$0\colon C \to \int\Pred_\square$ has a left adjoint~$Q\colon \int \Pred_\square \to C$.
Let~$X$ be some object in~$C$ with a predicate~$p$.
    Write~$\eta \colon (X,p) \to 0Q(X,p)$
        for the~$(X,p)$ component of the unit of the adjunction~$Q \dashv 0$.
By definition of maps in~$\int \Pred_\square$,
    we know~$1 \after \eta \leq p^\perp$.
We will show that~$\eta$ is a quotientfor~$p$.
To this end, assume~$f\colon X \to Y$ is some map in~$C$
    with~$1 \after f \leq p^\perp$.
Then~$f\colon (X,p) \to (Y,0)$ in~$\int \Pred_\square$.
By the universal mapping property,
    there is a unique map~$f'\colon Q(X,p) \to Y$ in~$C$
    with~$(0f')\after \eta = f$, i.e.~$f' \after \eta = f$.
This shows that~$\eta$ is indeed a quotient of~$p$.
Hence~$C$ has quotients.
\end{solution}
\begin{solution}{compr-grothendieck}%
Assume~$C$ is an effectus with comprehension.
We will show that~$1$ has a right adjoint by demonstrating
    the dual of the universal mapping property.
Pick for every object~$(X,p)$ in~$\int \Pred_\square$
    a comprehension~$\pi_p \colon \cmpr{X}{p} \to X$.
Note that for any map~$f$,
    we have~$1 \after f = p \after f$ if and only if~$p^\perp \after f = 0$.
    Thus~$(p^\perp \after \pi_p)^\perp = 1$, which shows
        that~$\pi_p$ is a map~$(\cmpr{X}{p}, 1) \to (X,p) $ in~$\int\Pred_\square$.
    We will use these maps as components for the counit of the adjunction.
Let~$f\colon (Y,1) \to (X,p)$ be some map in~$C$.
By  definition of maps in~$\int \Pred_\square$,
    we have~$1 \leq (p^\perp \after f)^\perp$,
    viz.~$1 \after f = p \after f$.
    Thus, by the universal property of~$\pi_p$,
    there is a unique map~$f'\colon Y \to \cmpr{X}{p}$ in~$C$
    with~$\pi_p \after f' = f$.
    Hence~$1 f'$ is the unique map in~$\int \Pred_\square$
        with~$\pi_p \after (1 f') = f$ in~$\int \Pred_\square$.
    We have shown that~$1$ has a right adjoint.

To prove the converse, assume~$C$ is an effectus
    where~$1\colon C \to \int \Pred_\square$
    has a right adjoint~$K\colon \int \Pred_\square \to C$.
Let~$X$ be some object of~$C$ with a predicate~$p$.
    Write~$\varepsilon\colon 1 K(X,p) \to (X,p)$ for the~$(X,p)$ component of
    the counit of the adjunction~$1 \dashv K$.
    By definition of maps in~$\int \Pred_\square$,
        we know~$1 \leq (p^\perp \after \varepsilon)^\perp$,
     hence~$1 \after \varepsilon = 1 \after \varepsilon$.
We will show that~$\varepsilon$ is a comprehension for~$p$.
To this end, assume~$f\colon Y \to X$ is some map in~$C$
    with~$1 \after f = p \after f$.
    Then~$1 = (p^\perp \after f)^\perp$
        and so~$f$ is a map~$(Y, 1) \to (X,p)$
        in~$\int \Pred_\square$.
By the dual of the universal mapping property,
    there is a unique map~$f'\colon Y \to K(X,p)$ in~$\int \Pred_\square$
    with~$\varepsilon \after (1f') = f$,
    i.e.~$\varepsilon \after f' = f$.
    This shows that~$\varepsilon$ is indeed a comprehension for~$p$.
    Hence~$C$ has comprehension.
\end{solution}
\begin{solution}{compr-basics}%
Let~$C$ be an effectus.
 \begin{enumerate}
    \item
Assume~$\pi \colon X \to Y$ is a comprehension for~$p$
    and~$\vartheta\colon Z \to X$ is an isomorphism.
We will show~$\pi \after \vartheta$ is a comprehension for~$p$
    as well.
To start, note~$1 \after \pi \after \vartheta  = p \after \pi \after \vartheta$.
Assume~$f\colon X' \to Y$ is any map with~$1 \after f = p\after f$.
Then, by the universal property of~$\pi$, there is a unique map~$f' \colon X' \to X$
    with~$\pi \after f' = f$.
         Then~$f = (\pi \after \vartheta) \after (\vartheta^{-1} \after f')$.
Assume~$g \colon X' \to X$ is a map with~$g = (\pi \after \vartheta) \after g$
    as well.
Then~$\vartheta \after g = f'$ by uniqueness of~$f'$ and so~$
    g = \vartheta^{-1} \after f'$.
We have shown that~$\pi \after \vartheta$ is a comprehension for~$p$ as well.
\item
Assume~$\pi_1 \colon X_1 \to Y$ and~$\pi_2 \colon X_2 \to Y$
    are comprehensions of~$p$.
By the universal property of~$\pi_1$, there is a unique map~$\vartheta_1\colon X_2 \to X_1$
    with~$\pi_2 = \pi_1 \after \vartheta_1$.
On the other side, by the universal property of~$\pi_2$,
     there is a unique map~$\vartheta_2\colon X_1 \to X_2$
    with~$\pi_1 = \pi_2 \after \vartheta_2$.
Then~$\pi_1 = \pi_2 \after \vartheta_2 = \pi_1 \after \vartheta_1 \after \vartheta_2$.
By the universal property of~$\pi_1$ again, the map~$\id \colon X_1 \to X_1$
    is the unique map with~$\pi_1 \after \id = \pi_1$
         and so~$\vartheta_1 \after \vartheta_2 = \id$.
    Similarly~$\vartheta_2 \after \vartheta_1 = \id$.
    Thus~$\vartheta_2$ is an isomorphism with~$\vartheta_1 = \pi_2 \after \vartheta_2$.
        It is the unique such isomorphism by the universal property of~$\pi_2$.
\item
Assume~$f\colon X \to Y$ is some map.  Trivially~$1 \after f = 1 \after f$ and~$1 \after \id = 1 \after \id$.
The map~$f$ itself is the unique map~$g$ with~$g \after \id  = f$
    and so~$\id$ is a comprehension for~$1$.
    Consequently any isomorphism is a comprehension for~$1$ by the first point.
\item
Assume~$Y$ is any object.
As~$0$ is a zero object, there is a unique map~$0_Y \colon 0 \to Y$.
Clearly~$1 \after 0_Y = 0 = 0 \after 0_Y$.
Assume~$f\colon X \to Y$ is any map with~$1 \after f = 0\after f$.
Then~$1 \after f= 0 \after f = 0$ and so~$f=0$.
The zero map is the unique map~$ X \to 0$
    and for that map, we have~$0_Y \after 0 = 0 = f$
         and so~$0_Y$ is a comprehension for~$0$.
\item
Assume~$\pi \after f_1 = \pi \after f_2$ for some comprehension~$\pi$
    of~$p$.
Then~$1 \after \pi \after f_1
         = p \after \pi \after f_1$
    and so by the universal property of~$\pi$,
         there is a unique map~$f$ with~$\pi \after f = \pi \after f_1$.
         Both~$f_1$ and~$f_2$ fit the bill, hence~$f_1 =f = f_2$.
         Thus~$\pi$ is monic.
\item
Assume~$\pi$ is a comprehension of~$p$.
         Then~$(p^\perp \after \pi) \ovee (p \after \pi) = 
                (p^\perp \ovee p) \after \pi = 1 \after \pi =p \after \pi$.
        Thus~$p^\perp \after \pi = 0$ by cancellation.
\end{enumerate}
\end{solution}
\spacingfix{}
\begin{solution}{compr-is-kernel}%
Let~$C$ be any effectus.
Assume~$\pi\colon X \to Y$ is a comprehension of~$p$.
We will show~$\pi$ is a categorical kernel~$p^\perp$.
To start, note~$p^\perp \after \pi = 0$ by~\sref{compr-basics}.
Now assume~$f\colon X' \to Y$ is some map with~$p^\perp \after f = 0$.
    Then~$p \after f = (p \after f) \ovee (p^\perp \after f)= 1 \after f $
    and so by the universal property of~$\pi$,
    there is a unique map~$f'\colon X' \to X$ with~$f = \pi \after f'$.
    Hence~$\pi$ is a kernel of~$p^\perp$.

To prove the converse, assume~$k\colon X \to Y$ is a kernel
        of~$p^\perp\colon Y \to 1$. 
We have~$1 \after k = (p \after p^\perp) \after k = (p \after k) \ovee (p^\perp \after k)
            = p \after k$.
Assume~$f\colon X' \to Y$ is some map with~$p \after f = 1 \after f$.
Then~$(p^\perp \after f) \ovee (p \after f) = 1 \after f = p \after f$
    and so by cancellation~$p^\perp \after f = 0$.
Hence by the universal property of~$k$
    there is a unique map~$f'\colon X' \to X$
        with~$k \after f' = f$.
    This shows that~$k$ is a comprehension of~$p$.
\end{solution}
\begin{solution}{im-ineq}%
We have~$(\IM f ) \after f \after g = 1 \after f \after g \leq 1 \after g$
    and so~$\IM f \after g \leq \IM f$.
    Consequently~$\IM f \after \alpha \leq \IM f
        = \IM f \after \alpha \after \alpha^{-1}
        \leq \IM f \after \alpha$ and so~$\IM f \after \alpha  = \IM f$.
\end{solution}
\begin{solution}{exc-quot-faithful}%
Assume~$\xi\colon X \to Y$ is some quotient and~$p$ is a predicate on~$Y$
    with~$p \after\xi = 0$.
    Then~$p \after \xi = 0 = 0 \after \xi$
        and so~$p = 0$ by epicity of~$\xi$.
Hence~$\xi$ is faithful.
\end{solution}
\begin{solution}{img-of-compr}%
Assume~$\floor{p} = p$.
    Then~$p = \floor{p} \equiv \IM \pi$ for some comprehension~$\pi$ of~$p$
        and so~$p$ is sharp.
To prove the converse, assume~$p$ is sharp;
    that is: $p = \IM f$ for some~$f$.
Pick any comprehension~$\pi$ of~$p$.
    Clearly~$p \after f = (\IM f) \after f = 1 \after f$
        and so there is a unique~$g$ with~$f = \pi \after g$.
    By~\sref{im-ineq}
        we find~$p = \IM f = \IM \pi \after g \leq \IM \pi \equiv \floor{p} \leq p$
        and so~$p = \floor{p}$.
\end{solution}
\begin{solution}{ceiling-within-ceiling}%
Note~$\ceil{p} \leq \ceil{q}$  whenever~$p \leq q$
    and so~$\ceil{p \after f} \leq \ceil{\ceil{p} \after f}
        \leq \ceil{\ceil{p \after f}} = \ceil{p \after f}$
        as~$p \leq \ceil{p}$
        and~$\ceil{p} \after f \leq \ceil{p \after f}$ by~\sref{ceiling-within-ceiling}.
    Thus~$\ceil{p \after f} = \ceil{\ceil{p} \after f}$.
\end{solution}
\begin{solution}{img-tupling}%
To start, note ~$[\IM f, \IM g]
        \langle f,g\rangle = ((\IM f) \after f) \ovee ((\IM g) \after g)
            = (1 \after f) \ovee (1 \after g) = 1 \after \langle f,g\rangle$
            by~\sref{eff-prod-rules}.
    Assume~$p\equiv[p_1,p_2]$ is some predicate with~$p \after \langle f,g \rangle
        = 1 \after \langle f, g \rangle$.
    Then~$(1 \after f) \ovee (1 \after g)
        = 1 \after \langle f, g \rangle 
        = p \after \langle f, g \rangle
        = (p_1 \after f) \ovee (p_2 \after g)
        \leq (p_1 \after f) \ovee (1 \after g)$
        and so by cancellation~$1 \after f \leq p_1 \after f \leq 1 \after f$,
        hence~$1 \after f = p_1 \after f$.
Thus~$\IM f \leq p_1$.
Similarly~$\IM g \leq p_2$.
    Thus~$[\IM f, \IM g] \leq [p_1, p_2] = p$.
We have shown~$\IM \langle f, g\rangle = [\IM f, \IM g]$.

Assume~$[p,q]$ is sharp.
    That is, there is some~$f$ with~$[p,q] = \IM f$.
    Then~$[p,q] = \IM f = \IM \langle \pproj_1 \after f, \pproj_2 \after f \rangle 
    = [\IM \pproj_1 \after f, \IM \pproj_2 \after f]$
    and so~$p = \IM \pproj_1 \after f$
    and~$q = \IM \pproj_2 \after f$,
    which shows~$p$ and~$q$ are sharp.

To prove the converse, assume~$p$ and~$q$ are sharp.
Pick~$f$ and~$g$ with~$p = \IM f$ and~$g = \IM g$.
Note~$f+g = \langle f \after \pproj_1 , g \after \pproj_2 \rangle$.
We have~$\IM f = \IM f \after \pproj_1 \after \kappa_1
                \leq \IM f \after \pproj_1  \leq \IM f$
                and so~$\IM f \after \pproj_1 = \IM f$.
Similarly~$\IM g \after \pproj_2 = \IM g$.
Thus~$\IM f+g =  \IM \langle f \after \pproj_1 , g \after \pproj_2 \rangle
    =  [\IM f \after \pproj_1, \IM g \after \pproj_2]
    = [\IM f, \IM g] = [p,q]$, which shows~$[p,q]$ is sharp.
\end{solution}
\begin{solution}{exc-cokernels}%
Let~$C$ be an effectus with comprehension and images.
   Assume~$f\colon X \to Y$ is a quotient of sharp~$s$.
Let~$\pi_s$ denote a comprehension of~$s$.
Clearly~$1 \after f \after \pi_s = s^\perp \after \pi_s = 0$.
Assume~$g\colon X \to Y'$ is some map with~$g \after \pi_s = 0$.
Then~$1 \after g \after \pi_s = 0$ and so~$1 \after g \leq \IMperp \pi_s = s^\perp$.
Thus by the universal property of~$f$ as a quotient,
    there is a unique map~$g'\colon Y \to Y'$
        with~$g = g' \after f$.
        This shows~$f$ is a cokernel of~$\pi_s$.

To prove the converse, assume~$f$ is a cokernel of~$\pi_s$.
By assumption~$f \after \pi_s = 0$.
Thus~$1 \after f \after \pi_s = 0$.
Hence~$1 \after f \leq \IMperp \pi_s = s^\perp$.
Let~$g\colon X \to Y'$ be a map with~$1 \after g \leq s^\perp$.
Then~$1 \after g \after \pi_s \leq s^\perp \after \pi_s = 0$
    and so~$g \after \pi_s = 0$.
    Thus by the universal property of~$f$ as a cokernel,
        there is a unique map~$g'\colon Y\to Y'$
        with~$g = g' \after f$.
This shows~$f$ is a quotient of~$s$.
\end{solution}
\begin{solution}{exc-diam-order-pres}%
Assume~$s \leq t$.
Then~$s \after f \leq t \after f$,
    hence~$(t \after f) ^\perp \leq (s \after f)^\perp$ so~$\floor{(t \after f)^\perp} \leq \floor{(s \after f)^\perp}$.
    Thus~$f^\diamond(s) \equiv \ceil{s \after f}
        \equiv \floor{(s \after f)^\perp}^\perp \leq \floor{(t \after f)^\perp}^\perp \equiv
            f^\diamond(t)$.
        This shows~$f^\diamond$ is order preserving.
Furthermore~$t^\perp \leq s^\perp$,
    so~$f^\diamond(t^\perp) \leq f^\diamond(s^\perp)$,
    hence~$f^\BOX (s) \equiv f^\diamond(s^\perp)^\perp
        \leq f^\diamond(t^\perp)^\perp \equiv f^\BOX(t)$,
        which shows that~$f^\BOX$ is order preserving.
\end{solution}
\begin{solution}{order-adj-basics}%
Let~$f\colon X \to Y$ be some map in a~$\diamond$-effectus.
\begin{enumerate}
\item
Assume~$s \leq t$.
        Clearly~$f_\diamond (t) \leq f_\diamond(t)$
            so~$s \leq t \leq f^\BOX(f_\diamond(t))$,
            hence~$f_\diamond(s) \leq f_\diamond(t)$,
            which shows that~$f_\diamond$ is order preserving.
\item
Assume~$D \subseteq \SPred X$ is a set of sharp predicates with a supremum~$\sup D$.
        We will show~$f_\diamond(\sup D)$ is the supremum of~$f_\diamond(D)$.
For any~$d \in D$, we have~$d \leq \sup D$
    and so~$f_\diamond (d) \leq f_\diamond(\sup D)$ by the previous point.
Assume~$x$ is some sharp predicate of~$Y$
    with~$f_\diamond(d) \leq x$ for all~$d \in D$ as well.
Then~$d \leq f^\BOX(x)$ for all~$d \in D$
    and so~$\sup D \leq f^\BOX(x)$.
Hence~$f_\diamond(\sup D) \leq x$, which shows
    that~$\sup f_\diamond (D) = f_\diamond(\sup D)$.
\item
Assume~$D \subseteq \SPred Y$ is a set of sharp predicates with an infimum~$\inf D$.
        We will show~$f^\BOX(\inf D)$ is the infimum of~$f^\BOX(D)$.
For any~$d \in D$, we have~$\inf D \leq d$
        and so~$f^\BOX(\inf D) \leq f^\BOX(d)$
        by~\sref{exc-diam-order-pres}.
Assume~$x$ is some sharp predicate of~$X$
    with~$x \leq f^\BOX(d)$ for all~$d \in D$ as well.
    Then~$f_\diamond(x) \leq d$ for all~$d \in D$
            and so~$f_\diamond(x) \leq \inf D$,
hence~$x \leq f^\BOX(\inf D)$.
        This shows that~$\inf f^\BOX(D) = f^\BOX(\inf D)$.
\item
    For for any~$D \subseteq \SPred Y$ with supremum,
        we have
        \begin{align*}
        f^\diamond(\sup D)
            & \ \equiv\ f^\BOX((\sup D)^\perp)^\perp \\
            &\  =\ f^\BOX(\inf_{d \in D} ( d^\perp))^\perp \\
            &\ =\ \bigl(\inf_{d \in D}f^\BOX(  d^\perp)\bigr)^\perp \\
            &\ =\ \sup_{d \in D}f^\BOX(  d^\perp)^\perp \\
            &\ =\ \sup_{d \in D}f^\diamond(d)
        \end{align*}
    and so~$f^\diamond$ preserves suprema.
\item
For any~$s \in \SPred X$
        we have~$f_\diamond (s) \leq f_\diamond(s)$
        and so~$s \leq f^\BOX (f_\diamond(s))$.
Similarly, for any~$t \in \SPred Y$
    we have~$f_\diamond (f^\BOX (t)) \leq t$.
Hence~$ f_\diamond (f^\BOX (f_\diamond(s)))
                \leq f_\diamond(s)
                \leq f_\diamond(f^\BOX(f_\diamond(s)))$,
                so~$f_\diamond = f_\diamond \after f^\BOX \after f_\diamond$.
\item
We have~$ f^\BOX(f_\diamond(f^\BOX(t)))
            \leq f^\BOX(t) \leq f^\BOX(f_\diamond(f^\BOX(t)))$
            and so~$f^\BOX = f^\BOX \after f_\diamond \after f^\BOX$.
\end{enumerate}
\end{solution}
\spacingfix{}
\begin{solution}{diamond-equiv-equiv}%
Let~$f,g\colon X \to Y$  be some maps in a~$\diamond$-effectus.
Assume~$f^\diamond = g^\diamond$ and~$t \in \SPred X$.
    Trivially~$g_\diamond(t) \leq (g_\diamond(t)^\perp)^\perp$
        and so~$f^\diamond(g_\diamond(t)^\perp) =
        g^\diamond(g_\diamond(t)^\perp) \leq t^\perp$.
Thus~$f_\diamond(t) \leq g_\diamond(t)^{\perp\perp} = g_\diamond(t)$.
Reasoning in the same way with~$f$ and~$g$ swapped,
        we see~$g_\diamond(t) \leq f_\diamond(t)$ and so~$f_\diamond = g_\diamond$.

To prove the converse, assume~$f_\diamond = g_\diamond$ and~$s \in \SPred Y$.
    Trivially~$g^\diamond(s) \leq (g^\diamond(s)^\perp)^\perp$
        and so~$f_\diamond(g^\diamond(s)^\perp) =
        g_\diamond(g^\diamond(s)^\perp) \leq t^\perp$.
Thus~$f^\diamond(t) \leq g^\diamond(t)^{\perp\perp} = g^\diamond(t)$.
Reasoning in the same way with~$f$ and~$g$ swapped,
        we see~$g^\diamond(t) \leq f^\diamond(t)$ and so~$f^\diamond = g^\diamond$.
\end{solution}
\begin{solution}{spred-sup}%
As~$\xi_\diamond$ is left adjoint to~$\xi^\BOX$,
    we have~$t \leq \xi^\BOX (\xi_\diamond(t))$.
    Furthermore~$s = \ceil{s^\perp}^\perp = \ceil{1 \after \xi}^\perp
        = \xi^\BOX(0) \leq
        \xi^\BOX (\xi_\diamond(t))$.
Thus~$\xi^\BOX (\xi_\diamond(t))$ is an upper bound of~$s$ and~$t$.
We have to show it is the least sharp upper bound.
So assume~$r$ is any sharp predicate with~$r \geq s$ and~$r \geq t$.
Then~$1 \after \xi_r = r^\perp \leq s^\perp$
    and so there is a (unique) map~$h$ with~$\xi_r = h \after \xi$. We have
\begin{equation*}
    (\xi_r)^\BOX
     \ = \ \xi^\BOX \after h^\BOX
     \ = \ \xi^\BOX  \after \xi_\diamond \after \xi^\BOX\after h^\BOX
     \ = \ \xi^\BOX  \after \xi_\diamond \after (\xi_r)^\BOX
\end{equation*}
    and so~$r
    = (\xi_r)^\BOX(0)
    = \xi^\BOX (\xi_\diamond((\xi_r)^\BOX(0)))
    = \xi^\BOX (\xi_\diamond(r))
    \geq \xi^\BOX (\xi_\diamond(t)) $.
\end{solution}
\begin{solution}{exc-diamond-adj}%
For the first point, assume~$f^\diamond = g_\diamond$
    for some~$f : X \leftrightarrows Y : g$.
    Trivially~$g^\diamond(s) \leq (g^\diamond(s)^\perp)^\perp$ for any~
    for any~$s \in \SPred X$ and
    so~$ f^\diamond(g^\diamond(s)^\perp) =
    g_\diamond(g^\diamond(s)^\perp) \leq s^\perp$.
    Hence~$f_\diamond(s) \leq g^\diamond(s)^{\perp\perp} = g^\diamond(s)$.
With essentially the same argument
    (swapping  both~$f$ and~$g$ and the up/down position of the diamonds),
    we see~$g^\diamond(s) \leq  f_\diamond(s)$,
    hence~$g^\diamond = f_\diamond$.
Swapping~$f$ and~$g$ again, gives the argument for the converse.

For the second point, assume~$f$ and~$g$ are~$\diamond$-adjoint.
    Then~$\IM f = f_\diamond(1) = g^\diamond(1) = \ceil{1 \after g}$,
    as desired.
\end{solution}
\begin{solution}{diamond-squares}%
For the first point, assume~$f$ is~$\diamond$-self-adjoint.
Then~$(f \after f)_\diamond
        = f_\diamond \after f_\diamond
        = f^\diamond \after f^\diamond
        = (f \after f)^\diamond $
        and so~$f\after f$ is~$\diamond$-self-adjoint as well.

For the second point, assume~$f$ is~$\diamond$-positive.
Then~$f = g \after g$ for some~$\diamond$-self-adjoint~$g$.
Hence~$f$ is~$\diamond$-self-adjoint by the previous point.

For the third point, assume~$f$ is~$\diamond$-positive
    and~$f \after f$ is pure.
By the previous point~$f$ is~$\diamond$-self-adjoint
        and so~$f\after f$ is~$\diamond$-self-adjoint by the first point.
    Thus~$f \after f$ is~$\diamond$-positive.
\end{solution}
\begin{solution}{sharp-ceil}%
Assume~$f\colon X \to Y$ is a sharp map and~$p$ is a predicate on~$Y$.
Using~\sref{ceiling-within-ceiling}
    and~\sref{img-of-compr}
    we see~$\ceil{p \after f} = \ceil{\ceil{p} \after f} = \ceil{p} \after f$.
For the converse assume~$\ceil{p} \after f = \ceil{p \after f}$
    for every predicate~$p$ on~$Y$.
Assume~$s$ is any sharp predicate on~$Y$.
    Then~$\ceil{s \after f} = \ceil{s} \after f = s \after f$ and so~$s \after f$
        is sharp.  Thus~$f$ is sharp.
\end{solution}
\begin{solution}{andthen-square-rule}%
    Note~$1 \after \asrt_p \after \asrt_p = p \after \asrt_p \equiv \andthen{p}{p}$.
By \sref{diamond-squares} the map~$\asrt_p \after \asrt_p$
    is~$\diamond$-positive
    and thus~$\asrt_p \after \asrt_p = \asrt_{\andthen{p}{p}}$
    by uniqueness of~$\diamond$-positive maps.
\end{solution}
\begin{solution}{asrt-absorp-rule}%
Assume~$f\colon X\to Y$ is a map in a~$\&$-effectus
    with~$s \in \SPred Y$ and~$t \in \SPred X$.

Assume~$\IM f \leq s$.
Then~$1\after f = s \after f$.
Thus there is a~$g$ with~$f = \pi_s \after g$.
Then~$\asrt_s \after f = \pi_s \after \zeta_s \after\pi_s \after g
        = \pi_s \after g = f$.
For the converse, assume~$\asrt_s \after f = f$.
Then~$1 \after f = 1 \after \asrt_s \after f =s \after f$
    and so~$\IM f \leq s$.
We have shown the first equivalence.

For the second equivalence, assume~$1 \after f \leq t$.
    Then~$f = h \after \zeta_t$ for a (unique) map~$h$.
So~$f \after \asrt_t = h \after \zeta_t \after \pi_t \after \zeta_t
        = h \after \zeta_t = f$.
For the converse, assume~$f = f \after \asrt_t$.
Then~$1 \after f = 1 \after f \after \asrt_t \leq 1 \after \asrt_t = t$.
\end{solution}
\begin{solution}{simple-andthen-absorption}%
Clearly~$p \leq s$ if and only if~$1 \after p \leq s$.
Even more trivially,  $\andthen{s}{p} = p$
        if and only if~$p \after \asrt_s \equiv \andthen{s}{p} = p$.
Thus the result follows by applying the second equivalence
    of~\sref{asrt-absorp-rule} with~$f=p$ and~$t=s$.
\end{solution}
\begin{solution}{exc-prod-sharp-maps}%
Let~$f \colon X \to A$ and~$g\colon X \to B$
    be maps in a~$\&$-effectus.
   Assume~$\langle f, g \rangle$ is sharp.
Pick any sharp predicate~$s$ on~$A$.
By~\sref{img-tupling} we know~$[s,0]$ is sharp
    and so~$[s,0] \after \langle f, g \rangle = (s \after f) \ovee (0 \after g) = s \after f$
    is sharp as well.
This shows~$f$ is sharp.
With a similar argument we see~$g$ is sharp.
For the converse, assume~$f$ and~$g$ are sharp.
    Assume~$[s,t]$ is some sharp predicate on~$A+B$.
    By~\sref{img-tupling} we know~$s$ and~$t$ are sharp.
    Thus~$s \after f$ and~$t \after g$ are sharp.
    Consequently~$[s,t] \after \langle f, g \rangle = (s \after f) \ovee (t \after g)$
    is sharp as well by~\sref{diamond-oml}.
    Thus~$\langle f,g\rangle$ is sharp, as desired.
\end{solution}
\begin{solution}{zeta-through-asrt}%
Assume~$C$ is a~$\&$-effectus
    where every predicate has a square root
    and where~$\pi_s$ is~$\diamond$-adjoint to~$\zeta_s$.
Let~$p$ be a predicate and~$s$ be another predicate on the same object that is sharp.
    Pick a predicate~$q$ with~$\andthen{q}{q}=p$.
    The map~$ f\equiv \pi_s \after \asrt_q \after \zeta_s$ is pure and~$\diamond$-self-adjoint:
\begin{equation*}
    f^\diamond
    \ =\ (\zeta_s)^\diamond \after (\asrt_q)^\diamond \after (\pi_s)^\diamond
    \ =\ (\pi_s)_\diamond \after (\asrt_q)_\diamond \after (\zeta_s)_\diamond
    \ =\ f_\diamond.
\end{equation*}
Thus~$f \after f$ is~$\diamond$-positive.
By~\sref{andthen-square-rule}, we get~$f \after f = \pi_s \after \asrt_q \after \zeta_s \after \pi_s \after \asrt_q \after \zeta_s
    = \pi_s \after \asrt_{\andthen{q}{q}} \after \zeta_s
    = \pi_s \after \asrt_{p} \after \zeta_s$
    and so~$1 \after f \after f = 1 \after \pi_s \after \asrt_p \after
        \zeta_s = p \after \zeta_s$.
    By uniqueness of~$\diamond$-positive maps,
    we get~$\pi_s \after \asrt_{p} \after \zeta_s = f \after f 
        = \asrt_{1 \after f \after f} = \asrt_{p \after \zeta_s}$.
Thus~$\asrt_p \after \zeta_s =
    \zeta_s \after \pi_s \after \asrt_{p} \after \zeta_s
        = \zeta_s \after \asrt_{p \after \zeta_s}$.
\end{solution}
\begin{solution}{dagger-prime-basics}%
By~\sref{asrt-absorp-rule} we have
 \begin{equation*}
    \asrt_p \ = \ 
     \asrt_p \after \asrt_{\ceil{p}}
     \ = \ \pi_{\ceil{p}} \after \id \after \zeta_{\ceil{p}} \after \asrt_p
 \end{equation*}
    and so
\begin{equation*}
    \asrt^\dagger_p \ =\  \asrt_p \after \pi_{\ceil{p}} \after \id^{-1} \after \zeta_{\ceil{p}}
        \ =\  \asrt_p \after \asrt_{\ceil{p}} \ =\  \asrt_p.
\end{equation*}
    Next, note~$\zeta_1 = \pi_1 = \asrt_1 = \id$
    and so~$\pi_s = \pi_s \after \id \after \zeta_1 \after \asrt_1$,
    hence
\begin{equation*}
    \pi_s^\dagger 
       \ =\ \asrt_1 \after \pi_1 \after \id^{-1} \after \zeta_s
       \ =\ \zeta_s.
\end{equation*}
Also~$\zeta_s = \pi_1 \after \id \after \zeta_s \after \asrt_1$
    and so
\begin{equation*}
    \zeta_s^\dagger \ = \ \asrt_1 \after \pi_s \after \id^{-1} \after \zeta_1 \ = \ \pi_s.
\end{equation*}
    Finally, note~$\alpha = \asrt_1 = \pi_1 \after \alpha \after \zeta_1$
        for any isomorphism~$\alpha$ and so
\begin{equation*}
    \alpha^\dagger \ =  \ \asrt_1 \after \pi_1 \after \alpha^{-1} \after \zeta_{1} \ = \ \alpha^{-1}.
\end{equation*}
\end{solution}
\spacingfix{}
\begin{solution}{standard-form-pristine}%
Assume~$h$ is a pristine map in a~$\&$-effectus.
By~\sref{standard-form-map} we know there is an isomorphism~$\alpha$
    with~$h = \pi_{\IM h} \after \alpha \after \zeta_{\ceil{1 \after h}} \after \asrt_{1 \after h}$.
As~$1 \after h$ is sharp,
    we have~$\zeta_{\ceil{1 \after h}} \after \asrt_{1 \after h}
        = \zeta_{1 \after h} \after \pi_{1 \after h} \after \zeta_{1 \after h}
        = \zeta_{1 \after h}$
        and so~$
    h = \pi_{\IM h} \after \alpha \after \zeta_{1 \after h}$.
\end{solution}
\begin{solution}{asrt-pristine-reverse}%
Assume~$h \equiv \pi_{\IM h} \after \alpha \after \zeta_{1 \after h}$
    is a pristine map in a~$\dagger'$-effectus, cf.~\sref{standard-form-pristine}.
\begin{enumerate}
\item
Clearly~$h =\pi_{\IM h} \after \alpha \after \zeta_{\ceil{1 \after h}} \after \asrt_{1 \after h}$
as~$1 \after h$ is sharp and by~\sref{asrt-absorp-rule}.
Thus
\begin{equation*}
    h^\dagger \ = \ \asrt_{1 \after h} \after \pi_{\ceil{1 \after h}}
        \after \alpha^{-1} \after \zeta_{\IM h}
        \ = \ \pi_{1 \after h} \after \alpha^{-1} \after \zeta_{\IM h}
\end{equation*}
        using other equivalence of~\sref{asrt-absorp-rule}.
\item
Note that by the previous point~$1 \after h^\dagger = \IM h$
    and~$\IM h^\dagger = 1\after h$, hence
        by a second appeal to the previous point:
\begin{equation*}
    h^{\dagger\dagger}
    \ = \  \pi_{1 \after h^\dagger} \after (\alpha^{-1})^{-1} \after \zeta_{\IM h^\dagger}
    \ = \  \pi_{\IM h} \after \alpha \after \zeta_{1 \after h} \ = \ h.
\end{equation*}
\item
By the first point~$ h^\dagger \after h
     =  \pi_{1 \after h} \after \alpha^{-1} \after \zeta_{\IM h}
        \after \pi_{\IM h} \after \alpha \after \zeta_{1 \after h}
     =  \pi_{1 \after h} \after \zeta_{1 \after h}
     =  \asrt_{1 \after h}$ as promised.
\item
    Clearly~$p \after h^\dagger \leq 1 \after h^\dagger = \IM h$.
\item
    Assume~$p \leq 1 \after h$.
Then~$p \after \asrt_{1 \after h} = p$ by \sref{asrt-absorp-rule} and
        so~
    \begin{equation*}
        \asrt_{p \after h^\dagger} \after h
                \ \stackrel{\mathclap{\sref{pristine-asrt}}}{=}\
                h \after \asrt_{p \after h^\dagger \after h}
                \ \stackrel{\mathclap{\text{pt.~3}}}{=}\  
                h \after \asrt_{p \after \asrt_{1\after h}}
                \ = \ h \after \asrt_p,
    \end{equation*}
        as desired.
\end{enumerate}
\end{solution}
\spacingfix{}
\begin{solution}{dagger-iso-alpha2}%
    The map~$
    \alpha \after \zeta_{\ceil{t \after \asrt_p \after \pi_{\ceil{p}}}}
        \after \varphi^{-1} \after \zeta_s
        $
        is a quotient for~$\ceil{t \after f^\dagger}$:
\begin{align*}
    1 \after \alpha \after \zeta_{\ceil{t \after \asrt_p \after \pi_{\ceil{p}}}}
        \after \varphi^{-1} \after \zeta_s
        & \ = \ 
     \ceil{t \after \asrt_p \after \pi_{\ceil{p}}}
        \after \varphi^{-1} \after \zeta_s \\
        & \ = \ 
     \ceil{t \after \asrt_p \after \pi_{\ceil{p}}
        \after \varphi^{-1} \after \zeta_s} \\
        & \ = \ 
     \ceil{t \after f^\dagger}.
\end{align*}
Furthermore
\begin{align*}
    & \alpha \after \zeta_{\ceil{t \after \asrt_p \after \pi_{\ceil{p}}}}
    \after \varphi^{-1} \after \zeta_s \after \pi_{\ceil{t \after f^\dagger}} \\
    & \qquad \ \stackrel{\mathclap{\eqref{dagger-iso-alpha}}}{=} \ 
    \alpha \after \zeta_{\ceil{t \after \asrt_p \after \pi_{\ceil{p}}}}
    \after \varphi^{-1} \after \zeta_s \after \pi_s \after \varphi
        \after \pi_{\ceil{t \after \asrt_p \after \pi_{\ceil{p}}}} \after \alpha^{-1}\\
      &\qquad \ = \ \id.
\end{align*}
    Thus~$\alpha \after \zeta_{\ceil{t \after \asrt_p \after \pi_{\ceil{p}}}}
        \after \varphi^{-1} \after \zeta_s$
        is the unique quotient corresponding to the
        comprehension~$\pi_{\ceil{t \after f^\dagger}}$.
    Hence~$\alpha \after \zeta_{\ceil{t \after \asrt_p \after \pi_{\ceil{p}}}}
        \after \varphi^{-1} \after \zeta_s
        = \zeta_{\ceil{t \after f^\dagger}}$, as desired.
\end{solution}
\begin{solution}{dagger-iso-zeta2}%
    Note that~$\pi_{\ceil{q}} \after \psi^{-1}
        = \pi_{\ceil{q}} \after \psi^{-1} \after \zeta_{1} \after \asrt_1$
        and so by definition of the dagger, we have~$(\pi_{\ceil{q}} \after \psi^{-1})^\dagger
            = \asrt_1 \after \pi_1 \after (\psi^{-1})^{-1} \after \zeta_{\ceil{q}}
            = \psi\after \zeta_{\ceil{q}}$.
    Furthermore~$p \after \pi_t \leq 1 = 1 \after \pi_{\ceil{q}} \after \psi^{-1}$.
Thus we can apply the fourth point of~\sref{asrt-pristine-reverse}:
\begin{align*}
    \pi_{\ceil{q}} \after \psi^{-1} \after \asrt_{p \after \pi_t}
    & \ = \ \asrt_{p \after \pi_t \after (\pi_{\ceil{q}} \after \psi^{-1})^{-1}}
    \after \pi_{\ceil{q}} \after \psi^{-1} \\
    & \ = \ \asrt_{p \after \pi_t \after \psi \after \zeta{\ceil{q}}} \after \pi_{\ceil{q}} \after \psi^{-1} \\
    & \ = \ \asrt_{p \after k} \after \pi_{\ceil{q}} \after \psi^{-1},
\end{align*}
which is the desired equality.
\end{solution}
\begin{solution}{exc-cvn-no-dilations}%
Assume~$\scrA$ is a commutative von Neumann algebra
    with effect~$a \in \scrA$.
The standard corner~$h_a\colon \scrA \to \ceil{a} \scrA \ceil{a}$
    for~$a$ (see~\sref{standard-corner-dils})
    is given by~$b \mapsto \ceil{a} b \ceil{a}$.
If~$\scrA$ is commutative,
    then~$h_a(bc) = \ceil{a} bc \ceil{a}
            = \ceil{a} b \ceil{a} \ceil{a} c \ceil{a} = h_a(b) h_a(c)$
            and so~$h_a$ is an nmiu-map.
Every corner is a standard corner up-to an isomorphism,
    thus every corner on a commutative von Neumann algebra
    is an nmiu-map as well.
If~$\op\CvN$ were to have dilations,
    then any ncpu-map would be the composition of a corner and an nmiu-map.
We just saw that such corners are nmiu and so any ncp-map would be an nmiu-map
    as well, which is absurd.  Thus~~$\op\CvN$ does not have dilations
        in the sense of~\sref{dfn-eff-dilations}.
\end{solution}
\begin{solution}{exc-purec-no-biproduct}%
   We will first show that for any non-zero pure map~$f\colon \scrA \to \C$,
        there is a Hilbert space~$\scrH$, an element~$x \in \scrH$,
    von Neumann algebra~$\scrC$ and isomorphism~$\varphi\colon \scrB(\scrH)\oplus \scrC \to \scrA$
    with~$f (\varphi (T, c))= \langle x, T x\rangle$.
By definition of purity, we know~$f = c \after h$
    for some corner~$h$ and filter~$c$ for~$f(1)$.
Clearly~$c(1)=f(1) \neq 0$ as~$f\neq 0$ and
    so~$c\colon \C \to \C$ is simply given by~$c(1)\id$,
    hence~$f = f(1) h$ for a corner~$h\colon \scrA \to \C$
        of~$p \in \scrA$.
    As~$p \scrA p \cong \C$ both have only two projections
        and so~$p$ is a mimimal projection.
    In particular~$0$ is the only central element strictly below~$p$,
        thus~$\cceil{p}\scrA$ is a factor.
Let~$\varrho\colon \cceil{p}\scrA \to \scrB(\scrH)$, $y\in \scrH$
    be a GNS-representation of the restriction of~$h$ to~$\cceil{p}\scrA$
    with~$y\in\scrH$ such that~$h(\cceil{p}a) = \langle y, \varrho(\cceil{p}a) y \rangle$.
    This restriction of~$h$ is still a corner, hence pure
    and so by~\sref{paschke-pure} we know~$\varrho$ is surjective.
    As~$\cceil{p}\scrA$ is a factor, $\varrho$ must be injective as well.
    Thus~$\varrho$ is an nmiu-isomorphism.
    Hence~$\varphi\colon \scrB(\scrH) \oplus \scrB \to \scrA$
        defined by~$\varphi(T,c) \equiv \varrho^{-1}(T) + c$
                is an nmiu-isomorphism as well.
Define~$x \equiv \sqrt{f(1)}$. Then~$f(\varphi(T, c))
            = f(1) h(\varrho^{-1}(T) + c)
            = f(1) h(\varrho^{-1}(T))
            = f(1) \langle y, \varrho (\varrho^{-1}(T)) y \rangle
            = f(1) \langle y, T y \rangle
            =  \langle x, T x \rangle$ as desired.

We continue with the main task of this exercise.
Assume~$(\Pure \op\vN)^{\mathsf{op}}$ has a
        product~$ \C \xleftarrow{\pi_1} \scrA \xrightarrow{\pi_2} \C$
        of~$\C$ and~$\C$.
Note~$\pi_1 \neq 0$
    for otherwise~$\id_\C = \pi_1 \after h = 0$,
    where~$h\colon \C \to \scrA$ is the unique map given by the universal property
    of the product for~$\id_\C$ and~$\id_\C$.
    Similarly~$\pi_2 \neq 0$.

By assumption, there is a unqiue pure map~$h\colon \C \to \scrA$
    with~$\pi_1 \after h = \id$ and~$\pi_2 \after h = \id$.
With the reasoning at the start,
    we know there is a Hilbert space~$\scrH$, an element~ $z\in \scrH$,
    a von Neumann algebras~$\scrC$
    and an isomorphism~$\varphi\colon \scrB(\scrH)\oplus \scrC$
    with~$h^\dagger (\varphi(T,c)) = \langle z, Tz \rangle$.
    Note that~$T \mapsto \langle z, T z\rangle$
        is the dagger of the standard filter~$c_{\ketbra{z}{z}}$.
    Write~$p\colon \scrB(\scrH) \oplus \C \to \scrB(\scrH)$
        for the regular coprojection~$(T,c) \mapsto T$,
        which is corner with~$p^\dagger(T) = (T,0)$.
        Then~$h^\dagger \after \varphi = c^\dagger_{\ketbra{z}{z}} \after p$
            and so~$h = \varphi \after p^\dagger \after c_{\ketbra{z}{z}}$.
    By definition of~$h$, we have 
\begin{equation}\label{exc-224VI1}
    \id \ =\ \pi_1 \after h = \pi_1 \after \varphi \after p^\dagger \after c_{\ketbra{z}{z}}.
\end{equation}
Reasoning the same way as we did for~$f$ at the start
    of this exercise,
    we know that~$\ceil{\pi_1 \after \varphi}$ is a minimal projection.
    If~$\ceil{\pi_1 \after \varphi} \leq (0,1)$,
        then~$\pi_1 \after \varphi \after p^\dagger = 0$
        and so~$\id = 0$ by~\eqref{exc-224VI1}, which is absurd.
    Thus~$\ceil{\pi_1 \after \varphi} = (p_1,0)$
    for some minimal projection~$p_1$ of~$\scrB(\scrH)$.
Thus~$\pi_1 (\varphi(T,c)) = \langle x, T x\rangle$ for some~$x \in \scrH$.
Similarly~$\pi_2 (\varphi(T,c)) = \langle y, T y\rangle$ for some~$y \in \scrH$.

Note that~$
\pi_i \after \varphi \after (\id \oplus \vartheta)
=\pi_i \after \varphi$
for any map~$\vartheta\colon \scrC \to \scrC$ and~$i=1,2$.
Hence~$\id \after \pi_i \after \varphi \after (\id \oplus\vartheta) \after \varphi^{-1} \after h$
    for~$i=1,2$ and so~$(\id\oplus\vartheta) \after \varphi^{-1}\after h = \varphi^{-1}\after h$
    for every~$\vartheta$ by uniqueness of~$h$,
    which is absurd unless~$\scrC = 0$.

Suppose~$\dim \scrH \geq 3$.
Then there is some non-zero~$z_0 \in \scrH$
    with~$z_0 \perp x$ and~$x_0 \perp y$.
So~$\pi_i \after \varphi \after c_{\ketbra{z_0}{z_0}} = 0$ for~$i=1,2$.
By the presumed universal property of~$\scrA$,
    there is a unique map~$h\colon \C \to \scrA$
        with~$\pi_1 \after h = \pi_1 \after \varphi \after c_{\ketbra{z_0}{z_0}}$
        and~$\pi_2 \after h = \pi_2 \after \varphi \after c_{\ketbra{z_0}{z_0}}$
Both~$\varphi \after c_{\ketbra{z_0}{z_0}}$ and~$0$ fit the bill,
    so~$\varphi \after c_{\ketbra{z_0}{z_0}} = 0$,
    which is absured. Thus~$\dim \scrH \leq 2$.

Clearly~$\dim \scrH \neq 0$, for otherwise there would be a unqiue pure map~$\C \to \C$, which is absurd.
If~$\dim \scrH = 1$, then~$x = \lambda y$ for some~$\lambda \in \C, \lambda \neq 0$
    and so~$\pi_1 = \lambda \pi_2$,
    hence for any two pure maps~$f_1,f_2\colon \scrB \to \C$ we would have~$f_1 = \lambda f_2$,
        which is absurd.

Reducing clutter a bit,
    we may assume without loss of generality~$\scrA = M_2$, $\pi_1(T) = \langle x, Tx\rangle$
    and~$\pi_2(T) = \langle y, T y\rangle$ for some~$x,y \in \C^2$.
Write~$f_1,f_2\colon \C^2\to\C$
    for the maps~$f_1(\lambda,\mu) = \lambda$
    and~$f_2(\lambda,\mu) = \mu$.
By the presumed universal property of~$M_2$ as product,
    there is a unique pure map~$f\colon \C^2 \to \scrA$
    with~$\pi_1 \after f = f_1$ and~$\pi_2 \after f = f_2$.
From the definition of purity,
    it follows that~$f = c \after h$
        for some corner~$h\colon \C^2 \to \ceil{h(1)} M_2 \ceil{h(1)}$
        and the standard filter~$c \colon \ceil{h(1)} M_2 \ceil{h(1)} \to M_2$.
Note~$\ceil{h} \C^2 \ceil{h} \cong \ceil{h(1)} M_2 \ceil{h(1)}$.
The former is commutative, so~$\ceil{h(1)} \neq 1$.
If~$\ceil{h(1)} = 0$, then~$h = 0$ and so~$f_1=0$, quod non.
Thus~$\ceil{h(1)} M_2 \ceil{h(1)} \cong \C$ and so
    either~$\ceil{h} = (1,0)$
    or~$\ceil{h} = (0,1)$.
In the first case~$1 = f_2 (0,1)= \pi_2 (h(0,1)) = \pi_2(0) = 0$,
    which is a contradiction.
The second case leads to a similar contradiction.
So there is no product of~$\C$ and~$\C$ in~$(\Pure \op\vN)^{\mathsf{op}}$.
\end{solution}
\begin{solution}{exc-purec-equal}%
   Assume~$\xi\colon \scrC \to M_4$ is some filter
    with~$\ad_\sigma \after \xi = \xi$.
    Note~$\sigma^* = \sigma$.
    Assume~$c \in \scrC$.
    Then~$\xi(c) = \ad_\sigma(\xi (c)) = \sigma \xi(c) \sigma$.
Let~$b \in M_4$ be any element with~$\sigma b \sigma = b$.
    Then~$
    (p_\scrS + p_\scrA) b (p_\scrS +p_\scrA)
        = b = \sigma b \sigma
  =  (\sigma p_\scrS + \sigma p_\scrA) b (p_\scrS  \sigma +p_\scrA \sigma)
  =  ( p_\scrS -  p_\scrA) b (p_\scrS  -p_\scrA )$
    and so~$p_\scrA b p_\scrS + p_\scrS b p_\scrA = 0$.
As~$p_\scrA$ and~$p_\scrS$ have orthogonal ranges,
    we see~$p_\scrA b p_\scrS = 0 = p_\scrS b p_\scrA$.
Hence~$b = p_\scrS b p_\scrS + p_\scrA b p_\scrA$.
Furthermore~$b p_\scrS =
    p_\scrS b p_\scrS$.
If~$b$ is additionally self adjoint, then~$p_\scrS b = 
    p_\scrS b^* p_\scrS =
    p_\scrS b p_\scrS = p_\scrS b$.
Thus~$p_\scrS \xi(c) = \xi(c) p_\scrS$ for any~$c \in \scrC$
There is some isomorphism~$\vartheta \colon \scrC \to \ceil{\xi(1)} M_4 \ceil{\xi(1)}$
    with~$\xi(c) = \sqrt{\xi(1)} \vartheta(c) \sqrt{\xi(1)}$
    and so~$ p_\scrS \sqrt{\xi(1)} b \sqrt{\xi(1)} = \sqrt{\xi(1)} b \sqrt{\xi(1)} p_\scrS$
    for any~$b \in M_4$.
    From this it follows (using the pseudoinverse of~$\sqrt{\xi(1)}$)
    that~$
    p_\scrS \ceil{\xi(1)} b  \ceil{\xi(1)} =
     \ceil{\xi(1)} b  \ceil{\xi(1)}  p_\scrS$.
    Thus~$\ceil{\xi(1)}p_\scrS\ceil{\xi(1)}$ is central in~$\ceil{\xi(1)}M_4 \ceil{\xi(1)}$.
So either~$\ceil{\xi(1)} p_\scrS \ceil{\xi(1)} = 0$
    or~$\ceil{\xi(1)} p_\scrS \ceil{\xi(1)} = \ceil{\xi(1)}$.
    In the first case, we have~$\ceil{\xi(1)} p_\scrS = 0$
        hence~$\xi(c) = p_\scrA \xi(c) p_\scrA + p_\scrS \xi(c) \ceil{\xi(1)} p_\scrS
        = p_\scrA \xi(c) p_\scrA \leq p_\scrA$ for all~$c \in \scrC $.
    With similar reasoning we see that~$\xi(c) \leq p_\scrS$ for all~$c \in \scrC$
        in the second case.

Assume there is an equalizer~$e\colon E \to M_4$
    of~$\id_{M_4}$ and~$\ad_\sigma$ in~$(\Pure \op\vN)^{\mathsf{op}}$.
    By definition of purity~$e = c \after h$ for some filter~$c$ and comprehension~$h$.
As~$e$ is supposed to be an equalizer of~$\id$ and~$\ad_\sigma$,
    we have~$\ad_\sigma \after c \after h = \ad_\sigma \after e = e = c \after h$,
    hence~$\ad_\sigma \after c = c$ by the surjectivity of the corner~$h$.
We have just shown that then
        either~$c(b) \leq p_\scrA$ for all~$b$ or
        or~$c(b) \leq p_\scrS$ for all~$b$.
    Hence~$h(e) \leq p_\scrA$ for all~$e \in E$
        or~$h(e) \leq p_\scrS$ for all~$e \in E$.
    In the first case, consider the map~$\ad_{e^\dagger_\scrS} \colon M_3 \to M_4$.
    Clearly~$\ad_\sigma \after \ad_{e^\dagger_\scrS}
     = \ad_{e^\dagger_\scrS}$
     and so there must be a (unique) map~$h\colon M_3 \to E$
     with~$e \after h = \ad_{e^\dagger_\scrS}$.
    However~$p_\scrA \geq e(h(T)) = \ad_{e^\dagger_\scrS}(T) \leq p_\scrS$
        for all~$T \in M_3$
        and so~$\ad_{e^\dagger_\scrS} = 0$, quod non.
    We reach a contradiction in the other case in a similar fashion.
    Hence~$\Pure \op\vN$ has no coequalizers.
\end{solution}

% vim: se ft=tex.latex :

\oldchapter{Errata}\label{errata}

These are errata to the printed version of this thesis,
    which have been corrected in this PDF.
Additional errata will be published
    on the arXiv under~\href{https://arxiv.org/abs/1803.01911}{1803.01911}.

\begin{erratum}{physics-stinespring}%
The second part of the exercise starts with
    ``Conclude that any quantum channel~$\Phi\colon \scrB(\scrK) \to \scrB(\scrH)$ \ldots''.
This should have read
    ``Conclude that any quantum channel~$\Phi$ from~$\scrH$ to itself \ldots''
\end{erratum}
\begin{erratum}{ess-uniq-pur}%
The third-to-last sentence reads
``Derive from the latter
that for each~$y \in \scrK'$ and rank-one projector~$e \in \scrH$,
there is a~$y' \in \scrK'$
with~$U_0 (e \otimes y) = e \otimes y'$.''
This should have been
``Derive from the latter
that for each~$y \in \scrK'$ and unit vector~$e \in \scrH$,
there is a~$y' \in \scrK'$
with~$U_0 (e \otimes y) = e \otimes y'$.''.
\end{erratum}
\begin{erratum}{moved-dfn-selfdual}%
    A pre-Hilbert~$\scrB$-module~$X$ is \Define{self dual}
    if every bounded~$\scrB$-linear~$\tau \colon X \to \scrB$
    is of the form~$\tau(x) = \langle t, x\rangle$ for some~$t \in X$.
In print the condition that~$\tau$ should be bounded
        was accidentally lost when it was moved from~\sref{dils-selfdual}.
\end{erratum}
\begin{erratum}{dils-uniform-spaces-basics}%
The assumption  that the net~$x_\alpha$ is Cauchy in point 3
    is superfluous.
\end{erratum}
\begin{erratum}{hilbmod-adjoint-exists}%
The map~$T$ should be presumed to be bounded.
\end{erratum}
\begin{erratum}{hilmod-fixed-on-V}%
The C$^*$-algebra~$\scrB$ isn't assumed to be a von Neumann algebra,
    which it should have been (for otherwise~\sref{dils-completion} wouldn't
    be applicable.)
\end{erratum}
\begin{erratum}{err159IV}%
The displayed
    inequality
\begin{equation*}
    |f(T - p_S T p_S)| \ \leq \  |f(\, (1-p_S) T | +| p_S T (1-p_S) \,)|
\end{equation*}
    should have read
    $|f(T - p_S T p_S)| \ \leq \  |f(\, (1-p_S) T \,) | +| f(\, p_S T (1-p_S) \,)|$.
\end{erratum}
\begin{erratum}{dils-stand-filter}%
The standard filter~$c_b$ for ~$b \in \scrB$
    is a map~``$\ceil{b} \scrB \ceil{b} \to \scrB$''
    instead of a map~``$\scrB \to \ceil{b} \scrB \ceil{b}$''.
\end{erratum}
\begin{erratum}{quotient-basics}%
Not every zero map is a quotient; only those into~$0$.
\end{erratum}
\begin{erratum}{compr-basics}%
Not every zero map is a comprehension; only those from~$0$.
\end{erratum}
\begin{erratum}{asrt-pristine-reverse}%
Also in the last two points of the exercise,
    the map~$h$ should be presumed to be pristine.
\end{erratum}
\begin{erratum}{exc-purec-no-biproduct}%
The condition
    ``$f \after \varphi = \langle x, (\,\cdot\,) x\rangle$''
        is a typo.
        What was meant is~``$f (\varphi(T, c)) = \langle x, T x\rangle $
        for all~$c \in \scrC$
        and~$T \in \scrB(\scrH)$''.
The same mistake was made writing
    ``$\pi_1 \after \varphi = \langle x, (\,\cdot\,) x\rangle$'' and
    ``$\pi_2 \after \varphi = \langle y, (\,\cdot\,) y\rangle$''.
\end{erratum}

% vim: se ft=tex.latex :

% Trick arXiv into running latex 4 times
\message{Rerun to get citations correct.}

\end{document}